\documentclass[12pt]{article}
\usepackage{amsmath,amssymb, amsthm, amsopn, amsfonts,epsfig,graphics,relsize,exscale}{}
\usepackage{hyperref}

\newcommand{\field}[1]{\mathbb{#1}} \newcommand{\rz}{\field{R}}
\newcommand{\cz}{\field{C}} \newcommand{\nz}{\field{N}}
\newcommand{\zz}{\field{Z}} \newcommand{\tz}{\field{T}}

\newcommand{\sz}{\field{S}}

\newcommand{\Id}{{\textrm{Id}}}
\newcommand{\Imag}{{\textrm{Im\,}}}

\DeclareMathOperator\Real{\rz e\,}
\DeclareMathOperator\supp{\textrm{supp}}
\DeclareMathOperator\sign{\textrm{sign}\,}
\DeclareMathOperator\DiffP{\textrm{DiffP}}
\newcommand{\Hess}{{\text{Hess\,}}}
\newcommand{\Ker}{{\text{Ker\,}}}
\newcommand{\Ran}{{\text{Ran\,}}}
\newcommand{\ad}{{\text{ad\,}}}
\newcommand{\ccap}{\mathop{\cap}}
\newcommand{\ccup}{\mathop{\cup}}
\newtheorem{theorem}{Theorem}[section]
\newtheorem{lemma}[theorem]{Lemma}

\newtheorem{proposition}[theorem]{Proposition}
\newtheorem{definition}[theorem]{Definition}
\newtheorem{remark}[theorem]{Remark}
\newtheorem{corollary}[theorem]{Corollary}

\newtheorem{assumption}{Hypothesis}

\def\dim{\textrm{dim\;}}
\title{Boundary conditions and subelliptic estimates for geometric
  Kramers-Fokker-Planck operators on manifolds with boundaries
}
\author{
  F.~Nier
\thanks{IRMAR, UMR-CNRS 6625, Universit{\'e} de Rennes 1, Campus
    de Beaulieu, F35042 Rennes Cedex. \textit{email:
      francis.nier@univ-rennes1.fr} 
}
}
\begin{document}
%
%
%
%
\ifx\figforTeXisloaded\relax \else\global\let\figforTeXisloaded=\relax\fi
\message{version 1.8.4}
\catcode`\@=11
\ifx\ctr@ln@m\undefined\else%
    \immediate\write16{*** Fig4TeX WARNING : \string\ctr@ln@m\space already defined.}\fi
\def\ctr@ln@m#1{\ifx#1\undefined\else%
    \immediate\write16{*** Fig4TeX WARNING : \string#1 already defined.}\fi}
\ctr@ln@m\ctr@ld@f
\def\ctr@ld@f#1#2{\ctr@ln@m#2#1#2}
\ctr@ld@f\def\ctr@ln@w#1#2{\ctr@ln@m#2\csname#1\endcsname#2}
{\catcode`\/=0 \catcode`/\=12 /ctr@ld@f/gdef/BS@{\}}
\ctr@ld@f\def\ctr@lcsn@m#1{\expandafter\ifx\csname#1\endcsname\relax\else%
    \immediate\write16{*** Fig4TeX WARNING : \BS@\expandafter\string#1\space already defined.}\fi}
\ctr@ld@f\edef\colonc@tcode{\the\catcode`\:}
\ctr@ld@f\edef\semicolonc@tcode{\the\catcode`\;}
\ctr@ld@f\def\t@stc@tcodech@nge{{\let\c@tcodech@nged=\z@%
    \ifnum\colonc@tcode=\the\catcode`\:\else\let\c@tcodech@nged=\@ne\fi%
    \ifnum\semicolonc@tcode=\the\catcode`\;\else\let\c@tcodech@nged=\@ne\fi%
    \ifx\c@tcodech@nged\@ne%
    \immediate\write16{}
    \immediate\write16{!!!=============================================================!!!}
    \immediate\write16{ Fig4TeX WARNING :}
    \immediate\write16{ The category code of some characters has been changed, which will}
    \immediate\write16{ result in an error (message "Runaway argument?").}
    \immediate\write16{ This probably comes from another package that changed the category}
    \immediate\write16{ code after Fig4TeX was loaded. If that proves to be exact, the}
    \immediate\write16{ solution is to exchange the loading commands on top of your file}
    \immediate\write16{ so that Fig4TeX is loaded last. For example, in LaTeX, we should}
    \immediate\write16{ say :}
    \immediate\write16{\BS@ usepackage[french]{babel}}
    \immediate\write16{\BS@ usepackage{fig4tex}}
    \immediate\write16{!!!=============================================================!!!}
    \immediate\write16{}
    \fi}}
\ctr@ld@f\def\FigforTeX{F\kern-.05em i\kern-.05em g\kern-.1em\raise-.14em\hbox{4}\kern-.19em\TeX}
\ctr@ln@w{newdimen}\epsil@n\epsil@n=0.00005pt
\ctr@ln@w{newdimen}\Cepsil@n\Cepsil@n=0.005pt
\ctr@ln@w{newdimen}\dcq@\dcq@=254pt
\ctr@ln@w{newdimen}\PI@\PI@=3.141592pt
\ctr@ln@w{newdimen}\DemiPI@deg\DemiPI@deg=90pt
\ctr@ln@w{newdimen}\PI@deg\PI@deg=180pt
\ctr@ln@w{newdimen}\DePI@deg\DePI@deg=360pt
\ctr@ld@f\chardef\t@n=10
\ctr@ld@f\chardef\c@nt=100
\ctr@ld@f\chardef\@lxxiv=74
\ctr@ld@f\chardef\@xci=91
\ctr@ld@f\mathchardef\@nMnCQn=9949
\ctr@ld@f\chardef\@vi=6
\ctr@ld@f\chardef\@xxx=30
\ctr@ld@f\chardef\@lvi=56
\ctr@ld@f\chardef\@@lxxi=71
\ctr@ld@f\chardef\@lxxxv=85
\ctr@ld@f\mathchardef\@@mmmmlxviii=4068
\ctr@ld@f\mathchardef\@ccclx=360
\ctr@ld@f\mathchardef\@dccxx=720
\ctr@ln@w{newcount}\p@rtent \ctr@ln@w{newcount}\f@ctech \ctr@ln@w{newcount}\result@tent
\ctr@ln@w{newdimen}\v@lmin \ctr@ln@w{newdimen}\v@lmax \ctr@ln@w{newdimen}\v@leur
\ctr@ln@w{newdimen}\result@t\ctr@ln@w{newdimen}\result@@t
\ctr@ln@w{newdimen}\mili@u \ctr@ln@w{newdimen}\c@rre \ctr@ln@w{newdimen}\delt@
\ctr@ld@f\def\degT@rd{0.017453 }  
\ctr@ld@f\def\rdT@deg{57.295779 } 
\ctr@ln@m\v@leurseule
{\catcode`p=12 \catcode`t=12 \gdef\v@leurseule#1pt{#1}}
\ctr@ld@f\def\repdecn@mb#1{\expandafter\v@leurseule\the#1\space}
\ctr@ld@f\def\arct@n#1(#2,#3){{\v@lmin=#2\v@lmax=#3%
    \maxim@m{\mili@u}{-\v@lmin}{\v@lmin}\maxim@m{\c@rre}{-\v@lmax}{\v@lmax}%
    \delt@=\mili@u\m@ech\mili@u%
    \ifdim\c@rre>\@nMnCQn\mili@u\divide\v@lmax\tw@\c@lATAN\v@leur(\z@,\v@lmax)
    \else%
    \maxim@m{\mili@u}{-\v@lmin}{\v@lmin}\maxim@m{\c@rre}{-\v@lmax}{\v@lmax}%
    \m@ech\c@rre%
    \ifdim\mili@u>\@nMnCQn\c@rre\divide\v@lmin\tw@
    \maxim@m{\mili@u}{-\v@lmin}{\v@lmin}\c@lATAN\v@leur(\mili@u,\z@)%
    \else\c@lATAN\v@leur(\delt@,\v@lmax)\fi\fi%
    \ifdim\v@lmin<\z@\v@leur=-\v@leur\ifdim\v@lmax<\z@\advance\v@leur-\PI@%
    \else\advance\v@leur\PI@\fi\fi%
    \global\result@t=\v@leur}#1=\result@t}
\ctr@ld@f\def\m@ech#1{\ifdim#1>1.646pt\divide\mili@u\t@n\divide\c@rre\t@n\m@ech#1\fi}
\ctr@ld@f\def\c@lATAN#1(#2,#3){{\v@lmin=#2\v@lmax=#3\v@leur=\z@\delt@=\tw@ pt%
    \un@iter{0.785398}{\v@lmax<}%
    \un@iter{0.463648}{\v@lmax<}%
    \un@iter{0.244979}{\v@lmax<}%
    \un@iter{0.124355}{\v@lmax<}%
    \un@iter{0.062419}{\v@lmax<}%
    \un@iter{0.031240}{\v@lmax<}%
    \un@iter{0.015624}{\v@lmax<}%
    \un@iter{0.007812}{\v@lmax<}%
    \un@iter{0.003906}{\v@lmax<}%
    \un@iter{0.001953}{\v@lmax<}%
    \un@iter{0.000976}{\v@lmax<}%
    \un@iter{0.000488}{\v@lmax<}%
    \un@iter{0.000244}{\v@lmax<}%
    \un@iter{0.000122}{\v@lmax<}%
    \un@iter{0.000061}{\v@lmax<}%
    \un@iter{0.000030}{\v@lmax<}%
    \un@iter{0.000015}{\v@lmax<}%
    \global\result@t=\v@leur}#1=\result@t}
\ctr@ld@f\def\un@iter#1#2{%
    \divide\delt@\tw@\edef\dpmn@{\repdecn@mb{\delt@}}%
    \mili@u=\v@lmin%
    \ifdim#2\z@%
      \advance\v@lmin-\dpmn@\v@lmax\advance\v@lmax\dpmn@\mili@u%
      \advance\v@leur-#1pt%
    \else%
      \advance\v@lmin\dpmn@\v@lmax\advance\v@lmax-\dpmn@\mili@u%
      \advance\v@leur#1pt%
    \fi}
\ctr@ld@f\def\c@ssin#1#2#3{\expandafter\ifx\csname COS@\number#3\endcsname\relax\c@lCS{#3pt}%
    \expandafter\xdef\csname COS@\number#3\endcsname{\repdecn@mb\result@t}%
    \expandafter\xdef\csname SIN@\number#3\endcsname{\repdecn@mb\result@@t}\fi%
    \edef#1{\csname COS@\number#3\endcsname}\edef#2{\csname SIN@\number#3\endcsname}}
\ctr@ld@f\def\c@lCS#1{{\mili@u=#1\p@rtent=\@ne%
    \relax\ifdim\mili@u<\z@\red@ng<-\else\red@ng>+\fi\f@ctech=\p@rtent%
    \relax\ifdim\mili@u<\z@\mili@u=-\mili@u\f@ctech=-\f@ctech\fi\c@@lCS}}
\ctr@ld@f\def\c@@lCS{\v@lmin=\mili@u\c@rre=-\mili@u\advance\c@rre\DemiPI@deg\v@lmax=\c@rre%
    \mili@u\@@lxxi\mili@u\divide\mili@u\@@mmmmlxviii%
    \edef\v@larg{\repdecn@mb{\mili@u}}\mili@u=-\v@larg\mili@u%
    \edef\v@lmxde{\repdecn@mb{\mili@u}}%
    \c@rre\@@lxxi\c@rre\divide\c@rre\@@mmmmlxviii%
    \edef\v@largC{\repdecn@mb{\c@rre}}\c@rre=-\v@largC\c@rre%
    \edef\v@lmxdeC{\repdecn@mb{\c@rre}}%
    \fctc@s\mili@u\v@lmin\global\result@t\p@rtent\v@leur%
    \let\t@mp=\v@larg\let\v@larg=\v@largC\let\v@largC=\t@mp%
    \let\t@mp=\v@lmxde\let\v@lmxde=\v@lmxdeC\let\v@lmxdeC=\t@mp%
    \fctc@s\c@rre\v@lmax\global\result@@t\f@ctech\v@leur}
\ctr@ld@f\def\fctc@s#1#2{\v@leur=#1\relax\ifdim#2<\@lxxxv\p@\cosser@h\else\sinser@t\fi}
\ctr@ld@f\def\cosser@h{\advance\v@leur\@lvi\p@\divide\v@leur\@lvi%
    \v@leur=\v@lmxde\v@leur\advance\v@leur\@xxx\p@%
    \v@leur=\v@lmxde\v@leur\advance\v@leur\@ccclx\p@%
    \v@leur=\v@lmxde\v@leur\advance\v@leur\@dccxx\p@\divide\v@leur\@dccxx}
\ctr@ld@f\def\sinser@t{\v@leur=\v@lmxdeC\p@\advance\v@leur\@vi\p@%
    \v@leur=\v@largC\v@leur\divide\v@leur\@vi}
\ctr@ld@f\def\red@ng#1#2{\relax\ifdim\mili@u#1#2\DemiPI@deg\advance\mili@u#2-\PI@deg%
    \p@rtent=-\p@rtent\red@ng#1#2\fi}
\ctr@ld@f\def\pr@c@lCS#1#2#3{\ctr@lcsn@m{COS@\number#3 }%
    \expandafter\xdef\csname COS@\number#3\endcsname{#1}%
    \expandafter\xdef\csname SIN@\number#3\endcsname{#2}}
\pr@c@lCS{1}{0}{0}
\pr@c@lCS{0.7071}{0.7071}{45}\pr@c@lCS{0.7071}{-0.7071}{-45}
\pr@c@lCS{0}{1}{90}          \pr@c@lCS{0}{-1}{-90}
\pr@c@lCS{-1}{0}{180}        \pr@c@lCS{-1}{0}{-180}
\pr@c@lCS{0}{-1}{270}        \pr@c@lCS{0}{1}{-270}
\ctr@ld@f\def\invers@#1#2{{\v@leur=#2\maxim@m{\v@lmax}{-\v@leur}{\v@leur}%
    \f@ctech=\@ne\m@inv@rs%
    \multiply\v@leur\f@ctech\edef\v@lv@leur{\repdecn@mb{\v@leur}}%
    \p@rtentiere{\p@rtent}{\v@leur}\v@lmin=\p@\divide\v@lmin\p@rtent%
    \inv@rs@\multiply\v@lmax\f@ctech\global\result@t=\v@lmax}#1=\result@t}
\ctr@ld@f\def\m@inv@rs{\ifdim\v@lmax<\p@\multiply\v@lmax\t@n\multiply\f@ctech\t@n\m@inv@rs\fi}
\ctr@ld@f\def\inv@rs@{\v@lmax=-\v@lmin\v@lmax=\v@lv@leur\v@lmax%
    \advance\v@lmax\tw@ pt\v@lmax=\repdecn@mb{\v@lmin}\v@lmax%
    \delt@=\v@lmax\advance\delt@-\v@lmin\ifdim\delt@<\z@\delt@=-\delt@\fi%
    \ifdim\delt@>\epsil@n\v@lmin=\v@lmax\inv@rs@\fi}
\ctr@ld@f\def\minim@m#1#2#3{\relax\ifdim#2<#3#1=#2\else#1=#3\fi}
\ctr@ld@f\def\maxim@m#1#2#3{\relax\ifdim#2>#3#1=#2\else#1=#3\fi}
\ctr@ld@f\def\p@rtentiere#1#2{#1=#2\divide#1by65536 }
\ctr@ld@f\def\r@undint#1#2{{\v@leur=#2\divide\v@leur\t@n\p@rtentiere{\p@rtent}{\v@leur}%
    \v@leur=\p@rtent pt\global\result@t=\t@n\v@leur}#1=\result@t}
\ctr@ld@f\def\sqrt@#1#2{{\v@leur=#2%
    \minim@m{\v@lmin}{\p@}{\v@leur}\maxim@m{\v@lmax}{\p@}{\v@leur}%
    \f@ctech=\@ne\m@sqrt@\sqrt@@%
    \mili@u=\v@lmin\advance\mili@u\v@lmax\divide\mili@u\tw@\multiply\mili@u\f@ctech%
    \global\result@t=\mili@u}#1=\result@t}
\ctr@ld@f\def\m@sqrt@{\ifdim\v@leur>\dcq@\divide\v@leur\c@nt\v@lmax=\v@leur%
    \multiply\f@ctech\t@n\m@sqrt@\fi}
\ctr@ld@f\def\sqrt@@{\mili@u=\v@lmin\advance\mili@u\v@lmax\divide\mili@u\tw@%
    \c@rre=\repdecn@mb{\mili@u}\mili@u%
    \ifdim\c@rre<\v@leur\v@lmin=\mili@u\else\v@lmax=\mili@u\fi%
    \delt@=\v@lmax\advance\delt@-\v@lmin\ifdim\delt@>\epsil@n\sqrt@@\fi}
\ctr@ld@f\def\extrairelepremi@r#1\de#2{\expandafter\lepremi@r#2@#1#2}
\ctr@ld@f\def\lepremi@r#1,#2@#3#4{\def#3{#1}\def#4{#2}\ignorespaces}
\ctr@ld@f\def\@cfor#1:=#2\do#3{%
  \edef\@fortemp{#2}%
  \ifx\@fortemp\empty\else\@cforloop#2,\@nil,\@nil\@@#1{#3}\fi}
\ctr@ln@m\@nextwhile
\ctr@ld@f\def\@cforloop#1,#2\@@#3#4{%
  \def#3{#1}%
  \ifx#3\Fig@nnil\let\@nextwhile=\Fig@fornoop\else#4\relax\let\@nextwhile=\@cforloop\fi%
  \@nextwhile#2\@@#3{#4}}

\ctr@ld@f\def\@ecfor#1:=#2\do#3{%
  \def\@@cfor{\@cfor#1:=}%
  \edef\@@@cfor{#2}%
  \expandafter\@@cfor\@@@cfor\do{#3}}
\ctr@ld@f\def\Fig@nnil{\@nil}
\ctr@ld@f\def\Fig@fornoop#1\@@#2#3{}
\ctr@ln@m\list@@rg
\ctr@ld@f\def\trtlis@rg#1#2{\def\list@@rg{#1}%
    \@ecfor\p@rv@l:=\list@@rg\do{\expandafter#2\p@rv@l|}}
\ctr@ln@w{newbox}\b@xvisu
\ctr@ln@w{newtoks}\let@xte
\ctr@ln@w{newif}\ifitis@K
\ctr@ln@w{newcount}\s@mme
\ctr@ln@w{newcount}\l@mbd@un \ctr@ln@w{newcount}\l@mbd@de
\ctr@ln@w{newcount}\superc@ntr@l\superc@ntr@l=\@ne        
\ctr@ln@w{newcount}\typec@ntr@l\typec@ntr@l=\superc@ntr@l 
\ctr@ln@w{newdimen}\v@lX  \ctr@ln@w{newdimen}\v@lY  \ctr@ln@w{newdimen}\v@lZ
\ctr@ln@w{newdimen}\v@lXa \ctr@ln@w{newdimen}\v@lYa \ctr@ln@w{newdimen}\v@lZa
\ctr@ln@w{newdimen}\unit@\unit@=\p@ 
\ctr@ld@f\def\unit@util{pt}
\ctr@ld@f\def\ptT@ptps{0.996264 }
\ctr@ld@f\def\ptpsT@pt{1.00375 }
\ctr@ld@f\def\ptT@unit@{1} 
\ctr@ld@f\def\setunit@#1{\def\unit@util{#1}\setunit@@#1:\invers@{\result@t}{\unit@}%
    \edef\ptT@unit@{\repdecn@mb\result@t}}
\ctr@ld@f\def\setunit@@#1#2:{\ifcat#1a\unit@=\@ne#1#2\else\unit@=#1#2\fi}
\ctr@ld@f\def\d@fm@cdim#1#2{{\v@leur=#2\v@leur=\ptT@unit@\v@leur\xdef#1{\repdecn@mb\v@leur}}}
\ctr@ln@w{newif}\ifBdingB@x\BdingB@xtrue
\ctr@ln@w{newdimen}\c@@rdXmin \ctr@ln@w{newdimen}\c@@rdYmin  
\ctr@ln@w{newdimen}\c@@rdXmax \ctr@ln@w{newdimen}\c@@rdYmax
\ctr@ld@f\def\b@undb@x#1#2{\ifBdingB@x%
    \relax\ifdim#1<\c@@rdXmin\global\c@@rdXmin=#1\fi%
    \relax\ifdim#2<\c@@rdYmin\global\c@@rdYmin=#2\fi%
    \relax\ifdim#1>\c@@rdXmax\global\c@@rdXmax=#1\fi%
    \relax\ifdim#2>\c@@rdYmax\global\c@@rdYmax=#2\fi\fi}
\ctr@ld@f\def\b@undb@xP#1{{\Figg@tXY{#1}\b@undb@x{\v@lX}{\v@lY}}}
\ctr@ld@f\def\ellBB@x#1;#2,#3(#4,#5,#6){{\s@uvc@ntr@l\et@tellBB@x%
    \setc@ntr@l{2}\figptell-2::#1;#2,#3(#4,#6)\b@undb@xP{-2}%
    \figptell-2::#1;#2,#3(#5,#6)\b@undb@xP{-2}%
    \c@ssin{\C@}{\S@}{#6}\v@lmin=\C@ pt\v@lmax=\S@ pt%
    \mili@u=#3\v@lmin\delt@=#2\v@lmax\arct@n\v@leur(\delt@,\mili@u)%
    \mili@u=-#3\v@lmax\delt@=#2\v@lmin\arct@n\c@rre(\delt@,\mili@u)%
    \v@leur=\rdT@deg\v@leur\advance\v@leur-\DePI@deg%
    \c@rre=\rdT@deg\c@rre\advance\c@rre-\DePI@deg%
    \v@lmin=#4pt\v@lmax=#5pt%
    \loop\ifdim\v@leur<\v@lmax\ifdim\v@leur>\v@lmin%
    \edef\@ngle{\repdecn@mb\v@leur}\figptell-2::#1;#2,#3(\@ngle,#6)%
    \b@undb@xP{-2}\fi\advance\v@leur\PI@deg\repeat%
    \loop\ifdim\c@rre<\v@lmax\ifdim\c@rre>\v@lmin%
    \edef\@ngle{\repdecn@mb\c@rre}\figptell-2::#1;#2,#3(\@ngle,#6)%
    \b@undb@xP{-2}\fi\advance\c@rre\PI@deg\repeat%
    \resetc@ntr@l\et@tellBB@x}\ignorespaces}
\ctr@ld@f\def\initb@undb@x{\c@@rdXmin=\maxdimen\c@@rdYmin=\maxdimen%
    \c@@rdXmax=-\maxdimen\c@@rdYmax=-\maxdimen}
\ctr@ld@f\def\c@ntr@lnum#1{%
    \relax\ifnum\typec@ntr@l=\@ne%
    \ifnum#1<\z@%
    \immediate\write16{*** Forbidden point number (#1). Abort.}\end\fi\fi%
    \set@bjc@de{#1}}
\ctr@ln@m\objc@de
\ctr@ld@f\def\set@bjc@de#1{\edef\objc@de{@BJ\ifnum#1<\z@ M\romannumeral-#1\else\romannumeral#1\fi}}
\s@mme=\m@ne\loop\ifnum\s@mme>-19
  \set@bjc@de{\s@mme}\ctr@lcsn@m\objc@de\ctr@lcsn@m{\objc@de T}
\advance\s@mme\m@ne\repeat
\s@mme=\@ne\loop\ifnum\s@mme<6
  \set@bjc@de{\s@mme}\ctr@lcsn@m\objc@de\ctr@lcsn@m{\objc@de T}
\advance\s@mme\@ne\repeat
\ctr@ld@f\def\setc@ntr@l#1{\ifnum\superc@ntr@l>#1\typec@ntr@l=\superc@ntr@l%
    \else\typec@ntr@l=#1\fi}
\ctr@ld@f\def\resetc@ntr@l#1{\global\superc@ntr@l=#1\setc@ntr@l{#1}}
\ctr@ld@f\def\s@uvc@ntr@l#1{\edef#1{\the\superc@ntr@l}}
\ctr@ln@m\c@lproscal
\ctr@ld@f\def\c@lproscalDD#1[#2,#3]{{\Figg@tXY{#2}%
    \edef\Xu@{\repdecn@mb{\v@lX}}\edef\Yu@{\repdecn@mb{\v@lY}}\Figg@tXY{#3}%
    \global\result@t=\Xu@\v@lX\global\advance\result@t\Yu@\v@lY}#1=\result@t}
\ctr@ld@f\def\c@lproscalTD#1[#2,#3]{{\Figg@tXY{#2}\edef\Xu@{\repdecn@mb{\v@lX}}%
    \edef\Yu@{\repdecn@mb{\v@lY}}\edef\Zu@{\repdecn@mb{\v@lZ}}%
    \Figg@tXY{#3}\global\result@t=\Xu@\v@lX\global\advance\result@t\Yu@\v@lY%
    \global\advance\result@t\Zu@\v@lZ}#1=\result@t}
\ctr@ld@f\def\c@lprovec#1{%
    \det@rmC\v@lZa(\v@lX,\v@lY,\v@lmin,\v@lmax)%
    \det@rmC\v@lXa(\v@lY,\v@lZ,\v@lmax,\v@leur)%
    \det@rmC\v@lYa(\v@lZ,\v@lX,\v@leur,\v@lmin)%
    \Figv@ctCreg#1(\v@lXa,\v@lYa,\v@lZa)}
\ctr@ld@f\def\det@rm#1[#2,#3]{{\Figg@tXY{#2}\Figg@tXYa{#3}%
    \delt@=\repdecn@mb{\v@lX}\v@lYa\advance\delt@-\repdecn@mb{\v@lY}\v@lXa%
    \global\result@t=\delt@}#1=\result@t}
\ctr@ld@f\def\det@rmC#1(#2,#3,#4,#5){{\global\result@t=\repdecn@mb{#2}#5%
    \global\advance\result@t-\repdecn@mb{#3}#4}#1=\result@t}
\ctr@ld@f\def\getredf@ctDD#1(#2,#3){{\maxim@m{\v@lXa}{-#2}{#2}\maxim@m{\v@lYa}{-#3}{#3}%
    \maxim@m{\v@lXa}{\v@lXa}{\v@lYa}
    \ifdim\v@lXa>\@xci pt\divide\v@lXa\@xci%
    \p@rtentiere{\p@rtent}{\v@lXa}\advance\p@rtent\@ne\else\p@rtent=\@ne\fi%
    \global\result@tent=\p@rtent}#1=\result@tent\ignorespaces}
\ctr@ld@f\def\getredf@ctTD#1(#2,#3,#4){{\maxim@m{\v@lXa}{-#2}{#2}\maxim@m{\v@lYa}{-#3}{#3}%
    \maxim@m{\v@lZa}{-#4}{#4}\maxim@m{\v@lXa}{\v@lXa}{\v@lYa}%
    \maxim@m{\v@lXa}{\v@lXa}{\v@lZa}
    \ifdim\v@lXa>\@lxxiv pt\divide\v@lXa\@lxxiv%
    \p@rtentiere{\p@rtent}{\v@lXa}\advance\p@rtent\@ne\else\p@rtent=\@ne\fi%
    \global\result@tent=\p@rtent}#1=\result@tent\ignorespaces}
\ctr@ld@f\def\FigptintercircB@zDD#1:#2:#3,#4[#5,#6,#7,#8]{{\s@uvc@ntr@l\et@tfigptintercircB@zDD%
    \setc@ntr@l{2}\figvectPDD-1[#5,#8]\Figg@tXY{-1}\getredf@ctDD\f@ctech(\v@lX,\v@lY)%
    \mili@u=#4\unit@\divide\mili@u\f@ctech\c@rre=\repdecn@mb{\mili@u}\mili@u%
    \figptBezierDD-5::#3[#5,#6,#7,#8]%
    \v@lmin=#3\p@\v@lmax=\v@lmin\advance\v@lmax0.1\p@%
    \loop\edef\T@{\repdecn@mb{\v@lmax}}\figptBezierDD-2::\T@[#5,#6,#7,#8]%
    \figvectPDD-1[-5,-2]\n@rmeucCDD{\delt@}{-1}\ifdim\delt@<\c@rre\v@lmin=\v@lmax%
    \advance\v@lmax0.1\p@\repeat%
    \loop\mili@u=\v@lmin\advance\mili@u\v@lmax%
    \divide\mili@u\tw@\edef\T@{\repdecn@mb{\mili@u}}\figptBezierDD-2::\T@[#5,#6,#7,#8]%
    \figvectPDD-1[-5,-2]\n@rmeucCDD{\delt@}{-1}\ifdim\delt@>\c@rre\v@lmax=\mili@u%
    \else\v@lmin=\mili@u\fi\v@leur=\v@lmax\advance\v@leur-\v@lmin%
    \ifdim\v@leur>\epsil@n\repeat\figptcopyDD#1:#2/-2/%
    \resetc@ntr@l\et@tfigptintercircB@zDD}\ignorespaces}
\ctr@ln@m\figptinterlines
\ctr@ld@f\def\inters@cDD#1:#2[#3,#4;#5,#6]{{\s@uvc@ntr@l\et@tinters@cDD%
    \setc@ntr@l{2}\vecunit@{-1}{#4}\vecunit@{-2}{#6}%
    \Figg@tXY{-1}\setc@ntr@l{1}\Figg@tXYa{#3}%
    \edef\A@{\repdecn@mb{\v@lX}}\edef\B@{\repdecn@mb{\v@lY}}%
    \v@lmin=\B@\v@lXa\advance\v@lmin-\A@\v@lYa%
    \Figg@tXYa{#5}\setc@ntr@l{2}\Figg@tXY{-2}%
    \edef\C@{\repdecn@mb{\v@lX}}\edef\D@{\repdecn@mb{\v@lY}}%
    \v@lmax=\D@\v@lXa\advance\v@lmax-\C@\v@lYa%
    \delt@=\A@\v@lY\advance\delt@-\B@\v@lX%
    \invers@{\v@leur}{\delt@}\edef\v@ldelta{\repdecn@mb{\v@leur}}%
    \v@lXa=\A@\v@lmax\advance\v@lXa-\C@\v@lmin%
    \v@lYa=\B@\v@lmax\advance\v@lYa-\D@\v@lmin%
    \v@lXa=\v@ldelta\v@lXa\v@lYa=\v@ldelta\v@lYa%
    \setc@ntr@l{1}\Figp@intregDD#1:{#2}(\v@lXa,\v@lYa)%
    \resetc@ntr@l\et@tinters@cDD}\ignorespaces}
\ctr@ld@f\def\inters@cTD#1:#2[#3,#4;#5,#6]{{\s@uvc@ntr@l\et@tinters@cTD%
    \setc@ntr@l{2}\figvectNVTD-1[#4,#6]\figvectNVTD-2[#6,-1]\figvectPTD-1[#3,#5]%
    \r@pPSTD\v@leur[-2,-1,#4]\edef\v@lcoef{\repdecn@mb{\v@leur}}%
    \figpttraTD#1:{#2}=#3/\v@lcoef,#4/\resetc@ntr@l\et@tinters@cTD}\ignorespaces}
\ctr@ld@f\def\r@pPSTD#1[#2,#3,#4]{{\Figg@tXY{#2}\edef\Xu@{\repdecn@mb{\v@lX}}%
    \edef\Yu@{\repdecn@mb{\v@lY}}\edef\Zu@{\repdecn@mb{\v@lZ}}%
    \Figg@tXY{#3}\v@lmin=\Xu@\v@lX\advance\v@lmin\Yu@\v@lY\advance\v@lmin\Zu@\v@lZ%
    \Figg@tXY{#4}\v@lmax=\Xu@\v@lX\advance\v@lmax\Yu@\v@lY\advance\v@lmax\Zu@\v@lZ%
    \invers@{\v@leur}{\v@lmax}\global\result@t=\repdecn@mb{\v@leur}\v@lmin}%
    #1=\result@t}
\ctr@ln@m\n@rminf
\ctr@ld@f\def\n@rminfDD#1#2{{\Figg@tXY{#2}\maxim@m{\v@lX}{\v@lX}{-\v@lX}%
    \maxim@m{\v@lY}{\v@lY}{-\v@lY}\maxim@m{\global\result@t}{\v@lX}{\v@lY}}%
    #1=\result@t}
\ctr@ld@f\def\n@rminfTD#1#2{{\Figg@tXY{#2}\maxim@m{\v@lX}{\v@lX}{-\v@lX}%
    \maxim@m{\v@lY}{\v@lY}{-\v@lY}\maxim@m{\v@lZ}{\v@lZ}{-\v@lZ}%
    \maxim@m{\v@lX}{\v@lX}{\v@lY}\maxim@m{\global\result@t}{\v@lX}{\v@lZ}}%
    #1=\result@t}
\ctr@ld@f\def\n@rmeucCDD#1#2{\Figg@tXY{#2}\divide\v@lX\f@ctech\divide\v@lY\f@ctech%
    #1=\repdecn@mb{\v@lX}\v@lX\v@lX=\repdecn@mb{\v@lY}\v@lY\advance#1\v@lX}
\ctr@ld@f\def\n@rmeucCTD#1#2{\Figg@tXY{#2}%
    \divide\v@lX\f@ctech\divide\v@lY\f@ctech\divide\v@lZ\f@ctech%
    #1=\repdecn@mb{\v@lX}\v@lX\v@lX=\repdecn@mb{\v@lY}\v@lY\advance#1\v@lX%
    \v@lX=\repdecn@mb{\v@lZ}\v@lZ\advance#1\v@lX}
\ctr@ln@m\n@rmeucSV
\ctr@ld@f\def\n@rmeucSVDD#1#2{{\Figg@tXY{#2}%
    \v@lXa=\repdecn@mb{\v@lX}\v@lX\v@lYa=\repdecn@mb{\v@lY}\v@lY%
    \advance\v@lXa\v@lYa\sqrt@{\global\result@t}{\v@lXa}}#1=\result@t}
\ctr@ld@f\def\n@rmeucSVTD#1#2{{\Figg@tXY{#2}\v@lXa=\repdecn@mb{\v@lX}\v@lX%
    \v@lYa=\repdecn@mb{\v@lY}\v@lY\v@lZa=\repdecn@mb{\v@lZ}\v@lZ%
    \advance\v@lXa\v@lYa\advance\v@lXa\v@lZa\sqrt@{\global\result@t}{\v@lXa}}#1=\result@t}
\ctr@ln@m\n@rmeuc
\ctr@ld@f\def\n@rmeucDD#1#2{{\Figg@tXY{#2}\getredf@ctDD\f@ctech(\v@lX,\v@lY)%
    \divide\v@lX\f@ctech\divide\v@lY\f@ctech%
    \v@lXa=\repdecn@mb{\v@lX}\v@lX\v@lYa=\repdecn@mb{\v@lY}\v@lY%
    \advance\v@lXa\v@lYa\sqrt@{\global\result@t}{\v@lXa}%
    \global\multiply\result@t\f@ctech}#1=\result@t}
\ctr@ld@f\def\n@rmeucTD#1#2{{\Figg@tXY{#2}\getredf@ctTD\f@ctech(\v@lX,\v@lY,\v@lZ)%
    \divide\v@lX\f@ctech\divide\v@lY\f@ctech\divide\v@lZ\f@ctech%
    \v@lXa=\repdecn@mb{\v@lX}\v@lX%
    \v@lYa=\repdecn@mb{\v@lY}\v@lY\v@lZa=\repdecn@mb{\v@lZ}\v@lZ%
    \advance\v@lXa\v@lYa\advance\v@lXa\v@lZa\sqrt@{\global\result@t}{\v@lXa}%
    \global\multiply\result@t\f@ctech}#1=\result@t}
\ctr@ln@m\vecunit@
\ctr@ld@f\def\vecunit@DD#1#2{{\Figg@tXY{#2}\getredf@ctDD\f@ctech(\v@lX,\v@lY)%
    \divide\v@lX\f@ctech\divide\v@lY\f@ctech%
    \Figv@ctCreg#1(\v@lX,\v@lY)\n@rmeucSV{\v@lYa}{#1}%
    \invers@{\v@lXa}{\v@lYa}\edef\v@lv@lXa{\repdecn@mb{\v@lXa}}%
    \v@lX=\v@lv@lXa\v@lX\v@lY=\v@lv@lXa\v@lY%
    \Figv@ctCreg#1(\v@lX,\v@lY)\multiply\v@lYa\f@ctech\global\result@t=\v@lYa}}
\ctr@ld@f\def\vecunit@TD#1#2{{\Figg@tXY{#2}\getredf@ctTD\f@ctech(\v@lX,\v@lY,\v@lZ)%
    \divide\v@lX\f@ctech\divide\v@lY\f@ctech\divide\v@lZ\f@ctech%
    \Figv@ctCreg#1(\v@lX,\v@lY,\v@lZ)\n@rmeucSV{\v@lYa}{#1}%
    \invers@{\v@lXa}{\v@lYa}\edef\v@lv@lXa{\repdecn@mb{\v@lXa}}%
    \v@lX=\v@lv@lXa\v@lX\v@lY=\v@lv@lXa\v@lY\v@lZ=\v@lv@lXa\v@lZ%
    \Figv@ctCreg#1(\v@lX,\v@lY,\v@lZ)\multiply\v@lYa\f@ctech\global\result@t=\v@lYa}}
\ctr@ld@f\def\vecunitC@TD[#1,#2]{\Figg@tXYa{#1}\Figg@tXY{#2}%
    \advance\v@lX-\v@lXa\advance\v@lY-\v@lYa\advance\v@lZ-\v@lZa\c@lvecunitTD}
\ctr@ld@f\def\vecunitCV@TD#1{\Figg@tXY{#1}\c@lvecunitTD}
\ctr@ld@f\def\c@lvecunitTD{\getredf@ctTD\f@ctech(\v@lX,\v@lY,\v@lZ)%
    \divide\v@lX\f@ctech\divide\v@lY\f@ctech\divide\v@lZ\f@ctech%
    \v@lXa=\repdecn@mb{\v@lX}\v@lX%
    \v@lYa=\repdecn@mb{\v@lY}\v@lY\v@lZa=\repdecn@mb{\v@lZ}\v@lZ%
    \advance\v@lXa\v@lYa\advance\v@lXa\v@lZa\sqrt@{\v@lYa}{\v@lXa}%
    \invers@{\v@lXa}{\v@lYa}\edef\v@lv@lXa{\repdecn@mb{\v@lXa}}%
    \v@lX=\v@lv@lXa\v@lX\v@lY=\v@lv@lXa\v@lY\v@lZ=\v@lv@lXa\v@lZ}
\ctr@ln@m\figgetangle
\ctr@ld@f\def\figgetangleDD#1[#2,#3,#4]{\ifps@cri{\s@uvc@ntr@l\et@tfiggetangleDD\setc@ntr@l{2}%
    \figvectPDD-1[#2,#3]\figvectPDD-2[#2,#4]\vecunit@{-1}{-1}%
    \c@lproscalDD\delt@[-2,-1]\figvectNVDD-1[-1]\c@lproscalDD\v@leur[-2,-1]%
    \arct@n\v@lmax(\delt@,\v@leur)\v@lmax=\rdT@deg\v@lmax%
    \ifdim\v@lmax<\z@\advance\v@lmax\DePI@deg\fi\xdef#1{\repdecn@mb{\v@lmax}}%
    \resetc@ntr@l\et@tfiggetangleDD}\ignorespaces\fi}
\ctr@ld@f\def\figgetangleTD#1[#2,#3,#4,#5]{\ifps@cri{\s@uvc@ntr@l\et@tfiggetangleTD\setc@ntr@l{2}%
    \figvectPTD-1[#2,#3]\figvectPTD-2[#2,#5]\figvectNVTD-3[-1,-2]%
    \figvectPTD-2[#2,#4]\figvectNVTD-4[-3,-1]%
    \vecunit@{-1}{-1}\c@lproscalTD\delt@[-2,-1]\c@lproscalTD\v@leur[-2,-4]%
    \arct@n\v@lmax(\delt@,\v@leur)\v@lmax=\rdT@deg\v@lmax%
    \ifdim\v@lmax<\z@\advance\v@lmax\DePI@deg\fi\xdef#1{\repdecn@mb{\v@lmax}}%
    \resetc@ntr@l\et@tfiggetangleTD}\ignorespaces\fi}
\ctr@ld@f\def\figgetdist#1[#2,#3]{\ifps@cri{\s@uvc@ntr@l\et@tfiggetdist\setc@ntr@l{2}%
    \figvectP-1[#2,#3]\n@rmeuc{\v@lX}{-1}\v@lX=\ptT@unit@\v@lX\xdef#1{\repdecn@mb{\v@lX}}%
    \resetc@ntr@l\et@tfiggetdist}\ignorespaces\fi}
\ctr@ld@f\def\Figg@tT#1{\c@ntr@lnum{#1}%
    {\expandafter\expandafter\expandafter\extr@ctT\csname\objc@de\endcsname:%
     \ifnum\B@@ltxt=\z@\ptn@me{#1}\else\csname\objc@de T\endcsname\fi}}
\ctr@ld@f\def\extr@ctT#1,#2,#3/#4:{\def\B@@ltxt{#3}}
\ctr@ld@f\def\Figg@tXY#1{\c@ntr@lnum{#1}%
    \expandafter\expandafter\expandafter\extr@ctC\csname\objc@de\endcsname:}
\ctr@ln@m\extr@ctC
\ctr@ld@f\def\extr@ctCDD#1/#2,#3,#4:{\v@lX=#2\v@lY=#3}
\ctr@ld@f\def\extr@ctCTD#1/#2,#3,#4:{\v@lX=#2\v@lY=#3\v@lZ=#4}
\ctr@ld@f\def\Figg@tXYa#1{\c@ntr@lnum{#1}%
    \expandafter\expandafter\expandafter\extr@ctCa\csname\objc@de\endcsname:}
\ctr@ln@m\extr@ctCa
\ctr@ld@f\def\extr@ctCaDD#1/#2,#3,#4:{\v@lXa=#2\v@lYa=#3}
\ctr@ld@f\def\extr@ctCaTD#1/#2,#3,#4:{\v@lXa=#2\v@lYa=#3\v@lZa=#4}
\ctr@ln@m\t@xt@
\ctr@ld@f\def\figinit#1{\t@stc@tcodech@nge\initpr@lim\Figinit@#1,:\initpss@ttings\ignorespaces}
\ctr@ld@f\def\Figinit@#1,#2:{\setunit@{#1}\def\t@xt@{#2}\ifx\t@xt@\empty\else\Figinit@@#2:\fi}
\ctr@ld@f\def\Figinit@@#1#2:{\if#12 \else\Figs@tproj{#1}\initTD@\fi}
\ctr@ln@w{newif}\ifTr@isDim
\ctr@ld@f\def\UnD@fined{UNDEFINED}
\ctr@ld@f\def\ifundefined#1{\expandafter\ifx\csname#1\endcsname\relax}
\ctr@ln@m\@utoFN
\ctr@ln@m\@utoFInDone
\ctr@ln@m\disob@unit
\ctr@ld@f\def\initpr@lim{\initb@undb@x\figsetmark{}\figsetptname{$A_{##1}$}\def\Sc@leFact{1}%
    \initDD@\figsetroundcoord{yes}\ps@critrue\expandafter\setupd@te\defaultupdate:%
    \edef\disob@unit{\UnD@fined}\edef\t@rgetpt{\UnD@fined}\gdef\@utoFInDone{1}\gdef\@utoFN{0}}
\ctr@ld@f\def\initDD@{\Tr@isDimfalse%
    \ifPDFm@ke%
     \let\Ps@rcerc=\Ps@rcercBz%
     \let\Ps@rell=\Ps@rellBz%
    \fi
    \let\c@lDCUn=\c@lDCUnDD%
    \let\c@lDCDeux=\c@lDCDeuxDD%
    \let\c@ldefproj=\relax%
    \let\c@lproscal=\c@lproscalDD%
    \let\c@lprojSP=\relax%
    \let\extr@ctC=\extr@ctCDD%
    \let\extr@ctCa=\extr@ctCaDD%
    \let\extr@ctCF=\extr@ctCFDD%
    \let\Figp@intreg=\Figp@intregDD%
    \let\Figpts@xes=\Figpts@xesDD%
    \let\n@rmeucSV=\n@rmeucSVDD\let\n@rmeuc=\n@rmeucDD\let\n@rminf=\n@rminfDD%
    \let\pr@dMatV=\pr@dMatVDD%
    \let\ps@xes=\ps@xesDD%
    \let\vecunit@=\vecunit@DD%
    \let\figcoord=\figcoordDD%
    \let\figgetangle=\figgetangleDD%
    \let\figpt=\figptDD%
    \let\figptBezier=\figptBezierDD%
    \let\figptbary=\figptbaryDD%
    \let\figptcirc=\figptcircDD%
    \let\figptcircumcenter=\figptcircumcenterDD%
    \let\figptcopy=\figptcopyDD%
    \let\figptcurvcenter=\figptcurvcenterDD%
    \let\figptell=\figptellDD%
    \let\figptendnormal=\figptendnormalDD%
    \let\figptinterlineplane=\figptinterlineplaneDD%
    \let\figptinterlines=\inters@cDD%
    \let\figptorthocenter=\figptorthocenterDD%
    \let\figptorthoprojline=\figptorthoprojlineDD%
    \let\figptorthoprojplane=\figptorthoprojplaneDD%
    \let\figptrot=\figptrotDD%
    \let\figptscontrol=\figptscontrolDD%
    \let\figptsintercirc=\figptsintercircDD%
    \let\figptsinterlinell=\figptsinterlinellDD%
    \let\figptsorthoprojline=\figptsorthoprojlineDD%
    \let\figptorthoprojplane=\figptorthoprojplaneDD%
    \let\figptsrot=\figptsrotDD%
    \let\figptssym=\figptssymDD%
    \let\figptstra=\figptstraDD%
    \let\figptsym=\figptsymDD%
    \let\figpttraC=\figpttraCDD%
    \let\figpttra=\figpttraDD%
    \let\figptvisilimSL=\figptvisilimSLDD%
    \let\figsetobdist=\figsetobdistDD%
    \let\figsettarget=\figsettargetDD%
    \let\figsetview=\figsetviewDD%
    \let\figvectDBezier=\figvectDBezierDD%
    \let\figvectN=\figvectNDD%
    \let\figvectNV=\figvectNVDD%
    \let\figvectP=\figvectPDD%
    \let\figvectU=\figvectUDD%
    \let\psarccircP=\psarccircPDD%
    \let\psarccirc=\psarccircDD%
    \let\psarcell=\psarcellDD%
    \let\psarcellPA=\psarcellPADD%
    \let\psarrowBezier=\psarrowBezierDD%
    \let\psarrowcircP=\psarrowcircPDD%
    \let\psarrowcirc=\psarrowcircDD%
    \let\psarrowhead=\psarrowheadDD%
    \let\psarrow=\psarrowDD%
    \let\psBezier=\psBezierDD%
    \let\pscirc=\pscircDD%
    \let\pscurve=\pscurveDD%
    \let\psnormal=\psnormalDD%
    }
\ctr@ld@f\def\initTD@{\Tr@isDimtrue\initb@undb@xTD\newt@rgetptfalse\newdis@bfalse%
    \let\c@lDCUn=\c@lDCUnTD%
    \let\c@lDCDeux=\c@lDCDeuxTD%
    \let\c@ldefproj=\c@ldefprojTD%
    \let\c@lproscal=\c@lproscalTD%
    \let\extr@ctC=\extr@ctCTD%
    \let\extr@ctCa=\extr@ctCaTD%
    \let\extr@ctCF=\extr@ctCFTD%
    \let\Figp@intreg=\Figp@intregTD%
    \let\Figpts@xes=\Figpts@xesTD%
    \let\n@rmeucSV=\n@rmeucSVTD\let\n@rmeuc=\n@rmeucTD\let\n@rminf=\n@rminfTD%
    \let\pr@dMatV=\pr@dMatVTD%
    \let\ps@xes=\ps@xesTD%
    \let\vecunit@=\vecunit@TD%
    \let\figcoord=\figcoordTD%
    \let\figgetangle=\figgetangleTD%
    \let\figpt=\figptTD%
    \let\figptBezier=\figptBezierTD%
    \let\figptbary=\figptbaryTD%
    \let\figptcirc=\figptcircTD%
    \let\figptcircumcenter=\figptcircumcenterTD%
    \let\figptcopy=\figptcopyTD%
    \let\figptcurvcenter=\figptcurvcenterTD%
    \let\figptinterlineplane=\figptinterlineplaneTD%
    \let\figptinterlines=\inters@cTD%
    \let\figptorthocenter=\figptorthocenterTD%
    \let\figptorthoprojline=\figptorthoprojlineTD%
    \let\figptorthoprojplane=\figptorthoprojplaneTD%
    \let\figptrot=\figptrotTD%
    \let\figptscontrol=\figptscontrolTD%
    \let\figptsintercirc=\figptsintercircTD%
    \let\figptsorthoprojline=\figptsorthoprojlineTD%
    \let\figptsorthoprojplane=\figptsorthoprojplaneTD%
    \let\figptsrot=\figptsrotTD%
    \let\figptssym=\figptssymTD%
    \let\figptstra=\figptstraTD%
    \let\figptsym=\figptsymTD%
    \let\figpttraC=\figpttraCTD%
    \let\figpttra=\figpttraTD%
    \let\figptvisilimSL=\figptvisilimSLTD%
    \let\figsetobdist=\figsetobdistTD%
    \let\figsettarget=\figsettargetTD%
    \let\figsetview=\figsetviewTD%
    \let\figvectDBezier=\figvectDBezierTD%
    \let\figvectN=\figvectNTD%
    \let\figvectNV=\figvectNVTD%
    \let\figvectP=\figvectPTD%
    \let\figvectU=\figvectUTD%
    \let\psarccircP=\psarccircPTD%
    \let\psarccirc=\psarccircTD%
    \let\psarcell=\psarcellTD%
    \let\psarcellPA=\psarcellPATD%
    \let\psarrowBezier=\psarrowBezierTD%
    \let\psarrowcircP=\psarrowcircPTD%
    \let\psarrowcirc=\psarrowcircTD%
    \let\psarrowhead=\psarrowheadTD%
    \let\psarrow=\psarrowTD%
    \let\psBezier=\psBezierTD%
    \let\pscirc=\pscircTD%
    \let\pscurve=\pscurveTD%
    }
\ctr@ld@f\def\un@v@ilable#1{\immediate\write16{*** The macro #1 is not available in the current context.}}
\ctr@ld@f\def\figinsert#1{{\def\t@xt@{#1}\relax%
    \ifx\t@xt@\empty\ifnum\@utoFInDone>\z@\Figinsert@\DefGIfilen@me,:\fi%
    \else\expandafter\FiginsertNu@#1 :\fi}\ignorespaces}
\ctr@ld@f\def\FiginsertNu@#1 #2:{\def\t@xt@{#1}\relax\ifx\t@xt@\empty\def\t@xt@{#2}%
    \ifx\t@xt@\empty\ifnum\@utoFInDone>\z@\Figinsert@\DefGIfilen@me,:\fi%
    \else\FiginsertNu@#2:\fi\else\expandafter\FiginsertNd@#1 #2:\fi}
\ctr@ld@f\def\FiginsertNd@#1#2:{\ifcat#1a\Figinsert@#1#2,:\else%
    \ifnum\@utoFInDone>\z@\Figinsert@\DefGIfilen@me,#1#2,:\fi\fi}
\ctr@ln@m\Sc@leFact
\ctr@ld@f\def\Figinsert@#1,#2:{\def\t@xt@{#2}\ifx\t@xt@\empty\xdef\Sc@leFact{1}\else%
    \X@rgdeux@#2\xdef\Sc@leFact{\@rgdeux}\fi%
    \Figdisc@rdLTS{#1}{\t@xt@}\@psfgetbb{\t@xt@}%
    \v@lX=\@psfllx\p@\v@lX=\ptpsT@pt\v@lX\v@lX=\Sc@leFact\v@lX%
    \v@lY=\@psflly\p@\v@lY=\ptpsT@pt\v@lY\v@lY=\Sc@leFact\v@lY%
    \b@undb@x{\v@lX}{\v@lY}%
    \v@lX=\@psfurx\p@\v@lX=\ptpsT@pt\v@lX\v@lX=\Sc@leFact\v@lX%
    \v@lY=\@psfury\p@\v@lY=\ptpsT@pt\v@lY\v@lY=\Sc@leFact\v@lY%
    \b@undb@x{\v@lX}{\v@lY}%
    \ifPDFm@ke\Figinclud@PDF{\t@xt@}{\Sc@leFact}\else%
    \v@lX=\c@nt pt\v@lX=\Sc@leFact\v@lX\edef\F@ct{\repdecn@mb{\v@lX}}%
    \ifx\TeXturesonMacOSltX\special{postscriptfile #1 vscale=\F@ct\space hscale=\F@ct}%
    \else\includegraphics{#1}\fi\fi%
    \message{[\t@xt@]}\ignorespaces}
\ctr@ld@f\def\Figdisc@rdLTS#1#2{\expandafter\Figdisc@rdLTS@#1 :#2}
\ctr@ld@f\def\Figdisc@rdLTS@#1 #2:#3{\def#3{#1}\relax\ifx#3\empty\expandafter\Figdisc@rdLTS@#2:#3\fi}
\ctr@ld@f\def\figinsertE#1{\FiginsertE@#1,:\ignorespaces}
\ctr@ld@f\def\FiginsertE@#1,#2:{{\def\t@xt@{#2}\ifx\t@xt@\empty\xdef\Sc@leFact{1}\else%
    \X@rgdeux@#2\xdef\Sc@leFact{\@rgdeux}\fi%
    \Figdisc@rdLTS{#1}{\t@xt@}\pdfximage{\t@xt@}%
    \setbox\Gb@x=\hbox{\pdfrefximage\pdflastximage}%
    \v@lX=\z@\v@lY=-\Sc@leFact\dp\Gb@x\b@undb@x{\v@lX}{\v@lY}%
    \advance\v@lX\Sc@leFact\wd\Gb@x\advance\v@lY\Sc@leFact\dp\Gb@x%
    \advance\v@lY\Sc@leFact\ht\Gb@x\b@undb@x{\v@lX}{\v@lY}%
    \v@lX=\Sc@leFact\wd\Gb@x\pdfximage width \v@lX {\t@xt@}%
    \rlap{\pdfrefximage\pdflastximage}\message{[\t@xt@]}}\ignorespaces}
\ctr@ld@f\def\X@rgdeux@#1,{\edef\@rgdeux{#1}}
\ctr@ln@m\figpt
\ctr@ld@f\def\figptDD#1:#2(#3,#4){\ifps@cri\c@ntr@lnum{#1}%
    {\v@lX=#3\unit@\v@lY=#4\unit@\Fig@dmpt{#2}{\z@}}\ignorespaces\fi}
\ctr@ld@f\def\Fig@dmpt#1#2{\def\t@xt@{#1}\ifx\t@xt@\empty\def\B@@ltxt{\z@}%
    \else\expandafter\gdef\csname\objc@de T\endcsname{#1}\def\B@@ltxt{\@ne}\fi%
    \expandafter\xdef\csname\objc@de\endcsname{\ifitis@vect@r\C@dCl@svect%
    \else\C@dCl@spt\fi,\z@,\B@@ltxt/\the\v@lX,\the\v@lY,#2}}
\ctr@ld@f\def\C@dCl@spt{P}
\ctr@ld@f\def\C@dCl@svect{V}
\ctr@ln@m\c@@rdYZ
\ctr@ln@m\c@@rdY
\ctr@ld@f\def\figptTD#1:#2(#3,#4){\ifps@cri\c@ntr@lnum{#1}%
    \def\c@@rdYZ{#4,0,0}\extrairelepremi@r\c@@rdY\de\c@@rdYZ%
    \extrairelepremi@r\c@@rdZ\de\c@@rdYZ%
    {\v@lX=#3\unit@\v@lY=\c@@rdY\unit@\v@lZ=\c@@rdZ\unit@\Fig@dmpt{#2}{\the\v@lZ}%
    \b@undb@xTD{\v@lX}{\v@lY}{\v@lZ}}\ignorespaces\fi}
\ctr@ln@m\Figp@intreg
\ctr@ld@f\def\Figp@intregDD#1:#2(#3,#4){\c@ntr@lnum{#1}%
    {\result@t=#4\v@lX=#3\v@lY=\result@t\Fig@dmpt{#2}{\z@}}\ignorespaces}
\ctr@ld@f\def\Figp@intregTD#1:#2(#3,#4){\c@ntr@lnum{#1}%
    \def\c@@rdYZ{#4,\z@,\z@}\extrairelepremi@r\c@@rdY\de\c@@rdYZ%
    \extrairelepremi@r\c@@rdZ\de\c@@rdYZ%
    {\v@lX=#3\v@lY=\c@@rdY\v@lZ=\c@@rdZ\Fig@dmpt{#2}{\the\v@lZ}%
    \b@undb@xTD{\v@lX}{\v@lY}{\v@lZ}}\ignorespaces}
\ctr@ln@m\figptBezier
\ctr@ld@f\def\figptBezierDD#1:#2:#3[#4,#5,#6,#7]{\ifps@cri{\s@uvc@ntr@l\et@tfigptBezierDD%
    \FigptBezier@#3[#4,#5,#6,#7]\Figp@intregDD#1:{#2}(\v@lX,\v@lY)%
    \resetc@ntr@l\et@tfigptBezierDD}\ignorespaces\fi}
\ctr@ld@f\def\figptBezierTD#1:#2:#3[#4,#5,#6,#7]{\ifps@cri{\s@uvc@ntr@l\et@tfigptBezierTD%
    \FigptBezier@#3[#4,#5,#6,#7]\Figp@intregTD#1:{#2}(\v@lX,\v@lY,\v@lZ)%
    \resetc@ntr@l\et@tfigptBezierTD}\ignorespaces\fi}
\ctr@ld@f\def\FigptBezier@#1[#2,#3,#4,#5]{\setc@ntr@l{2}%
    \edef\T@{#1}\v@leur=\p@\advance\v@leur-#1pt\edef\UNmT@{\repdecn@mb{\v@leur}}%
    \figptcopy-4:/#2/\figptcopy-3:/#3/\figptcopy-2:/#4/\figptcopy-1:/#5/%
    \l@mbd@un=-4 \l@mbd@de=-\thr@@\p@rtent=\m@ne\c@lDecast%
    \l@mbd@un=-4 \l@mbd@de=-\thr@@\p@rtent=-\tw@\c@lDecast%
    \l@mbd@un=-4 \l@mbd@de=-\thr@@\p@rtent=-\thr@@\c@lDecast\Figg@tXY{-4}}
\ctr@ln@m\c@lDCUn
\ctr@ld@f\def\c@lDCUnDD#1#2{\Figg@tXY{#1}\v@lX=\UNmT@\v@lX\v@lY=\UNmT@\v@lY%
    \Figg@tXYa{#2}\advance\v@lX\T@\v@lXa\advance\v@lY\T@\v@lYa%
    \Figp@intregDD#1:(\v@lX,\v@lY)}
\ctr@ld@f\def\c@lDCUnTD#1#2{\Figg@tXY{#1}\v@lX=\UNmT@\v@lX\v@lY=\UNmT@\v@lY\v@lZ=\UNmT@\v@lZ%
    \Figg@tXYa{#2}\advance\v@lX\T@\v@lXa\advance\v@lY\T@\v@lYa\advance\v@lZ\T@\v@lZa%
    \Figp@intregTD#1:(\v@lX,\v@lY,\v@lZ)}
\ctr@ld@f\def\c@lDecast{\relax\ifnum\l@mbd@un<\p@rtent\c@lDCUn{\l@mbd@un}{\l@mbd@de}%
    \advance\l@mbd@un\@ne\advance\l@mbd@de\@ne\c@lDecast\fi}
\ctr@ld@f\def\figptmap#1:#2=#3/#4/#5/{\ifps@cri{\s@uvc@ntr@l\et@tfigptmap%
    \setc@ntr@l{2}\figvectP-1[#4,#3]\Figg@tXY{-1}%
    \pr@dMatV/#5/\figpttra#1:{#2}=#4/1,-1/%
    \resetc@ntr@l\et@tfigptmap}\ignorespaces\fi}
\ctr@ln@m\pr@dMatV
\ctr@ld@f\def\pr@dMatVDD/#1,#2;#3,#4/{\v@lXa=#1\v@lX\advance\v@lXa#2\v@lY%
    \v@lYa=#3\v@lX\advance\v@lYa#4\v@lY\Figv@ctCreg-1(\v@lXa,\v@lYa)}
\ctr@ld@f\def\pr@dMatVTD/#1,#2,#3;#4,#5,#6;#7,#8,#9/{%
    \v@lXa=#1\v@lX\advance\v@lXa#2\v@lY\advance\v@lXa#3\v@lZ%
    \v@lYa=#4\v@lX\advance\v@lYa#5\v@lY\advance\v@lYa#6\v@lZ%
    \v@lZa=#7\v@lX\advance\v@lZa#8\v@lY\advance\v@lZa#9\v@lZ%
    \Figv@ctCreg-1(\v@lXa,\v@lYa,\v@lZa)}
\ctr@ln@m\figptbary
\ctr@ld@f\def\figptbaryDD#1:#2[#3;#4]{\ifps@cri{\edef\list@num{#3}\extrairelepremi@r\p@int\de\list@num%
    \s@mme=\z@\@ecfor\c@ef:=#4\do{\advance\s@mme\c@ef}%
    \edef\listec@ef{#4,0}\extrairelepremi@r\c@ef\de\listec@ef%
    \Figg@tXY{\p@int}\divide\v@lX\s@mme\divide\v@lY\s@mme%
    \multiply\v@lX\c@ef\multiply\v@lY\c@ef%
    \@ecfor\p@int:=\list@num\do{\extrairelepremi@r\c@ef\de\listec@ef%
           \Figg@tXYa{\p@int}\divide\v@lXa\s@mme\divide\v@lYa\s@mme%
           \multiply\v@lXa\c@ef\multiply\v@lYa\c@ef%
           \advance\v@lX\v@lXa\advance\v@lY\v@lYa}%
    \Figp@intregDD#1:{#2}(\v@lX,\v@lY)}\ignorespaces\fi}
\ctr@ld@f\def\figptbaryTD#1:#2[#3;#4]{\ifps@cri{\edef\list@num{#3}\extrairelepremi@r\p@int\de\list@num%
    \s@mme=\z@\@ecfor\c@ef:=#4\do{\advance\s@mme\c@ef}%
    \edef\listec@ef{#4,0}\extrairelepremi@r\c@ef\de\listec@ef%
    \Figg@tXY{\p@int}\divide\v@lX\s@mme\divide\v@lY\s@mme\divide\v@lZ\s@mme%
    \multiply\v@lX\c@ef\multiply\v@lY\c@ef\multiply\v@lZ\c@ef%
    \@ecfor\p@int:=\list@num\do{\extrairelepremi@r\c@ef\de\listec@ef%
           \Figg@tXYa{\p@int}\divide\v@lXa\s@mme\divide\v@lYa\s@mme\divide\v@lZa\s@mme%
           \multiply\v@lXa\c@ef\multiply\v@lYa\c@ef\multiply\v@lZa\c@ef%
           \advance\v@lX\v@lXa\advance\v@lY\v@lYa\advance\v@lZ\v@lZa}%
    \Figp@intregTD#1:{#2}(\v@lX,\v@lY,\v@lZ)}\ignorespaces\fi}
\ctr@ld@f\def\figptbaryR#1:#2[#3;#4]{\ifps@cri{%
    \v@leur=\z@\@ecfor\c@ef:=#4\do{\maxim@m{\v@lmax}{\c@ef pt}{-\c@ef pt}%
    \ifdim\v@lmax>\v@leur\v@leur=\v@lmax\fi}%
    \ifdim\v@leur<\p@\f@ctech=\@M\else\ifdim\v@leur<\t@n\p@\f@ctech=\@m\else%
    \ifdim\v@leur<\c@nt\p@\f@ctech=\c@nt\else\ifdim\v@leur<\@m\p@\f@ctech=\t@n\else%
    \f@ctech=\@ne\fi\fi\fi\fi%
    \def\listec@ef{0}%
    \@ecfor\c@ef:=#4\do{\sc@lec@nvRI{\c@ef pt}\edef\listec@ef{\listec@ef,\the\s@mme}}%
    \extrairelepremi@r\c@ef\de\listec@ef\figptbary#1:#2[#3;\listec@ef]}\ignorespaces\fi}
\ctr@ld@f\def\sc@lec@nvRI#1{\v@leur=#1\p@rtentiere{\s@mme}{\v@leur}\advance\v@leur-\s@mme\p@%
    \multiply\v@leur\f@ctech\p@rtentiere{\p@rtent}{\v@leur}%
    \multiply\s@mme\f@ctech\advance\s@mme\p@rtent}
\ctr@ln@m\figptcirc
\ctr@ld@f\def\figptcircDD#1:#2:#3;#4(#5){\ifps@cri{\s@uvc@ntr@l\et@tfigptcircDD%
    \c@lptellDD#1:{#2}:#3;#4,#4(#5)\resetc@ntr@l\et@tfigptcircDD}\ignorespaces\fi}
\ctr@ld@f\def\figptcircTD#1:#2:#3,#4,#5;#6(#7){\ifps@cri{\s@uvc@ntr@l\et@tfigptcircTD%
    \setc@ntr@l{2}\c@lExtAxes#3,#4,#5(#6)\figptellP#1:{#2}:#3,-4,-5(#7)%
    \resetc@ntr@l\et@tfigptcircTD}\ignorespaces\fi}
\ctr@ln@m\figptcircumcenter
\ctr@ld@f\def\figptcircumcenterDD#1:#2[#3,#4,#5]{\ifps@cri{\s@uvc@ntr@l\et@tfigptcircumcenterDD%
    \setc@ntr@l{2}\figvectNDD-5[#3,#4]\figptbaryDD-3:[#3,#4;1,1]%
                  \figvectNDD-6[#4,#5]\figptbaryDD-4:[#4,#5;1,1]%
    \resetc@ntr@l{2}\inters@cDD#1:{#2}[-3,-5;-4,-6]%
    \resetc@ntr@l\et@tfigptcircumcenterDD}\ignorespaces\fi}
\ctr@ld@f\def\figptcircumcenterTD#1:#2[#3,#4,#5]{\ifps@cri{\s@uvc@ntr@l\et@tfigptcircumcenterTD%
    \setc@ntr@l{2}\figvectNTD-1[#3,#4,#5]%
    \figvectPTD-3[#3,#4]\figvectNVTD-5[-1,-3]\figptbaryTD-3:[#3,#4;1,1]%
    \figvectPTD-4[#4,#5]\figvectNVTD-6[-1,-4]\figptbaryTD-4:[#4,#5;1,1]%
    \resetc@ntr@l{2}\inters@cTD#1:{#2}[-3,-5;-4,-6]%
    \resetc@ntr@l\et@tfigptcircumcenterTD}\ignorespaces\fi}
\ctr@ln@m\figptcopy
\ctr@ld@f\def\figptcopyDD#1:#2/#3/{\ifps@cri{\Figg@tXY{#3}%
    \Figp@intregDD#1:{#2}(\v@lX,\v@lY)}\ignorespaces\fi}
\ctr@ld@f\def\figptcopyTD#1:#2/#3/{\ifps@cri{\Figg@tXY{#3}%
    \Figp@intregTD#1:{#2}(\v@lX,\v@lY,\v@lZ)}\ignorespaces\fi}
\ctr@ln@m\figptcurvcenter
\ctr@ld@f\def\figptcurvcenterDD#1:#2:#3[#4,#5,#6,#7]{\ifps@cri{\s@uvc@ntr@l\et@tfigptcurvcenterDD%
    \setc@ntr@l{2}\c@lcurvradDD#3[#4,#5,#6,#7]\edef\Sprim@{\repdecn@mb{\result@t}}%
    \figptBezierDD-1::#3[#4,#5,#6,#7]\figpttraDD#1:{#2}=-1/\Sprim@,-5/%
    \resetc@ntr@l\et@tfigptcurvcenterDD}\ignorespaces\fi}
\ctr@ld@f\def\figptcurvcenterTD#1:#2:#3[#4,#5,#6,#7]{\ifps@cri{\s@uvc@ntr@l\et@tfigptcurvcenterTD%
    \setc@ntr@l{2}\figvectDBezierTD -5:1,#3[#4,#5,#6,#7]%
    \figvectDBezierTD -6:2,#3[#4,#5,#6,#7]\vecunit@TD{-5}{-5}%
    \edef\Sprim@{\repdecn@mb{\result@t}}\figvectNVTD-1[-6,-5]%
    \figvectNVTD-5[-5,-1]\c@lproscalTD\v@leur[-6,-5]%
    \invers@{\v@leur}{\v@leur}\v@leur=\Sprim@\v@leur\v@leur=\Sprim@\v@leur%
    \figptBezierTD-1::#3[#4,#5,#6,#7]\edef\Sprim@{\repdecn@mb{\v@leur}}%
    \figpttraTD#1:{#2}=-1/\Sprim@,-5/\resetc@ntr@l\et@tfigptcurvcenterTD}\ignorespaces\fi}
\ctr@ld@f\def\c@lcurvradDD#1[#2,#3,#4,#5]{{\figvectDBezierDD -5:1,#1[#2,#3,#4,#5]%
    \figvectDBezierDD -6:2,#1[#2,#3,#4,#5]\vecunit@DD{-5}{-5}%
    \edef\Sprim@{\repdecn@mb{\result@t}}\figvectNVDD-5[-5]\c@lproscalDD\v@leur[-6,-5]%
    \invers@{\v@leur}{\v@leur}\v@leur=\Sprim@\v@leur\v@leur=\Sprim@\v@leur%
    \global\result@t=\v@leur}}
\ctr@ln@m\figptell
\ctr@ld@f\def\figptellDD#1:#2:#3;#4,#5(#6,#7){\ifps@cri{\s@uvc@ntr@l\et@tfigptell%
    \c@lptellDD#1::#3;#4,#5(#6)\figptrotDD#1:{#2}=#1/#3,#7/%
    \resetc@ntr@l\et@tfigptell}\ignorespaces\fi}
\ctr@ld@f\def\c@lptellDD#1:#2:#3;#4,#5(#6){\c@ssin{\C@}{\S@}{#6}\v@lmin=\C@ pt\v@lmax=\S@ pt%
    \v@lmin=#4\v@lmin\v@lmax=#5\v@lmax%
    \edef\Xc@mp{\repdecn@mb{\v@lmin}}\edef\Yc@mp{\repdecn@mb{\v@lmax}}%
    \setc@ntr@l{2}\figvectC-1(\Xc@mp,\Yc@mp)\figpttraDD#1:{#2}=#3/1,-1/}
\ctr@ld@f\def\figptellP#1:#2:#3,#4,#5(#6){\ifps@cri{\s@uvc@ntr@l\et@tfigptellP%
    \setc@ntr@l{2}\figvectP-1[#3,#4]\figvectP-2[#3,#5]%
    \v@leur=#6pt\c@lptellP{#3}{-1}{-2}\figptcopy#1:{#2}/-3/%
    \resetc@ntr@l\et@tfigptellP}\ignorespaces\fi}
\ctr@ln@m\@ngle
\ctr@ld@f\def\c@lptellP#1#2#3{\edef\@ngle{\repdecn@mb\v@leur}\c@ssin{\C@}{\S@}{\@ngle}%
    \figpttra-3:=#1/\C@,#2/\figpttra-3:=-3/\S@,#3/}
\ctr@ln@m\figptendnormal
\ctr@ld@f\def\figptendnormalDD#1:#2:#3,#4[#5,#6]{\ifps@cri{\s@uvc@ntr@l\et@tfigptendnormal%
    \Figg@tXYa{#5}\Figg@tXY{#6}%
    \advance\v@lX-\v@lXa\advance\v@lY-\v@lYa%
    \setc@ntr@l{2}\Figv@ctCreg-1(\v@lX,\v@lY)\vecunit@{-1}{-1}\Figg@tXY{-1}%
    \delt@=#3\unit@\maxim@m{\delt@}{\delt@}{-\delt@}\edef\l@ngueur{\repdecn@mb{\delt@}}%
    \v@lX=\l@ngueur\v@lX\v@lY=\l@ngueur\v@lY%
    \delt@=\p@\advance\delt@-#4pt\edef\l@ngueur{\repdecn@mb{\delt@}}%
    \figptbaryR-1:[#5,#6;#4,\l@ngueur]\Figg@tXYa{-1}%
    \advance\v@lXa\v@lY\advance\v@lYa-\v@lX%
    \setc@ntr@l{1}\Figp@intregDD#1:{#2}(\v@lXa,\v@lYa)\resetc@ntr@l\et@tfigptendnormal}%
    \ignorespaces\fi}
\ctr@ld@f\def\figptexcenter#1:#2[#3,#4,#5]{\ifps@cri{\let@xte={-}%
    \Figptexinsc@nter#1:#2[#3,#4,#5]}\ignorespaces\fi}
\ctr@ld@f\def\figptincenter#1:#2[#3,#4,#5]{\ifps@cri{\let@xte={}%
    \Figptexinsc@nter#1:#2[#3,#4,#5]}\ignorespaces\fi}
\ctr@ld@f\let\figptinscribedcenter=\figptincenter
\ctr@ld@f\def\Figptexinsc@nter#1:#2[#3,#4,#5]{%
    \figgetdist\LA@[#4,#5]\figgetdist\LB@[#3,#5]\figgetdist\LC@[#3,#4]%
    \figptbaryR#1:{#2}[#3,#4,#5;\the\let@xte\LA@,\LB@,\LC@]}
\ctr@ln@m\figptinterlineplane
\ctr@ld@f\def\figptinterlineplaneDD{\un@v@ilable{figptinterlineplane}}
\ctr@ld@f\def\figptinterlineplaneTD#1:#2[#3,#4;#5,#6]{\ifps@cri{\s@uvc@ntr@l\et@tfigptinterlineplane%
    \setc@ntr@l{2}\figvectPTD-1[#3,#5]\vecunit@TD{-2}{#6}%
    \r@pPSTD\v@leur[-2,-1,#4]\edef\v@lcoef{\repdecn@mb{\v@leur}}%
    \figpttraTD#1:{#2}=#3/\v@lcoef,#4/\resetc@ntr@l\et@tfigptinterlineplane}\ignorespaces\fi}
\ctr@ln@m\figptorthocenter
\ctr@ld@f\def\figptorthocenterDD#1:#2[#3,#4,#5]{\ifps@cri{\s@uvc@ntr@l\et@tfigptorthocenterDD%
    \setc@ntr@l{2}\figvectNDD-3[#3,#4]\figvectNDD-4[#4,#5]%
    \resetc@ntr@l{2}\inters@cDD#1:{#2}[#5,-3;#3,-4]%
    \resetc@ntr@l\et@tfigptorthocenterDD}\ignorespaces\fi}
\ctr@ld@f\def\figptorthocenterTD#1:#2[#3,#4,#5]{\ifps@cri{\s@uvc@ntr@l\et@tfigptorthocenterTD%
    \setc@ntr@l{2}\figvectNTD-1[#3,#4,#5]%
    \figvectPTD-2[#3,#4]\figvectNVTD-3[-1,-2]%
    \figvectPTD-2[#4,#5]\figvectNVTD-4[-1,-2]%
    \resetc@ntr@l{2}\inters@cTD#1:{#2}[#5,-3;#3,-4]%
    \resetc@ntr@l\et@tfigptorthocenterTD}\ignorespaces\fi}
\ctr@ln@m\figptorthoprojline
\ctr@ld@f\def\figptorthoprojlineDD#1:#2=#3/#4,#5/{\ifps@cri{\s@uvc@ntr@l\et@tfigptorthoprojlineDD%
    \setc@ntr@l{2}\figvectPDD-3[#4,#5]\figvectNVDD-4[-3]\resetc@ntr@l{2}%
    \inters@cDD#1:{#2}[#3,-4;#4,-3]\resetc@ntr@l\et@tfigptorthoprojlineDD}\ignorespaces\fi}
\ctr@ld@f\def\figptorthoprojlineTD#1:#2=#3/#4,#5/{\ifps@cri{\s@uvc@ntr@l\et@tfigptorthoprojlineTD%
    \setc@ntr@l{2}\figvectPTD-1[#4,#3]\figvectPTD-2[#4,#5]\vecunit@TD{-2}{-2}%
    \c@lproscalTD\v@leur[-1,-2]\edef\v@lcoef{\repdecn@mb{\v@leur}}%
    \figpttraTD#1:{#2}=#4/\v@lcoef,-2/\resetc@ntr@l\et@tfigptorthoprojlineTD}\ignorespaces\fi}
\ctr@ln@m\figptorthoprojplane
\ctr@ld@f\def\figptorthoprojplaneDD{\un@v@ilable{figptorthoprojplane}}
\ctr@ld@f\def\figptorthoprojplaneTD#1:#2=#3/#4,#5/{\ifps@cri{\s@uvc@ntr@l\et@tfigptorthoprojplane%
    \setc@ntr@l{2}\figvectPTD-1[#3,#4]\vecunit@TD{-2}{#5}%
    \c@lproscalTD\v@leur[-1,-2]\edef\v@lcoef{\repdecn@mb{\v@leur}}%
    \figpttraTD#1:{#2}=#3/\v@lcoef,-2/\resetc@ntr@l\et@tfigptorthoprojplane}\ignorespaces\fi}
\ctr@ld@f\def\figpthom#1:#2=#3/#4,#5/{\ifps@cri{\s@uvc@ntr@l\et@tfigpthom%
    \setc@ntr@l{2}\figvectP-1[#4,#3]\figpttra#1:{#2}=#4/#5,-1/%
    \resetc@ntr@l\et@tfigpthom}\ignorespaces\fi}
\ctr@ln@m\figptrot
\ctr@ld@f\def\figptrotDD#1:#2=#3/#4,#5/{\ifps@cri{\s@uvc@ntr@l\et@tfigptrotDD%
    \c@ssin{\C@}{\S@}{#5}\setc@ntr@l{2}\figvectPDD-1[#4,#3]\Figg@tXY{-1}%
    \v@lXa=\C@\v@lX\advance\v@lXa-\S@\v@lY%
    \v@lYa=\S@\v@lX\advance\v@lYa\C@\v@lY%
    \Figv@ctCreg-1(\v@lXa,\v@lYa)\figpttraDD#1:{#2}=#4/1,-1/%
    \resetc@ntr@l\et@tfigptrotDD}\ignorespaces\fi}
\ctr@ld@f\def\figptrotTD#1:#2=#3/#4,#5,#6/{\ifps@cri{\s@uvc@ntr@l\et@tfigptrotTD%
    \c@ssin{\C@}{\S@}{#5}%
    \setc@ntr@l{2}\figptorthoprojplaneTD-3:=#4/#3,#6/\figvectPTD-2[-3,#3]%
    \n@rmeucTD\v@leur{-2}\ifdim\v@leur<\Cepsil@n\Figg@tXYa{#3}\else%
    \edef\v@lcoef{\repdecn@mb{\v@leur}}\figvectNVTD-1[#6,-2]%
    \Figg@tXYa{-1}\v@lXa=\v@lcoef\v@lXa\v@lYa=\v@lcoef\v@lYa\v@lZa=\v@lcoef\v@lZa%
    \v@lXa=\S@\v@lXa\v@lYa=\S@\v@lYa\v@lZa=\S@\v@lZa\Figg@tXY{-2}%
    \advance\v@lXa\C@\v@lX\advance\v@lYa\C@\v@lY\advance\v@lZa\C@\v@lZ%
    \Figg@tXY{-3}\advance\v@lXa\v@lX\advance\v@lYa\v@lY\advance\v@lZa\v@lZ\fi%
    \Figp@intregTD#1:{#2}(\v@lXa,\v@lYa,\v@lZa)\resetc@ntr@l\et@tfigptrotTD}\ignorespaces\fi}
\ctr@ln@m\figptsym
\ctr@ld@f\def\figptsymDD#1:#2=#3/#4,#5/{\ifps@cri{\s@uvc@ntr@l\et@tfigptsymDD%
    \resetc@ntr@l{2}\figptorthoprojlineDD-5:=#3/#4,#5/\figvectPDD-2[#3,-5]%
    \figpttraDD#1:{#2}=#3/2,-2/\resetc@ntr@l\et@tfigptsymDD}\ignorespaces\fi}
\ctr@ld@f\def\figptsymTD#1:#2=#3/#4,#5/{\ifps@cri{\s@uvc@ntr@l\et@tfigptsymTD%
    \resetc@ntr@l{2}\figptorthoprojplaneTD-3:=#3/#4,#5/\figvectPTD-2[#3,-3]%
    \figpttraTD#1:{#2}=#3/2,-2/\resetc@ntr@l\et@tfigptsymTD}\ignorespaces\fi}
\ctr@ln@m\figpttra
\ctr@ld@f\def\figpttraDD#1:#2=#3/#4,#5/{\ifps@cri{\Figg@tXYa{#5}\v@lXa=#4\v@lXa\v@lYa=#4\v@lYa%
    \Figg@tXY{#3}\advance\v@lX\v@lXa\advance\v@lY\v@lYa%
    \Figp@intregDD#1:{#2}(\v@lX,\v@lY)}\ignorespaces\fi}
\ctr@ld@f\def\figpttraTD#1:#2=#3/#4,#5/{\ifps@cri{\Figg@tXYa{#5}\v@lXa=#4\v@lXa\v@lYa=#4\v@lYa%
    \v@lZa=#4\v@lZa\Figg@tXY{#3}\advance\v@lX\v@lXa\advance\v@lY\v@lYa%
    \advance\v@lZ\v@lZa\Figp@intregTD#1:{#2}(\v@lX,\v@lY,\v@lZ)}\ignorespaces\fi}
\ctr@ln@m\figpttraC
\ctr@ld@f\def\figpttraCDD#1:#2=#3/#4,#5/{\ifps@cri{\v@lXa=#4\unit@\v@lYa=#5\unit@%
    \Figg@tXY{#3}\advance\v@lX\v@lXa\advance\v@lY\v@lYa%
    \Figp@intregDD#1:{#2}(\v@lX,\v@lY)}\ignorespaces\fi}
\ctr@ld@f\def\figpttraCTD#1:#2=#3/#4,#5,#6/{\ifps@cri{\v@lXa=#4\unit@\v@lYa=#5\unit@\v@lZa=#6\unit@%
    \Figg@tXY{#3}\advance\v@lX\v@lXa\advance\v@lY\v@lYa\advance\v@lZ\v@lZa%
    \Figp@intregTD#1:{#2}(\v@lX,\v@lY,\v@lZ)}\ignorespaces\fi}
\ctr@ld@f\def\figptsaxes#1:#2(#3){\ifps@cri{\an@lys@xes#3,:\ifx\t@xt@\empty%
    \ifTr@isDim\Figpts@xes#1:#2(0,#3,0,#3,0,#3)\else\Figpts@xes#1:#2(0,#3,0,#3)\fi%
    \else\Figpts@xes#1:#2(#3)\fi}\ignorespaces\fi}
\ctr@ln@m\Figpts@xes
\ctr@ld@f\def\Figpts@xesDD#1:#2(#3,#4,#5,#6){%
    \s@mme=#1\figpttraC\the\s@mme:$x$=#2/#4,0/%
    \advance\s@mme\@ne\figpttraC\the\s@mme:$y$=#2/0,#6/}
\ctr@ld@f\def\Figpts@xesTD#1:#2(#3,#4,#5,#6,#7,#8){%
    \s@mme=#1\figpttraC\the\s@mme:$x$=#2/#4,0,0/%
    \advance\s@mme\@ne\figpttraC\the\s@mme:$y$=#2/0,#6,0/%
    \advance\s@mme\@ne\figpttraC\the\s@mme:$z$=#2/0,0,#8/}
\ctr@ld@f\def\figptsmap#1=#2/#3/#4/{\ifps@cri{\s@uvc@ntr@l\et@tfigptsmap%
    \setc@ntr@l{2}\def\list@num{#2}\s@mme=#1%
    \@ecfor\p@int:=\list@num\do{\figvectP-1[#3,\p@int]\Figg@tXY{-1}%
    \pr@dMatV/#4/\figpttra\the\s@mme:=#3/1,-1/\advance\s@mme\@ne}%
    \resetc@ntr@l\et@tfigptsmap}\ignorespaces\fi}
\ctr@ln@m\figptscontrol
\ctr@ld@f\def\figptscontrolDD#1[#2,#3,#4,#5]{\ifps@cri{\s@uvc@ntr@l\et@tfigptscontrolDD\setc@ntr@l{2}%
    \v@lX=\z@\v@lY=\z@\Figtr@nptDD{-5}{#2}\Figtr@nptDD{2}{#5}%
    \divide\v@lX\@vi\divide\v@lY\@vi%
    \Figtr@nptDD{3}{#3}\Figtr@nptDD{-1.5}{#4}\Figp@intregDD-1:(\v@lX,\v@lY)%
    \v@lX=\z@\v@lY=\z@\Figtr@nptDD{2}{#2}\Figtr@nptDD{-5}{#5}%
    \divide\v@lX\@vi\divide\v@lY\@vi\Figtr@nptDD{-1.5}{#3}\Figtr@nptDD{3}{#4}%
    \s@mme=#1\advance\s@mme\@ne\Figp@intregDD\the\s@mme:(\v@lX,\v@lY)%
    \figptcopyDD#1:/-1/\resetc@ntr@l\et@tfigptscontrolDD}\ignorespaces\fi}
\ctr@ld@f\def\figptscontrolTD#1[#2,#3,#4,#5]{\ifps@cri{\s@uvc@ntr@l\et@tfigptscontrolTD\setc@ntr@l{2}%
    \v@lX=\z@\v@lY=\z@\v@lZ=\z@\Figtr@nptTD{-5}{#2}\Figtr@nptTD{2}{#5}%
    \divide\v@lX\@vi\divide\v@lY\@vi\divide\v@lZ\@vi%
    \Figtr@nptTD{3}{#3}\Figtr@nptTD{-1.5}{#4}\Figp@intregTD-1:(\v@lX,\v@lY,\v@lZ)%
    \v@lX=\z@\v@lY=\z@\v@lZ=\z@\Figtr@nptTD{2}{#2}\Figtr@nptTD{-5}{#5}%
    \divide\v@lX\@vi\divide\v@lY\@vi\divide\v@lZ\@vi\Figtr@nptTD{-1.5}{#3}\Figtr@nptTD{3}{#4}%
    \s@mme=#1\advance\s@mme\@ne\Figp@intregTD\the\s@mme:(\v@lX,\v@lY,\v@lZ)%
    \figptcopyTD#1:/-1/\resetc@ntr@l\et@tfigptscontrolTD}\ignorespaces\fi}
\ctr@ld@f\def\Figtr@nptDD#1#2{\Figg@tXYa{#2}\v@lXa=#1\v@lXa\v@lYa=#1\v@lYa%
    \advance\v@lX\v@lXa\advance\v@lY\v@lYa}
\ctr@ld@f\def\Figtr@nptTD#1#2{\Figg@tXYa{#2}\v@lXa=#1\v@lXa\v@lYa=#1\v@lYa\v@lZa=#1\v@lZa%
    \advance\v@lX\v@lXa\advance\v@lY\v@lYa\advance\v@lZ\v@lZa}
\ctr@ld@f\def\figptscontrolcurve#1,#2[#3]{\ifps@cri{\s@uvc@ntr@l\et@tfigptscontrolcurve%
    \def\list@num{#3}\extrairelepremi@r\Ak@\de\list@num%
    \extrairelepremi@r\Ai@\de\list@num\extrairelepremi@r\Aj@\de\list@num%
    \s@mme=#1\figptcopy\the\s@mme:/\Ai@/%
    \setc@ntr@l{2}\figvectP -1[\Ak@,\Aj@]%
    \@ecfor\Ak@:=\list@num\do{\advance\s@mme\@ne\figpttra\the\s@mme:=\Ai@/\curv@roundness,-1/%
       \figvectP -1[\Ai@,\Ak@]\advance\s@mme\@ne\figpttra\the\s@mme:=\Aj@/-\curv@roundness,-1/%
       \advance\s@mme\@ne\figptcopy\the\s@mme:/\Aj@/%
       \edef\Ai@{\Aj@}\edef\Aj@{\Ak@}}\advance\s@mme-#1\divide\s@mme\thr@@%
       \xdef#2{\the\s@mme}%
    \resetc@ntr@l\et@tfigptscontrolcurve}\ignorespaces\fi}
\ctr@ln@m\figptsintercirc
\ctr@ld@f\def\figptsintercircDD#1[#2,#3;#4,#5]{\ifps@cri{\s@uvc@ntr@l\et@tfigptsintercircDD%
    \setc@ntr@l{2}\let\c@lNVintc=\c@lNVintcDD\Figptsintercirc@#1[#2,#3;#4,#5]%
    \resetc@ntr@l\et@tfigptsintercircDD}\ignorespaces\fi}
\ctr@ld@f\def\figptsintercircTD#1[#2,#3;#4,#5;#6]{\ifps@cri{\s@uvc@ntr@l\et@tfigptsintercircTD%
    \setc@ntr@l{2}\let\c@lNVintc=\c@lNVintcTD\vecunitC@TD[#2,#6]%
    \Figv@ctCreg-3(\v@lX,\v@lY,\v@lZ)\Figptsintercirc@#1[#2,#3;#4,#5]%
    \resetc@ntr@l\et@tfigptsintercircTD}\ignorespaces\fi}
\ctr@ld@f\def\Figptsintercirc@#1[#2,#3;#4,#5]{\figvectP-1[#2,#4]%
    \vecunit@{-1}{-1}\delt@=\result@t\f@ctech=\result@tent%
    \s@mme=#1\advance\s@mme\@ne\figptcopy#1:/#2/\figptcopy\the\s@mme:/#4/%
    \ifdim\delt@=\z@\else%
    \v@lmin=#3\unit@\v@lmax=#5\unit@\v@leur=\v@lmin\advance\v@leur\v@lmax%
    \ifdim\v@leur>\delt@%
    \v@leur=\v@lmin\advance\v@leur-\v@lmax\maxim@m{\v@leur}{\v@leur}{-\v@leur}%
    \ifdim\v@leur<\delt@%
    \divide\v@lmin\f@ctech\divide\v@lmax\f@ctech\divide\delt@\f@ctech%
    \v@lmin=\repdecn@mb{\v@lmin}\v@lmin\v@lmax=\repdecn@mb{\v@lmax}\v@lmax%
    \invers@{\v@leur}{\delt@}\advance\v@lmax-\v@lmin%
    \v@lmax=-\repdecn@mb{\v@leur}\v@lmax\advance\delt@\v@lmax\delt@=.5\delt@%
    \v@lmax=\delt@\multiply\v@lmax\f@ctech%
    \edef\t@ille{\repdecn@mb{\v@lmax}}\figpttra-2:=#2/\t@ille,-1/%
    \delt@=\repdecn@mb{\delt@}\delt@\advance\v@lmin-\delt@%
    \sqrt@{\v@leur}{\v@lmin}\multiply\v@leur\f@ctech\edef\t@ille{\repdecn@mb{\v@leur}}%
    \c@lNVintc\figpttra#1:=-2/-\t@ille,-1/\figpttra\the\s@mme:=-2/\t@ille,-1/\fi\fi\fi}
\ctr@ld@f\def\c@lNVintcDD{\Figg@tXY{-1}\Figv@ctCreg-1(-\v@lY,\v@lX)} 
\ctr@ld@f\def\c@lNVintcTD{{\Figg@tXY{-3}\v@lmin=\v@lX\v@lmax=\v@lY\v@leur=\v@lZ%
    \Figg@tXY{-1}\c@lprovec{-3}\vecunit@{-3}{-3}
    \Figg@tXY{-1}\v@lmin=\v@lX\v@lmax=\v@lY%
    \v@leur=\v@lZ\Figg@tXY{-3}\c@lprovec{-1}}} 
\ctr@ln@m\figptsinterlinell
\ctr@ld@f\def\figptsinterlinellDD#1[#2,#3,#4,#5;#6,#7]{\ifps@cri{\s@uvc@ntr@l\et@tfigptsinterlinellDD%
    \figptcopy#1:/#6/\s@mme=#1\advance\s@mme\@ne\figptcopy\the\s@mme:/#7/%
    \v@lmin=#3\unit@\v@lmax=#4\unit@
    \setc@ntr@l{2}\figptbaryDD-4:[#6,#7;1,1]\figptsrotDD-3=-4,#7/#2,-#5/
    \Figg@tXY{-3}\Figg@tXYa{#2}\advance\v@lX-\v@lXa\advance\v@lY-\v@lYa
    \figvectP-1[-3,-2]\Figg@tXYa{-1}\figvectP-3[-4,#7]\Figptsint@rLE{#1}
    \resetc@ntr@l\et@tfigptsinterlinellDD}\ignorespaces\fi}
\ctr@ld@f\def\figptsinterlinellP#1[#2,#3,#4;#5,#6]{\ifps@cri{\s@uvc@ntr@l\et@tfigptsinterlinellP%
    \figptcopy#1:/#5/\s@mme=#1\advance\s@mme\@ne\figptcopy\the\s@mme:/#6/\setc@ntr@l{2}%
    \figvectP-1[#2,#3]\vecunit@{-1}{-1}\v@lmin=\result@t
    \figvectP-2[#2,#4]\vecunit@{-2}{-2}\v@lmax=\result@t
    \figptbary-4:[#5,#6;1,1]
    \figvectP-3[#2,-4]\c@lproscal\v@lX[-3,-1]\c@lproscal\v@lY[-3,-2]
    \figvectP-3[-4,#6]\c@lproscal\v@lXa[-3,-1]\c@lproscal\v@lYa[-3,-2]
    \Figptsint@rLE{#1}\resetc@ntr@l\et@tfigptsinterlinellP}\ignorespaces\fi}
\ctr@ld@f\def\Figptsint@rLE#1{%
    \getredf@ctDD\f@ctech(\v@lmin,\v@lmax)%
    \getredf@ctDD\p@rtent(\v@lX,\v@lY)\ifnum\p@rtent>\f@ctech\f@ctech=\p@rtent\fi%
    \getredf@ctDD\p@rtent(\v@lXa,\v@lYa)\ifnum\p@rtent>\f@ctech\f@ctech=\p@rtent\fi%
    \divide\v@lmin\f@ctech\divide\v@lmax\f@ctech\divide\v@lX\f@ctech\divide\v@lY\f@ctech%
    \divide\v@lXa\f@ctech\divide\v@lYa\f@ctech%
    \c@rre=\repdecn@mb\v@lXa\v@lmax\mili@u=\repdecn@mb\v@lYa\v@lmin%
    \getredf@ctDD\f@ctech(\c@rre,\mili@u)%
    \c@rre=\repdecn@mb\v@lX\v@lmax\mili@u=\repdecn@mb\v@lY\v@lmin%
    \getredf@ctDD\p@rtent(\c@rre,\mili@u)\ifnum\p@rtent>\f@ctech\f@ctech=\p@rtent\fi%
    \divide\v@lmin\f@ctech\divide\v@lmax\f@ctech\divide\v@lX\f@ctech\divide\v@lY\f@ctech%
    \divide\v@lXa\f@ctech\divide\v@lYa\f@ctech%
    \v@lmin=\repdecn@mb{\v@lmin}\v@lmin\v@lmax=\repdecn@mb{\v@lmax}\v@lmax%
    \edef\G@xde{\repdecn@mb\v@lmin}\edef\P@xde{\repdecn@mb\v@lmax}%
    \c@rre=-\v@lmax\v@leur=\repdecn@mb\v@lY\v@lY\advance\c@rre\v@leur\c@rre=\G@xde\c@rre%
    \v@leur=\repdecn@mb\v@lX\v@lX\v@leur=\P@xde\v@leur\advance\c@rre\v@leur
    \v@lmin=\repdecn@mb\v@lYa\v@lmin\v@lmax=\repdecn@mb\v@lXa\v@lmax%
    \mili@u=\repdecn@mb\v@lX\v@lmax\advance\mili@u\repdecn@mb\v@lY\v@lmin
    \v@lmax=\repdecn@mb\v@lXa\v@lmax\advance\v@lmax\repdecn@mb\v@lYa\v@lmin
    \ifdim\v@lmax>\epsil@n%
    \maxim@m{\v@leur}{\c@rre}{-\c@rre}\maxim@m{\v@lmin}{\mili@u}{-\mili@u}%
    \maxim@m{\v@leur}{\v@leur}{\v@lmin}\maxim@m{\v@lmin}{\v@lmax}{-\v@lmax}%
    \maxim@m{\v@leur}{\v@leur}{\v@lmin}\p@rtentiere{\p@rtent}{\v@leur}\advance\p@rtent\@ne%
    \divide\c@rre\p@rtent\divide\mili@u\p@rtent\divide\v@lmax\p@rtent%
    \delt@=\repdecn@mb{\mili@u}\mili@u\v@leur=\repdecn@mb{\v@lmax}\c@rre%
    \advance\delt@-\v@leur\ifdim\delt@<\z@\else\sqrt@\delt@\delt@%
    \invers@\v@lmax\v@lmax\edef\Uns@rAp{\repdecn@mb\v@lmax}%
    \v@leur=-\mili@u\advance\v@leur-\delt@\v@leur=\Uns@rAp\v@leur%
    \edef\t@ille{\repdecn@mb\v@leur}\figpttra#1:=-4/\t@ille,-3/\s@mme=#1\advance\s@mme\@ne%
    \v@leur=-\mili@u\advance\v@leur\delt@\v@leur=\Uns@rAp\v@leur%
    \edef\t@ille{\repdecn@mb\v@leur}\figpttra\the\s@mme:=-4/\t@ille,-3/\fi\fi}
\ctr@ln@m\figptsorthoprojline
\ctr@ld@f\def\figptsorthoprojlineDD#1=#2/#3,#4/{\ifps@cri{\s@uvc@ntr@l\et@tfigptsorthoprojlineDD%
    \setc@ntr@l{2}\figvectPDD-3[#3,#4]\figvectNVDD-4[-3]\resetc@ntr@l{2}%
    \def\list@num{#2}\s@mme=#1\@ecfor\p@int:=\list@num\do{%
    \inters@cDD\the\s@mme:[\p@int,-4;#3,-3]\advance\s@mme\@ne}%
    \resetc@ntr@l\et@tfigptsorthoprojlineDD}\ignorespaces\fi}
\ctr@ld@f\def\figptsorthoprojlineTD#1=#2/#3,#4/{\ifps@cri{\s@uvc@ntr@l\et@tfigptsorthoprojlineTD%
    \setc@ntr@l{2}\figvectPTD-2[#3,#4]\vecunit@TD{-2}{-2}%
    \def\list@num{#2}\s@mme=#1\@ecfor\p@int:=\list@num\do{%
    \figvectPTD-1[#3,\p@int]\c@lproscalTD\v@leur[-1,-2]%
    \edef\v@lcoef{\repdecn@mb{\v@leur}}\figpttraTD\the\s@mme:=#3/\v@lcoef,-2/%
    \advance\s@mme\@ne}\resetc@ntr@l\et@tfigptsorthoprojlineTD}\ignorespaces\fi}
\ctr@ln@m\figptsorthoprojplane
\ctr@ld@f\def\figptsorthoprojplaneDD{\un@v@ilable{figptsorthoprojplane}}
\ctr@ld@f\def\figptsorthoprojplaneTD#1=#2/#3,#4/{\ifps@cri{\s@uvc@ntr@l\et@tfigptsorthoprojplane%
    \setc@ntr@l{2}\vecunit@TD{-2}{#4}%
    \def\list@num{#2}\s@mme=#1\@ecfor\p@int:=\list@num\do{\figvectPTD-1[\p@int,#3]%
    \c@lproscalTD\v@leur[-1,-2]\edef\v@lcoef{\repdecn@mb{\v@leur}}%
    \figpttraTD\the\s@mme:=\p@int/\v@lcoef,-2/\advance\s@mme\@ne}%
    \resetc@ntr@l\et@tfigptsorthoprojplane}\ignorespaces\fi}
\ctr@ld@f\def\figptshom#1=#2/#3,#4/{\ifps@cri{\s@uvc@ntr@l\et@tfigptshom%
    \setc@ntr@l{2}\def\list@num{#2}\s@mme=#1%
    \@ecfor\p@int:=\list@num\do{\figvectP-1[#3,\p@int]%
    \figpttra\the\s@mme:=#3/#4,-1/\advance\s@mme\@ne}%
    \resetc@ntr@l\et@tfigptshom}\ignorespaces\fi}
\ctr@ln@m\figptsrot
\ctr@ld@f\def\figptsrotDD#1=#2/#3,#4/{\ifps@cri{\s@uvc@ntr@l\et@tfigptsrotDD%
    \c@ssin{\C@}{\S@}{#4}\setc@ntr@l{2}\def\list@num{#2}\s@mme=#1%
    \@ecfor\p@int:=\list@num\do{\figvectPDD-1[#3,\p@int]\Figg@tXY{-1}%
    \v@lXa=\C@\v@lX\advance\v@lXa-\S@\v@lY%
    \v@lYa=\S@\v@lX\advance\v@lYa\C@\v@lY%
    \Figv@ctCreg-1(\v@lXa,\v@lYa)\figpttraDD\the\s@mme:=#3/1,-1/\advance\s@mme\@ne}%
    \resetc@ntr@l\et@tfigptsrotDD}\ignorespaces\fi}
\ctr@ld@f\def\figptsrotTD#1=#2/#3,#4,#5/{\ifps@cri{\s@uvc@ntr@l\et@tfigptsrotTD%
    \c@ssin{\C@}{\S@}{#4}%
    \setc@ntr@l{2}\def\list@num{#2}\s@mme=#1%
    \@ecfor\p@int:=\list@num\do{\figptorthoprojplaneTD-3:=#3/\p@int,#5/%
    \figvectPTD-2[-3,\p@int]%
    \figvectNVTD-1[#5,-2]\n@rmeucTD\v@leur{-2}\edef\v@lcoef{\repdecn@mb{\v@leur}}%
    \Figg@tXYa{-1}\v@lXa=\v@lcoef\v@lXa\v@lYa=\v@lcoef\v@lYa\v@lZa=\v@lcoef\v@lZa%
    \v@lXa=\S@\v@lXa\v@lYa=\S@\v@lYa\v@lZa=\S@\v@lZa\Figg@tXY{-2}%
    \advance\v@lXa\C@\v@lX\advance\v@lYa\C@\v@lY\advance\v@lZa\C@\v@lZ%
    \Figg@tXY{-3}\advance\v@lXa\v@lX\advance\v@lYa\v@lY\advance\v@lZa\v@lZ%
    \Figp@intregTD\the\s@mme:(\v@lXa,\v@lYa,\v@lZa)\advance\s@mme\@ne}%
    \resetc@ntr@l\et@tfigptsrotTD}\ignorespaces\fi}
\ctr@ln@m\figptssym
\ctr@ld@f\def\figptssymDD#1=#2/#3,#4/{\ifps@cri{\s@uvc@ntr@l\et@tfigptssymDD%
    \setc@ntr@l{2}\figvectPDD-3[#3,#4]\Figg@tXY{-3}\Figv@ctCreg-4(-\v@lY,\v@lX)%
    \resetc@ntr@l{2}\def\list@num{#2}\s@mme=#1%
    \@ecfor\p@int:=\list@num\do{\inters@cDD-5:[#3,-3;\p@int,-4]\figvectPDD-2[\p@int,-5]%
    \figpttraDD\the\s@mme:=\p@int/2,-2/\advance\s@mme\@ne}%
    \resetc@ntr@l\et@tfigptssymDD}\ignorespaces\fi}
\ctr@ld@f\def\figptssymTD#1=#2/#3,#4/{\ifps@cri{\s@uvc@ntr@l\et@tfigptssymTD%
    \setc@ntr@l{2}\vecunit@TD{-2}{#4}\def\list@num{#2}\s@mme=#1%
    \@ecfor\p@int:=\list@num\do{\figvectPTD-1[\p@int,#3]%
    \c@lproscalTD\v@leur[-1,-2]\v@leur=2\v@leur\edef\v@lcoef{\repdecn@mb{\v@leur}}%
    \figpttraTD\the\s@mme:=\p@int/\v@lcoef,-2/\advance\s@mme\@ne}%
    \resetc@ntr@l\et@tfigptssymTD}\ignorespaces\fi}
\ctr@ln@m\figptstra
\ctr@ld@f\def\figptstraDD#1=#2/#3,#4/{\ifps@cri{\Figg@tXYa{#4}\v@lXa=#3\v@lXa\v@lYa=#3\v@lYa%
    \def\list@num{#2}\s@mme=#1\@ecfor\p@int:=\list@num\do{\Figg@tXY{\p@int}%
    \advance\v@lX\v@lXa\advance\v@lY\v@lYa%
    \Figp@intregDD\the\s@mme:(\v@lX,\v@lY)\advance\s@mme\@ne}}\ignorespaces\fi}
\ctr@ld@f\def\figptstraTD#1=#2/#3,#4/{\ifps@cri{\Figg@tXYa{#4}\v@lXa=#3\v@lXa\v@lYa=#3\v@lYa%
    \v@lZa=#3\v@lZa\def\list@num{#2}\s@mme=#1\@ecfor\p@int:=\list@num\do{\Figg@tXY{\p@int}%
    \advance\v@lX\v@lXa\advance\v@lY\v@lYa\advance\v@lZ\v@lZa%
    \Figp@intregTD\the\s@mme:(\v@lX,\v@lY,\v@lZ)\advance\s@mme\@ne}}\ignorespaces\fi}
\ctr@ln@m\figptvisilimSL
\ctr@ld@f\def\figptvisilimSLDD{\un@v@ilable{figptvisilimSL}}
\ctr@ld@f\def\figptvisilimSLTD#1:#2[#3,#4;#5,#6]{\ifps@cri{\s@uvc@ntr@l\et@tfigptvisilimSLTD%
    \setc@ntr@l{2}\figvectP-1[#3,#4]\n@rminf{\delt@}{-1}%
    \ifcase\curr@ntproj\v@lX=\cxa@\p@\v@lY=-\p@\v@lZ=\cxb@\p@
    \Figv@ctCreg-2(\v@lX,\v@lY,\v@lZ)\figvectP-3[#5,#6]\figvectNV-1[-2,-3]%
    \or\figvectP-1[#5,#6]\vecunitCV@TD{-1}\v@lmin=\v@lX\v@lmax=\v@lY
    \v@leur=\v@lZ\v@lX=\cza@\p@\v@lY=\czb@\p@\v@lZ=\czc@\p@\c@lprovec{-1}%
    \or\c@ley@pt{-2}\figvectN-1[#5,#6,-2]\fi
    \edef\Ai@{#3}\edef\Aj@{#4}\figvectP-2[#5,\Ai@]\c@lproscal\v@leur[-1,-2]%
    \ifdim\v@leur>\z@\p@rtent=\@ne\else\p@rtent=\m@ne\fi%
    \figvectP-2[#5,\Aj@]\c@lproscal\v@leur[-1,-2]%
    \ifdim\p@rtent\v@leur>\z@\figptcopy#1:#2/#3/%
    \message{*** \BS@ figptvisilimSL: points are on the same side.}\else%
    \figptcopy-3:/#3/\figptcopy-4:/#4/%
    \loop\figptbary-5:[-3,-4;1,1]\figvectP-2[#5,-5]\c@lproscal\v@leur[-1,-2]%
    \ifdim\p@rtent\v@leur>\z@\figptcopy-3:/-5/\else\figptcopy-4:/-5/\fi%
    \divide\delt@\tw@\ifdim\delt@>\epsil@n\repeat%
    \figptbary#1:#2[-3,-4;1,1]\fi\resetc@ntr@l\et@tfigptvisilimSLTD}\ignorespaces\fi}
\ctr@ld@f\def\c@ley@pt#1{\t@stp@r\ifitis@K\v@lX=\cza@\p@\v@lY=\czb@\p@\v@lZ=\czc@\p@%
    \Figv@ctCreg-1(\v@lX,\v@lY,\v@lZ)\Figp@intreg-2:(\wd\Bt@rget,\ht\Bt@rget,\dp\Bt@rget)%
    \figpttra#1:=-2/-\disob@intern,-1/\else\end\fi}
\ctr@ld@f\def\t@stp@r{\itis@Ktrue\ifnewt@rgetpt\else\itis@Kfalse%
    \message{*** \BS@ figptvisilimXX: target point undefined.}\fi\ifnewdis@b\else%
    \itis@Kfalse\message{*** \BS@ figptvisilimXX: observation distance undefined.}\fi%
    \ifitis@K\else\message{*** This macro must be called after \BS@ psbeginfig or after
    having set the missing parameter(s) with \BS@ figset proj()}\fi}
\ctr@ld@f\def\figscan#1(#2,#3){{\s@uvc@ntr@l\et@tfigscan\@psfgetbb{#1}\if@psfbbfound\else%
    \def\@psfllx{0}\def\@psflly{20}\def\@psfurx{540}\def\@psfury{640}\fi\figscan@{#2}{#3}%
    \resetc@ntr@l\et@tfigscan}\ignorespaces}
\ctr@ld@f\def\figscan@#1#2{%
    \unit@=\@ne bp\setc@ntr@l{2}\figsetmark{}%
    \def\minst@p{20pt}%
    \v@lX=\@psfllx\p@\v@lX=\Sc@leFact\v@lX\r@undint\v@lX\v@lX%
    \v@lY=\@psflly\p@\v@lY=\Sc@leFact\v@lY\ifdim\v@lY>\z@\r@undint\v@lY\v@lY\fi%
    \delt@=\@psfury\p@\delt@=\Sc@leFact\delt@%
    \advance\delt@-\v@lY\v@lXa=\@psfurx\p@\v@lXa=\Sc@leFact\v@lXa\v@leur=\minst@p%
    \edef\valv@lY{\repdecn@mb{\v@lY}}\edef\LgTr@it{\the\delt@}%
    \loop\ifdim\v@lX<\v@lXa\edef\valv@lX{\repdecn@mb{\v@lX}}%
    \figptDD -1:(\valv@lX,\valv@lY)\figwriten -1:\hbox{\vrule height\LgTr@it}(0)%
    \ifdim\v@leur<\minst@p\else\figsetmark{\raise-8bp\hbox{$\scriptscriptstyle\triangle$}}%
    \figwrites -1:\@ffichnb{0}{\valv@lX}(6)\v@leur=\z@\figsetmark{}\fi%
    \advance\v@leur#1pt\advance\v@lX#1pt\repeat%
    \def\minst@p{10pt}%
    \v@lX=\@psfllx\p@\v@lX=\Sc@leFact\v@lX\ifdim\v@lX>\z@\r@undint\v@lX\v@lX\fi%
    \v@lY=\@psflly\p@\v@lY=\Sc@leFact\v@lY\r@undint\v@lY\v@lY%
    \delt@=\@psfurx\p@\delt@=\Sc@leFact\delt@%
    \advance\delt@-\v@lX\v@lYa=\@psfury\p@\v@lYa=\Sc@leFact\v@lYa\v@leur=\minst@p%
    \edef\valv@lX{\repdecn@mb{\v@lX}}\edef\LgTr@it{\the\delt@}%
    \loop\ifdim\v@lY<\v@lYa\edef\valv@lY{\repdecn@mb{\v@lY}}%
    \figptDD -1:(\valv@lX,\valv@lY)\figwritee -1:\vbox{\hrule width\LgTr@it}(0)%
    \ifdim\v@leur<\minst@p\else\figsetmark{$\triangleright$\kern4bp}%
    \figwritew -1:\@ffichnb{0}{\valv@lY}(6)\v@leur=\z@\figsetmark{}\fi%
    \advance\v@leur#2pt\advance\v@lY#2pt\repeat}
\ctr@ld@f\let\figscanI=\figscan
\ctr@ld@f\def\figscan@E#1(#2,#3){{\s@uvc@ntr@l\et@tfigscan@E%
    \Figdisc@rdLTS{#1}{\t@xt@}\pdfximage{\t@xt@}%
    \setbox\Gb@x=\hbox{\pdfrefximage\pdflastximage}%
    \edef\@psfllx{0}\v@lY=-\dp\Gb@x\edef\@psflly{\repdecn@mb{\v@lY}}%
    \edef\@psfurx{\repdecn@mb{\wd\Gb@x}}%
    \v@lY=\dp\Gb@x\advance\v@lY\ht\Gb@x\edef\@psfury{\repdecn@mb{\v@lY}}%
    \figscan@{#2}{#3}\resetc@ntr@l\et@tfigscan@E}\ignorespaces}
\ctr@ld@f\def\figshowpts[#1,#2]{{\figsetmark{$\bullet$}\figsetptname{\bf ##1}%
    \p@rtent=#2\relax\ifnum\p@rtent<\z@\p@rtent=\z@\fi%
    \s@mme=#1\relax\ifnum\s@mme<\z@\s@mme=\z@\fi%
    \loop\ifnum\s@mme<\p@rtent\pt@rvect{\s@mme}%
    \ifitis@K\figwriten{\the\s@mme}:(4pt)\fi\advance\s@mme\@ne\repeat%
    \pt@rvect{\s@mme}\ifitis@K\figwriten{\the\s@mme}:(4pt)\fi}\ignorespaces}
\ctr@ld@f\def\pt@rvect#1{\set@bjc@de{#1}%
    \expandafter\expandafter\expandafter\inqpt@rvec\csname\objc@de\endcsname:}
\ctr@ld@f\def\inqpt@rvec#1#2:{\if#1\C@dCl@spt\itis@Ktrue\else\itis@Kfalse\fi}
\ctr@ld@f\def\figshowsettings{{%
    \immediate\write16{====================================================================}%
    \immediate\write16{ Current settings about:}%
    \immediate\write16{ --- GENERAL ---}%
    \immediate\write16{Scale factor and Unit = \unit@util\space (\the\unit@)
     \space -> \BS@ figinit{ScaleFactorUnit}}%
    \immediate\write16{Update mode = \ifpsupdatem@de yes\else no\fi
     \space-> \BS@ psset(update=yes/no) or \BS@ pssetdefault(update=yes/no)}%
    \immediate\write16{ --- PRINTING ---}%
    \immediate\write16{Implicit point name = \ptn@me{i} \space-> \BS@ figsetptname{Name}}%
    \immediate\write16{Point marker = \the\c@nsymb \space -> \BS@ figsetmark{Mark}}%
    \immediate\write16{Print rounded coordinates = \ifr@undcoord yes\else no\fi
     \space-> \BS@ figsetroundcoord{yes/no}}%
    \immediate\write16{ --- GRAPHICAL (general) ---}%
    \immediate\write16{First-level (or primary) settings:}%
    \immediate\write16{ Color = \curr@ntcolor \space-> \BS@ psset(color=ColorDefinition)}%
    \immediate\write16{ Filling mode = \iffillm@de yes\else no\fi
     \space-> \BS@ psset(fillmode=yes/no)}%
    \immediate\write16{ Line join = \curr@ntjoin \space-> \BS@ psset(join=miter/round/bevel)}%
    \immediate\write16{ Line style = \curr@ntdash \space-> \BS@ psset(dash=Index/Pattern)}%
    \immediate\write16{ Line width = \curr@ntwidth
     \space-> \BS@ psset(width=real in PostScript units)}%
    \immediate\write16{Second-level (or secondary) settings:}%
    \immediate\write16{ Color = \sec@ndcolor \space-> \BS@ psset second(color=ColorDefinition)}%
    \immediate\write16{ Line style = \curr@ntseconddash
     \space-> \BS@ psset second(dash=Index/Pattern)}%
    \immediate\write16{ Line width = \curr@ntsecondwidth
     \space-> \BS@ psset second(width=real in PostScript units)}%
    \immediate\write16{Third-level (or ternary) settings:}%
    \immediate\write16{ Color = \th@rdcolor \space-> \BS@ psset third(color=ColorDefinition)}%
    \immediate\write16{ --- GRAPHICAL (specific) ---}%
    \immediate\write16{Arrow-head:}%
    \immediate\write16{ (half-)Angle = \@rrowheadangle
     \space-> \BS@ psset arrowhead(angle=real in degrees)}%
    \immediate\write16{ Filling mode = \if@rrowhfill yes\else no\fi
     \space-> \BS@ psset arrowhead(fillmode=yes/no)}%
    \immediate\write16{ "Outside" = \if@rrowhout yes\else no\fi
     \space-> \BS@ psset arrowhead(out=yes/no)}%
    \immediate\write16{ Length = \@rrowheadlength
     \if@rrowratio\space(not active)\else\space(active)\fi
     \space-> \BS@ psset arrowhead(length=real in user coord.)}%
    \immediate\write16{ Ratio = \@rrowheadratio
     \if@rrowratio\space(active)\else\space(not active)\fi
     \space-> \BS@ psset arrowhead(ratio=real in [0,1])}%
    \immediate\write16{Curve: Roundness = \curv@roundness
     \space-> \BS@ psset curve(roundness=real in [0,0.5])}%
    \immediate\write16{Mesh: Diagonal = \c@ntrolmesh
     \space-> \BS@ psset mesh(diag=integer in {-1,0,1})}%
    \immediate\write16{Flow chart:}%
    \immediate\write16{ Arrow position = \@rrowp@s
     \space-> \BS@ psset flowchart(arrowposition=real in [0,1])}%
    \immediate\write16{ Arrow reference point = \ifcase\@rrowr@fpt start\else end\fi
     \space-> \BS@ psset flowchart(arrowrefpt = start/end)}%
    \immediate\write16{ Line type = \ifcase\fclin@typ@ curve\else polygon\fi
     \space-> \BS@ psset flowchart(line=polygon/curve)}%
    \immediate\write16{ Padding = (\Xp@dd, \Yp@dd)
     \space-> \BS@ psset flowchart(padding = real in user coord.)}%
    \immediate\write16{\space\space\space\space(or
     \BS@ psset flowchart(xpadding=real, ypadding=real) )}%
    \immediate\write16{ Radius = \fclin@r@d
     \space-> \BS@ psset flowchart(radius=positive real in user coord.)}%
    \immediate\write16{ Shape = \fcsh@pe
     \space-> \BS@ psset flowchart(shape = rectangle, ellipse or lozenge)}%
    \immediate\write16{ Thickness = \thickn@ss
     \space-> \BS@ psset flowchart(thickness = real in user coord.)}%
    \ifTr@isDim%
    \immediate\write16{ --- 3D to 2D PROJECTION ---}%
    \immediate\write16{Projection : \typ@proj \space-> \BS@ figinit{ScaleFactorUnit, ProjType}}%
    \immediate\write16{Longitude (psi) = \v@lPsi \space-> \BS@ figset proj(psi=real in degrees)}%
    \ifcase\curr@ntproj\immediate\write16{Depth coeff. (Lambda)
     \space = \v@lTheta \space-> \BS@ figset proj(lambda=real in [0,1])}%
    \else\immediate\write16{Latitude (theta)
     \space = \v@lTheta \space-> \BS@ figset proj(theta=real in degrees)}%
    \fi%
    \ifnum\curr@ntproj=\tw@%
    \immediate\write16{Observation distance = \disob@unit
     \space-> \BS@ figset proj(dist=real in user coord.)}%
    \immediate\write16{Target point = \t@rgetpt \space-> \BS@ figset proj(targetpt=pt number)}%
     \v@lX=\ptT@unit@\wd\Bt@rget\v@lY=\ptT@unit@\ht\Bt@rget\v@lZ=\ptT@unit@\dp\Bt@rget%
    \immediate\write16{ Its coordinates are
     (\repdecn@mb{\v@lX}, \repdecn@mb{\v@lY}, \repdecn@mb{\v@lZ})}%
    \fi%
    \fi%
    \immediate\write16{====================================================================}%
    \ignorespaces}}
\ctr@ln@w{newif}\ifitis@vect@r
\ctr@ld@f\def\figvectC#1(#2,#3){{\itis@vect@rtrue\figpt#1:(#2,#3)}\ignorespaces}
\ctr@ld@f\def\Figv@ctCreg#1(#2,#3){{\itis@vect@rtrue\Figp@intreg#1:(#2,#3)}\ignorespaces}
\ctr@ln@m\figvectDBezier
\ctr@ld@f\def\figvectDBezierDD#1:#2,#3[#4,#5,#6,#7]{\ifps@cri{\s@uvc@ntr@l\et@tfigvectDBezierDD%
    \FigvectDBezier@#2,#3[#4,#5,#6,#7]\v@lX=\c@ef\v@lX\v@lY=\c@ef\v@lY%
    \Figv@ctCreg#1(\v@lX,\v@lY)\resetc@ntr@l\et@tfigvectDBezierDD}\ignorespaces\fi}
\ctr@ld@f\def\figvectDBezierTD#1:#2,#3[#4,#5,#6,#7]{\ifps@cri{\s@uvc@ntr@l\et@tfigvectDBezierTD%
    \FigvectDBezier@#2,#3[#4,#5,#6,#7]\v@lX=\c@ef\v@lX\v@lY=\c@ef\v@lY\v@lZ=\c@ef\v@lZ%
    \Figv@ctCreg#1(\v@lX,\v@lY,\v@lZ)\resetc@ntr@l\et@tfigvectDBezierTD}\ignorespaces\fi}
\ctr@ld@f\def\FigvectDBezier@#1,#2[#3,#4,#5,#6]{\setc@ntr@l{2}%
    \edef\T@{#2}\v@leur=\p@\advance\v@leur-#2pt\edef\UNmT@{\repdecn@mb{\v@leur}}%
    \ifnum#1=\tw@\def\c@ef{6}\else\def\c@ef{3}\fi%
    \figptcopy-4:/#3/\figptcopy-3:/#4/\figptcopy-2:/#5/\figptcopy-1:/#6/%
    \l@mbd@un=-4 \l@mbd@de=-\thr@@\p@rtent=\m@ne\c@lDecast%
    \ifnum#1=\tw@\c@lDCDeux{-4}{-3}\c@lDCDeux{-3}{-2}\c@lDCDeux{-4}{-3}\else%
    \l@mbd@un=-4 \l@mbd@de=-\thr@@\p@rtent=-\tw@\c@lDecast%
    \c@lDCDeux{-4}{-3}\fi\Figg@tXY{-4}}
\ctr@ln@m\c@lDCDeux
\ctr@ld@f\def\c@lDCDeuxDD#1#2{\Figg@tXY{#2}\Figg@tXYa{#1}%
    \advance\v@lX-\v@lXa\advance\v@lY-\v@lYa\Figp@intregDD#1:(\v@lX,\v@lY)}
\ctr@ld@f\def\c@lDCDeuxTD#1#2{\Figg@tXY{#2}\Figg@tXYa{#1}\advance\v@lX-\v@lXa%
    \advance\v@lY-\v@lYa\advance\v@lZ-\v@lZa\Figp@intregTD#1:(\v@lX,\v@lY,\v@lZ)}
\ctr@ln@m\figvectN
\ctr@ld@f\def\figvectNDD#1[#2,#3]{\ifps@cri{\Figg@tXYa{#2}\Figg@tXY{#3}%
    \advance\v@lX-\v@lXa\advance\v@lY-\v@lYa%
    \Figv@ctCreg#1(-\v@lY,\v@lX)}\ignorespaces\fi}
\ctr@ld@f\def\figvectNTD#1[#2,#3,#4]{\ifps@cri{\vecunitC@TD[#2,#4]\v@lmin=\v@lX\v@lmax=\v@lY%
    \v@leur=\v@lZ\vecunitC@TD[#2,#3]\c@lprovec{#1}}\ignorespaces\fi}
\ctr@ln@m\figvectNV
\ctr@ld@f\def\figvectNVDD#1[#2]{\ifps@cri{\Figg@tXY{#2}\Figv@ctCreg#1(-\v@lY,\v@lX)}\ignorespaces\fi}
\ctr@ld@f\def\figvectNVTD#1[#2,#3]{\ifps@cri{\vecunitCV@TD{#3}\v@lmin=\v@lX\v@lmax=\v@lY%
    \v@leur=\v@lZ\vecunitCV@TD{#2}\c@lprovec{#1}}\ignorespaces\fi}
\ctr@ln@m\figvectP
\ctr@ld@f\def\figvectPDD#1[#2,#3]{\ifps@cri{\Figg@tXYa{#2}\Figg@tXY{#3}%
    \advance\v@lX-\v@lXa\advance\v@lY-\v@lYa%
    \Figv@ctCreg#1(\v@lX,\v@lY)}\ignorespaces\fi}
\ctr@ld@f\def\figvectPTD#1[#2,#3]{\ifps@cri{\Figg@tXYa{#2}\Figg@tXY{#3}%
    \advance\v@lX-\v@lXa\advance\v@lY-\v@lYa\advance\v@lZ-\v@lZa%
    \Figv@ctCreg#1(\v@lX,\v@lY,\v@lZ)}\ignorespaces\fi}
\ctr@ln@m\figvectU
\ctr@ld@f\def\figvectUDD#1[#2]{\ifps@cri{\n@rmeuc\v@leur{#2}\invers@\v@leur\v@leur%
    \delt@=\repdecn@mb{\v@leur}\unit@\edef\v@ldelt@{\repdecn@mb{\delt@}}%
    \Figg@tXY{#2}\v@lX=\v@ldelt@\v@lX\v@lY=\v@ldelt@\v@lY%
    \Figv@ctCreg#1(\v@lX,\v@lY)}\ignorespaces\fi}
\ctr@ld@f\def\figvectUTD#1[#2]{\ifps@cri{\n@rmeuc\v@leur{#2}\invers@\v@leur\v@leur%
    \delt@=\repdecn@mb{\v@leur}\unit@\edef\v@ldelt@{\repdecn@mb{\delt@}}%
    \Figg@tXY{#2}\v@lX=\v@ldelt@\v@lX\v@lY=\v@ldelt@\v@lY\v@lZ=\v@ldelt@\v@lZ%
    \Figv@ctCreg#1(\v@lX,\v@lY,\v@lZ)}\ignorespaces\fi}
\ctr@ld@f\def\figvisu#1#2#3{\c@ldefproj\initb@undb@x\xdef\figforTeXFigno{\figforTeXnextFigno}%
    \s@mme=\figforTeXnextFigno\advance\s@mme\@ne\xdef\figforTeXnextFigno{\number\s@mme}%
    \setbox\b@xvisu=\hbox{\ifnum\@utoFN>\z@\figinsert{}\gdef\@utoFInDone{0}\fi\ignorespaces#3}%
    \gdef\@utoFInDone{1}\gdef\@utoFN{0}%
    \v@lXa=-\c@@rdYmin\v@lYa=\c@@rdYmax\advance\v@lYa-\c@@rdYmin%
    \v@lX=\c@@rdXmax\advance\v@lX-\c@@rdXmin%
    \setbox#1=\hbox{#2}\v@lY=-\v@lX\maxim@m{\v@lX}{\v@lX}{\wd#1}%
    \advance\v@lY\v@lX\divide\v@lY\tw@\advance\v@lY-\c@@rdXmin%
    \setbox#1=\vbox{\parindent0mm\hsize=\v@lX\vskip\v@lYa%
    \rlap{\hskip\v@lY\smash{\raise\v@lXa\box\b@xvisu}}%
    \def\t@xt@{#2}\ifx\t@xt@\empty\else\medskip\centerline{#2}\fi}\wd#1=\v@lX}
\ctr@ld@f\def\figDecrementFigno{{\xdef\figforTeXnextFigno{\figforTeXFigno}%
    \s@mme=\figforTeXFigno\advance\s@mme\m@ne\xdef\figforTeXFigno{\number\s@mme}}}
\ctr@ln@w{newbox}\Bt@rget\setbox\Bt@rget=\null
\ctr@ln@w{newbox}\BminTD@\setbox\BminTD@=\null
\ctr@ln@w{newbox}\BmaxTD@\setbox\BmaxTD@=\null
\ctr@ln@w{newif}\ifnewt@rgetpt\ctr@ln@w{newif}\ifnewdis@b
\ctr@ld@f\def\b@undb@xTD#1#2#3{%
    \relax\ifdim#1<\wd\BminTD@\global\wd\BminTD@=#1\fi%
    \relax\ifdim#2<\ht\BminTD@\global\ht\BminTD@=#2\fi%
    \relax\ifdim#3<\dp\BminTD@\global\dp\BminTD@=#3\fi%
    \relax\ifdim#1>\wd\BmaxTD@\global\wd\BmaxTD@=#1\fi%
    \relax\ifdim#2>\ht\BmaxTD@\global\ht\BmaxTD@=#2\fi%
    \relax\ifdim#3>\dp\BmaxTD@\global\dp\BmaxTD@=#3\fi}
\ctr@ld@f\def\c@ldefdisob{{\ifdim\wd\BminTD@<\maxdimen\v@leur=\wd\BmaxTD@\advance\v@leur-\wd\BminTD@%
    \delt@=\ht\BmaxTD@\advance\delt@-\ht\BminTD@\maxim@m{\v@leur}{\v@leur}{\delt@}%
    \delt@=\dp\BmaxTD@\advance\delt@-\dp\BminTD@\maxim@m{\v@leur}{\v@leur}{\delt@}%
    \v@leur=5\v@leur\else\v@leur=800pt\fi\c@ldefdisob@{\v@leur}}}
\ctr@ln@m\disob@intern
\ctr@ln@m\disob@
\ctr@ln@m\divf@ctproj
\ctr@ld@f\def\c@ldefdisob@#1{{\v@leur=#1\ifdim\v@leur<\p@\v@leur=800pt\fi%
    \xdef\disob@intern{\repdecn@mb{\v@leur}}%
    \delt@=\ptT@unit@\v@leur\xdef\disob@unit{\repdecn@mb{\delt@}}%
    \f@ctech=\@ne\loop\ifdim\v@leur>\t@n pt\divide\v@leur\t@n\multiply\f@ctech\t@n\repeat%
    \xdef\disob@{\repdecn@mb{\v@leur}}\xdef\divf@ctproj{\the\f@ctech}}%
    \global\newdis@btrue}
\ctr@ln@m\t@rgetpt
\ctr@ld@f\def\c@ldeft@rgetpt{\newt@rgetpttrue\def\t@rgetpt{CenterBoundBox}{%
    \delt@=\wd\BmaxTD@\advance\delt@-\wd\BminTD@\divide\delt@\tw@%
    \v@leur=\wd\BminTD@\advance\v@leur\delt@\global\wd\Bt@rget=\v@leur%
    \delt@=\ht\BmaxTD@\advance\delt@-\ht\BminTD@\divide\delt@\tw@%
    \v@leur=\ht\BminTD@\advance\v@leur\delt@\global\ht\Bt@rget=\v@leur%
    \delt@=\dp\BmaxTD@\advance\delt@-\dp\BminTD@\divide\delt@\tw@%
    \v@leur=\dp\BminTD@\advance\v@leur\delt@\global\dp\Bt@rget=\v@leur}}
\ctr@ln@m\c@ldefproj
\ctr@ld@f\def\c@ldefprojTD{\ifnewt@rgetpt\else\c@ldeft@rgetpt\fi\ifnewdis@b\else\c@ldefdisob\fi}
\ctr@ld@f\def\c@lprojcav{
    \v@lZa=\cxa@\v@lY\advance\v@lX\v@lZa%
    \v@lZa=\cxb@\v@lY\v@lY=\v@lZ\advance\v@lY\v@lZa\ignorespaces}
\ctr@ln@m\v@lcoef
\ctr@ld@f\def\c@lprojrea{
    \advance\v@lX-\wd\Bt@rget\advance\v@lY-\ht\Bt@rget\advance\v@lZ-\dp\Bt@rget%
    \v@lZa=\cza@\v@lX\advance\v@lZa\czb@\v@lY\advance\v@lZa\czc@\v@lZ%
    \divide\v@lZa\divf@ctproj\advance\v@lZa\disob@ pt\invers@{\v@lZa}{\v@lZa}%
    \v@lZa=\disob@\v@lZa\edef\v@lcoef{\repdecn@mb{\v@lZa}}%
    \v@lXa=\cxa@\v@lX\advance\v@lXa\cxb@\v@lY\v@lXa=\v@lcoef\v@lXa%
    \v@lY=\cyb@\v@lY\advance\v@lY\cya@\v@lX\advance\v@lY\cyc@\v@lZ%
    \v@lY=\v@lcoef\v@lY\v@lX=\v@lXa\ignorespaces}
\ctr@ld@f\def\c@lprojort{
    \v@lXa=\cxa@\v@lX\advance\v@lXa\cxb@\v@lY%
    \v@lY=\cyb@\v@lY\advance\v@lY\cya@\v@lX\advance\v@lY\cyc@\v@lZ%
    \v@lX=\v@lXa\ignorespaces}
\ctr@ld@f\def\Figptpr@j#1:#2/#3/{{\Figg@tXY{#3}\superc@lprojSP%
    \Figp@intregDD#1:{#2}(\v@lX,\v@lY)}\ignorespaces}
\ctr@ln@m\figsetobdist
\ctr@ld@f\def\figsetobdistDD{\un@v@ilable{figsetobdist}}
\ctr@ld@f\def\figsetobdistTD(#1){{\ifcurr@ntPS%
    \immediate\write16{*** \BS@ figsetobdist is ignored inside a
     \BS@ psbeginfig-\BS@ psendfig block.}%
    \else\v@leur=#1\unit@\c@ldefdisob@{\v@leur}\fi}\ignorespaces}
\ctr@ln@m\c@lprojSP
\ctr@ln@m\curr@ntproj
\ctr@ln@m\typ@proj
\ctr@ln@m\superc@lprojSP
\ctr@ld@f\def\Figs@tproj#1{%
    \if#13 \d@faultproj\else\if#1c\d@faultproj%
    \else\if#1o\xdef\curr@ntproj{1}\xdef\typ@proj{orthogonal}%
         \figsetviewTD(\def@ultpsi,\def@ulttheta)%
         \global\let\c@lprojSP=\c@lprojort\global\let\superc@lprojSP=\c@lprojort%
    \else\if#1r\xdef\curr@ntproj{2}\xdef\typ@proj{realistic}%
         \figsetviewTD(\def@ultpsi,\def@ulttheta)%
         \global\let\c@lprojSP=\c@lprojrea\global\let\superc@lprojSP=\c@lprojrea%
    \else\d@faultproj\message{*** Unknown projection. Cavalier projection assumed.}%
    \fi\fi\fi\fi}
\ctr@ld@f\def\d@faultproj{\xdef\curr@ntproj{0}\xdef\typ@proj{cavalier}\figsetviewTD(\def@ultpsi,0.5)%
         \global\let\c@lprojSP=\c@lprojcav\global\let\superc@lprojSP=\c@lprojcav}
\ctr@ln@m\figsettarget
\ctr@ld@f\def\figsettargetDD{\un@v@ilable{figsettarget}}
\ctr@ld@f\def\figsettargetTD[#1]{{\ifcurr@ntPS%
    \immediate\write16{*** \BS@ figsettarget is ignored inside a
     \BS@ psbeginfig-\BS@ psendfig block.}%
    \else\global\newt@rgetpttrue\xdef\t@rgetpt{#1}\Figg@tXY{#1}\global\wd\Bt@rget=\v@lX%
    \global\ht\Bt@rget=\v@lY\global\dp\Bt@rget=\v@lZ\fi}\ignorespaces}
\ctr@ln@m\figsetview
\ctr@ld@f\def\figsetviewDD{\un@v@ilable{figsetview}}
\ctr@ld@f\def\figsetviewTD(#1){\ifcurr@ntPS%
     \immediate\write16{*** \BS@ figsetview is ignored inside a
     \BS@ psbeginfig-\BS@ psendfig block.}\else\Figsetview@#1,:\fi\ignorespaces}
\ctr@ld@f\def\Figsetview@#1,#2:{{\xdef\v@lPsi{#1}\def\t@xt@{#2}%
    \ifx\t@xt@\empty\def\@rgdeux{\v@lTheta}\else\X@rgdeux@#2\fi%
    \c@ssin{\costhet@}{\sinthet@}{#1}\v@lmin=\costhet@ pt\v@lmax=\sinthet@ pt%
    \ifcase\curr@ntproj%
    \v@leur=\@rgdeux\v@lmin\xdef\cxa@{\repdecn@mb{\v@leur}}%
    \v@leur=\@rgdeux\v@lmax\xdef\cxb@{\repdecn@mb{\v@leur}}\v@leur=\@rgdeux pt%
    \relax\ifdim\v@leur>\p@\message{*** Lambda too large ! See \BS@ figset proj() !}\fi%
    \else%
    \v@lmax=-\v@lmax\xdef\cxa@{\repdecn@mb{\v@lmax}}\xdef\cxb@{\costhet@}%
    \ifx\t@xt@\empty\edef\@rgdeux{\def@ulttheta}\fi\c@ssin{\C@}{\S@}{\@rgdeux}%
    \v@lmax=-\S@ pt%
    \v@leur=\v@lmax\v@leur=\costhet@\v@leur\xdef\cya@{\repdecn@mb{\v@leur}}%
    \v@leur=\v@lmax\v@leur=\sinthet@\v@leur\xdef\cyb@{\repdecn@mb{\v@leur}}%
    \xdef\cyc@{\C@}\v@lmin=-\C@ pt%
    \v@leur=\v@lmin\v@leur=\costhet@\v@leur\xdef\cza@{\repdecn@mb{\v@leur}}%
    \v@leur=\v@lmin\v@leur=\sinthet@\v@leur\xdef\czb@{\repdecn@mb{\v@leur}}%
    \xdef\czc@{\repdecn@mb{\v@lmax}}\fi%
    \xdef\v@lTheta{\@rgdeux}}}
\ctr@ld@f\def\def@ultpsi{40}
\ctr@ld@f\def\def@ulttheta{25}
\ctr@ln@m\l@debut
\ctr@ln@m\n@mref
\ctr@ld@f\def\figset#1(#2){\def\t@xt@{#1}\ifx\t@xt@\empty\trtlis@rg{#2}{\Figsetwr@te}
    \else\keln@mde#1|%
    \def\n@mref{pr}\ifx\l@debut\n@mref\ifcurr@ntPS
     \immediate\write16{*** \BS@ figset proj(...) is ignored inside a
     \BS@ psbeginfig-\BS@ psendfig block.}\else\trtlis@rg{#2}{\Figsetpr@j}\fi\else%
    \def\n@mref{wr}\ifx\l@debut\n@mref\trtlis@rg{#2}{\Figsetwr@te}\else
    \immediate\write16{*** Unknown keyword: \BS@ figset #1(...)}%
    \fi\fi\fi\ignorespaces}
\ctr@ld@f\def\Figsetpr@j#1=#2|{\keln@mtr#1|%
    \def\n@mref{dep}\ifx\l@debut\n@mref\Figsetd@p{#2}\else
    \def\n@mref{dis}\ifx\l@debut\n@mref%
     \ifnum\curr@ntproj=\tw@\figsetobdist(#2)\else\Figset@rr\fi\else
    \def\n@mref{lam}\ifx\l@debut\n@mref\Figsetd@p{#2}\else
    \def\n@mref{lat}\ifx\l@debut\n@mref\Figsetth@{#2}\else
    \def\n@mref{lon}\ifx\l@debut\n@mref\figsetview(#2)\else
    \def\n@mref{psi}\ifx\l@debut\n@mref\figsetview(#2)\else
    \def\n@mref{tar}\ifx\l@debut\n@mref%
     \ifnum\curr@ntproj=\tw@\figsettarget[#2]\else\Figset@rr\fi\else
    \def\n@mref{the}\ifx\l@debut\n@mref\Figsetth@{#2}\else
    \immediate\write16{*** Unknown attribute: \BS@ figset proj(..., #1=...).}%
    \fi\fi\fi\fi\fi\fi\fi\fi}
\ctr@ld@f\def\Figsetd@p#1{\ifnum\curr@ntproj=\z@\figsetview(\v@lPsi,#1)\else\Figset@rr\fi}
\ctr@ld@f\def\Figsetth@#1{\ifnum\curr@ntproj=\z@\Figset@rr\else\figsetview(\v@lPsi,#1)\fi}
\ctr@ld@f\def\Figset@rr{\message{*** \BS@ figset proj(): Attribute "\n@mref" ignored, incompatible
    with current projection}}
\ctr@ld@f\def\initb@undb@xTD{\wd\BminTD@=\maxdimen\ht\BminTD@=\maxdimen\dp\BminTD@=\maxdimen%
    \wd\BmaxTD@=-\maxdimen\ht\BmaxTD@=-\maxdimen\dp\BmaxTD@=-\maxdimen}
\ctr@ln@w{newbox}\Gb@x      
\ctr@ln@w{newbox}\Gb@xSC    
\ctr@ln@w{newtoks}\c@nsymb  
\ctr@ln@w{newif}\ifr@undcoord\ctr@ln@w{newif}\ifunitpr@sent
\ctr@ld@f\def\unssqrttw@{0.707106 }
\ctr@ld@f\def\figAst{\raise-1.15ex\hbox{$\ast$}}
\ctr@ld@f\def\figBullet{\raise-1.15ex\hbox{$\bullet$}}
\ctr@ld@f\def\figCirc{\raise-1.15ex\hbox{$\circ$}}
\ctr@ld@f\def\figDiamond{\raise-1.15ex\hbox{$\diamond$}}%
\ctr@ld@f\def\boxit#1#2{\leavevmode\hbox{\vrule\vbox{\hrule\vglue#1%
    \vtop{\hbox{\kern#1{#2}\kern#1}\vglue#1\hrule}}\vrule}}
\ctr@ld@f\def\centertext#1#2{\vbox{\hsize#1\parindent0cm%
    \leftskip=0pt plus 1fil\rightskip=0pt plus 1fil\parfillskip=0pt{#2}}}
\ctr@ld@f\def\lefttext#1#2{\vbox{\hsize#1\parindent0cm\rightskip=0pt plus 1fil#2}}
\ctr@ld@f\def\c@nterpt{\ignorespaces%
    \kern-.5\wd\Gb@xSC%
    \raise-.5\ht\Gb@xSC\rlap{\hbox{\raise.5\dp\Gb@xSC\hbox{\copy\Gb@xSC}}}%
    \kern .5\wd\Gb@xSC\ignorespaces}
\ctr@ld@f\def\b@undb@xSC#1#2{{\v@lXa=#1\v@lYa=#2%
    \v@leur=\ht\Gb@xSC\advance\v@leur\dp\Gb@xSC%
    \advance\v@lXa-.5\wd\Gb@xSC\advance\v@lYa-.5\v@leur\b@undb@x{\v@lXa}{\v@lYa}%
    \advance\v@lXa\wd\Gb@xSC\advance\v@lYa\v@leur\b@undb@x{\v@lXa}{\v@lYa}}}
\ctr@ln@m\Dist@n
\ctr@ln@m\l@suite
\ctr@ld@f\def\@keldist#1#2{\edef\Dist@n{#2}\y@tiunit{\Dist@n}%
    \ifunitpr@sent#1=\Dist@n\else#1=\Dist@n\unit@\fi}
\ctr@ld@f\def\y@tiunit#1{\unitpr@sentfalse\expandafter\y@tiunit@#1:}
\ctr@ld@f\def\y@tiunit@#1#2:{\ifcat#1a\unitpr@senttrue\else\def\l@suite{#2}%
    \ifx\l@suite\empty\else\y@tiunit@#2:\fi\fi}
\ctr@ln@m\figcoord
\ctr@ld@f\def\figcoordDD#1{{\v@lX=\ptT@unit@\v@lX\v@lY=\ptT@unit@\v@lY%
    \ifr@undcoord\ifcase#1\v@leur=0.5pt\or\v@leur=0.05pt\or\v@leur=0.005pt%
    \or\v@leur=0.0005pt\else\v@leur=\z@\fi%
    \ifdim\v@lX<\z@\advance\v@lX-\v@leur\else\advance\v@lX\v@leur\fi%
    \ifdim\v@lY<\z@\advance\v@lY-\v@leur\else\advance\v@lY\v@leur\fi\fi%
    (\@ffichnb{#1}{\repdecn@mb{\v@lX}},\ifmmode\else\thinspace\fi%
    \@ffichnb{#1}{\repdecn@mb{\v@lY}})}}
\ctr@ld@f\def\@ffichnb#1#2{{\def\@@ffich{\@ffich#1(}\edef\n@mbre{#2}%
    \expandafter\@@ffich\n@mbre)}}
\ctr@ld@f\def\@ffich#1(#2.#3){{#2\ifnum#1>\z@.\fi\def\dig@ts{#3}\s@mme=\z@%
    \loop\ifnum\s@mme<#1\expandafter\@ffichdec\dig@ts:\advance\s@mme\@ne\repeat}}
\ctr@ld@f\def\@ffichdec#1#2:{\relax#1\def\dig@ts{#20}}
\ctr@ld@f\def\figcoordTD#1{{\v@lX=\ptT@unit@\v@lX\v@lY=\ptT@unit@\v@lY\v@lZ=\ptT@unit@\v@lZ%
    \ifr@undcoord\ifcase#1\v@leur=0.5pt\or\v@leur=0.05pt\or\v@leur=0.005pt%
    \or\v@leur=0.0005pt\else\v@leur=\z@\fi%
    \ifdim\v@lX<\z@\advance\v@lX-\v@leur\else\advance\v@lX\v@leur\fi%
    \ifdim\v@lY<\z@\advance\v@lY-\v@leur\else\advance\v@lY\v@leur\fi%
    \ifdim\v@lZ<\z@\advance\v@lZ-\v@leur\else\advance\v@lZ\v@leur\fi\fi%
    (\@ffichnb{#1}{\repdecn@mb{\v@lX}},\ifmmode\else\thinspace\fi%
     \@ffichnb{#1}{\repdecn@mb{\v@lY}},\ifmmode\else\thinspace\fi%
     \@ffichnb{#1}{\repdecn@mb{\v@lZ}})}}
\ctr@ld@f\def\figsetroundcoord#1{\expandafter\Figsetr@undcoord#1:\ignorespaces}
\ctr@ld@f\def\Figsetr@undcoord#1#2:{\if#1n\r@undcoordfalse\else\r@undcoordtrue\fi}
\ctr@ld@f\def\Figsetwr@te#1=#2|{\keln@mun#1|%
    \def\n@mref{m}\ifx\l@debut\n@mref\figsetmark{#2}\else
    \immediate\write16{*** Unknown attribute: \BS@ figset (..., #1=...)}%
    \fi}
\ctr@ld@f\def\figsetmark#1{\c@nsymb={#1}\setbox\Gb@xSC=\hbox{\the\c@nsymb}\ignorespaces}
\ctr@ln@m\ptn@me
\ctr@ld@f\def\figsetptname#1{\def\ptn@me##1{#1}\ignorespaces}
\ctr@ld@f\def\FigWrit@L#1:#2(#3,#4){\ignorespaces\@keldist\v@leur{#3}\@keldist\delt@{#4}%
    \C@rp@r@m\def\list@num{#1}\@ecfor\p@int:=\list@num\do{\FigWrit@pt{\p@int}{#2}}}
\ctr@ld@f\def\FigWrit@pt#1#2{\FigWp@r@m{#1}{#2}\Vc@rrect\figWp@si%
    \ifdim\wd\Gb@xSC>\z@\b@undb@xSC{\v@lX}{\v@lY}\fi\figWBB@x}
\ctr@ld@f\def\FigWp@r@m#1#2{\Figg@tXY{#1}%
    \setbox\Gb@x=\hbox{\def\t@xt@{#2}\ifx\t@xt@\empty\Figg@tT{#1}\else#2\fi}\c@lprojSP}
\ctr@ld@f\let\Vc@rrect=\relax
\ctr@ld@f\let\C@rp@r@m=\relax
\ctr@ld@f\def\figwrite[#1]#2{{\ignorespaces\def\list@num{#1}\@ecfor\p@int:=\list@num\do{%
    \setbox\Gb@x=\hbox{\def\t@xt@{#2}\ifx\t@xt@\empty\Figg@tT{\p@int}\else#2\fi}%
    \Figwrit@{\p@int}}}\ignorespaces}
\ctr@ld@f\def\Figwrit@#1{\Figg@tXY{#1}\c@lprojSP%
    \rlap{\kern\v@lX\raise\v@lY\hbox{\unhcopy\Gb@x}}\v@leur=\v@lY%
    \advance\v@lY\ht\Gb@x\b@undb@x{\v@lX}{\v@lY}\advance\v@lX\wd\Gb@x%
    \v@lY=\v@leur\advance\v@lY-\dp\Gb@x\b@undb@x{\v@lX}{\v@lY}}
\ctr@ld@f\def\figwritec[#1]#2{{\ignorespaces\def\list@num{#1}%
    \@ecfor\p@int:=\list@num\do{\Figwrit@c{\p@int}{#2}}}\ignorespaces}
\ctr@ld@f\def\Figwrit@c#1#2{\FigWp@r@m{#1}{#2}%
    \rlap{\kern\v@lX\raise\v@lY\hbox{\rlap{\kern-.5\wd\Gb@x%
    \raise-.5\ht\Gb@x\hbox{\raise.5\dp\Gb@x\hbox{\unhcopy\Gb@x}}}}}%
    \v@leur=\ht\Gb@x\advance\v@leur\dp\Gb@x%
    \advance\v@lX-.5\wd\Gb@x\advance\v@lY-.5\v@leur\b@undb@x{\v@lX}{\v@lY}%
    \advance\v@lX\wd\Gb@x\advance\v@lY\v@leur\b@undb@x{\v@lX}{\v@lY}}
\ctr@ld@f\def\figwritep[#1]{{\ignorespaces\def\list@num{#1}\setbox\Gb@x=\hbox{\c@nterpt}%
    \@ecfor\p@int:=\list@num\do{\Figwrit@{\p@int}}}\ignorespaces}
\ctr@ld@f\def\figwritew#1:#2(#3){\figwritegcw#1:{#2}(#3,0pt)}
\ctr@ld@f\def\figwritee#1:#2(#3){\figwritegce#1:{#2}(#3,0pt)}
\ctr@ld@f\def\figwriten#1:#2(#3){{\def\Vc@rrect{\v@lZ=\v@leur\advance\v@lZ\dp\Gb@x}%
    \Figwrit@NS#1:{#2}(#3)}\ignorespaces}
\ctr@ld@f\def\figwrites#1:#2(#3){{\def\Vc@rrect{\v@lZ=-\v@leur\advance\v@lZ-\ht\Gb@x}%
    \Figwrit@NS#1:{#2}(#3)}\ignorespaces}
\ctr@ld@f\def\Figwrit@NS#1:#2(#3){\let\figWp@si=\FigWp@siNS\let\figWBB@x=\FigWBB@xNS%
    \FigWrit@L#1:{#2}(#3,0pt)}
\ctr@ld@f\def\FigWp@siNS{\rlap{\kern\v@lX\raise\v@lY\hbox{\rlap{\kern-.5\wd\Gb@x%
    \raise\v@lZ\hbox{\unhcopy\Gb@x}}\c@nterpt}}}
\ctr@ld@f\def\FigWBB@xNS{\advance\v@lY\v@lZ%
    \advance\v@lY-\dp\Gb@x\advance\v@lX-.5\wd\Gb@x\b@undb@x{\v@lX}{\v@lY}%
    \advance\v@lY\ht\Gb@x\advance\v@lY\dp\Gb@x%
    \advance\v@lX\wd\Gb@x\b@undb@x{\v@lX}{\v@lY}}
\ctr@ld@f\def\figwritenw#1:#2(#3){{\let\figWp@si=\FigWp@sigW\let\figWBB@x=\FigWBB@xgWE%
    \def\C@rp@r@m{\v@leur=\unssqrttw@\v@leur\delt@=\v@leur%
    \ifdim\delt@=\z@\delt@=\epsil@n\fi}\let@xte={-}\FigWrit@L#1:{#2}(#3,0pt)}\ignorespaces}
\ctr@ld@f\def\figwritesw#1:#2(#3){{\let\figWp@si=\FigWp@sigW\let\figWBB@x=\FigWBB@xgWE%
    \def\C@rp@r@m{\v@leur=\unssqrttw@\v@leur\delt@=-\v@leur%
    \ifdim\delt@=\z@\delt@=-\epsil@n\fi}\let@xte={-}\FigWrit@L#1:{#2}(#3,0pt)}\ignorespaces}
\ctr@ld@f\def\figwritene#1:#2(#3){{\let\figWp@si=\FigWp@sigE\let\figWBB@x=\FigWBB@xgWE%
    \def\C@rp@r@m{\v@leur=\unssqrttw@\v@leur\delt@=\v@leur%
    \ifdim\delt@=\z@\delt@=\epsil@n\fi}\let@xte={}\FigWrit@L#1:{#2}(#3,0pt)}\ignorespaces}
\ctr@ld@f\def\figwritese#1:#2(#3){{\let\figWp@si=\FigWp@sigE\let\figWBB@x=\FigWBB@xgWE%
    \def\C@rp@r@m{\v@leur=\unssqrttw@\v@leur\delt@=-\v@leur%
    \ifdim\delt@=\z@\delt@=-\epsil@n\fi}\let@xte={}\FigWrit@L#1:{#2}(#3,0pt)}\ignorespaces}
\ctr@ld@f\def\figwritegw#1:#2(#3,#4){{\let\figWp@si=\FigWp@sigW\let\figWBB@x=\FigWBB@xgWE%
    \let@xte={-}\FigWrit@L#1:{#2}(#3,#4)}\ignorespaces}
\ctr@ld@f\def\figwritege#1:#2(#3,#4){{\let\figWp@si=\FigWp@sigE\let\figWBB@x=\FigWBB@xgWE%
    \let@xte={}\FigWrit@L#1:{#2}(#3,#4)}\ignorespaces}
\ctr@ld@f\def\FigWp@sigW{\v@lXa=\z@\v@lYa=\ht\Gb@x\advance\v@lYa\dp\Gb@x%
    \ifdim\delt@>\z@\relax%
    \rlap{\kern\v@lX\raise\v@lY\hbox{\rlap{\kern-\wd\Gb@x\kern-\v@leur%
          \raise\delt@\hbox{\raise\dp\Gb@x\hbox{\unhcopy\Gb@x}}}\c@nterpt}}%
    \else\ifdim\delt@<\z@\relax\v@lYa=-\v@lYa%
    \rlap{\kern\v@lX\raise\v@lY\hbox{\rlap{\kern-\wd\Gb@x\kern-\v@leur%
          \raise\delt@\hbox{\raise-\ht\Gb@x\hbox{\unhcopy\Gb@x}}}\c@nterpt}}%
    \else\v@lXa=-.5\v@lYa%
    \rlap{\kern\v@lX\raise\v@lY\hbox{\rlap{\kern-\wd\Gb@x\kern-\v@leur%
          \raise-.5\ht\Gb@x\hbox{\raise.5\dp\Gb@x\hbox{\unhcopy\Gb@x}}}\c@nterpt}}%
    \fi\fi}
\ctr@ld@f\def\FigWp@sigE{\v@lXa=\z@\v@lYa=\ht\Gb@x\advance\v@lYa\dp\Gb@x%
    \ifdim\delt@>\z@\relax%
    \rlap{\kern\v@lX\raise\v@lY\hbox{\c@nterpt\kern\v@leur%
          \raise\delt@\hbox{\raise\dp\Gb@x\hbox{\unhcopy\Gb@x}}}}%
    \else\ifdim\delt@<\z@\relax\v@lYa=-\v@lYa%
    \rlap{\kern\v@lX\raise\v@lY\hbox{\c@nterpt\kern\v@leur%
          \raise\delt@\hbox{\raise-\ht\Gb@x\hbox{\unhcopy\Gb@x}}}}%
    \else\v@lXa=-.5\v@lYa%
    \rlap{\kern\v@lX\raise\v@lY\hbox{\c@nterpt\kern\v@leur%
          \raise-.5\ht\Gb@x\hbox{\raise.5\dp\Gb@x\hbox{\unhcopy\Gb@x}}}}%
    \fi\fi}
\ctr@ld@f\def\FigWBB@xgWE{\advance\v@lY\delt@%
    \advance\v@lX\the\let@xte\v@leur\advance\v@lY\v@lXa\b@undb@x{\v@lX}{\v@lY}%
    \advance\v@lX\the\let@xte\wd\Gb@x\advance\v@lY\v@lYa\b@undb@x{\v@lX}{\v@lY}}
\ctr@ld@f\def\figwritegcw#1:#2(#3,#4){{\let\figWp@si=\FigWp@sigcW\let\figWBB@x=\FigWBB@xgcWE%
    \let@xte={-}\FigWrit@L#1:{#2}(#3,#4)}\ignorespaces}
\ctr@ld@f\def\figwritegce#1:#2(#3,#4){{\let\figWp@si=\FigWp@sigcE\let\figWBB@x=\FigWBB@xgcWE%
    \let@xte={}\FigWrit@L#1:{#2}(#3,#4)}\ignorespaces}
\ctr@ld@f\def\FigWp@sigcW{\rlap{\kern\v@lX\raise\v@lY\hbox{\rlap{\kern-\wd\Gb@x\kern-\v@leur%
     \raise-.5\ht\Gb@x\hbox{\raise\delt@\hbox{\raise.5\dp\Gb@x\hbox{\unhcopy\Gb@x}}}}%
     \c@nterpt}}}
\ctr@ld@f\def\FigWp@sigcE{\rlap{\kern\v@lX\raise\v@lY\hbox{\c@nterpt\kern\v@leur%
    \raise-.5\ht\Gb@x\hbox{\raise\delt@\hbox{\raise.5\dp\Gb@x\hbox{\unhcopy\Gb@x}}}}}}
\ctr@ld@f\def\FigWBB@xgcWE{\v@lZ=\ht\Gb@x\advance\v@lZ\dp\Gb@x%
    \advance\v@lX\the\let@xte\v@leur\advance\v@lY\delt@\advance\v@lY.5\v@lZ%
    \b@undb@x{\v@lX}{\v@lY}%
    \advance\v@lX\the\let@xte\wd\Gb@x\advance\v@lY-\v@lZ\b@undb@x{\v@lX}{\v@lY}}
\ctr@ld@f\def\figwritebn#1:#2(#3){{\def\Vc@rrect{\v@lZ=\v@leur}\Figwrit@NS#1:{#2}(#3)}\ignorespaces}
\ctr@ld@f\def\figwritebs#1:#2(#3){{\def\Vc@rrect{\v@lZ=-\v@leur}\Figwrit@NS#1:{#2}(#3)}\ignorespaces}
\ctr@ld@f\def\figwritebw#1:#2(#3){{\let\figWp@si=\FigWp@sibW\let\figWBB@x=\FigWBB@xbWE%
    \let@xte={-}\FigWrit@L#1:{#2}(#3,0pt)}\ignorespaces}
\ctr@ld@f\def\figwritebe#1:#2(#3){{\let\figWp@si=\FigWp@sibE\let\figWBB@x=\FigWBB@xbWE%
    \let@xte={}\FigWrit@L#1:{#2}(#3,0pt)}\ignorespaces}
\ctr@ld@f\def\FigWp@sibW{\rlap{\kern\v@lX\raise\v@lY\hbox{\rlap{\kern-\wd\Gb@x\kern-\v@leur%
          \hbox{\unhcopy\Gb@x}}\c@nterpt}}}
\ctr@ld@f\def\FigWp@sibE{\rlap{\kern\v@lX\raise\v@lY\hbox{\c@nterpt\kern\v@leur%
          \hbox{\unhcopy\Gb@x}}}}
\ctr@ld@f\def\FigWBB@xbWE{\v@lZ=\ht\Gb@x\advance\v@lZ\dp\Gb@x%
    \advance\v@lX\the\let@xte\v@leur\advance\v@lY\ht\Gb@x\b@undb@x{\v@lX}{\v@lY}%
    \advance\v@lX\the\let@xte\wd\Gb@x\advance\v@lY-\v@lZ\b@undb@x{\v@lX}{\v@lY}}
\ctr@ln@w{newread}\frf@g  \ctr@ln@w{newwrite}\fwf@g
\ctr@ln@w{newif}\ifcurr@ntPS
\ctr@ln@w{newif}\ifps@cri
\ctr@ln@w{newif}\ifUse@llipse
\ctr@ln@w{newif}\ifpsdebugmode \psdebugmodefalse
\ctr@ln@w{newif}\ifPDFm@ke
\ifx\pdfliteral\undefined\else\ifnum\pdfoutput>\z@\PDFm@ketrue\fi\fi
\ctr@ld@f\def\initPDF@rDVI{%
\ifPDFm@ke
 \let\figscan=\figscan@E
 \let\newGr@FN=\newGr@FNPDF
 \ctr@ld@f\def\c@mcurveto{c}
 \ctr@ld@f\def\c@mfill{f}
 \ctr@ld@f\def\c@mgsave{q}
 \ctr@ld@f\def\c@mgrestore{Q}
 \ctr@ld@f\def\c@mlineto{l}
 \ctr@ld@f\def\c@mmoveto{m}
 \ctr@ld@f\def\c@msetgray{g}     \ctr@ld@f\def\c@msetgrayStroke{G}
 \ctr@ld@f\def\c@msetcmykcolor{k}\ctr@ld@f\def\c@msetcmykcolorStroke{K}
 \ctr@ld@f\def\c@msetrgbcolor{rg}\ctr@ld@f\def\c@msetrgbcolorStroke{RG}
 \ctr@ld@f\def\d@fprimarC@lor{\curr@ntcolor\space\curr@ntcolorc@md%
               \space\curr@ntcolor\space\curr@ntcolorc@mdStroke}
 \ctr@ld@f\def\d@fsecondC@lor{\sec@ndcolor\space\sec@ndcolorc@md%
               \space\sec@ndcolor\space\sec@ndcolorc@mdStroke}
 \ctr@ld@f\def\d@fthirdC@lor{\th@rdcolor\space\th@rdcolorc@md%
              \space\th@rdcolor\space\th@rdcolorc@mdStroke}
 \ctr@ld@f\def\c@msetdash{d}
 \ctr@ld@f\def\c@msetlinejoin{j}
 \ctr@ld@f\def\c@msetlinewidth{w}
 \ctr@ld@f\def\f@gclosestroke{\immediate\write\fwf@g{s}}
 \ctr@ld@f\def\f@gfill{\immediate\write\fwf@g{\fillc@md}}
 \ctr@ld@f\def\f@gnewpath{}
 \ctr@ld@f\def\f@gstroke{\immediate\write\fwf@g{S}}
\else
 \let\figinsertE=\figinsert
 \let\newGr@FN=\newGr@FNDVI
 \ctr@ld@f\def\c@mcurveto{curveto}
 \ctr@ld@f\def\c@mfill{fill}
 \ctr@ld@f\def\c@mgsave{gsave}
 \ctr@ld@f\def\c@mgrestore{grestore}
 \ctr@ld@f\def\c@mlineto{lineto}
 \ctr@ld@f\def\c@mmoveto{moveto}
 \ctr@ld@f\def\c@msetgray{setgray}          \ctr@ld@f\def\c@msetgrayStroke{}
 \ctr@ld@f\def\c@msetcmykcolor{setcmykcolor}\ctr@ld@f\def\c@msetcmykcolorStroke{}
 \ctr@ld@f\def\c@msetrgbcolor{setrgbcolor}  \ctr@ld@f\def\c@msetrgbcolorStroke{}
 \ctr@ld@f\def\d@fprimarC@lor{\curr@ntcolor\space\curr@ntcolorc@md}
 \ctr@ld@f\def\d@fsecondC@lor{\sec@ndcolor\space\sec@ndcolorc@md}
 \ctr@ld@f\def\d@fthirdC@lor{\th@rdcolor\space\th@rdcolorc@md}
 \ctr@ld@f\def\c@msetdash{setdash}
 \ctr@ld@f\def\c@msetlinejoin{setlinejoin}
 \ctr@ld@f\def\c@msetlinewidth{setlinewidth}
 \ctr@ld@f\def\f@gclosestroke{\immediate\write\fwf@g{closepath\space stroke}}
 \ctr@ld@f\def\f@gfill{\immediate\write\fwf@g{\fillc@md}}
 \ctr@ld@f\def\f@gnewpath{\immediate\write\fwf@g{newpath}}
 \ctr@ld@f\def\f@gstroke{\immediate\write\fwf@g{stroke}}
\fi}
\ctr@ld@f\def\c@pypsfile#1#2{\c@pyfil@{\immediate\write#1}{#2}}
\ctr@ld@f\def\Figinclud@PDF#1#2{\openin\frf@g=#1\pdfliteral{q #2 0 0 #2 0 0 cm}%
    \c@pyfil@{\pdfliteral}{\frf@g}\pdfliteral{Q}\closein\frf@g}
\ctr@ln@w{newif}\ifmored@ta
\ctr@ln@m\bl@nkline
\ctr@ld@f\def\c@pyfil@#1#2{\def\bl@nkline{\par}{\catcode`\%=12
    \loop\ifeof#2\mored@tafalse\else\mored@tatrue\immediate\read#2 to\tr@c
    \ifx\tr@c\bl@nkline\else#1{\tr@c}\fi\fi\ifmored@ta\repeat}}
\ctr@ld@f\def\keln@mun#1#2|{\def\l@debut{#1}\def\l@suite{#2}}
\ctr@ld@f\def\keln@mde#1#2#3|{\def\l@debut{#1#2}\def\l@suite{#3}}
\ctr@ld@f\def\keln@mtr#1#2#3#4|{\def\l@debut{#1#2#3}\def\l@suite{#4}}
\ctr@ld@f\def\keln@mqu#1#2#3#4#5|{\def\l@debut{#1#2#3#4}\def\l@suite{#5}}
\ctr@ld@f\let\@psffilein=\frf@g 
\ctr@ln@w{newif}\if@psffileok    
\ctr@ln@w{newif}\if@psfbbfound   
\ctr@ln@w{newif}\if@psfverbose   
\@psfverbosetrue
\ctr@ln@m\@psfllx \ctr@ln@m\@psflly
\ctr@ln@m\@psfurx \ctr@ln@m\@psfury
\ctr@ln@m\resetcolonc@tcode
\ctr@ld@f\def\@psfgetbb#1{\global\@psfbbfoundfalse%
\global\def\@psfllx{0}\global\def\@psflly{0}%
\global\def\@psfurx{30}\global\def\@psfury{30}%
\openin\@psffilein=#1\relax
\ifeof\@psffilein\errmessage{I couldn't open #1, will ignore it}\else
   \edef\resetcolonc@tcode{\catcode`\noexpand\:\the\catcode`\:\relax}%
   {\@psffileoktrue \chardef\other=12
    \def\do##1{\catcode`##1=\other}\dospecials \catcode`\ =10 \resetcolonc@tcode
    \loop
       \read\@psffilein to \@psffileline
       \ifeof\@psffilein\@psffileokfalse\else
          \expandafter\@psfaux\@psffileline:. \\%
       \fi
   \if@psffileok\repeat
   \if@psfbbfound\else
    \if@psfverbose\message{No bounding box comment in #1; using defaults}\fi\fi
   }\closein\@psffilein\fi}%
\ctr@ln@m\@psfbblit
\ctr@ln@m\@psfpercent
{\catcode`\%=12 \global\let\@psfpercent=
\ctr@ln@m\@psfaux
\long\def\@psfaux#1#2:#3\\{\ifx#1\@psfpercent
   \def\testit{#2}\ifx\testit\@psfbblit
      \@psfgrab #3 . . . \\%
      \@psffileokfalse
      \global\@psfbbfoundtrue
   \fi\else\ifx#1\par\else\@psffileokfalse\fi\fi}%
\ctr@ld@f\def\@psfempty{}%
\ctr@ld@f\def\@psfgrab #1 #2 #3 #4 #5\\{%
\global\def\@psfllx{#1}\ifx\@psfllx\@psfempty
      \@psfgrab #2 #3 #4 #5 .\\\else
   \global\def\@psflly{#2}%
   \global\def\@psfurx{#3}\global\def\@psfury{#4}\fi}%
\ctr@ld@f\def\PSwrit@cmd#1#2#3{{\Figg@tXY{#1}\c@lprojSP\b@undb@x{\v@lX}{\v@lY}%
    \v@lX=\ptT@ptps\v@lX\v@lY=\ptT@ptps\v@lY%
    \immediate\write#3{\repdecn@mb{\v@lX}\space\repdecn@mb{\v@lY}\space#2}}}
\ctr@ld@f\def\PSwrit@cmdS#1#2#3#4#5{{\Figg@tXY{#1}\c@lprojSP\b@undb@x{\v@lX}{\v@lY}%
    \global\result@t=\v@lX\global\result@@t=\v@lY%
    \v@lX=\ptT@ptps\v@lX\v@lY=\ptT@ptps\v@lY%
    \immediate\write#3{\repdecn@mb{\v@lX}\space\repdecn@mb{\v@lY}\space#2}}%
    \edef#4{\the\result@t}\edef#5{\the\result@@t}}
\ctr@ld@f\def\psaltitude#1[#2,#3,#4]{{\ifcurr@ntPS\ifps@cri%
    \PSc@mment{psaltitude Square Dim=#1, Triangle=[#2 / #3,#4]}%
    \s@uvc@ntr@l\et@tpsaltitude\resetc@ntr@l{2}\figptorthoprojline-5:=#2/#3,#4/%
    \figvectP -1[#3,#4]\n@rminf{\v@leur}{-1}\vecunit@{-3}{-1}%
    \figvectP -1[-5,#3]\n@rminf{\v@lmin}{-1}\figvectP -2[-5,#4]\n@rminf{\v@lmax}{-2}%
    \ifdim\v@lmin<\v@lmax\s@mme=#3\else\v@lmax=\v@lmin\s@mme=#4\fi%
    \figvectP -4[-5,#2]\vecunit@{-4}{-4}\delt@=#1\unit@%
    \edef\t@ille{\repdecn@mb{\delt@}}\figpttra-1:=-5/\t@ille,-3/%
    \figptstra-3=-5,-1/\t@ille,-4/\psline[#2,-5]\s@uvdash{\typ@dash}%
    \pssetdash{\defaultdash}\psline[-1,-2,-3]\pssetdash{\typ@dash}%
    \ifdim\v@leur<\v@lmax\Pss@tsecondSt\psline[-5,\the\s@mme]\Psrest@reSt\fi%
    \PSc@mment{End psaltitude}\resetc@ntr@l\et@tpsaltitude\fi\fi}}
\ctr@ld@f\def\Ps@rcerc#1;#2(#3,#4){\ellBB@x#1;#2,#2(#3,#4,0)%
    \f@gnewpath{\delt@=#2\unit@\delt@=\ptT@ptps\delt@%
    \BdingB@xfalse%
    \PSwrit@cmd{#1}{\repdecn@mb{\delt@}\space #3\space #4\space arc}{\fwf@g}}}
\ctr@ln@m\psarccirc
\ctr@ld@f\def\psarccircDD#1;#2(#3,#4){\ifcurr@ntPS\ifps@cri%
    \PSc@mment{psarccircDD Center=#1 ; Radius=#2 (Ang1=#3, Ang2=#4)}%
    \iffillm@de\Ps@rcerc#1;#2(#3,#4)%
    \f@gfill%
    \else\Ps@rcerc#1;#2(#3,#4)\f@gstroke\fi%
    \PSc@mment{End psarccircDD}\fi\fi}
\ctr@ld@f\def\psarccircTD#1,#2,#3;#4(#5,#6){{\ifcurr@ntPS\ifps@cri\s@uvc@ntr@l\et@tpsarccircTD%
    \PSc@mment{psarccircTD Center=#1,P1=#2,P2=#3 ; Radius=#4 (Ang1=#5, Ang2=#6)}%
    \setc@ntr@l{2}\c@lExtAxes#1,#2,#3(#4)\psarcellPATD#1,-4,-5(#5,#6)%
    \PSc@mment{End psarccircTD}\resetc@ntr@l\et@tpsarccircTD\fi\fi}}
\ctr@ld@f\def\c@lExtAxes#1,#2,#3(#4){%
    \figvectPTD-5[#1,#2]\vecunit@{-5}{-5}\figvectNTD-4[#1,#2,#3]\vecunit@{-4}{-4}%
    \figvectNVTD-3[-4,-5]\delt@=#4\unit@\edef\r@yon{\repdecn@mb{\delt@}}%
    \figpttra-4:=#1/\r@yon,-5/\figpttra-5:=#1/\r@yon,-3/}
\ctr@ln@m\psarccircP
\ctr@ld@f\def\psarccircPDD#1;#2[#3,#4]{{\ifcurr@ntPS\ifps@cri\s@uvc@ntr@l\et@tpsarccircPDD%
    \PSc@mment{psarccircPDD Center=#1; Radius=#2, [P1=#3, P2=#4]}%
    \Ps@ngleparam#1;#2[#3,#4]\ifdim\v@lmin>\v@lmax\advance\v@lmax\DePI@deg\fi%
    \edef\@ngdeb{\repdecn@mb{\v@lmin}}\edef\@ngfin{\repdecn@mb{\v@lmax}}%
    \psarccirc#1;\r@dius(\@ngdeb,\@ngfin)%
    \PSc@mment{End psarccircPDD}\resetc@ntr@l\et@tpsarccircPDD\fi\fi}}
\ctr@ld@f\def\psarccircPTD#1;#2[#3,#4,#5]{{\ifcurr@ntPS\ifps@cri\s@uvc@ntr@l\et@tpsarccircPTD%
    \PSc@mment{psarccircPTD Center=#1; Radius=#2, [P1=#3, P2=#4, P3=#5]}%
    \setc@ntr@l{2}\c@lExtAxes#1,#3,#5(#2)\psarcellPP#1,-4,-5[#3,#4]%
    \PSc@mment{End psarccircPTD}\resetc@ntr@l\et@tpsarccircPTD\fi\fi}}
\ctr@ld@f\def\Ps@ngleparam#1;#2[#3,#4]{\setc@ntr@l{2}%
    \figvectPDD-1[#1,#3]\vecunit@{-1}{-1}\Figg@tXY{-1}\arct@n\v@lmin(\v@lX,\v@lY)%
    \figvectPDD-2[#1,#4]\vecunit@{-2}{-2}\Figg@tXY{-2}\arct@n\v@lmax(\v@lX,\v@lY)%
    \v@lmin=\rdT@deg\v@lmin\v@lmax=\rdT@deg\v@lmax%
    \v@leur=#2pt\maxim@m{\mili@u}{-\v@leur}{\v@leur}%
    \edef\r@dius{\repdecn@mb{\mili@u}}}
\ctr@ld@f\def\Ps@rcercBz#1;#2(#3,#4){\Ps@rellBz#1;#2,#2(#3,#4,0)}
\ctr@ld@f\def\Ps@rellBz#1;#2,#3(#4,#5,#6){%
    \ellBB@x#1;#2,#3(#4,#5,#6)\BdingB@xfalse%
    \c@lNbarcs{#4}{#5}\v@leur=#4pt\setc@ntr@l{2}\figptell-13::#1;#2,#3(#4,#6)%
    \f@gnewpath\PSwrit@cmd{-13}{\c@mmoveto}{\fwf@g}%
    \s@mme=\z@\bcl@rellBz#1;#2,#3(#6)\BdingB@xtrue}
\ctr@ld@f\def\bcl@rellBz#1;#2,#3(#4){\relax%
    \ifnum\s@mme<\p@rtent\advance\s@mme\@ne%
    \advance\v@leur\delt@\edef\@ngle{\repdecn@mb\v@leur}\figptell-14::#1;#2,#3(\@ngle,#4)%
    \advance\v@leur\delt@\edef\@ngle{\repdecn@mb\v@leur}\figptell-15::#1;#2,#3(\@ngle,#4)%
    \advance\v@leur\delt@\edef\@ngle{\repdecn@mb\v@leur}\figptell-16::#1;#2,#3(\@ngle,#4)%
    \figptscontrolDD-18[-13,-14,-15,-16]%
    \PSwrit@cmd{-18}{}{\fwf@g}\PSwrit@cmd{-17}{}{\fwf@g}%
    \PSwrit@cmd{-16}{\c@mcurveto}{\fwf@g}%
    \figptcopyDD-13:/-16/\bcl@rellBz#1;#2,#3(#4)\fi}
\ctr@ld@f\def\Ps@rell#1;#2,#3(#4,#5,#6){\ellBB@x#1;#2,#3(#4,#5,#6)%
    \f@gnewpath{\v@lmin=#2\unit@\v@lmin=\ptT@ptps\v@lmin%
    \v@lmax=#3\unit@\v@lmax=\ptT@ptps\v@lmax\BdingB@xfalse%
    \PSwrit@cmd{#1}%
    {#6\space\repdecn@mb{\v@lmin}\space\repdecn@mb{\v@lmax}\space #4\space #5\space ellipse}{\fwf@g}}%
    \global\Use@llipsetrue}
\ctr@ln@m\psarcell
\ctr@ld@f\def\psarcellDD#1;#2,#3(#4,#5,#6){{\ifcurr@ntPS\ifps@cri%
    \PSc@mment{psarcellDD Center=#1 ; XRad=#2, YRad=#3 (Ang1=#4, Ang2=#5, Inclination=#6)}%
    \iffillm@de\Ps@rell#1;#2,#3(#4,#5,#6)%
    \f@gfill%
    \else\Ps@rell#1;#2,#3(#4,#5,#6)\f@gstroke\fi%
    \PSc@mment{End psarcellDD}\fi\fi}}
\ctr@ld@f\def\psarcellTD#1;#2,#3(#4,#5,#6){{\ifcurr@ntPS\ifps@cri\s@uvc@ntr@l\et@tpsarcellTD%
    \PSc@mment{psarcellTD Center=#1 ; XRad=#2, YRad=#3 (Ang1=#4, Ang2=#5, Inclination=#6)}%
    \setc@ntr@l{2}\figpttraC -8:=#1/#2,0,0/\figpttraC -7:=#1/0,#3,0/%
    \figvectC -4(0,0,1)\figptsrot -8=-8,-7/#1,#6,-4/\psarcellPATD#1,-8,-7(#4,#5)%
    \PSc@mment{End psarcellTD}\resetc@ntr@l\et@tpsarcellTD\fi\fi}}
\ctr@ln@m\psarcellPA
\ctr@ld@f\def\psarcellPADD#1,#2,#3(#4,#5){{\ifcurr@ntPS\ifps@cri\s@uvc@ntr@l\et@tpsarcellPADD%
    \PSc@mment{psarcellPADD Center=#1,PtAxis1=#2,PtAxis2=#3 (Ang1=#4, Ang2=#5)}%
    \setc@ntr@l{2}\figvectPDD-1[#1,#2]\vecunit@DD{-1}{-1}\v@lX=\ptT@unit@\result@t%
    \edef\XR@d{\repdecn@mb{\v@lX}}\Figg@tXY{-1}\arct@n\v@lmin(\v@lX,\v@lY)%
    \v@lmin=\rdT@deg\v@lmin\edef\Inclin@{\repdecn@mb{\v@lmin}}%
    \figgetdist\YR@d[#1,#3]\psarcellDD#1;\XR@d,\YR@d(#4,#5,\Inclin@)%
    \PSc@mment{End psarcellPADD}\resetc@ntr@l\et@tpsarcellPADD\fi\fi}}
\ctr@ld@f\def\psarcellPATD#1,#2,#3(#4,#5){{\ifcurr@ntPS\ifps@cri\s@uvc@ntr@l\et@tpsarcellPATD%
    \PSc@mment{psarcellPATD Center=#1,PtAxis1=#2,PtAxis2=#3 (Ang1=#4, Ang2=#5)}%
    \iffillm@de\Ps@rellPATD#1,#2,#3(#4,#5)%
    \f@gfill%
    \else\Ps@rellPATD#1,#2,#3(#4,#5)\f@gstroke\fi%
    \PSc@mment{End psarcellPATD}\resetc@ntr@l\et@tpsarcellPATD\fi\fi}}
\ctr@ld@f\def\Ps@rellPATD#1,#2,#3(#4,#5){\let\c@lprojSP=\relax%
    \setc@ntr@l{2}\figvectPTD-1[#1,#2]\figvectPTD-2[#1,#3]\c@lNbarcs{#4}{#5}%
    \v@leur=#4pt\c@lptellP{#1}{-1}{-2}\Figptpr@j-5:/-3/%
    \f@gnewpath\PSwrit@cmdS{-5}{\c@mmoveto}{\fwf@g}{\X@un}{\Y@un}%
    \edef\C@nt@r{#1}\s@mme=\z@\bcl@rellPATD}
\ctr@ld@f\def\bcl@rellPATD{\relax%
    \ifnum\s@mme<\p@rtent\advance\s@mme\@ne%
    \advance\v@leur\delt@\c@lptellP{\C@nt@r}{-1}{-2}\Figptpr@j-4:/-3/%
    \advance\v@leur\delt@\c@lptellP{\C@nt@r}{-1}{-2}\Figptpr@j-6:/-3/%
    \advance\v@leur\delt@\c@lptellP{\C@nt@r}{-1}{-2}\Figptpr@j-3:/-3/%
    \v@lX=\z@\v@lY=\z@\Figtr@nptDD{-5}{-5}\Figtr@nptDD{2}{-3}%
    \divide\v@lX\@vi\divide\v@lY\@vi%
    \Figtr@nptDD{3}{-4}\Figtr@nptDD{-1.5}{-6}\v@lmin=\v@lX\v@lmax=\v@lY%
    \v@lX=\z@\v@lY=\z@\Figtr@nptDD{2}{-5}\Figtr@nptDD{-5}{-3}%
    \divide\v@lX\@vi\divide\v@lY\@vi\Figtr@nptDD{-1.5}{-4}\Figtr@nptDD{3}{-6}%
    \BdingB@xfalse%
    \Figp@intregDD-4:(\v@lmin,\v@lmax)\PSwrit@cmdS{-4}{}{\fwf@g}{\X@de}{\Y@de}%
    \Figp@intregDD-4:(\v@lX,\v@lY)\PSwrit@cmdS{-4}{}{\fwf@g}{\X@tr}{\Y@tr}%
    \BdingB@xtrue\PSwrit@cmdS{-3}{\c@mcurveto}{\fwf@g}{\X@qu}{\Y@qu}%
    \B@zierBB@x{1}{\Y@un}(\X@un,\X@de,\X@tr,\X@qu)%
    \B@zierBB@x{2}{\X@un}(\Y@un,\Y@de,\Y@tr,\Y@qu)%
    \edef\X@un{\X@qu}\edef\Y@un{\Y@qu}\figptcopyDD-5:/-3/\bcl@rellPATD\fi}
\ctr@ld@f\def\c@lNbarcs#1#2{%
    \delt@=#2pt\advance\delt@-#1pt\maxim@m{\v@lmax}{\delt@}{-\delt@}%
    \v@leur=\v@lmax\divide\v@leur45 \p@rtentiere{\p@rtent}{\v@leur}\advance\p@rtent\@ne%
    \s@mme=\p@rtent\multiply\s@mme\thr@@\divide\delt@\s@mme}
\ctr@ld@f\def\psarcellPP#1,#2,#3[#4,#5]{{\ifcurr@ntPS\ifps@cri\s@uvc@ntr@l\et@tpsarcellPP%
    \PSc@mment{psarcellPP Center=#1,PtAxis1=#2,PtAxis2=#3 [Point1=#4, Point2=#5]}%
    \setc@ntr@l{2}\figvectP-2[#1,#3]\vecunit@{-2}{-2}\v@lmin=\result@t%
    \invers@{\v@lmax}{\v@lmin}%
    \figvectP-1[#1,#2]\vecunit@{-1}{-1}\v@leur=\result@t%
    \v@leur=\repdecn@mb{\v@lmax}\v@leur\edef\AsB@{\repdecn@mb{\v@leur}}
    \c@lAngle{#1}{#4}{\v@lmin}\edef\@ngdeb{\repdecn@mb{\v@lmin}}%
    \c@lAngle{#1}{#5}{\v@lmax}\ifdim\v@lmin>\v@lmax\advance\v@lmax\DePI@deg\fi%
    \edef\@ngfin{\repdecn@mb{\v@lmax}}\psarcellPA#1,#2,#3(\@ngdeb,\@ngfin)%
    \PSc@mment{End psarcellPP}\resetc@ntr@l\et@tpsarcellPP\fi\fi}}
\ctr@ld@f\def\c@lAngle#1#2#3{\figvectP-3[#1,#2]%
    \c@lproscal\delt@[-3,-1]\c@lproscal\v@leur[-3,-2]%
    \v@leur=\AsB@\v@leur\arct@n#3(\delt@,\v@leur)#3=\rdT@deg#3}
\ctr@ln@w{newif}\if@rrowratio\@rrowratiotrue
\ctr@ln@w{newif}\if@rrowhfill
\ctr@ln@w{newif}\if@rrowhout
\ctr@ld@f\def\Psset@rrowhe@d#1=#2|{\keln@mun#1|%
    \def\n@mref{a}\ifx\l@debut\n@mref\pssetarrowheadangle{#2}\else
    \def\n@mref{f}\ifx\l@debut\n@mref\pssetarrowheadfill{#2}\else
    \def\n@mref{l}\ifx\l@debut\n@mref\pssetarrowheadlength{#2}\else
    \def\n@mref{o}\ifx\l@debut\n@mref\pssetarrowheadout{#2}\else
    \def\n@mref{r}\ifx\l@debut\n@mref\pssetarrowheadratio{#2}\else
    \immediate\write16{*** Unknown attribute: \BS@ psset arrowhead(..., #1=...)}%
    \fi\fi\fi\fi\fi}
\ctr@ln@m\@rrowheadangle
\ctr@ln@m\C@AHANG \ctr@ln@m\S@AHANG \ctr@ln@m\UNSS@N
\ctr@ld@f\def\pssetarrowheadangle#1{\edef\@rrowheadangle{#1}{\c@ssin{\C@}{\S@}{#1}%
    \xdef\C@AHANG{\C@}\xdef\S@AHANG{\S@}\v@lmax=\S@ pt%
    \invers@{\v@leur}{\v@lmax}\maxim@m{\v@leur}{\v@leur}{-\v@leur}%
    \xdef\UNSS@N{\the\v@leur}}}
\ctr@ld@f\def\pssetarrowheadfill#1{\expandafter\set@rrowhfill#1:}
\ctr@ld@f\def\set@rrowhfill#1#2:{\if#1n\@rrowhfillfalse\else\@rrowhfilltrue\fi}
\ctr@ld@f\def\pssetarrowheadout#1{\expandafter\set@rrowhout#1:}
\ctr@ld@f\def\set@rrowhout#1#2:{\if#1n\@rrowhoutfalse\else\@rrowhouttrue\fi}
\ctr@ln@m\@rrowheadlength
\ctr@ld@f\def\pssetarrowheadlength#1{\edef\@rrowheadlength{#1}\@rrowratiofalse}
\ctr@ln@m\@rrowheadratio
\ctr@ld@f\def\pssetarrowheadratio#1{\edef\@rrowheadratio{#1}\@rrowratiotrue}
\ctr@ln@m\defaultarrowheadlength
\ctr@ld@f\def\psresetarrowhead{%
    \pssetarrowheadangle{\defaultarrowheadangle}%
    \pssetarrowheadfill{\defaultarrowheadfill}%
    \pssetarrowheadout{\defaultarrowheadout}%
    \pssetarrowheadratio{\defaultarrowheadratio}%
    \d@fm@cdim\defaultarrowheadlength{\defaulth@rdahlength}
    \pssetarrowheadlength{\defaultarrowheadlength}}
\ctr@ld@f\def\defaultarrowheadratio{0.1}
\ctr@ld@f\def\defaultarrowheadangle{20}
\ctr@ld@f\def\defaultarrowheadfill{no}
\ctr@ld@f\def\defaultarrowheadout{no}
\ctr@ld@f\def\defaulth@rdahlength{8pt}
\ctr@ln@m\psarrow
\ctr@ld@f\def\psarrowDD[#1,#2]{{\ifcurr@ntPS\ifps@cri\s@uvc@ntr@l\et@tpsarrow%
    \PSc@mment{psarrowDD [Pt1,Pt2]=[#1,#2]}\pssetfillmode{no}%
    \psarrowheadDD[#1,#2]\setc@ntr@l{2}\psline[#1,-3]%
    \PSc@mment{End psarrowDD}\resetc@ntr@l\et@tpsarrow\fi\fi}}
\ctr@ld@f\def\psarrowTD[#1,#2]{{\ifcurr@ntPS\ifps@cri\s@uvc@ntr@l\et@tpsarrowTD%
    \PSc@mment{psarrowTD [Pt1,Pt2]=[#1,#2]}\resetc@ntr@l{2}%
    \Figptpr@j-5:/#1/\Figptpr@j-6:/#2/\let\c@lprojSP=\relax\psarrowDD[-5,-6]%
    \PSc@mment{End psarrowTD}\resetc@ntr@l\et@tpsarrowTD\fi\fi}}
\ctr@ln@m\psarrowhead
\ctr@ld@f\def\psarrowheadDD[#1,#2]{{\ifcurr@ntPS\ifps@cri\s@uvc@ntr@l\et@tpsarrowheadDD%
    \if@rrowhfill\def\@hangle{-\@rrowheadangle}\else\def\@hangle{\@rrowheadangle}\fi%
    \if@rrowratio%
    \if@rrowhout\def\@hratio{-\@rrowheadratio}\else\def\@hratio{\@rrowheadratio}\fi%
    \PSc@mment{psarrowheadDD Ratio=\@hratio, Angle=\@hangle, [Pt1,Pt2]=[#1,#2]}%
    \Ps@rrowhead\@hratio,\@hangle[#1,#2]%
    \else%
    \if@rrowhout\def\@hlength{-\@rrowheadlength}\else\def\@hlength{\@rrowheadlength}\fi%
    \PSc@mment{psarrowheadDD Length=\@hlength, Angle=\@hangle, [Pt1,Pt2]=[#1,#2]}%
    \Ps@rrowheadfd\@hlength,\@hangle[#1,#2]%
    \fi%
    \PSc@mment{End psarrowheadDD}\resetc@ntr@l\et@tpsarrowheadDD\fi\fi}}
\ctr@ld@f\def\psarrowheadTD[#1,#2]{{\ifcurr@ntPS\ifps@cri\s@uvc@ntr@l\et@tpsarrowheadTD%
    \PSc@mment{psarrowheadTD [Pt1,Pt2]=[#1,#2]}\resetc@ntr@l{2}%
    \Figptpr@j-5:/#1/\Figptpr@j-6:/#2/\let\c@lprojSP=\relax\psarrowheadDD[-5,-6]%
    \PSc@mment{End psarrowheadTD}\resetc@ntr@l\et@tpsarrowheadTD\fi\fi}}
\ctr@ld@f\def\Ps@rrowhead#1,#2[#3,#4]{\v@leur=#1\p@\maxim@m{\v@leur}{\v@leur}{-\v@leur}%
    \ifdim\v@leur>\Cepsil@n{
    \PSc@mment{ps@rrowhead Ratio=#1, Angle=#2, [Pt1,Pt2]=[#3,#4]}\v@leur=\UNSS@N%
    \v@leur=\curr@ntwidth\v@leur\v@leur=\ptpsT@pt\v@leur\delt@=.5\v@leur
    \setc@ntr@l{2}\figvectPDD-3[#4,#3]%
    \Figg@tXY{-3}\v@lX=#1\v@lX\v@lY=#1\v@lY\Figv@ctCreg-3(\v@lX,\v@lY)%
    \vecunit@{-4}{-3}\mili@u=\result@t%
    \ifdim#2pt>\z@\v@lXa=-\C@AHANG\delt@%
     \edef\c@ef{\repdecn@mb{\v@lXa}}\figpttraDD-3:=-3/\c@ef,-4/\fi%
    \edef\c@ef{\repdecn@mb{\delt@}}%
    \v@lXa=\mili@u\v@lXa=\C@AHANG\v@lXa%
    \v@lYa=\ptpsT@pt\p@\v@lYa=\curr@ntwidth\v@lYa\v@lYa=\sDcc@ngle\v@lYa%
    \advance\v@lXa-\v@lYa\gdef\sDcc@ngle{0}%
    \ifdim\v@lXa>\v@leur\edef\c@efendpt{\repdecn@mb{\v@leur}}%
    \else\edef\c@efendpt{\repdecn@mb{\v@lXa}}\fi%
    \Figg@tXY{-3}\v@lmin=\v@lX\v@lmax=\v@lY%
    \v@lXa=\C@AHANG\v@lmin\v@lYa=\S@AHANG\v@lmax\advance\v@lXa\v@lYa%
    \v@lYa=-\S@AHANG\v@lmin\v@lX=\C@AHANG\v@lmax\advance\v@lYa\v@lX%
    \setc@ntr@l{1}\Figg@tXY{#4}\advance\v@lX\v@lXa\advance\v@lY\v@lYa%
    \setc@ntr@l{2}\Figp@intregDD-2:(\v@lX,\v@lY)%
    \v@lXa=\C@AHANG\v@lmin\v@lYa=-\S@AHANG\v@lmax\advance\v@lXa\v@lYa%
    \v@lYa=\S@AHANG\v@lmin\v@lX=\C@AHANG\v@lmax\advance\v@lYa\v@lX%
    \setc@ntr@l{1}\Figg@tXY{#4}\advance\v@lX\v@lXa\advance\v@lY\v@lYa%
    \setc@ntr@l{2}\Figp@intregDD-1:(\v@lX,\v@lY)%
    \ifdim#2pt<\z@\fillm@detrue\psline[-2,#4,-1]
    \else\figptstraDD-3=#4,-2,-1/\c@ef,-4/\psline[-2,-3,-1]\fi
    \ifdim#1pt>\z@\figpttraDD-3:=#4/\c@efendpt,-4/\else\figptcopyDD-3:/#4/\fi%
    \PSc@mment{End ps@rrowhead}}\fi}
\ctr@ld@f\def\sDcc@ngle{0}
\ctr@ld@f\def\Ps@rrowheadfd#1,#2[#3,#4]{{%
    \PSc@mment{ps@rrowheadfd Length=#1, Angle=#2, [Pt1,Pt2]=[#3,#4]}%
    \setc@ntr@l{2}\figvectPDD-1[#3,#4]\n@rmeucDD{\v@leur}{-1}\v@leur=\ptT@unit@\v@leur%
    \invers@{\v@leur}{\v@leur}\v@leur=#1\v@leur\edef\R@tio{\repdecn@mb{\v@leur}}%
    \Ps@rrowhead\R@tio,#2[#3,#4]\PSc@mment{End ps@rrowheadfd}}}
\ctr@ln@m\psarrowBezier
\ctr@ld@f\def\psarrowBezierDD[#1,#2,#3,#4]{{\ifcurr@ntPS\ifps@cri\s@uvc@ntr@l\et@tpsarrowBezierDD%
    \PSc@mment{psarrowBezierDD Control points=#1,#2,#3,#4}\setc@ntr@l{2}%
    \if@rrowratio\c@larclengthDD\v@leur,10[#1,#2,#3,#4]\else\v@leur=\z@\fi%
    \Ps@rrowB@zDD\v@leur[#1,#2,#3,#4]%
    \PSc@mment{End psarrowBezierDD}\resetc@ntr@l\et@tpsarrowBezierDD\fi\fi}}
\ctr@ld@f\def\psarrowBezierTD[#1,#2,#3,#4]{{\ifcurr@ntPS\ifps@cri\s@uvc@ntr@l\et@tpsarrowBezierTD%
    \PSc@mment{psarrowBezierTD Control points=#1,#2,#3,#4}\resetc@ntr@l{2}%
    \Figptpr@j-7:/#1/\Figptpr@j-8:/#2/\Figptpr@j-9:/#3/\Figptpr@j-10:/#4/%
    \let\c@lprojSP=\relax\ifnum\curr@ntproj<\tw@\psarrowBezierDD[-7,-8,-9,-10]%
    \else\f@gnewpath\PSwrit@cmd{-7}{\c@mmoveto}{\fwf@g}%
    \if@rrowratio\c@larclengthDD\mili@u,10[-7,-8,-9,-10]\else\mili@u=\z@\fi%
    \p@rtent=\NBz@rcs\advance\p@rtent\m@ne\subB@zierTD\p@rtent[#1,#2,#3,#4]%
    \f@gstroke%
    \advance\v@lmin\p@rtent\delt@
    \v@leur=\v@lmin\advance\v@leur0.33333 \delt@\edef\unti@rs{\repdecn@mb{\v@leur}}%
    \v@leur=\v@lmin\advance\v@leur0.66666 \delt@\edef\deti@rs{\repdecn@mb{\v@leur}}%
    \figptcopyDD-8:/-10/\c@lsubBzarc\unti@rs,\deti@rs[#1,#2,#3,#4]%
    \figptcopyDD-8:/-4/\figptcopyDD-9:/-3/\Ps@rrowB@zDD\mili@u[-7,-8,-9,-10]\fi%
    \PSc@mment{End psarrowBezierTD}\resetc@ntr@l\et@tpsarrowBezierTD\fi\fi}}
\ctr@ld@f\def\c@larclengthDD#1,#2[#3,#4,#5,#6]{{\p@rtent=#2\figptcopyDD-5:/#3/%
    \delt@=\p@\divide\delt@\p@rtent\c@rre=\z@\v@leur=\z@\s@mme=\z@%
    \loop\ifnum\s@mme<\p@rtent\advance\s@mme\@ne\advance\v@leur\delt@%
    \edef\T@{\repdecn@mb{\v@leur}}\figptBezierDD-6::\T@[#3,#4,#5,#6]%
    \figvectPDD-1[-5,-6]\n@rmeucDD{\mili@u}{-1}\advance\c@rre\mili@u%
    \figptcopyDD-5:/-6/\repeat\global\result@t=\ptT@unit@\c@rre}#1=\result@t}
\ctr@ld@f\def\Ps@rrowB@zDD#1[#2,#3,#4,#5]{{\pssetfillmode{no}%
    \if@rrowratio\delt@=\@rrowheadratio#1\else\delt@=\@rrowheadlength pt\fi%
    \v@leur=\C@AHANG\delt@\edef\R@dius{\repdecn@mb{\v@leur}}%
    \FigptintercircB@zDD-5::0,\R@dius[#5,#4,#3,#2]%
    \pssetarrowheadlength{\repdecn@mb{\delt@}}\psarrowheadDD[-5,#5]%
    \let\n@rmeuc=\n@rmeucDD\figgetdist\R@dius[#5,-3]%
    \FigptintercircB@zDD-6::0,\R@dius[#5,#4,#3,#2]%
    \figptBezierDD-5::0.33333[#5,#4,#3,#2]\figptBezierDD-3::0.66666[#5,#4,#3,#2]%
    \figptscontrolDD-5[-6,-5,-3,#2]\psBezierDD1[-6,-5,-4,#2]}}
\ctr@ln@m\psarrowcirc
\ctr@ld@f\def\psarrowcircDD#1;#2(#3,#4){{\ifcurr@ntPS\ifps@cri\s@uvc@ntr@l\et@tpsarrowcircDD%
    \PSc@mment{psarrowcircDD Center=#1 ; Radius=#2 (Ang1=#3,Ang2=#4)}%
    \pssetfillmode{no}\Pscirc@rrowhead#1;#2(#3,#4)%
    \setc@ntr@l{2}\figvectPDD -4[#1,-3]\vecunit@{-4}{-4}%
    \Figg@tXY{-4}\arct@n\v@lmin(\v@lX,\v@lY)%
    \v@lmin=\rdT@deg\v@lmin\v@leur=#4pt\advance\v@leur-\v@lmin%
    \maxim@m{\v@leur}{\v@leur}{-\v@leur}%
    \ifdim\v@leur>\DemiPI@deg\relax\ifdim\v@lmin<#4pt\advance\v@lmin\DePI@deg%
    \else\advance\v@lmin-\DePI@deg\fi\fi\edef\ar@ngle{\repdecn@mb{\v@lmin}}%
    \ifdim#3pt<#4pt\psarccirc#1;#2(#3,\ar@ngle)\else\psarccirc#1;#2(\ar@ngle,#3)\fi%
    \PSc@mment{End psarrowcircDD}\resetc@ntr@l\et@tpsarrowcircDD\fi\fi}}
\ctr@ld@f\def\psarrowcircTD#1,#2,#3;#4(#5,#6){{\ifcurr@ntPS\ifps@cri\s@uvc@ntr@l\et@tpsarrowcircTD%
    \PSc@mment{psarrowcircTD Center=#1,P1=#2,P2=#3 ; Radius=#4 (Ang1=#5, Ang2=#6)}%
    \resetc@ntr@l{2}\c@lExtAxes#1,#2,#3(#4)\let\c@lprojSP=\relax%
    \figvectPTD-11[#1,-4]\figvectPTD-12[#1,-5]\c@lNbarcs{#5}{#6}%
    \if@rrowratio\v@lmax=\degT@rd\v@lmax\edef\D@lpha{\repdecn@mb{\v@lmax}}\fi%
    \advance\p@rtent\m@ne\mili@u=\z@%
    \v@leur=#5pt\c@lptellP{#1}{-11}{-12}\Figptpr@j-9:/-3/%
    \f@gnewpath\PSwrit@cmdS{-9}{\c@mmoveto}{\fwf@g}{\X@un}{\Y@un}%
    \edef\C@nt@r{#1}\s@mme=\z@\bcl@rcircTD\f@gstroke%
    \advance\v@leur\delt@\c@lptellP{#1}{-11}{-12}\Figptpr@j-5:/-3/%
    \advance\v@leur\delt@\c@lptellP{#1}{-11}{-12}\Figptpr@j-6:/-3/%
    \advance\v@leur\delt@\c@lptellP{#1}{-11}{-12}\Figptpr@j-10:/-3/%
    \figptscontrolDD-8[-9,-5,-6,-10]%
    \if@rrowratio\c@lcurvradDD0.5[-9,-8,-7,-10]\advance\mili@u\result@t%
    \maxim@m{\mili@u}{\mili@u}{-\mili@u}\mili@u=\ptT@unit@\mili@u%
    \mili@u=\D@lpha\mili@u\advance\p@rtent\@ne\divide\mili@u\p@rtent\fi%
    \Ps@rrowB@zDD\mili@u[-9,-8,-7,-10]%
    \PSc@mment{End psarrowcircTD}\resetc@ntr@l\et@tpsarrowcircTD\fi\fi}}
\ctr@ld@f\def\bcl@rcircTD{\relax%
    \ifnum\s@mme<\p@rtent\advance\s@mme\@ne%
    \advance\v@leur\delt@\c@lptellP{\C@nt@r}{-11}{-12}\Figptpr@j-5:/-3/%
    \advance\v@leur\delt@\c@lptellP{\C@nt@r}{-11}{-12}\Figptpr@j-6:/-3/%
    \advance\v@leur\delt@\c@lptellP{\C@nt@r}{-11}{-12}\Figptpr@j-10:/-3/%
    \figptscontrolDD-8[-9,-5,-6,-10]\BdingB@xfalse%
    \PSwrit@cmdS{-8}{}{\fwf@g}{\X@de}{\Y@de}\PSwrit@cmdS{-7}{}{\fwf@g}{\X@tr}{\Y@tr}%
    \BdingB@xtrue\PSwrit@cmdS{-10}{\c@mcurveto}{\fwf@g}{\X@qu}{\Y@qu}%
    \if@rrowratio\c@lcurvradDD0.5[-9,-8,-7,-10]\advance\mili@u\result@t\fi%
    \B@zierBB@x{1}{\Y@un}(\X@un,\X@de,\X@tr,\X@qu)%
    \B@zierBB@x{2}{\X@un}(\Y@un,\Y@de,\Y@tr,\Y@qu)%
    \edef\X@un{\X@qu}\edef\Y@un{\Y@qu}\figptcopyDD-9:/-10/\bcl@rcircTD\fi}
\ctr@ld@f\def\Pscirc@rrowhead#1;#2(#3,#4){{%
    \PSc@mment{pscirc@rrowhead Center=#1 ; Radius=#2 (Ang1=#3,Ang2=#4)}%
    \v@leur=#2\unit@\edef\s@glen{\repdecn@mb{\v@leur}}\v@lY=\z@\v@lX=\v@leur%
    \resetc@ntr@l{2}\Figv@ctCreg-3(\v@lX,\v@lY)\figpttraDD-5:=#1/1,-3/%
    \figptrotDD-5:=-5/#1,#4/%
    \figvectPDD-3[#1,-5]\Figg@tXY{-3}\v@leur=\v@lX%
    \ifdim#3pt<#4pt\v@lX=\v@lY\v@lY=-\v@leur\else\v@lX=-\v@lY\v@lY=\v@leur\fi%
    \Figv@ctCreg-3(\v@lX,\v@lY)\vecunit@{-3}{-3}%
    \if@rrowratio\v@leur=#4pt\advance\v@leur-#3pt\maxim@m{\mili@u}{-\v@leur}{\v@leur}%
    \mili@u=\degT@rd\mili@u\v@leur=\s@glen\mili@u\edef\s@glen{\repdecn@mb{\v@leur}}%
    \mili@u=#2\mili@u\mili@u=\@rrowheadratio\mili@u\else\mili@u=\@rrowheadlength pt\fi%
    \figpttraDD-6:=-5/\s@glen,-3/\v@leur=#2pt\v@leur=2\v@leur%
    \invers@{\v@leur}{\v@leur}\c@rre=\repdecn@mb{\v@leur}\mili@u
    \mili@u=\c@rre\mili@u=\repdecn@mb{\c@rre}\mili@u%
    \v@leur=\p@\advance\v@leur-\mili@u
    \invers@{\mili@u}{2\v@leur}\delt@=\c@rre\delt@=\repdecn@mb{\mili@u}\delt@%
    \xdef\sDcc@ngle{\repdecn@mb{\delt@}}
    \sqrt@{\mili@u}{\v@leur}\arct@n\v@leur(\mili@u,\c@rre)%
    \v@leur=\rdT@deg\v@leur
    \ifdim#3pt<#4pt\v@leur=-\v@leur\fi%
    \if@rrowhout\v@leur=-\v@leur\fi\edef\cor@ngle{\repdecn@mb{\v@leur}}%
    \figptrotDD-6:=-6/-5,\cor@ngle/\psarrowheadDD[-6,-5]%
    \PSc@mment{End pscirc@rrowhead}}}
\ctr@ln@m\psarrowcircP
\ctr@ld@f\def\psarrowcircPDD#1;#2[#3,#4]{{\ifcurr@ntPS\ifps@cri%
    \PSc@mment{psarrowcircPDD Center=#1; Radius=#2, [P1=#3,P2=#4]}%
    \s@uvc@ntr@l\et@tpsarrowcircPDD\Ps@ngleparam#1;#2[#3,#4]%
    \ifdim\v@leur>\z@\ifdim\v@lmin>\v@lmax\advance\v@lmax\DePI@deg\fi%
    \else\ifdim\v@lmin<\v@lmax\advance\v@lmin\DePI@deg\fi\fi%
    \edef\@ngdeb{\repdecn@mb{\v@lmin}}\edef\@ngfin{\repdecn@mb{\v@lmax}}%
    \psarrowcirc#1;\r@dius(\@ngdeb,\@ngfin)%
    \PSc@mment{End psarrowcircPDD}\resetc@ntr@l\et@tpsarrowcircPDD\fi\fi}}
\ctr@ld@f\def\psarrowcircPTD#1;#2[#3,#4,#5]{{\ifcurr@ntPS\ifps@cri\s@uvc@ntr@l\et@tpsarrowcircPTD%
    \PSc@mment{psarrowcircPTD Center=#1; Radius=#2, [P1=#3,P2=#4,P3=#5]}%
    \figgetangleTD\@ngfin[#1,#3,#4,#5]\v@leur=#2pt%
    \maxim@m{\mili@u}{-\v@leur}{\v@leur}\edef\r@dius{\repdecn@mb{\mili@u}}%
    \ifdim\v@leur<\z@\v@lmax=\@ngfin pt\advance\v@lmax-\DePI@deg%
    \edef\@ngfin{\repdecn@mb{\v@lmax}}\fi\psarrowcircTD#1,#3,#5;\r@dius(0,\@ngfin)%
    \PSc@mment{End psarrowcircPTD}\resetc@ntr@l\et@tpsarrowcircPTD\fi\fi}}
\ctr@ld@f\def\psaxes#1(#2){{\ifcurr@ntPS\ifps@cri\s@uvc@ntr@l\et@tpsaxes%
    \PSc@mment{psaxes Origin=#1 Range=(#2)}\an@lys@xes#2,:\resetc@ntr@l{2}%
    \ifx\t@xt@\empty\ifTr@isDim\ps@xes#1(0,#2,0,#2,0,#2)\else\ps@xes#1(0,#2,0,#2)\fi%
    \else\ps@xes#1(#2)\fi\PSc@mment{End psaxes}\resetc@ntr@l\et@tpsaxes\fi\fi}}
\ctr@ld@f\def\an@lys@xes#1,#2:{\def\t@xt@{#2}}
\ctr@ln@m\ps@xes
\ctr@ld@f\def\ps@xesDD#1(#2,#3,#4,#5){%
    \figpttraC-5:=#1/#2,0/\figpttraC-6:=#1/#3,0/\psarrowDD[-5,-6]%
    \figpttraC-5:=#1/0,#4/\figpttraC-6:=#1/0,#5/\psarrowDD[-5,-6]}
\ctr@ld@f\def\ps@xesTD#1(#2,#3,#4,#5,#6,#7){%
    \figpttraC-7:=#1/#2,0,0/\figpttraC-8:=#1/#3,0,0/\psarrowTD[-7,-8]%
    \figpttraC-7:=#1/0,#4,0/\figpttraC-8:=#1/0,#5,0/\psarrowTD[-7,-8]%
    \figpttraC-7:=#1/0,0,#6/\figpttraC-8:=#1/0,0,#7/\psarrowTD[-7,-8]}
\ctr@ln@m\newGr@FN
\ctr@ld@f\def\newGr@FNPDF#1{\s@mme=\Gr@FNb\advance\s@mme\@ne\xdef\Gr@FNb{\number\s@mme}}
\ctr@ld@f\def\newGr@FNDVI#1{\newGr@FNPDF{}\xdef#1{\jobname GI\Gr@FNb.anx}}
\ctr@ld@f\def\psbeginfig#1{\newGr@FN\DefGIfilen@me\gdef\@utoFN{0}%
    \def\t@xt@{#1}\relax\ifx\t@xt@\empty\psupdatem@detrue%
    \gdef\@utoFN{1}\Psb@ginfig\DefGIfilen@me\else\expandafter\Psb@ginfigNu@#1 :\fi}
\ctr@ld@f\def\Psb@ginfigNu@#1 #2:{\def\t@xt@{#1}\relax\ifx\t@xt@\empty\def\t@xt@{#2}%
    \ifx\t@xt@\empty\psupdatem@detrue\gdef\@utoFN{1}\Psb@ginfig\DefGIfilen@me%
    \else\Psb@ginfigNu@#2:\fi\else\Psb@ginfig{#1}\fi}
\ctr@ln@m\PSfilen@me \ctr@ln@m\auxfilen@me
\ctr@ld@f\def\Psb@ginfig#1{\ifcurr@ntPS\else%
    \edef\PSfilen@me{#1}\edef\auxfilen@me{\jobname.anx}%
    \ifpsupdatem@de\ps@critrue\else\openin\frf@g=\PSfilen@me\relax%
    \ifeof\frf@g\ps@critrue\else\ps@crifalse\fi\closein\frf@g\fi%
    \curr@ntPStrue\c@ldefproj\expandafter\setupd@te\defaultupdate:%
    \ifps@cri\initb@undb@x%
    \immediate\openout\fwf@g=\auxfilen@me\initpss@ttings\fi%
    \fi}
\ctr@ld@f\def\Gr@FNb{0}
\ctr@ld@f\def\figforTeXFileno{\Gr@FNb}
\ctr@ld@f\def\figforTeXFigno{0 }
\ctr@ld@f\def\figforTeXnextFigno{1 }
\ctr@ld@f\edef\DefGIfilen@me{\jobname GI.anx}
\ctr@ld@f\def\initpss@ttings{\psreset{arrowhead,curve,first,flowchart,mesh,second,third}%
    \Use@llipsefalse}
\ctr@ld@f\def\B@zierBB@x#1#2(#3,#4,#5,#6){{\c@rre=\t@n\epsil@n
    \v@lmax=#4\advance\v@lmax-#5\v@lmax=\thr@@\v@lmax\advance\v@lmax#6\advance\v@lmax-#3%
    \mili@u=#4\mili@u=-\tw@\mili@u\advance\mili@u#3\advance\mili@u#5%
    \v@lmin=#4\advance\v@lmin-#3\maxim@m{\v@leur}{-\v@lmax}{\v@lmax}%
    \maxim@m{\delt@}{-\mili@u}{\mili@u}\maxim@m{\v@leur}{\v@leur}{\delt@}%
    \maxim@m{\delt@}{-\v@lmin}{\v@lmin}\maxim@m{\v@leur}{\v@leur}{\delt@}%
    \ifdim\v@leur>\c@rre\invers@{\v@leur}{\v@leur}\edef\Uns@rM@x{\repdecn@mb{\v@leur}}%
    \v@lmax=\Uns@rM@x\v@lmax\mili@u=\Uns@rM@x\mili@u\v@lmin=\Uns@rM@x\v@lmin%
    \maxim@m{\v@leur}{-\v@lmax}{\v@lmax}\ifdim\v@leur<\c@rre%
    \maxim@m{\v@leur}{-\mili@u}{\mili@u}\ifdim\v@leur<\c@rre\else%
    \invers@{\mili@u}{\mili@u}\v@leur=-0.5\v@lmin%
    \v@leur=\repdecn@mb{\mili@u}\v@leur\m@jBBB@x{\v@leur}{#1}{#2}(#3,#4,#5,#6)\fi%
    \else\delt@=\repdecn@mb{\mili@u}\mili@u\v@leur=\repdecn@mb{\v@lmax}\v@lmin%
    \advance\delt@-\v@leur\ifdim\delt@<\z@\else\invers@{\v@lmax}{\v@lmax}%
    \edef\Uns@rAp{\repdecn@mb{\v@lmax}}\sqrt@{\delt@}{\delt@}%
    \v@leur=-\mili@u\advance\v@leur\delt@\v@leur=\Uns@rAp\v@leur%
    \m@jBBB@x{\v@leur}{#1}{#2}(#3,#4,#5,#6)%
    \v@leur=-\mili@u\advance\v@leur-\delt@\v@leur=\Uns@rAp\v@leur%
    \m@jBBB@x{\v@leur}{#1}{#2}(#3,#4,#5,#6)\fi\fi\fi}}
\ctr@ld@f\def\m@jBBB@x#1#2#3(#4,#5,#6,#7){{\relax\ifdim#1>\z@\ifdim#1<\p@%
    \edef\T@{\repdecn@mb{#1}}\v@lX=\p@\advance\v@lX-#1\edef\UNmT@{\repdecn@mb{\v@lX}}%
    \v@lX=#4\v@lY=#5\v@lZ=#6\v@lXa=#7\v@lX=\UNmT@\v@lX\advance\v@lX\T@\v@lY%
    \v@lY=\UNmT@\v@lY\advance\v@lY\T@\v@lZ\v@lZ=\UNmT@\v@lZ\advance\v@lZ\T@\v@lXa%
    \v@lX=\UNmT@\v@lX\advance\v@lX\T@\v@lY\v@lY=\UNmT@\v@lY\advance\v@lY\T@\v@lZ%
    \v@lX=\UNmT@\v@lX\advance\v@lX\T@\v@lY%
    \ifcase#2\or\v@lY=#3\or\v@lY=\v@lX\v@lX=#3\fi\b@undb@x{\v@lX}{\v@lY}\fi\fi}}
\ctr@ld@f\def\PsB@zier#1[#2]{{\f@gnewpath%
    \s@mme=\z@\def\list@num{#2,0}\extrairelepremi@r\p@int\de\list@num%
    \PSwrit@cmdS{\p@int}{\c@mmoveto}{\fwf@g}{\X@un}{\Y@un}\p@rtent=#1\bclB@zier}}
\ctr@ld@f\def\bclB@zier{\relax%
    \ifnum\s@mme<\p@rtent\advance\s@mme\@ne\BdingB@xfalse%
    \extrairelepremi@r\p@int\de\list@num\PSwrit@cmdS{\p@int}{}{\fwf@g}{\X@de}{\Y@de}%
    \extrairelepremi@r\p@int\de\list@num\PSwrit@cmdS{\p@int}{}{\fwf@g}{\X@tr}{\Y@tr}%
    \BdingB@xtrue%
    \extrairelepremi@r\p@int\de\list@num\PSwrit@cmdS{\p@int}{\c@mcurveto}{\fwf@g}{\X@qu}{\Y@qu}%
    \B@zierBB@x{1}{\Y@un}(\X@un,\X@de,\X@tr,\X@qu)%
    \B@zierBB@x{2}{\X@un}(\Y@un,\Y@de,\Y@tr,\Y@qu)%
    \edef\X@un{\X@qu}\edef\Y@un{\Y@qu}\bclB@zier\fi}
\ctr@ln@m\psBezier
\ctr@ld@f\def\psBezierDD#1[#2]{\ifcurr@ntPS\ifps@cri%
    \PSc@mment{psBezierDD N arcs=#1, Control points=#2}%
    \iffillm@de\PsB@zier#1[#2]%
    \f@gfill%
    \else\PsB@zier#1[#2]\f@gstroke\fi%
    \PSc@mment{End psBezierDD}\fi\fi}
\ctr@ln@m\et@tpsBezierTD
\ctr@ld@f\def\psBezierTD#1[#2]{\ifcurr@ntPS\ifps@cri\s@uvc@ntr@l\et@tpsBezierTD%
    \PSc@mment{psBezierTD N arcs=#1, Control points=#2}%
    \iffillm@de\PsB@zierTD#1[#2]%
    \f@gfill%
    \else\PsB@zierTD#1[#2]\f@gstroke\fi%
    \PSc@mment{End psBezierTD}\resetc@ntr@l\et@tpsBezierTD\fi\fi}
\ctr@ld@f\def\PsB@zierTD#1[#2]{\ifnum\curr@ntproj<\tw@\PsB@zier#1[#2]\else\PsB@zier@TD#1[#2]\fi}
\ctr@ld@f\def\PsB@zier@TD#1[#2]{{\f@gnewpath%
    \s@mme=\z@\def\list@num{#2,0}\extrairelepremi@r\p@int\de\list@num%
    \let\c@lprojSP=\relax\setc@ntr@l{2}\Figptpr@j-7:/\p@int/%
    \PSwrit@cmd{-7}{\c@mmoveto}{\fwf@g}%
    \loop\ifnum\s@mme<#1\advance\s@mme\@ne\extrairelepremi@r\p@intun\de\list@num%
    \extrairelepremi@r\p@intde\de\list@num\extrairelepremi@r\p@inttr\de\list@num%
    \subB@zierTD\NBz@rcs[\p@int,\p@intun,\p@intde,\p@inttr]\edef\p@int{\p@inttr}\repeat}}
\ctr@ld@f\def\subB@zierTD#1[#2,#3,#4,#5]{\delt@=\p@\divide\delt@\NBz@rcs\v@lmin=\z@%
    {\Figg@tXY{-7}\edef\X@un{\the\v@lX}\edef\Y@un{\the\v@lY}%
    \s@mme=\z@\loop\ifnum\s@mme<#1\advance\s@mme\@ne%
    \v@leur=\v@lmin\advance\v@leur0.33333 \delt@\edef\unti@rs{\repdecn@mb{\v@leur}}%
    \v@leur=\v@lmin\advance\v@leur0.66666 \delt@\edef\deti@rs{\repdecn@mb{\v@leur}}%
    \advance\v@lmin\delt@\edef\trti@rs{\repdecn@mb{\v@lmin}}%
    \figptBezierTD-8::\trti@rs[#2,#3,#4,#5]\Figptpr@j-8:/-8/%
    \c@lsubBzarc\unti@rs,\deti@rs[#2,#3,#4,#5]\BdingB@xfalse%
    \PSwrit@cmdS{-4}{}{\fwf@g}{\X@de}{\Y@de}\PSwrit@cmdS{-3}{}{\fwf@g}{\X@tr}{\Y@tr}%
    \BdingB@xtrue\PSwrit@cmdS{-8}{\c@mcurveto}{\fwf@g}{\X@qu}{\Y@qu}%
    \B@zierBB@x{1}{\Y@un}(\X@un,\X@de,\X@tr,\X@qu)%
    \B@zierBB@x{2}{\X@un}(\Y@un,\Y@de,\Y@tr,\Y@qu)%
    \edef\X@un{\X@qu}\edef\Y@un{\Y@qu}\figptcopyDD-7:/-8/\repeat}}
\ctr@ld@f\def\NBz@rcs{2}
\ctr@ld@f\def\c@lsubBzarc#1,#2[#3,#4,#5,#6]{\figptBezierTD-5::#1[#3,#4,#5,#6]%
    \figptBezierTD-6::#2[#3,#4,#5,#6]\Figptpr@j-4:/-5/\Figptpr@j-5:/-6/%
    \figptscontrolDD-4[-7,-4,-5,-8]}
\ctr@ln@m\pscirc
\ctr@ld@f\def\pscircDD#1(#2){\ifcurr@ntPS\ifps@cri\PSc@mment{pscircDD Center=#1 (Radius=#2)}%
    \psarccircDD#1;#2(0,360)\PSc@mment{End pscircDD}\fi\fi}
\ctr@ld@f\def\pscircTD#1,#2,#3(#4){\ifcurr@ntPS\ifps@cri%
    \PSc@mment{pscircTD Center=#1,P1=#2,P2=#3 (Radius=#4)}%
    \psarccircTD#1,#2,#3;#4(0,360)\PSc@mment{End pscircTD}\fi\fi}
\ctr@ln@m\p@urcent
{\catcode`\%=12\gdef\p@urcent{
\ctr@ld@f\def\PSc@mment#1{\ifpsdebugmode\immediate\write\fwf@g{\p@urcent\space#1}\fi}
\ctr@ln@m\acc@louv \ctr@ln@m\acc@lfer
{\catcode`\[=1\catcode`\{=12\gdef\acc@louv[{}}
{\catcode`\]=2\catcode`\}=12\gdef\acc@lfer{}]]
\ctr@ld@f\def\PSdict@{\ifUse@llipse%
    \immediate\write\fwf@g{/ellipsedict 9 dict def ellipsedict /mtrx matrix put}%
    \immediate\write\fwf@g{/ellipse \acc@louv ellipsedict begin}%
    \immediate\write\fwf@g{ /endangle exch def /startangle exch def}%
    \immediate\write\fwf@g{ /yrad exch def /xrad exch def}%
    \immediate\write\fwf@g{ /rotangle exch def /y exch def /x exch def}%
    \immediate\write\fwf@g{ /savematrix mtrx currentmatrix def}%
    \immediate\write\fwf@g{ x y translate rotangle rotate xrad yrad scale}%
    \immediate\write\fwf@g{ 0 0 1 startangle endangle arc}%
    \immediate\write\fwf@g{ savematrix setmatrix end\acc@lfer def}%
    \fi\PShe@der{EndProlog}}
\ctr@ld@f\def\Pssetc@rve#1=#2|{\keln@mun#1|%
    \def\n@mref{r}\ifx\l@debut\n@mref\pssetroundness{#2}\else
    \immediate\write16{*** Unknown attribute: \BS@ psset curve(..., #1=...)}%
    \fi}
\ctr@ln@m\curv@roundness
\ctr@ld@f\def\pssetroundness#1{\edef\curv@roundness{#1}}
\ctr@ld@f\def\defaultroundness{0.2} 
\ctr@ln@m\pscurve
\ctr@ld@f\def\pscurveDD[#1]{{\ifcurr@ntPS\ifps@cri\PSc@mment{pscurveDD Points=#1}%
    \s@uvc@ntr@l\et@tpscurveDD%
    \iffillm@de\Psc@rveDD\curv@roundness[#1]%
    \f@gfill%
    \else\Psc@rveDD\curv@roundness[#1]\f@gstroke\fi%
    \PSc@mment{End pscurveDD}\resetc@ntr@l\et@tpscurveDD\fi\fi}}
\ctr@ld@f\def\pscurveTD[#1]{{\ifcurr@ntPS\ifps@cri%
    \PSc@mment{pscurveTD Points=#1}\s@uvc@ntr@l\et@tpscurveTD\let\c@lprojSP=\relax%
    \iffillm@de\Psc@rveTD\curv@roundness[#1]%
    \f@gfill%
    \else\Psc@rveTD\curv@roundness[#1]\f@gstroke\fi%
    \PSc@mment{End pscurveTD}\resetc@ntr@l\et@tpscurveTD\fi\fi}}
\ctr@ld@f\def\Psc@rveDD#1[#2]{%
    \def\list@num{#2}\extrairelepremi@r\Ak@\de\list@num%
    \extrairelepremi@r\Ai@\de\list@num\extrairelepremi@r\Aj@\de\list@num%
    \f@gnewpath\PSwrit@cmdS{\Ai@}{\c@mmoveto}{\fwf@g}{\X@un}{\Y@un}%
    \setc@ntr@l{2}\figvectPDD -1[\Ak@,\Aj@]%
    \@ecfor\Ak@:=\list@num\do{\figpttraDD-2:=\Ai@/#1,-1/\BdingB@xfalse%
       \PSwrit@cmdS{-2}{}{\fwf@g}{\X@de}{\Y@de}%
       \figvectPDD -1[\Ai@,\Ak@]\figpttraDD-2:=\Aj@/-#1,-1/%
       \PSwrit@cmdS{-2}{}{\fwf@g}{\X@tr}{\Y@tr}\BdingB@xtrue%
       \PSwrit@cmdS{\Aj@}{\c@mcurveto}{\fwf@g}{\X@qu}{\Y@qu}%
       \B@zierBB@x{1}{\Y@un}(\X@un,\X@de,\X@tr,\X@qu)%
       \B@zierBB@x{2}{\X@un}(\Y@un,\Y@de,\Y@tr,\Y@qu)%
       \edef\X@un{\X@qu}\edef\Y@un{\Y@qu}\edef\Ai@{\Aj@}\edef\Aj@{\Ak@}}}
\ctr@ld@f\def\Psc@rveTD#1[#2]{\ifnum\curr@ntproj<\tw@\Psc@rvePPTD#1[#2]\else\Psc@rveCPTD#1[#2]\fi}
\ctr@ld@f\def\Psc@rvePPTD#1[#2]{\setc@ntr@l{2}%
    \def\list@num{#2}\extrairelepremi@r\Ak@\de\list@num\Figptpr@j-5:/\Ak@/%
    \extrairelepremi@r\Ai@\de\list@num\Figptpr@j-3:/\Ai@/%
    \extrairelepremi@r\Aj@\de\list@num\Figptpr@j-4:/\Aj@/%
    \f@gnewpath\PSwrit@cmdS{-3}{\c@mmoveto}{\fwf@g}{\X@un}{\Y@un}%
    \figvectPDD -1[-5,-4]%
    \@ecfor\Ak@:=\list@num\do{\Figptpr@j-5:/\Ak@/\figpttraDD-2:=-3/#1,-1/%
       \BdingB@xfalse\PSwrit@cmdS{-2}{}{\fwf@g}{\X@de}{\Y@de}%
       \figvectPDD -1[-3,-5]\figpttraDD-2:=-4/-#1,-1/%
       \PSwrit@cmdS{-2}{}{\fwf@g}{\X@tr}{\Y@tr}\BdingB@xtrue%
       \PSwrit@cmdS{-4}{\c@mcurveto}{\fwf@g}{\X@qu}{\Y@qu}%
       \B@zierBB@x{1}{\Y@un}(\X@un,\X@de,\X@tr,\X@qu)%
       \B@zierBB@x{2}{\X@un}(\Y@un,\Y@de,\Y@tr,\Y@qu)%
       \edef\X@un{\X@qu}\edef\Y@un{\Y@qu}\figptcopyDD-3:/-4/\figptcopyDD-4:/-5/}}
\ctr@ld@f\def\Psc@rveCPTD#1[#2]{\setc@ntr@l{2}%
    \def\list@num{#2}\extrairelepremi@r\Ak@\de\list@num%
    \extrairelepremi@r\Ai@\de\list@num\extrairelepremi@r\Aj@\de\list@num%
    \Figptpr@j-7:/\Ai@/%
    \f@gnewpath\PSwrit@cmd{-7}{\c@mmoveto}{\fwf@g}%
    \figvectPTD -9[\Ak@,\Aj@]%
    \@ecfor\Ak@:=\list@num\do{\figpttraTD-10:=\Ai@/#1,-9/%
       \figvectPTD -9[\Ai@,\Ak@]\figpttraTD-11:=\Aj@/-#1,-9/%
       \subB@zierTD\NBz@rcs[\Ai@,-10,-11,\Aj@]\edef\Ai@{\Aj@}\edef\Aj@{\Ak@}}}
\ctr@ld@f\def\psendfig{\ifcurr@ntPS\ifps@cri\immediate\closeout\fwf@g%
    \immediate\openout\fwf@g=\PSfilen@me\relax%
    \ifPDFm@ke\PSBdingB@x\else%
    \immediate\write\fwf@g{\p@urcent\string!PS-Adobe-2.0 EPSF-2.0}%
    \PShe@der{Creator\string: TeX (fig4tex.tex)}%
    \PShe@der{Title\string: \PSfilen@me}%
    \PShe@der{CreationDate\string: \the\day/\the\month/\the\year}%
    \PSBdingB@x%
    \PShe@der{EndComments}\PSdict@\fi%
    \immediate\write\fwf@g{\c@mgsave}%
    \openin\frf@g=\auxfilen@me\c@pypsfile\fwf@g\frf@g\closein\frf@g%
    \immediate\write\fwf@g{\c@mgrestore}%
    \PSc@mment{End of file.}\immediate\closeout\fwf@g%
    \immediate\openout\fwf@g=\auxfilen@me\immediate\closeout\fwf@g%
    \immediate\write16{File \PSfilen@me\space created.}\fi\fi\curr@ntPSfalse\ps@critrue}
\ctr@ld@f\def\PShe@der#1{\immediate\write\fwf@g{\p@urcent\p@urcent#1}}
\ctr@ld@f\def\PSBdingB@x{{\v@lX=\ptT@ptps\c@@rdXmin\v@lY=\ptT@ptps\c@@rdYmin%
     \v@lXa=\ptT@ptps\c@@rdXmax\v@lYa=\ptT@ptps\c@@rdYmax%
     \PShe@der{BoundingBox\string: \repdecn@mb{\v@lX}\space\repdecn@mb{\v@lY}%
     \space\repdecn@mb{\v@lXa}\space\repdecn@mb{\v@lYa}}}}
\ctr@ld@f\def\psfcconnect[#1]{{\ifcurr@ntPS\ifps@cri\PSc@mment{psfcconnect Points=#1}%
    \pssetfillmode{no}\s@uvc@ntr@l\et@tpsfcconnect\resetc@ntr@l{2}%
    \fcc@nnect@[#1]\resetc@ntr@l\et@tpsfcconnect\PSc@mment{End psfcconnect}\fi\fi}}
\ctr@ld@f\def\fcc@nnect@[#1]{\let\N@rm=\n@rmeucDD\def\list@num{#1}%
    \extrairelepremi@r\Ai@\de\list@num\edef\pr@m{\Ai@}\v@leur=\z@\p@rtent=\@ne\c@llgtot%
    \ifcase\fclin@typ@\edef\list@num{[\pr@m,#1,\Ai@}\expandafter\pscurve\list@num]%
    \else\ifdim\fclin@r@d\p@>\z@\Pslin@conge[#1]\else\psline[#1]\fi\fi%
    \v@leur=\@rrowp@s\v@leur\edef\list@num{#1,\Ai@,0}%
    \extrairelepremi@r\Ai@\de\list@num\mili@u=\epsil@n\c@llgpart%
    \advance\mili@u-\epsil@n\advance\mili@u-\delt@\advance\v@leur-\mili@u%
    \ifcase\fclin@typ@\invers@\mili@u\delt@%
    \ifnum\@rrowr@fpt>\z@\advance\delt@-\v@leur\v@leur=\delt@\fi%
    \v@leur=\repdecn@mb\v@leur\mili@u\edef\v@lt{\repdecn@mb\v@leur}%
    \extrairelepremi@r\Ak@\de\list@num%
    \figvectPDD-1[\pr@m,\Aj@]\figpttraDD-6:=\Ai@/\curv@roundness,-1/%
    \figvectPDD-1[\Ak@,\Ai@]\figpttraDD-7:=\Aj@/\curv@roundness,-1/%
    \delt@=\@rrowheadlength\p@\delt@=\C@AHANG\delt@\edef\R@dius{\repdecn@mb{\delt@}}%
    \ifcase\@rrowr@fpt%
    \FigptintercircB@zDD-8::\v@lt,\R@dius[\Ai@,-6,-7,\Aj@]\psarrowheadDD[-5,-8]\else%
    \FigptintercircB@zDD-8::\v@lt,\R@dius[\Aj@,-7,-6,\Ai@]\psarrowheadDD[-8,-5]\fi%
    \else\advance\delt@-\v@leur%
    \p@rtentiere{\p@rtent}{\delt@}\edef\C@efun{\the\p@rtent}%
    \p@rtentiere{\p@rtent}{\v@leur}\edef\C@efde{\the\p@rtent}%
    \figptbaryDD-5:[\Ai@,\Aj@;\C@efun,\C@efde]\ifcase\@rrowr@fpt%
    \delt@=\@rrowheadlength\unit@\delt@=\C@AHANG\delt@\edef\t@ille{\repdecn@mb{\delt@}}%
    \figvectPDD-2[\Ai@,\Aj@]\vecunit@{-2}{-2}\figpttraDD-5:=-5/\t@ille,-2/\fi%
    \psarrowheadDD[\Ai@,-5]\fi}
\ctr@ld@f\def\c@llgtot{\@ecfor\Aj@:=\list@num\do{\figvectP-1[\Ai@,\Aj@]\N@rm\delt@{-1}%
    \advance\v@leur\delt@\advance\p@rtent\@ne\edef\Ai@{\Aj@}}}
\ctr@ld@f\def\c@llgpart{\extrairelepremi@r\Aj@\de\list@num\figvectP-1[\Ai@,\Aj@]\N@rm\delt@{-1}%
    \advance\mili@u\delt@\ifdim\mili@u<\v@leur\edef\pr@m{\Ai@}\edef\Ai@{\Aj@}\c@llgpart\fi}
\ctr@ld@f\def\Pslin@conge[#1]{\ifnum\p@rtent>\tw@{\def\list@num{#1}%
    \extrairelepremi@r\Ai@\de\list@num\extrairelepremi@r\Aj@\de\list@num%
    \figptcopy-6:/\Ai@/\figvectP-3[\Ai@,\Aj@]\vecunit@{-3}{-3}\v@lmax=\result@t%
    \@ecfor\Ak@:=\list@num\do{\figvectP-4[\Aj@,\Ak@]\vecunit@{-4}{-4}%
    \minim@m\v@lmin\v@lmax\result@t\v@lmax=\result@t%
    \det@rm\delt@[-3,-4]\maxim@m\mili@u{\delt@}{-\delt@}\ifdim\mili@u>\Cepsil@n%
    \ifdim\delt@>\z@\figgetangleDD\Angl@[\Aj@,\Ak@,\Ai@]\else%
    \figgetangleDD\Angl@[\Aj@,\Ai@,\Ak@]\fi%
    \v@leur=\PI@deg\advance\v@leur-\Angl@\p@\divide\v@leur\tw@%
    \edef\Angl@{\repdecn@mb\v@leur}\c@ssin{\C@}{\S@}{\Angl@}\v@leur=\fclin@r@d\unit@%
    \v@leur=\S@\v@leur\mili@u=\C@\p@\invers@\mili@u\mili@u%
    \v@leur=\repdecn@mb{\mili@u}\v@leur%
    \minim@m\v@leur\v@leur\v@lmin\edef\t@ille{\repdecn@mb{\v@leur}}%
    \figpttra-5:=\Aj@/-\t@ille,-3/\psline[-6,-5]\figpttra-6:=\Aj@/\t@ille,-4/%
    \figvectNVDD-3[-3]\figvectNVDD-8[-4]\inters@cDD-7:[-5,-3;-6,-8]%
    \ifdim\delt@>\z@\psarccircP-7;\fclin@r@d[-5,-6]\else\psarccircP-7;\fclin@r@d[-6,-5]\fi%
    \else\psline[-6,\Aj@]\figptcopy-6:/\Aj@/\fi
    \edef\Ai@{\Aj@}\edef\Aj@{\Ak@}\figptcopy-3:/-4/}\psline[-6,\Aj@]}\else\psline[#1]\fi}
\ctr@ld@f\def\psfcnode[#1]#2{{\ifcurr@ntPS\ifps@cri\PSc@mment{psfcnode Points=#1}%
    \s@uvc@ntr@l\et@tpsfcnode\resetc@ntr@l{2}%
    \def\t@xt@{#2}\ifx\t@xt@\empty\def\g@tt@xt{\setbox\Gb@x=\hbox{\Figg@tT{\p@int}}}%
    \else\def\g@tt@xt{\setbox\Gb@x=\hbox{#2}}\fi%
    \v@lmin=\h@rdfcXp@dd\advance\v@lmin\Xp@dd\unit@\multiply\v@lmin\tw@%
    \v@lmax=\h@rdfcYp@dd\advance\v@lmax\Yp@dd\unit@\multiply\v@lmax\tw@%
    \Figv@ctCreg-8(\unit@,-\unit@)\def\list@num{#1}%
    \delt@=\curr@ntwidth bp\divide\delt@\tw@%
    \fcn@de\PSc@mment{End psfcnode}\resetc@ntr@l\et@tpsfcnode\fi\fi}}
\ctr@ld@f\def\d@butn@de{\g@tt@xt\v@lX=\wd\Gb@x%
    \v@lY=\ht\Gb@x\advance\v@lY\dp\Gb@x\advance\v@lX\v@lmin\advance\v@lY\v@lmax}
\ctr@ld@f\def\fcn@deE{%
    \@ecfor\p@int:=\list@num\do{\d@butn@de\v@lX=\unssqrttw@\v@lX\v@lY=\unssqrttw@\v@lY%
    \ifdim\thickn@ss\p@>\z@
    \v@lXa=\v@lX\advance\v@lXa\delt@\v@lXa=\ptT@unit@\v@lXa\edef\XR@d{\repdecn@mb\v@lXa}%
    \v@lYa=\v@lY\advance\v@lYa\delt@\v@lYa=\ptT@unit@\v@lYa\edef\YR@d{\repdecn@mb\v@lYa}%
    \arct@n\v@leur(\v@lXa,\v@lYa)\v@leur=\rdT@deg\v@leur\edef\@nglde{\repdecn@mb\v@leur}%
    {\c@lptellDD-2::\p@int;\XR@d,\YR@d(\@nglde)}
    \advance\v@leur-\PI@deg\edef\@nglun{\repdecn@mb\v@leur}%
    {\c@lptellDD-3::\p@int;\XR@d,\YR@d(\@nglun)}%
    \figptstra-6=-3,-2,\p@int/\thickn@ss,-8/\pssetfillmode{yes}\us@secondC@lor%
    \psline[-2,-3,-6,-5]\psarcell-4;\XR@d,\YR@d(\@nglun,\@nglde,0)\fi
    \v@lX=\ptT@unit@\v@lX\v@lY=\ptT@unit@\v@lY%
    \edef\XR@d{\repdecn@mb\v@lX}\edef\YR@d{\repdecn@mb\v@lY}%
    \pssetfillmode{yes}\us@thirdC@lor\psarcell\p@int;\XR@d,\YR@d(0,360,0)%
    \pssetfillmode{no}\us@primarC@lor\psarcell\p@int;\XR@d,\YR@d(0,360,0)}}
\ctr@ld@f\def\fcn@deL{\delt@=\ptT@unit@\delt@\edef\t@ille{\repdecn@mb\delt@}%
    \@ecfor\p@int:=\list@num\do{\Figg@tXYa{\p@int}\d@butn@de%
    \ifdim\v@lX>\v@lY\itis@Ktrue\else\itis@Kfalse\fi%
    \advance\v@lXa-\v@lX\Figp@intreg-1:(\v@lXa,\v@lYa)%
    \advance\v@lXa\v@lX\advance\v@lYa-\v@lY\Figp@intreg-2:(\v@lXa,\v@lYa)%
    \advance\v@lXa\v@lX\advance\v@lYa\v@lY\Figp@intreg-3:(\v@lXa,\v@lYa)%
    \advance\v@lXa-\v@lX\advance\v@lYa\v@lY\Figp@intreg-4:(\v@lXa,\v@lYa)%
    \ifdim\thickn@ss\p@>\z@\Figg@tXYa{\p@int}\pssetfillmode{yes}\us@secondC@lor
    \c@lpt@xt{-1}{-4}\c@lpt@xt@\v@lXa\v@lYa\v@lX\v@lY\c@rre\delt@%
    \Figp@intregDD-9:(\v@lZ,\v@lYa)\Figp@intregDD-11:(\v@lZa,\v@lYa)%
    \c@lpt@xt{-4}{-3}\c@lpt@xt@\v@lYa\v@lXa\v@lY\v@lX\delt@\c@rre%
    \Figp@intregDD-12:(\v@lXa,\v@lZ)\Figp@intregDD-10:(\v@lXa,\v@lZa)%
    \ifitis@K\figptstra-7=-9,-10,-11/\thickn@ss,-8/\psline[-9,-11,-5,-6,-7]\else%
    \figptstra-7=-10,-11,-12/\thickn@ss,-8/\psline[-10,-12,-5,-6,-7]\fi\fi
    \pssetfillmode{yes}\us@thirdC@lor\psline[-1,-2,-3,-4]%
    \pssetfillmode{no}\us@primarC@lor\psline[-1,-2,-3,-4,-1]}}
\ctr@ld@f\def\c@lpt@xt#1#2{\figvectN-7[#1,#2]\vecunit@{-7}{-7}\figpttra-5:=#1/\t@ille,-7/%
    \figvectP-7[#1,#2]\Figg@tXY{-7}\c@rre=\v@lX\delt@=\v@lY\Figg@tXY{-5}}
\ctr@ld@f\def\c@lpt@xt@#1#2#3#4#5#6{\v@lZ=#6\invers@{\v@lZ}{\v@lZ}\v@leur=\repdecn@mb{#5}\v@lZ%
    \v@lZ=#2\advance\v@lZ-#4\mili@u=\repdecn@mb{\v@leur}\v@lZ%
    \v@lZ=#3\advance\v@lZ\mili@u\v@lZa=-\v@lZ\advance\v@lZa\tw@#1}
\ctr@ld@f\def\fcn@deR{\@ecfor\p@int:=\list@num\do{\Figg@tXYa{\p@int}\d@butn@de%
    \advance\v@lXa-0.5\v@lX\advance\v@lYa-0.5\v@lY\Figp@intreg-1:(\v@lXa,\v@lYa)%
    \advance\v@lXa\v@lX\Figp@intreg-2:(\v@lXa,\v@lYa)%
    \advance\v@lYa\v@lY\Figp@intreg-3:(\v@lXa,\v@lYa)%
    \advance\v@lXa-\v@lX\Figp@intreg-4:(\v@lXa,\v@lYa)%
    \ifdim\thickn@ss\p@>\z@\pssetfillmode{yes}\us@secondC@lor
    \Figv@ctCreg-5(-\delt@,-\delt@)\figpttra-9:=-1/1,-5/%
    \Figv@ctCreg-5(\delt@,-\delt@)\figpttra-10:=-2/1,-5/%
    \Figv@ctCreg-5(\delt@,\delt@)\figpttra-11:=-3/1,-5/%
    \figptstra-7=-9,-10,-11/\thickn@ss,-8/\psline[-9,-11,-5,-6,-7]\fi
    \pssetfillmode{yes}\us@thirdC@lor\psline[-1,-2,-3,-4]%
    \pssetfillmode{no}\us@primarC@lor\psline[-1,-2,-3,-4,-1]}}
\ctr@ln@m\@rrowp@s
\ctr@ln@m\Xp@dd     \ctr@ln@m\Yp@dd
\ctr@ln@m\fclin@r@d \ctr@ln@m\thickn@ss
\ctr@ld@f\def\Pssetfl@wchart#1=#2|{\keln@mtr#1|%
    \def\n@mref{arr}\ifx\l@debut\n@mref\expandafter\keln@mtr\l@suite|%
     \def\n@mref{owp}\ifx\l@debut\n@mref\edef\@rrowp@s{#2}\else
     \def\n@mref{owr}\ifx\l@debut\n@mref\setfcr@fpt#2|\else
     \immediate\write16{*** Unknown attribute: \BS@ psset flowchart(..., #1=...)}%
     \fi\fi\else%
    \def\n@mref{lin}\ifx\l@debut\n@mref\setfccurv@#2|\else
    \def\n@mref{pad}\ifx\l@debut\n@mref\edef\Xp@dd{#2}\edef\Yp@dd{#2}\else
    \def\n@mref{rad}\ifx\l@debut\n@mref\edef\fclin@r@d{#2}\else
    \def\n@mref{sha}\ifx\l@debut\n@mref\setfcshap@#2|\else
    \def\n@mref{thi}\ifx\l@debut\n@mref\edef\thickn@ss{#2}\else
    \def\n@mref{xpa}\ifx\l@debut\n@mref\edef\Xp@dd{#2}\else
    \def\n@mref{ypa}\ifx\l@debut\n@mref\edef\Yp@dd{#2}\else
    \immediate\write16{*** Unknown attribute: \BS@ psset flowchart(..., #1=...)}%
    \fi\fi\fi\fi\fi\fi\fi\fi}
\ctr@ln@m\@rrowr@fpt \ctr@ln@m\fclin@typ@
\ctr@ld@f\def\setfcr@fpt#1#2|{\if#1e\def\@rrowr@fpt{1}\else\def\@rrowr@fpt{0}\fi}
\ctr@ld@f\def\setfccurv@#1#2|{\if#1c\def\fclin@typ@{0}\else\def\fclin@typ@{1}\fi}
\ctr@ln@m\h@rdfcXp@dd \ctr@ln@m\h@rdfcYp@dd
\ctr@ln@m\fcn@de \ctr@ln@m\fcsh@pe
\ctr@ld@f\def\setfcshap@#1#2|{%
    \if#1e\let\fcn@de=\fcn@deE\def\h@rdfcXp@dd{4pt}\def\h@rdfcYp@dd{4pt}%
     \edef\fcsh@pe{ellipse}\else%
    \if#1l\let\fcn@de=\fcn@deL\def\h@rdfcXp@dd{4pt}\def\h@rdfcYp@dd{4pt}%
     \edef\fcsh@pe{lozenge}\else%
          \let\fcn@de=\fcn@deR\def\h@rdfcXp@dd{6pt}\def\h@rdfcYp@dd{6pt}%
     \edef\fcsh@pe{rectangle}\fi\fi}
\ctr@ld@f\def\psline[#1]{{\ifcurr@ntPS\ifps@cri\PSc@mment{psline Points=#1}%
    \let\pslign@=\Pslign@P\Pslin@{#1}\PSc@mment{End psline}\fi\fi}}
\ctr@ld@f\def\pslineF#1{{\ifcurr@ntPS\ifps@cri\PSc@mment{pslineF Filename=#1}%
    \let\pslign@=\Pslign@F\Pslin@{#1}\PSc@mment{End pslineF}\fi\fi}}
\ctr@ld@f\def\pslineC(#1){{\ifcurr@ntPS\ifps@cri\PSc@mment{pslineC}%
    \let\pslign@=\Pslign@C\Pslin@{#1}\PSc@mment{End pslineC}\fi\fi}}
\ctr@ld@f\def\Pslin@#1{\iffillm@de\pslign@{#1}%
    \f@gfill%
    \else\pslign@{#1}\ifx\derp@int\premp@int%
    \f@gclosestroke%
    \else\f@gstroke\fi\fi}
\ctr@ld@f\def\Pslign@P#1{\def\list@num{#1}\extrairelepremi@r\p@int\de\list@num%
    \edef\premp@int{\p@int}\f@gnewpath%
    \PSwrit@cmd{\p@int}{\c@mmoveto}{\fwf@g}%
    \@ecfor\p@int:=\list@num\do{\PSwrit@cmd{\p@int}{\c@mlineto}{\fwf@g}%
    \edef\derp@int{\p@int}}}
\ctr@ld@f\def\Pslign@F#1{\s@uvc@ntr@l\et@tPslign@F\setc@ntr@l{2}\openin\frf@g=#1\relax%
    \ifeof\frf@g\message{*** File #1 not found !}\end\else%
    \read\frf@g to\tr@c\edef\premp@int{\tr@c}\expandafter\extr@ctCF\tr@c:%
    \f@gnewpath\PSwrit@cmd{-1}{\c@mmoveto}{\fwf@g}%
    \loop\read\frf@g to\tr@c\ifeof\frf@g\mored@tafalse\else\mored@tatrue\fi%
    \ifmored@ta\expandafter\extr@ctCF\tr@c:\PSwrit@cmd{-1}{\c@mlineto}{\fwf@g}%
    \edef\derp@int{\tr@c}\repeat\fi\closein\frf@g\resetc@ntr@l\et@tPslign@F}
\ctr@ln@m\extr@ctCF
\ctr@ld@f\def\extr@ctCFDD#1 #2:{\v@lX=#1\unit@\v@lY=#2\unit@\Figp@intregDD-1:(\v@lX,\v@lY)}
\ctr@ld@f\def\extr@ctCFTD#1 #2 #3:{\v@lX=#1\unit@\v@lY=#2\unit@\v@lZ=#3\unit@%
    \Figp@intregTD-1:(\v@lX,\v@lY,\v@lZ)}
\ctr@ld@f\def\Pslign@C#1{\s@uvc@ntr@l\et@tPslign@C\setc@ntr@l{2}%
    \def\list@num{#1}\extrairelepremi@r\p@int\de\list@num%
    \edef\premp@int{\p@int}\f@gnewpath%
    \expandafter\Pslign@C@\p@int:\PSwrit@cmd{-1}{\c@mmoveto}{\fwf@g}%
    \@ecfor\p@int:=\list@num\do{\expandafter\Pslign@C@\p@int:%
    \PSwrit@cmd{-1}{\c@mlineto}{\fwf@g}\edef\derp@int{\p@int}}%
    \resetc@ntr@l\et@tPslign@C}
\ctr@ld@f\def\Pslign@C@#1 #2:{{\def\t@xt@{#1}\ifx\t@xt@\empty\Pslign@C@#2:
    \else\extr@ctCF#1 #2:\fi}}
\ctr@ln@m\c@ntrolmesh
\ctr@ld@f\def\Pssetm@sh#1=#2|{\keln@mun#1|%
    \def\n@mref{d}\ifx\l@debut\n@mref\pssetmeshdiag{#2}\else
    \immediate\write16{*** Unknown attribute: \BS@ psset mesh(..., #1=...)}%
    \fi}
\ctr@ld@f\def\pssetmeshdiag#1{\edef\c@ntrolmesh{#1}}
\ctr@ld@f\def\defaultmeshdiag{0}    
\ctr@ld@f\def\psmesh#1,#2[#3,#4,#5,#6]{{\ifcurr@ntPS\ifps@cri%
    \PSc@mment{psmesh N1=#1, N2=#2, Quadrangle=[#3,#4,#5,#6]}%
    \s@uvc@ntr@l\et@tpsmesh\Pss@tsecondSt\setc@ntr@l{2}%
    \ifnum#1>\@ne\Psmeshp@rt#1[#3,#4,#5,#6]\fi%
    \ifnum#2>\@ne\Psmeshp@rt#2[#4,#5,#6,#3]\fi%
    \ifnum\c@ntrolmesh>\z@\Psmeshdi@g#1,#2[#3,#4,#5,#6]\fi%
    \ifnum\c@ntrolmesh<\z@\Psmeshdi@g#2,#1[#4,#5,#6,#3]\fi\Psrest@reSt%
    \psline[#3,#4,#5,#6,#3]\PSc@mment{End psmesh}\resetc@ntr@l\et@tpsmesh\fi\fi}}
\ctr@ld@f\def\Psmeshp@rt#1[#2,#3,#4,#5]{{\l@mbd@un=\@ne\l@mbd@de=#1\loop%
    \ifnum\l@mbd@un<#1\advance\l@mbd@de\m@ne\figptbary-1:[#2,#3;\l@mbd@de,\l@mbd@un]%
    \figptbary-2:[#5,#4;\l@mbd@de,\l@mbd@un]\psline[-1,-2]\advance\l@mbd@un\@ne\repeat}}
\ctr@ld@f\def\Psmeshdi@g#1,#2[#3,#4,#5,#6]{\figptcopy-2:/#3/\figptcopy-3:/#6/%
    \l@mbd@un=\z@\l@mbd@de=#1\loop\ifnum\l@mbd@un<#1%
    \advance\l@mbd@un\@ne\advance\l@mbd@de\m@ne\figptcopy-1:/-2/\figptcopy-4:/-3/%
    \figptbary-2:[#3,#4;\l@mbd@de,\l@mbd@un]%
    \figptbary-3:[#6,#5;\l@mbd@de,\l@mbd@un]\Psmeshdi@gp@rt#2[-1,-2,-3,-4]\repeat}
\ctr@ld@f\def\Psmeshdi@gp@rt#1[#2,#3,#4,#5]{{\l@mbd@un=\z@\l@mbd@de=#1\loop%
    \ifnum\l@mbd@un<#1\figptbary-5:[#2,#5;\l@mbd@de,\l@mbd@un]%
    \advance\l@mbd@de\m@ne\advance\l@mbd@un\@ne%
    \figptbary-6:[#3,#4;\l@mbd@de,\l@mbd@un]\psline[-5,-6]\repeat}}
\ctr@ln@m\psnormal
\ctr@ld@f\def\psnormalDD#1,#2[#3,#4]{{\ifcurr@ntPS\ifps@cri%
    \PSc@mment{psnormal Length=#1, Lambda=#2 [Pt1,Pt2]=[#3,#4]}%
    \s@uvc@ntr@l\et@tpsnormal\resetc@ntr@l{2}\figptendnormal-6::#1,#2[#3,#4]%
    \figptcopyDD-5:/-1/\psarrow[-5,-6]%
    \PSc@mment{End psnormal}\resetc@ntr@l\et@tpsnormal\fi\fi}}
\ctr@ld@f\def\psreset#1{\trtlis@rg{#1}{\Psreset@}}
\ctr@ld@f\def\Psreset@#1|{\keln@mde#1|%
    \def\n@mref{ar}\ifx\l@debut\n@mref\psresetarrowhead\else
    \def\n@mref{cu}\ifx\l@debut\n@mref\psset curve(roundness=\defaultroundness)\else
    \def\n@mref{fi}\ifx\l@debut\n@mref\psset (color=\defaultcolor,dash=\defaultdash,%
         fill=\defaultfill,join=\defaultjoin,width=\defaultwidth)\else
    \def\n@mref{fl}\ifx\l@debut\n@mref\psset flowchart(arrowp=\defaultfcarrowposition,%
 arrowr=\defaultfcarrowrefpt,line=\defaultfcline,xpadd=\defaultfcxpadding,%
 ypadd=\defaultfcypadding,radius=\defaultfcradius,shape=\defaultfcshape,%
 thick=\defaultfcthickness)\else
    \def\n@mref{me}\ifx\l@debut\n@mref\psset mesh(diag=\defaultmeshdiag)\else
    \def\n@mref{se}\ifx\l@debut\n@mref\psresetsecondsettings\else
    \def\n@mref{th}\ifx\l@debut\n@mref\psset third(color=\defaultthirdcolor)\else
    \immediate\write16{*** Unknown keyword #1 (\BS@ psreset).}%
    \fi\fi\fi\fi\fi\fi\fi}
\ctr@ld@f\def\psset#1(#2){\def\t@xt@{#1}\ifx\t@xt@\empty\trtlis@rg{#2}{\Pssetf@rst}
    \else\keln@mde#1|%
    \def\n@mref{ar}\ifx\l@debut\n@mref\trtlis@rg{#2}{\Psset@rrowhe@d}\else
    \def\n@mref{cu}\ifx\l@debut\n@mref\trtlis@rg{#2}{\Pssetc@rve}\else
    \def\n@mref{fi}\ifx\l@debut\n@mref\trtlis@rg{#2}{\Pssetf@rst}\else
    \def\n@mref{fl}\ifx\l@debut\n@mref\trtlis@rg{#2}{\Pssetfl@wchart}\else
    \def\n@mref{me}\ifx\l@debut\n@mref\trtlis@rg{#2}{\Pssetm@sh}\else
    \def\n@mref{se}\ifx\l@debut\n@mref\trtlis@rg{#2}{\Pssets@cond}\else
    \def\n@mref{th}\ifx\l@debut\n@mref\trtlis@rg{#2}{\Pssetth@rd}\else
    \immediate\write16{*** Unknown keyword: \BS@ psset #1(...)}%
    \fi\fi\fi\fi\fi\fi\fi\fi}
\ctr@ld@f\def\pssetdefault#1(#2){\ifcurr@ntPS\immediate\write16{*** \BS@ pssetdefault is ignored
    inside a \BS@ psbeginfig-\BS@ psendfig block.}%
    \immediate\write16{*** It must be called before \BS@ psbeginfig.}\else%
    \def\t@xt@{#1}\ifx\t@xt@\empty\trtlis@rg{#2}{\Pssd@f@rst}\else\keln@mde#1|%
    \def\n@mref{ar}\ifx\l@debut\n@mref\trtlis@rg{#2}{\Pssd@@rrowhe@d}\else
    \def\n@mref{cu}\ifx\l@debut\n@mref\trtlis@rg{#2}{\Pssd@c@rve}\else
    \def\n@mref{fi}\ifx\l@debut\n@mref\trtlis@rg{#2}{\Pssd@f@rst}\else
    \def\n@mref{fl}\ifx\l@debut\n@mref\trtlis@rg{#2}{\Pssd@fl@wchart}\else
    \def\n@mref{me}\ifx\l@debut\n@mref\trtlis@rg{#2}{\Pssd@m@sh}\else
    \def\n@mref{se}\ifx\l@debut\n@mref\trtlis@rg{#2}{\Pssd@s@cond}\else
    \def\n@mref{th}\ifx\l@debut\n@mref\trtlis@rg{#2}{\Pssd@th@rd}\else
    \immediate\write16{*** Unknown keyword: \BS@ pssetdefault #1(...)}%
    \fi\fi\fi\fi\fi\fi\fi\fi\initpss@ttings\fi}
\ctr@ld@f\def\Pssd@f@rst#1=#2|{\keln@mun#1|%
    \def\n@mref{c}\ifx\l@debut\n@mref\edef\defaultcolor{#2}\else
    \def\n@mref{d}\ifx\l@debut\n@mref\edef\defaultdash{#2}\else
    \def\n@mref{f}\ifx\l@debut\n@mref\edef\defaultfill{#2}\else
    \def\n@mref{j}\ifx\l@debut\n@mref\edef\defaultjoin{#2}\else
    \def\n@mref{u}\ifx\l@debut\n@mref\edef\defaultupdate{#2}\pssetupdate{#2}\else
    \def\n@mref{w}\ifx\l@debut\n@mref\edef\defaultwidth{#2}\else
    \immediate\write16{*** Unknown attribute: \BS@ pssetdefault (..., #1=...)}%
    \fi\fi\fi\fi\fi\fi}
\ctr@ld@f\def\Pssd@@rrowhe@d#1=#2|{\keln@mun#1|%
    \def\n@mref{a}\ifx\l@debut\n@mref\edef\defaultarrowheadangle{#2}\else
    \def\n@mref{f}\ifx\l@debut\n@mref\edef\defaultarrowheadangle{#2}\else
    \def\n@mref{l}\ifx\l@debut\n@mref\y@tiunit{#2}\ifunitpr@sent%
     \edef\defaulth@rdahlength{#2}\else\edef\defaulth@rdahlength{#2pt}%
     \message{*** \BS@ pssetdefault (..., #1=#2, ...) : unit is missing, pt is assumed.}%
     \fi\else
    \def\n@mref{o}\ifx\l@debut\n@mref\edef\defaultarrowheadout{#2}\else
    \def\n@mref{r}\ifx\l@debut\n@mref\edef\defaultarrowheadratio{#2}\else
    \immediate\write16{*** Unknown attribute: \BS@ pssetdefault arrowhead(..., #1=...)}%
    \fi\fi\fi\fi\fi}
\ctr@ld@f\def\Pssd@c@rve#1=#2|{\keln@mun#1|%
    \def\n@mref{r}\ifx\l@debut\n@mref\edef\defaultroundness{#2}\else%
    \immediate\write16{*** Unknown attribute: \BS@ pssetdefault curve(..., #1=...)}%
    \fi}
\ctr@ld@f\def\Pssd@fl@wchart#1=#2|{\keln@mtr#1|%
    \def\n@mref{arr}\ifx\l@debut\n@mref\expandafter\keln@mtr\l@suite|%
     \def\n@mref{owp}\ifx\l@debut\n@mref\edef\defaultfcarrowposition{#2}\else
     \def\n@mref{owr}\ifx\l@debut\n@mref\edef\defaultfcarrowrefpt{#2}\else
     \immediate\write16{*** Unknown attribute: \BS@ pssetdefault flowchart(..., #1=...)}%
     \fi\fi\else%
    \def\n@mref{lin}\ifx\l@debut\n@mref\edef\defaultfcline{#2}\else
    \def\n@mref{pad}\ifx\l@debut\n@mref\edef\defaultfcxpadding{#2}%
                    \edef\defaultfcypadding{#2}\else
    \def\n@mref{rad}\ifx\l@debut\n@mref\edef\defaultfcradius{#2}\else
    \def\n@mref{sha}\ifx\l@debut\n@mref\edef\defaultfcshape{#2}\else
    \def\n@mref{thi}\ifx\l@debut\n@mref\edef\defaultfcthickness{#2}\else
    \def\n@mref{xpa}\ifx\l@debut\n@mref\edef\defaultfcxpadding{#2}\else
    \def\n@mref{ypa}\ifx\l@debut\n@mref\edef\defaultfcypadding{#2}\else
    \immediate\write16{*** Unknown attribute: \BS@ pssetdefault flowchart(..., #1=...)}%
    \fi\fi\fi\fi\fi\fi\fi\fi}
\ctr@ld@f\def\defaultfcarrowposition{0.5}
\ctr@ld@f\def\defaultfcarrowrefpt{start}
\ctr@ld@f\def\defaultfcline{polygon}
\ctr@ld@f\def\defaultfcradius{0}
\ctr@ld@f\def\defaultfcshape{rectangle}
\ctr@ld@f\def\defaultfcthickness{0}
\ctr@ld@f\def\defaultfcxpadding{0}
\ctr@ld@f\def\defaultfcypadding{0}
\ctr@ld@f\def\Pssd@m@sh#1=#2|{\keln@mun#1|%
    \def\n@mref{d}\ifx\l@debut\n@mref\edef\defaultmeshdiag{#2}\else%
    \immediate\write16{*** Unknown attribute: \BS@ pssetdefault mesh(..., #1=...)}%
    \fi}
\ctr@ld@f\def\Pssd@s@cond#1=#2|{\keln@mun#1|%
    \def\n@mref{c}\ifx\l@debut\n@mref\edef\defaultsecondcolor{#2}\else%
    \def\n@mref{d}\ifx\l@debut\n@mref\edef\defaultseconddash{#2}\else%
    \def\n@mref{w}\ifx\l@debut\n@mref\edef\defaultsecondwidth{#2}\else%
    \immediate\write16{*** Unknown attribute: \BS@ pssetdefault second(..., #1=...)}%
    \fi\fi\fi}
\ctr@ld@f\def\Pssd@th@rd#1=#2|{\keln@mun#1|%
    \def\n@mref{c}\ifx\l@debut\n@mref\edef\defaultthirdcolor{#2}\else%
    \immediate\write16{*** Unknown attribute: \BS@ pssetdefault third(..., #1=...)}%
    \fi}
\ctr@ln@w{newif}\iffillm@de
\ctr@ld@f\def\pssetfillmode#1{\expandafter\setfillm@de#1:}
\ctr@ld@f\def\setfillm@de#1#2:{\if#1n\fillm@defalse\else\fillm@detrue\fi}
\ctr@ld@f\def\defaultfill{no}     
\ctr@ln@w{newif}\ifpsupdatem@de
\ctr@ld@f\def\pssetupdate#1{\ifcurr@ntPS\immediate\write16{*** \BS@ pssetupdate is ignored inside a
     \BS@ psbeginfig-\BS@ psendfig block.}%
    \immediate\write16{*** It must be called before \BS@ psbeginfig.}%
    \else\expandafter\setupd@te#1:\fi}
\ctr@ld@f\def\setupd@te#1#2:{\if#1n\psupdatem@defalse\else\psupdatem@detrue\fi}
\ctr@ld@f\def\defaultupdate{no}     
\ctr@ln@m\curr@ntcolor \ctr@ln@m\curr@ntcolorc@md
\ctr@ld@f\def\Pssetc@lor#1{\ifps@cri\result@tent=\@ne\expandafter\c@lnbV@l#1 :%
    \def\curr@ntcolor{}\def\curr@ntcolorc@md{}%
    \ifcase\result@tent\or\pssetgray{#1}\or\or\pssetrgb{#1}\or\pssetcmyk{#1}\fi\fi}
\ctr@ln@m\curr@ntcolorc@mdStroke
\ctr@ld@f\def\pssetcmyk#1{\ifps@cri\def\curr@ntcolor{#1}\def\curr@ntcolorc@md{\c@msetcmykcolor}%
    \def\curr@ntcolorc@mdStroke{\c@msetcmykcolorStroke}%
    \ifcurr@ntPS\PSc@mment{pssetcmyk Color=#1}\us@primarC@lor\fi\fi}
\ctr@ld@f\def\pssetrgb#1{\ifps@cri\def\curr@ntcolor{#1}\def\curr@ntcolorc@md{\c@msetrgbcolor}%
    \def\curr@ntcolorc@mdStroke{\c@msetrgbcolorStroke}%
    \ifcurr@ntPS\PSc@mment{pssetrgb Color=#1}\us@primarC@lor\fi\fi}
\ctr@ld@f\def\pssetgray#1{\ifps@cri\def\curr@ntcolor{#1}\def\curr@ntcolorc@md{\c@msetgray}%
    \def\curr@ntcolorc@mdStroke{\c@msetgrayStroke}%
    \ifcurr@ntPS\PSc@mment{pssetgray Gray level=#1}\us@primarC@lor\fi\fi}
\ctr@ln@m\fillc@md
\ctr@ld@f\def\us@primarC@lor{\immediate\write\fwf@g{\d@fprimarC@lor}%
    \let\fillc@md=\prfillc@md}
\ctr@ld@f\def\prfillc@md{\d@fprimarC@lor\space\c@mfill}
\ctr@ld@f\def\defaultcolor{0}       
\ctr@ld@f\def\c@lnbV@l#1 #2:{\def\t@xt@{#1}\relax\ifx\t@xt@\empty\c@lnbV@l#2:
    \else\c@lnbV@l@#1 #2:\fi}
\ctr@ld@f\def\c@lnbV@l@#1 #2:{\def\t@xt@{#2}\ifx\t@xt@\empty%
    \def\t@xt@{#1}\ifx\t@xt@\empty\advance\result@tent\m@ne\fi
    \else\advance\result@tent\@ne\c@lnbV@l@#2:\fi}
\ctr@ld@f\def\Blackcmyk{0 0 0 1}
\ctr@ld@f\def\Whitecmyk{0 0 0 0}
\ctr@ld@f\def\Cyancmyk{1 0 0 0}
\ctr@ld@f\def\Magentacmyk{0 1 0 0}
\ctr@ld@f\def\Yellowcmyk{0 0 1 0}
\ctr@ld@f\def\Redcmyk{0 1 1 0}
\ctr@ld@f\def\Greencmyk{1 0 1 0}
\ctr@ld@f\def\Bluecmyk{1 1 0 0}
\ctr@ld@f\def\Graycmyk{0 0 0 0.50}
\ctr@ld@f\def\BrickRedcmyk{0 0.89 0.94 0.28} 
\ctr@ld@f\def\Browncmyk{0 0.81 1 0.60} 
\ctr@ld@f\def\ForestGreencmyk{0.91 0 0.88 0.12} 
\ctr@ld@f\def\Goldenrodcmyk{ 0 0.10 0.84 0} 
\ctr@ld@f\def\Marooncmyk{0 0.87 0.68 0.32} 
\ctr@ld@f\def\Orangecmyk{0 0.61 0.87 0} 
\ctr@ld@f\def\Purplecmyk{0.45 0.86 0 0} 
\ctr@ld@f\def\RoyalBluecmyk{1. 0.50 0 0} 
\ctr@ld@f\def\Violetcmyk{0.79 0.88 0 0} 
\ctr@ld@f\def\Blackrgb{0 0 0}
\ctr@ld@f\def\Whitergb{1 1 1}
\ctr@ld@f\def\Redrgb{1 0 0}
\ctr@ld@f\def\Greenrgb{0 1 0}
\ctr@ld@f\def\Bluergb{0 0 1}
\ctr@ld@f\def\Cyanrgb{0 1 1}
\ctr@ld@f\def\Magentargb{1 0 1}
\ctr@ld@f\def\Yellowrgb{1 1 0}
\ctr@ld@f\def\Grayrgb{0.5 0.5 0.5}
\ctr@ld@f\def\Chocolatergb{0.824 0.412 0.118}
\ctr@ld@f\def\DarkGoldenrodrgb{0.722 0.525 0.043}
\ctr@ld@f\def\DarkOrangergb{1 0.549 0}
\ctr@ld@f\def\Firebrickrgb{0.698 0.133 0.133}
\ctr@ld@f\def\ForestGreenrgb{0.133 0.545 0.133}
\ctr@ld@f\def\Goldrgb{1 0.843 0}
\ctr@ld@f\def\HotPinkrgb{1 0.412 0.706}
\ctr@ld@f\def\Maroonrgb{0.690 0.188 0.376}
\ctr@ld@f\def\Pinkrgb{1 0.753 0.796}
\ctr@ld@f\def\RoyalBluergb{0.255 0.412 0.882}
\ctr@ld@f\def\Pssetf@rst#1=#2|{\keln@mun#1|%
    \def\n@mref{c}\ifx\l@debut\n@mref\Pssetc@lor{#2}\else
    \def\n@mref{d}\ifx\l@debut\n@mref\pssetdash{#2}\else
    \def\n@mref{f}\ifx\l@debut\n@mref\pssetfillmode{#2}\else
    \def\n@mref{j}\ifx\l@debut\n@mref\pssetjoin{#2}\else
    \def\n@mref{u}\ifx\l@debut\n@mref\pssetupdate{#2}\else
    \def\n@mref{w}\ifx\l@debut\n@mref\pssetwidth{#2}\else
    \immediate\write16{*** Unknown attribute: \BS@ psset (..., #1=...)}%
    \fi\fi\fi\fi\fi\fi}
\ctr@ln@m\curr@ntdash
\ctr@ld@f\def\s@uvdash#1{\edef#1{\curr@ntdash}}
\ctr@ld@f\def\defaultdash{1}        
\ctr@ld@f\def\pssetdash#1{\ifps@cri\edef\curr@ntdash{#1}\ifcurr@ntPS\expandafter\Pssetd@sh#1 :\fi\fi}
\ctr@ld@f\def\Pssetd@shI#1{\PSc@mment{pssetdash Index=#1}\ifcase#1%
    \or\immediate\write\fwf@g{[] 0 \c@msetdash}
    \or\immediate\write\fwf@g{[6 2] 0 \c@msetdash}
    \or\immediate\write\fwf@g{[4 2] 0 \c@msetdash}
    \or\immediate\write\fwf@g{[2 2] 0 \c@msetdash}
    \or\immediate\write\fwf@g{[1 2] 0 \c@msetdash}
    \or\immediate\write\fwf@g{[2 4] 0 \c@msetdash}
    \or\immediate\write\fwf@g{[3 5] 0 \c@msetdash}
    \or\immediate\write\fwf@g{[3 3] 0 \c@msetdash}
    \or\immediate\write\fwf@g{[3 5 1 5] 0 \c@msetdash}
    \or\immediate\write\fwf@g{[6 4 2 4] 0 \c@msetdash}
    \fi}
\ctr@ld@f\def\Pssetd@sh#1 #2:{{\def\t@xt@{#1}\ifx\t@xt@\empty\Pssetd@sh#2:
    \else\def\t@xt@{#2}\ifx\t@xt@\empty\Pssetd@shI{#1}\else\s@mme=\@ne\def\debutp@t{#1}%
    \an@lysd@sh#2:\ifodd\s@mme\edef\debutp@t{\debutp@t\space\finp@t}\def\finp@t{0}\fi%
    \PSc@mment{pssetdash Pattern=#1 #2}%
    \immediate\write\fwf@g{[\debutp@t] \finp@t\space\c@msetdash}\fi\fi}}
\ctr@ld@f\def\an@lysd@sh#1 #2:{\def\t@xt@{#2}\ifx\t@xt@\empty\def\finp@t{#1}\else%
    \edef\debutp@t{\debutp@t\space#1}\advance\s@mme\@ne\an@lysd@sh#2:\fi}
\ctr@ln@m\curr@ntwidth
\ctr@ld@f\def\s@uvwidth#1{\edef#1{\curr@ntwidth}}
\ctr@ld@f\def\defaultwidth{0.4}     
\ctr@ld@f\def\pssetwidth#1{\ifps@cri\edef\curr@ntwidth{#1}\ifcurr@ntPS%
    \PSc@mment{pssetwidth Width=#1}\immediate\write\fwf@g{#1 \c@msetlinewidth}\fi\fi}
\ctr@ln@m\curr@ntjoin
\ctr@ld@f\def\pssetjoin#1{\ifps@cri\edef\curr@ntjoin{#1}\ifcurr@ntPS\expandafter\Pssetj@in#1:\fi\fi}
\ctr@ld@f\def\Pssetj@in#1#2:{\PSc@mment{pssetjoin join=#1}%
    \if#1r\def\t@xt@{1}\else\if#1b\def\t@xt@{2}\else\def\t@xt@{0}\fi\fi%
    \immediate\write\fwf@g{\t@xt@\space\c@msetlinejoin}}
\ctr@ld@f\def\defaultjoin{miter}   
\ctr@ld@f\def\Pssets@cond#1=#2|{\keln@mun#1|%
    \def\n@mref{c}\ifx\l@debut\n@mref\Pssets@condcolor{#2}\else%
    \def\n@mref{d}\ifx\l@debut\n@mref\pssetseconddash{#2}\else%
    \def\n@mref{w}\ifx\l@debut\n@mref\pssetsecondwidth{#2}\else%
    \immediate\write16{*** Unknown attribute: \BS@ psset second(..., #1=...)}%
    \fi\fi\fi}
\ctr@ln@m\curr@ntseconddash
\ctr@ld@f\def\pssetseconddash#1{\edef\curr@ntseconddash{#1}}
\ctr@ld@f\def\defaultseconddash{4}  
\ctr@ln@m\curr@ntsecondwidth
\ctr@ld@f\def\pssetsecondwidth#1{\edef\curr@ntsecondwidth{#1}}
\ctr@ld@f\edef\defaultsecondwidth{\defaultwidth} 
\ctr@ld@f\def\psresetsecondsettings{%
    \pssetseconddash{\defaultseconddash}\pssetsecondwidth{\defaultsecondwidth}%
    \Pssets@condcolor{\defaultsecondcolor}}
\ctr@ln@m\sec@ndcolor \ctr@ln@m\sec@ndcolorc@md
\ctr@ld@f\def\Pssets@condcolor#1{\ifps@cri\result@tent=\@ne\expandafter\c@lnbV@l#1 :%
    \def\sec@ndcolor{}\def\sec@ndcolorc@md{}%
    \ifcase\result@tent\or\pssetsecondgray{#1}\or\or\pssetsecondrgb{#1}%
    \or\pssetsecondcmyk{#1}\fi\fi}
\ctr@ln@m\sec@ndcolorc@mdStroke
\ctr@ld@f\def\pssetsecondcmyk#1{\def\sec@ndcolor{#1}\def\sec@ndcolorc@md{\c@msetcmykcolor}%
    \def\sec@ndcolorc@mdStroke{\c@msetcmykcolorStroke}}
\ctr@ld@f\def\pssetsecondrgb#1{\def\sec@ndcolor{#1}\def\sec@ndcolorc@md{\c@msetrgbcolor}%
    \def\sec@ndcolorc@mdStroke{\c@msetrgbcolorStroke}}
\ctr@ld@f\def\pssetsecondgray#1{\def\sec@ndcolor{#1}\def\sec@ndcolorc@md{\c@msetgray}%
    \def\sec@ndcolorc@mdStroke{\c@msetgrayStroke}}
\ctr@ld@f\def\us@secondC@lor{\immediate\write\fwf@g{\d@fsecondC@lor}%
    \let\fillc@md=\sdfillc@md}
\ctr@ld@f\def\sdfillc@md{\d@fsecondC@lor\space\c@mfill}
\ctr@ld@f\edef\defaultsecondcolor{\defaultcolor} 
\ctr@ld@f\def\Pss@tsecondSt{%
    \s@uvdash{\typ@dash}\pssetdash{\curr@ntseconddash}%
    \s@uvwidth{\typ@width}\pssetwidth{\curr@ntsecondwidth}\us@secondC@lor}
\ctr@ld@f\def\Psrest@reSt{\pssetwidth{\typ@width}\pssetdash{\typ@dash}\us@primarC@lor}
\ctr@ld@f\def\Pssetth@rd#1=#2|{\keln@mun#1|%
    \def\n@mref{c}\ifx\l@debut\n@mref\Pssetth@rdcolor{#2}\else%
    \immediate\write16{*** Unknown attribute: \BS@ psset third(..., #1=...)}%
    \fi}
\ctr@ln@m\th@rdcolor \ctr@ln@m\th@rdcolorc@md
\ctr@ld@f\def\Pssetth@rdcolor#1{\ifps@cri\result@tent=\@ne\expandafter\c@lnbV@l#1 :%
    \def\th@rdcolor{}\def\th@rdcolorc@md{}%
    \ifcase\result@tent\or\Pssetth@rdgray{#1}\or\or\Pssetth@rdrgb{#1}%
    \or\Pssetth@rdcmyk{#1}\fi\fi}
\ctr@ln@m\th@rdcolorc@mdStroke
\ctr@ld@f\def\Pssetth@rdcmyk#1{\def\th@rdcolor{#1}\def\th@rdcolorc@md{\c@msetcmykcolor}%
    \def\th@rdcolorc@mdStroke{\c@msetcmykcolorStroke}}
\ctr@ld@f\def\Pssetth@rdrgb#1{\def\th@rdcolor{#1}\def\th@rdcolorc@md{\c@msetrgbcolor}%
    \def\th@rdcolorc@mdStroke{\c@msetrgbcolorStroke}}
\ctr@ld@f\def\Pssetth@rdgray#1{\def\th@rdcolor{#1}\def\th@rdcolorc@md{\c@msetgray}%
    \def\th@rdcolorc@mdStroke{\c@msetgrayStroke}}
\ctr@ld@f\def\us@thirdC@lor{\immediate\write\fwf@g{\d@fthirdC@lor}%
    \let\fillc@md=\thfillc@md}
\ctr@ld@f\def\thfillc@md{\d@fthirdC@lor\space\c@mfill}
\ctr@ld@f\def\defaultthirdcolor{1}  
\ctr@ld@f\def\pstrimesh#1[#2,#3,#4]{{\ifcurr@ntPS\ifps@cri%
    \PSc@mment{pstrimesh Type=#1, Triangle=[#2,#3,#4]}%
    \s@uvc@ntr@l\et@tpstrimesh\ifnum#1>\@ne\Pss@tsecondSt\setc@ntr@l{2}%
    \Pstrimeshp@rt#1[#2,#3,#4]\Pstrimeshp@rt#1[#3,#4,#2]%
    \Pstrimeshp@rt#1[#4,#2,#3]\Psrest@reSt\fi\psline[#2,#3,#4,#2]%
    \PSc@mment{End pstrimesh}\resetc@ntr@l\et@tpstrimesh\fi\fi}}
\ctr@ld@f\def\Pstrimeshp@rt#1[#2,#3,#4]{{\l@mbd@un=\@ne\l@mbd@de=#1\loop\ifnum\l@mbd@de>\@ne%
    \advance\l@mbd@de\m@ne\figptbary-1:[#2,#3;\l@mbd@de,\l@mbd@un]%
    \figptbary-2:[#2,#4;\l@mbd@de,\l@mbd@un]\psline[-1,-2]%
    \advance\l@mbd@un\@ne\repeat}}
\initpr@lim\initpss@ttings\initPDF@rDVI
\ctr@ln@w{newbox}\figBoxA
\ctr@ln@w{newbox}\figBoxB
\ctr@ln@w{newbox}\figBoxC
\catcode`\@=12

\bibliographystyle{plain}


\maketitle 
\begin{abstract}
This article is concerned with the maximal accretive realizations of
geometric Kramers-Fokker-Planck operators on manifolds with
boundaries. A general class of boundary conditions is
introduced which ensures the maximal accretivity and some global
subelliptic estimates. Those estimates imply nice spectral properties
as well as exponential decay properties for the associated semigroup.
Admissible boundary conditions cover a wide range of applications
for the usual scalar Kramer-Fokker-Planck equation or Bismut's
hypoelliptic laplacian.
\end{abstract}
\noindent\textbf{MSC2010: 35H10, 35H20, 58J32, 58J50, 58J65, 60J65}\\
\noindent\textbf{Keywords:} Kramers-Fokker-Planck equation, Langevin
process, hypoelliptic Laplacian, boundary value problem, subelliptic
estimates.
\tableofcontents{}
\section{Introduction}
\subsection{Motivations}
\label{se.motiv}
A few years ago, while I was visiting the WIAS in Berlin, Holger
Stephan asked me about the relevant boundary conditions for the
Kramers-Fokker-Planck operator
$$
p.\partial_{q}-(\partial_{q}V(q)).\partial_{p}+\frac{-\Delta_{p}+|p|^{2}}{2}
$$
in $\Omega\times \rz^{d}_{p}$ when $\Omega\subset \rz^{d}_{q}$ is a
bounded regular domain. At that time I was aware of the previous works
developed within the mathematical analysis of kinetic models (see
a.e. \cite{DeMG}\cite{Car}\cite{Luc}) but I noticed, and I have kept in mind
since, that the weak
formulation (of time-dependent problems) were far from providing the
accurate semigroup or resolvent information which was accessible in
the boundaryless case (see e.g. \cite{HerNi}\cite{HelNi}). After this
various things occurred:
\begin{itemize}
\item We realized recently with D.~Le~Peutrec and C.~Viterbo in
  \cite{LNV} that the introduction of artificial boundary value
  problems was an important step in the accurate spectral analysis of
  Witten Laplacians acting on $p$-forms in the low temperature limit.
Actually the harmonic forms of these artificial boundary value
problems, 
after Witten's deformation, 
play in the case of $p$-forms, $p>1$\,, the role of the truncated
version $\chi(q)e^{-V(q)/h}$ of the explicit equilibrium
density $e^{-V(q)/h}$ used in \cite{HKN} for the case $p=0$\,. 
The Kramers-Fokker-Planck operator is
actually the $p=0$ version of Bismut's hypoelliptic Laplacian. We
may hope to extend the accurate spectral analysis of the case $p=0$
done by H{\'e}rau-Hitrik-Sj{\"o}strand in \cite{HHS2}, if we understand the
suitable realizations of the hypoelliptic Laplacian on a manifold with boundary.
\item While we were working with T.~Leli{\`e}vre about quasi-stationnary
  distributions for the Einstein-Smoluchowski case (SDE in a gradient
  field),  with the help
  of boundary Witten Laplacians acting on functions ($p=0$) and
  $1$-forms, T.~Leli{\`e}vre repeatedly insisted that the Langevin process
  was a more natural process for molecular dynamics. The
  quasi-stationnary distributions are especially used to develop or
  study efficients algorithms for molecular dynamics. The PDE
  formulation of the Langevin process is the Kramers-Fokker-Planck
  equation and quasi-stationnary distributions especially make sense on
  bounded domains (in the position variable $q$).
\item A bit more than one year ago, we had discussions with F.~H{\'e}rau and
  D.~Le~Peutrec who elaborated an approach of those boundary value
  problems
relying on the weak formulation of \cite{Car} combined with the
so-called hypocoercive techniques proposed by C.~Villani in
\cite{Vil}. They were embarrassed with the definition of traces and I
suggested that the problem should be reconsidered from the beginning,
by trying to mimic what is done for elliptic operators with the
introduction of Calderon projectors. We continued on our different
ways while keeping in touch. Their point of view may be efficient for
nonlinear problems.
\item In 2007, G.~Lebeau besides his work \cite{BiLe} in collaboration
  with J.M.~Bismut about the hypoelliptic Laplacian, proved maximal
  estimates for  what he called the geometric 
(Kramers)-Fokker-Planck operator. This is actually the scalar
principal part of the hypoelliptic Laplacian. Those
  estimates have not been, seemingly,  employed seriously  up to now.
We shall see that they are instrumental in absorbing some
singular perturbations due to the curvature of the boundary.
\end{itemize}
This article presents the functional analysis of geometric
Kramers-Fokker-Planck operators with a rather wide and natural class of
boundary conditions. 
About the application to Bismut's hypoelliptic Laplacian
and if one compares with the program which has been
achieved for the accurate asymptotic and spectral analysis of 
boundary Witten Laplacians in \cite{HeNi} and \cite{Lep3} and its
applications in \cite{LNV} and \cite{LeNi}, it is only the
beginning. Nothing is said about the asymptotic analysis with respect
to small parameters, nor about the supersymmetric arguments so
effective when dealing with the Witten Laplacian. In the case
of the Kramers-Fokker-Planck equation two parameters can be
introduced with independent or correlated asymptotics, the
temperature and the friction (see for example \cite{Ris}\cite{HerNi}).\\
Other boundary conditions which are proposed as examples of
application, are heuristically introduced  by completing the Langevin
process with a jump process when particles hit the boundary. Except in
the case of specular reflection studied in \cite{Lap}\cite{BoJa} or
exotic one-dimensional problems with non elastic boundary conditions
 in \cite{Ber}, no definite result
seems available for justifying these heuristic arguments. Our
subelliptic (regularity and decay) estimates of the corresponding
boundary value problem, may help to get a better understanding.\\
Finally one may wonder whether it is necessary to work in the general
framework of riemannian geometry. Our approach passes through the
local reduction to straight half-space problem where the coordinate
changes prevent us from sticking to the euclidean case. 
Moreover, even when the metric is flat (the Riemann
curvature tensor vanishes), the extrinsic curvature of the boundary
raises crucial difficulties in the analysis. Working with a general
metric on a riemannian manifold with boundary (in the $q$ variable)
does not add any serious difficulties.
I nevertheless chose to put the stress on a presentation  in
coordinate systems, in order to make the essential points of the
 analysis more obvious, and to avoid confusing the possibly 
 non familiar reader with the concise but nevertheless subtle
 notations of intrinsic geometry.  Such
 formulations occur when necessary in the end of the article, for
 example while applying the general framework to Bismut's hypoelliptic Laplacian.
\subsection{The problem}
\label{se.prob}
We shall consider the geometric Kramers-Fokker-Planck equation set on
$X=T^{*}Q$ when $\overline{Q}=Q\cup \partial Q$ is a $d$-dimensional 
compact riemannian manifold with boundary or a
compact perturbation of the euclidean half-space
$\overline{\rz^{2d}_{-}}=(-\infty,0]\times\rz^{d-1}\times
\rz^{d}$\,. When $q\in \partial Q$\,, the fiber
$T^{*}_{q}Q$ is the direct sum
 $T^{*}_{q}\partial Q\oplus N^{*}_{q,\partial Q}Q$
where $N^{*}_{q,\partial Q}Q$ is the conormal fiber at $q$\,.\\
In a neighborhood $U$ of $q\in \overline{Q}$\,, position coordinates are denoted by 
$(q^{1},\ldots,q^{d})$\,.  We shall use Einstein's convention of
up and down repeated indices. An element of the fiber $p\in T^{*}_{q}Q$\,,
$q\in U$\,, is written $p=p_{i}dq^{i}$ and $(q^{1},\ldots,q^{d},
p_{1},\ldots, p_{d})$ are symplectic coordinates in $U\times\rz^{d}\sim
T^{*}U\subset T^{*}Q$\,.\\
The metric on $\overline{Q}$\,, i.e. on the tangent fiber bundle
$\pi_{TQ}:TQ\to Q$\,, is denoted by
$g(q)=g^{T}(q)=\left(g_{ij}(q)\right)_{1\leq i,j\leq d}$ or
$g=g_{ij}(q)dq^{i}dq^{j}$\,, and its dual metric on the cotangent
bundle $\pi_{T^{*}Q}:T^{*}Q \to Q$ is
$g^{-1}(q)=(g^{ij}(q))_{1\leq i,j\leq d}$\,.
The cotangent bundle $X=T^{*}Q$ or
$\overline{X}=T^{*}\overline{Q}=X\sqcup \partial X$\,, viewed as a manifold,
 is endowed with the metric $g\oplus
g^{-1}$\,. Actually for every $x\in \overline{X}$\,, the tangent space
$T_{x}X=T_{x}T^{*}Q$ is decomposed into an horizontal and vertical
component $T_{x}T^{*}Q=(T_{x}T^{*}Q)^{H}\oplus (T_{x}T^{*}Q)^{V}$ (see
Section~\ref{se.geoKFP} for details) which specifies the orthogonal decomposition
of the metric $g\oplus g^{-1}$\,.
The corresponding volume on $X=T^{*}Q$ 
coincides with the symplectic volume form
 $\frac{(-1)^{d}}{d!}(dp_{i}\wedge dq^{i})^{\wedge d}$ on
$T^{*}Q$\,,
 and the integration measure is,
locally in $T^{*}U$ with symplectic coordinates, the
$\rz^{2d}$-Lebesgue's measure simply denoted by $dqdp$\,. The $L^{2}$-space on $X$
will be denoted $L^{2}(X,dqdp)$ and the scalar product (extended
as a duality product between distributions and test functions) and the
norm are
$$
\langle u\,,\,v\rangle=\int_{X}\overline{u(q,p)}{v(q,p)}~dqdp\quad,\quad
\|u\|^{2}=\langle u\,,\,v\rangle\,.
$$
We shall consider also $\mathfrak{f}-$valued functions on $\overline{X}$\,, where
$\mathfrak{f}$ is a complex Hilbert-space\footnote{All our Hilbert
  spaces are assumed separable}. The above scalar (duality) product has to be replaced by
$$
\langle u\,,\,v\rangle=\int_{X}\langle
u(q,p)\,,\,v(q,p)\rangle_{\mathfrak{f}}~dqdp\,,
$$
and the corresponding $L^{2}$-space will be denoted
$L^{2}(X,dqdp;\mathfrak{f})=L^{2}(X,dqdp)\otimes
\mathfrak{f}$ while the space of compactly supported regular sections on
$\overline{X}$ (resp. $X$) is denoted by
$\mathcal{C}^{\infty}_{0}(\overline{X};\mathfrak{f})$
(resp. $\mathcal{C}^{\infty}_{0}(X;\mathfrak{f})$)\,.
Hilbert-completed tensor products  of two Hilbert spaces
$\mathfrak{f}_{1}\,, \mathfrak{f}_{2}\,,$ will be
  denoted by $\mathfrak{f}_{1}\otimes\mathfrak{f}_{2}$\,. 
When necessary the algebraic tensor product
  (or other completions) will be specified by a notation like
  $\mathfrak{f}_{1}\otimes^{alg}\mathfrak{f}_{2}$\,.\\
In a vertical fiber $T^{*}_{q}\overline{Q}$ endowed with the scalar product
$g^{-1}(q)$\,, the length $|p|_{q}$\,, the vertical Laplacian
$\Delta_{p}$ and the harmonic oscillator hamiltonian are defined by
\begin{eqnarray*}
&&
|p|_{q}^{2}=g^{ij}(q)p_{i}p_{j}=p^{T}g^{-1}(q)p\,,
\\
&&
\Delta_{p}=\partial_{p_{i}}g_{ij}(q)\partial_{p_{j}}\quad,\quad
\mathcal{O}_{g}=\frac{-\Delta_{p}+|p|^{2}_{q}}{2}\,.
\end{eqnarray*}
With the energy 
$\mathcal{E}(q,p)=\frac{1}{2}|p|_{q}^{2}$ defined on $T^{*}Q$\,, we
associate the Hamiltonian vector field
$$
\mathcal{Y}_{\mathcal{E}}=g^{ij}(q)p_{i}\partial_{q^{j}}-\frac{1}{2}\partial_{q^{k}}g^{ij}(q)p_{i}p_{j}\partial_{p_{k}}\,,
$$
The geometric Kramers-Fokker-Planck operator that we consider is the scalar
(even when $\mathfrak{f}\neq\cz$) differential operator 
$$
P_{\pm,Q,g}=\pm \mathcal{Y}_{\mathcal{E}}+\mathcal{\mathcal{O}}_{g}\,,
$$
well-defined on $\mathcal{C}^{\infty}_{0}(\overline{X};\mathfrak{f})$\,.
Natural functional spaces are associated with the vertical operator
$\mathcal{O}_{g}$\,.
Along a vertical fiber they are the Sobolev spaces associated with the
harmonic oscillator hamiltonian
$\mathcal{O}_{g}(q)=\frac{-\Delta_{p}+|p|_{q}^{2}}{2}$:
\begin{equation}
  \label{eq.Hsq}
\mathcal{H}^{s}(q)=\left\{u\in
  \mathcal{S}'(\rz^{d},dp;\mathfrak{f})\,,\quad
  (\frac{d}{2}+\mathcal{O}_{g}(q))^{s/2}u\in L^{2}(\rz^{d},dp;\mathfrak{f})\right\}\,.
\end{equation}
The distributional space $\mathcal{H}^{s}(q)$ coincides with the same
space given by the euclidean distance on 
$\rz^{d}_{p}$\,, owing to the global ellipticity of
$(\frac{d}{2}+\mathcal{O}_{g(q)})$ (see \cite{Helgl}\cite{HormIII}-Chap~18). 
Only the corresponding norm depends on $g(q)$\,, via a simple linear
change of the variable $p$\,, and
this provides a Hermitian bundle structure on $\pi_{\mathcal{H}^{s}}:\mathcal{H}^{s}
\to \overline{Q}$\,, $\pi_{\mathcal{H}^{s}}^{-1}(q)=\mathcal{H}^{s}(q)$\,.
Global norms on $L^{2}(Q;\mathcal{H}^{1})$ and $L^{2}(Q;\mathcal{H}^{2})$
are given by
\begin{eqnarray*}
  && \|u\|_{L^{2}(Q;\mathcal{H}^{1})}^{2}=\frac{1}{2}\int_{X}
  \|g^{1/2}(q)\partial_{p}u(q,p)\|_{\mathfrak{f}}^{2}
+\|g(q)^{-1/2}pu(q,p)\|_{\mathfrak{f}}^{2}~dqdp
+d\|u\|^{2}\\
&&\|u\|_{L^{2}(Q;\mathcal{H}^{2})}^{2}=\|(\frac{d}{2}+\mathcal{O}_{g})u\|^{2}\,.
\end{eqnarray*}
The standard Sobolev spaces on $\overline{Q}$ are denoted by
$H^{s}(\overline{Q})$\,, $s\in\rz$ and
the corresponding spaces of $\mathcal{H}^{s'}$-sections  are denoted by
$H^{s}(\overline{Q};\mathcal{H}^{s'})$\,,
$s,s'\in\rz$\,.
According to \cite{ChPi}, the spaces $H^{s}(\overline{Q})$ are locally
identified with
$H^{s}(\overline{\rz^{d}_{-}})=H^{s}((-\infty,0]\times \rz^{d-1})$\,,
which is the set of $u\in
\mathcal{D}'(\rz^{d}_{-})=\mathcal{D}'((-\infty,0)\times \rz^{d-1})$
such that $u=\tilde{u}\big|_{\rz^{d}_{-}}$ with $\tilde{u}\in
H^{s}(\rz^{d})$\,. 
The definitions of $\mathcal{C}^{\infty}_{0}(\overline{Q})$ and
$\mathcal{C}^{\infty}_{0}(\overline{X};\mathfrak{f})$ follow the same rule.
Similarly we shall use the notations
 $\mathcal{S}(\overline{\rz^{d}_{-}})$\,,
 $\mathcal{S}(\overline{\rz^{2d}_{-}})$\,.
\\
In order to study closed (and maximal accretive realizations) of
$P_{\pm,Q,g}$ we need to specify boundary conditions. In a neighborhood
$U\subset \overline{Q}$ of $q\in \partial Q$\,, coordinates can be chosen so that the
metric $g$ equals
\begin{equation}
  \label{eq.defgm}
g(q)=
\begin{pmatrix}
  1&0\\
0& m_{ij}(q^{1},q')
\end{pmatrix}\quad,\quad
g_{ij}(q)dq^{i}dq^{j}=(dq^{1})^{2}+m_{ij}(q^{1},q')d{q'}^{i}d{q'}^{j}\,.
\end{equation}
The corresponding coordinates on $T^{*}U$ are $(q^{1},q',p_{1},p')$
and $(q^{'},p_{1})$ is a coordinate system for the conormal bundle
$N^{*}_{\partial Q}\overline{Q}$\,. 
The symplectic volume element is 
$$
dqdp=dq_{1}dp_{1}dq'dp'\,.
$$
On $\partial X=T^{*}Q\big|_{\partial Q}$ we shall use the measure
$|p_{1}|dp_{1}dq'dp'$ and the corresponding $L^{2}$-space will be
denoted $L^{2}(\partial X,|p_{1}|dq'dp;\mathfrak{f})$
with the norm locally defined by
$$
\|\gamma\|_{L^{2}(T^{*}_{Q'}U,|p_{1}|dq'dp;\mathfrak{f})}^{2}=\int_{T^{*}_{Q'}U} \|\gamma(q,p)\|_{\mathfrak{f}}^{2}~|p_{1}|dq'dp\,,
$$
This definition does not depend on the choice of coordinates $(q^{1},\ldots,q^{d})$ for the
decomposition \eqref{eq.defgm} summarized as $g=1\oplus^{\perp}m$\,. 
Polar coordinates in the momentum
variable
can be defined globally as well. 
By writing the momentum $p\in T^{*}_{q}Q$\,, $p=r\omega$ with
$r=|p|_{q}$ and $\omega \in S^{*}_{q}Q$\,, the space $L^{2}(\partial
X, |p_{1}|dq'dp;\mathfrak{f})$ equals
\begin{eqnarray*}
L^{2}(T^{*}_{\partial Q}Q,|p_{1}|dq'dp;\mathfrak{f})&=&L^{2}((0,+\infty),r^{d-1}dr;
L^{2}(S^{*}_{\partial Q}Q, |\omega_{1}|dq'd\omega;\mathfrak{f}))\\
&=&
L^{2}(\partial Q\times (0,+\infty),r^{d-1}drdq';
L^{2}(S^{*}_{q}Q, |\omega_{1}|d\omega;\mathfrak{f}))\,.
\end{eqnarray*}
Once the decomposition $g=1\oplus^{\perp}m$ is assumed, the
 last line defines a space of $L^{2}$-sections of a Hilbert
bundle on $\partial Q\times (0,+\infty)$\,.\\
On those coordinates, the mappings 
\begin{eqnarray*}
&&(q',p_{1})\to (q',-p_{1})\\
\text{(resp.)}&&(q',p_{1},p')\to (q',-p_{1},p')  
\end{eqnarray*}
defines an involution on the conormal bundle $N^{*}_{\partial
  Q}Q$ (resp. on $T^{*}_{\partial Q}Q$)\,. On $\partial X=T^{*}_{\partial Q}Q$
this involution preserves the energy shells $\left\{|p|_{q}^{2}=C\right\}$\,.
The function $\sign(p_{1})$ is well-defined on  $N^{*}_{\partial Q}Q$ 
 and $T^{*}_{\partial Q}Q$\,, $p_{1}>0$ corresponding to the exterior
 conormal orientation\,.
When the space $\mathfrak{f}$ is endowed with a unitary involution $j$\,,
the mapping
$$
v(q',p_{1},p')\to jv(q',-p_{1},p')
$$
defines a
unitary involution on $L^{2}(\partial X,|p_{1}|dq'dp;\mathfrak{f})$\,.
The even and odd part of $\gamma\in
L^{2}(\partial X,|p_{1}|dq'dp;\mathfrak{f})$ are given by
\begin{eqnarray}
 \label{eq.intgev}
&&\gamma_{ev}(q',p_{1},p')=[\Pi_{ev}\gamma] (q',p)=\frac{\gamma(q',p_{1},p')+
  j\gamma(q',-p_{1},p')}{2}\,,\\
\label{eq.intgodd}&&
\gamma_{odd}(q',p_{1},p')=[\Pi_{odd}\gamma](q',p)=
\frac{\gamma(q',p_{1},p')-
  j\gamma(q',-p_{1},p')}{2}\,.
\end{eqnarray}
The operator $\Pi_{ev}$ and $\Pi_{odd}=1-\Pi_{ev}$ are pointwise
operations in $(q',|p|_{q})$ which can be written as operator valued
multiplications by
$\Pi_{odd,ev}(q',|p|_{q})$
 and they are orthogonal
projections in 
$$
L^{2}(\partial X,|p_{1}|dq'dp;\mathfrak{f})=
L^{2}(\partial Q\times (0,+\infty),r^{d-1}drdq';
L^{2}(S^{*}_{q}Q, |\omega_{1}|d\omega;\mathfrak{f}))\,.
$$
We shall also use a bounded accretive operator $A$ on $L^{2}(\partial
X,|p_{1}|dq'dp;\mathfrak{f})$ which is also diagonal in the variable
$(q',|p|_{q})$\,.
We shall assume that there exists $\|A\|>0$ and $c_{A}>0$ such that
\begin{eqnarray}
\label{eq.intArq1}
  && \hspace{-1cm}\|A(q,r)\|_{\mathcal{L}(L^{2}(S^{*}_{q}Q,
    |\omega_{1}|d\omega ;\mathfrak{f}))}\leq
  \|A\|\quad\text{for~a.e.}~(q,r)\in \partial Q\times \rz_{+}\,,\\
&&
\label{eq.intArq0}
\left[A(q,r)\,,\, \Pi_{ev}(q,r)\right]=0\quad \text{for~a.e.}~(q,r)\in \partial Q\times \rz_{+}\,,\\
\label{eq.intArq2}
\text{with~either}&& \min \sigma(\Real A(q,r))\geq
c_{A}>0\quad\text{for~a.e.}~(q,r)\in \partial Q\times \rz_{+}\,,\\ 
\label{eq.intArq3}
\text{or}&& A(q,r)=0 \quad\text{for~a.e.}~(q,r)\in\partial Q\times \rz_{+}\,.
\end{eqnarray}
Our aim is to study the operator $K_{\pm,A,g}$ defined in
$L^{2}(X,dqdp;\mathfrak{f})$ with the domain $D(K_{\pm,A,g})$ 
characterized by
\begin{eqnarray}
\label{eq.intDKA1}
  && u\in L^{2}(Q,dq;\mathcal{H}^{1})\quad,\quad P_{\pm,Q,g}u\in
  L^{2}(X,dqdp; \mathfrak{f})\,,\\
\label{eq.intDKA2}
&& \forall R>0\,, 1_{[0,R]}(|p|_{q})\gamma u \in L^{2}(\partial X,|p_{1}|dq'dp; \mathfrak{f})\,,\\
\label{eq.intDKA3}
&& \gamma_{odd}u=\pm\sign(p_{1})A\gamma_{ev}u\,.
\end{eqnarray}
For the case $A=0$\,, it is interesting to introduce the set
$\mathcal{D}(\overline{X},j)$ of regular functions which satsify
\begin{eqnarray}
  \label{eq.DXj1}
  && u\in
  \mathcal{C}^{\infty}_{0}(\overline{X};\mathfrak{f})\quad,\quad
\gamma_{odd}u=0\,, i.e.~\gamma u(q,-p)=j\gamma u(q,p)\\
\text{and}\label{eq.DXj2}
&& \partial_{q^{1}}u=\mathcal{O}(|q^{1}|^{\infty})\,,
\end{eqnarray}
where the last line is the local writing in a coordinate systems such
that $g=1\oplus^{\perp}m$\,.\\
\textbf{Convention:} We keep the letter $P$\,, often with additional informational
indices, for differential operators acting on $\mathcal{C}^{\infty}$-functions or
distributions. When $P$ denotes a Kramers-Fokker-Planck operator,
 the letter $K$ will be used for closed (and
actually maximal accretive) realizations, also parametrized by $A$
(and $j$)\,,
 of $P$ in $L^{2}(X,dqdp;\mathfrak{f})$\,.
\subsection{Main results}
\label{se.typres}
Although the analysis mixes the two cases $A=0$ and $A\neq 0$\,, we
separate them here for the sake of clarity.
\begin{theorem}
\label{th.main0} 
Let $Q$ be a riemannian compact manifold with boundary or a compact
perturbation of the euclidean half-space
$\overline{\rz^{d}_{-}}=(-\infty,0]\times \rz^{d-1}$\,.
The operator
$K_{\pm,0,g}-\frac{d}{2}=P_{\pm,g}-\frac{d}{2}$ with the domain defined by
\eqref{eq.intDKA1}\eqref{eq.intDKA2}\eqref{eq.intDKA3} with $A=0$
 is maximal accretive.\\
It satisfies the integration by part identity
\begin{eqnarray*}
 \Real\langle u\,,\, (K_{\pm,0,g}+\frac{d}{2})u\rangle
=
\langle u\,,\,(\frac{d}{2}+\mathcal{O}_{g})u\rangle
=\|u\|_{L^{2}(Q;\mathcal{H}^{1})}^{2}
\end{eqnarray*}
for all $u\in D(K_{\pm,0,g})$\,.\\
There exists $C>0$\,,
 such that
\begin{multline*}
\langle \lambda\rangle^{1/4}\|u\|+\langle
\lambda\rangle^{\frac{1}{8}}\|u\|_{L^{2}(Q;\mathcal{H}^{1})}
+\|u\|_{H^{\frac{1}{3}}(\overline{Q};\mathcal{H}^{0})}
\\
+\langle\lambda\rangle^{\frac{1}{4}}
\|(1+|p|_{q})^{-1}\gamma u\|_{L^{2}(\partial X,|p_{1}|dq'dp;\mathfrak{f})}
\leq C\|(K_{\pm,0,g}-i\lambda)u\|
\end{multline*}
for all $u\in D(K_{\pm,0,g})$ and all $\lambda\in\rz$\,.\\
When $\Phi\in \mathcal{C}^{\infty}_{b}([0,+\infty))$\,, is such that
  $\Phi(0)=0$\,,
 there exist a constant $C'>0$
 independent of $\Phi$ and constant $C_{\Phi}>0$ 
such that
$$
\|\Phi(d_{g}(q,\partial Q))\mathcal{O}_{g}u\|\leq 
C'\|\Phi\|_{L^{\infty}}\|(K_{\pm,0,g}-i\lambda)u\|+C_{\Phi}\|u\|\,,
$$
for all $u\in D(K_{\pm,0,g})$ and all $\lambda\in\rz$\,.\\
The adjoint $K_{\pm,0,g}^{*}$ equals $K_{\mp,0,g}$\,.\\
Finally the set $\mathcal{D}(\overline{X},j)$ defined by
\eqref{eq.DXj1}\eqref{eq.DXj2} is dense in $D(K_{\pm,0,g})$ endowed with its graph norm.
\end{theorem}
In the case when $A\neq 0$ and by assuming
\eqref{eq.intArq1}\eqref{eq.intArq0}\eqref{eq.intArq2}\eqref{eq.intArq3},
we have better estimates of the trace, but no density theorem in general.
\begin{theorem}
\label{th.mainA} 
Let $Q$ be a riemannian compact manifold with boundary or a compact
perturbation of the euclidean half-space
$\overline{\rz^{d}_{-}}=(-\infty,0]\times \rz^{d-1}$\,.
Assume that the operator $A\in \mathcal{L}(L^{2}(\partial
X,|p_{1}|dq'dp;\mathfrak{f}))$ satisfies \eqref{eq.intArq1}\eqref{eq.intArq0}\eqref{eq.intArq2}\eqref{eq.intArq3}.\\
The operator
$K_{\pm,A,g}-\frac{d}{2}=P_{\pm,g}-\frac{d}{2}$ with the domain defined by
\eqref{eq.intDKA1}\eqref{eq.intDKA2}\eqref{eq.intDKA3}
 is maximal accretive.\\
It satisfies the integration by part identity
\begin{eqnarray*}
 \Real\langle u\,,\, (K_{\pm,A,g}+\frac{d}{2})u\rangle
=
\|u\|_{L^{2}(Q;\mathcal{H}^{1})}^{2}+\Real\langle
\gamma_{ev}u\,,\,A\gamma_{ev}u\rangle_{L^{2}(\partial X,|p_{1}|dq'dp;\mathfrak{f})}\,,
\end{eqnarray*}
for all $u\in D(K_{\pm,A,g})$\,.\\
There exists $C>0$ and for any $t\in [0,\frac{1}{18})$ there exists $C_{t}>0$
 such that
\begin{multline*}
\langle \lambda\rangle^{1/4}\|u\|+\langle
\lambda\rangle^{\frac{1}{8}}\|u\|_{L^{2}(Q;\mathcal{H}^{1})}
+C_{t}^{-1}\langle \lambda\rangle^{\frac{1}{8}}\|u\|_{H^{t}(\overline{Q};\mathcal{H}^{0})}
\\
+\langle\lambda\rangle^{\frac{1}{8}}
\|\gamma u\|_{L^{2}(\partial X,|p_{1}|dq'dp;\mathfrak{f})}
\leq C\|(K_{\pm,A,g}-i\lambda)u\|
\end{multline*}
for all $u\in D(K_{\pm,A,g})$ and all $\lambda\in\rz$\,.\\
When $\Phi\in \mathcal{C}^{\infty}_{b}([0,+\infty))$\,, is such that
  $\Phi(0)=0$\,,
 there exist a constant $C'>0$
 independent of $\Phi$ and constant $C_{\Phi}>0$ 
such that
$$
\|\Phi(d_{g}(q,\partial Q))\mathcal{O}_{g}u\|\leq 
C'\|\Phi\|_{L^{\infty}}\|(K_{\pm,A,g}-i\lambda)u\|+C_{\Phi}\|u\|\,,
$$
for all $u\in D(K_{\pm,A,g})$ and all $\lambda\in\rz$\,.\\
The adjoint $K_{\pm,A,g}^{*}$ equals $K_{\mp,A^{*},g}$\,.
\end{theorem}

Specific cases, in particular the one dimensional case and the flat
multidimensional case, will be studied with weaker assumptions or
stronger results. Other results are gathered in
Section~\ref{se.variation}.
Various applications are listed in Section~\ref{se.appli}.

\subsection{Guidelines for reading this text}
Although this text is rather long the strategy is really a classical
one for boundary value problems (see e.g. \cite{ChPi}\cite{BdM}\cite{HormIII}-Chap~20):
\begin{enumerate}
\item The first step (see especially \cite{BdM}) consists in a full
  understanding of the simplest one-dimensional problem.
\item In the second step, separation of variable arguments are
  introduced in order to treat straight half-space problems.
\item The last step, is devoted to the local reduction of the general
  problem to the straight half-space problem, by checking that the
  correction terms due to the change of coordinates can be considered
  in some sense as perturbative terms.
\end{enumerate}
Once this is said, one has to face two difficulties:
\begin{description}
\item[a)] The one-dimensional boundary value problem is a two-dimensional
  problem with $(q^{1},p_{1})\in (-\infty,0]\times \rz$\,, with
  $p_{1}$-dependent coefficients. Moreover it rather looks like a corner
  problem because at $q^{1}=0$ the cases $p_{1}>0$ and $p_{1}<0$ are
  discontinuously partitionned. Singularities actually occur at
  $p_{1}=0$ and have to be handled with weighted $L^{2}$-norms.

\vspace{1cm}
\figinit{1.pt}

\figpt 0: (-200,0)
\figpt 1: $q^{1}$(40,0)
\figpt 2: (0,-80)
\figpt 3:$p_{1}$ (0,80)
\figpt 4: ${(0,0)}$ (0,0)
\figpt 5: $X$ (-120,-30)
\figpt 6: $\partial X$ (0,-40)
\figpt 10: (-35,20)
\figpt 11: (-5,20)
\figpt 12: (-65,40)
\figpt 13: (-5,40)
\figpt 14: (-95,60)
\figpt 15: (-5,60)
\figpt 20: (-35,-20)
\figpt 21: (-5,-20)
\figpt 22: (-65,-40)
\figpt 23: (-5,-40)
\figpt 24: (-95,-60)
\figpt 25: (-5,-60)

\psbeginfig{}
\psset(width=2)
\psarrow[0,1]
\psarrow[2,3]
\psset(width=1)
\psarrow[10,11]
\psarrow[12,13]
\psarrow[14,15]
\psarrow[21,20]
\psarrow[23,22]
\psarrow[25,24]
\psendfig
\figvisu{\figBoxA}{\vbox{
    \begin{center}
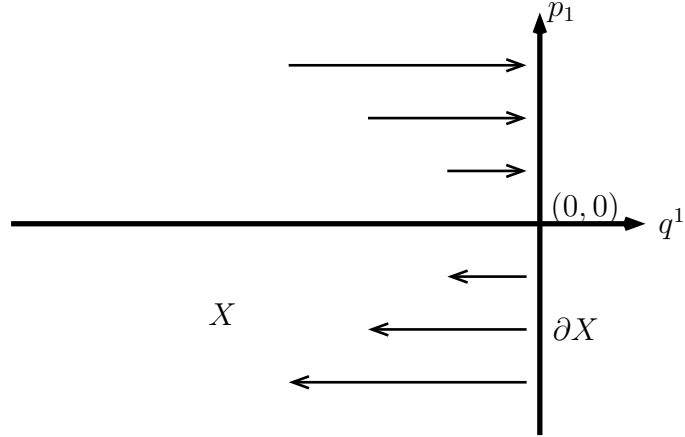
\bf Fig.1: The boundary $\partial X=\left\{q^{1}=0\right\}$ and the
vector field $p_{1}\partial_{q^{1}}$ are represented. For the
absorbing case, the boundary condition says $\gamma u(p_{1})=0$ for
$p_{1}<0$ and corresponds to the case ($j=1$ and $A=1$).
\end{center}
}}{
\figwritee 1:(1)
\figwritene 4:(0.2)
\figwritee  3:(3)
\figwrites 5:(0)
\figwritee 6:(1)
}
\centerline{\box\figBoxA}
\item[b)] In Step~3, the extrinsic curvature of the boundary brings
  a singular perturbation : after a change of
  coordinates  the  corresponding perturbative terms are not negligible
as compared to the regularity and decay estimates obtained in Step~2. 
In the end, the subelliptic estimates are deteriorated by this
curvature effect. Here is the geometric reason: When one
considers the geodesic flow on $\overline{Q}=Q\sqcup\partial Q$\,,
that is the flow of the hamiltonian vector field $\mathcal{Y}$ on
$T^{*}Q$ or $S^{*}Q$ completed by specular reflection at the boundary,
one has to face the well known problem of glancing rays. There are the two
categories of gliding and grazing rays (see fig.2 below) which prevent
from smooth symplectic reductions to half-space problems (no
$\mathcal{C}^{\infty}$ solutions to the eikonal equations for example).
Those problems have been widely studied with the propagation of
singuralities for
 the wave equation (see
e.g. \cite{AnMe}\cite{Tay1}\cite{Tay2}\cite{MeSj1}\cite{MeSj2}). Here
the question is: To what extent the dissipative term
$\mathcal{O}_{g}=\frac{-\Delta_{p}+|p|^{2}_{q}}{2}$ and the
hypoelliptic bracketing (with the hamiltonian vector field
$\mathcal{Y}$) allow to absorb those regularity problems~?

\figinit{1.pt}

\figpt 1: (-42,5)
\figpt 2: (-45,20)
\figpt 3: (-55,40)
\figpt 6: (-75,60)
\figpt 11: (-42,-5)
\figpt 12: (-45,-20)
\figpt 13: (-55,-40)
\figpt 16: (-75,-60)
\figpt 17:(-85,-55)
\figpt 20: $Q$(-100,-10)
\figpt 21: $\partial Q$(-40,0)

\figpttraC101:=1/150,0/
\figpttraC102:=2/150,0/
\figpttraC103:=3/150,0/
\figpttraC106:=6/150,0/
\figpttraC111:=11/150,0/
\figpttraC112:=12/150,0/
\figpttraC113:=13/150,0/
\figpttraC116:=16/150,0/
\figpt 117: (110,-60)
\figpt 118:(110,60)
\figpt 120: $Q$ (120,10)
\figpt 121: $\partial Q$ (110,0)
\psbeginfig{}
\psset(width=2)
\pscurve[16,16,13,12,11,1,2,3,6,6]
\pscurve[116,116,113,112,111,101,102,103,106,106]
\psset(width=1)
\psline[17,13,1,3,6]
\psarrow[117,118]
\psendfig
\figvisu{\figBoxA}{\vbox{
    \begin{center}
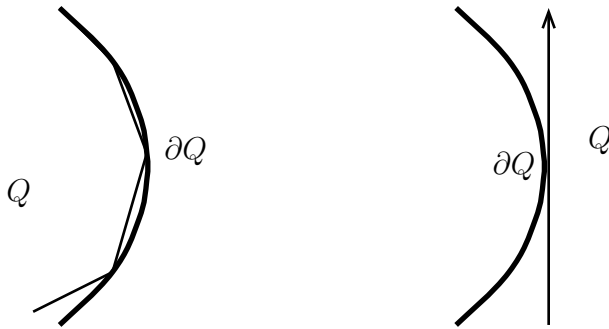
\bf Fig.2: The left picture show a (approximately) gliding ray and the
right one a grazing ray.
\end{center}
}}{
\figwritee 20:(1)
\figwritene 21:(1)
\figwritee 120:(1)
\figwritew 121:(1)
}
\centerline{\box\figBoxA}
\end{description}

A deductive reading is possible starting from Section~\ref{se.model}
to Section~\ref{se.appli}, by referring occasionally to the rather
elementary
 arguments gathered in Appendix~\ref{se.kfpline} and Appendix~\ref{se.partunit}.
According to what the reader is looking for, below are some details
about the various Sections.
\begin{itemize}
\item The reader who wants to know the consequences and applications
  of Theorem~\ref{th.main0} and Theorem~\ref{th.mainA} can go directly
  to Section~\ref{se.variation} and
  Section~\ref{se.appli}. Section~\ref{se.variation} contains
  corollaries about the spectral properties of $K_{\pm, A,g}$ and the
  exponential decay properties of $e^{-tK_{\pm,A,g}}$\,. There are also
  extensions to the case when a potential $V(q)$ is added to the
  kinetic energy $\mathcal{E}(q,p)=\frac{|p|_{q}^{2}}{2}$\,, or to the case when
  $Q\times\mathfrak{f}$ is replaced by some fiber
  bundle. Section~\ref{se.appli} lists various natural boundary
  condition operators $A$ for the scalar case, with interpretations in
  terms of the Langevin stochastic process, and then considers
  specific cases for Bismut's hypoelliptic Laplacian.
\item Appendix~\ref{se.kfpline} recalls the maximal (i.e. with optimal
  exponent) subelliptic estimates in the translation invariant
  case. This is carried out with elementary arguments: Fourier
  transform, harmonic oscillator Hamiltonian $-\Delta_{p}+|p|^{2}$ and
  the one-dimensional complex Airy operator
  $-\partial_{x}^{2}+ix$\,. For the freshman in subelliptic estimates,
  this can be a good starting point.
\item Section~\ref{se.model} provides a thorough study of the
  one-dimensional problem $\overline{Q}=(-\infty,0]$ endowed with the
  euclidean metric $g=(dq^{1})^{2}$\,. The singularity at $p_{1}=0$ is
  solved by introducing a quantization $S$ of the function
  $\sign(p_{1})$ related to some kind of Fourier series in the
  $p_{1}$-variable. The ``Fourier'' basis is a set of eigenvectors for
the compact self-adjoint operator $(\frac{1}{2}+\mathcal{O}_{g})^{-1}p_{1}$
acting in $\mathcal{H}^{1}$\,. The interior problem 
$(p_{1}\partial_{q^{1}}+\frac{1}{2}+\mathcal{O}_{g})u=f$\,, the
condition $\gamma u\in L^{2}(\rz,|p_{1}|dp_{1})$ as
well as the boundary condtions
$\gamma_{odd}u=\sign(p_{1})A\gamma_{odd}u$ are trivialized in terms of
those Fourier series. A Calderon projector is introduced and general
boundary value problems can be studied.
\item Section~\ref{se.cusp} gathers several functional analysis
 properties for semigroup   generators which satisy some subelliptic
 estimates.
 Like in \cite{HerNi} (see also \cite{EcHa}) those generators have
 good resolvent estimates in some cuspidal domain of the complex
 plane,
 and they interpolate between the
 notion of sectorial operators and the one of general maximal
 accretive operators. We shall speak of ``cuspidal semigroup'' or
 ``cuspidal operator'' for the maximal accretive generator. 
Although this does not provide very good subelliptic estimates, this
property is stable by tensorization, which allow separation of
variable arguments for the straight half-space problem. 
Perturbative results are also provided.
\item With Appendix~\ref{se.kfpline}, Section~\ref{se.model}
  and Section~\ref{se.cusp} only, the reader will be convinced that
  subelliptic estimates are reasonable for boundary value problem
 for  which the differential operator and 
the boundary conditions allow the separation of the variables
$(q^{1},p_{1})$
and $(q',p')$\,.
This is reconsidered in Section~\ref{se.insvar} with the specific case
$A=0$ and $A=1$\,. For $A=1$ non homogeneous boundary value problems
are studied with the help of the one-dimensional case and the
variational argument proposed in kinetic theory (see
\cite{Car}\cite{Luc}) and relying on the  weak
formulations of \cite{Lio}. For the case $A=0$\,, the problem is extended to the
whole space with the reflection $(q^{1},p_{1})\to (-q^{1},-p_{1})$\,.
Both cases are useful and their nice properties compensate the fact
that we do not have accurate enough information about the Calderon projector for
Kramers-Fokker-Planck operators in $\overline{\rz^{2d}_{-}}$ when
$d>1$\,.
\item Half-space problems with a general boundary operator $A$ are
  treated in Section~\ref{se.genBC} after a reduction to the boundary
  with the help of the case $A=1$\,. General subelliptic estimates are
  obtained by using the reflection principle $(q^{1},p_{1})\to
  (-q^{1},-p_{1})$ and the nice properties of the case $A=0$\,.
\item The geometric analysis really starts in Section~\ref{se.geoKFP}
  where the maximal subelliptic estimates of Lebeau in
  \cite{Leb2} are written in the way which is used afterwards.
It concerns the case of manifolds without boundary, but the 
difficulty of a non vanishing curvature is recalled after
\cite{Leb1}\cite{Leb2} and explained in view of boundary value problems.
Spatial partitions of unities and dyadic
  partitions of unities in the momentum variables are also introduced
  in this section. Repeatedly used easy formulas for partitions of unities
  are recalled in Appendix~\ref{se.partunit}.
\item The most technical part is in Section~\ref{se.multi}: We study
  the case of cylinders $\overline{Q}=\times (-\infty,0]\times Q'$
  where $Q'$ is compact and the metric $g$ has the canonical form
  \eqref{eq.defgm}. 
Accurate subelliptic estimates for the case  when
$\partial_{q^{1}}g\equiv 0$ (straight cylinder)
are deduced from the general functional framework of
Section~\ref{se.genBC}. With a dyadic partition of unity in the
momentum $p$\,, the large momentum analysis is replaced by some kind
of semiclassical asymptotics on a compact set in $(q,p)$\,.
The more general case $\partial_{q^{1}}g\neq 0$ is sent to the case
$\partial_{q^{1}}g\equiv 0$ with the help of a non symplectic transformation on
$X=T^{*}Q$ which  is the identity along the boundary $\partial X=\left\{(q^{1},q',p_{1},p')\,,
  q^{1}=0\right\}$\,.
 A delicate use of the second resolvent formula for the semiclassical
 problem then
allows to aborb the perturbation which encodes the curvature effect.
 After gluing all the dyadic pieces in
$p$ and then using a spatial partition of unity (the local form
$A=A(q,|p|_{q})$ is used for both steps) Theorem~\ref{th.main0} and
Theorem~\ref{th.mainA} are proved.
\end{itemize}

\section{One dimensional model problem}
\label{se.model}
As we learn from the general theory of boundary value problems  for
linear PDEs, relying on the construction of Calderon's projector,  as
presented in \cite{ChPi}\cite{BdM}\cite{HormIII}-Chap~20,  the key point
is a good understanding of half-line one dimensional model problems
with constant coefficients. Here the one dimensional case is actually
a bidimensional problem with $p$-dependent coefficients.
The main ingredient of this section is the introduction of adapted
``Fourier series'' in the $p$-variable which allow a thourough study
of Calderon's projector and general boundary value problems.
\subsection{Presentation}
\label{se.presentmodel}
We consider the simple case when
$\overline{Q}=\overline{\rz_{-}}=(-\infty,0]$ is endowed with the
euclidean metric. The cotangent bundle is $X=T^{*}Q=\rz^{2}_{-}$\,,
$\overline{X}=T^{*}\overline{Q}=\overline{\rz^{2}_{-}}$\,.
\\
The Kramers-Fokker-Planck operator is simply given by
\begin{eqnarray*}
  && \mathcal{Y}_{\mathcal{E}}=p\partial_{q}\quad,\quad
  \mathcal{O}=\frac{1}{2}(-\partial_{p}^{2}+p ^{2})\,,\\
\text{and}&&
P=p\partial_{q}+\frac{-\partial_{p}^{2}+p^{2}}{2}\quad\text{on}\quad
\rz^{2}_{-}=\left\{(q,p)\in \rz^{2}\,,\, q<0\right\}\,.
\end{eqnarray*}
The space
$\mathcal{H}^{s}(q)=(\frac{1}{2}+\mathcal{O})^{-s/2}L^{2}(\rz^{d};\mathfrak{f})$
and its norm do not depend on $q\in \overline{\rz_{-}}$
 and we simply write $\mathcal{H}^{s}(q)=\mathcal{H}^{s}$\,. The Sobolev space
$H^{s}(\overline{\rz_{-}})$ is the usual one and the notations,
$L^{2}(\rz_{-};\mathcal{H}^{s})$ is better written here $L^{2}(\rz_{-},dq;\mathcal{H}^{s})$ and
$H^{s}(\overline{\rz_{-}};\mathcal{H}^{s})$ keeps the same meaning as
in the introduction.\\
We want to understand the boundary conditions along $\{q=0\}$, which ensure
the maximal accretivity of the associated closed operators.
When $u\in \mathcal{S}(\overline{\rz^{2}_{-}};\mathfrak{f})$ one computes
\begin{equation*}
\Real\langle u\,,\,
Pu\rangle=\frac{1}{2}
\int_{\rz^{2}_{-}}|\partial_{p}u|_{\mathfrak{f}}^{2}+|pu|_{\mathfrak{f}}^{2}~dqdp
+\frac 1 2\int_{\rz}p|u(0,p)|_{\mathfrak{f}}^{2}~dp\,.
\end{equation*}
For $u\in \mathcal{S}(\overline{\rz^{2}_{-}};\mathfrak{f})$ (or possibly less
regular $u$), a
trace $\gamma u (p)=u(0,p)$ is decomposed according to the general
definitions \eqref{eq.intgev}\eqref{eq.intgodd} into
\begin{eqnarray*}
 \gamma_{ev}u(p)=\frac{1}{2}\left[\gamma u(p)+j\gamma u(-p)\right]
&\quad\text{and}\quad&
 \gamma_{odd}u(p)=\frac{1}{2}\left[\gamma u(p)-j\gamma u(-p)\right]\;,\\
\gamma_{+}u(p)=\gamma u(p)1_{(0,+\infty)}(p)
&\quad\text{and}\quad&
\gamma_{-}u(p)=\gamma u(p)1_{(-\infty,0)}(p)\;,
\end{eqnarray*}
where $j$ is the unitary involution on $\mathfrak{f}$\,. Notice the
relations
\begin{eqnarray*}
&&
[\gamma_{ev}u+\sign(p)\gamma_{odd}u](p)=\gamma_{+}u(p)+j(\gamma_{+}u)(-p)\,,\\
\text{and}&&
[\gamma_{ev}u-\sign(p)\gamma_{odd}u](p)=j(\gamma_{-}u)(-p)+\gamma_{-}u(p)\,.
\end{eqnarray*}
With those notations the boundary term can be written
\begin{eqnarray*}
&&
\frac 1 2
\int_{\rz}p|u(0,p)|_{\mathfrak{f}}^{2}~dp=\frac 1 2 \int_{\rz}\left[|\gamma_{+}u(p)|_{\mathfrak{f}}^{2}-|\gamma_{-}u(p)|_{\mathfrak{f}}^{2}\right]~|p|dp\\
&&\qquad=
\frac{1}{4}\int_{\rz}|\gamma_{ev}u+\sign(p)\gamma_{odd}u|_{\mathfrak{f}}^{2}-
|\gamma_{ev}u-\sign(p)\gamma_{odd}u|_{\mathfrak{f}}^{2}~|p|dp
\\
&&\qquad=
\Real\int_{\rz}\langle \gamma_{ev}u(p)\,,\,\sign(p)\gamma_{odd}u(p)\rangle_{\mathfrak{f}}~|p|dp\,
\\
&&\qquad
=\Real\langle \gamma_{ev}u\,,\, \sign(p)\gamma_{odd}u\rangle_{L^{2}(\rz,|p|dp;\mathfrak{f})}\,.
\end{eqnarray*}
A natural space for the traces is
$L^{2}(\rz,|p|dp;\mathfrak{f})$ 
 and when the
boundary conditions are given by a linear relation of the form
$$
\gamma_{odd}u=\sign(p)\times (A\gamma_{ev}u)\,,
$$
the accretivity of the operator $A$ is a necessary condition for
having an accretive realization of $P$\,.  We shall further assume
that $A$ is maximal accretive and commutes with the orthogonal
projection $\Pi_{ev}:L^{2}(\rz,|p|dp;\mathfrak{f})\ni \gamma\mapsto
\Pi_{ev}\gamma=\gamma_{ev}$\,.  Then we consider the operator $K_{A}$
given by
\begin{eqnarray}
\label{eq.defDKA1d}
& D(K_{A})=\left\{u\in L^{2}(\rz_{-};\mathcal{H}^{1})\,,~Pu\in
    L^{2}(\rz^{2}_{-};\mathfrak{f})\,,~
    \gamma_{odd}u=\sign(p)A\gamma_{ev}u\right\}\,,\\
\label{eq.defKA1d}
& \forall u\in
D(K_{A})\,,~K_{A}u=Pu=(p\partial_{q}+\frac{-\partial_{p}^{2}+p^{2}}{2})u\,.
\end{eqnarray}
Note that the conditions occuring in $D(K_{A})$ are the minimal ones
to give a meaning to the previous calculations.\\
Guided by the hyperbolic nature of $p\partial_{q}$\,, the kinetic
theory (see \cite{Bar}\cite{Car}\cite{Luc}) 
more often formulates the boundary conditions in term of 
$\gamma_{+}u$ (outflow) and $\gamma_{-}u$ (inflow). There are two
fundamental examples with $j=\Id_{\mathfrak{f}}$
\begin{itemize}
\item \underline{Specular reflection:} It is usually written
  $\gamma_{-}u=\gamma_{+}u$ and it reads now
$$
\gamma_{odd}u=0\quad,\quad A=0\,.
$$
\item \underline{Absorbing boundary:} It is usually written
  $\gamma_{-}u=0$ and it reads now
$$
\gamma_{odd}u=\sign(p)\gamma_{ev}u\quad,\quad A=\Id.
$$
\end{itemize}
Both example fulfill the maximal accretivity and commutation conditions.
\subsection{Results}
\label{se.res1D}
The general results hold here for possibly unbounded maximal accretive
operators $(A,D(A))$ in $L^{2}(\rz,|p|dp;\mathfrak{f})$\,. The
commutation with $\Pi_{ev}$ is then defined by 
$\Pi_{ev}e^{-tA}=e^{-tA}\Pi_{ev}$ for all $t\geq 0$\,.
\begin{theorem}
\label{th.maxaccKA}
  Assume that $(A,D(A))$  is maximal accretive on
  $L^{2}(\rz,|p|dp;\mathfrak{f})$ and commutes with $\Pi_{ev}$\,. Then the operator
$K_{A}-\frac{1}{2}$ defined by \eqref{eq.defDKA1d}\eqref{eq.defKA1d}
is densely defined and maximal accretive in
$L^{2}(\rz_{-}^{2},dqdp;\mathfrak{f})$\,.\\
For every $u\in D(K_{A})$\,, the traces $\gamma_{ev}u=\Pi_{ev}\gamma u$
and $\gamma_{odd}u=\Pi_{odd}\gamma u$ are well defined with 
\begin{eqnarray*}
&& \|\gamma_{odd}u\|_{L^{2}(\rz,|p|dp;\mathfrak{f})}\leq
C\|u\|_{D(K_{A})}=C\left[\|u\|+\|K_{A}u\|\right]\,,
\\
&&
\|\gamma_{ev}u\|_{L^{2}(\rz,|p|dp;\mathfrak{f})}+
\|A\gamma_{ev}u\|_{L^{2}(\rz,|p|dp;\mathfrak{f})}
=
\|\gamma_{ev}u\|_{D(A)}\leq C\left[\|u\|+\|K_{A}u\|\right]\,.
\end{eqnarray*}
Any $u\in D(K_{A})$ satisfies the integration
by part identity
\begin{equation}
  \label{eq.intbypartTh1}
\Real\langle \gamma_{ev}u\,,\,
A\gamma_{ev}u\rangle_{L^{2}(\rz,|p|dp;\mathfrak{f})}
+ \|u\|_{L^{2}(\rz_{-},dq;\mathcal{H}^{1})}^{2}=\Real\langle u\,,\,
((\frac{1}{2}+K_{A})u\rangle\,.
\end{equation}
\end{theorem}
The above result is completed by the two following propositions.
\begin{proposition}
\label{pr.extH-1H1}
Under the hypothesis of Theorem~\ref{th.maxaccKA}\,, the boundary
value problem
$$
(P-z)u=f\quad,\quad \gamma_{odd}u=\sign(p)A\gamma_{ev}u\,,
$$
admits a unique solution in $L^{2}(\rz_{-},dq;\mathcal{H}^{1})$ when
$f\in L^{2}(\rz_{-},dq;\mathcal{H}^{-1})$ and $\Real z<
\frac{1}{2}$\,. This solution belongs to
$\mathcal{C}^{0}_{b}((-\infty,0];L^{2}(\rz,|p|dp;\mathfrak{f}))$ and
$\gamma_{ev}u\in D(A)$\,. This provides a unique continuous extension of the
resolvent $(K_{A}-z)^{-1}:L^{2}(\rz_{-},dq;\mathcal{H}^{-1})\to
L^{2}(\rz_{-},dq;\mathcal{H}^{1})$ when $\Real z<\frac{1}{2}$\,.
\end{proposition}
\begin{proposition}
\label{pr.adj1} Under the hypothesis of Theorem~\ref{th.maxaccKA}\,,
the adjoint of $K_{A}^{*}$ of $(K_{A},D(K_{A}))$ is given by
\begin{eqnarray*}
&&D(K_{A}^{*})=\left\{u\in L^{2}(\rz_{-},dq;\mathcal{H}^{1})\,,
  \begin{array}[c]{l}
    P_{-}u\in
    L^{2}(\rz^{2}_{-},dqdp;\mathfrak{f})\,,\\
\gamma_{odd}u=-\sign(p)A^{*}\gamma_{ev}u
  \end{array}
    \right\}\,,\\
&&\forall u\in D(K_{A}^{*})\,,\quad
K_{A}^{*}u=P_{-}u=(-p\partial_{q}+\frac{-\partial_{p}^{2}+p^{2}}{2})u\,. 
\end{eqnarray*}
\end{proposition}
\subsection{Fourier series in $\mathcal{H}^{1}$ and
  $L^{2}(\rz,|p|dp)$}
\label{se.Fseries}
Some results of this section  will be stated with $\mathfrak{f}=\cz$
and $j=\Id$\,. In general we shall use the orthogonal decomposition
\begin{eqnarray}
\nonumber
&&\mathfrak{f}=\mathfrak{f}_{ev}\mathop{\oplus}^{\perp}\mathfrak{f}_{odd}\,, \\
\text{with}
&&
\label{eq.foddev}
\mathfrak{f}_{ev}=\ker(j-\Id)\quad,\quad \mathfrak{f}_{odd}=\ker(j+\Id)\,.
\end{eqnarray}
The harmonic oscillator hamiltonian,
$\mathcal{O}=\frac{1}{2}(-\partial_{p}^{2}+p^{2})\otimes\Id_{\mathfrak{f}}$
satisfies
\begin{eqnarray*}
&&  \mathcal{O}-\frac{1}{2}=a^{*}a\otimes\Id_{\mathfrak{f}}\quad,\quad
  \mathcal{O}+\frac{1}{2}=aa^{*}\otimes\Id_{\mathfrak{f}}=(1+a^{*}a)\otimes\Id_{\mathfrak{f}}\,,\\
\text{with}&& a=\frac{1}{\sqrt{2}}(\partial_{p}+p)\quad,\quad a^{*}=\frac{1}{\sqrt{2}}(-\partial_{p}+p)\,.
\end{eqnarray*}
The Hermite functions are normalized as
$$
\varphi_{0}=\pi^{-1/4}e^{-\frac{p^{2}}{2}}\quad,\quad
\varphi_{n}=(n!)^{-1/2}(a^{*})^{n}\varphi_{0}\,, 
$$
and form an orthonormal family of $L^{2}(\rz,dp;\cz)$ of eigenvectors of
$a^{*}a$\,. Note
\begin{eqnarray*}
  &&(\frac{1}{2}+\mathcal{O})^{s}=\sum_{n\in\nz^{*}}n^{s}|\left(\varphi_{n-1}\rangle\langle\varphi_{n-1}|\otimes\Id_{\mathfrak{f}}\right)\,\\
\text{and}&&
\|u\|_{\mathcal{H}^{s}}^{2}=\sum_{n\in\nz^{*}}n^{s}|u_{n}|_{\mathfrak{f}}^{2}\quad,\quad 
u(p)=\sum_{n\in\nz^{*}}u_{n}\varphi_{n-1}(p)\,.
\end{eqnarray*}
While studying the maximal accretivity of $K_{A}$ one has to study the
equation
$$
(\frac{1}{2}+P)u=f\quad,\quad 
\gamma_{odd}u=\sign p \times A\gamma_{ev}u\quad,\quad u\in L^{2}(\rz_{-},dq; \mathcal{H}^{1})\,.
$$
For $f\in L^{2}(\rz_{-},dq;\mathcal{H}^{-1})$ and by setting
$f=(\frac{1}{2}+\mathcal{O})\check{f}$\,,
 it becomes
$$
(\frac{1}{2}+\mathcal{O})^{-1}p\partial_{q}u=\check{f}\quad,\quad
\gamma_{odd}u=\sign p\times A \gamma_{ev}u
$$
with now $\check{f}\in L^{2}(\rz_{-},dq;\mathcal{H}^{1})$\,.
The clue is the spectral analysis of the operator $(\frac{1}{2}+\mathcal{O})^{-1}p$ in
$\mathcal{H}^{1}$\,.

\subsubsection{Spectral resolution of $(\frac{1}{2}+\mathcal{O})^{-1}p$ on $\mathcal{H}^{1}$}
\label{se.spectrvH0}
The first result diagonalizes $(\frac{1}{2}+\mathcal{O})^{-1}p$ when
$\mathfrak{f}=\cz$\,.
\begin{proposition}
\label{pr.vH0}
Assume $\mathfrak{f}=\cz$\,.
The operator $(\frac{1}{2}+\mathcal{O})^{-1}p=(1+a^{*}a)^{-1}p$ is self-adjoint and compact in
$\mathcal{H}^{1}$\,.
Its spectrum equals $\pm (2\nz^{*})^{-1/2}$ and all
 eigenvalues are simple.
For $\nu\in \pm (2\nz^{*})^{-1/2}$\,, a normalized eigenvector
can be chosen as
\begin{equation*}
e_{\nu}(p)=i^{\frac{\sign(\nu)}{2\nu^{2}}}\nu\varphi_{[\frac{1}{2\nu^{2}}-1]}(p-\frac{1}{\nu})\,.
\end{equation*}
With this choice $e_{\nu}(-p)=e_{-\nu}(p)$ so that
$e_{\nu}+e_{-\nu}$ is even and $e_{\nu}-e_{-\nu}$ is odd.
\end{proposition}
\begin{proof}
The operator $(\frac{1}{2}+\mathcal{O})^{1/2}:\mathcal{H}^{1}\to
L^{2}(\rz,dp)$ 
is unitary and $(\frac{1}{2}+\mathcal{O})^{-1}p$ is unitarily
equivalent to $(\frac{1}{2}+\mathcal{O})^{-1/2}p(\frac{1}{2}+\mathcal{O})^{-1/2}$\,. As
an operator lying in $\textrm{OpS}(\langle p,\eta\rangle^{-1},
\frac{dp^{2}}{\langle
  p\rangle^{2}}+\frac{d\eta^{2}}{\langle\eta\rangle^{2}})$  (see
\cite{HormIII}-Chap~18), it is
bounded and compact on $L^{2}(\rz,dp)$ and clearly self-adjoint.
The spectral equation writes
$$
(\frac{1}{2}+\mathcal{O})^{-1}pe=\nu e\quad \text{in}~\mathcal{S}'(\rz)\,,
$$
which is equivalent to ($\nu=0$ cannot happen)
$$
\left[-\partial_{p}^{2}+(p-\frac{1}{\nu})^{2}+1-
\frac{1}{\nu^{2}}\right]e=0\,.
$$
After setting $e(p)=\varphi(p-\frac{1}{\nu})$\,, this gives
$(aa^{*}-\frac{1}{2\nu^{2}})\varphi=0$ which is possible only when
$$
\frac{1}{2\nu^{2}}=n\in 1+ \nz=\nz^{*}\quad\text{and}\quad \varphi\in \cz \varphi_{n-1}\,.
$$
For $\nu=\pm(2n)^{-1/2}$ with $n\in\nz^{*}$, we take 
$$
e_{\nu}(p)=\frac{e^{i\alpha_{\nu}}}
{\|\varphi_{n-1}(.-\frac{1}{\nu})\|_{\mathcal{H}^{1}}}
\varphi_{n-1}(.-\frac{1}{\nu})\,,
$$
with $\alpha_{\nu}\in \rz$\,. The norm
$\|\varphi_{n-1}(.-\frac{1}{\nu})\|_{\mathcal{H}^{1}}$
 is computed by
 \begin{eqnarray*}
   \|\varphi_{n-1}(.-\frac{1}{\nu})\|_{\mathcal{H}^{1}}^{2}
&=&
\langle
     \varphi_{n-1}(.-\frac{1}{\nu})\,,\, (\frac{1}{2}+\mathcal{O})
     \varphi_{n-1}(.-\frac{1}{\nu})
\rangle
\\
&=& \langle
     \varphi_{n-1}\,,\, \frac{(-\partial_{p}^{2}+(p+\frac{1}{\nu})^{2}+1)}{2}
     \varphi_{n-1} \rangle
\\
&=& \langle
     \varphi_{n-1}\,,\, (\frac{1}{2}+\mathcal{O})
     \varphi_{n-1}\rangle
+
\frac{1}{2\nu^{2}}=\frac{1}{\nu^{2}}\,,
 \end{eqnarray*}
where we have used $\langle \varphi_{n-1}\,,\,
p\varphi_{n-1}\rangle=0$\,. Therefore 
$$
\frac{1}{\|\varphi_{n-1}(.-\frac{1}{\nu})\|_{\mathcal{H}^{1}}}=|\nu|\,,
$$
and we can take
$$
e_{\nu}(p)=e^{i\beta_{\nu}}\nu
\varphi_{n-1}(p-\frac{1}{\nu})\,,\quad \beta_{\nu}\in\rz\,.
$$
If one wants the two functions $e_{\nu}+e_{-\nu}$ and $e_{\nu}-e_{-\nu}$ to
be even or odd, the equality $e_{\nu}(-p)=e_{-\nu}(p)$
enforces
$$
\beta_{\nu}-\beta_{-\nu}=n\pi\quad\mod(2\pi)\,,
$$
and $\beta_{\nu}=\sign(\nu)\frac{\pi}{4\nu^{2}} \mod (2\pi)$ works.
\end{proof}

With this spectral resolution one can define the following objects for
general $(\mathfrak{f},j)$\,. 
\begin{definition}
The self-adjoint operator
$(\frac{1}{2}+\mathcal{O})^{-1}p=
[(1+a^{*}a)^{-1}p]\otimes \Id_{\mathfrak{f}}$ in $\mathcal{H}^{1}$ is
denoted by $A_{0}$\,.\\
  For $s\in \rz$\,, the Hilbert space
  $\mathcal{D}_{s}=A_{0}^{s}\mathcal{H}^{1}$\,, equivalently defined by
$$
\mathcal{D}_{s}=\left\{u=\sum_{\nu\in \pm(
    2\nz^{*})^{-1/2}}e_{\nu}\otimes u_{\nu}~, \sum_{\nu\in  \pm(2\nz^{*})^{-1/2}}|\nu|^{-2s}|u_{\nu}|_{\mathfrak{f}}^{2}<+\infty\right\}\,,
$$
is endowed
with the scalar product
$$
\langle u\,,\, u'\rangle_{\mathcal{D}_{s}}=
\sum_{\nu\in  \pm(2\nz^{*})^{-1/2}}|\nu|^{-2s}\langle u_{\nu}\,,\,u'_{\nu}\rangle_{\mathfrak{f}}\,.
$$
The space $\mathcal{D}_{\infty}=\ccap_{s\in\rz}\mathcal{D}_{s}$ is a
Fr{\'e}chet space with the family of norms
$(\|~\|_{\mathcal{D}_{n}})_{n\in\nz}$ dense in $\mathcal{H}^{1}$ and
its dual is $\mathcal{D}_{-\infty}=\ccup_{s\in\rz}\mathcal{D}_{s}$\,.\\
In those spaces, the operator $S$ is defined by
$$
S\left(\sum_{\nu\in\pm(2\nz^{*})^{-1/2}}e_{\nu}\otimes
    u_{\nu}\right)=
\sum_{\nu\in\pm(2\nz^{*})^{-1/2}}\sign(\nu)e_{\nu}\otimes u_{\nu}\,.
$$
For $\nu\in (2\nz^{*})^{-1/2}$\,, the following notations will be used
    \begin{eqnarray*}
&&      V_{\nu}=
   \left(\cz e_{\nu}\mathop{\oplus}^{\perp}\cz e_{-\nu}\right)\otimes
   \mathfrak{f}
=V_{\nu,ev}\mathop{\oplus}^{\perp}V_{\nu,odd}\,,\\
&&
e_{\nu,ev}=\frac{1}{\sqrt{2}}(e_{\nu}+e_{-\nu})\quad,\quad
e_{\nu, odd}= \frac{1}{\sqrt{2}}(e_{\nu}-e_{-\nu})\,,\\
&& V_{\nu,ev}=(\cz e_{\nu,ev})\otimes\mathfrak{f}_{ev}
\mathop{\oplus}^{\perp}
(\cz e_{\nu,odd})\otimes\mathfrak{f}_{odd}\,,\\
&& V_{\nu,odd}=(\cz e_{\nu,odd})\otimes\mathfrak{f}_{ev}
\mathop{\oplus}^{\perp}(\cz e_{\nu,ev})\otimes\mathfrak{f}_{odd}
\,,
    \end{eqnarray*}
where $\mathfrak{f}_{ev,odd}$ are defined by \eqref{eq.foddev}.
\end{definition}
Here is a list of obvious properties
\begin{itemize}
\item $\mathcal{D}_{s}=\mathop{\bigoplus}_{\nu\in (2\nz^{*})^{-1/2}}^{\perp}V_{\nu}$ and 
 $\|u\|_{\mathcal{D}_{s}}^{2}=\sum_{\nu\in
   (2\nz^{*})^{-1/2}}\nu^{-2s}(|u_{\nu}|_{\mathfrak{f}}^{2}
+|u_{-\nu}|_{\mathfrak{f}}^{2})$\,.
\item $S$ is a unitary self-adjoint operator in any $\mathcal{D}_{s}$
  such that $S^{2}=S$\,. The orthogonal projections $\frac{1+S}{2}=1_{\rz_{+}}(S)$
  and $\frac{1-S}{2}=1_{\rz_{-}}(S)$ are given by
  \begin{eqnarray*}
    \frac{1+S}{2}\big(\sum_{\nu\in \pm (2\nz^{*})^{-\frac 1
          2}}e_{\nu}\otimes u_{\nu}\big)=
\sum_{\nu\in +(2\nz^{*})^{-\frac 1 2}}e_{\nu}\otimes u_{\nu}\,,\\
\frac{1-S}{2}\big(\sum_{\nu\in \pm (2\nz^{*})^{-\frac 1
          2}}e_{\nu}\otimes u_{\nu}\big)
=
\sum_{\nu\in -(2\nz^{*})^{-\frac 1 2}}e_{\nu}\otimes u_{\nu}\,.
  \end{eqnarray*}
Note also the relations
\begin{equation}
  \label{eq.Seevodd}
S(e_{\nu,ev}\otimes u_{\nu})=e_{\nu,odd}\otimes u_{\nu}\quad,\quad
S(e_{\nu,odd}\otimes u_{\nu})=e_{\nu,ev}\otimes u_{\nu}\,.  
\end{equation}
\item With $\mathcal{D}_{0}=\mathcal{H}^{1}$\,,
  $\mathcal{D}_{s}=D(|A_{0}|^{-s})$ when $s>0$ while, for $s<0$\,, $\mathcal{D}_{s}$ is the
  completion of $\mathcal{D}_{0}=\mathcal{H}^{1}$ (or
  $\mathcal{S}(\rz)$) for the norm
  $\|u\|_{\mathcal{D}_{s}}=
\||A_{0}|^{-s}u\|_{\mathcal{H}^{1}}$\,. In particular $\mathcal{D}_{-1}$ is
the completion of $\mathcal{S}(\rz)$ for the norm 
$$
\||A_{0}|u\|_{\mathcal{H}^{1}}=
\|A_{0}u\|_{\mathcal{H}^{1}}=\|(\frac{1}{2}+\mathcal{O})pu\|_{\mathcal{H}^{1}}=\|pu\|_{\mathcal{H}^{-1}}\,.
$$
Note that $\mathcal{D}_{-1}$ is not embedded in $\mathcal{D}'(\rz)$\,,
only in $\mathcal{D}'(\rz\setminus\left\{0\right\})$\,, for
it contains functions behaving like $\frac{1}{p}$ around $p=0$\,.
\item When $u$ belongs to $ \mathcal{D}_{s}\subset  \mathcal{D}'(\rz\setminus\left\{0\right\})$
 for some $s\in\rz$\,, the even  (resp. odd) part,
  $\Pi_{ev}u=\frac{1}{2}(u(p)+ju(-p))$
  (resp. $\Pi_{odd}u=\frac{1}{2}(u(p)-ju(-p))$),
is the orthogonal projection of $u$ onto 
$\mathop{\oplus}_{\nu\in (2\nz^{*})^{-1/2}}^{\perp}V_{\nu,ev}$ 
(resp. $\mathop{\oplus}_{\nu\in  (2\nz^{*})^{-1/2}}^{\perp}V_{\nu,odd}$).
More precisely, for 
\begin{eqnarray*}
&&u=\sum_{\nu\in \pm(2\nz^{*})^{-1/2}}e_{\nu}\otimes u_{\nu}
=
\sum_{\nu\in\pm(2\nz^{*})^{-1/2}}e_{\nu}\otimes(\frac{1+j}{2}u_{\nu}+\frac{1-j}{2}u_{\nu})\,,
\\
&&\text{with}\quad \frac{1+j}{2}u_{\nu}\in \mathfrak{f}_{ev}
\quad,\quad\frac{1-j}{2}u_{\nu}\in \mathfrak{f}_{odd}\,,
\end{eqnarray*}
the orthogonal decomposition is 
\begin{multline}
  \label{eq.orthevodd}
u=\underbrace{\sum_{\nu>0}e_{\nu,ev}\otimes(\frac{1+j}{2})\frac{u_{\nu}+u_{-\nu}}{\sqrt{2}}
+
e_{\nu,odd}\otimes(\frac{1-j}{2})\frac{u_{\nu}-u_{-\nu}}{\sqrt{2}}}_{\Pi_{ev}u}
\\
+
\underbrace{\sum_{\nu>0}e_{\nu,odd}\otimes(\frac{1+j}{2})\frac{u_{\nu}-u_{-\nu}}{\sqrt{2}}
+e_{\nu,ev}\otimes(\frac{1-j}{2})\frac{u_{\nu}+u_{-\nu}}{\sqrt{2}}
}_{\Pi_{odd}u}\,.
\end{multline}
This formula and \eqref{eq.Seevodd} imply $S\Pi_{ev}=\Pi_{odd}S$\,.
\item For $u,u'\in \mathcal{S}(\rz;\mathfrak{f})$\,, the following identites hold:
  \begin{multline*}
    \langle u\,,\,\sign(p)u'\rangle_{L^{2}(\rz,|p|dp;\mathfrak{f})}=\langle u\,,\,
    pu'\rangle
=\langle u\,,\, (\frac{1}{2}+\mathcal{O})(\frac{1}{2}+\mathcal{O})^{-1}pu'\rangle
\\=
\langle u\,,\, A_{0}u'\rangle_{\mathcal{H}^{1}}=\langle u\,,\,
Su\rangle_{\mathcal{D}_{-\frac 1 2}}\,.
  \end{multline*}
\end{itemize}

\subsubsection{An interpolation result}
\label{se.interp}

The main result of this section is the
\begin{proposition}
\label{pr.equalspaces}
The interpolated Hilbert space $\mathcal{D}_{-\frac 1
  2}=[\mathcal{H}^{1}\,,\, A_{0}^{-1}\mathcal{H}^{1}]_{\frac 1 2}$ is nothing but
$$
\mathcal{D}_{-\frac 1 2}=L^{2}(\rz, |p|dp;\mathfrak{f})\,.
$$
Moreover there is a bounded invertible positive operator $M$ on
$\mathcal{D}_{-\frac 1 2}$  (which has the same properties on
$L^{2}(\rz,|p|dp)$) such that
\begin{eqnarray*}
  &&\forall u,u'\in L^{2}(\rz,|p|dp;\mathfrak{f})\,,\quad \langle
  u\,,\,u'\rangle_{L^{2}(\rz,|p|dp;\mathfrak{f})}=\langle u\,,\, M
  u'\rangle_{\mathcal{D}_{-\frac 12}}\,,\\
&& \forall u,u'\in \mathcal{D}_{-\frac 1 2}\,,\quad \langle u\,,\,
u'\rangle_{\mathcal{D}_{-\frac 1 2}}=
\langle u\,,\, M^{-1}u'\rangle_{L^{2}(\rz,|p|dp;\mathfrak{f})}\,,\\
&& S= M\circ \sign p=\sign p \circ M^{-1}\,.
\end{eqnarray*}
The projections $\Pi_{ev},\Pi_{odd}$ given by $(\Pi_{ev}u)(p)=\frac{1}{2}[u(p)+ju(-p)]$ and
$(\Pi_{odd})u(p)=\frac{1}{2}\left[u(p)-ju(-p)\right]$ are
orthogonal for both scalar products and commute with $M$\,.
\end{proposition}
\begin{proof}
By functional calculus the Hilbert space $\mathcal{D}_{-\frac 1 2}$ is
the interpolation (complex or real interpolation are equivalent in
this case see \cite{BeLo}) of
$\left[\mathcal{D}_{0},\mathcal{D}_{-1}\right]_{\frac 1 2}$\,. Moreover
$\mathcal{D}_{0}=\mathcal{H}^{1}$ while $\mathcal{D}_{-1}=A_{0}^{-1}\mathcal{H}^{1}$ is the
completion of $\mathcal{S}(\rz;\mathfrak{f})$ for the norm
$\|p.\|_{\mathcal{H}^{-1}}$\,.\\
\noindent\textbf{a)} Consider the multiplication operator $u\mapsto
\frac{p}{\sqrt{1+p^{2}}}u$\,. It belongs to $\textrm{OpS}(1,
\frac{dp^{2}}{\langle
  p\rangle^{2}}+\frac{d\eta^{2}}{\langle\eta\rangle^{2}})$ and 
$(\frac{1}{2}+\mathcal{O})^{\frac 1
  2}(\frac{p}{\sqrt{1+p^{2}}}\times)(\frac{1}{2}+\mathcal{O})^{-\frac 1 2}$ belongs to 
$\textrm{OpS}(1,
\frac{dp^{2}}{\langle
  p\rangle^{2}}+\frac{d\eta^{2}}{\langle\eta\rangle^{2}})$\,. 
We deduce
$$
\|\frac{p}{\sqrt{1+p^{2}}}u\|_{\mathcal{H}^{1}}\leq
C\|u\|_{\mathcal{H}^{1}}\leq C\|u\|_{\mathcal{D}_{0}}\,,
$$
for all $u\in \mathcal{H}^{1}$\,.\\
Similarly $(\frac{1}{2}+\mathcal{O})^{-\frac 1
  2}\frac{1}{\sqrt{1+p^{2}}}
(\frac{1}{2}+\mathcal{O})^{\frac 1 2}\in \textrm{OpS}(1,
\frac{dp^{2}}{\langle
  p\rangle^{2}}+\frac{d\eta^{2}}{\langle\eta\rangle^{2}})$ implies
\begin{multline*}
  \|\frac{p}{\sqrt{1+p^{2}}}u\|_{\mathcal{H}^{-1}}\leq 
\|(\frac{1}{2}+\mathcal{O})^{-\frac 1 2}\left(\frac{1}{\sqrt{1+p^{2}}}
\times\right)(\frac{1}{2}+\mathcal{O})^{\frac 1 2}(\frac{1}{2}+\mathcal{O})^{-\frac 1 2}pu\|\\
\leq C\|(\frac{1}{2}+\mathcal{O})^{-\frac 1 2}pu\|
\leq
  C\|pu\|_{\mathcal{H}^{-1}}\leq C\|u\|_{\mathcal{D}_{-1}}\,,
\end{multline*}
for all $u\in \mathcal{S}(\rz;\mathfrak{f})$ and by density for all $u\in
\mathcal{D}_{-1}$\,.
By interpolation the mapping $\Phi_{2}:u\mapsto \frac{p}{\sqrt{1+p^{2}}}u$ is continuous
from $\mathcal{D}_{-\frac 1
  2}=\left[\mathcal{D}_{0}\,,\,\mathcal{D}_{-1}\right]_{\frac 1 2}$ into 
$\left[\mathcal{H}^{1}\,,\,\mathcal{H}^{-1}\right]_{\frac 1
  2}=L^{2}(\rz, dp;\mathfrak{f})$ with norm less than $C$\,.\\
\noindent\textbf{b)} Consider the multiplication operator $u\mapsto
\frac{1}{\sqrt{1+p^{2}}}u$\,.
With the same arguments as above, which say now that $(\frac{1}{2}+\mathcal{O})^{\frac 1
  2}(\frac{1}{\sqrt{1+p^{2}}}\times)(\frac{1}{2}+\mathcal{O})^{-\frac 1 2}$
and $(\frac{1}{2}+\mathcal{O})^{-\frac 1
  2}(\frac{1}{\sqrt{1+p^{2}}}\times)(\frac{1}{2}+\mathcal{O})^{\frac 1 2}$
belong to $\textrm{OpS}(1,
\frac{dp^{2}}{\langle
  p\rangle^{2}}+\frac{d\eta^{2}}{\langle\eta\rangle^{2}})$\,, we
deduce
\begin{eqnarray*}
  &&\|\frac{1}{\sqrt{1+p^{2}}}u\|_{\mathcal{D}_{0}}=\|\frac{1}{\sqrt{1+p^{2}}}u\|_{\mathcal{H}^{1}}\leq
  C'\|u\|_{\mathcal{H}^{1}}~,~
  \forall u\in \mathcal{H}^{1}\,,\\
\text{and}\hspace{-0.5cm}&&
\|\frac{1}{\sqrt{1+p^{2}}}u\|_{\mathcal{D}_{-1}}=\|(\frac{1}{2}+\mathcal{O})^{-1/2}\frac{1}{\sqrt{1+p^{2}}}u\|\leq
C'\|u\|_{\mathcal{H}^{-1}}~,~\forall u\in \mathcal{H}^{-1}\,.
\end{eqnarray*}
By interpolation the multiplication operator $\Phi_{1}:u\mapsto \frac{1}{\sqrt{1+p^{2}}}u$ is
continuous from 
$\left[\mathcal{H}^{1}\,,\, \mathcal{H}^{-1}\right]_{\frac 1
  2}=L^{2}(\rz,dp;\mathfrak{f})$ into $\mathcal{D}_{-\frac 1 2}$\,.\\
\noindent\textbf{c)} The mapping $\Phi_{1}:L^{2}(\rz,dp;\mathfrak{f})\to
\mathcal{D}_{-\frac 1 2}$ is one to one.  Indeed assume $\Phi_{1}(u)=0$
in $\mathcal{D}_{-\frac 1 2}$\,. Then $\Phi_{2}\circ
\Phi_{1}(u)=\frac{p}{1+p^{2}}u=0$ in $L^{2}(\rz,dp;\mathfrak{f})$\,, which implies
$u=0$ in $L^{2}(\rz,dp;\mathfrak{f})$\,.\\
\noindent\textbf{d)} The mapping $\Phi_{2}:\mathcal{D}_{-\frac 1 2}\to
L^{2}(\rz,dp;\mathfrak{f})$ is one to one. The proof is a little more delicate than
c).
For $\varphi\in \mathcal{S}(\rz;\mathfrak{f})$
the equalities
\begin{equation}
  \label{eq.idenSD12}
\langle \varphi\,,\, Su\rangle_{\mathcal{D}_{-\frac 1 2}}
=\int_{\rz}\langle\sqrt{1+p^{2}}\varphi(p)\,,\,\frac{p}{\sqrt{1+p^{2}}}u(p)\rangle_{\mathfrak{f}}~dp
=\langle \sqrt{1+p^{2}}\varphi\,,\, \Phi_{2}u\rangle
\end{equation}
holds for any  $u\in \mathcal{S}(\rz;\mathfrak{f})$ and by density
 for any  $u\in \mathcal{D}_{-\frac 1 2}$\,.
If $\Phi_{2}(u)=0$ in $L^{2}(\rz,dp;\mathfrak{f})$ for some $u\in
\mathcal{D}_{-\frac 1 2}$\,, the previous identity gives
$$
\forall \varphi\in \mathcal{S}(\rz;\mathfrak{f})\,,\quad \langle \varphi\,,\,
Su\rangle_{\mathcal{D}_{-\frac 1 2}}=0\,,
$$
which implies $Su=0$ and therefore $u=0$\,.\\
\noindent\textbf{e)} The identity \eqref{eq.idenSD12} holds for
$\varphi=\Phi_{1}u_{1}=\frac{1}{\sqrt{1+p^{2}}}u_{1}$ with $u_{1}\in
\mathcal{S}(\rz;\mathfrak{f})$ and $u=u_{2}\in \mathcal{D}_{-\frac 1 2}$\,. We
obtain
$$
\langle \Phi_{1}u_{1}\,,\, Su_{2}\rangle_{\mathcal{D}_{-\frac 1 2}}
=\langle \sqrt{1+p^{2}}\,\Phi_{1}u_{1}\,,\, \Phi_{2}u_{2}\rangle
=\langle u_{1}\,,\, \Phi_{2}u_{2}\rangle
$$
which extends by continuity to any $u_{1}\in L^{2}(\rz,dp;\mathfrak{f})$ and any
$u_{2}\in \mathcal{D}_{-\frac 1 2}$\,. We deduce
\begin{equation}
  \label{eq.adjPhi}
  \Phi_{1}^{*}=\Phi_{2}S\quad\text{and}\quad \Phi_{2}^{*}=S\Phi_{1}\,.
\end{equation}
\noindent\textbf{f)} 
Call $\tilde{\Phi}_{1}$ (resp. $\tilde{\Phi}_{2}$) the bijection from
 $L^{2}(\rz,dp;\mathfrak{f})$ onto
 $L^{2}(\rz,(1+p^{2})dp;\mathfrak{f})$ (resp. from
 $L^{2}(\rz,\frac{p^{2}}{1+p^{2}}dp;\mathfrak{f})$ onto $L^{2}(\rz,dp;\mathfrak{f})$) given by
$\tilde{\Phi}_{1}u=\frac{1}{\sqrt{1+p^{2}}}u$
(resp. $\tilde{\Phi}_{2}u=\frac{|p|}{\sqrt{1+p^{2}}}u$)\,. One has
\begin{eqnarray*}
&&
\Phi_{1}=j_{1}\circ \tilde{\Phi}_{1}\quad\text{and}\quad
\Phi_{2}=\tilde{\Phi}_{2}\circ(\sign(p)\times)\circ j_{2}\,,
\\
\text{with}&& j_{1}u_{1}=u_{1}\;\text{for}\; u_{1}\in
L^{2}(\rz,(1+p^{2})dp;\mathfrak{f})\quad\text{and}\quad
j_{2}u_{2}=u_{2}\;\text{for}\;u_{2}\in \mathcal{D}_{-\frac 1 2}\,. 
\end{eqnarray*}
Moreover the equality 
$$
\langle j_{1}\tilde{\Phi}_{1}u_{1}\,,\,
Su_{2}\rangle_{\mathcal{D}_{-\frac 1 2}}
=\langle u_{1}\,,\, \tilde{\Phi}_{2}\sign(p)j_{2}u_{2}\rangle
$$
valid for any $u_{1},u_{2}\in \mathcal{S}(\rz;\mathfrak{f})$  can be written, with
$u_{1}'=\tilde{\Phi}_{1}(u_{1})=\frac{u_{1}}{\sqrt{1+p^{2}}}$,
\begin{eqnarray*}
\langle Sj_{1}u_{1}'\,,\, u_{2}\rangle_{\mathcal{D}_{-\frac 1 2}}
&=&
\langle u_{1}\,,\,
\tilde{\Phi}_{2}\sign(p)j_{2}u_{2}\rangle
\\
&=&\int_{\rz}\langle\sqrt{1+p^{2}}\,u_{1}'(p)\,,\,
\frac{|p|}{\sqrt{1+p^{2}}}
\sign(p)(j_{2}u_{2})(p)\rangle_{\mathfrak{f}}~dp
\\
&=& \int_{\rz}
\langle u'_{1}(p)\,,\,\sign(p)(j_{2}u_{2})(p)\rangle_{\mathfrak{f}}~|p|dp\,.
\end{eqnarray*}
where $\int_{\rz}\langle u(p)\,,\,v(p)\rangle_{\mathfrak{f}}~|p|dp$ can be interpreted as a
duality product between
$L^{2}(\rz,(1+p^{2})dp;\mathfrak{f})$ and $L^{2}(\rz, \frac{p^{2}}{1+p^{2}}dp;\mathfrak{f})$\,.
After identifying $j_{k}u_{k}=u_{k}$ (and dropping the $'$)\,, it is better written
\begin{multline}
\label{eq.Suuusignu}
\langle Su_{1}\,,\, u_{2}\rangle_{\mathcal{D}_{-\frac 1
    2}}=\int_{\rz}\langle u_{1}\,,\,(\sign(p)u_{2}\rangle_{\mathfrak{f}}~|p|dp
\\
=\langle u_{1}\,,\,
\sign(p)u_{2}\rangle_{L^{2}(\rz,(1+|p|^{2})dp;\mathfrak{f})\,,\,
  L^{2}(\rz,\frac{p^{2}}{1+p^{2}}dp,\mathfrak{f})}\,,
\end{multline}
for all $u_{1}\in L^{2}(\rz,(1+p^{2})dp;\mathfrak{f})\subset \mathcal{D}_{-\frac 1
  2}$ and all $u_{2}\in \mathcal{D}_{-\frac 1 2}$\,. 
This implies that
$$
L^{2}(\rz,(1+p^{2})dp;\mathfrak{f})\stackrel{Sj_{1}}{\longrightarrow}\mathcal{D}_{-\frac 1 2}
\stackrel{\sign(p)j_{2}=(Sj_{1})^{*}}{\longrightarrow}L^{2}(\rz,\frac{p^{2}}{1+p^{2}}dp;\mathfrak{f})
$$
is a Hilbert triple. With the closed quadratic form $q(u)=\|u\|_{L^{2}(\rz,
  (1+p^{2})dp;\mathfrak{f})}^{2}+\|u\|_{\mathcal{D}_{-\frac 1 2}}^{2}$\,, we can
associate a non negative self-adjoint operator $B_{0}\geq 1$ such that 
$L^{2}(\rz,(1+p^{2})dp;\mathfrak{f})= B_{0}^{-1/2}\mathcal{D}_{-\frac 1 2}$ and
$L^{2}(\rz,\frac{p^{2}}{1+p^{2}}dp;\mathfrak{f})=B_{0}^{1/2}\mathcal{D}_{-\frac 1
  2}$\,. By the functional calculus for $B_{0}$ and complex
interpolation,
 we deduce
$$
\mathcal{D}_{-\frac 1 2}=\left[L^{2}(\rz,(1+p^{2})dp;\mathfrak{f})\,,\, L^{2}(\rz,
  \frac{p^{2}}{1+p^{2}}dp;\mathfrak{f})\right]_{\frac 1 2}= L^{2}(\rz,|p|dp;\mathfrak{f})\,,
$$
with equivalent Hilbert norms.\\
\noindent\textbf{g)} Let us specify the equivalence of the Hilbert
scalar products. The Hilbert space $\mathcal{D}_{-\frac 1
  2}=L^{2}(\rz,|p|dp;\mathfrak{f})$ is now endowed with two scalar products
$$
\langle u_{1}\,,\, u_{2}\rangle_{1}=\langle u_{1}\,,\,
u_{2}\rangle_{\mathcal{D}_{-\frac 1 2}}
\quad\text{and}\quad
\langle u_{1}\,,\, u_{2}\rangle_{2}=\int_{\rz}\langle u_{1}(p)\,,\,
u_{2}(p)\rangle_{\mathfrak{f}}~|p|dp\,.
$$
The identity \eqref{eq.Suuusignu} implies now
$$
\forall u_{1},u_{2}\in \mathcal{D}_{-\frac 1 2}\,,\quad
\langle Su_{1}\,, u_{2}\rangle_{1}=\langle u_{1}\,,\, \sign(p)u_{2}\rangle_{2}\,,
$$
where $S$ is unitary for $\langle~\,,\,~\rangle_{1}$ and $\sign(p)$ is
unitary
for $\langle~\,,\,~\rangle_{2}$\,. We deduce that $M=\sign(p)^{*1}S$ and
$M'=S^{*2}\sign(p)$ are
positive bounded operators (for the respective scalar products) such
that
$$
\langle u_{1}\,,M u_{2}\rangle_{1}=\langle u_{1}\,,\,u_{2}\rangle_{2}
\quad\text{and}\quad 
\langle u_{1}\,,\, u_{2}\rangle_{1}=\langle u_{1}\,,\, M'u_{2}\rangle_{2}\,.
$$
Hence $M'=M^{-1}$ and $M$ is a bounded invertible positive
operator for both scalar products.\\
The writing
$$
\langle u_{1}\,,\, Su_{2}\rangle_{1}=\langle Su_{1}\,,\,u_{2}\rangle=\langle u_{1}\,,\,
\sign(p)u_{2}\rangle_{2}=\langle u_{1}\,,\, M\sign(p)u_{2}\rangle_{1}\,,
$$
implies $S=M\sign(p)$ which can be combined with the relations
$S=S^{-1}=S^{*1}$ and $\sign(p)=\sign(p)^{-1}=\sign(p)^{*2}$\,.\\
\noindent\textbf{h)} It remains to prove the orthogonality of
$\Pi_{ev}(u)$ and $\Pi_{odd}(u)$ for both scalar products and the
commutation with $M$\,.
For the second one it results from $\int_{\rz}f(p)|p|dp=0$ for every
odd element of $L^{1}(\rz,|p|dp)$\,. For the first scalar product, we
have already checked in \eqref{eq.orthevodd} that $\Pi_{ev}$ and
$\Pi_{odd}$ are orthogonal projections in $\mathcal{D}^{s}$\,.
The relation
$\sign(p)\circ\Pi_{ev}=\Pi_{odd}\circ\sign(p)$ comes from the
definition \eqref{eq.intgev}\eqref{eq.intgodd} while
\eqref{eq.orthevodd}
led to $S\circ\Pi_{ev}=\Pi_{odd}\circ S$\,. We obtain
$$
M\circ \Pi_{ev}=
S\circ\sign(p)\circ\Pi_{ev}=S\circ\Pi_{odd}\circ\sign(p)
=\Pi_{ev}\circ S\circ\sign(p)=\Pi_{ev}\circ M\,.
$$
\end{proof}
We end this section with an additional lemma, related with the
previous result: It ensures that some maximal accretivity is preserved
when the scalar product is changed.
\begin{lemma}
\label{le.equaccr}
  Consider a complex Hilbert space $\mathcal{D}$ endowed with two equivalent scalar
  products $\langle~,~\rangle_{1}$ and $\langle~,~\rangle_{2}$, like 
in Proposition~\ref{pr.equalspaces}, and let $M$ be a positive bounded
invertible operator (in $(\mathcal{D},\langle~,~\rangle_{1})$ or  $(\mathcal{D},\langle~,~\rangle_{2})$) such that
$$
\forall u,u'\in \mathcal{D}\,,\quad \langle u\,,\,
Mu'\rangle_{1}=\langle u\,,\, u'\rangle_{2}\,.
$$
Then a densely defined operator $(A, D(A))$ is maximal accretive in
$(\mathcal{D}, \langle~,~\rangle_{2})$ iff $(MA, D(A))$ is
maximal accretive in $(\mathcal{D}, \langle~,~\rangle_{1})$\,.
\end{lemma}
\begin{proof}
Assume $(A,D(A))$ maximal accretive in $(\mathcal{D}, \langle~,~\rangle_{2})$\,.\\
 \noindent\textbf{a)} Then $(MA, D(A))$ is clearly accretive owing to
$$
\forall u\in D(A)\,, \quad\Real\langle u\,,\, MA u \rangle_{1}=\Real
\langle u\,,\, Au\rangle_{2}\geq 0\,.
$$
\noindent\textbf{b)} Let us check that $\Id +MA$ is invertible. Our
assumptions ensure that $\varepsilon\Id+ A$ is invertible for
$\varepsilon>0$ in
$\mathcal{D}$ and therefore it is the case also for
$\varepsilon M+MA$\,. Hence $\varepsilon M+MA$ which is accretive in
$(\mathcal{D}, \langle~,~\rangle_{1})$ is maximal
accretive. Therefore $B_{\varepsilon}=\Id+\varepsilon M+MA$ with
$D(B_{\varepsilon})=D(A)$ is invertible for $\varepsilon>0$ with the
estimate
$$
\forall u\in D(A)\,,\quad
\|u\|_{1}\|(\Id+\varepsilon M+MA)\|_{1}\geq \Real \langle u\,,\,
(\Id+\varepsilon M+MA)u\rangle_{1}\geq \|u\|_{1}^{2}\,,
$$
which means 
$$
\|B_{\varepsilon}^{-1}\|_{\mathcal{L}(\mathcal{D},\|~\|_{1})}\leq 1\,.
$$
The second resolvent formula
$$
(\Id+MA)^{-1}=B_{\varepsilon}^{-1}\left[\Id-\varepsilon MB_{\varepsilon}^{-1}\right]^{-1}\,,
$$
and choosing $\varepsilon< \|M\|_{\mathcal{L}(\mathcal{D},\|~\|_{1})}^{-1}$ imply
that $\Id+MA$ is invertible with $\Ran(I+MA)^{-1}=D(B_{\varepsilon})=D(A)$\,.\\

The two points $a)$ and $b)$ prove the maximal accretivity of $MA$ in
$(\mathcal{D}, \langle~,~\rangle_{1})$ when $(A,D(A))$ is
maximal accretive in $(\mathcal{D},
\langle~,~\rangle_{2})$\,. The converse statement comes from 
$A=M^{-1}MA$\,.
\end{proof}

\subsection{System of ODE and boundary value problem}
\label{systODEBVP}
The spectral resolution of $(\frac{1}{2} + \mathcal{O})^{-1}p$ and the
interpolation result of Proposition~\ref{pr.equalspaces} reduce
boundary value problems for $P$ on $\overline{\rz^{2}_{-}}$ to an infinite
system of ODE's in $\overline{\rz_{-}}$\,.
\subsubsection{Rewriting $Pu=f$}
\label{se.rewriting}
Consider the two cases $I=\rz_{-}$ and $I=\rz$\,.
When $u\in L^{2}(I,dq;\mathcal{H}^{1})$ and $f\in
L^{2}(I,dq;\mathcal{H}^{-1})$ the relation $Pu=f$ can be written
\begin{eqnarray}
  \label{eq.maxaccr1}
&&(\frac{1}{2}+\mathcal{O})^{-1}p\partial_{q}u+u=\check{f}\,,\\
\nonumber\text{with}&&
\check{f}=(\frac{1}{2}+\mathcal{O})^{-1}f\in
L^{2}(I,dq;\mathcal{H}^{1})
=L^{2}(I,dq;\mathcal{D}_{0})\,.
\end{eqnarray}
Meanwhile the boundary condition, in the case $I=\rz_{-}$\,, is written for traces in
$L^{2}(\rz,|p|dp;\mathfrak{f})=\mathcal{D}_{-\frac 1 2}$\,. 
Hence we can use the basis $\left(e_{\nu}\right)_{\nu\in
  \pm(2\nz^{*})^{-\frac 1 2}}$ and use
$\mathcal{D}_{s}=\mathop{\oplus}_{\nu\in \pm(2\nz)^{-\frac 1 2}}^{\perp}(\cz
  e_{\nu})\otimes \mathfrak{f}$ for all $s\in\rz$\,. By writing
$$
u(q,p)=\sum_{\nu\in \pm(2\nz^{*})^{-\frac 1
    2}}u_{\nu}(q)e_{\nu}\,,\quad u_{\nu}(q)\in \mathfrak{f}~\text{a.e}
$$
the squared $L^{2}(\rz_{-},dq;\mathcal{D}_{s})$-norm is nothing but
$$
\|u\|_{L^{2}(\rz_{-},dq;\mathcal{D}_{s})}^{2}=\sum_{\nu\in \pm(2\nz^{*})^{-\frac
  1 2}}|\nu|^{2s}\int_{-\infty}^{0}|u_{\nu}(q)|_{\mathfrak{f}}^{2}~dq
$$
while the trace, when defined in $\mathcal{D}_{-\frac 1 2}$ equals
$u(0,p)=\sum_{\nu\in \pm (2\nz^{*})^{-\frac 1
    2}}u_{\nu}(0)e_{\nu}$ with the squared-norm
$$
\|u(0,.)\|_{\mathcal{D}_{-\frac 1 2}}^{2}
=\sum_{\nu\in \pm (2\nz^{*})^{-\frac 1 2 }}|\nu||u_{\nu}(0)|_{\mathfrak{f}}^{2}\,.
$$
Using the same decomposition for $\check{f}\in
L^{2}(\rz_{-},dq;\mathcal{D}_{0})$\,, the equation \eqref{eq.maxaccr1}
becomes
$$
\forall \nu\in \pm (2\nz^{*})^{-\frac 1 2}\,,\quad
\nu\partial_{q}u_{\nu}+u_{\nu}=\check{f}_{\nu}\quad\text{in}~I\,.
$$
When $I=\rz$\,, the only possible solution is given by
\begin{eqnarray*}
 \nu>0\,,&& u_{\nu}(q)= \frac{1}{\nu}\int_{-\infty}^{q}e^{-\frac{q-s}{\nu}}\check{f}_{\nu}(s)~ds
=(\frac{e^{-\frac{.}{\nu}}}{\nu}1_{\rz_{+}})*\check{f}_{\nu}(q)\\
\nu<0\,,&&
u_{\nu}(q)=\frac{1}{\nu}\int_{0}^{q}e^{-\frac{(q-s)}{\nu}}\check{f}_{\nu}(s)~ds
=
\left(\frac{1}{\nu}e^{-\frac{.}{\nu}}1_{\rz_{-}}\right)*\check{f}_{\nu}(q)\,.
\end{eqnarray*}
Actually this formula provides a solution. Owing to the next Lemma.
\begin{lemma}
\label{le.estL1}
For any $s\in\rz$ the operator $E$ defined by 
$$
E(\sum_{\nu\in \pm(\nz)^{-\frac 1 2}}\check{f}_{\nu}(q)e_{\nu})=
\sum_{\nu\in \pm(\nz)^{-\frac 1 2}}
\left[\left(\frac{e^{-\frac{.}{\nu}}}{\nu}1_{\sign(\nu)\rz_{+}}\right)*\check{f}_{\nu}\right](q)e_{\nu}\,,
$$
is a contraction on $L^{2}(\rz,dq;\mathcal{D}_{s})$\,.
\end{lemma}
\begin{proof}
  It suffices to notice that
$$
\left\|\left(\frac{e^{-\frac{.}{\nu}}}{\nu}1_{\sign(\nu)\rz_{+}}\right)\right\|_{L^{1}(\rz,dq)}=1
$$
\end{proof}
When $I=\rz_{-}$\,, we distinguish the two cases $\nu>0$ and $\nu<0$~:
\begin{description}
\item[$\nu>0$:] The solution  $u_{\nu}\in
  L^{2}(\rz_{-},dq;\mathfrak{f})$ to
  $(\nu\partial_{q}+1)u_{\nu}=\check{f}_{\nu}$ is still
  \begin{equation}
\label{eq.solnu+}
u_{\nu}(q)=\frac{1}{\nu}\int_{-\infty}^{q}e^{-\frac{q-s}{\nu}}\check{f}_{\nu}(s)~ds
=(\frac{e^{-\frac{.}{\nu}}}{\nu}1_{\rz_{+}})*(\check{f}_{\nu}1_{\rz_{-}})(q)\,,
\end{equation}
which contains
\begin{equation}
\label{eq.sol0nu+}
u_{\nu}(0)=\int_{-\infty}^{0}\frac{e^{\frac{s}{\nu}}}{\nu}\check{f}_{\nu}(s)~ds\,.
\end{equation}
\item[$\nu<0$:] The possible solutions $u_{\nu}\in
  L^{2}(\rz_{-},dq;\mathfrak{f})$ to
  $(\nu\partial_{q}+1)u_{\nu}=\check{f}_{\nu}$ are given by
  \begin{eqnarray}
\nonumber
u_{\nu}(q)&=&
u_{\nu}(0)e^{-\frac{q}{\nu}}+\frac{1}{\nu}\int_{0}^{q}e^{-\frac{(q-s)}{\nu}}\check{f}_{\nu}(s)~ds\\
\label{eq.solnu-}
&=&u_{\nu}(0)e^{-\frac{q}{\nu}}+
\left(\frac{1}{\nu}e^{-\frac{.}{\nu}}1_{\rz_{-}}\right)*(\check{f}_{\nu}1_{\rz_{-}})(q)\,.
\end{eqnarray}
\end{description}

\subsubsection{Trace theorem and integration by parts}
\label{se.intbypart}

We next prove a trace theorem and an integration by part formula adapted to the Fourier
series in $\mathcal{D}_{s}$\,.
\begin{definition}
  \label{de.EI}
For an open interval $I$ of $\rz$\,,
$I=(a,b)$ with
$-\infty\leq a <b\leq +\infty$\,, one sets 
$$
\mathcal{E}_{I}(\mathfrak{f})=\left\{u\in
  L^{2}(I,dq;\mathcal{H}^{1})\,,~p\partial_{q}u \in
  L^{2}(I,dq;\mathcal{H}^{-1})\right\}\,,
$$
and its norm is given by
$$
\|u\|_{\mathcal{E}_{I}(\mathfrak{f})}^{2}=\|u\|_{L^{2}(I,dq;\mathcal{H}^{1})}^{2}
+\|p\partial_{q}u\|_{L^{2}(I,dq;\mathcal{H}^{-1})}^{2}\,.
$$
\end{definition}
Remember that  $\mathcal{H}^{s}$ is a space of  $\mathfrak{f}$-valued distributions\,.
\begin{proposition}
\label{pr.intbypart} 
The following properties hold:
\begin{itemize}
\item $\mathcal{E}_{I}(\mathfrak{f})$ is continuously embedded in
  $\mathcal{C}^{0}_{b}(\overline{I};L^{2}(\rz,|p|dp;\mathfrak{f}))=\mathcal{C}^{0}_{b}(\overline{I};\mathcal{D}_{-\frac
  1 2})$\,;
\item $\mathcal{S}(\overline{I}\times \rz;\mathfrak{f})$ is dense in
  $\mathcal{E}_{I}(\mathfrak{f})$ and if 
  $a=-\infty$ (resp. $b=+\infty$)  the norm $\|u(q,.)\|_{L^{2}(\rz,
    |p|dp;\mathfrak{f})}$ goes to $0$ as $q\to -\infty$ (resp. $q\to +\infty$)\,.
\item  Any $u\in \mathcal{E}_{I}(\mathfrak{f})$ fulfills the integration by parts formula
  \begin{equation}
    \label{eq.IPP}
2\Real \langle u\,,\,
p\partial_{q}u\rangle=1_{\rz}(b)\int_{\rz}|u(b,p)|_{\mathfrak{f}}^{2}pdp -
1_{\rz}(a)\int_{\rz}|u(a,p)|_{\mathfrak{f}}^{2}pdp\,.
\end{equation}
\item When $a>-\infty$ (resp. $b<+\infty$) the trace map
  $\gamma_{\bullet}:\mathcal{E}_{I}(\mathfrak{f})\to
  \gamma_{\bullet}u=u(\bullet,p)\in L^{2}(\rz,
  |p|dp;\mathfrak{f})=\mathcal{D}_{-\frac 1 2}$  (with
  $\bullet=a$ or $b$ respectively) is surjective, with a continuous left-inverse.
\end{itemize}
\end{proposition}
\begin{proof}
  After noting that $p\partial_{q}u\in L^{2}(I,dq;\mathcal{H}^{-1})$ is
  equivalent to $(\frac{1}{2}+ \mathcal{O})^{-1}p\partial_{q}u\in
  L^{2}(I,dq;\mathcal{H}^{1})$\,, introduce the orthonormal basis
  $\left(e_{\nu}\right)_{\nu\in \pm (2\nz^{*})^{-\frac 1 2}}$\,. The space $\mathcal{E}_{I}(\mathfrak{f})$ is the set of series 
$u=\sum_{\nu\in \pm (2\nz^{*})^{-\frac 1 2}}u_{\nu}(q)e_{\nu}$\,,
$u_{\nu}(q)\in \mathfrak{f}$ a.e.,
with the norm
$$
\|u\|_{\mathcal{E}_{I}(\mathfrak{f})}^{2}=\sum_{\nu\in\pm
  (2\nz^{*})^{-\frac 1 2}}
\int_{a}^{b}|u_{\nu}(q)|_{\mathfrak{f}}^{2}+|\nu|^{2}|\partial_{q}u_{\nu}(q)|_{\mathfrak{f}}^{2}~dq\,.
$$
Every $u_{\nu}$ belongs to $\mathcal{C}^{0}_{b}(I;\mathfrak{f})$ and goes to $0$ at
$\infty$ when $a=-\infty$ or $b=+\infty$ with
\begin{eqnarray*}
  |\nu||\sup_{q\in I}u_{\nu}(q)|_{\mathfrak{f}}^{2}
&\leq&
C\|u_{\nu}\|_{L^{2}(I,dq;\mathfrak{f})}\left[\|u_{\nu}\|_{L^{2}(I,dq;\mathfrak{f})}+
\|\nu\partial_{q}u_{\nu}\|_{
L^{2}(I,dq;\mathfrak{f})}\right]
\\
&\leq& C'
\left[
\|u_{\nu}\|_{L^{2}(I,dq;\mathfrak{f})}^{2}+|\nu|^{2}\|\partial_{q}u_{\nu}\|_{L^{2}(I,dq;\mathfrak{f})}^{2}
\right]
\leq C'
\|u\|_{\mathcal{E}_{I}(\mathfrak{f})}^{2}\,.
\end{eqnarray*}
The dominated convergence theorem applied
to the series 
$$
\sum_{\nu\in\pm (2\nz^{*})^{-\frac 1
    2}}|\nu|\left|u_{\nu}(q)-u_{\nu}(q_{0})\right|_{\mathfrak{f}}^{2}
\quad
\text{as} \quad
q\to q_{0}\,,
$$ 
with the upper bound
$$
|\nu||u_{\nu}(q)-u_{\nu}(q_{0})|_{\mathfrak{f}}^{2}\leq 2C'
\left[
\|u_{\nu}\|_{L^{2}(I,dq;\mathfrak{f})}^{2}+|\nu|^{2}\|\partial_{q}u_{\nu}\|_{L^{2}(I,dq;\mathfrak{f})}^{2}
\right]\,,
$$
provides the $\mathcal{D}_{-\frac 1 2}$-continuity w.r.t $q$\,.\\
For the density of $\mathcal{S}(\overline{I}\times \rz;\mathfrak{f})$\,, it suffices to
approximate $\displaystyle u=\sum_{\nu \in \pm (2\nz^{*})^{-\frac 1 2}} 
u_{\nu}(q)e_{\nu}$ by a finite sum $u^{N}=\sum_{|\nu|\geq  N^{-1}}
u_{\nu}^{N}(q)e_{\nu}$, with $u_{\nu}^{N}\in
\mathcal{S}(\overline{I};\mathfrak{f})$\,,
\begin{align*}
  &
\|u_{\nu}^{N}\|_{L^{2}(I,dq;\mathfrak{f})}^{2}
+|\nu|^{2}\|\partial_{q}u_{\nu}^{N}\|_{L^{2}(I,dq;\mathfrak{f})}^{2}
\leq 2
\|u_{\nu}\|_{L^{2}(I,dq;\mathfrak{f})}^{2}
+2|\nu|^{2}\|\partial_{q}u_{\nu}\|_{L^{2}(I,dq;\mathfrak{f})}^{2}
\,,\\
\text{and}&
\lim_{N\to\infty}\|u_{\nu}^{N}-u_{\nu}\|_{L^{2}(I,dq;\mathfrak{f})}^{2}
+
|\nu|^{2}\|\partial_{q}(u_{\nu}^{N}-u_{\nu})\|_{L^{2}(I,dq;\mathfrak{f})}^{2}
=0
\end{align*}
for all $\nu\in \pm(2\nz^{*})^{-\frac 1 2}$\,.\\
The integration by parts formula \eqref{eq.IPP} is true when $u\in
\mathcal{S}(\overline{I}\times\rz;\mathfrak{f})$ and all its terms are
continuous on $\mathcal{E}_{I}(\mathfrak{f})$\,.\\
For the surjectivity, the translation invariance allows to assume  $b=0$\,. Take
$\gamma=\sum_{\nu\in (2\nz^{*})^{-1/2}}\gamma_{\nu}e_{\nu}$ with
$\|\gamma\|_{\mathcal{D}_{-\frac 1 2}}^{2}=\sum_{\nu\in
  (2\nz^{*})^{-1/2}}|\nu||\gamma_{\nu}|_{\mathfrak{f}}^{2}<\infty$\,. Choose
a non negative cut-off function $\chi\in
\mathcal{C}^{\infty}_{0}((a,0])$ such that $\chi\equiv 1$ in the
neighborhood of $0$\,.
We use the presentation of Subsection~\ref{se.rewriting} and write
$Pu=f$ in the form $((\frac{1}{2}+\mathcal{O})^{-1}p\partial_{q}+1)u=\check{f}=(\frac{1}{2}+\mathcal{O})^{-1}f$\,.\\
For
$\nu<0$\,, take 
$$
u_{\nu}=\gamma_{\nu}\chi(\frac{q}{\nu})
\quad\text{and}\quad
\check{f}_\nu=(\nu\partial_{q}u_{\nu}+u_{\nu})=
\gamma_{\nu}(\chi'(\frac{q}{\nu})+\chi(\frac{q}{\nu}))\,.
$$
For $\nu>0$\,, take
$\check{f}_{\nu}(q)=2\gamma_{\nu}\chi(\frac{q}{\nu})$ and
\begin{eqnarray*}
u_{\nu}(q)&=&\frac{\gamma_{\nu}}{\nu\int_{-\infty}^{0}e^{s}\chi(s)~ds}\int_{-\infty}^{q}e^{-\frac{q-s}{\nu}}\chi(\frac{s}{\nu})~ds
\\
&=&
\frac{1}{\int_{-\infty}^{0}e^{s}\chi(s)~ds}
\left(\frac{e^{-\frac{.}{\nu}}}{\nu}1_{\rz_{+}}\right)
*(\chi(\frac{.}{\nu})1_{\rz_{-}})\gamma_{\nu}\,.
\end{eqnarray*}
With
$\left\|\left(\frac{e^{-\frac{.}{\nu}}}{\nu}1_{\rz_{+}}\right)\right\|_{L^{1}}=1$
  and
$\|\gamma_{\nu}\phi(\frac{.}{\nu})\|_{L^{2}(\rz_{-},dq;\mathfrak{f})}^{2}=|\nu|
|\gamma_{\nu}|_{\mathfrak{f}}^{2}\int_{-\infty}^{0}|\phi(s)|^{2}~ds$ with
$\phi=\chi$ or $\phi=\chi'+\chi$\,,
we deduce
$u, \check{f}\in L^{2}(\rz_{-},dq;\mathcal{H}^{1})$
and therefore $u\in  L^{2}(\rz_{-},dq;\mathcal{H}^{1})$
 and $Pu=f=(\frac{1}{2}+\mathcal{O})\check{f}\in
L^{2}(\rz_{-},dq;\mathcal{H}^{-1})$\,.
Meanwhile the integral
$\int_{-\infty}^{0}e^{\frac{s}{\nu}}\chi(s)~ds=\nu\int_{-\infty}^{0}e^{s}\chi(s)~ds$ for $\nu>0$\,,
leads to $u(q=0)=\gamma$\,. Our construction also ensures $\supp
u\subset \left\{q\in\supp
  \chi\right\}\subset\left\{q\in(a,b=0]\right\}$\,.\\
The surjectivity for the trace at $a$ when
$a>-\infty$ is deduced by the symmetry $(q,p)\to (a+b-q,-p)$
when $-\infty<a<b<+\infty$ or $(q,p)\to (2a-q,-p)$ when $-\infty<a$\,, $b=+\infty$\,.
\end{proof}

\subsubsection{Calderon projector and boundary value problem}
\label{se.caldproj}

The operator $\tilde{P}=\frac{1}{2}+P$ is a local (differential operator) on $\rz^{2}$
and we  construct the associated Calderon projector at $q=0$ for
boundary value problems on $\overline{\rz^{2}_{-}}$\,. For Cauchy
problems for ODE's with constant coefficients, the Calderon  projector
selects the relevant exponentially decaying modes. We recall and apply
the general definition for differential operators.
Let us first check the invertibility of $\tilde{P}=\frac{1}{2}+P$ on the whole space
$\rz^{2}$\,. 
\begin{proposition}
\label{pr.invertR}
  The operator $\tilde{P}=\frac{1}{2}+P=p\partial_{q}+(\frac{1}{2}+\mathcal{O})$ defines
 an isomorphism from $\mathcal{E}_{\rz}(\mathfrak{f})$ to
 $L^{2}(\rz,dq;\mathcal{H}^{-1})$\,.
\end{proposition}
\begin{proof}
When $u\in \mathcal{E}_{\rz}(\mathfrak{f})$\,, $(\frac{1}{2}+\mathcal{O})u\in
L^{2}(\rz,dq;\mathcal{H}^{-1})$ and 
\begin{eqnarray*}
\|\tilde{P}u\|_{L^{2}(\rz,dq;\mathcal{H}^{-1})}
&=&\|p\partial_{q}u
+(\frac{1}{2}+\mathcal{O})u\|_{L^{2}(\rz,dq;\mathcal{H}^{-1})}
\\
&\leq&
\|p\partial_{q}u\|_{L^{2}(\rz,dq;\mathcal{H}^{-1})}
+\|(\frac{1}{2}+\mathcal{O})u\|_{L^{2}(\rz,dq;\mathcal{H}^{-1})}
\leq 2\|u\|_{\mathcal{E}_{\rz}(\mathfrak{f})}\,.
\end{eqnarray*}
When $f\in L^{2}(\rz,dq;\mathcal{H}^{-1}\otimes \mathfrak{f})$\,, the
rewriting of
$\tilde{P}u=f$
as
$(\frac{1}{2}+\mathcal{O})^{-1}p\partial_{q}u+u=\check{f}=(\frac{1}{2}+\mathcal{O})^{-1}f\in
L^{2}(\rz,dq;\mathcal{H}^{-1}\otimes \mathfrak{f})$ like in
Subsection~\ref{se.rewriting}, combined with Lemma~\ref{le.estL1} (case
$s=0$)\,, provides a solution $u\in L^{2}(\rz,dq;\mathcal{H}^{1}\otimes
\mathfrak{f})$ to $\tilde{P}u=f$\,.\\
If there are two solutions $u_{1},u_{2}\in \mathcal{E}_{\rz}(\mathfrak{f})$ to
$\tilde{P}u=f$\,, the difference $u_{2}-u_{1}$ belongs to $\mathcal{E}_{\rz}(\mathfrak{f})$
and solves $p\partial_{q}(u_{2}-u_{1})=-(\frac{1}{2}+\mathcal{O})(u_{2}-u_{1})$\,.
The integration by parts \eqref{eq.IPP} with $a=-\infty$ and
$b=+\infty$ gives
$$
-\Real\langle (u_{2}-u_{1})\,,(\frac{1}{2}+\mathcal{O})(u_{2}-u_{1})\rangle=
\Real\langle (u_{2}-u_{1})\,,\, p\partial_{q}(u_{2}-u_{1}) \rangle=0\,.
$$
This means $\|u_{2}-u_{1}\|_{L^{2}(\rz,dq;\mathcal{H}^{1})}=0$ and $u_{2}=u_{1}$\,.
\end{proof}
\begin{definition}
  When $\mathfrak{g}$ is a separable Hilbert space, we call $r_{-}$ the
  restriction operator
$$
r_{-}:L^{2}(\rz,dq;\mathfrak{g})\to L^{2}(\rz_{-},dq;\mathfrak{g})\quad
r_{-}u=u\big|_{\rz_{-}}\,,
$$
and $e_{-}$ the extension by $0$\,, $e_{-}=r_{-}^{*}$\,,
$$
e_{-}:L^{2}(\rz_{-},dq;\mathfrak{g})\to
L^{2}(\rz,dq;\mathfrak{g})\quad,\quad e_{-}u=u\times 1_{\rz_{-}}\,.
$$
\end{definition}
\begin{proposition}
\label{pr.cald}
For $\tilde{P}=\frac{1}{2}+P$\,,
the expression
$$
\mathcal{K}\gamma=u-r_{-}{\tilde{P}}^{-1}e_{-}\tilde{P}u\quad \text{for}~\gamma_{0}u=\gamma\,,
$$
with $u\in \mathcal{E}_{\rz_{-}}(\mathfrak{f})$ defines
a continuous operator from $L^{2}(\rz,|p|dp;\mathfrak{f})$ to
$\mathcal{E}_{\rz_{-}}(\mathfrak{f})$ such that $\tilde{P}(\mathcal{K}\gamma)=0$ for all
$\gamma\in L^{2}(\rz,|p|dp;\mathfrak{f})$\,.\\
The operator $C_{0}=\gamma_{0}\circ \mathcal{K}$ is the orthogonal projection
$\frac{1-S}{2}=1_{\rz_{-}}(S)$ on
$L^{2}(\rz,|p|dp;\mathfrak{f})=\mathcal{D}_{-\frac 1 2}$\,.
\end{proposition}
\begin{proof}
 We know that $\gamma_{0}:u\in \mathcal{E}_{\rz_{-}}(\mathfrak{f})\to \gamma_{0}u\in
 L^{2}(\rz,|p|dp;\mathfrak{f})$ is surjective.\\
For every $u\in \mathcal{E}_{\rz_{-}}(\mathfrak{f})$\,, the continuity of
$\tilde{P}:\mathcal{E}_{\rz_{-}}(\mathfrak{f})\to L^{2}(\rz_{-},dq;\mathcal{H}^{-1})$ and Proposition~\ref{pr.invertR} ensure that 
$$
v_{u}=u-r_{-}{\tilde{P}}^{-1}e_{-}\tilde{P}u\,,
$$
belongs to $\mathcal{E}_{\rz_{-}}(\mathfrak{f})$ and solves $\tilde{P}v_{u}=0$ for
$q<0$\,.\\
For $\gamma\in L^{2}(\rz,|p|dp;\mathfrak{f})$ assume that $u_{1}\in
\mathcal{E}_{\rz_{-}}(\mathfrak{f})$ and $u_{2}\in \mathcal{E}_{\rz_{-}}(\mathfrak{f})$ satisfy
$\gamma_{0}u_{2}=\gamma_{0}u_{1}=\gamma$\,. The difference
$u_{2}-u_{1}$ has a null trace at $q=0$ and since $\tilde{P}$ is a first order
differential operator in $q$\,, one gets the equality
$e_{-}Pu=P(e_{-}u)$\,. The invertibility of $\tilde{P}:\mathcal{E}_{\rz}(\mathfrak{f})\to
L^{2}(\rz_{-},dq;\mathcal{H}^{-1})$ (already used to
define $v_{u}$) provides the nullity 
$$
v_{u_{2}}-v_{u_{1}}=u_{2}-u_{1}-r_{-}\tilde{P}^{-1}\tilde{P}e_{-}(u_{2}-u_{})=0\quad\text{for}~q<0\,.
$$
The possibility to choose $u\in \mathcal{E}_{\rz_{-}}(\mathfrak{f})$ such that
$\gamma_{0}u=\gamma$ and
$\|u\|_{\mathcal{E}_{\rz_{-}}(\mathfrak{f})}\leq
C\|\gamma\|_{L^{2}(\rz,|p|dp;\mathfrak{f})}$
according to Proposition~\ref{pr.intbypart}, combined with
$\|v_{u}\|_{\mathcal{E}_{\rz_{-}}(\mathfrak{f})}\leq
C'\|u\|_{\mathcal{E}_{\rz_{-}}(\mathfrak{f})}$\,, implies the continuity of
$\mathcal{K}$\,.\\
With the rewriting of $\tilde{P}u=f$ as
$(\frac{1}{2}+\mathcal{O})^{-1}\tilde{P}u=\check{f}=(\frac{1}{2}+\mathcal{O})^{-1}f$\,,
$\mathcal{K}\gamma$ equals
$$
\mathcal{K}\gamma=u-r_{-}\tilde{P}^{-1}e_{-}\tilde{P}u=u-r_{-}{\tilde{P}}^{-1}(\frac{1}{2}+\mathcal{O})e_{-}(\frac{1}{2}+\mathcal{O})^{-1}\tilde{P}u\,.
$$
After introducing the orthogonal decomposition $((\cz e_{\nu})\otimes\mathfrak{f})_{\nu\in
  \pm(2\nz^{*})^{-\frac 1 2}}$ of $\mathcal{H}^{1}$ and
$\mathcal{D}_{-\frac 1 2}=L^{2}(\rz,|p|dp)$\,, it becomes
$$
(\mathcal{K}\gamma)_{\nu}(q)=u_{\nu}(q)-\left[r_{-}(\nu\partial_{q}+1)^{-1}e_{-}(\nu\partial_{q}+1)u_{\nu}\right](q)\,,\quad
\forall \nu\in \pm(2\nz^{*})^{-\frac 1 2}\,.
$$

All the $u_{\nu}$'s belong to $H^{1}_{loc}(\rz_{-})$ and the equality
$$
e_{-}(\nu\partial_{q}+1)u_{\nu}=(\nu\partial_{q}+1)e_{-}u_{\nu}+\nu u_{\nu}(0)\delta_{0}\,,
$$
holds in the distributional sense in $\rz$\,. 
The general solution on $\rz$ to $(\nu\partial_{q}+1)v=\nu
u_{\nu}(0)\delta_{0}$ in $\rz$ is
$$
1_{\rz_{-}}(q)e^{-\frac{q}{\nu}}C+(u_{\nu}(0)+C)1_{\rz_{+}}(q)e^{-\frac{q}{\nu}}\,,
\quad C\in\mathfrak{f}\,.
$$
Thus, the only possible
value for
$(\nu\partial_{q}+1)^{-1}e_{-}(\nu\partial_{q}+1)u_{\nu}$  is given by
\begin{eqnarray*}
 &u_{\nu}(q)1_{\rz_{-}}(q)+u_{\nu}(0)e^{-\frac{q}{\nu}}1_{\rz_{+}}(q)\,,&\text{for}~\nu>0\,,\\
&u_{\nu}(q)1_{\rz_{-}}(q)-u_{\nu}(0)1_{\rz_{-}}(q)e^{-\frac{q}{\nu}}\,,&\text{for}~\nu<0\,.
\end{eqnarray*}
We finally obtain after taking the restriction to $\rz_{-}$\,,
$$
(K\gamma)_{\nu}(q)=
\left\{
  \begin{array}[c]{ll}
    0&\text{for}~\nu>0\,,\\
    u_{\nu}(0)1_{\rz_{-}}(q)e^{-\frac{q}{\nu}}=\gamma_{\nu}1_{\rz_{-}}(q)e^{-\frac{q}{\nu}}&\text{for}~\nu<0\,.
  \end{array}
\right.
$$
 Taking the
trace at $q=0$ yields
$$
(C_{0}\gamma)_{\nu}=\gamma_{\nu}1_{\rz_{-}}(\nu)=\left(1_{\rz_{-}}(S)\gamma\right)_{\nu}\,,\quad
\forall \nu\in \pm (2\nz^{*})^{-\frac 1 2}\,.
$$
\end{proof}
\begin{remark}
  In general the local nature of $\tilde{P}$ suffices to prove that the
  Calderon projector $C_{0}$ fulfills $C_{0}\circ C_{0}=C_{0}$\,. Here
  (actually all the construction was made for this) we identify
  directly $C_{0}=1_{\rz_{-}}(S)=\frac{1-S}{2}$\,. The range of $C_{0}$ is nothing
  but the kernel of $1_{\rz_{+}}(S)=\frac{1+S}{2}$\,.\\
The operator $\mathcal{K}$ is usually called a Poisson operator.
\end{remark}
The Calderon projector is used to study the well-posedness of boundary value
problems (see \cite{ChPi}\cite{BdM}\cite{HormIII}-Chap~20 or \cite{Tay}). This
is summarized in the following straightforward proposition.
\begin{proposition}
\label{pr.BVPCald}
Set $\tilde{P}=\frac{1}{2}+P$\,. Let  $\mathcal{T}$
   be a Hilbert spaces  and let
 $T:\mathcal{D}_{-\frac 1 2}\to \mathcal{T}$ be a
 continuous operator. 
Consider the
  boundary value problem
  \begin{equation}
    \label{eq.BVPCald}
 \left\{
  \begin{array}[c]{ll}
    \tilde{P}u=f\in L^{2}(\rz_{-},dq; \mathcal{H}^{-1})&
    \\
T\gamma_{0} u=f_{\partial}\in \mathcal{T}\,,
  \end{array}
\right.
\quad u\in \mathcal{E}_{\rz_{-}}(\mathfrak{f})\,.
\end{equation}
Call $T_{-}$ the restriction $T\big|_{\Ran C_{0}}=T\big|_{\ker(1+S)}$\,.
When $T_{-}$ is injective the boundary value problem \eqref{eq.BVPCald} admits at most one
solution in $\mathcal{E}_{\rz}(\mathfrak{f})$\,.
When $T_{-}$ is surjective the boundary value problem \eqref{eq.BVPCald}
has solutions for any $f_{\partial}\in \mathcal{T}$ and $f\in
L^{2}(\rz,dq;\mathcal{H}^{-1})$\,.
When $T_{-}$ is an isomorphism, then \eqref{eq.BVPCald} admits a unique
solution 
\begin{eqnarray*}
&&u=r_{-}\tilde{P}^{-1}e_{-}f+ \mathcal{K}T_{-}^{-1}[f_{\partial}-T\gamma_{0}(r_{-}\tilde{P}^{-1}e_{-}f)]\,.
\\
\text{with}&&
\|u\|\leq C\left[\|f\|_{L^{2}(\rz_{-}, \mathcal{H}^{-1})}+\|f_{\partial}\|_{\mathcal{T}}\right]\,.
\end{eqnarray*}
\end{proposition}
Since $\Ran C_{0}$ is the kernel of the simple operator $\frac{(1+S)}{2}=1_{\rz_{+}}(S)$\,, one
can study the well posedness of the boundary value problem by
considering $T:\mathcal{D}_{-\frac 1
  2} =L^{2}(\rz,|p|dp;\mathfrak{f}) \to \mathcal{T}$ and by studying whether the
system
\begin{equation}
  \label{eq.BVPCald2}
\left\{
  \begin{array}[c]{l}
    T\gamma=f_{\partial}'\in \mathcal{T}\\
    (1+S)\gamma=0
  \end{array}
\right.
\end{equation}
admits a unique solution $\gamma\in L^{2}(\rz,|p|dp;\mathfrak{f})=\mathcal{D}_{-\frac 1
  2}$ for any $f_{\partial}'\in \mathcal{T}$\,.
A typical example is $\mathcal{T}=\Ran C_{0}$ and
$T=C_{0}=\frac{1-S}{2}=1_{\rz_{-}}(S)$\,. 
The boundary value problem
$$
\left\{
  \begin{array}[c]{l}
    \tilde{P}u=f\,,\\
   1_{\rz_{-}}(S)\gamma_{0} u=f_{\partial}\,,
  \end{array}
\right.
$$
is well-posed in $\mathcal{E}_{\rz_{-}}(\mathfrak{f})$ for all  $f\in
L^{2}(\rz_{-},dq;\mathcal{H}^{-1})$ and all $f_{\partial}\in \Ran
C_{0}=\Ran 1_{\rz_{-}}(S)$\,.\\
A solution $u\in \mathcal{E}_{\rz_{-}}(\mathfrak{f})$ to $\tilde{P}u=f\in
L^{2}(\rz_{-},dq;\mathcal{H}^{-1})$ coincides with the
restriction to $\rz_{-}$ of ${\tilde{P}}^{-1}e_{-}f={\tilde{P}}^{-1}(f 1_{\rz_{-}}(q))$ if
and  only if
$C_{0}u=0$\,,
that is $1_{\rz_{-}}(S)u=0$\,.

\subsection{Maximal accretivity}
\label{se.maxacc1D}
In order to prove Theorem~\ref{th.maxaccKA}, we study the boundary value problem 
$$
Pu=f
\quad,\quad
\gamma_{odd}u=\sign(p)\times A\gamma_{ev}u\,,
$$
in a larger framework by considering $f\in
L^{2}(\rz_{-},dq;\mathcal{H}^{-1})$ and $u\in
L^{2}(\rz_{-},dq;\mathcal{H}^{1})$ better suited for functional analysis
arguments.\\
Remember $\gamma u= \gamma_{0}u=u\big|_{q=0}$ and $\gamma_{ev,
  odd}u=\Pi_{ev,odd}\gamma u$\,.

\subsubsection{Boundary value problem related with $A$}
\label{se.sysODE}
Although they are equal, the Hilbert spaces $L^{2}(\rz,|p|dp,\mathfrak{f})$ and
$\mathcal{D}_{-\frac{1}{2}}$ are endowed with two different scalar products.
When $(A,D(A))$ is maximal accretive in
$L^{2}(\rz,|p|dp;\mathfrak{f})$ commuting with $\Pi_{ev}$\,, the same
holds for $(MA,D(MA))$ in $D_{-\frac{1}{2}}$\,, with $M=S\circ \sign(p)$
and $D(MA)=D(A)$\,,
by  Lemma~\ref{le.equaccr} and the commutation $M\Pi_{ev}=\Pi_{ev}M$\,.
We shall first work in $\mathcal{D}_{-\frac{1}{2}}$\,. For the
corresponding scalar product, the adjoint 
$((MA)^{*},D((MA)^{*}))$ is also maximal accretive owing to
Hille-Yosida theorem and
$\|(\lambda+(MA)^{*})^{-1}\|_{\mathcal{L}(\mathcal{D}_{-\frac 1 2})}=\|(\lambda+MA)^{-1}\|_{\mathcal{L}(\mathcal{D}_{-\frac 1 2})}\leq
\frac{1}{\lambda}$ for $\lambda>0$\,. The operator $(1+MA)$ (resp. $1+MA^{*}$) defines
isomorphisms in the diagram
\begin{eqnarray*}
&&  D(A)\stackrel{1+MA}{\to}\mathcal{D}_{-\frac{1}{2}}\stackrel{1+MA}{\to}D((MA)^{*})'\\
\text{resp.}&&
  D((MA)^{*})\stackrel{1+(MA)^{*}}{\to}\mathcal{D}_{-\frac{1}{2}}
\stackrel{1+(MA)^{*}}{\to}D(MA)'\,,
\end{eqnarray*}
where $D((MA)^{*})'$ (resp. $D(MA)'$) is the dual of $D((MA)^{*})$
(resp. $D(MA)$)\,. This will be more systematically used in
Section~\ref{se.insvar}.
Those isomorphisms and the commutation $\Pi_{ev}MA=MA\Pi_{ev}$ ensure
$E=\Pi_{ev}E\oplus^{\perp} \Pi_{odd}E$ for $E=D(MA), D((MA)^{*}),$ 
$D(MA)',
D((MA)^{*})'$\,, endowed with the suitable scalar product.

\begin{proposition}
\label{pr.ODE}
Keep the notation $\tilde{P}=\frac{1}{2}+P$\,. 
Then for any $f\in
L^{2}(\rz_{-},dq;\mathcal{H}^{-1})$ and any $f_{\partial}\in
\Pi_{ev}D((MA)^{*})'$\,, 
the boundary value problem
\begin{equation}
\label{eq.ODE}
\tilde{P}u=f\quad,\quad S\gamma_{odd}u-MA\gamma_{ev}u=f_{\partial}\,,
\end{equation}
admits a unique solution in $u\in L^{2}(\rz_{-},dq;\mathcal{H}^{1})$\,,
which belongs to $\mathcal{E}_{\rz_{-}}(\mathfrak{f})\subset
\mathcal{C}^{0}_{b}((-\infty,0]; L^{2}(\rz, 
|p|dp;\mathfrak{f})$ and 
satisfies
$$
\|u\|_{\mathcal{E}_{\rz_{-}}(\mathfrak{f})}
\leq
C\left[\|f\|_{L^{2}(\rz_{-},dq;\mathcal{H}^{-1})}
+\|(1+MA)^{-1}f_{\partial}\|_{\mathcal{D}_{-\frac{1}{2}}}\right]\,.
$$
\end{proposition}
\begin{proof}
The conditions $u\in L^{2}(\rz_{-},dq;\mathcal{H}^{1})$ and
$\tilde{P}u\in L^{2}(\rz_{-},dq;\mathcal{H}^{-1})$ imply $u\in
\mathcal{E}_{\rz_{-}}(\mathfrak{f})$ and $\gamma u\in
\mathcal{D}_{-\frac{1}{2}}$\,.
Since $S\Pi_{odd}\gamma u=\Pi_{ev}S\gamma u$ and $MA\Pi_{ev}\gamma
u\in \Pi_{ev}D((MA)^{*})'$\,, the boundary condition makes sense for
$f_{\partial}\in \Pi_{ev}D((MA)^{*})'$\,.
According to Proposition~\ref{pr.BVPCald} and the reformulation
\eqref{eq.BVPCald2}, the boundary value problem \eqref{eq.ODE} is well
posed if 
the system
$$
\left\{
  \begin{array}[c]{l}
    S\gamma_{odd}-MA\gamma_{ev}=f_{\partial}'~\,,\quad\text{with}~\gamma_{ev,odd}=\Pi_{ev,odd}\gamma\,,\\
   (1+S)\gamma=0
  \end{array}
\right.
$$
admits a unique solution $\gamma\in \mathcal{D}_{-\frac 1 2}$ 
for all $f_{\partial}'\in \Pi_{ev}D((MA)^{*})'$\,.
After \eqref{eq.Seevodd}\eqref{eq.orthevodd} we checked
$S\Pi_{ev}=\Pi_{odd}S$ and the above system is equivalent to
$$
\left\{
  \begin{array}[c]{l}
    S\gamma_{odd}-MA\gamma_{ev}=f_{\partial}'~,\\
   \Pi_{ev}(1+S)\gamma=0\,,\\
\Pi_{odd}(1+S)\gamma=0
  \end{array}
\right.
\Leftrightarrow
\left\{
  \begin{array}[c]{l}
   (\Id +MA)\gamma_{ev}=-f_{\partial}'~,\\
   \gamma_{odd}=-S\gamma_{ev}\,,
  \end{array}
\right.
$$
which admits a unique solution $\gamma\in \mathcal{D}_{-\frac{1}{2}}$
such that $\|\gamma\|_{\mathcal{D}_{-\frac{1}{2}}}\leq
C_{1}\|f_{\partial}'\|_{D((MA)^{*})'}$\,. The final estimate is also
provided by Proposition~\ref{pr.BVPCald}.
\end{proof}
We now translate the previous result in the more usual structure on
 $L^{2}(\rz,|p|dp;\mathfrak{f})$\,. 
By making use of the integration by part formula~\eqref{eq.IPP}, the
estimate of $\|u\|_{L^{2}(\rz_{-},dq;\mathcal{H}^{1})}$ will
also be made more accurate.
\begin{proposition}
\label{pr.ODEbis} Set $\tilde{P}=(\frac{1}{2}+P)$ and assume $(A,D(A))$ maximal accretive in
$L^{2}(\rz,|p|dp;\mathfrak{f})$ and commuting with $\Pi_{ev}$\,.
For any $f\in L^{2}(\rz_{-},dq;\mathcal{H}^{-1})$ and 
$f_{\partial}\in \Pi_{odd}L^{2}(\rz,|p|dp;\mathfrak{f})$
the  boundary value problem
\begin{equation}
  \label{eq.ODEbis}
\tilde{P}u=f\quad,\quad \gamma_{odd}u-\sign(p)A\gamma_{ev}u=f_{\partial}  
\end{equation}
admits a unique solution  
$u\in L^{2}(\rz_{-},dq;\mathcal{H}^{1})$ which belongs to $\mathcal{E}_{\rz_{-}}(\mathfrak{f})\subset
\mathcal{C}^{0}_{b}((-\infty,0]; L^{2}(\rz, 
|p|dp;\mathfrak{f})$ with $\gamma_{ev}u\in D(A)$\,.
This solution satisfies
\begin{multline*}
\Real\langle \gamma_{ev}u\,,\,
A\gamma_{ev}u\rangle_{L^{2}(\rz,|p|dp;\mathfrak{f})}
+\|u\|_{L^{2}(\rz_{-},dq;\mathcal{H}^{1})}^{2}
=\Real\langle
f\,,\, u\rangle
\\
-\Real\int_{\rz}\langle f_{\partial}(p),\gamma_{ev}u(p)\rangle_{\mathfrak{f}}~pdp\,.
\end{multline*}
When $f_{\partial}=0$ this implies 
$$
\|u\|_{L^{2}(\rz_{-},dq;\mathcal{H}^{1})}\leq
\|f\|_{L^{2}(\rz_{-},dq;\mathcal{H}^{-1})}\,.
$$
\end{proposition}
\begin{proof}
The boundary value problem \eqref{eq.ODE} with $f_{\partial}\in
\Pi_{ev}\mathcal{D}_{-\frac{1}{2}}$ admits a unique solution such that
$\gamma_{ev}u\in \mathcal{D}_{-\frac{1}{2}}$ and 
$$
MA\gamma_{ev}u=f_{\partial}-S\gamma_{odd}u \in \mathcal{D}_{-\frac{1}{2}}\,.
$$
Therefore $\gamma_{ev}u\in D(MA)=D(A)$ and, with $M=S\circ \sign(p)$\,, the boundary value problem
\eqref{eq.ODE} is equivalent to
$$
\tilde{P}u=f\quad,\quad \gamma_{odd}u-\sign(p)A\gamma_{ev}u=Sf_{\partial}\,,
$$
which is \eqref{eq.ODEbis} after changing $Sf_{\partial}$ into
$f_{\partial}$\,.\\
We can now use the integration by
parts formula \eqref{eq.IPP} which gives
\begin{eqnarray*}
0 \leq \Real\langle \gamma_{ev}u\,,\,
A\gamma_{ev}u\rangle_{L^{2}(\rz,|p|dp;\mathfrak{f})}
&&=
+\Real\langle \gamma_{ev}u\,,\,
\sign(p)\gamma_{odd}u\rangle_{L^{2}(\rz,|p|dp; \mathfrak{f})}
\\
&&
\qquad -\Real\langle \gamma_{ev}u\,,\, \sign(p)f_{\partial}\rangle_{L^{2}(\rz,|p|dp; \mathfrak{f})}
\\
&&\hspace{-3cm}=
\frac 1 2 \int_{\rz}|u(0,p)|_{\mathfrak{f}}^{2}~pdp
-\Real\int_{\rz}\langle f_{\partial}(p)\,,\,\gamma_{ev}u(p)\rangle_{\mathfrak{f}}~pdp
\\
&&
\hspace{-3cm}=
\Real\langle u\,, p\partial_{q}u\rangle 
-\Real\int_{\rz}\langle f_{\partial}(p)\,,\,\gamma_{ev}u(p)\rangle_{\mathfrak{f}}~pdp
\\
&&
\hspace{-3cm}=
\Real\langle u\,,
Pu\rangle-\|u\|_{L^{2}(\rz,dp;\mathcal{H}^{1})}^{2}
-\Real\int_{\rz}\overline{f_{\partial}(p)}\gamma_{ev}u(p)~pdp\,.
\end{eqnarray*}
\end{proof}
\begin{remark}
  Note that only the even part $A_{ev}=\Pi_{ev}A\Pi_{ev}$ appears in
  the assumption and all the analysis. So only this restriction can be
  considered  with an arbitrary maximal accretive extension on
  $\Ran\Pi_{odd}$\,.
 In applications, it is easier to consider
  a realization of $(A,D(A))$ in $L^{2}(\rz,|p|dp;\mathfrak{f})$
  without the parity in the domain definition and just check $A\Pi_{ev}=\Pi_{ev}A$\,.
\end{remark}
\subsubsection{Maximal accretivity of $K_{A}$}
\label{se.domKA}
We end here the proof of Theorem~\ref{th.maxaccKA}.
\begin{proof}
  The domain $D(K_{A})$ contains
  $\mathcal{C}^{\infty}_{0}((-\infty,0)\times
  (\rz\setminus\left\{0\right\});\mathfrak{f})$\,. It is dense in
  $L^{2}(\rz_{-}^{2},dqdp;\mathfrak{f})$\,.\\
From Proposition~\ref{pr.ODEbis}  applied with $f_{\partial}=0$\,, we know that
$\|u\|_{E_{A}}=\|(\frac{1}{2}+P)u\|_{L^{2}(\rz_{-},dq;\mathcal{H}^{-1})}$ is a norm on
 $$
E_{A}=\left\{u\in
  L^{2}(\rz_{-},dq;\mathcal{H}^{1})\,, Pu\in
  L^{2}(\rz_{-},dq;\mathcal{H}^{-1})\,, \gamma_{odd}u=\sign(p)A
  \gamma_{ev}u\right\}
\,,
$$
which is made a Banach space
isomorphic via $(\frac{1}{2}+P)$ to $L^{2}(\rz,dq;\mathcal{H}^{-1})$ and continuously
embedded in 
\begin{eqnarray*}
  &&\mathcal{E}_{\rz_{-}}(\mathfrak{f})\subset 
L^{2}(\rz_{-},dq;\mathcal{H}^{1})\cap
\mathcal{C}^{0}_{b}((-\infty,0];L^{2}(\rz,|p|dp;\mathfrak{f})\\
\text{with}&&
(u\in E_{A})\Rightarrow (\gamma_{ev}u \in D(A))\,.
\end{eqnarray*}
Since $L^{2}(\rz^{2}_{-},dqdp;\mathfrak{f})$ is complete and continuously embedded in
$L^{2}(\rz_{-},dq;\mathcal{H}^{-1})$\,, the domain $D(K_{A})$ endowed with the graph norm
$\|u\|_{D(K_{A})}=\|u\|+\|(\frac{1}{2}+P)u\|$ is complete, which means
that $K_{A}$ is closed.\\
By Proposition~\ref{pr.ODEbis} we also know that $\frac{1}{2}+P:E_{A}\to
L^{2}(\rz,dq;\mathcal{H}^{-1})$ is surjective, which implies $\Ran
K_{A}=L^{2}(\rz_{-}^{2},dqdp;\mathfrak{f})$\,.\\
For the accretivity and the integration by part formula
\eqref{eq.intbypartTh1}, 
simply use Proposition~\ref{pr.ODEbis} with $f_{\partial}=0$\,:
\begin{equation*}
 \Real\langle u\,,\, (K_{A}-\frac{1}{2})u\rangle
=
\Real \langle \gamma_{ev}u\,,\, A\gamma_{ev}u\rangle_{L^{2}(\rz,|p|dp;\mathfrak{f})}
+ 
\|u\|_{L^{2}(\rz_{-},dq;\mathcal{H}^{1})}^{2}-\|u\|^{2}\geq 0\,.
\end{equation*}
\end{proof}

\subsection{Extension of the resolvent and adjoint}
\label{se.extresadj}
We end the proofs of Proposition~\ref{pr.extH-1H1} and
Proposition~\ref{pr.adj1}.
\begin{proof}[Proof of Proposition~\ref{pr.extH-1H1}]
For $f\in L^{2}(\rz_{-},dq;\mathcal{H}^{1})$ and $z\in\cz$  a solution
$u\in L^{2}(\rz_{-},dq;\mathcal{H}^{1})$
$$
(P-z)u=f\quad,\quad \gamma_{odd}f=\sign(p)A\gamma_{ev}u
$$
satisfies $p\partial_{q}u=f+zu-\mathcal{O}u\in
L^{2}(\rz_{-},dq;\mathcal{H}^{-1})$\,. It necessarily belongs to
$\mathcal{E}_{\rz_{-}}(\mathfrak{f})$ and satisfies $\gamma_{ev}u\in
D(A)$\,.\\
\noindent\textbf{Uniqueness:} If there are two solutions $u_{1},u_{2}$, the
difference $u=u_{2}-u_{1}$\,, solves $(P-z)u=0\in
L^{2}(\rz^{2}_{-},dqdp;\mathfrak{f})$\,. When $\Real z<\frac{1}{2}$\,,
it means $u=(K_{A}-z)^{-1}0=0$\,.\\
\noindent\textbf{Existence:} From the integration by part of
Theorem~\ref{th.maxaccKA}
an element $u\in D(K_{A})$ satisfies
\begin{eqnarray*}
  \Real\langle u\,,\, (K_{A}-z)u\rangle
&\geq&\|u\|_{L^{2}(\rz_{-},dq;\mathcal{H}^{1})}^{2}-(\frac{1}{2}+\Real
z)\|u\|^{2}
\\
&\geq&(\frac{1}{2}-\Real z)\|u\|^{2}\,,
\end{eqnarray*}
and for $\Real z<\frac{1}{2}$
$$
\left(1+\frac{1}{\frac{1}{2}-\Real
    z}\right)\|u\|_{L^{2}(\rz_{-},dq;\mathcal{H}^{1})}\|(K_{A}-z)u\|_{L^{2}(\rz_{-},dq;\mathcal{H}^{-1})}
\geq \|u\|_{L^{2}(\rz_{-},dq;\mathcal{H}^{1})}^{2}\,.
$$
Hence $(K_{A}-z)^{-1}$ has a unique continuous extension from
$L^{2}(\rz_{-},dq;\mathcal{H}^{-1})$ to
$L^{2}(\rz_{-},dq;\mathcal{H}^{1})$\,.\\
For $f\in L^{2}(\rz_{-},dq;\mathcal{H}^{-1})$ and $\Real
z<\frac{1}{2}$\,, 
the function $u=(K_{A}-z)^{-1}f$ satisfies $(P-z)u=f$ in
$\mathcal{D}'(\rz^{2}_{-};\mathfrak{f})$\,. It belongs to
$\mathcal{E}_{\rz_{-}}(\mathfrak{f})$ and the traces
$\gamma_{ev,odd}u$ are well defined. When
$f=\lim_{n\to\infty}f_{n}$ in $L^{2}(\rz_{-},dq;\mathcal{H}^{-1})$
with $f_{n}\in L^{2}(\rz^{2}_{-},dqdp;\mathfrak{f})$\,, the sequence
$u_{n}=(K_{A}-z)^{-1}f_{n}$ satisfies
\begin{eqnarray*}
&&\gamma_{ev}u_{n}\in D(A)\quad,\quad
  \gamma_{odd}u_{n}=\sign(p)A\gamma_{ev}u_{n}\,,\\
&&\lim_{n\to
  \infty}\|\gamma_{ev}u-\gamma_{ev}u_{n}\|_{L^{2}(\rz,|p|dp;\mathfrak{f})}=0\,,\\
&&
\lim_{n\to
  \infty}\|\sign(p)\gamma_{odd}u-A\gamma_{ev}u_{n}\|_{L^{2}(\rz,|p|dp;\mathfrak{f})}=0\,.
\end{eqnarray*}
Hence $\gamma_{ev}u\in D(A)$ and $\gamma_{odd}u=\sign(p)A\gamma_{ev}u$\,.
\end{proof}
\begin{proof}[Proof of Proposition~\ref{pr.adj1}]
By changing $p$ into $-p$ and $A$ into $A^{*}$\,, the operator
$K_{-,A^{*}}$ defined by 
\begin{eqnarray*}
&D(K_{-,A^{*}})=\left\{u\in L^{2}(\rz_{-},dq;\mathcal{H}^{1})\,,~
  \begin{array}[c]{l}
P_{-}u\in
    L^{2}(\rz^{2}_{-},dqdp;\mathfrak{f})\,,\\
    \gamma_{odd}u=-\sign(p)A^{*}\gamma_{ev}u 
  \end{array}
\right\}\,,\\
&\forall u\in D(K_{A}^{*})\,,\quad
K_{A}^{*}u=P_{-}u=(-p\partial_{q}+\frac{-\partial_{p}^{2}+p^{2}}{2})u\,, 
\end{eqnarray*}
is maximal accretive. Hence it suffices to prove $K_{-,A^{*}}\subset
K_{A}^{*}$\,, which means
\begin{equation}
  \label{eq.adj1}
\forall u\in D(K_{-,A^{*}})\,,\, \forall v\in D(K_{A})\,,\quad
\langle K_{-,A^{*}}u\,,\, v\rangle=\langle u\,,\, K_{A}v\rangle\,.
\end{equation}
Both $u\in D(K_{-,A^{*}})$ and $v\in D(K_{A})$ belong to
$\mathcal{E}_{\rz_{-}}(\mathfrak{f})$ and the polarized integration by
parts formula of Proposition~\ref{pr.intbypart} is
\begin{multline*}
\langle p\partial_{q}u\,,\,v\rangle+\langle u\,,\,
p\partial_{q}v\rangle
=\langle u\,,\,\sign(p)v\rangle_{L^{2}(\rz,|p|dp;\mathfrak{f})}\\
=\langle
\sign(p)\gamma_{odd}u\,,\,\gamma_{ev}v\rangle_{L^{2}(\rz,|p|dp;\mathfrak{f})}
+
\langle
\gamma_{ev}u\,,\,\sign(p)\gamma_{odd}v\rangle_{L^{2}(\rz,|p|dp;\mathfrak{f})}
\\
=\langle
-A^{*}\gamma_{ev}u\,,\,\gamma_{ev}v\rangle_{L^{2}(\rz,|p|dp;\mathfrak{f})}
+
\langle
\gamma_{ev}u\,,\,A\gamma_{ev}v\rangle_{L^{2}(\rz,|p|dp;\mathfrak{f})}=0\,,
\end{multline*}
which yields \eqref{eq.adj1}.
\end{proof}

\section{Cuspidal semigroups}
\label{se.cusp}
From the functional analysis point of view, the Kramers-Fokker-Planck
operators have some specific properties, which were already used in
\cite{HerNi}\cite{EcHa} and \cite{HelNi} in the case without boundary.
 Those properties behave well 
after tensorisation and are stable under relevant perturbations.
\subsection{Definition and first properties}
\label{se.defcusp}
We work here in a complex Hilbert space $\mathcal{H}$ with the scalar product
$\langle\,,\,\rangle$ and norm $\|~\|$\,.
A maximal accretive operator $(K,D(K))$ is a densely defined operator
such 
$$
\forall u\in D(K)\,,~\Real \langle u\,,\, Ku\rangle\geq 0\,,
$$
and $1+K:D(K)\to \mathcal{H}$ is a bijection. (see
\cite{Bre}\cite{ReSi}\cite{EnNa}). This is equivalent the fact that
$(e^{-tK})_{t\geq 0}$ is a strongly continuous semigroup of
contractions. This property passes to the adjoint $(K^{*},D(K^{*}))$
(use $\|(\lambda+K^{*})^{-1}\|=\|(\lambda+K)^{-1}\|\leq
\frac{1}{\lambda}$ and Hille-Yosida theorem). The operator $(1+K)$
(resp. $1+K^{*}$) defines isomorphisms in the diagram
\begin{eqnarray*}
  && D(K)\stackrel{1+K}{\to}\mathcal{H}\stackrel{1+K}{\to}
  D(K^{*})'\,,\\
&&  D(K^{*})\stackrel{1+K^{*}}{\to}\mathcal{H}\stackrel{1+K^{*}}{\to}
  D(K)'\,,
\end{eqnarray*}
where $D(K^{*})'$ (resp. $D(K)'$) is the dual of $D(K)$
(resp. $D(K^{*})$) and the density of $D(K)$ (resp. $D(K^{*})$)
provides
the embedding $\mathcal{H}\subset D(K)'$
(resp. $\mathcal{H}\subset D(K^{*})'$)\,.
The operator $(1+K^{*})(1+K)$ is the self-adjoint positive operator
defined as a Friedrichs extension with
$$
D((1+K^{*})(1+K))=\left\{u\in
  D(K)\,, (1+K)^{*}(1+K)u\in \mathcal{H}\right\}\,.
$$
The modulus $|1+K|$ is the square root,
$|1+K|=\sqrt{(1+K^{*})(1+K)}$\,.\\
An accretive operator of which the closure is maximal accretive will
be said essentially maximal accretive (see \cite{HelNi}).
\begin{definition}
Let $(\Lambda,D(\Lambda))$ be a positive self-adjoint operator such
that $\Lambda\geq 1$ and let $r$ belong to $(0,1]$\,.
  A maximal accretive operator $(K,D(K))$ is said $(\Lambda,r)$-cuspidal
 if $D(\Lambda)$ is dense in $D(K)$
  endowed with the graph norm and if there exists a constant $C>0$
  such that
  \begin{eqnarray}
\label{eq.KuLu} 
 &&  \forall u\in D(\Lambda)\,,\quad \|Ku\|\leq C\|\Lambda u\|\,,\\
\label{eq.LruKu}
&& \forall u\in D(\Lambda)\,, \forall \lambda\in\rz\,,\quad
\|\Lambda^{r}u\|\leq C\left[\|(K-i\lambda)u\|+\|u\|\right]\,.
  \end{eqnarray}
\end{definition}
\begin{definition}
\label{de.pseudsp}
  We say that a maximal accretive operator $(K,D(K))$ is
  $r$-pseudospectral 
for some  $r\in (0,1]$\,, if there
  exists $C_{K}>0$ such that 
$$
\forall \lambda\in\rz,\quad \|(-1+i\lambda-K)^{-1}\|\leq C_{K}\langle \lambda\rangle^{-r}\,.
$$ 
\end{definition}
An easy reformulation of this definition is given by the
\begin{proposition}
\label{pr.pseudspeq}
A maximal accretive operator $(K,D(K))$ is $r$-pseudospectral
with $r\in(0,1]$\,, iff there exists $C_{K}'>0$
such that
\begin{itemize}
\item the spectrum of $K$ is contained in  $S_{K}\cap \left\{\Real z\geq
  0\right\}$ with
$$
S_{K}=\left\{z\in\cz\,, |z+1|\leq C_{K}'|\Real z+ 1|^{1/r}\,, \Real
  z\geq -1\right\}\,;
$$
\item for $z\not\in S_{K}$ with $\Real z \geq -1$\,, the resolvent
norm is estimated by
$$
\|(z-K)^{-1}\|\leq C'_{K}\langle z\rangle^{-r}\,.
$$
\end{itemize}
\end{proposition}
\begin{proof}
  The ``if'' part is obvious. For the ``only if'' part it suffices to
  use the first resolvent identity. More precisely write
$$
(K-\mu-i\lambda)=(K+1-i\lambda)\left[1-(1+\mu)(K+1-i\lambda)^{-1}\right]\,.
$$
 When
$\langle\lambda\rangle^{r}>2C_{K}|1+\mu|$ the right-hand side is
invertible and  $\|(K-\mu-i\lambda)^{-1}\|\leq 2\|(K+1-i\lambda)^{-1}\|$\,.
\end{proof}
Note that the case $r=1$ corresponds to the case of sectorial
operators.\\
Among other consequences of these properties, we shall prove that they
are essentially equivalent (with some loss in the exponent)\,.
\begin{theorem}
\label{th.cusppsd}
  For a maximal accretive operator $(K,D(K))$ the following statements
  satisfy
$$
(i)\Rightarrow (ii)\Rightarrow (iii)\,.
$$
\begin{description}
\item[(i)] $(K,D(K))$ is $(\Lambda,r)$-cuspidal, $r\in(0,1]$\,.
\item[(ii)] $(K,D(K))$ is $r$-pseudospectral.
\item[(iii)] $(K,D(K))$ is $(|1+K|,r')$-cuspidal for $r'< \frac{r}{2-r}$\,.
\end{description}
\end{theorem}
The first result taken from \cite{HerNi} is about $(i)\Rightarrow (ii)$\,.
\begin{proposition}
\label{pr.pseudsp}
When $(K,D(K))$ is $(\Lambda,r)$-cuspidal with exponent $r\in
(0,1]$\,, it is $r$-pseudospectral with the same exponent.
\end{proposition}
\begin{proof}
For $u\in D(K)$\,, we start from
$$
2\|(K-z)u\|^{2}+2|\Real z|^{2}\|u\|^{2}\geq \|(K-\Imag
z)u\|^{2}\,,
$$
which gives
$$
2\|(K-z)u\|^{2}+(1+2|\Real z|^{2})\|u\|^{2}\geq \|(K-\Imag
z)u\|^{2}+\|u\|^{2}\geq \frac{1}{2C}\|\Lambda^{r}u\|^{2}\,.
$$
The condition \eqref{eq.KuLu} provides the operator inequality
$1\leq (1+K)^{*}(1+K)\leq \Lambda^{2}$\,. The operator monotonicity of
$x\to x^{r}$ for $r\in [0,1]$ implies $1\leq
\left[(1+K^{*})(1+K)\right]^{r}\leq \Lambda^{2r}$ and
$$
\|(K-z)u\|^{2}+(1+\Real z)^{2}\|u\|^{2}\geq
\frac{1}{8C}\left[\|(K-z)u\|^{2}+\langle u\,,\, [(1+K^{*})(1+K)]^{r}u\rangle\right]\,.
$$
Lemma~B.1 of \cite{HerNi} then says
$$
4\left[\|(K-z)u\|^{2}+\langle u\,,\,
  [(1+K^{*})(1+K)]^{r}u\rangle\right]\geq |z+1|^{2r}\|u\|^{2}\,,
\forall u\in D(K)\,,
$$
as soon as $\Real z\geq -1$\,.
We have proved
$$
\forall u\in D(K)\,, \|(K-z)u\|^{2}\geq
\left[\frac{1}{32C}|z+1|^{2r}-(1+\Real z)^{2}\right]\|u\|^{2}\,,
$$
for all $z\in\cz$ such that $\Real z\geq -1$\,.
Any $z=(-1+i\lambda)$ with $\lambda\in\rz$  belongs to the
resolvent set of $K$ and we get the inequality
$$
\|(K+1-i\lambda)^{-1}\|\leq 8\sqrt{C}\langle \lambda\rangle^{-r}\,.
$$
\end{proof}
Remember the formula 
$$
e^{-tK}=\frac{1}{2i\pi}\int_{+i\infty}^{-i\infty}e^{-tz}(z-K)^{-1}~dz\,,
$$
 which can be understood for general maximal accretive operators as
 the  limit
$$
\forall \psi\in D(K)\,,~e^{-tK}\psi=\lim_{\varepsilon \to 0^{+}}\lim_{k\to
  \infty}\frac{1}{2i\pi}\int_{-\varepsilon +i\infty}^{-\varepsilon
  -i\infty}\frac{e^{-tz}}{1+\frac{z}{k}}(z-K)^{-1}\psi~dz\,,
$$
extended to any $\psi\in \mathcal{H}$ by density.\\
Proposition~\ref{pr.pseudspeq} allows  a contour deformation 
which provides a norm convergent integral for $t>0$\,.

\begin{proposition}
\label{pr.contour} Assume that $(K,D(K))$ is $r$-pseudospectral
with exponent
$r\in (0,1]$\,. Let $C_{K}'$
be the constant of Proposition~\ref{pr.pseudspeq} and let $\Gamma_{K}$
be the contour $\left\{|\Imag z|=2C_{K}|1+\Real z|^{\frac{1}{r}}\,,
  \Real z\geq -1\right\}$ oriented from $+i\infty$ to $-i\infty$\,.
Then for any $t>0$
\begin{equation}
  \label{eq.intGK}
e^{-tK}=\frac{1}{2i\pi}\int_{\Gamma_{K}}e^{-tz}(z-K)^{-1}~dz
\end{equation}
where the right-hand side is a norm convergent integral in $\mathcal{L}(\mathcal{H})$\,.
\end{proposition}
\begin{proof}
  For $\psi\in D(K)$\,, $t>0$ and $k\geq 2$ the function
  $\frac{e^{-tz}}{(1+\frac{z}{k})}(z-K)^{-1}\psi$ is a
  holomorphic function of $z$ in $\left\{z\in \cz\setminus S_{K}\,,
    \Real z >-2\right\}$ with
$$
\|\frac{e^{-t z}}{(1+\frac{z}{k})}(z-K)^{-1}\psi\|\leq
C_{k}\frac{e^{-t\Real z}}{\langle z\rangle^{2}}\|(1+K)\psi\|\,,
$$
when $\Real z\geq -1$\,, $z\in \cz\setminus S_{K}$\,. The contour
integral
$$
\int_{-\varepsilon +i\infty}^{-\varepsilon-i\infty}
\frac{e^{-tz}}{1+\frac{z}{k}}(z-K)^{-1}\psi~dz\,,
$$
for any $\varepsilon\in (0,1)$\,, can thus be deformed into
$$
\int_{\Gamma_{K}}\frac{e^{-tz}}{1+\frac{z}{k}}(z-K)^{-1}\psi~dz\,.
$$
With the inequalities
\begin{eqnarray}
\label{eq.z-K}
  && \|(z-K)^{-1}\|\leq C_{K}'\langle z\rangle^{-r}\leq C_{K}'\langle
  \Imag z\rangle^{-r}\,,\\
\label{eq.expdz}
&& |e^{-tz}|=e^{-t\Real z}\leq e^{t}e^{-t\frac{|\Imag
    z|^{r}}{2C_{K}'}}\;,\; |dz|=\left[1+\mathcal{O}(|\Imag z|^{2(r-1)}) \right]|d~\Imag z|
\end{eqnarray}
valid for all $z\in \Gamma_{K}$\,, we can take the limit as $k\to
\infty$ for any fixed $t>0$\,, in the integral $\int_{\Gamma_{K}}$\,.
 The convergence holds for any $\psi\in \mathcal{H}$ and the
integral \eqref{eq.intGK} is a norm convergent integral.
\end{proof}
A classical results (see for instance \cite{EnNa}) provides the
equivalence between:
\begin{itemize}
\item $(e^{-z A})_{z\in C\cup\left\{0\right\}}$ is a bounded analytic
  semigroup for some open convex cone $C$\,.
\item $A$ is sectorial (see \cite{EnNa} for a general definition).
\item For all positive $\tau$ the estimate $\|(-\tau+is -A)^{-1}\|\leq
  \frac{C}{|s|}$ holds for all $s\neq 0$\,.
\item The quantity $\sup_{t>0}\|tAe^{-tA}\|$ is finite. 
\end{itemize}
Cuspidal semigroup, and this is a key idea of \cite{HerNi}, is a
fractional version of the above notions. Actually
Proposition~\ref{pr.pseudspeq} says that $(1+K)$ is sectorial when
$r=1$\,.
 In this direction, the next
result completes Proposition~\ref{pr.contour}.

\begin{proposition}
\label{pr.fractps} 
If $(K,D(K))$ is $r$-pseudospectral with exponent $r\in
(0,1)$\,, then
\begin{equation}
\label{eq.frtpsK}
\sup_{t> 0}\|t^{\frac{2}{r}-1}(1+K)e^{-t(1+K)}\|<+\infty\,.
\end{equation}
If $(K,D(K))$ is $(\Lambda,r)$-cuspidal with $r\in(0,1)$\,, then
\begin{equation}
\label{eq.frtpsL}
\sup_{t> 0}\|t^{\frac{2}{r}-1}\Lambda^{r}e^{-t(1+K)}\|<+\infty\,.
\end{equation}
\end{proposition}
\begin{proof}
For $u\in D(K)$ and $t>0$\,,  the formula \eqref{eq.intGK} gives
$$
e^{-tK}Ku=\frac{1}{2i\pi}\int_{\Gamma_{K}}e^{-tz}(z-K)^{-1}Ku~dz=
\frac{1}{2i\pi}\int_{\Gamma_{K}}e^{-tz}z(z-K)^{-1}u~dz\,,
$$
because $(z-K)^{-1}Ku=-u+z(z-K)^{-1}u$\,. From the inequalities
\eqref{eq.z-K} and \eqref{eq.expdz} and $|z|\leq C''(1+|\Imag z|)$ along
$\Gamma_{K}$\,, we deduce
$$
\|e^{-tK}Ku\|\leq
C^{(3)}e^{t}\int_{0}^{+\infty}(1+\lambda)e^{-t\lambda^{r}}\lambda^{-r}~d\lambda
\leq C^{(4)}e^{t}\left[t^{1-\frac{1}{r}}+t^{1-\frac{2}{r}}\right]
$$
The density of $D(K)$ yields
$$
\forall t>0\,,\quad e^{-t}\|Ke^{-t(1+K)}\|\leq
C^{(4)}\left[t^{1-\frac{1}{r}}+t^{1-\frac{2}{r}}\right]e^{-t}\leq C^{(5)}t^{1-\frac{2}{r}}\,.
$$
while we know $e^{-t}\|e^{-t(1+K)}\|\leq e^{-2t}\leq
C_{r}t^{1-\frac{2}{r}}$\,.
We have proved that 
$\sup_{t\geq 0}\|t^{\frac{2}{r}-1}(2+K)e^{-t(2+K)}\|$
is finite. Replacing $K$ with $2K$\,, which is also $r$-pseudospectral 
owing to Proposition~\ref{pr.pseudspeq}, 
and $2t$ by $t$\,, finishes the proof of
\eqref{eq.frtpsK}.\\
When $(K,D(K))$ is $(\Lambda,r)$-cuspidal with $r<1$\,, then it is
$r$-pseudospectral
according to Proposition~\ref{pr.pseudsp}.
The second estimate then comes from
$$
\|\Lambda^{r}u\|\leq C\left[\|Ku\|+\|u\|\right]\leq C'\|(1+K)u\|\,.
$$
\end{proof}
Below is the converse implication of the second statement of Proposition~\ref{pr.fractps}.
\begin{proposition}
\label{pr.frtpsconv}
Let $(K,D(K))$ be a maximal accretive operator.
Assume that \eqref{eq.frtpsL} is finite for $r\in (0,1]$ and assume
the existence of an
operator $\Lambda\geq 1$ such that $D(\Lambda)$ is dense in $D(K)$ and
$$
\forall u\in D(\Lambda)\,,\quad\|Ku\|\leq \|\Lambda u\|\,.
$$
Then the operator $(K,D(K))$ is $(\Lambda,\theta r)$-cuspidal 
 $\theta\in (0,\frac{r}{2-r})\subset(0,r)$: there
exists $C_{\theta}>0$ such that
$$
\forall u\in D(K)\,,\forall \lambda\in\rz\,,\quad \|\Lambda^{\theta r}u\|\leq C_{\theta}\left[\|(K-i\lambda)u\|+\|u\|\right]\,.
$$
\end{proposition}
\begin{proof}
  For $t>0$, $\lambda\in \rz$ and $v\in \mathcal{H}$\,, the estimate
$$
\|\Lambda^{r}e^{-t(1+K-i\lambda)}v\|\leq C_{0}t^{1-\frac{2}{r}}\|v\|\,,
$$
is interpolated into
$$
\|\Lambda^{\theta r}e^{-t(1+K-i\lambda)}v\|\leq
C_{0}^{\theta}t^{\theta\frac{r-2}{r}}\|v\|\,,
$$
for any $\theta\in [0,1]$\,. Then the right-hand side of 
$$
\|\Lambda^{\theta r}e^{-t(K+2-i\lambda)}v\|\leq C_{0}^{\theta}e^{-t}t^{\theta\frac{r-2}{r}}\|v\|
$$
is integrable on $(0,+\infty)$ as soon as $\theta\frac{r-2}{r}>-1$\,.
From
$$
(2+K-i\lambda)^{-1}v=\int_{0}^{+\infty}e^{-t(2+K-i\lambda)}v~dt\,,
$$
we deduce $(2+K-i\lambda)^{-1}v\in D(\Lambda^{r\theta})$ and
$$
\|\Lambda^{r\theta}(2+K-i\lambda)^{-1}v\|\leq C_{\theta}\|v\|\,.
$$
Setting $u=(2+K-i\lambda)^{-1}v$ gives
$$
\|\Lambda^{r\theta }u\|\leq C_{\theta}\|2u +(K-i\lambda)u\|\leq C\left[\|(K-i\lambda)u\|+\|u\|\right]\,,
$$
for all $u\in D(K)$\,, as soon as $\theta\in (0,\frac{r}{2-r})$\,.
\end{proof}
\begin{proof}[Proof of $(ii)\Rightarrow (iii)$ in
  Theorem~\ref{th.cusppsd}.] 
Firstly, $D(|1+K|)$ is the form domain of $(1+K^{*})(1+K)$ and equals
$D(K)$\,, with
$$
\forall u\in D(K)\,\quad \|Ku\|^{2}\leq \|(1+K)u\|^{2}=\||1+K|u\|^{2}\,.
$$ 
The end is a variant of the previous argument. 
The relation \eqref{eq.frtpsK} can be written 
$$
\||1+K|e^{-t(1+K)}\|\leq Ct^{1-\frac{2}{r}}\,,
$$
and leads to
$$
\forall \lambda\in\rz\,,
\forall v\in \mathcal{H}\,, 
\||1+K|^{\theta}e^{-t(2+K-i\lambda)}\|\leq C_{\theta}e^{-t}t^{\theta\frac{r-2}{r}}\|v\|\,.
$$
By taking $\theta<\frac{r}{2-r}$\,, we deduce $(2+K-i\lambda)^{-1}v\in
D(|1+K|^{\theta})$ and 
$$
\||1+K|^{\theta}(2+K-i\lambda)^{-1}v\|\leq C_{\theta}'\|v\|\,.
$$
With $u=(2+K-i\lambda)^{-1}v$ this means
$$
\forall \lambda\in\rz\,, \forall u\in D(K)\,,\quad 
\||1+K|^{\theta}u\|\leq C_{\theta}''\left[\|(K-i\lambda)u\|+\|u\|\right]\,.
$$
\end{proof}
\subsection{Perturbation}
\label{se.perturbcusp}
The Theorem~X.50 of \cite{ReSi} says that a pair of accretive operators
$A$ and $C$  defined on a same dense domain $D\subset \mathcal{H}$\,, such that
$$
\forall u\in D\,,\quad \|(A-C)u\|\leq
a\left[\|Au\|+\|Cu\|\right]+b\|u\|\,,
$$
for some  fixed $a<1$ and $b>0$\,, satisfy~:
\begin{itemize}
\item the closures $\overline{A}$ and $\overline{C}$ have the same domain;
\item $\overline{A}$ is maximal accretive if and only if
  $\overline{C}$ is.
\end{itemize}
Of course this applies to $C=A+B$ when $B$ is a relatively bounded
perturbation
$$
\forall u\in D(A)\,,\quad \|Bu\|\leq a\|Au\|+b\|u\|\,.
$$
Below is the cuspidal version, which requires a uniform control with
respect to $i\lambda\in i\rz$\,.
\begin{proposition}
\label{pr.perturb}
Let $(K,D(K))$ be a $(\Lambda,r)$-cuspidal operator with exponent $r\in
(0,1]$\,. Assume that $B$ is a relatively bounded perturbation of $K$
such that
$$
\forall u\in D(K)\,,\forall \lambda\in\rz\,,\quad \|Bu\|\leq a\|(K-i\lambda)u\|+ b\|u\|
$$
with some fixed $a<1$ and $b>0$\,.
If $(K+B)$ is accretive on $D(K)$\,, then $(K+B,D(K))$ is
$(\Lambda,r)$-cuspidal\,.\\
The same statement holds for $r$-pseudospectral operators with
exponent $r\in (0,1]$\,.
\end{proposition}
\begin{proof}
The Theorem~X.50 of \cite{ReSi} says that $(K+B,D(K))$ is maximal
accretive.\\
For the inequalities, write simply for $u\in D(\Lambda)\subset
D(K)$ and $\lambda\in\rz$\,,
\begin{eqnarray*}
\|(K+B)u\|&&\leq \|Ku\|+\|Bu\|\leq (1+a)\|Ku\|+b\|u\|\leq
[C(1+a)+b]\|\Lambda u\|\,,
\\
\|(K+B-i\lambda)u\|&&\geq \|(K-i\lambda)u\|-\|Bu\|\geq
(1-a)\|(K-i\lambda)u\|-
b\|u\| 
\\
&&\geq (1-a)C^{-1}\|\Lambda^{r}u\|-(b+1)\|u\|\,.
\end{eqnarray*}
This yields
$$
\forall \lambda\in\rz\,,\, \forall u\in D(\Lambda)\,,
\|Ku\|\leq C'\|\Lambda u\|\quad,\quad \|\Lambda^{r}u\|\leq C'\left[\|(K-i\lambda)u\|+\|u\|\right]\,,
$$
with $C'=\max\left\{C(1+a)+b, \frac{C(b+1)}{1-a}\right\}$\,.\\
When $(K,D(K))$ is $r$-pseudospectral, the same lower bound
$\|(K+B-i\lambda)u\|\geq
(1-a)\|(K-i\lambda)u\|-
b\|u\|$ allows to conclude that $(K+B,D(K))$ is $r$-pseudospectral.
\end{proof}
\begin{corollary}
\label{co.perturb}
A maximal accretive operator $(K,D(K))$ is $(\Lambda,r)$-cuspidal 
 if and only if there exists a constant $C\in\rz$ such
that $(C+K,D(K))$ is $(\Lambda,r)$-cuspidal.\\
If $(K,D(K))$ is $(\Lambda,r)$-cuspidal with exponent $r\in (0,1]$ and $B$ is
defined on $D(K)$ with
$$
\forall u\in D(K)\,, \|Bu\|^{2}\leq C_{B}\Real\langle u\,,\, (1+K)u\rangle
$$
and $K+B$ accretive. Then $(K+B,D(K))$ is $(\Lambda,r)$-cuspidal\,.\\
The same statements hold for $r$-pseudospectral  operators.
\end{corollary}
\begin{proof}
  The first result comes readily from 
$$\|-Cu\|\leq
  0\times \|(K-i\lambda)u\|+|C|\times\|u\|\,.
$$
The second one comes from
\begin{eqnarray*}
\|Bu\|^{2}
&\leq& C_{B}\Real\langle u\,,\, (K-i\lambda)u\rangle\leq
C_{B}\|u\|\|(K-i\lambda)u\|\\
&\leq& \left[\frac{1}{2}\|(K-i\lambda)u\|+C_{B}\|u\|\right]^{2}\,.
\end{eqnarray*}
\end{proof}
\subsection{Tensorisation}
\label{se.tensor}
Take two Hilbert spaces $\mathcal{H}_{1}$ and
$\mathcal{H}_{2}$ and consider the Hilbert tensor product
$\mathcal{H}=\mathcal{H}_{1}\otimes \mathcal{H}_{2}$\,.
We first recall a result about maximal accretivity.
\begin{proposition}
\label{pr.maxaccK1K2}
 Assume that $(K_{1},D(K_{1}))$ and $(K_{2},D(K_{2}))$ are maximal accretive 
in $\mathcal{H}_{1}$ and $\mathcal{H}_{2}$\,, respectively.
Then the closure of $(K_{1}\otimes \Id_{\mathcal{H}_{2}}+\Id_{\mathcal{H}_{1}}\otimes
K_{2})$ initially defined on the algebraic tensor product
$D(K_{1})\stackrel{alg}{\otimes} D(K_{2})$\,, is maximal accretive\,.  
\end{proposition}
Then we will prove the cuspidal version.
\begin{proposition}
\label{pr.cuspK1K2}
  Assume that $(K_{1},D(K_{1}))$ and $(K_{2},D(K_{2}))$ are
  respectively  $(\Lambda_{1},r_{1})$ and $(\Lambda_{2},r_{2})$-cuspidal
  operators 
in $\mathcal{H}_{1}$ and $\mathcal{H}_{2}$\,.
Then the closure of $(K_{1}\otimes \Id_{\mathcal{H}_{2}}+\Id_{\mathcal{H}_{1}}\otimes
K_{2})$ initially defined on the algebraic tensor product
$D(K_{1})\stackrel{alg}{\otimes} D(K_{2})$\,, is
$(\Lambda,r)$-cuspidal with
\begin{eqnarray*}
  &&\Lambda=\Lambda_{1}\otimes \Id +\Id\otimes \Lambda_{2}\,,\\
\text{and}&&
r<\min\{\frac{r_{1}^{2}}{8-2r_{1}},\frac{r_{2}^{2}}{8-2r_{2}}\}\,.
\end{eqnarray*}
\end{proposition}
  \begin{proof}[Proof of Proposition~\ref{pr.maxaccK1K2} (for the sake of completeness).]
Consider the strongly semigroup of contractions $(e^{-tK_{1}}\otimes
e^{-tK_{2}})_{t\geq 0}$ on $\mathcal{H}=\mathcal{H}_{1}\otimes
\mathcal{H}_{2}$
 and call $(K,D(K))$ its maximal accretive generator.
The operator $(K,D(K))$ is the closure of $(K_{1}\otimes \Id_{\mathcal{H}_{2}}+\Id_{\mathcal{H}_{1}}\otimes
K_{2})$ defined on 
$D=D(K_{1})\stackrel{alg}{\otimes} D(K_{2})$~: Indeed,
for every pair $(\varphi_{1},\varphi_{2})\in D(K_{1})\times
D(K_{2})$\,, $\rz_{+}\ni t\mapsto e^{-tK}(\varphi_{1}\otimes
\varphi_{2})=(e^{-tK_{1}}\varphi_{1})\otimes (e^{-tK_{2}}\varphi_{2})$
is a $\mathcal{C}^{1}$-function. This implies $D\subset D(K)$\,.
For $s>0$\,, the expression
$$
A_{s}\varphi=\frac{1}{s}\int_{0}^{s}e^{-t(1+K)}\varphi~dt\,,
$$
defines a continuous operator from $\mathcal{H}$ to $D(K)$\,.
For $\varphi\in D(K)$\,, the relation
$$
(1+K)\varphi-(1+K)A_{s}\varphi=(1+K)\varphi+\frac{1}{s}\left(e^{-s(K+1)}\varphi-\varphi\right)
$$
implies
$$
\|\varphi-A_{s}\varphi\|_{D(K)}=\|(1+K)(\varphi-A_{s}\varphi)\|_{\mathcal{H}}\stackrel{s\to
0^{+}}{\to}0\,.
$$
Since $D(K_{j})$ is dense in $\mathcal{H}_{j}$ for $j=1,2$\,,
$\varphi$ can be approximated in $\mathcal{H}$ by an element
$\sum_{j=1}^{J}\varphi_{1,j}\otimes \varphi_{2,h}$ of $D$
with
$$
\|A_{s}\left(\varphi-\sum_{j=1}^{J}\varphi_{1,j}\otimes
  \varphi_{2,h}\right)\|_{D(K)}
\leq C_{s}\|\varphi-\sum_{j=1}^{J}\varphi_{1,j}\otimes
  \varphi_{2,h}\|_{\mathcal{H}}\,.
$$ 
Now the integral
$$
A_{s}\left(\sum_{j=1}^{J}\varphi_{1,j}\otimes
  \varphi_{2,h}\right)
=\sum_{j=1}^{J}\frac{1}{s}\int_{0}^{s}(e^{-tK_{1}}\varphi_{1,j})\otimes
(e^{-tK_{2}}\varphi_{2,j})~dt
$$
can be approximated  in $D(K)$ by a Riemann sum which belongs to
$D$\,.
Hence $D$ is dense in $D(K)$ endowed with the graph norm
$\|u\|_{D(K)}=\|(1+K)u\|_{\mathcal{H}}$\,.
It is a core for $(K,D(K))$\,.
\end{proof}
\begin{proof}[Proof of Proposition~\ref{pr.cuspK1K2}]
We know from Proposition~\ref{pr.maxaccK1K2} that 
$K=K_{1}\otimes \Id +  \Id \otimes K_{2}$ is essentially maximal
accretive on $D(K_{1})\otimes^{alg} D(K_{2})$\,.\\
Let us prove the cuspidal property. The operator $\Lambda_{j}\geq 1$,
$j=1,2$\,,
 is a self-adjoint operator in $\mathcal{H}_{j}$
 with $D(\Lambda_{j})$ dense in $D(K_{j})$ and such that
$$
\|K_{j}u\|_{\mathcal{H}_{j}}\leq C_{j}\|\Lambda_{j}u\|_{\mathcal{H}_{j}}\quad,\quad
\|\Lambda_{j}^{r_{j}}u\|_{\mathcal{H}_{j}}\leq C_{j}\left[\left\|(K_{j}-i\lambda)u\right\|_{\mathcal{H}_{j}}+\|u\|_{\mathcal{H}_{j}}\right]\,,
$$
hold for any $u\in D(\Lambda_{j})$ and $\lambda\in\rz$\,. Call
$\Lambda=\Lambda_{1}\otimes \Id +\Id\otimes \Lambda_{2}$ the
essentially self-adjoint operator on $D(\Lambda_{1})\otimes^{alg}
D(\Lambda_{2})$\,. Its domain $D(\Lambda)$ is dense in $D(K)$\,, and
$$
\forall u\in D(\Lambda)\,,\quad
\|Ku\|\leq C\|\Lambda u\|\,,
$$
is obtained after taking first $u\in D(\Lambda_{1})\otimes^{alg}
D(\Lambda_{2})$\,.\\
The inequality \eqref{eq.frtpsL} applied for $j=1,2$ gives the
uniform bound
$$
\|t^{\frac{2}{r_{j}}-1}\Lambda_{j}^{r_{j}}e^{-t(1+K_{j})}\|_{\mathcal{L}(\mathcal{H}_{j})}\leq C_{j}'\,.
$$
This yields the uniform bound
$$
\|(\Lambda_{1}^{r_{1}}+\Lambda_{2}^{r_{2}})e^{-t(2+K)}\|_{\mathcal{L}(\mathcal{H})}\leq
\max\{C_{1}'t^{1-\frac{2}{r_{1}}}, C_{2}'t^{1-\frac{2}{r_{2}}}\}e^{-t}
$$
We take $\varrho=\min\left\{r_{1},r_{2}\right\}$ and the functional calculus
of commuting self-adjoint operators gives
\begin{eqnarray*}
  &&
\Lambda^{2\varrho}\leq C_{r_{1},r_{2}}(\Lambda_{1}^{r_{1}}+\Lambda_{2}^{r_{2}})^{2}\,,\\
\text{and}&&
\|\Lambda^{\varrho}e^{-t(2+K)}\|\leq C_{r_{1},r_{2}}' t^{2-\frac{4}{\varrho}}e^{-t}\,,
\end{eqnarray*}
By interpolation we deduce
$$
\sup_{t>0}\|t^{\frac{2}{\varrho}-1}\Lambda^{\frac{\varrho}{2}}e^{-t(2+K)}\|<+\infty\,.
$$
Proposition~\ref{pr.frtpsconv} implies that $(K,D(K))$ is
$(\Lambda,r)$-cuspidal for any $r$ in $(0,\frac{\varrho^{2}}{8-2\varrho})$\,. 
\end{proof}

\section{Separation of variables}
\label{se.insvar}

The results of Subsection~\ref{se.tensor} provide a strategy for
studying the maximal accretivity and cuspidal property for 
operators which make possible a complete separation of
variables. Half-space problems associated with differential operators
with separated variables can be reduced to the one dimensional half-line
problem when the boundary conditions agree with this separation of variables.
We focus on the cases when the boundary
condition is given by $\gamma_{odd}u=\nu\sign(p)\gamma_{ev}u$ with
$\nu\in \left\{0,1\right\}$\,.
 
Although Proposition~\ref{pr.cuspK1K2} provides the guideline,
we do not apply it naively. Instead, we develop an
analysis of abstract one-dimensional problems which will lead in the
next section to more accurate results.
In particular, we specify the domains
of the operators which are implicit Proposition~\ref{pr.maxaccK1K2}
and Proposition~\ref{pr.cuspK1K2}.
Then  we solve inhomogeneous
boundary value problems by using variational
arguments inspired by 
\cite{Lio}\cite{Bar}\cite{Luc}\cite{Car} in the case $\nu=1$\,.

\subsection{Some notations}
\label{se.notasep}

In this section, we assume the operators $(L_{\pm},D(L_{\pm}))$ to be maximal
accretive in the separable Hilbert space $\mathfrak{L}$ with $L_{\pm}^{*}=L_{\mp}$\,. For $I=\rz_{-}$ or $I=\rz$\,, we
shall consider the operator
$$
P_{\pm}^{L}=\pm
p.\partial_{q}+\frac{-\partial_{p}^{2}+p^{2}+1}{2}+L_{\pm}=\pm
p\partial_{q}+
\frac{1}{2}+\mathcal{O}+L_{\pm}
$$
defined with the proper domain, containing
$\mathcal{C}^{\infty}_{0}(I\times \rz;D(L_{\pm}))$ and specified below, in
$L^{2}(I\times\rz,dqdp;\mathfrak{L})=L^{2}(I\times \rz,dqdp)\otimes
\mathfrak{L}$\,. The formal adjoint of $P_{\pm}^{L}$ is $P_{\mp}^{L}$\,.
 Both operators $1+L_{\pm}$ are
isomorphisms
 in the diagram
$$
D(L_{\pm})\stackrel{1+L_{\pm}}{\to}\mathfrak{L}\stackrel{1+L_{\pm}}{\to}
  D(L_{\pm}^{*})'=D(L_{\mp})'\,.
$$
where $D(L_{\pm}^{*})'$  is the dual of $D(L_{\pm}^{*})$
and the density of $D(L_{\pm})$ 
provides
the embedding $\mathfrak{L}\subset D(L^{*}_{\pm})'$\,. 
In particular
$P_{\pm}^{L}$ defines a continuous operator from $L^{2}(I\times
\rz,dqdp;\mathfrak{L})$ to $\mathcal{D}'(I\times
\rz; D(L_{\mp})')$ for any open interval $I\subset\rz$\,.\\
The densely defined quadratic form $\Real\langle u\,,\, L_{\pm} u\rangle$
defines a self-adjoint operator (take the Friedrichs extension see
\cite{ReSi}-Theorem~X.23). This operator is (abusively) denoted by $\Real
L$\,. For any $u\in D(L_{+})\cap D(L_{-})$\,, $\Real
Lu$ equals $\frac{L_{+}+L_{-}}{2}u$\,.
The notation $\Real L$ is especially justified under the following
assumption.
\begin{assumption}
\label{hyp.L1}
The intersection $D(L_{\pm})\cap D(L_{-})$ is dense in $D(L_{\pm})$
endowed with its graph norm and
 $\Real L$ is essentially self-adjoint on $D(L_{+})\cap D(L_{-})$\,.
\end{assumption}
Under this assumption, $D(L_{+})$, $D(L_{-})$ and $D(L_{+})\cap D(L_{-})$ are  densely included
in $D((\Real L)^{1/2})$\,, which is the quadratic form domain of
$\Real L$\,, and we have the embeddings
$$
D(L_{\pm})\subset D((\Real L)^{1/2})\subset \mathfrak{L}\subset
  D((\Real L)^{1/2})'\subset D(L_{\mp})'\,.
$$
Contrary to the case when $(L_{+},D(L_{+}))$ is a non negative self-adjoint
operator, there is no reason to assume in general 
that $L_{+}$ is continuous from $D((\Real L)^{1/2})$
to $D((\Real L)^{1/2})'$\,. This would mean that $\Imag
L=\frac{L_{+}-L_{-}}{2i}$ is estimated  in terms of $\Real L$\,.\\
According to this multiplicity of spaces, the notations of
Section~\ref{se.model} will be adapted with various choices for the
Hilbert space $\mathfrak{f}$\,. This will be specified in every case.
For example $\mathcal{H}^{s}$\,, $s\in\rz$\,, will follow the definition
$\mathcal{H}^{s}=(\frac{1}{2}+\mathcal{O})^{-s/2}L^{2}(\rz,dp;\cz)$ of
Section~\ref{se.presentmodel} with $\mathfrak{f}=\cz$\,.\\
We shall use the following space
\begin{equation}
  \label{eq.defH1L}
\mathcal{H}^{1,L}=\left\{u\in L^{2}(\rz,dp;\mathfrak{L})\,,
  \|u\|_{\mathcal{H}^{1}}^{2}+\langle u\,,\,
  (\Real L)u\rangle_{\mathfrak{L}} < +\infty\right\}\,,
\end{equation}
endowed with its natural norm and scalar product. Its dual is denoted
by $\mathcal{H}^{-1, L}$\,. When $L$ is bounded in $\mathfrak{L}$\,,
it is nothing but
$\mathcal{H}^{1,0}=\mathcal{H}^{1}\otimes\mathfrak{L}$ and $\mathcal{H}^{-1,0}=\mathcal{H}^{-1}\otimes\mathfrak{L}$\,.
For an interval $I$ of $\rz$
\begin{multline*}
L^{2}(I,dq;\mathcal{H}^{1}\otimes D(L_{\pm}))
\subset
L^{2}(I,dq;\mathcal{H}^{1, L})\subset L^{2}(I\times \rz,dqdp;\mathfrak{L})
\\
\subset L^{2}(I,dq;\mathcal{H}^{-1, L})\subset
L^{2}(I,dq;\mathcal{H}^{-1}\otimes D(L_{\mp})')\,.
\end{multline*}
The odd and even part of elements of $L^{2}(\rz,|p|dp;\mathfrak{L})$
involves an involution $j$ acting on $\mathfrak{L}$\,.
 In the end, this
construction will be applied  with
$\mathfrak{L}=\mathfrak{L}_{1}\otimes \mathfrak{f}_{1}$ with $L=L_{1}\otimes
\Id$ and $j=\Id\otimes j_{1}$\,.
\begin{definition}
  \label{de.invoL}
The Hilbert space $\mathfrak{L}$ is endowed with a unitary involution
$j$ which commutes with $e^{-tL_{\pm}}$ for all $t\geq 0$\,.
In $L^{2}(\rz,|p|dp; \mathfrak{L})$ or
$\mathcal{D}'(\rz^{*}; \mathfrak{f})$ with $\mathfrak{f}\in
\left\{D(L_{\pm}),\mathfrak{L}, D(L_{\pm})'\right\}$\,,
 the even and odd part are given by
\begin{eqnarray}
 \label{eq.gev}
&&\gamma_{ev}(p)=[\Pi_{ev}\gamma] (p)=\frac{\gamma(p)+
  j\gamma(-p)}{2}\,,\\
\label{eq.godd}&&
\gamma_{odd}(p)=[\Pi_{odd}\gamma](p)=
\frac{\gamma(p)-
  j\gamma(-p)}{2}\,.
\end{eqnarray}
The operators $\Pi_{+}$ and $\Pi_{-}$ are given by
\begin{equation}
  \label{eq.defPi+-}
\Pi_{+}=\Pi_{ev}+\sign(p)\Pi_{odd}\quad,\quad
\Pi_{-}=\Pi_{ev}-\sign(p)\Pi_{odd}\,,
\end{equation}
\end{definition}
The projections $\Pi_{ev}$ and $\Pi_{odd}$ are orthogonal in
$L^{2}(\rz,|p|dp;\mathfrak{L})$ with
\begin{eqnarray}
\nonumber
 &&\Pi_{ev}^{*}=\Pi_{ev}=(1-\Pi_{odd})\;,\\
&&
\label{eq.piose}
  \sign(p)\circ\Pi_{ev}=\Pi_{odd}\circ\sign(v)\,,\\
&&
\nonumber
\int_{\rz}\langle\gamma(p)\,,\,\gamma'(p)\rangle_{\mathfrak{L}}~pdp
=
\langle \gamma\,,\,
\sign(p)\gamma'\rangle_{L^{2}(\rz,|p|dp;\mathfrak{L})}
\\
\label{eq.intoddev}
&&
\quad=\langle \gamma_{ev}\,,\,
\sign(p)\gamma_{odd}'\rangle_{L^{2}(\rz,|p|dp;\mathfrak{L})}
+\langle \gamma_{odd}\,,\,
\sign(p)\gamma_{ev}'\rangle_{L^{2}(\rz,|p|dp;\mathfrak{L})}\,.
\end{eqnarray}
Moreover the two operators $\Pi_{+}$ and $\Pi_{-}$
are projections owing to $\sign(p)\Pi_{ev}=\Pi_{odd}\sign(p)$\,, with
the same range $\Ran \Pi_{+}=\Ran \Pi_{-}=\Ran\Pi_{ev}$\,.
With $\Pi_{ev}=\frac{\Pi_{+}+\Pi_{-}}{2}$ and
$\sign(p)\Pi_{odd}=\frac{\Pi_{+}-\Pi_{-}}{2}$\,,
the quantity
\eqref{eq.intoddev} also equals
\begin{equation}
  \label{eq.int+-}
\frac{1}{2}\left[\langle \Pi_{+}\gamma\,,\,
  \Pi_{+}\gamma'\rangle_{L^{2}(\rz,|p|dp;\mathfrak{L})}
-
\langle \Pi_{-}\gamma\,,\,
 \Pi_{-} \gamma'\rangle_{L^{2}(\rz,|p|dp;\mathfrak{L})}
\right]\,.
\end{equation}
The projection $\Pi_{+}$ (resp. $\Pi_{-}$) has the same kernel as $1_{\rz_{+}}(p)$
 (resp. $1_{\rz_{-}}(p)$)
according to the explicit formulas:
$$
\Pi_{+}\gamma=(1_{\rz_{+}}(p)+j1_{\rz_{-}}(p)) \gamma (|p|)\quad,\quad
\Pi_{-}\gamma=(j1_{\rz_{+}}(p)+1_{\rz_{-}}(p)) \gamma(-|p|)\,.
$$
The mapping $\Pi_{+}\times\Pi_{-}: L^{2}(\rz,|p|dp;\mathfrak{L})\to
(\Ran \Pi_{ev})^{2}$ is an isomorphism
\begin{eqnarray}
  \label{eq.isom}
&& (\Pi_{+},\Pi_{-}):L^{2}(\rz,|p|dp;\mathfrak{L})\leftrightarrows
(\Ran \Pi_{ev})^{2}\\
\nonumber
\text{with}&& 
L^{2}(\rz,|p|dp;\mathfrak{L})=\ker(\Pi_{+})\mathop{\oplus}^{\perp} \ker(\Pi_{-})\,,
\\
\nonumber
\text{and}
&&  (\Pi_{+},\Pi_{-})\gamma_{ev}=(\gamma_{ev},\gamma_{ev})\,,\\
\nonumber
&&(\Pi_{+},\Pi_{-})\gamma_{odd}=(\sign(p)\gamma_{odd},-\sign(p)\gamma_{odd})\,.
\end{eqnarray}
\subsection{Traces and integration by parts}
\label{se.traces}
For a given interval, we shall consider the space
$$
\mathcal{E}_{I}^{L,\pm}=\left\{u\in L^{2}(I,dq;\mathcal{H}^{1, L})\,,\,
    P^{L}_{\pm}u\in L^{2}(I,dq;\mathcal{H}^{-1, L})\right\}\,,
$$
endowed with the norm
$$
\|u\|_{\mathcal{E}_{I}^{L,\pm}}=\|u\|_{L^{2}(I,dq;\mathcal{H}^{1, L})}
+\|P_{\pm}^{L}u\|_{L^{2}(I,dq;\mathcal{H}^{-1, L})}\,.
$$
\begin{proposition}
\label{pr.EIDL}
Under Hypothesis~\ref{hyp.L1} the space
 $\mathcal{E}_{I}^{L,\pm}$ is continuously
  embedded in the space 
  $\mathcal{E}_{I}(D(L_{\mp})')$ of Definition~\ref{de.EI} with
  $\mathfrak{f}=D(L_{\mp})'$\,,
and therefore in
$\mathcal{C}^{0}_{b}(\overline{I};L^{2}(\rz,|p|dp; D(L_{\mp})')$\,.
\end{proposition}
\begin{proof}
From the embeddings 
$\mathcal{H}^{1, L}\subset
\mathcal{H}^{1, 0}=\mathcal{H}^{1}\otimes
\mathfrak{L}\subset\mathcal{H}^{1}\otimes D(L_{\mp})'$ and
$\mathcal{H}^{-1, L}\subset \mathcal{H}^{-1}\otimes D(L_{\mp})'$\,, we
deduce
\begin{eqnarray*}
&&u\in L^{2}(I,dq;\mathcal{H}^{1}\otimes\mathfrak{L})\subset
L^{2}(I,dq;\mathcal{H}^{1}\otimes D(L_{\mp})')\,,\\
\text{and}&&
p\partial_{q}u=\pm\left[P^{L}_{\pm}u-(\frac{1}{2}+\mathcal{O})u-L_{\pm}u\right]\in
L^{2}(I,dq;\mathcal{\mathcal{H}}^{-1}\otimes D(L_{\mp})')\,,
\end{eqnarray*}
because $u\in L^{2}(I,dq;\mathcal{H}^{1}\otimes\mathfrak{L})$ implies
$L_{\pm}u\in  L^{2}(I,dq;\mathcal{H}^{1}\otimes D(L_{\mp})')$ and
$(\frac{1}{2}+\mathcal{O})u\in L^{2}(I,dq; \mathcal{H}^{-1}\otimes \mathfrak{L})$\,.
We conclude by referring to Proposition~\ref{pr.intbypart} with
$\mathfrak{f}$ replaced by $D(L_{\mp})'$\,.
\end{proof}
\begin{proposition}
\label{pr.EIDLint}
Assume Hypothesis~\ref{hyp.L1} and for $I=\rz_{-}$ set $\gamma
u=u(q=0)$ for $u\in \mathcal{E}_{\rz_{-}}^{L,\pm}$\,.
The integration by part formula
\begin{eqnarray*}
&& \langle v\,,\,P^{L}_{+}u\rangle
=\lim_{\varepsilon\to 0^{+}}
  \langle \gamma v_{\varepsilon}\,, \sign(p)\gamma u_{\varepsilon}\rangle_{L^{2}(\rz,|p|dp;\mathfrak{L})}
+\langle P^{L}_{-}v\,,\, u\rangle\,,\\
&& v_{\varepsilon}=(1+\varepsilon L_{-})^{-1}v\quad,\quad
u_{\varepsilon}=(1+\varepsilon L_{+})^{-1}u\,,
\end{eqnarray*}
holds for pairs $(u,v)\in \mathcal{E}^{L,+}_{\rz_{-}}\times
\mathcal{C}^{\infty}_{0}((-\infty,0]\times\rz;D(L_{-}))$ 
and for pairs $(u,v)$ which satisfy
\begin{itemize}
\item $u$ and $v$ belong to $L^{2}(\rz_{-},dq,\mathcal{H}^{1,L})$\,;
\item $P^{L}_{+}u$ and $P^{L}_{-}v$ belong to
  $L^{2}(\rz_{-},dq;\mathcal{H}^{-1}\otimes\mathfrak{L})$\,.
\end{itemize}
\end{proposition}
\begin{proof}
In both cases, Proposition~\ref{pr.EIDL} ensures $u\in
\mathcal{E}_{\rz_{-}}(D(L_{-})')$ and $v\in
\mathcal{E}_{\rz_{-}}(D(L_{+})')$\,.\\
When $(u,v)\in \mathcal{E}_{\rz_{-}}^{L,+}\times
\mathcal{C}^{\infty}_{0}((-\infty,0]\times \rz;D(L_{-})))$\,, the integration by part
\begin{multline*}
\langle v\,,\, p\partial_{q}u+(\frac{1}{2}+\mathcal{O})u\rangle=\int_{\rz}\langle
v(0,p)\;,\; u(0,p)\rangle_{D(L_{-})\,,\, D(L_{-})'}~pdp
\\
+\langle (-p\partial_{q}v+\frac{1}{2}+\mathcal{O})v\,,\,u\rangle
\end{multline*}
makes sense, while the other term of $P^{L}_{+}u$ gives the integrable quantities
\begin{equation*}
\langle v(q,p)\,,\, L_{+}u(q,p)\rangle_{D(L_{-})\,,\, D(L_{-})'}
=
\langle L_{-}v(q,p)\,,\, u(q,p)\rangle\,.
\end{equation*}
In the second case $u_{\varepsilon}=(1+\varepsilon L_{+})^{-1}u$  and
$v_{\varepsilon}=(1+\varepsilon L_{-})^{-1}v$ both belong to
$\mathcal{E}_{\rz_{-}}(\mathfrak{L})$
and we can apply the polarized integration by part of
Proposition~\ref{pr.intbypart}
$$
\langle p\partial_{q}u_{\varepsilon}\,,\, v_{\varepsilon}\rangle+
\langle u_{\varepsilon}\,,\, p\partial_{q}v_{\varepsilon}\rangle
=\int_{\rz}\langle u_{\varepsilon}(0,p)\,,\,
v_{\varepsilon}(0,p)\rangle_{\mathfrak{L}}~pdp\,.
$$
The term with $(\frac{1}{2}+\mathcal{O})$ involves only the duality between
$\mathcal{H}^{1}$ and $\mathcal{H}^{-1}$ while the terms with $L_{+}$ and
$L_{-}$ is concerned with $L^{2}(\rz^{2}_{-},dqdp;\mathfrak{L})$ functions.
We obtain
\begin{multline*}
\langle (1+\varepsilon L_{-})^{-1}v\,,\,(1+\varepsilon
L_{+})^{-1}P^{L}_{+}u\rangle
-
\langle (1+\varepsilon L_{-})^{-1}P^{L}_{-}v\,,\,
(1+\varepsilon L_{+})^{-1}u\rangle\\
=\langle (1+\varepsilon L_{-})^{-1}\gamma v\,,\,
\sign(p)(1+\varepsilon
L_{+})^{-1}\gamma u\rangle_{L^{2}(\rz,|p|dp;\mathfrak{L})}
\end{multline*}
Our assumptions were made so that  the left-hand side, and
therefore the right-hand side, converge. In
particular $P^{L}_{+}u\,,\, P^{L}_{-}v \in
L^{2}(\rz_{-},dq;\mathcal{H}^{-1}\otimes 
\mathfrak{L})$ is used with $\mathop{\mathrm{s-lim}}_{\varepsilon\to
  0}(1+\varepsilon L_{\pm})^{-1}=1$ in $\mathfrak{L}$ and $\|(1+\varepsilon
L_{\pm})^{-1}\|_{\mathcal{L}(\mathfrak{L})}\leq 1$\,.
\end{proof}
\begin{remark}
 The last argument cannot be used for general $(u,v)\in
 \mathcal{E}_{\rz_{-}}^{L,+}\times \mathcal{E}_{\rz_{-}}^{L,-}$
 because we do 
  not control the effect of the regularization $(1+\varepsilon
  L_{\pm})^{-1}$ on $D((\Real L)^{1/2})$ and $\mathcal{H}^{\pm 1, L}$\,.
\end{remark}
\begin{proposition}
\label{pr.EIL2}
Assume Hypothesis~\ref{hyp.L1}. Take $I=(a,b)$ with $-\infty\leq a< b\leq +\infty$\,.
Assume $u\in \mathcal{E}_{I}^{L,\pm}$ with $P^{L}_{\pm}u\in
L^{2}(I,dq;\mathcal{H}^{-1}\otimes \mathfrak{L})$ and $\gamma_{a}u=u(a,.)\in
L^{2}(\rz,|p|dp; \mathfrak{L})$ if $-\infty<a$ and $\gamma_{b}u=u(b,.)\in
L^{2}(\rz,|p|dp;\mathfrak{L})$ if $b<+\infty$\,.
Then the following inequality holds
\begin{multline*}
2 \|u\|_{L^{2}(I,dq;\mathcal{H}^{1,L})}^{2}\leq 2 \Real\langle
u\,,\,P^{L}_{\pm}u\rangle
\pm 1_{\rz}(a)\int_{\rz}|u(a,p)|_{\mathfrak{L}}^{2}~pdp
\\
\mp 1_{\rz}(b)\int_{\rz}|u(b,p)|_{\mathfrak{L}}^{2}~pdp\,.
\end{multline*}
\end{proposition}
\begin{proof}
With an obvious change of signs, it
 suffices to consider the $+$ case\,.
  From Proposition~\ref{pr.EIDL} we know $u\in
  \mathcal{E}_{I}(D(L_{-})')$ and the traces are well defined. Set
  $f=P_{+}^{L}u$\,, $f_{\partial}^{a}=\gamma_{a}u$ and $f_{\partial}^{b}=\gamma_{b}u$\,.
For $\varepsilon>0$ take $u_{\varepsilon}=(1+\varepsilon L_{+})^{-1}u$ and
set $f_{\varepsilon}=(1+\varepsilon L_{+})^{-1}f$\,,
$f_{\partial,\varepsilon}^{a,b}=(1+\varepsilon L_{+})^{-1}f_{\partial}^{a,b}$\,.
Then $u_{\varepsilon}$ belongs to $L^{2}(I,dq;\mathcal{H}^{1}\otimes
D(L_{+}))\subset L^{2}(I,dq;\mathcal{H}^{1}\otimes\mathfrak{L})$ and satisfies
\begin{eqnarray*}
  &&
  p\partial_{q}u_{\varepsilon}=f_{\varepsilon}-(\frac{1}{2}+\mathcal{O})u_{\varepsilon}-L_{+}u_{\varepsilon}
\quad\in
 L^{2}(I,dq;\mathcal{H}^{-1}\otimes\mathfrak{L})\,,\\
&&
\gamma_{a,b}u_{\varepsilon}=f_{\partial,\varepsilon}^{a,b}\in L^{2}(\rz,|p|dp;\mathfrak{L})\,.
\end{eqnarray*}
Hence the integration by parts formula of
Proposition~\ref{pr.intbypart} can be used with
$\mathfrak{f}=\mathfrak{L}$:
\begin{eqnarray*}
1_{\rz}(b)\int_{\rz}|f_{\partial,\varepsilon}^{b}(p)|_{\mathfrak{L}}^{2}~pdp
-1_{\rz}(a)\int_{\rz}|f_{\partial,\varepsilon}^{a}(p)|_{\mathfrak{L}}^{2}~pdp
&&=2\Real\langle u_{\varepsilon}\,,\,
p\partial_{q}u_{\varepsilon}\rangle\\
&&\hspace{-2cm}=2\Real\langle u_{\varepsilon}\,, f_{\varepsilon}\rangle
-2\|u_{\varepsilon}\|_{L^{2}(I,dq;\mathcal{H}^{1,L})}^{2}\,.
\end{eqnarray*}
With $\|(1+\varepsilon L_{+})^{-1}\|_{\mathcal{L}(\mathfrak{L})}\leq 1$
and $\mathcal{H}^{1,0}=\mathcal{H}^{1}\otimes \mathfrak{L}$\,,
we deduce from this the inequalities
\begin{eqnarray*}
\|u_{\varepsilon}\|_{L^{2}(I,dq;\mathcal{H}^{1,0})}^{2}+\|u_{\varepsilon}\|_{L^{2}(I,dq;\mathcal{H}^{1,L})}^{2}
&\leq&
2\|u_{\varepsilon}\|_{L^{2}(I,dq;\mathcal{H}^{1,L})}^{2}
\\
&&\hspace{-3cm}\leq
\|u_{\varepsilon}\|_{L^{2}(I,dq;\mathcal{H}^{1,0})}\|f\|_{L^{2}(I,dq;\mathcal{H}^{-1,0})}
+1_{\rz}(b)\|f_{\partial}^{b}\|_{L^{2}(\rz,|p|dp;\mathfrak{L})}^{2}
\\
&&\hspace{3cm}
+1_{\rz}(a)\|f_{\partial}^{a}\|_{L^{2}(\rz,|p|dp;\mathfrak{L})}^{2}\,,
\end{eqnarray*}
and a uniform bound for
$\|u_{\varepsilon}\|_{L^{2}(I,dq;\mathcal{H}^{1,L})}$\,.\\
As $\varepsilon\to 0$ the operators $(1+\varepsilon L_{+})^{-1}$\,,
(resp. $\Id_{\mathfrak{g}}\otimes (1+\varepsilon L_{+})^{-1}$) converges
 in the strong operator topology to $\Id_{\mathfrak{L}}$
 (resp. $\Id_{\mathfrak{g}\otimes \mathfrak{L}}$).
Therefore we know
\begin{eqnarray*}
  &&\lim_{\varepsilon\to
    0}\|f_{\partial,\varepsilon}^{a,b}-f_{\partial}^{a,b}\|_{L^{2}(\rz,|p|dp;\mathfrak{L})}=0\,,\\
&& \lim_{\varepsilon\to
  0}\|u_{\varepsilon}-u\|_{L^{2}(I,dq;\mathcal{H}^{1,0})}=0\,,\\
&& \lim_{\varepsilon\to 0}\|f_{\varepsilon}-f\|_{L^{2}(I,dq;\mathcal{H}^{-1,0})}=0\,,
\end{eqnarray*}
while $u_{\varepsilon}$ converges weakly to some $v$ in
$L^{2}(I,dq;\mathcal{H}^{1,L})$\,.
Since $\mathcal{H}^{-1,0}=(\mathcal{H}^{1,0})'$ is continuously embedded in
$\mathcal{H}^{-1,L}=(\mathcal{H}^{1,L})'$\,, the second limit implies $v=u$ and
$$
\|u\|_{L^{2}(I,dq;\mathcal{H}^{1,L})}^{2}\leq \liminf_{\varepsilon\to
  0}\|u_{\varepsilon}\|_{L^{2}(I,dq;\mathcal{H}^{1,L})}^{2}\,.
$$
This ends the proof.
\end{proof}
\begin{remark}
  Again the strong limit argument $\lim_{\varepsilon\to
    0}(1+\varepsilon L_{+})^{-1}f=f$ works in $\mathfrak{L}$ but not in
  $D((\Real L)^{1/2})$\,, $\mathcal{H}^{1,L}$ or $\mathcal{H}^{-1,L}$\,.
\end{remark}
\subsection{Identifying the domains}
\label{se.fullsep}
We shall first consider the whole space problem on $\rz^{2}$ and then
the
half-space problem in $\rz_{-}\times \rz$\,.\\
In the whole space problem, $(K^{L}_{\pm},D(K^{L}_{\pm}))$ will denote the
maximal accretive realization defined in
$L^{2}(\rz^{2},dqdp;
\mathfrak{L})=
L^{2}(\rz^{2},dqdp)\otimes\mathfrak{L}$ like in
Subsection~\ref{se.tensor},  with
$\mathcal{H}_{1}=L^{2}(\rz^{2},dqdp)$\,, $\mathcal{H}_{2}=\mathfrak{L}$\,,
$K_{1}=K_{\pm}+\frac{1}{2}=\pm p\partial_{q}+(\frac{1}{2}+\mathcal{O})$ and
$K_{2}=L_{\pm}$ (the domains $D(K_{\pm})$ and $D(L_{\pm})$
are known).\\
For the half-space problem only the boundary conditions
$$
\gamma_{odd}u=\pm\sign(p)\nu\gamma_{ev}u\quad\text{with}\quad\nu\in \left\{0,1\right\}
$$
are considered. The corresponding maximal accretive realization of
$\pm p.\partial_{q}+\frac{1}{2}+\mathcal{O}$ in $L^{2}(\rz^{2}_{-},dqdp)$ defined in
Subsection~\ref{se.maxacc1D} is denoted by $(K_{\pm, \nu}+\frac{1}{2},D(K_{\pm, \nu}))$\,.
The maximal accretive operator $(K_{\pm, \nu}^{L},D(K_{\pm, \nu}^{L}))$  is  the
maximal accretive realization defined in
$L^{2}(\rz^{2}_{-},dqdp; \mathfrak{L})=
L^{2}(\rz^{2}_{-},dqdp)\otimes\mathfrak{L}$ like in
Subsection~\ref{se.tensor},  with
$\mathcal{H}_{1}=L^{2}(\rz^{2}_{-},dqdp)$\,, $\mathcal{H}_{2}=\mathfrak{L}$\,,
$K_{1}=K_{\pm,\nu}+\frac{1}{2}$ and $K_{2}=L_{\pm}$\,.
\begin{proposition}
\label{pr.idenDK}
With Hypothesis~\ref{hyp.L1}, the domain $D(K^{L}_{\pm})$ equals
$$
\left\{u\in L^{2}(\rz,dq;\mathcal{H}^{1,L})\,,\,
  P^{L}_{\pm}u\in L^{2}(\rz^{2},dqdp; \mathfrak{L})\right\}\,.
$$
The relation $K^{L}_{\pm}u=P^{L}_{\pm}u$ and the integration by part equality
$$
\|u\|_{L^{2}(\rz,dq;\mathcal{H}^{1,L})}^{2}=\Real\langle u\,,\,K^{L}_{\pm}u\rangle
$$
hold for any $u\in D(K^{L}_{\pm})$\,.\\
The adjoint of $K_{\pm}^{L}$ is $K_{\mp}^{L}$\,.
\end{proposition}
\begin{proof}
From Theorem~\ref{th.app}\,, we know that
$K_{\pm}^{*}=K_{\mp}$\,.\\
The maximal accretive operator $K_{\pm}^{L}$ being defined as the
generator of \\
$(e^{-t(K_{\pm}+\frac{1}{2})}e^{-tL_{\pm}})_{t\geq 0}$\,, let us first
check that 
$K_{\pm}^{L}u=P_{\pm}^{L}u$ when $u\in
D(K_{\pm}^{L})$\,. For such a $u\in D(K_{\pm}^{L})$\,, 
 $e^{-tK_{\pm}^{L}}u=e^{-t(K_{\pm}+\frac{1}{2})}e^{-tL_{\pm}}u$  defines an
 $L^{2}(\rz^{2}, dqdp; \mathfrak{L})$-valued 
$\mathcal{C}^{1}$ function on $[0,+\infty)$\,.
From the equality
$$
(e^{-tK_{\pm}}e^{-tL_{\pm}})^{*}=e^{-tK_{\pm}^{*}}e^{-tL_{\pm}^{*}}=e^{-tK_{\mp}}e^{-tL_{\mp}}\,.
$$
we deduce that $(K_{\pm}^{L})^{*}=K_{\mp}^{L}$ and we know 
that $\mathcal{C}^{\infty}_{0}(\rz^{2})\otimes^{alg} D(L_{-})$\,,
and therefore $\mathcal{C}^{\infty}_{0}(\rz^{2};D(L_{-}))$, is contained in
$D(K_{\mp}^{L})=D((K_{\pm}^{L})^{*})$\,. Looking at the derivative at $t=0$
of
$$
\langle e^{-t(K_{\mp}+\frac{1}{2})}e^{-tL_{\mp}}v\,, u\rangle=\langle v\,,\,
e^{-tK_{\pm}^{L}}u\rangle\,,
$$
we obtain
$$
\langle ((K_{\mp}+\frac{1}{2})\otimes \Id +\Id\otimes L_{\mp})v\,,\,
u\rangle=\langle v\,,\, K_{\pm}^{L}u\rangle\,,
$$
for all $v\in\mathcal{C}^{\infty}_{0}(\rz^{2};D(L_{\mp}))$\,.
But this means exactly
$$
K_{\pm}^{L}u=(\frac{1}{2}+K_{\pm})u+L_{\pm}u=P_{\pm}^{L}u\quad\text{in}~\mathcal{D}'(\rz^{2};D(L_{\mp})')\supset L^{2}(\rz^{2},dqdp;\mathfrak{L})\,.
$$
Proposition~\ref{pr.EIL2} applied with $a=-\infty$\,, $b=+\infty$ and
 provides the inequality
$$
\|u\|_{L^{2}(\rz,dq;\mathcal{H}^{1,L})}^{2}\leq
\Real\langle u\,,\, K^{L}_{\pm}u\rangle\quad,\quad \forall u\in D(K^{L}_{\pm})\,.
$$
With $\|u\|_{L^{2}(\rz,dq;\mathcal{H}^{1,L})}\geq
\|u\|$\,, this implies $\|K_{\pm}^{L}u\|\geq
\|u\|_{L^{2}(\rz,dq;\mathcal{H}^{1,L})}$ and
 $D(K_{\pm}^{L})$ is continuously embedded in
$L^{2}(\rz,dq; \mathcal{H}^{1,L})$\,.\\
Conversely, if $v$ belongs to $L^{2}(\rz,dq;\mathcal{H}^{1,L})$ with $P_{\pm}^{L}v=f\in
L^{2}(\rz^{2},dqdp; \mathfrak{L})$\,, there exists $u\in
D(K_{\pm}^{L})$ such that $K_{\pm}^{L}u=f$ because $(K_{\pm}^{L}-1,D(K_{\pm}^{L}))$ is maximal
accretive.
The difference $w=u-v$ belongs to $L^{2}(\rz,dq;\mathcal{H}^{1,L})$
with $P^{L}w=0$\,. Proposition~\ref{pr.EIL2} then implies
$$
\|w\|_{L^{2}(\rz,dq;\mathcal{H}^{1,L})}^{2}\leq 0\,.
$$
Therefore $v=u$ belongs to $D(K_{\pm}^{L})$\,.\\
The equality
$$
\Real\langle u\,,\, K_{\pm}^{L}u\rangle=\Real\langle u\,,\, (\frac{1}{2}+K_{\pm})u\rangle + 
\Real\langle u\,,\, L_{\pm} u\rangle=\|u\|_{L^{2}(\rz,dq;\mathcal{H}^{1,L})}^{2}\,,
$$
holds  when $u\in D(K_{\pm})\otimes^{alg}D(L_{\pm})$\,. The algebraic tensor product
$D(K_{\pm})\otimes^{alg}D(L_{\pm})$ is dense in $D(K_{\pm}^{L})$ according to
Proposition~\ref{pr.maxaccK1K2} while both sides are continuous on
$D(K_{\pm}^{L})$\,. This proves the equality for all $u\in
D(K_{\pm}^{L})$\,.
\end{proof}
\begin{proposition}
\label{pr.idenDKAId}
Under Hypotheses~\ref{hyp.L1} 
and with $\nu\in \left\{0,1\right\}$\,, the domain $D(K^{L}_{\pm,\nu})$
is nothing but
$$
\left\{u\in
  L^{2}(\rz_{-},dq;\mathcal{H}^{1,L})\,,\,
  \begin{array}[c]{l}
 P^{L}_{\pm}u\in
  L^{2}(\rz^{2}_{-},dqdp;\mathfrak{L})\\
\gamma_{odd}u=\pm\sign(p) \nu\gamma_{ev}u
\end{array}
\right\}\,.
$$
Moreover the relation
$$
K^{L}_{\pm,\nu}u=P_{\pm}^{L}u\,,
$$
and
 the integration by parts identity
$$
\|u\|_{L^{2}(\rz_{-},dq;\mathcal{H}^{1,L})}^{2}+
\|\gamma_{odd}u\|_{L^{2}(\rz,|p|dp; \mathfrak{L})}^{2}
=
\Real\langle u\,,\, K^{L}_{\pm,\nu}u\rangle
$$
holds for all $u\in D(K^{L}_{\pm,\nu})$\,.\\
Finally the adjoint of $K^{L}_{\pm,\nu}$ is $K^{L}_{\mp,\nu}$\,.
\end{proposition}
\begin{remark}
The boundary condition $\gamma_{odd}u=\pm\sign(p)\nu\gamma_{ev}u$ makes
sense because the trace $\gamma u$ is well-defined in
$L^{2}(\rz,|p|dp; D(L_{\mp})')$ when $u\in
\mathcal{E}^{L,\pm}_{\rz_{-}}$ and a fortiori when $u\in
L^{2}(\rz_{-},dq;\mathcal{H}^{1,L})$ and $P_{\pm}^{L}u\in
L^{2}(\rz^{2}_{-},dqdp;\mathfrak{L})$\,. 
For $\nu=0$\,, the integration by parts inequality 
says nothing about $\gamma_{ev}u$\,. This will be studied later in 
Section~\ref{se.genBC}.
\end{remark}
\begin{proof}
We proved in Proposition~\ref{pr.adj1} that the adjoint of
$K_{\pm,\nu}$ is $K_{\mp,\nu}$ for $\nu=\left\{0,1\right\}$\,.\\
The proof of $(K_{\pm,\nu}^{L})^{*}=K_{\mp,\nu}^{L}$ and $K^{L}_{\pm,\nu}u=P_{\pm}^{L}u$ when $u\in
D(K^{L}_{\pm,\nu})$ is the same as in
Proposition~\ref{pr.idenDK}\,, with
$\mathcal{C}^{\infty}_{0}(\rz^{2}_{-};D(L_{\mp}))\subset D((K_{\pm,\nu}^{L})^{*})$\,.\\
For $u\in D(K_{\pm,\nu})\otimes^{alg}D(L_{\pm})$ the integration by part formula
$$
\|u\|_{L^{2}(\rz_{-},dq;\mathcal{H}^{1}\otimes\mathfrak{L})}^{2}+\Real\langle\gamma_{ev}u\,,\,
\nu\gamma_{ev}u\rangle_{L^{2}(\rz,|p|dp;\mathfrak{L})}=\Real\langle
u\,,\, K^{L}_{\pm,\nu}u\rangle-\Real\langle u\,,\, L_{\pm} u\rangle
$$
comes from Theorem~\ref{th.maxaccKA} or
Proposition~\ref{pr.intbypart}, applied with
$\mathfrak{f}$ replaced by $\mathfrak{L}$  and 
$$
(\frac{1}{2}+K_{\pm,\nu})u=K^{L}_{\pm, \nu}u-
L_{\pm}u\quad\text{in}~L^{2}(\rz^{2}_{-},dqdp;\mathfrak{L})\,.
$$
We know also, for such a $u$\,, $\gamma u \in
L^{2}(\rz,|p|dp; \mathfrak{L})$ and $\gamma_{odd}u
=\pm\sign(p)\nu\gamma_{ev}u$\,.
With $\|u\|_{L^{2}(\rz_{-},dq;\mathcal{H}^{1}\otimes\mathfrak{L})}\geq
\|u\|$, we deduce
$$
\forall u\in D(K_{\pm,\nu})\otimes^{alg}D(L_{\pm})\,,\quad
\|u\|_{L^{2}(\rz_{-},dq;\mathcal{H}^{1,L})}+\|\gamma_{odd}
u\|_{L^{2}(\rz,|p|dp;\mathfrak{L})}\leq 8\|K^{L}_{\pm,\nu}u\|\,.
$$
Since $D(K_{\pm,\nu})\otimes^{alg}D(L_{\pm})$ is dense in $D(K^{L}_{\pm,\nu})$\,, 
we conclude that for all $u\in D(K^{L}_{\pm,\nu})$\,,
$\gamma_{odd} u=\pm\sign(p)\nu\gamma_{ev}u$ belongs to
$L^{2}(\rz,|p|dp;\mathfrak{L})$
and $u$ belongs
to $L^{2}(\rz_{-},dq;\mathcal{H}^{1,L})$\,.\\
With the same density argument, the equality
$$
\|u\|_{L^{2}(\rz_{-},dq;\mathcal{H}^{1,L})}^{2}+\|\gamma_{odd}u\|^{2}_{L^{2}(\rz,|p|dp;\mathfrak{L})}
=\langle u\,,\, K^{L}_{\pm,\nu}u\rangle
$$
can be extended
to any $u\in D(K^{L}_{\pm,\nu})$\,.\\
Let us complete the identification of $D(K^{L}_{\pm,\nu})$\,.
Since $K^{L}_{\pm,\nu}-1$ is maximal accretive the equation
$K^{L}_{\pm,\nu}u=f$ admits a unique solution $u\in
D(K^{L}_{\pm,\nu})$ for any $f\in
L^{2}(\rz^{2}_{-},dqdp;\mathfrak{L})$\,. 
If $v$ belongs to
$L^{2}(\rz_{-},dq;\mathcal{H}^{1,L})$ and satisfies
\begin{eqnarray*}
  && P_{\pm}^{L}v=f\, \in L^{2}(\rz^{2}_{-},dqdp; \mathfrak{L})\,,\\
&& \gamma_{odd}v=\pm\sign(p)\nu\gamma_{ev}v\quad(\text{in}~L^{2}(\rz,|p|dp;D(L_{\mp})')\,,
\end{eqnarray*}
solve $K^{L}_{\pm,\nu}u=f$\,. The difference $w=u-v$
belongs to $L^{2}(\rz_{-},dq;\mathcal{H}^{1,L})$ and satisfies
$$
P_{\pm}^{L}w=0\quad,\quad \gamma_{odd}w=\pm\sign(p)\nu\gamma_{ev}w\,. 
$$
This implies $w\in
\mathcal{E}_{\rz_{-}}^{L,\pm}\subset\mathcal{E}_{\rz_{-}}(D(L_{\mp})')$
according to Proposition~\ref{pr.EIDL}. The regularisation $w'=(1+L)^{-1}w$ belongs to
$\mathcal{E}_{\rz_{-}}(\mathfrak{L})$
and has a trace $\gamma w'\in L^{2}(\rz,|p|dp;\mathfrak{L})$\,.
It satisfies $P_{\pm}^{L}w'=0$\,, $w'\in
L^{2}(\rz_{-};\mathcal{H}^{1,L)})$ and $\gamma w'\in L^{2}(\rz,|p|dp;
\mathfrak{L})$ with
$\gamma_{odd}w'=\pm\sign(p)\nu\gamma_{ev}w'$\,.
The integration by part inequality of
Proposition~\ref{pr.EIL2} applies to $w'$:
$$
\|w'\|_{L^{2}(\rz_{-};\mathcal{H}^{1,L})}^{2}+\nu\langle \gamma_{ev}w'\,,
\gamma_{ev}w'\rangle\leq 0\,,
$$
which yields $w'=0$\,, $w=0$ and  $v=u\in D(K^{L}_{\pm,\nu})$\,.
\end{proof}
\begin{proposition}
  \label{pr.extH-1} 
Assume  Hypotheses~\ref{hyp.L1}
  and take $\lambda\in\rz$\,.
  \begin{itemize}
  \item The resolvent $(K^{L}_{\pm}-i\lambda)^{-1}$ admits a
    continuous extension which sends
    $L^{2}(\rz,dq;\mathcal{H}^{-1,L})$ into $L^{2}(\rz,dq;\mathcal{H}^{1,L})$ and the equation
$$
(P^{L}_{\pm}-i\lambda) u=f
$$
has a unique solution $u\in L^{2}(\rz,dq;\mathcal{H}^{1,L})$ when $f\in
L^{2}(\rz,dq;\mathcal{H}^{-1,L})$\,.
\item For $\nu\in \left\{0,1\right\}$\,, 
the resolvent $(K^{L}_{\pm, \nu}-i\lambda)^{-1}$ admits a continuous extension
which sends $L^{2}(\rz_{-},dq;\mathcal{H}^{-1,L})$ into 
$L^{2}(\rz_{-},dq;\mathcal{H}^{1,L})$ and the equation 
$$
(P^{L}_{\pm}-i\lambda)u=f\quad,\quad \gamma_{odd}u=\pm\sign(p)\nu\gamma_{ev}u\,,
$$
has a unique solution $u\in L^{2}(\rz_{-},dq;\mathcal{H}^{1,L})$ when $f\in
L^{2}(\rz_{-},dq;\mathcal{H}^{-1,L})$\,, with
$$
\|u\|_{L^{2}(\rz_{-},dq;\mathcal{H}^{1,L})}^{2}+\|\gamma_{odd}u\|_{L^{2}(\rz,|p|dp;
  \mathfrak{L})}^{2}=\Real\langle
f\,,\,u\rangle\,.
$$ 
\end{itemize}
\end{proposition}
\begin{proof}
  The proof is the same in the two cases and even simpler when there is no
  boundary term in the first case. We focus on the second case.\\
For $f\in L^{2}(\rz^{2}_{-},dq; \mathfrak{L})$ and
$u=(K^{L}_{\pm,\nu}-i\lambda)^{-1}f\in D(K^{L}_{\pm,\nu})$\,, the inequalities
\begin{eqnarray*}
  \|u\|_{L^{2}(\rz_{-},dq;\mathcal{H}^{1,L})}^{2}+\|\gamma_{odd}u\|_{L^{2}(\rz,|p|dp;
    \mathfrak{L})}^{2}
&=&
 \Real\langle u\,,\, (K^{L}_{\pm,\nu}-i\lambda)u\rangle\\
&=&\Real \langle u\,,\, f\rangle \\
&\leq& \|u\|_{L^{2}(\rz_{-},dq;\mathcal{H}^{1,L})}\|f\|_{L^{2}(\rz_{-},dq;\mathcal{H}^{-1,L})}\,,
\end{eqnarray*}
 allows to extend $(K_{\mp,\nu}^{L}-i\lambda)^{-1}f$ to $f\in
 L^{2}(\rz_{-},dq;\mathcal{H}^{-1,L})$\,. Approximating a general $f\in
 L^{2}(\rz_{-},dq;\mathcal{H}^{-1,L})$ by $f_{n}\in
 L^{2}(\rz^{2}_{-};\mathfrak{L})$ implies that the above equalities are still
 valid for
$f\in
 L^{2}(\rz_{-},dq;\mathcal{H}^{-1,L})$\,.\\
Assume now that  $v\in L^{2}(\rz_{-},dq;\mathcal{H}^{1,L})$ solves
$(P_{\pm}^{L}-i\lambda)v=f\in L^{2}(\rz_{-},dq;\mathcal{H}^{-1,L})$ and
$\gamma_{odd}v=\pm\sign(p)\nu\gamma_{ev}v$ in $L^{2}(\rz,|p|dp; D(L_{\mp})')$\,.
Then the difference $w=v-(K^{L}_{\pm,\nu}-i\lambda)^{-1}f$ belongs to
$L^{2}(\rz_{-};\mathcal{H}^{1,L})$ and satisfies  $P^{L}w=i\lambda
w\in L^{2}(\rz_{-}^{2},dqdp; \mathfrak{L})$\,,
$\gamma_{odd}w=\pm\sign(p)\nu\gamma_{ev}w$\,. It belongs to
$D(K^{L}_{\pm,\nu})$ and solves $(K^{L}_{\pm,\nu}-i\lambda)w=0$\,. This implies
$w=0$ and the solution is unique.
\end{proof}
\begin{definition}
\label{de.Ktext}
  We keep the notations $(K^{L}_{\pm}-i\lambda)^{-1}f$ when $f\in
  L^{2}(\rz,dq;\mathcal{H}^{-1,L})$ and $(K^{L}_{\pm,\nu}-i\lambda)^{-1}f$
  when 
$f\in L^{2}(\rz_{-},dq;\mathcal{H}^{-1,L})$\,.
\end{definition}

\subsection{Inhomogeneous boundary value problems}
\label{se.inhBVPL}
In this paragraph, we study 
an inhomogeneous boundary value problems, naturally associated with
$K_{+,1}^{L}$ and $K^{L}_{+,0}$\,. The control of $\gamma_{ev}u$
for $K^{L}_{\pm,1}$ allows  the variational arguments of
\cite{Luc}\cite{Car}\cite{Lio}. A general a priori estimate is deduced.
\begin{proposition}
\label{pr.inhBVPL1}
Assume Hypothesis~\ref{hyp.L1}, $f_{\partial}=\Pi_{\mp}L^{2}(\rz,|p|dp;
\mathfrak{L})=\Ran \Pi_{ev}$\,, 
$f\in L^{2}(\rz_{-},dq;\mathcal{H}^{-1,L})$ and $\lambda\in\rz$\,. Then the boundary value problem 
\begin{equation}
  \label{eq.inhBVPL}
(P^{L}_{\pm}-i\lambda)u=f\quad,\quad \Pi_{\mp}\gamma u=\gamma_{ev}u\mp\sign(p)\gamma_{odd}u=f_{\partial}\,,
\end{equation}
admits a unique solution in $L^{2}(\rz_{-},dq;\mathcal{H}^{1,L})$\,,
which is characterized by
\begin{equation}
  \label{eq.charac}
\forall \varphi \in D(K^{L}_{-,1})\,,\quad
\langle u\,,\,
(K^{L}_{-,1}+i\lambda)\varphi\rangle=
\langle u\,,\, f\rangle
+\frac{1}{2}\langle f_{\partial}\,,\,
\Pi_{\pm}\gamma\varphi\rangle_{L^{2}(\rz,|p dp|;\mathfrak{L})}\,.
\end{equation}
It satisfies the integration by part identity

\begin{equation}
  \label{eq.ippinhBVPL}
\frac{1}{4}\|\Pi_{\pm}\gamma u\|_{L^{2}(\rz,|p|dp;
  \mathfrak{L})}^{2}+
\|u\|_{L^{2}(\rz_{-},dq;\mathcal{H}^{1,L})}^{2}
=
\frac{1}{4}\|f_{\partial}\|_{L^{2}(\rz,|p|dp;
  \mathfrak{L})}^{2}
+
\Real\langle f\,,\,u\rangle
\,.
\end{equation}
\end{proposition}
\begin{proof}
In this proof it is important to distinguish $P^{L}_{+}$ and
$P^{L}_{-}$\,. While focusing on $P_{+}^{L}$ (the other case is
deduced by a transposition of $\left\{+,-\right\}$)\,,
we will work with $P_{+}^{L}$ and $P_{-}^{L}$\,.\\
  \textbf{Uniqueness:} If there are two solutions $u_{1},u_{2}\in
  L^{2}(\rz_{-},dq;\mathcal{H}^{1,L})$ then the difference  $w$ belongs to
$\mathcal{E}^{L,+}_{\rz_{-}}$ and solves $(P^{L}_{+}-i\lambda)w=0$
with
$\gamma_{odd}w=\sign(p)\gamma_{ev}w\in
L^{2}(\rz,|p|dp;D(L_{-})')$\,.  
This means exactly 
$w=(K^{L}_{+, 1}-i\lambda)^{-1}0=0$ and $u_{2}=u_{1}$\,.\\
\noindent\textbf{Existence when $f_{\partial}=0$:} Simply apply
Proposition~\ref{pr.extH-1} and $u=(K^{L}_{+,1}-i\lambda)^{-1}f$
according to Definition~\ref{de.Ktext}.
Because $(K_{+,1}^{L})=K_{-,1}^{L}$\,, the solution is characterized
by
$$
\forall \varphi\in D(K_{-,1}^{L})\,,\quad 
\langle u\,,\, (K_{-,1}^{L}+i\lambda)\varphi\rangle
=\langle f\,,\, \varphi\rangle\,.
$$\\
\noindent\textbf{Existence when $f=0$:}
 We adapt the approach of
\cite{Luc}\cite{Car}.
Let $V$ be the space $L^{2}(\rz_{-},dq;\mathcal{H}^{1,L})$ and let $Y$
be the domain $D(K^{L}_{-,1})$ in
$L^{2}(\rz^{2}_{-},dqdp;\mathfrak{L})$:
We are in the case $\nu=+1$ and
$Y$ is the set of $\varphi\in
  L^{2}(\rz_{-},dq;\mathcal{H}^{1,L})$ which satisfy
  \begin{eqnarray*}
    &&(P^{L}_{-}+i\lambda)\varphi\in L^{2}(\rz_{-}^{2},dqdp;
    \mathfrak{L})\,,\\
&& \gamma \varphi \in L^{2}(\rz,|p|dp;\mathfrak{L})\quad,\quad
\Pi_{+}\gamma\varphi=0\,\quad(\text{or}\quad
\gamma_{ev}\varphi+\sign(p)\gamma_{odd}\varphi=0)\,.
  \end{eqnarray*}
The space $Y$ is
endowed with the norm
$$
\|\varphi\|_{Y}=\left(\|\varphi\|_{L^{2}(\rz_{-},dq;\mathcal{H}^{1,L})}^{2}+
\|\gamma_{ev}\varphi\|_{L^{2}(\rz,|p|dp;\mathfrak{L})}^{2}\right)^{\frac 1 2}\,.
$$
The space $(Y,\|~\|_{Y})$ is continuously embedded in
$V=L^{2}(\rz_{-},dq;\mathcal{H}^{1,L})$ and contains
$
\left\{\varphi\in \mathcal{C}^{\infty}_{0}(\overline{\rz^{2}_{-}};
  D(L_{-}))\,, 1_{\rz_{+}}(p)\varphi(0,p)\equiv 0\right\}$\,, which is
dense in $V$\,.
Let $a:V\times Y\to \rz$ be the real bilinear form given by
$$
a(u,\varphi)=\Real
\langle u\,,\, (P^{L}_{-}+i\lambda)\varphi\rangle
\,,
$$
For any fixed $\varphi\in Y$ the real linear map
$V\ni u\mapsto a(u,\varphi)\in\rz$ is continuous.
For $\varphi\in Y=D(K^{L}_{-,1})$\,,
Proposition~\ref{pr.idenDKAId} implies that $a(\varphi,\varphi)$ equals
\begin{multline*}
a(\varphi,\varphi)
=\Real \langle \varphi\,,\, (K^{L}_{-,1}+i\lambda)\varphi\rangle
=\|\gamma_{ev}\varphi\|_{L^{2}(\rz,|p|dp; \mathfrak{L})}^{2}
+\|\varphi\|_{L^{2}(\rz_{-},dq;\mathcal{H}^{1,L})}^{2}
\end{multline*}
and $a$ is coercive on $(Y,\|~\|_{Y})$\,.
Finally the linear form $\ell: F\to \rz$ given by
$$
\ell(\varphi)=-\Real\langle f_{\partial}\,,\,\sign(p)
\gamma\varphi \rangle_{L^{2}(\rz,|p|dp; \mathfrak{L})}
=\frac{1}{2}\Real
\langle \Pi_{-}f_{\partial}\,,\, \Pi_{-}\gamma\varphi\rangle_{L^{2}(\rz,|p|dp;
  \mathfrak{L})}\,,
$$
where we use \eqref{eq.int+-} and $\Pi_{+}\gamma\varphi=0$\,,
 is
continuous.\\
 By Lions'Theorem (see
\cite{Lio}) there exists $u\in V$ such that 
$$
\forall \varphi\in Y\,,\quad a(u,\varphi)=
\frac{1}{2}\Real\langle \Pi_{-}f_{\partial}\,,\,
\Pi_{-}\gamma\varphi\rangle_{L^{2}(\rz,|p|dp; \mathfrak{L})}\,.
$$
By taking test functions $\varphi\in \mathcal{C}^{\infty}_{0}(\rz^{2}_{-};
D(L_{-}))\subset Y$ 
supported away from $\left\{q=0\right\}$\,, the integration by part of
Proposition~\ref{pr.EIDLint} (first case) implies
$$
(P^{L}_{+}-i\lambda)u=0\quad\text{in}~\mathcal{D}'(\rz^{2}_{-};
D(L_{-})')\supset L^{2}(\rz^{2}_{-},dqdp;\mathfrak{L})\,.
$$
Now $u$ and any $\varphi\in Y$ fulfill the assumptions of the second
case of Proposition~\ref{pr.EIDLint}. 
We obtain
\begin{eqnarray*}
\frac{1}{2}
\Real\langle \Pi_{-}f_{\partial}\,,\, \Pi_{-}\varphi\rangle_{L^{2}(\rz,|p|dp;\mathfrak{L})}
&&=a(u,\varphi)\\
&&\hspace{-4cm}=0-\lim_{\varepsilon\to 0^{+}}\Real\langle (1+\varepsilon
L_{+})^{-1}\gamma u\,,\, \sign(p)(1+\varepsilon
L_{-})^{-1}\gamma\varphi\rangle_{L^{2}(\rz,|p|dp;\mathfrak{L})}\,,\\
&&\hspace{-4cm}
=\lim_{\varepsilon\to 0^{+}}\frac{1}{2}\Real\langle (1+\varepsilon
L_{+})^{-1}\Pi_{-}\gamma u\,,\, (1+\varepsilon
L_{-})^{-1}\Pi_{-}\gamma\varphi\rangle_{L^{2}(\rz,|p|dp;\mathfrak{L})}\,.
\end{eqnarray*}
By taking an arbitrary $\gamma \varphi\in
\mathcal{C}^{\infty}_{0}((-\infty,0);D(L_{-}))$\,, the right-hand side
converges to 
$$
\frac{1}{2}\langle \Pi_{-}\gamma u\,,\,
\Pi_{-}\gamma\varphi\rangle_{L^{2}(\rz,|p|dp;D(L_{-})'), L^{2}(\rz,|p|dp;D(L_{-}))}\,.
$$
Since the equality holds for all
$$
\varphi\in\Pi_{-}\mathcal{C}^{\infty}_{0}((-\infty,0);D(L_{-}))=
\Pi_{-}\mathcal{C}^{\infty}_{0}(\rz^{*};D(L_{-}))
=\Pi_{ev}\mathcal{C}^{\infty}_{0}(\rz^{*};D(L_{-}))\,,
$$ 
this proves
$$
\Pi_{-}\gamma u=\Pi_{-}f_{\partial}=f_{\partial}\,.
$$
The solution $u$ is characterized by 
$$
\forall \varphi \in D(K^{L}_{-,1})\,,\quad
\langle u\,,\,
(K^{L}_{-,1}+i\lambda)\varphi\rangle=\frac{1}{2}\langle f_{\partial}\,,\,
\Pi_{-}\gamma\varphi\rangle_{L^{2}(\rz,|p dp|;\mathfrak{L})}\,,
$$
and \eqref{eq.charac} is obtained by adding the two cases $f=0$ and
$f_{\partial}=0$\,.\\
\noindent\textbf{Estimate:} We prove now that the solution to
\eqref{eq.inhBVPL} satisfies
\begin{equation}
  \label{eq.estiminhBVPL}
\|u\|_{L^{2}(\rz_{-},dq;\mathcal{H}^{1,L})}+\|\gamma u\|_{L^{2}(\rz,|p|dp; \mathfrak{L})}
\leq C\left[\|f\|_{L^{2}(\rz_{-},dq;\mathcal{H}^{-1,L})}+\|f_{\partial}\|_{L^{2}(\rz,|p|dp;
  \mathfrak{L})}
\right]\,,
\end{equation}
where $C>0$ is independent of $\lambda\in\rz$\,.
The solution $u$ to \eqref{eq.inhBVPL} is the sum
$$
u=u_{1}+u_{0}=(K_{+,1}^{L}-i\lambda)^{-1}f +u_{0}\,,
$$
where $u_{0}$ solves the boundary value problem \eqref{eq.inhBVPL}
with $f=0$\,.\\
Since $f\in L^{2}(\rz_{-},dq;\mathcal{H}^{-1,L})$\,,
Proposition~\ref{pr.extH-1}
gives
\begin{multline*}
\|u_{1}\|_{L^{2}(\rz_{-},dq;\mathcal{H}^{1,L})}^{2}+\|\gamma_{odd}u_{1}\|_{L^{2}(\rz,|p|dp;
  \mathfrak{L})}^{2}=\Real\langle u_{1}\,,f\,
\rangle_{L^{2}(\rz_{-},dq;\mathcal{H}^{1,L})\,, L^{2}(\rz_{-},dq;\mathcal{H}^{-1,L})}
\\
\leq
\|u_{1}\|_{L^{2}(\rz_{-},dq;\mathcal{H}^{1,L})}\|f\|_{L^{2}(\rz_{-},dq;\mathcal{H}^{-1,L})}\,,
\end{multline*}
where the boundary condition
$\gamma_{odd}u_{1}=\sign(p)\gamma_{ev}u_{1}$ implies
$$
\|\gamma_{odd}u_{1}\|_{L^{2}(\rz,|p|dp;\mathfrak{L})}^{2}=\frac{1}{2}
\|\gamma u_{1}\|_{L^{2}(\rz,|p|dp;\mathfrak{L})}^{2}\,.
$$
We deduce 
$$
\|u_{1}\|_{L^{2}(\rz_{-},dq;\mathcal{H}^{1,L})}+\|\gamma
u\|_{L^{2}(\rz,|p|dp;\mathfrak{L})}
\leq 2\|f\|_{L^{2}(\rz,dq;\mathcal{H}^{-1,L})}\,.
$$
The function $u_{0}\in L^{2}(\rz_{-},dq;\mathcal{H}^{1,L})$ solves
\eqref{eq.inhBVPL}
with $f=0$\,. Therefore,
it belongs to $\mathcal{E}^{L,+}_{\rz_{-}}$ and has a trace $\gamma u\in
L^{2}(\rz,|p|dp; D(L_{-})')$ according to
Proposition~\ref{pr.EIDL}. 
We follow the regularization scheme of Proposition~\ref{pr.EIL2}.
For $\varepsilon>0$ the function
$u_{\varepsilon}=(1+\varepsilon L_{+})^{-1}u_{0}$ belongs to
$L^{2}(\rz_{-},dq; \mathcal{H}^{1,L})$ and satisfies
$$
 (P^{L}_{+}-i\lambda)u_{\varepsilon}=0\quad,\quad \gamma
  u_{\varepsilon}\in L^{2}(\rz,|p|dp ;\mathfrak{L})\,.
$$
Hence the integration by part inequality of Proposition~\ref{pr.EIL2}
applies and gives
\begin{eqnarray*}
2\|u_{\varepsilon}\|_{L^{2}(\rz_{-},dq;\mathcal{H}^{1,L})}^{2}&\leq&
-\int_{\rz}|\gamma u_{\varepsilon}(p)|_{\mathfrak{L}}^{2}~pdp
\\
&\leq&
-2\Real\langle \gamma_{ev}u_{\varepsilon}\,,\,
\sign(p)\gamma_{odd}u_{\varepsilon}\rangle_{L^{2}(\rz,|p|dp;
  \mathfrak{L})}
\\
&\leq&
-\frac{1}{2}\|\Pi_{+}\gamma u_{\varepsilon}\|_{L^{2}(\rz,|p|dp;
  \mathfrak{L})}^{2}
+\frac{1}{2}\|f_{\partial}\|^{2}_{L^{2}(\rz,|p|dp;
  \mathfrak{L})}\,,
\end{eqnarray*}
where \eqref{eq.int+-} and $\|(1+\varepsilon L_{+})^{-1}f_{\partial}\|_{L^{2}(\rz,|p|dp;
  \mathfrak{L})}\leq \|f\|_{L^{2}(\rz,|p|dp; \mathfrak{L})}$ were
used for the last inequality.\\
With
\begin{eqnarray*}
\|\gamma u_{\varepsilon}\|_{L^{2}(\rz,|p|dp;
  \mathfrak{L})}
&\leq &
\frac{\|\Pi_{+}\gamma u_{\varepsilon}\|_{L^{2}(\rz,|p|dp;
  \mathfrak{L})}+\|\Pi_{-}\gamma u_{\varepsilon}\|_{L^{2}(\rz,|p|dp;
  \mathfrak{L})}}{2}
\\
&\leq &\|f_{\partial}\|_{L^{2}(\rz,|p|dp;\mathfrak{L})}\,,
\end{eqnarray*}
we infer the uniform bound
$$
\|u_{\varepsilon}\|_{L^{2}(\rz_{-},dq;\mathcal{H}^{1,L})}+\|\gamma
u_{\varepsilon}\|_{L^{2}(\rz,|p|dp;
  \mathfrak{L})}
\leq \frac{3}{2}\|f_{\partial}\|_{L^{2}(\rz,|p|dp;
  \mathfrak{L})}\,,
$$
while we know the strong convergences $\lim_{\varepsilon\to 0}u_{\varepsilon}=u_{0}$ in
$L^{2}(\rz^{2}_{-},dqdp;\mathfrak{L})$ and $\lim_{\varepsilon\to 0}\gamma u_{\varepsilon}=
\gamma u_{0}$ in 
$L^{2}(\rz,|p|dp; D(L_{-})')$\,.
Therefore $u_{0}\in L^{2}(\rz_{-},dq;\mathcal{H}^{1,L})$ has a trace $\gamma
u_{0}\in L^{2}(\rz,|p|dp; \mathfrak{L})$ and
\begin{align*}
\|u_{0}\|_{L^{2}(\rz_{-},dq;\mathcal{H}^{1,L})}
+\|\gamma u_{0}\|_{L^{2}(\rz,|p|dp;
  \mathfrak{L})}
&\leq
 \liminf_{\varepsilon\to 0^{+}}
\left[
\|u_{\varepsilon}\|_{L^{2}(\rz_{-},dq;\mathcal{H}^{1,L})}
\right.\\
&
\quad \phantom{\liminf_{\varepsilon\to 0^{+}}\left[\right.}
\left.+\|\gamma
u_{\varepsilon}\|_{L^{2}(\rz,|p|dp;
  \mathfrak{L})}\right]
\\
&\leq \frac{3}{2} \|f_{\partial}\|_{L^{2}(\rz,|p|dp;
  \mathfrak{L})}\,.
\end{align*}
\noindent\textbf{Integration by part identity:} By
\eqref{eq.estiminhBVPL} all the terms of \eqref{eq.ippinhBVPL} are
continuous with respect to $f\in
L^{2}(\rz_{-},dq;\mathcal{H}^{-1,L})$\,.
We can therefore assume $f\in L^{2}(\rz^{2}_{-},dqdp;\mathfrak{L})$\,.
We set $f_{\varepsilon}=(1+\varepsilon L)^{-1}f$ and
$f_{\partial,\varepsilon}=(1+\varepsilon L)^{-1}f_{\partial}$\,.
Then $u_{\varepsilon}=(1+\varepsilon L)u\in
L^{2}(\rz_{-},dq;\mathcal{H}^{1,L})$ is the unique solution to 
$$
(P_{+}^{L}-i\lambda)u_{\varepsilon}=f_{\varepsilon}\quad,\quad
\Pi_{-}\gamma u_{\varepsilon}=f_{\partial,\varepsilon}\,.
$$ 
This implies
$$
p\partial_{q}u_{\varepsilon}=f_{\varepsilon}+i\lambda u_{\varepsilon}
-(\frac{1}{2}+\mathcal{O})u_{\varepsilon}-Lu_{\varepsilon} \quad\in
L^{2}(\rz_{-},dq;\mathcal{H}^{-1}\otimes \mathfrak{L})\,,
$$
while
$$
u_{\varepsilon}\in L^{2}(\rz_{-},dq;\mathcal{H}^{1}\otimes
\mathfrak{L})
\quad\text{and}
\quad \gamma u_{\varepsilon}=(1+\varepsilon L)^{-1}\gamma u\in L^{2}(\rz,|p|dp;\mathfrak{L})\,.
$$
Proposition~\ref{pr.intbypart} applied with
$\mathfrak{f}=\mathfrak{L}$\,, $a=-\infty$ and $b=0$\,, says
\begin{align*}
2\Real \langle u_{\varepsilon}\,, p\partial_{q}u_{\varepsilon}\rangle
&=\int_{\rz}|u(0,p)|_{\mathfrak{L}}^{2}~pdp
=2\Real\langle
\gamma_{ev}u_{\varepsilon}\,,\,\sign(p)\gamma_{odd}u_{\varepsilon}
\rangle_{L^{2}(\rz,|p|dp;\mathfrak{L})}
\\
&=\frac{1}{2}\left[
\|\Pi_{+}\gamma u_{\varepsilon}\|_{L^{2}(\rz,|p|dp;\mathfrak{L})}^{2}
-
\|\Pi_{-}\gamma u_{\varepsilon}\|_{L^{2}(\rz,|p|dp;\mathfrak{L})}^{2}
\right]\,.
\end{align*}
We obtain
$$
\Real\langle u_{\varepsilon}\,,\, f_{\varepsilon}\rangle 
+\frac{1}{4}\|f_{\partial,\varepsilon}\|_{L^{2}(\rz,|p|dp;\mathfrak{L})}^{2}
=\frac{1}{4}\|\Pi_{+}\gamma
u_{\varepsilon}\|_{L^{2}(\rz,|p|dp;\mathfrak{L})}^{2}
+\|u_{\varepsilon}\|_{L^{2}(\rz_{-},dq,\mathcal{H}^{1,L})}^{2}\,.
$$
With $f\in L^{2}(\rz^{2}_{-},dqdp;\mathfrak{L})$\,,
taking the limit as $\varepsilon \to 0$ gives
\begin{eqnarray*}
  &&
  \lim_{\varepsilon}\|f-f_{\varepsilon}\|_{L^{2}(\rz_{-},dq;\mathcal{H}^{-1,L})}\leq
  \lim_{\varepsilon \to 0}\|f-f_{\varepsilon}\|=0\,,\\
&&
\lim_{\varepsilon\to 0}\|f_{\partial}-f_{\partial,\varepsilon}\|_{L^{2}(\rz,|p|dp;\mathfrak{L})}=0\,.
\end{eqnarray*}
Owing to the continuity estimate \eqref{eq.estiminhBVPL} (this is the
important point), we deduce
\begin{eqnarray*}
  &&
  \lim_{\varepsilon}\|u-u_{\varepsilon}\|_{L^{2}(\rz_{-},dq;\mathcal{H}^{1,L})}=0\,,\\
&&
\lim_{\varepsilon\to 0}\|\gamma u-\gamma u_{\varepsilon}\|_{L^{2}(\rz,|p|dp;\mathfrak{L})}=0\,,
\end{eqnarray*}
and \eqref{eq.ippinhBVPL} is proved.
\end{proof}
An easy consequence of Proposition~\ref{pr.inhBVPL1} is the
following result.
\begin{proposition}
\label{pr.TrScal}  If $u\in L^{2}(\rz_{-},dq;\mathcal{H}^{1,L})$
solves $(P_{\pm}-i\lambda)u=0$ with $\gamma u\in
L^{2}(\rz,|p|dp;\mathfrak{L})$\,, then 
$$
\pm 2\Real\langle \gamma_{ev}u\,,\, \sign(p)\gamma_{odd}u \rangle =
-\|u\|_{L^{2}(\rz_{-},dq;\mathcal{H})^{1}}^{2}\leq 0\,.
$$
\end{proposition}
\begin{proof}
  The function $u$ solves \eqref{eq.inhBVPL} with
  $f_{\partial}=\Pi_{\mp}\gamma u$\,. Apply simply
  \eqref{eq.ippinhBVPL} by referring again to \eqref{eq.int+-}.
\end{proof}
\begin{definition}
\label{de.TC} With the inhomogeneous boundary value problem
\eqref{eq.inhBVPL} with $f=0$\,, $\lambda\in\rz$\,, $f_{\partial}\in
\Pi_{ev}L^{2}(\rz,|p|dp;\mathfrak{L})$ and $u\in
L^{2}(\rz_{-},dq;\mathcal{H}^{1,L})$ solution, we define the two
operators
\begin{eqnarray*}
  && R^{L}_{\pm}(\lambda): \Pi_{ev}L^{2}(\rz,|p|dp;\mathfrak{L})\to
  L^{2}(\rz_{-},dq;\mathcal{H}^{1,L})\quad,\quad R^{L}_{\pm}(\lambda)f_{\partial}=u\,,
\\
&& C^{L}_{\pm}(\lambda): \Pi_{ev}L^{2}(\rz,|p|dp; \mathfrak{L})\to \Pi_{ev}L^{2}(\rz,|p|dp;
\mathfrak{L})\quad,\quad C^{L}_{\pm}f_{\partial}= \Pi_{\pm}\gamma
u\,,
\end{eqnarray*}
with $C^{L}_{\pm}(\lambda)=\Pi_{\pm}\circ \gamma\circ R^{L}_{\pm}(\lambda)$\,.
\end{definition}

\begin{proposition}
\label{pr.imfermeecontr}
Assume Hypothesis~\ref{hyp.L1} and take $\lambda\in\rz$\,.\\
The operator $R^{L}_{\pm}(\lambda)$ is continuous and injective.\\
Its adjoint is 
$$
(R^{L}_{\pm}(\lambda))^{*}=\Pi_{\mp}\circ \gamma \circ
(K^{L}_{\mp}+i\lambda)^{-1}:
L^{2}(\rz_{-},dq;\mathcal{H}^{-1,L})\to
\Pi_{ev}L^{2}(\rz,|p|dp;\mathfrak{L})
$$
and its range is dense in $\Pi_{ev}L^{2}(\rz,|p|dp;\mathfrak{L})$.\\
The operator $C^{L}_{\pm}(\lambda)$ is a contraction of $\Pi_{ev}L^{2}(\rz,dp;\mathfrak{L})$:
$$
\forall f_{\partial}\in \Pi_{ev}L^{2}(\rz,|p|dp;\mathfrak{L})\,,\quad
\|C^{L}_{\pm}(\lambda)f_{\partial}\|_{L^{2}(\rz,|p|dp;\mathfrak{L})}\leq
\|f_{\partial}\|_{L^{2}(\rz,|p|dp;\mathfrak{L})}\,.
$$
\end{proposition}
\begin{proof}
The operators $R_{\pm}(\lambda)$ and $C_{\pm}(\lambda)$ are continuous
owing to \eqref{eq.ippinhBVPL}. This also proves that
$C_{\pm}(\lambda)$ is a contraction.\\
\noindent\textbf{Injectivity:} If $u=R^{L}_{\pm}(\lambda)f_{\partial}=0$ then
$\gamma u=0$ and $f_{\partial}=\Pi_{\mp}\gamma u=0$\,.\\
\noindent\textbf{Adjoint:} In the proof of
Proposition~\ref{pr.inhBVPL1}, $u=R^{L}_{\pm}f_{\partial}$ is
characterized by 
$$
\forall \varphi \in D(K^{L}_{\mp,1})\,,\quad
\langle u\,,\, (K^{L}_{\mp,1}+i\lambda)\varphi\rangle=
\frac{1}{2}\langle f_{\partial}\,,\, \Pi_{\mp}\varphi\rangle_{L^{2}(\rz,|p|dp;\mathfrak{L})}\,.
$$
But both sides extend by continuity to any $\varphi$ such that
$(K^{L}_{\mp}+i\lambda)\varphi\in
L^{2}(\rz_{-},dq;\mathcal{H}^{-1,L})$ according to
Proposition~\ref{pr.extH-1}. Hence we get
$$
\forall f\in L^{2}(\rz_{-},dq;\mathcal{H}^{-1,L})\,,\quad
\langle R^{L}_{\pm}f_{\partial}\,,\, f\rangle=\frac{1}{2}
\langle f_{\partial}\,,\, \Pi_{\mp}\gamma (K^{L}_{\mp}+i\lambda)^{-1}f\rangle_{L^{2}(\rz,|p|dp;\mathfrak{L})}\,.
$$
This proves $R^{L}_{\pm}(\lambda)^{*}=\frac{1}{2}\Pi_{\mp}\circ
\gamma\circ (K^{L}+i\lambda)^{-1}$ and since $R^{L}_{\pm}$ is
injective,
 $R^{L}_{\pm}(\lambda)^{*}$  has a dense range.
\end{proof}

One may wonder about $C_{\pm}(\lambda)$ being a strict contraction:
$$
\|C_{\pm}^{L}(\lambda)\|_{\mathcal{L}(\Pi_{ev}L^{2}(\rz,|p|dp;\mathfrak{L}))}<
1\,.
$$
The following statements are equivalent:
\begin{description}
\item[i)] $R_{\pm}^{L}(\lambda)$ has a closed range;
\item[ii)] $R_{\pm}^{L}(\lambda)^{*}$ is surjective; 
\item[iii)] $C_{\pm}^{L}(\lambda)$ is a strict contraction;
\end{description}
Proving one of them, for the specific case when $L_{\pm}$ is a
geometric Kramers-Fokker-Planck operator, would allow to consider more general boundary
conditions and to relax the gap assumption $c_{A}>0$
in~\eqref{eq.intArq2} for Theorem~\ref{th.mainA}.

\section{General boundary conditions for half-space problems}
\label{se.genBC}
The analysis of abstract half-space problems is enhanced by
introducing more general boundary conditions which possibly couple the
variables. 
We keep the notations introduced in Section~\ref{se.insvar}.\\
By referring to Appendix~\ref{se.kfpline}, the analysis developed here
proves Theorem~\ref{th.main0} and Theorem~\ref{th.mainA} when
$\overline{Q}=(-\infty,0]\times \tz^{d_{1}'}\times \rz^{d_{2}'}$ is
endowed with the euclidean metric.

\subsection{Assumptions for $L$ and $A$}
\label{se.assLA}

We keep the notations introduced in Subsection~\ref{se.notasep}.
We shall need some additional assumptions for the vertical maximal
accretive operators $(L_{\pm},D(L_{\pm}))$\,, $L_{-}=L_{+}^{*}$\,,
 defined in $\mathfrak{L}$\,. Those are
formulated in terms of the maximal accretive operators
$K_{\pm}^{L}=P_{\pm}^{L}=\pm p\partial_{q}+\frac{1}{2}+\mathcal{O}+L_{\pm}$\,, 
defined on the whole line
$q\in\rz$\,, 
with $D(K_{\pm}^{L})\subset L^{2}(\rz^{2},dqdp; \mathfrak{L})$\,. We
have checked in Proposition~\ref{pr.idenDK} the equality $(K_{\pm}^{L})^{*}=K_{\mp}^{L}$\,.\\
Remember that Hypothesis~\ref{hyp.L1} says that $D(L_{+})\cap D(L_{-})$ is
a core for $\Real(L)^{1/2}$\,. This is strengthened by
\begin{assumption}
\label{hyp.whole} The domain $D(L_{+})=D(L_{-})$ is a core for $\Real(L)$\,.
The domains $D(K_{\pm}^{L})$ and $D(K_{\mp}^{L})$ for the whole line problem
are equal, 
 and the two operators $K_{+}^{L}$ and $K_{-}^{L}$ fulfill the estimate
$$
\|(K_{\pm}^{L}-K_{\mp}^{L}-2i\lambda)u\|+
\|(K_{\pm}^{L}+K_{\mp}^{L})u\|+ \langle \lambda\rangle^{1/2}\|u\|
+\|u\|_{\mathcal{Q}_{0}}\leq C_{L}\|(K_{\pm}^{L}-i\lambda)u\|\,,
$$
for all
$\lambda\in\rz$  and all $u\in D(K_{\pm}^{L})$ and where
 $\mathcal{Q}_{0}$ is a Hilbert space embedded in
$L^{2}(\rz^{2},dqdp;\mathfrak{L})$\,.
\end{assumption}
In applications, the Hilbert space $\mathcal{Q}_{0}$ will encode the regularity 
properties with respect to $q$\,. A stronger version is also needed
for general boundary conditions.
\begin{assumption}
\label{hyp.Q}
There exists a Hilbert space $\mathcal{Q}$ embedded in
$L^{2}(\rz^{2},dqdp;\mathfrak{L})$
such that any solution $v\in L^{2}(\rz,dq;\mathcal{H}^{1,L})$ to 
\begin{equation}
  \label{eq.PLgd}
(P_{\pm}^{L}-i\lambda)v=
\gamma\delta_{0}(q)\quad\text{in}~\mathcal{S}'(\rz^{2};D(L_{\mp})')\,,
\end{equation}
with $\gamma\in L^{2}(\rz,\frac{dp}{|p|};\mathfrak{L})$\,, 
satisfies
\begin{equation}
  \label{eq.hypQ1}
\|v\|_{\mathcal{Q}}\leq
C_{\mathcal{Q},L}\left[\|\gamma\|_{L^{2}(\rz,\frac{dp}{|p|};\mathfrak{L})}+
\|v\|_{L^{2}(\rz,dq;\mathcal{H}^{1,L})}\right]\,.
\end{equation}
We also assume
\begin{equation}
  \label{eq.hypQ2}
\forall u\in D(K_{\pm}^{L})\,,\quad
\langle \lambda\rangle^{\frac{1}{4}}
\|u\|_{\mathcal{Q}}\leq C_{\mathcal{Q},L}\|(K_{\pm}^{L}-i\lambda)u\|\,.
\end{equation}
\end{assumption}
\noindent\textbf{Example:}
An easy example which fulfills all these assumptions
is the Kramers-Fokker-Planck operator $L_{\pm}=\pm
p'.\partial_{q'}+\frac{-\Delta_{p'}+|p'|^{2}}{2}$ defined in
$L^{2}(X',dq'dp')$ with $X'=T^{*}Q'$ and $Q'=\rz^{d_{1}'}\times
\tz^{d_{2}'}$ endowed in the euclidean metric. The operator
$K_{\pm}^{L}-\frac{1}{2}=p.\partial_{q}+\frac{-\Delta_{p}+|p|^{2}}{2}$
is then the Kramers-Fokker-Planck operator in
$L^{2}(T^{*}Q,dqdpdq'dp')$ and $Q=\rz^{d_{1}'+1}\times \tz^{d_{2}'}$
is endowed with the euclidean metric.
Hypothesis~\ref{hyp.whole} is thus verified with
$\mathcal{Q}_{0}=H^{\frac{2}{3}}(Q;\mathcal{H}^{0})=H^{\frac{2}{3}}(Q;L^{2}(\rz^{1+d_{1}'+d_{2}'},dp))$
according to Theorem~\ref{th.app} applied with $s=s'=0$\,.
Hypothesis~\ref{hyp.Q} is verified with
$\mathcal{Q}=H^{t}(Q;\mathcal{H}^{0})=H^{t}(Q;L^{2}(\rz^{1+d_{1}'+d_{2}'},dp))$ for any $t\in
[0,\frac{1}{9})$: the first statement and the estimate \eqref{eq.PLgd}
is provided by Proposition~\ref{pr.regQcomp}; the estimate
\eqref{eq.hypQ2} is inferred by interpolating the inequalities of
Theorem~\ref{th.app} with $s=s'=0$\,,
\begin{multline*}
  \|u\|_{H^{t}(Q;\mathcal{H}^{0})}^{2}\leq
  \|u\|_{H^{\frac{1}{3}}(Q;\mathcal{H}^{0})}^{2}
\leq C\|u\|_{H^{\frac{2}{3}}(Q;\mathcal{H}^{0})}\|u\|
\leq C'\langle \lambda\rangle^{-\frac{1}{2}}\|(K_{\pm}^{L}-i\lambda)u\|^{2}\,.
\end{multline*}
The case when $L_{\pm}$ is a geometric Kramers-Fokker-Planck operator
on a general riemannian compact manifold will be studied in Section~\ref{se.geoKFP}.\\
The general boundary conditions at $q=0$\,, are formulated in terms of an
operator $(A,D(A))$ defined in $L^{2}(\rz,|p|dp;\mathfrak{L})$\,.
They are written according to the introduction as
$$
\gamma_{odd}u=\sign(p)A\gamma_{ev}u\quad,\quad \gamma u(p)=u(0,p)\,.
$$ 
Remember the Definition~\ref{de.invoL} for the unitary 
involution $j$ on $\mathfrak{L}$\,, for $\gamma_{ev}=\Pi_{ev}\gamma$ and
$\gamma_{odd}=\Pi_{odd}\gamma$ for $\gamma \in
L^{2}(\rz,|p|dp;\mathfrak{L})$\,. All these operators $j, \Pi_{ev}$
and $\Pi_{odd}$ commute with $L_{\pm}$\,.
We need some compatibility  between $(A,D(A))$ and the involution $j$
via $\Pi_{ev,odd}$\,, and we focus on the case $D(A)=L^{2}(\rz,|p|dp;\mathfrak{L})$\,.
\begin{assumption}
\label{hyp.Abdd} The operator $A$\,, $L^{2}(\rz,|p|dp;
\mathfrak{L})$\,, is a bounded accretive operator 
which commutes with the  orthogonal projections 
$\Pi_{ev}$ and $\Pi_{odd}=1-\Pi_{ev}$\,.
The norm of $A$ in $\mathcal{L}(L^{2}(\rz,|p|dp;\mathfrak{L}))$ is
simply written $\|A\|$\,.\\
We assume additionally
\begin{description}
\item[either] $c_{A}=\min \sigma(\Real A)>0$\,;
\item[or] $A=0$\,.
\end{description}
\end{assumption}
Mixed boundary conditions with $A=A_{1}\oplus 0$ when all the
operators are block diagonal in $\mathfrak{L}=\mathfrak{L}_{1}\oplus
\mathfrak{L}_{0}$ will follow at once from the analysis of the distinct
two cases. This will cover all the applications that we have in mind. 
More generally vanishing $A$'s are not considered in this article.
\subsection{Maximal accretivity}
\label{se.maxaccA}
In this section we check the maximal accretivity of the realizations of
$$
P_{\pm}^{L}=\pm p\partial_{q}+\frac{1}{2}+\mathcal{O}+L_{\pm}
$$
on $\rz^{2}_{-}=\rz_{-}\times \rz$\,, with boundary conditions given by
$\gamma_{odd}u=\pm\sign(p)A\gamma_{ev}u$\,.\\
The case $A=0$ was treated in Proposition~\ref{pr.idenDKAId} without
specifying $\gamma u\in L^{2}(\rz,|p|dp;\mathfrak{L})$\,. This last
point will be done in Subsection~\ref{se.halfwhL}. We now focus on the
case $A\neq 0$ of Hypothesis~\ref{hyp.Abdd}.\\
The strategy is a reduction to the boundary $\left\{q=0\right\}$\,, of the equation
$$
P_{\pm}^{L}u=f\quad\text{in}~L^{2}(\rz^{2}_{-},dqdp ;\mathfrak{L})~,\quad
\gamma_{odd}u=\pm\sign(p)A\gamma_{ev}u\quad\text{in}~L^{2}(\rz,|p|dp;\mathfrak{L})\,,
$$
where the absence of a Calderon projector, will be compensated by the
nice properties of $K^{L}_{\pm,1}$ (Hypothesis~\ref{hyp.L1} for
$L_{\pm}$ suffices here) combined with
Hypothesis~\ref{hyp.Abdd}.
\begin{proposition}
\label{pr.maxaccKAL}
  Assume Hypothesis~\ref{hyp.L1} and Hypothesis~\ref{hyp.Abdd} with
  $A\neq 0$\,.
Then the operator
$K^{L}_{\pm, A}-1$ defined by 
\begin{eqnarray*}
  &
D(K^{L}_{\pm, A})=\left\{\begin{array}[c]{ll}
u\in L^{2}(\rz_{-},dq;\mathcal{H}^{1,L})\,,
&P_{\pm}^{L}u\in
    L^{2}(\rz^{2}_{-},dqdp;\mathfrak{L})  
\\
&\gamma_{ev} u\in L^{2}(\rz,|p|dp;\mathfrak{L})\\
&\gamma_{odd}u=\pm \sign(p) A\gamma_{ev}u
    \end{array}
\right\}\,,\\
&\forall u\in D(K^{L}_{\pm, A})\,,~K^{L}_{\pm, A}u=P_{\pm}^{L}u
=(\pm p.\partial_{q}+\frac{1}{2}+\mathcal{O}+L)u\,,
\end{eqnarray*}
is densely defined and maximal accretive in $L^{2}(\rz_{-}^{2},dqdp;
\mathfrak{L})$\,.\\
Any $u\in D(K^{L}_{\pm, A})$ satisfies the integration by part identity
\begin{equation}
  \label{eq.idenipp}
\|u\|_{L^{2}(\rz_{-},dq;\mathcal{H}^{1,L})}^{2}+ \Real\langle
\gamma_{ev}u\,,
A\gamma_{ev}u\rangle_{L^{2}(\rz,|p|dp;\mathfrak{L})}
= \Real \langle
u\,,\, K^{L}_{\pm, A}u\rangle\,.
\end{equation}
Moreover for any $\lambda\in\rz$ the operators
$(K^{L}_{\pm, A}-i\lambda)^{-1}$ and
$\gamma\circ(K^{L}_{\pm, A}-i\lambda)^{-1}$ have continuous extensions
from $L^{2}(\rz_{-},dq;\mathcal{H}^{-1,L})$ respectively to
$L^{2}(\rz_{-},dq;\mathcal{H}^{1,L})$ and
$L^{2}(\rz,|p|dp;\mathfrak{L})$\,, with $\lambda$-dependent
norms.\\
Finally the adjoint of $K^{L}_{\pm, A}$ is $K^{L}_{\mp, A^{*}}$\,.
\end{proposition}
 Below is a 
more convenient (and more usual) rewriting of the boundary conditions.
\begin{lemma}
\label{le.eqA+-}
Let $(A,D(A))$ be a maximal accretive operator in $L^{2}(\rz,|p|dp;\mathfrak{L})$\,.
 For $\gamma\in L^{2}(\rz,|p|dp;\mathfrak{L})$\,, the
two relations
\begin{eqnarray*}
  &&\Pi_{odd}\gamma=\pm\sign(p)A\Pi_{ev}\gamma  \quad (\Pi_{ev}\gamma\in
  D(A))\,,\\
\text{and}
&&
\Pi_{\mp}\gamma=\frac{(1-A)}{(1+A)}\Pi_{\pm}\gamma
\end{eqnarray*}
are equivalent.\\
When $A$ fulfills Hypothesis~\ref{hyp.Abdd}, 
$$
\|\frac{1-A}{1+A}\|_{\mathcal{L}(L^{2}(\rz,|p|dp;\mathfrak{L}))}\leq
\left(1+\frac{2c_{A}}{1+\|A\|^{2}}
\right)^{-\frac{1}{2}}<1\,.
$$
\end{lemma}
\begin{proof}
With 
$$
\Pi_{+}=\Pi_{ev}+\sign(p)\Pi_{odd}\quad,\quad \Pi_{-}=\Pi_{ev}-\sign(p)\Pi_{odd}\,,
$$
one gets
\begin{eqnarray*}
  \left(\Pi_{\mp}\gamma=\frac{1-A}{1+A}\Pi_{\pm}\gamma\right)
&\Leftrightarrow&
\left(
\frac{2A}{1+A}\Pi_{ev}\gamma=\pm\frac{1}{1+A}\sign(p)\Pi_{odd}\gamma
\right)\\
&\Leftrightarrow&
\left(
A\Pi_{ev}\gamma=\pm\sign(p)\Pi_{odd}\gamma~\text{in}~D(A^{*})'
\right)
\\
&\stackrel{\hspace{-2cm}\Pi_{odd}\gamma \in
  L^{2}(\rz,|p|dp;\mathfrak{L})}{\Leftrightarrow}&
\left(
  \begin{array}[c]{l}
    \Pi_{ev}\gamma\in D(A)\\
    \Pi_{odd}\gamma=\pm\sign(p)A\Pi_{ev}\gamma
  \end{array}
\right)\,.
\end{eqnarray*} 
When $A$ is bounded and $\Real A\geq c_{A}>0$\,, the inequalities
\begin{eqnarray*}
  \|(1+A)u\|_{L^{2}(\rz,|p|dp;\mathfrak{L})}^{2}-\|(1-A)u\|_{L^{2}(\rz,|p|dp;\mathfrak{L})}^{2}
&=&2\Real\langle
  u\,,\, Au\rangle_{L^{2}(\rz,|p|dp;\mathfrak{L})}
\\
&\geq& 2c_{A}^{2}\|u\|_{L^{2}(\rz,|p|dp;\mathfrak{L})}^{2}\\
(1+\|A\|^{2})\|u\|_{L^{2}(\rz,|p|dp;\mathfrak{L})}^{2}&\geq & \|(1-A)u\|_{L^{2}(\rz,|p|dp;\mathfrak{L})}^{2}\,,
\end{eqnarray*}
yields
$$
(1+\frac{2c_{A}}{1+\|A\|^{2}})\|(1-A)u\|_{L^{2}(\rz,|p|dp;\mathfrak{L})}^{2}
\leq \|(1+A)u\|_{L^{2}(\rz,|p|dp;\mathfrak{L})}^{2}\,.
$$
\end{proof}
\begin{proof}[Proof of Proposition~\ref{pr.maxaccKAL}]
The domain contains $\mathcal{C}^{\infty}_{0}(\rz^{2}_{-};D(L))$ which
is dense in $L^{2}(\rz^{2}_{-},dqdp;\mathfrak{L})$\,.\\
  When $u$ belongs to $D(K^{L}_{\pm, A})$\,, it satisfies all the
  assumptions of Proposition~\ref{pr.EIL2}. We obtain
  \begin{eqnarray*}
    \Real\langle u\,,\, K^{L}_{\pm,A}u\rangle
&\geq&
    \|u\|_{L^{2}(\rz_{-},dq;\mathcal{H}^{1,L})}^{2}
\pm\Real\langle \gamma_{ev}u\,,\,\sign(p)\gamma_{odd}u\rangle_{L^{2}(\rz,|p|dp;\mathfrak{L})}\\
&\geq&
    \|u\|_{L^{2}(\rz_{-},dq;\mathcal{H}^{1,L})}^{2}
+\Real\langle \gamma_{ev}u\,,\,
A\gamma_{ev}u\rangle_{L^{2}(\rz,|p|dp;\mathfrak{L})}\\
&\geq & \|u\|_{L^{2}(\rz_{-},dq;\mathcal{H}^{1,L})}^{2}\geq \|u\|^{2}\,.
  \end{eqnarray*}
and
$K^{L}_{\pm,A}-1$ is accretive\,.\\
Let $f$ belong to $L^{2}(\rz^{2}_{-},dqdp; \mathfrak{L})$\,, we want to
find $u\in D(K^{L}_{\pm,A})$ such that $(K^{L}_{\pm,A}-i\lambda) u=f$\,.
By Proposition~\ref{pr.idenDKAId},  $u_{1} =(K_{\pm,1}^{L}-i\lambda)^{-1}f$
belongs to 
$L^{2}(\rz_{-},dq;\mathcal{H}^{1,L})$
 with $\gamma u_{1}\in L^{2}(\rz,|p|dp;\mathfrak{L})$\,,
$\gamma_{ev} u_{1}=\pm\sign(p)\gamma_{odd}u_{1}$ or  $\Pi_{\mp}\gamma
u_{1}=0$\,.\\
Set $u=u_{1}+v$ and the equation becomes
\begin{eqnarray*}
  && (P_{\pm}^{L}-i\lambda)v=0\quad,\quad v\in
  L^{2}(\rz_{-},dq;\mathcal{H}^{-1,L})\,,\\
&& \Pi_{\mp}\gamma v=\Pi_{\mp}\gamma u= (\frac{1-A}{1+A})\Pi_{\pm}\gamma u
=\frac{1-A}{1+A}\Pi_{\pm}\gamma u_{1}+ \frac{1-A}{1+A}\Pi_{\pm}\gamma v\,.
\end{eqnarray*}
According to Definition~\ref{de.TC}, the problem is solved when
$\Pi_{\pm}\gamma v$ is a solution to
$$
\gamma=C^{L}_{\pm}(\lambda)(\frac{1-A}{1+A}\Pi_{\pm}\gamma u_{1}+
\frac{1-A}{1+A}\gamma)\quad,\quad \gamma\in \Pi_{ev}L^{2}(\rz,|p|dp;\mathfrak{L})\,.
$$
Owing to $\|C^{L}_{\pm}(\lambda)\|\leq 1$
and the estimate $\|\frac{1-A}{1+A}\|\leq
(1+\frac{2c_{A}}{1+\|A\|^{2}})^{-\frac{1}{2}}$ of Lemma~\ref{le.eqA+-}\,, the
above equation admits a unique solution $\gamma_{\pm}$\,, with 
$$
\|\gamma_{\pm}\|_{L^{2}(\rz,|p|dp;\mathfrak{L})}\leq
C_{A}\|\Pi_{\pm}\gamma u_{1}\|_{L^{2}(\rz,|p|dp;\mathfrak{L})}
\leq C_{A}'\|f\|_{L^{2}(\rz_{-},dq;\mathcal{H}^{-1,L})}\,.
$$
Take
$$
u=(K^{L}_{\pm,1}-i\lambda)^{-1}f +
R_{\pm}^{L}(\lambda)\left[\frac{1-A}{1+A}\Pi_{\pm}\gamma u_{1}+ \frac{1-A}{1+A}\gamma_{\pm}\right]\,.
$$
It belongs to $L^{2}(\rz_{-},dq;\mathcal{H}^{1,L})$\,, solves
$(P_{\pm}^{L}-i\lambda)u=f$\,. Its trace, as the sum of two terms, belongs
to $L^{2}(\rz,|p|dp;\mathfrak{L})$ and solves
$$
\Pi_{\mp}\gamma(u)=\frac{1-A}{1+A}\Pi_{\pm}\gamma u_{1}+
\frac{1-A}{1+A}\gamma_{\pm}
=\frac{1-A}{1+A}\Pi_{\pm}\gamma u\,.
$$
This ends the proof of the maximal accretivity of $K_{\pm,A}^{L}$\,.\\
The function $u=(K_{\pm,A}-i\lambda)^{-1}f$ solves the boundary value
problem \eqref{eq.inhBVPL} with $f_{\partial}=\Pi_{\mp}\gamma u$\,. 
Therefore
Proposition~\ref{pr.inhBVPL1} now provides the integration by part
equality \eqref{eq.idenipp}
\begin{eqnarray*}
  \Real\langle u\,,\, (K_{\pm, A}^{L}-i\lambda) u
\rangle=
\frac{1}{4}\|\Pi_{+}\gamma u\|_{L^{2}(\rz,|p|dp;\mathfrak{L})}^{2}
-\frac{1}{4}\|\Pi_{-}\gamma u\|_{L^{2}(\rz,|p|dp;\mathfrak{L})}^{2}&
\\
+\|u\|_{L^{2}(\rz_{-},dq;\mathcal{H}^{1,L})}^{2}&\,,
\end{eqnarray*}
after using \eqref{eq.int+-}.\\
This identity \eqref{eq.idenipp} combined with
Hypothesis~\ref{hyp.Abdd} ($A\neq 0$) implies
\begin{multline*}
  \frac{c_{A}}{2}\|\gamma_{ev}u\|^{2}+\frac{c_{A}}{2\|A\|^{2}}\|\gamma_{odd}u\|^{2}
+\|u\|_{L^{2}(\rz_{-},dq;\mathcal{H}^{1,L})}^{2}
\\
\leq
\|u\|_{L^{2}(\rz_{-},dq;\mathcal{H}^{1}_{L})}\|(K_{\pm,A}^{L}-i\lambda)u\|_{L^{2}(\rz_{-},dq;\mathcal{H}^{-1,L})}\,, 
\end{multline*}
which ensures the continuous extensions
\begin{eqnarray*}
  &&(K_{\pm,A}^{L}-i\lambda)^{-1}:L^{2}(\rz_{-},dq;\mathcal{H}^{-1,L})\to
L^{2}(\rz,dq;\mathcal{H}^{1,L})\\
&& \gamma(K_{\pm,A}^{L}-i\lambda)^{-1} :
L^{2}(\rz_{-},dq;\mathcal{H}^{-1,L})\to
L^{2}(\rz,|p|dp;\mathfrak{L})\,.
\end{eqnarray*}
Two functions $u\in D(K^{L}_{+,A})$ and
$v\in(K^{L}_{-,A^{*}})$ fulfill the conditions of
Proposition~\ref{pr.EIDLint} (second case) and we get
\begin{eqnarray*}
\langle v\,,\, K^{L}_{+,A}u\rangle
&=&\langle v\,,\, P_{+}^{L}u\rangle\\
&=&\langle \gamma v\,,\,\sign(p)\gamma
u\rangle_{L^{2}(\rz,|p|dp;\mathfrak{L})}
+\langle P^{L}_{-}v\,,\,u\rangle\\
&=&\langle
\gamma_{ev}v\,,\,\sign(p)\gamma_{odd}u\rangle_{L^{2}(\rz,|p|dp;\mathfrak{L})}
\\
&&\qquad+\langle\sign(p)\gamma_{odd}v\,,\,
\gamma_{ev}u\rangle_{L^{2}(\rz,|p|dp;\mathfrak{L})}
+\langle K^{L}_{-,A^{*}}v\,,\,u\rangle\\
&=&\langle-A^{*}\gamma_{ev}v\,,\,
\gamma_{ev}u\rangle_{L^{2}(\rz,|p|dp;\mathfrak{L})}
\\
&&\qquad
+\langle\gamma_{ev}v\,,\,
A\gamma_{ev}u\rangle_{L^{2}(\rz,|p|dp;\mathfrak{L})}
+\langle K^{L}_{-,A^{*}}v\,,\,u\rangle\\
&=&\langle K^{L}_{-,A^{*}}v\,,\,u\rangle\,.
\end{eqnarray*}
This proves $K^{L}_{-,A^{*}}\subset (K^{L}_{+,A})^{*}$ and
therefore the equality, because
$K^{L}_{-,A^{*}}$ is maximal accretive.
\end{proof}

\subsection{Half-space and whole space problem}
\label{se.halfwhL}
\begin{definition}
\label{de.sigma}
Let $\mathfrak{f}$ be one of the Hilbert spaces $\mathfrak{L}$\,,\,
$D(L_{\pm})$ or \, $D(L_{\pm})'$\,, which are endowed
with the involution $j$ (see Definition~\ref{de.invoL})\,. 
Let $\mathcal{H}^{s}$\,, $s\in\rz$\,, be
the  space 
$\mathcal{H}^{s}=\left\{u\in S'(\rz)\,, (\frac{1}{2}+\mathcal{O})^{s/2}u\in L^{2}(\rz,dp)\right\}$\,.\\
  The operator $\Sigma:L^{2}(\rz_{-},dq;\mathcal{H}^{s}\otimes
\mathfrak{f})\to L^{2}(\rz,dq;\mathcal{H}^{s}\otimes\mathfrak{f})$ is
  defined by
$$
\Sigma u (q,p)=
\left\{
\begin{array}[c]{cc}
  u(q,p)&\text{if}~q<0\,,\\
  ju(-q,-p)&\text{if}~q>0\,.
\end{array}
\right.
$$
The operator $\tilde{\Sigma}:L^{2}(\rz,dq;\mathcal{H}^{s}\otimes
\mathfrak{f})\to L^{2}(\rz,dq;\mathcal{H}^{s}\otimes
\mathfrak{f})$ is defined by
$$
\tilde{\Sigma}u(q,p)=j u(-q,-p)\,.
$$
\end{definition}

\begin{proposition}
 \label{pr.extsym} With Hypothesis~\ref{hyp.whole}, assume  that $u\in
 L^{2}(\rz_{-},dq;\mathcal{H}^{1,L})\subset
 L^{2}(\rz^{2}_{-},dqdp;\mathfrak{L})$
 solves $(P_{\pm}^{L}-i\lambda)u=f\in
 L^{2}(\rz_{-},dq;\mathcal{H}^{-1,L})\subset L^{2}(\rz_{-},dq;
 \mathcal{H}^{-1}\otimes D(L_{\mp})')$ with $\lambda\in\rz$\,. 
Then $\Sigma u$ solves
$$
(P_{\pm}^{L}-i\lambda)\Sigma u=\Sigma f \mp 2p(\gamma_{odd}u)\delta_{0}(q)\,,
$$
in $\mathcal{S}'(\rz^{2}; D(L_{\mp})')$\,.\\
Moreover if $v\in L^{2}(\rz,dq;\mathcal{H}^{1,L})$ solves 
$$
(P_{\pm}^{L}-i\lambda)v=\Sigma f \mp 2p(\gamma_{odd}u)\delta_{0}(q)
$$
in $\mathcal{S}'(\rz^{2}; D(L_{\mp})')$\,, then
$v=\Sigma u$\,.
\end{proposition}
\begin{proof}
When $u\in L^{2}(\rz_{-},dq;\mathcal{H}^{1,L})\subset
L^{2}(\rz^{2}_{-},dqdp;\mathfrak{L})$
 solves
$(P_{\pm}^{L}-i\lambda)u=f$ with $f\in L^{2}(\rz_{-},dq;\mathcal{H}^{-1,L})\subset
L^{2}(\rz^{2}_{-},dqdp;\mathfrak{L})$\,, 
it belongs to
$\mathcal{E}_{\rz_{-}}^{L,\pm}\subset\mathcal{E}_{\rz_{-}}(D(L_{\mp})')$ and solves
$$
(\pm p\partial_{q}+\frac{1}{2}+\mathcal{O}-i\lambda)u=f -
L_{\pm}u\quad\text{in}~L_{\pm}^{2}(\rz_{-}^{2},dqdp; D(L_{\mp})')\,,
$$
with $\gamma u\in L^{2}(\rz,|p|dp; D(L_{\mp})')$\,.\\
With $\left[(1+L_{\pm})^{-1}\,,\, j\right]=0$ and
$$
[\pm p\partial_{q}+\frac{1}{2}+\mathcal{O}][ju(-q,-p)]= j
[(\pm p\partial_{q}+\frac{1}{2}+\mathcal{O})u](-q,-p)\quad \text{in}~\rz^{2}\setminus\left\{q=0\right\}\,,
$$
we obtain
\begin{eqnarray*}
(\pm p\partial_{q}+\frac{1}{2}+\mathcal{O})(\Sigma u)
&=&
\Sigma f-L_{\pm}\Sigma u \pm p\left[\Sigma
  u(0^{+},p)-\Sigma u(0^{-},p)\right]\delta_{0}(q)
\\
&=&
\Sigma f-L_{\pm}\Sigma u \mp 2p\gamma_{odd}u\delta_{0}(q)\,,
\end{eqnarray*}
in $\mathcal{S}'(\rz^{2}; D(L^{*})')$\,, which is
what we seek.\\
When $v\in L^{2}(\rz,dq;\mathcal{H}^{1,L})$ is another solution to 
$(P_{\pm}^{L}-i\lambda)v=\Sigma f \mp 2p\gamma_{odd}u\delta_{0}(q)$\,, then the difference
$w=v-\Sigma u$ belongs to $L^{2}(\rz,dq;\mathcal{H}^{1,L})$ and solves
$(P^{L}_{\pm}-i\lambda)w=0$\,. Therefore it belongs to $D(K_{\pm}^{L})$ and solves
$K_{\pm}^{L}w=0$\,, which implies $w=0$\,.
\end{proof}
From the previous result, one infers the equivalence
\begin{equation}
  \label{eq.equiKL0}
(u=(K^{L}_{\pm,0}-i\lambda)^{-1}f)\Leftrightarrow (\Sigma u=(K_{\pm}^{L}-i\lambda)^{-1}\Sigma
f)
\end{equation}
which combined with Hypothesis~\ref{hyp.whole} provides good estimates
for $(K_{\pm, 0}^{L}-i\lambda)^{-1}$ with respect to
$\lambda\in\rz$\,.
Actually this proves $\gamma u\in L^{2}(\rz,|p|dp;\mathfrak{L})$ when
$u\in D(K_{\pm,0}^{L})$ and
it solves completely the case $A=0$\,.
\begin{proposition}
  \label{pr.trKL0} 
With Hypothesis~\ref{hyp.whole}, there exists a constant $C>0$ such
that
\begin{multline*}
\langle \lambda\rangle^{\frac{1}{2}}\|u\|+
\|(\frac{1}{2}+\mathcal{O}+\Real L)u\|+\|(\pm
p\partial_{q}+\frac{L_{\pm}-L_{\mp}}{2}-i\lambda)u\|
+
\|\Sigma u\|_{\mathcal{Q}_{0}}
\\
+\langle \lambda\rangle^{\frac{1}{4}}\|u\|_{L^{2}(\rz_{-},dq;\mathcal{H}^{1,L})}
+\langle \lambda\rangle^{\frac{1}{4}}
\|\gamma u\|_{L^{2}(\rz,|p|dp;\mathfrak{L})}
\leq C\|(K_{\pm,0}^{L}-i\lambda)u\|
\end{multline*}
holds for all $\lambda\in\rz$ and all $u\in D(K_{\pm,0}^{L})$\,.\\
If $\mathcal{D}$ is dense in $D(K^{L}_{\pm})\subset
L^{2}(\rz,dqdp;\mathfrak{L})$ endowed with the graph norm, with
$\tilde{\Sigma}\mathcal{D}=\mathcal{D}$\,, then the set $\left\{u\in
  L^{2}(\rz_{-}^{2},dqdp;\mathfrak{L})\,, \Sigma u\in
  \mathcal{D}\right\}$ is dense in $D(K_{\pm,0}^{L})$ endowed with its
graph norm.
\end{proposition}
\begin{proof}
With $2\|(K_{\pm,0}^{L}-i\lambda)u\|^{2}=\|(K_{\pm}-i\lambda)\Sigma u\|^{2}$\,,
  the upper bound of 
\begin{multline*}
\sqrt{2}\left[\langle \lambda\rangle\|u\|+\|(\frac{1}{2}+\mathcal{O}+\Real
L)u\|+\|(\pm p\partial_{q}+\frac{L_{\pm}-L_{\mp}}{2}-i\lambda)
u\|\right]\\
=
\langle \lambda\rangle\|\Sigma u\|+\|(\frac{1}{2}+\mathcal{O}+\Real
L)\Sigma u\|+\|(\pm p\partial_{q}+\frac{L_{\pm}-L_{\mp}}{2}-i\lambda)\Sigma u\|
\end{multline*}
and $\|\Sigma u\|_{\mathcal{Q}_{0}}$
are assumed in Hypothesis~\ref{hyp.whole}.\\
The upper bound for $\langle
\lambda\rangle^{\frac{1}{4}}\|u\|_{L^{2}(\rz_{-},dq;\mathcal{H}^{1,L})}$
comes from 
$$
\|v\|_{L^{2}(\rz_{-},dq;\mathcal{H}^{1,L})}^{2}\leq
\|v\|\|(\frac{1}{2}+\mathcal{O}+\Real L)v\|\,.
$$
The upper bound for $\langle
\lambda\rangle^{\frac{1}{4}}\|\gamma
u\|_{L^{2}(\rz,|p|dp;\mathfrak{L})}$ is proved in
Proposition~\ref{pr.trL} below.\\
For the density of $\left\{u\in
  L^{2}(\rz^{2}_{-},dqdp;\mathfrak{L})\,, \Sigma u\in
  \mathcal{D}\right\}$\,, we start from the relations
\begin{equation*}
 \tilde{\Sigma}\Sigma=\Sigma\quad,\quad
K_{\pm}^{L}\Sigma =\Sigma K_{\pm,0}^{L}\,.
\end{equation*}
When $u\in D(K_{\pm,0}^{L})$\,, the symmetrized function $\Sigma u$
belongs to $D(K_{\pm}^{L})$ and can be approximated by a sequence
$(u_{n})_{n\in\nz}$ in $\mathcal{D}$\,.\\
With $\tilde{\Sigma}\mathcal{D}=\mathcal{D}$\,, we deduce
\begin{eqnarray*}
  && 0=\lim_{n\to
    \infty}\|\frac{1+\tilde{\Sigma}}{2}u_{n}-\frac{1+\tilde{\Sigma}}{2}\Sigma
  u\|=\lim_{n\to
    \infty}\|\frac{1+\tilde{\Sigma}}{2}u_{n}-\Sigma
  u\|\,,\\
&& 
0=\lim_{n\to \infty}
\|\frac{1+\tilde{\Sigma}}{2}K_{\pm}^{L}u_{n}-\frac{1+\tilde{\Sigma}}{2}\Sigma
K_{\pm,0}^{L}u\|=\lim_{n\to
  \infty}\|\frac{1+\tilde{\Sigma}}{2}K_{\pm}^{L}u_{n}
-\Sigma K_{\pm,0}^{L}u\|\,.
\end{eqnarray*}
The function $v_{n}=1_{\rz_{-}}(q)\frac{1+\tilde{\Sigma}}{2}u_{n}$
satisfies
$$
\Sigma v_{n}=\frac{1+\tilde{\Sigma}}{2}u_{n}\in \mathcal{D}
$$
and the commutation
$P_{\pm}^{L}\tilde{\Sigma}=\tilde{\Sigma}P_{\pm}^{L}$ implies
$$
K_{\pm}^{L}\Sigma
v_{n}=K_{\pm}^{L}\left[\frac{1+\tilde{\Sigma}}{2}u_{n}\right]
=\frac{1+\tilde{\Sigma}}{2}K_{\pm}^{L}u_{n}\,.
$$
We deduce $v_{n}\in D(K_{\pm,0}^{L})$\,, $\Sigma
K_{\pm,0}^{L}v_{n}=\frac{1+\tilde{\Sigma}}{2}K_{\pm}^{L}u_{n}$
and
$$
\lim_{n\to \infty}\|v_{n}-u\|+\|K_{\pm,0}^{L}(v_{n}-u)\|=\frac{1}{2}
\lim_{n\to \infty}\|\Sigma v_{n}-\Sigma u\|+\|\Sigma K_{\pm,0}^{L} (v_{n}-u)\|=0\,,
$$
while we have checked $\Sigma v_{n}\in \mathcal{D}$ for all $n\in\nz$\,.
\end{proof}
With \eqref{eq.equiKL0} the estimate of $\langle \lambda\rangle^{\frac
1 4}\|\gamma u\|_{L^{2}(\rz,|p|dp;\mathfrak{L})}$
is a variant of Proposition~\ref{pr.trace1D} and
 the proof follows the same lines.
\begin{proposition}
\label{pr.trL}
Assume Hypothesis~\ref{hyp.whole}. Any
$u\in D(K_{\pm}^{L})$ has a trace $\gamma u=u(0,p)\in
L^{2}(\rz,|p|dp; \mathfrak{L})$ with the estimate
$$
\forall \lambda\in\rz\,,\quad
\langle
\lambda\rangle^{\frac{1}{2}}\|u(q=0,p)\|_{L^{2}(\rz,|p|dp;
  \mathfrak{L})}^{2}\leq C'\|(K_{\pm}^{L}-i\lambda)u\|^{2}\,.
$$
When $u\in L^{2}(\rz^{2},dqdp; \mathfrak{L})$ solves $(P_{\pm}^{L}-i\lambda)u=\gamma(p)\delta_{0}(q)$ in
$\mathcal{S}'(\rz^{2}; D(L^{*})')$
with 
$\gamma\in L^{2}(\rz,\frac{dp}{|p|}; \mathfrak{L})$ and $\lambda\in\rz$\,, then 
$$
\langle\lambda\rangle^{\frac 1 2}\|u\|^{2}\leq
C''\|\gamma\|_{L^{2}(\rz,\frac{dp}{|p|}; \mathfrak{L})}^{2}\,.
$$
\end{proposition}
\begin{proof}
  We first take $u\in D(K_{\pm})\otimes^{alg}D(L_{\pm})$ and $\lambda\in\rz$\,. We
  compute
  \begin{eqnarray*}
    \int_{\rz}|p||u(0,p)|_{\mathfrak{L}}^{2}~dp
&=&
\pm 2\Real\int_{-\infty}^{0}\langle (\pm p\partial_{q}+i\Imag
L-i\lambda)u\,,\sign(p)u\rangle_{\mathfrak{L}}~dp\\
&=& \pm 2\Real\langle (i\Imag K_{\pm}^{L}-i\lambda)u\,,\,
1_{\rz_{-}}(q)\sign(p)u\rangle\\
&\leq&2\|(i\Imag K_{\pm}^{L}-i\lambda)u\|\|u\|\leq C\langle
\lambda\rangle^{-\frac 1 2}\|(K_{\pm}^{L}-i\lambda)u\|^{2}\,.
  \end{eqnarray*}
The results extends to any $u\in D(K_{\pm}^{L})$ by density.\\
A similar result is valid for $(K_{\pm}^{L})^{*}=K_{\mp}^{L}$\,, while
$\mathcal{S}(\rz^{2}; D(L_{\mp}))$ is a core for $K_{\mp}^{L}$\,.
Thus assuming
$$
\forall v\in \mathcal{S}(\rz^{2}; D(L_{\mp}))\,, 
\langle u\,,\, (K_{\mp}^{L}+i\lambda)v \rangle=\int_{\rz}\langle
\gamma(p)\,,\, v(0,p)\rangle_{ \mathfrak{L}}~dp\,,
$$
leads to
\begin{eqnarray*}
\left|\langle u\,,\, (K_{\mp}^{L}+i\lambda)v\rangle\right|
&\leq&
\||p|^{-\frac 1 2}\gamma\|_{L^{2}(\rz,dp; \mathfrak{L})}\||p|^{\frac 1 2}\gamma
v\|_{L^{2}(\rz,dp; \mathfrak{L})}
\\
&\leq& \|\gamma\|_{L^{2}(\rz,\frac{dp}{|p|}; \mathfrak{L})}\|\gamma
v\|_{L^{2}(\rz,|p|dp; \mathfrak{L})}\\
&\leq & C\langle \lambda\rangle^{-\frac 1
  4}\|\gamma\|_{L^{2}(\rz,\frac{dp}{|p|}; \mathfrak{L})}\|(K_{\mp}^{L}+i\lambda)v\|\,, 
\end{eqnarray*}
for all $v\in \mathcal{S}(\rz^{2}; D(L_{\mp}))$\,.
Since $\mathcal{S}(\rz^{2}; D(L_{\mp}))$ is dense in
$D(K_{\mp}^{L})$ and $(K_{\mp}^{L}+i\lambda)$ is an
isomorphism from $D(K_{\mp}^{L})$ to $L^{2}(\rz^{2},dqdp;
\mathfrak{L})$\,, we deduce
$$
\forall f\in L^{2}(\rz^{2},dqdp; \mathfrak{L})\,,\quad
\left|\langle u\,,\, f\rangle\right|\leq C\langle
\lambda\rangle^{-\frac 1 4} \|\gamma\|_{L^{2}(\rz,\frac{dp}{|p|}; \mathfrak{L})}\|f\|\,,
$$
which proves the estimate $\langle
\lambda\rangle^{\frac{1}{2}}\|u\|^{2}\leq
C''\|\gamma\|^{2}_{L^{2}(\rz,\frac{dp}{|p|};\mathcal{L})}$\,.
\end{proof}
\begin{remark}
  \label{re.nogo}
For any $\lambda\in\rz$\,, the only possible solution to 
$$
(P_{\pm}-i\lambda)u= 0\quad,\quad \gamma u \in \Pi_{odd}L^{2}(\rz,|p|dp;\mathfrak{L})
$$ 
is the trivial solution $u=0$\,. Actually, 
if $u\in L^{2}(\rz_{-},dq;\mathcal{H}^{1,L})$ is a solution, it
has a trace $\gamma\in L^{2}(\rz,|p|dp;D(L_{\mp})')$ and the function
$\tilde{u}(q,p)=u(q,p)1_{\rz_{-}}(q)$ solves
$$
(P_{\pm}-i\lambda)\tilde{u}=-p \gamma \delta_{0}(q)\quad
\text{in}~\mathcal{S}'(\rz^{2};D(L_{\mp})')\,.
$$
By Proposition~\ref{pr.extsym}, we deduce $\Sigma u
=2\tilde{u}$ and $\tilde{u}\big|_{\left\{q>0\right\}}\equiv 0$ implies
$\Sigma u=0$ and $u=0$\,.
\end{remark}

\subsection{Resolvent estimates}
\label{se.resestimL}

\subsubsection{Trace estimates}
\label{se.trestL}
\begin{proposition}
\label{pr.esttrL}
  Assume Hypotheses
  \ref{hyp.whole} and \ref{hyp.Abdd} with $A\neq 0$\,.
There exists a constant $C_{1,L,A}>0$ such that the
  estimate
$$
\forall f\in
L^{2}(\rz^{2}_{-},dqdp; \mathfrak{L})\,,\quad
\|\gamma (K^{L}_{\pm,A}-i\lambda)^{-1}f\|_{L^{2}(\rz,|p|dp;
  \mathfrak{L})}\leq
 C_{1,L,A}\langle \lambda\rangle^{-\frac 1 4}\|f\|\,,
$$
holds uniformly w.r.t $\lambda\in\rz$\,.
\end{proposition}
\begin{proof}
  For $f\in L^{2}(\rz^{2}_{-},dqdp;\mathfrak{L})$\,, the function
  $u=(K^{L}_{\pm,A}-i\lambda)^{-1}f\in D(K^{L}_{\pm,A})$ 
 solves
 \begin{eqnarray*}
&&
(P_{\pm}^{L}-i\lambda)u=f\quad,\quad \gamma_{odd}u=f_{\partial}
\\
\text{with}
&&f_{\partial}=\gamma_{odd}u=\pm\sign(p)A\gamma_{ev}u\quad,\quad
\gamma u \in L^{2}(\rz,|p|dp;\mathfrak{L})\,.
\end{eqnarray*}
The difference $v=u-(K_{\pm,0}-i\lambda)^{-1}f$ satisfies $v\in
L^{2}(\rz_{-},dq;\mathcal{H}^{1,L})$\,, 
$\gamma v\in L^{2}(\rz,|p|dp;\mathfrak{L})$ and solves
$$
(P_{\pm}-i\lambda)v=0\quad,\quad \gamma_{odd}v=f_{\partial}\,.
$$
Proposition~\ref{pr.TrScal} then says
$$
\pm\Real\langle \gamma_{ev}v\,,\,
\sign(p)\gamma_{odd}v\rangle_{L^{2}(\rz,|p|dp;\mathfrak{L})}
=
\Real\langle \gamma_{ev}v\,,\, \pm\sign(p)f_{\partial}\rangle_{L^{2}(\rz,|p|dp;\mathfrak{L})}\leq 0\,.
$$
The even part of the trace at $q=0$ of $u=(K_{\pm,0}^{L}-i\lambda)^{-1}f+v$ equals
\begin{equation*}
\gamma_{ev}u=\gamma_{ev}(K^{L}_{\pm,0}-i\lambda)^{-1}f +\gamma_{ev}v\,,
\end{equation*}
Taking the $L^{2}(\rz, |p|dp;\mathfrak{L})$-scalar product 
with $\pm\sign(p)f_{\partial}=A\gamma_{ev}u$
says that the quantity
$$
\Real\langle 	A \gamma_{ev}u\,,\, \gamma_{ev}u\rangle_{L^{2}(\rz,|p|dp;\mathfrak{L})}
=
\Real\langle \pm \sign(p)f_{\partial}\,,\, \gamma_{ev}u\rangle_{L^{2}(\rz,
  |p|dp;\mathfrak{L})}
$$
cannot be larger than
\begin{multline*}
\Real\langle \pm\sign(p)f_{\partial}\,,\,
\gamma_{ev}(K_{\pm,0}^{L}-i\lambda)^{-1}f\rangle_{L^{2}(\rz,|p|dp;\mathfrak{L})}
\\
\leq \|\gamma_{odd}u\|_{L^{2}(\rz,|p|dp;\mathfrak{L})}
\|\gamma (K^{L}_{\pm,0}-i\lambda)^{-1}f\|_{L^{2}(\rz,|p|dp;\mathfrak{L})}\,.
\end{multline*}
With Hypothesis~\ref{hyp.Abdd} and $A\neq 0$, we deduce
\begin{eqnarray*}
  \frac{c_{A}}{\|A\|^{2}}\|\gamma_{odd}u\|^{2}_{L^{2}(\rz,|p|dp;\mathfrak{L})}&=&\frac{c_{A}}{\|A\|^{2}}\|A
  \gamma_{ev}u\|_{L^{2}(\rz,|p|dp;\otimes\mathfrak{L})}^{2}
\\
&\leq& \Real\langle
  A\gamma_{ev}u\,, \gamma_{ev}u\rangle_{L^{2}(\rz,|p|dp;\otimes\mathfrak{L})}
\\
&\leq&
 \|\gamma_{odd}u\|_{L^{2}(\rz,|p|dp;\otimes\mathfrak{L})}\|\gamma (K^{L}_{\pm,0}-i\lambda)^{-1}
 f
\|_{L^{2}(\rz,|p|dp;\mathfrak{L})}\,.
\end{eqnarray*}
We infer from this
\begin{eqnarray*}
  &&\|\gamma_{odd}u\|_{L^{2}(\rz,|p|dp;\mathfrak{L})}\leq 
\frac{\|A\|^{2}}{c_{A}}\|\gamma (K^{L}_{\pm,0}-i\lambda)^{-1}
f\|_{L^{2}(\rz,|p|dp;\mathfrak{L})}\,,\\
\text{and}
&&
\|\gamma_{ev}u\|_{L^{2}(\rz,|p|dp;\mathfrak{L})}\leq
\frac{\|A\|}{c_{A}}
\|\gamma (K^{L}_{\pm,0}-i\lambda)^{-1}
f\|_{L^{2}(\rz,|p|dp;\mathfrak{L})}\,,
\end{eqnarray*}
and
we conclude by referring to
Proposition~\ref{pr.trKL0}.
\end{proof}

\subsubsection{$L^2$-estimates}
\label{se.L2L}

\begin{proposition}
\label{pr.L2L} Assume Hypotheses~\ref{hyp.whole} and \ref{hyp.Abdd}
with $A\neq 0$. 
There exists a constant $C_{2,L,A}>0$ such that 
$$
\forall f\in L^{2}(\rz^{2}_{-},dqdp;
\mathfrak{L})\,,\quad
\|(K^{L}_{\pm,A}-i\lambda)^{-1}f\|\leq
C_{2,L,A}\langle\lambda\rangle^{-\frac 1 2}\|f\|\,,
$$
holds uniformly w.r.t $\lambda\in\rz$\,.
\end{proposition}
\begin{proof}
Consider $u=(K^{L}_{\pm,A}-i\lambda)^{-1}$ with $f\in
L^{2}(\rz^{2}_{-},dqdp;\mathfrak{L})$ and decompose it into
$$
u=(K^{L}_{\pm,0}-i\lambda)^{-1}f+v\,.
$$
The
first term satisfies
$$
\Sigma(K^{L}_{\pm,0}-i\lambda)^{-1}f=(K_{\pm}^{L}-i\lambda)(\Sigma f)
$$
according to \eqref{eq.equiKL0}.
We deduce
$$
\|(K^{L}_{\pm,0}-i\lambda)^{-1}f\|\leq
\frac{1}{\sqrt{2}}\|(K_{\pm}^{L}-i\lambda)^{-1}(\Sigma f)\|
\leq \frac{C}{\sqrt{2}}\langle \lambda\rangle^{-\frac{1}{2}}\|f\|\,.
$$
The second term $v$ belongs to $L^{2}(\rz_{-},dq;\mathcal{H}^{1})$
with $\gamma v\in L^{2}(\rz,|p|dp;\mathfrak{L})$ and solves
$$
(P_{\pm}-i\lambda)v=0\quad,\quad \gamma_{odd}v=\gamma_{odd}u\,.
$$
Proposition~\ref{pr.extsym} says that $\Sigma v$ solves
$$
(P_{\mp}-i\lambda)\Sigma v=\mp 2p (\gamma_{odd}u)\delta_{0}(q)
$$
and Proposition~\ref{pr.trL} implies
\begin{eqnarray*}
  \|v\|^{2}&\leq&\frac{1}{2}\|\Sigma v\|^{2}\leq \frac{C''}{2}\langle
  \lambda\rangle^{-\frac{1}{2}}
\|p\gamma_{odd}u\|_{L^{2}(\rz,\frac{dp}{|p|};\mathfrak{L})}^{2}\\
&\leq&\frac{C''}{2}\langle
  \lambda\rangle^{-\frac{1}{2}}
\|\gamma u\|_{L^{2}(\rz,|p|dp;\mathfrak{L})}^{2}
\leq \frac{C'' C'}{2}\langle
  \lambda\rangle^{-1}\|f\|^{2}\,.
\end{eqnarray*}
\end{proof}

\subsubsection{Regularity estimates}
\label{se.regestL}
\begin{proposition}
\label{pr.regL1}
  Assume Hypotheses~\ref{hyp.whole} and \ref{hyp.Abdd} with $A\neq 0$\,. 
There exists a constant $C_{3,A,L}>0$ such that
$$
\langle \lambda\rangle\|u\|^{2}+\langle \lambda\rangle^{\frac 1
  2}\|u\|_{L^{2}(\rz_{-},dq;\mathcal{H}^{1,L})}^{2}
+\langle\lambda\rangle^{\frac 1 2}\|\gamma u\|_{L^{2}(\rz,|p|dp; \mathfrak{L})}^{2}
\leq C_{3,A,L}\|(K^{L}_{\pm,A}-i\lambda)u\|^{2}\,,
$$
holds for all $u\in D(K^{L}_{\pm,A})$ and all $\lambda\in\rz$\,.\\
If $\Phi\in \mathcal{C}^{\infty}_{b}((-\infty,0])$ satisfies
$\Phi(0)=0$\,, there exists $C_{L,\Phi}>0$ such that 
$$
\|\Phi(q)[\frac{1}{2}+\mathcal{O}+\Real L]u\|\leq
C\|\Phi\|_{L^{\infty}}\|(K^{L}_{\pm, A}-i\lambda)u\|
+C_{L,\Phi}\|u\|\,,
$$
holds true when $\lambda\in\rz$\,, $u\in D(K^{L}_{\pm,A})$ and 
$C$ is the constant of Hypothesis~\ref{hyp.whole}.
\end{proposition}
\begin{proof}
  \noindent\textbf{a)} The bound for $\langle \lambda\rangle\|u\|^{2}$
  and $\langle \lambda\rangle^{1/2}\|\gamma_{odd}u\|_{L^{2}(\rz,|p|dp;\mathfrak{L})}^{2}$
  are given in Proposition~\ref{pr.L2L} and Proposition~\ref{pr.esttrL}.\\
\noindent\textbf{b)} For $u\in D(K^{L}_{A})$\,, the integration by
part indentity of Proposition~\ref{pr.maxaccKAL}\,, combined with
Hypothesis~\ref{hyp.Abdd} ($A\neq 0$)\,, gives
\begin{eqnarray*}
\|u\|\|(K^{L}_{\pm,A}-i\lambda)u\|
&\geq &\Real\langle u\,,\, K^{L}_{\pm,A}u\rangle
\\
&\geq
&\|u\|_{L^{2}(\rz_{-},dq;\mathcal{H}^{1,L})}^{2}
+\Real\langle \gamma_{ev}u\,,\, A\gamma_{ev}u\rangle_{L^{2}(\rz,|p|dp;
  \mathfrak{L})}^{2}\\
&\geq&
\|u\|_{L^{2}(\rz_{-},dq;\mathcal{H}^{1,L})}^{2}
+\frac{c_{A}}{2}\|\gamma_{ev}u\|_{L^{2}(\rz,|p|dp;
  \mathfrak{L})}^{2}\\
&&\hspace{4cm}
+\frac{c_{A}}{2\|A\|^{2}}\|\gamma_{odd}u\|_{L^{2}(\rz,|p|dp;
  \mathfrak{L})}^{2}\,.
\end{eqnarray*}
With \textbf{a)} or Proposition~\ref{pr.L2L}, we deduce
\begin{eqnarray*}
\|(K^{L}_{\pm,A}-i\lambda)u\|^{2}
&\geq& C_{2,A,L}^{-1}\langle \lambda\rangle^{1/2}\|u\|\|(K^{L}_{\pm,A}-i\lambda)uu\|
\\
&\geq&
C_{A,L}^{-1}\langle\lambda\rangle^{1/2}\left[\|u\|_{L^{2}(\rz_{-},dq;\mathcal{H}^{1,L})}^{2}
+\|\gamma u\|_{L^{2}(\rz,|p|dp; \mathfrak{L})}^{2}\right]\,.
\end{eqnarray*}
\noindent\textbf{c)} For  $\Phi\in
\mathcal{C}^{\infty}_{b}(\overline{\rz_{-}})$ such that $\Phi(0)=0$\,,
$\lambda\in\rz$ 
and $u\in D(K^{L}_{\pm,A})$ with $(K^{L}_{\pm,A}-i\lambda)u=f$\,, write:
$$
(P_{\pm}^{L}-i\lambda)(\Phi(q)u)=\Phi(q)(P_{\pm}^{L}-i\lambda)u\pm p\Phi'(q)u=\Phi(q)f\pm\Phi'(q)(pu)\,.
$$
The assumption $\Phi(0)=0$ implies 
\begin{eqnarray*}
  &&(\Phi(q)u) 1_{\rz_{-}}(q)\in
D(K^{L}_{\pm})\\
\text{and}&&
(K_{\pm}^{L}-i\lambda)(\Phi(q)u
1_{\rz_{-}}(q))=1_{\rz_{-}}(q)\left[\Phi(q)f\pm
\Phi'(q)(pu)\right]\,.
\end{eqnarray*}
Hypothesis~\ref{hyp.whole} implies
\begin{eqnarray*}
  \|(\frac{1}{2}+\mathcal{O}+\Real L)\Phi(q)u\|&=&
\|\frac{K_{\pm}^{L}+K_{\mp}^{L^{*}}}{2}(\Phi(q)u 1_{\rz_{-}}(q))\|\\
&\leq&
C\|1_{\rz_{-}}(q)\left[\Phi(q)f+\Phi'(q)(pu)\right]\|
\\
&\leq&
CM^{0}_{\Phi}\|(K^{L}_{\pm,A}-i\lambda)u\| + M^{1}_{\Phi}\|pu\|
\\
&\leq& CM^{0}_{\Phi}\|(K^{L}_{\pm,A}-i\lambda)u\| 
+ 2M^{1}_{\Phi}\|u\|_{L^{2}(\rz_{-},dq;\mathcal{H}^{1,L})}\,,
\end{eqnarray*}
where we have used $M^{k}_{\Phi}=\|\partial_{q}^{k}\Phi\|_{L^{\infty}}$\,.
Again from 
\begin{equation*}
  \|u\|\|(K^{L}_{\pm,A}-i\lambda)u\|
\geq 
\Real\langle u\,,\, K^{L}_{\pm,A}u\rangle
\geq \|u\|_{L^{2}(\rz_{-},dq;\mathcal{H}^{1,L})}^{2}\,,
\end{equation*}
we deduce
$$
 \|(\frac{1}{2}+\mathcal{O}+\Real L)\Phi(q)u\|
\leq 2CM^{0}_{\Phi}\|(K^{L}_{\pm,A}-i\lambda)u\|+ C_{L,\Phi}\|u\|\,,
$$
by choosing simply $C_{L,\Phi}=\frac{2(M^{1}_{\Phi})^{2}}{CM^{0}_{\Phi}}$\,.
\end{proof}
\begin{proposition}
\label{pr.regq} 
Assume Hypotheses~\ref{hyp.whole}, \ref{hyp.Abdd} with $A\neq 0$ and Hypothesis~\ref{hyp.Q}.
Then for any $u\in D(K^{L}_{\pm,A})$\,, $\Sigma u$ belongs to
$\mathcal{Q}$ with the estimate
$$
\forall \lambda\in\rz\,, \forall u\in D(K^{L}_{\pm,A})\,,\quad
\langle \lambda\rangle^{\frac{1}{4}}\|\Sigma u\|_{\mathcal{Q}}\leq C_{\mathcal{Q},A,L}\|(K^{L}_{\pm,A}-i\lambda)u\|
$$
\end{proposition}
\begin{proof}
  Again we write that $\Sigma u=\Sigma(K^{L}_{\pm,A}-i\lambda)^{-1}f$ satisfies
$$
(P_{\pm}^{L}-i\lambda)u=\Sigma f\mp 2p(\gamma_{odd}u)\delta_{0}(q)\,,
$$
with $2p(\gamma_{odd}u)\in
L^{2}(\rz,\frac{dp}{|p|};\mathfrak{L})$\,.
According to
Proposition~\ref{pr.trL}\,, we call $v\in L^{2}(\rz,dq;\mathcal{H}^{1,L})$ the solution to
\eqref{eq.PLgd} with $\gamma=\mp 2p(\gamma_{odd}u)$ so that
$$
\Sigma u=(K_{\pm}^{L}-i\lambda)^{-1}(\Sigma f)+v\,.
$$
The first assumed inequality \eqref{eq.hypQ1} in Hypothesis~\ref{hyp.Q} gives
$$
\|v\|_{\mathcal{Q}}\leq
C\|\gamma_{odd}u\|_{L^{2}(\rz,|p|dp;\mathfrak{L})}
$$
while Proposition~\ref{pr.esttrL} implies 
$$
\|\gamma_{odd}u\|_{L^{2}(\rz,|p|dp;\mathfrak{L})}\leq C_{1,L,A}\langle \lambda\rangle^{-\frac{1}{4}}\|f\|\,.
$$
The second assumed inequality \eqref{eq.hypQ2} in
Hypothesis~\ref{hyp.Q}\,,
$$
\|(K_{\pm}^{L}-i\lambda)^{-1}(\Sigma
f)\|_{\mathcal{Q}}
\leq C\langle \lambda\rangle^{-1/4}\|f\|\,,
$$
yields the result.
\end{proof}
\section{Geometric Kramers-Fokker-Planck operator}
\label{se.geoKFP}

In \cite{Leb1}cite{Leb2}, G. Lebeau developed the resolvent analysis for
geometric (Kramers)-Fokker-Planck operators, which is the scalar model
for the hypoelliptic Laplacian introduced by J.M.~Bismut (see 
\cite{BiLe}\cite{Bis05}\cite{BisSurv}) and which
will be considered in Section~\ref{se.hypolapl}. Lebeau in
\cite{Leb1}\cite{Leb2} studied the case of compact manifolds but
a partition of unity allows to extend it to the case of cylinders
$Q=\rz\times Q'$ when $Q'$ is a compact manifold. We will check that
Hypothesis~\ref{hyp.whole} and Hypothesis~\ref{hyp.Q} are satisfied by
the geometric Kramers-Fokker-Planck operator on $T^{*}Q$ when
$Q=\rz\times Q'$ is endowed with a specific splitted metric.
This will prepare the analysis of half-cylinders models for the
boundary case. We also review some partition of unity arguments which
will be used more extensively in Section~\ref{se.multi} and
extract the usefull information from the results of \cite{Leb2}.\\
For this presentation we stick to the case of the scalar geometric
Kramers-Fokker-Planck equations. Tensorizing with a Hilbert space
$\mathfrak{f}$ makes no difficulty in the end.

\subsection{Notations and the geometric KFP-operator}
\label{se.notgFKP}

We consider a $d$-dimensional Riemannian manifold $(Q,g)$ (possibly
with a boundary $\partial Q$) endowed with
the metric
$g(q)=(g_{ij}(q))_{1\leq i,j\leq d}$\,,
$$
g(q;T)=\sum_{i,j=1}^{d}g_{ij}(q)T^{i}T^{j}=g_{ij}(q)T^{i}T^{j}=T^{T}g(q)T\,,\quad
T\in T_{q}Q\,,
$$
with Einstein summation convention and the matricial writing for the
last right-hand side. The $q$ coordinates are written
$(q^{1}\,\ldots,q^{d})\in\rz^{d}$ and when we use the matricial
notation the vectors are column-vectors. The inverse tensor $g^{-1}(q)$ is
denoted by $g^{-1}(q)=(g^{ij}(q))_{1\leq i,j\leq d}$ with
$$
g^{ik}g_{kj}=\delta^{i}_{j}\,.
$$
With the metric $g$ are associated the Christoffel symbols
$$
\Gamma^{\ell}_{ij}=\frac{1}{2}g^{\ell k}\left\{\partial_{q^{i}}g_{kj}+\partial_{q^{j}}g_{ik}-\partial_{q^{k}}g_{ij}\right\}=\Gamma^{\ell}_{ji}\,,
$$
the Levi-Civita connection
\begin{equation}
  \label{eq.LevCiv}
\nabla_{\partial_{q^{i}}}\left(\partial_{q^{j}}\right)=\Gamma^{\ell}_{ij}\partial_{q^{\ell}}=(\Gamma_{i})^{\ell}_{j}\partial_{q^{j}}\,,
\end{equation}
with its adjoint version
\begin{equation}
  \label{eq.adjLevCiv}
  \nabla_{\partial_{q^{i}}}\left(dq^{j}\right)=-\Gamma^{j}_{i\ell}dq^{\ell}\,,
\end{equation}
and the Riemann curvature tensor, which is an $\text{End}(TQ)$-valued
two form\,,
\begin{eqnarray*}
  &&R=R_{jk}dq^{j}\wedge dq^{k}
\\
&& R_{jk}=\frac{\partial \Gamma_{j}}{\partial q^{k}}-\frac{\partial
  \Gamma_{k}}{\partial q^{j}}+\left[\Gamma_{j}\,,\, \Gamma_{k}\right]\,.
\end{eqnarray*}
The geodesic flow in the tangent space $TQ$ is obtained after
minimizing
$$
\frac{1}{2}\int_{0}^{1}g(q;\dot \sigma)~dt
$$
with respect to $\sigma\in \mathcal{C}^{2}([0,1];Q)$ and is locally
equivalent to 
$$
{\ddot\sigma}^{\ell}+\Gamma^{\ell}_{ij}\dot{\sigma}^{i}\dot{\sigma}^{j}=0\,.
$$
This is the Lagrangian version of the geodesic flow which has the
Hamiltonian counterpart (see \cite{AbMa} for example).
The cotangent bundle over $Q$ is denoted by $X=T^{*}Q$\,. The length
of a covector $p_{i}dq^{i}\in T^{*}_{q}Q$
 is given by $|p|_{q}^{2}=g^{ij}(q)p_{i}p_{j}=p^{T}g^{-1}(q)p$\,. The Hamiltonian version
 of the geodesic flow is given by the Hamiltonian vector field $\mathcal{Y}_{\mathcal{E}}$
 associated with the energy
 $\mathcal{E}(q,p)=\frac{1}{2}|p|_{q}^{2}$\,.
With the
symplectic form $dp_{i}\wedge dq^{i}$\,, it is given by
$$
\mathcal{Y}_{\mathcal{E}}=p^{T}g^{-1}(q)\partial_{q}-\frac{1}{2}\partial_{q}\left[p^{T}g^{-1}(q)p\right]^{T}\partial_{p}=
g^{ij}(q)p_{i}\partial_{q^{j}}-\frac{1}{2}[\partial_{q^{k}}g^{ij}(q)]p_{i}p_{j}\partial_{p_{k}}
\,.
$$
Owing to the relation 
$$
-\frac{1}{2}\partial_{q^{k}}g^{ij}p_{i}p_{j} = g^{i\ell}\Gamma^{j}_{k\ell}p_{i}p_{j}
$$
it also equals
 \begin{eqnarray*}
   &&\mathcal{Y}_{\mathcal{E}}=g^{ij}(q)p_{i}e_{j}=p^{T}g^{-1}(q)e\,,\\
\text{with}&&e_{j}=\partial_{q^{j}}+\Gamma^{\ell}_{i,j}p_{\ell}\partial_{p_{i}}\quad,\quad
e=
\begin{pmatrix}
  e_{1}\\\vdots\\e_{d}
\end{pmatrix}\,.
 \end{eqnarray*}
When $q\to \tilde{q}=\Phi(q)$ is a change of coordinate in a neighborhood of
$q_{0}\in Q$\,, the
corresponding symplectic change of coordinates in $X=T^{*}Q$ is given
by 
$$
(\tilde{q},\tilde{p})=(\Phi(q),[{}^{t}d\Phi(q)]^{-1}p)\,,
$$
while the change of metric is given by 
$$
\tilde{g}(\tilde{q})=\left[d\Phi(q)\right]^{-1,T}g(q)\left[d\Phi(q)\right]^{-1}\quad,\quad
\left(\tilde{g}(\tilde{q})\right)^{-1}=[d\Phi(q)]g(q)^{-1}[d\Phi(q)]^{T}\,.
$$
The function $\mathcal{E}(q,p)=|p|_{q}^{2}/2$ and the Hamiltonian
vector field, $\mathcal{Y}_{\mathcal{E}}$\,, on $T^{*}Q$ do not depend on the choice of coordinates
on $Q$\,.\\
The manifold $X$ is endowed with the riemannian metric
$g^{X}=g\mathop{\oplus}^{\perp}g^{-1}$ where the orthogonal
decomposition follows the decomposition into the horizontal and
vertical spaces of $T_{x}X=(T_{x}X)^{H}\oplus (T_{x}X)^{V}$ with
$x=(q,p)$\, $(T_{x}X)^{H}\sim T_{q}Q$ and $(T_{x}X)^{V}\sim
T^{*}_{q}Q$\,. More precisely at a point $x=(q,p)$\,, $g^{X}(x)$ is
determined  by
$$
g^{X}(x;e_{i},e_{j})=g_{ij}(q)\quad,\quad
g^{X}(x;\partial_{p_{i}},\partial_{p_{j}})=g^{ij}(q)\quad\text{and}\quad
g^{X}(x;e_{i},\partial_{p_{j}})=0\,.
$$
\begin{remark}
\label{re.sign}
  One may be puzzled by the sign in the expression
  $e_{j}=\partial_{q^{j}}+\Gamma^{\ell}_{ij}p_{\ell}\partial_{p_{j}}$
  which does not coincide with the general definition of horizontal
  vector fields on a vector bundle $\pi:F\to Q$ endowed with an affine
  connection $\nabla^{F}$\,. Actually we follow the convention of
  G.~Lebeau and the horizontal tangent vectors on
  the \underline{cotangent} bundle, $X=T^{*}Q$\,, are expressed with the Levi-Civita
  connection on the
  \underline{tangent} bundle $TQ$ formulated with the Christoffel symbols
  $\Gamma^{\ell}_{ij}$\,. Compare in particular the adjoint formula
  \eqref{eq.adjLevCiv} with  \eqref{eq.LevCiv}.
\end{remark} 
The
 vertical Laplacian is given by
 $\Delta_{p}=\partial_{p_{i}}g_{ij}(q)\partial_{p_{j}}=\partial_{p}^{T}g(q)\partial_{p}$\,,
 which has also an invariant meaning.
 \begin{definition}
\label{de.geKFP}
   The geometric Kramers-Fokker-Planck operator on $X=T^{*}Q$ is the
   differential operator
$$
P_{\pm,Q,g}=\pm \mathcal{Y}_{\mathcal{E}}+\frac{1}{2}\left[-\Delta_{p}+|p|_{q}^{2}\right]\,.
$$
The notation $\mathcal{O}_{Q,g}$ will be used for
the operator acting along the fiber $T^{*}_{q}Q$: 
$$
\frac{1}{2}\left[-\Delta_{p}+|p|_{q}^{2}\right]=
\frac{1}{2}\left[-\partial_{p}^{T}g(q)\partial_{p}+p^{T}g^{-1}(q)p\right]\,.
$$
\end{definition}
Note that the vertical operator $\mathcal{O}_{Q,g}$ is self-adjoint with
the domain $D(\mathcal{O}_{Q,g})=\left\{u\in L^{2}(X)\,,\,
  \mathcal{O}_{Q,g}u\in L^{2}(X)\right\}$ and it satisfies the
uniform lower bound $\mathcal{O}_{Q,g}\geq \frac{\dim Q}{2}$\,.\\
We recall the notations introduced after \eqref{eq.Hsq} in the introduction.
\begin{definition}
Let $(\overline{Q},g)$\,, $\overline{Q}=Q\sqcup \partial Q$\,, be a riemannian manifold
with a possibly  non empty boundary $\partial Q$.\\
For any $s'\in\rz$\,, the space
$$
\mathcal{H}^{s'}(q)=\left\{u\in
  \mathcal{S}'(\rz^{d},dp)\,,\quad
  (\frac{d}{2}+\mathcal{O}_{g(q)})^{s'/2}u\in
  L^{2}(\rz^{d},dp)\right\}\,,\quad q\in Q\,,
$$
defines a hermitian bundle $\pi_{\mathcal{H}^{s'}}:\mathcal{H}^{s'}\to \overline{Q}$\,.
The space of its $L^{2}$-sections
(resp. $H^{s}$-sections, $s\in\rz$) is denoted by $L^{2}(Q;\mathcal{H}^{s'})$
(resp. $H^{s}(\overline{Q};\mathcal{H}^{s'})$).\\
On $L^{2}(Q;\mathcal{H}^{s'})$ the scalar product is given by 
$$
\langle
u\,,\,v\rangle_{L^{2}(Q;\mathcal{H}^{s'})}=\langle u\,,\,
(\frac{d}{2}+\mathcal{O}_{g})^{s'/2}v\rangle
$$
When $Q$ is compact without boundary,
 the intersection $\ccap_{s,s'\in\rz}H^{s}(Q;\mathcal{H}^{s'})$ is
denoted by $\mathcal{S}(X)$ and endowed with its Fr{\'e}chet space
topology. Its dual is
$\mathcal{S}'(X)$\,.
\end{definition}
When $Q$ is compact without boundary, the space $\mathcal{S}(X)$ is
the space  of rapidly decaying
$\mathcal{C}^{\infty}$-functions (like in $\rz^{2d}$) and one easily
checks (use for example the dyadic partition of unity in the
$p$-variable recalled in Subsection~\ref{se.momentpart}):
\begin{eqnarray*}
  &&\mathcal{S}'(X)=\ccup_{s,s'\in\rz}H^{s}(Q;\mathcal{H}^{s'})\,,\\
&&\mathcal{C}^{\infty}_{0}(X)\subset\mathcal{S}(X)\subset 
H^{s_{1}}(Q;\mathcal{H}^{s_{1}'})\subset
H^{s_{2}}(Q;\mathcal{H}^{s_{2}'})\subset
\mathcal{S}'(X)\subset \mathcal{D}'(X)\,, 
\end{eqnarray*}
when
$s_{1}\geq s_{2}\,, s_{1}'\geq s_{2}'$\,, with dense and continuous embeddings.\\
We close this section by giving the explicit form of $P_{\pm,Q,g}$ in a
specific local coordinate system along an hypersurface $Q'$ of $\overline{Q}$
(which can be the boundary)\,.
In a neighborhood $\mathcal{V}_{q_{0}}$ of $q_{0}\in Q'$ small enough, one can find a
coordinate system $(q^{1},q')$ such that
\begin{eqnarray*}
  && q^{1}(q_{0})=0\quad,\quad q'(q_{0})=0\,,\\
&& \forall q\in \mathcal{V}_{q_{0}}\,,\quad\left(q\in
  Q'\right)\Leftrightarrow\left(q^{1}(q)=0\right)\,,\\
&& \forall q\in \mathcal{V}_{q_{0}}\,,\quad 
g(q)=
\begin{pmatrix}
  1&0\\
0& m_{ij}(q^{1},q')
\end{pmatrix}\,.
\end{eqnarray*}
In such a coordinate system the equalities
\begin{eqnarray*}
  && 2\mathcal{E}=|p|_{q}^{2}=|p_{1}|^{2}+m^{ij}(q^{1},q')p'_{i}p'_{j}\,,\\
&&
\Delta_{p}=\partial_{p_{1}}^{2}+\partial_{p'_{i}}m_{ij}(q^{1},q')\partial_{p'_{j}}\,,\\
&&
\mathcal{Y}_{\mathcal{E}}=p_{1}\partial_{q^{1}}
+ m^{ij}(q^{1},q')p_{i}'\partial_{q^{'j}}
-\frac{1}{2}[\partial_{q^{'k}}m^{ij}(q^{1},q')]p_{i}'p_{j}'\partial_{p_{k}'}\\
&&\hspace{4cm}-\frac{1}{2}[\partial_{q^{1}}m^{ij}(q^{1},q')]p_{i}'p_{j}'\partial_{p_{1}}\,,
\end{eqnarray*}
hold with the convention  $p'^{T}=
(p_{2},\ldots,p_{d})$\,,
$q'^{T}=(q^{2},\ldots,q^{d})$\,, and the corresponding summation rule.
One obtains
\begin{equation}
  \label{eq.gKFPmet1}
P_{\pm,Q,g}=\pm p_{1}\partial_{q^{1}}-\partial_{p_{1}}^{2}+|p_{1}|^{2}+
P_{\pm, Q',m(q^{1},.)}-\frac{1}{2}\partial_{q^{1}}m^{ij}(q^{1},q')p_{i}'p_{j}'\partial_{p_{1}}\,.
\end{equation}
This has some simple interpretation: when $Q'$ is a curved
hypersurface of $Q$\,, which can itself have some intrinsic
curvature\,, 
the block component $(m_{ij}(q^{1},q'))_{2\leq i,j\leq d}$ cannot be
made independent of $q_{1}$\,. In particular the curvature of $Q'$ has
two consequences on the dynamics, a centrifugal (or centripetal) force
and a Coriolis type force. The Coriolis force is a tangential
phenomenon which can be included in $P_{\pm,Q',g(q^{1},.)}$\,. Only the
centrifugal force $\partial_{q^{1}}g^{ij}(q^{1},q')p_{i}p_{j}\partial_{p_{1}}$
and the parametrization of the metric $m(q^{1},.)$ by $q^{1}$\,,
prevent from a separation of variables. This is why, even when $\overline{Q}$ is
a domain of the flat euclidean space, the curvature term has to be analyzed
accurately for boundary value problems. Even in this basic case, it is
useful to formulate the problem in a geometric setting. We shall refer to
the recent results by G.~Lebeau in \cite{Leb1}\cite{Leb2} and
Bismut-Lebeau in \cite{BiLe} where the problems raised by the
curvature term have been accurately analyzed.

\subsection{The result by G.~Lebeau}
\label{se.lebeau}

In \cite{Leb1}\cite{Leb2}, G.~Lebeau noticed that in local
coordinates, the differential operators  involved in
$P_{\pm, Q,g}-z=\pm \mathcal{Y}_{\mathcal{E}}+\mathcal{O}_{Q,g}-\Real z +i\Imag
z$ have a non trivial
homogeneity: The main reason is that a symplectic change of variable 
$(q',p')=(\phi(q), d\phi(q)^{-T}p)$\,, $A^{-T}=(A^{-1})^{T}$\,,
leads to
$dq'=d\phi(q)d_{q}$\,, $dp'=(d\phi)^{-T}dp+(d(d\phi^{-T})p)dq$ and
\begin{eqnarray*}
  &&\partial_{p'}=d\phi(q)\partial_{p}\\
&&\partial_{q'}=(d\phi(q))^{-T}\left[\partial_{q}
+(d(d\phi^{T}))(d\phi)^{-T}p)^{T}\partial_{p}\right]\,.
\end{eqnarray*}
Hence assuming the order $1$ for
$\partial_{q}$ in a invariant way imposes the order $1/2$ for $p$ and
$\partial_{p}$\,.\\
Below is the list of orders for useful operators in a coordinate system:
\begin{eqnarray*}
  && \partial_{q^{j}}\;:\;1\quad; \quad p_{i},\partial_{p_{i}}\;:\;
  \frac{1}{2}\quad;\quad e_{j}=\partial_{q^{j}}+\Gamma_{ij}^{\ell}p_{\ell}\partial_{p_{i}}
\;:\;1\,,\\
&&
\mathcal{Y}_{\mathcal{E}}=g^{ij}(q)p_{i}e_{j}\;,\; \Imag z\;:\;\frac{3}{2}\quad;\quad
 \mathcal{O}_{Q,g}=\frac{-\Delta_{p}+|p|^{2}_{q}}{2}\;,\;\Real z\;:\; 1\,,\\
&&\langle p\rangle\partial_{p_{i}}\;,\; e_{j}\;,\; \langle
p\rangle^{2}+|\Real z|+\frac{|\Imag z|}{\langle
  p\rangle}\;:\;1\,,\quad 
(\langle p\rangle^{2}=1+|p|_{q}^{2})\,.
\end{eqnarray*}
\begin{definition}
\label{de.Wgz}
Assume that $(Q,g)$ is compact or $Q=\rz^{d}$ with $g-\Id_{\rz^{d}}\in
\mathcal{C}^{\infty}_{0}(\rz^{d})$\,. Consider a finite atlas
$Q=\cup_{\ell}^{L}\Omega_{\ell}$ with $\overline{\Omega_{\ell}}$
compact when $\ell>1$ and let
$\sum_{\ell=1}^{L}\theta_{\ell}(q)\equiv 1$\,, be a partition of
unity subordinate to this atlas.\\
With the corresponding  coordinate system $(q^{i},p_{i})_{i=1,\ldots,d}$ on
each $T^{*}\Omega_{\ell}$\,, $\ell\in \left\{1,\ldots, L\right\}$\,, set
$$
\mathfrak{T}_{\ell}=\left\{\langle p\rangle\partial_{p_{i}}\;,\; e_{j}\;,\; \langle
p\rangle^{2}+|\Real z|+\frac{|\Imag z|}{\langle p\rangle}\right\}\,.
$$
For $n\in\nz$ and $z\in \cz$\,, the spaces $\mathcal{W}_{g,z}^{n}$ are defined by
induction according to
\begin{eqnarray*}
  &&\mathcal{W}^{0}_{g,z}=L^{2}(X,dqdp)=L^{2}(T^{*}Q,dqdp)\,,\\
&&(u\in \mathcal{W}^{n+1}_{g,z})\Leftrightarrow
\left(\forall \ell\in \left\{1,\ldots,L\right\}, \forall T\in
  \mathfrak{T}_{\ell}\,,~T[\theta_{\ell}(q)u]\in \mathcal{W}^{n}_{g,z}\right)\,,\\
&&
\|u\|_{\mathcal{W}^{n+1}_{g,z}}^{2}=\sum_{\ell=1}^{M}\sum_{T\in
  \mathfrak{T_{\ell}}}\|T[\theta_{\ell}(q)u]\|_{\mathcal{W}^{n}_{g,z}}^{2}\,. 
\end{eqnarray*}
The space $\mathcal{W}_{g,z}^{s}$ for $s\in\rz$ are then defined by duality and
interpolation.
\end{definition}
For a fixed $s\in\rz$\,, the space $\mathcal{W}_{g,z}^{s}$ does not depend on
$z\in\cz$ but its norm $\|~\|_{\mathcal{W}^{s}_{g,z}}$ does. The quantities
$\sum_{i=1}^{d}\|\langle
p\rangle\partial_{p_{i}}v\|_{\mathcal{W}^{n}_{g,z}}^{2}+\|e_{i}v\|_{\mathcal{W}^{n}_{g,z}}^{2}$
and
$\sum_{i=1}^{d}\|\langle
p\rangle\partial_{p_{i}}v\|_{\mathcal{W}^{n}_{g,z}}^{2}+\|\partial_{q^{i}}v\|_{\mathcal{W}^{n}_{g,z}}^{2}$
are equivalent but the covariant derivative $e_{i}$ is a more natural
object. Note in particular the relations
\begin{eqnarray*}
  &&[e_{k}\,, |p|^{2}_{q}]=0\quad,\quad
\left[e_{k}\,,\,e_{k'}\right]={}^{t}R_{kk'}(q)\in \text{End}(T_{q}^{*}Q)\,,\\
&& \left[e_{k}\,,\,
  \Delta_{p}\right]=(\partial_{q^{k}}g_{ij})\partial_{p_{i}}\partial_{p_{j}}
+\Gamma^{\alpha}_{\beta,k}g_{ij}\left[p_{\alpha}\partial_{p_{\beta}}\,,\,
\partial_{p_{i}}\partial_{p_{j}}\right]\\
&&\qquad\qquad = \left(\partial_{q^{k}}g_{ij}-\Gamma^{\alpha}_{ik}g_{\alpha
    j}-\Gamma^{\alpha}_{jk}g_{i\alpha}\right)\partial_{p_{i}}\partial_{p_{j}}=0\,,\\
\text{owing~to}&&
\Gamma^{\alpha}_{ik}g_{\alpha j}+\Gamma^{\alpha}_{jk}g_{i \alpha}=\partial_{q^{k}}g_{ij}\,.
\end{eqnarray*}
By introducing a dyadic partition of unity in the $p$-variable
(see\cite{Leb1}\cite{Leb2} or Subsection~\ref{se.momentpart}), the squared norm
$\|u\|_{\mathcal{W}^{s}_{g,z}}^{2}$  is expressed as a quadratic series of
parameter dependent usual Sobolev norms. We leave the reader to check
with this approach that the space $\mathcal{W}^{n}_{g,z}$ does not depend on the
 the atlas and the subordinate partition of unity.
Actually in \cite{Leb1}\cite{Leb2} those spaces are defined with some
parameter dependent pseudo-differential calculus on compact manifolds
after introducing a dyadic partition of unity in the momentum variable
$p$\,.
Such a definition was extended in \cite{Bisorb} for some locally
symmetric non compact space $Q$\,. We focus here on the compact case
and will later simply extend some estimates to the case of 
cylinders $\rz\times Q$\,.
Again because $Q$ is compact without boundary, one has
\begin{eqnarray*}
&&
\mathcal{S}(X)=\ccap_{s\in\rz}\mathcal{W}^{s}_{g,z}\quad,\quad
  \mathcal{S}'(X)=\ccup_{s\in\rz}\mathcal{W}^{s}_{g,z}\,,\\
&&
\mathcal{C}^{\infty}_{0}(X)\subset
\mathcal{S}(X)\subset \mathcal{W}^{s_{1}}_{g,z}\subset
\mathcal{W}^{s_{2}}_{g,z}\subset \mathcal{S}'(X)\subset \mathcal{D}'(X)\,,\quad s_{1}\geq s_{2}\,,
\end{eqnarray*}
with dense and continuous embeddings.
\newtheorem*{thmLeb}{Theorem 1.2 of \cite{Leb2}}
\begin{thmLeb}
\label{th.Leb}
The riemannian manifold $(Q,g)$ is assumed to be compact (without boundary).
There exists $\delta_{1}>0$ such that the following properties hold.\\
As an operator in $\mathcal{S}(X)$ or $\mathcal{S}'(X)$\,, $P_{\pm, Q,g}$
satisfies
$$
\sigma(P_{\pm, Q,g})\subset
\left\{z\in\cz, \Real z\geq \frac{\dim(Q)}{2}+\delta_{1}|\Imag
  z|^{1/2}\right\}\,.
$$
For all $u\in \mathcal{S}'(X)$ and all $z\in U_{\delta_{1}}=\left\{z\in\cz\,, \Real z <
  \delta_{1}|\Imag z|^{1/2}\right\}$\,, the condition 
$(P_{\pm, Q,g}-z)u\in
\mathcal{W}^{s}_{g,0}$\,, $s\in\rz$\,, implies $u\in
\mathcal{W}^{s+\frac{2}{3}}_{g,0}$\,, $\mathcal{O}_{g}u\in
\mathcal{W}^{s}_{g,z}$\,, $\mathcal{Y}_{\mathcal{E}}u\in \mathcal{W}^{s}_{g,z}$
with the estimate
\begin{multline}
  \label{eq.subellLeb}
C_{s}\|(P_{\pm,Q,g}-z)u\|_{\mathcal{W}^{s}_{g,z}}
\geq
\|\mathcal{O}_{Q,g}u\|_{\mathcal{W}^{s}_{g,z}}+\|(\pm\mathcal{Y}_{\mathcal{E}}-i\Imag
z)u\|_{\mathcal{W}^{s}_{g,z}}
\\
+(|\Real z| +|\Imag z|^{1/2})\|u\|_{\mathcal{W}^{s}_{g,z}}+\|u\|_{\mathcal{W}^{s+\frac{2}{3}}_{g,z}}
\end{multline}
for some constant $C_{s}$\,.
\end{thmLeb}
We shall only use those estimates with $s\in [-1,1]$ and make the
connection with more usual Sobolev estimates. Note in particular for
$s\in [0,1]$ the embeddings and norm estimates
\begin{eqnarray*}
  && \mathcal{W}^{s}_{g,z}\subset \mathcal{W}^{s}_{g,0}\subset
  H^{s}(Q;\mathcal{H}^{0})\,,\\
&& \frac{1}{C_{s}}\|u\|_{H^{s}(Q;\mathcal{H}^{0})}\leq
\|u\|_{\mathcal{W}^{s}_{g,0}}
\leq \|u\|_{\mathcal{W}^{s}_{g,z}}\,,\\
&& H^{-s}(Q;\mathcal{H}^{0})\subset \mathcal{W}_{g,0}^{-s}\subset
\mathcal{W}_{g,z}^{-s}\,,\\
&& \|u\|_{\mathcal{W}^{-s}_{g,z}}\leq
\|u\|_{\mathcal{W}^{-s}_{g,0}}\leq C_{s}\|u\|_{H^{-s}(Q;\mathcal{H}^{0})}\,,
\end{eqnarray*}
obtained by duality and interpolation from the obvious cases $s=0$ and $s=1$\,.
\begin{corollary}
\label{co.Leb} 
Assume that $(Q,g)$ is a compact riemannian manifold (without boundary).
Consider the operator $K_{\pm,Q,g}-\frac{\dim Q}{2}$ defined in
$L^{2}(X)=L^{2}(T^{*}Q,dqdp)$ by
\begin{eqnarray*}
&  D(K_{\pm,Q,g})=\left\{u\in L^{2}(Q;\mathcal{H}^{1})\,,~P_{\pm,
      Q,g}u\in L^{2}(X)\right\}\\
&\forall u\in D(K_{\pm,Q,g})\,,\quad K_{\pm,Q,g}u=P_{\pm,Q,g}u=(\pm \mathcal{Y}_{\mathcal{E}}+\mathcal{O}_{Q,g})u\,,
\end{eqnarray*}
is maximal accretive and its adjoint is given by
$K_{\pm,Q,g}^{*}=K_{\mp,Q,g}$\,.\\
The domains $D(K_{\pm,Q,g})$ and $D(K_{\mp,Q,g})=D(K_{\pm,Q,g}^{*})$
are equal with equivalent graph norms. 
The spaces $\mathcal{S}(X)$ and $\mathcal{C}^{\infty}_{0}(X)$ are
dense in $D(K_{\pm,Q,g})$ endowed with its graph norm.\\
The estimates 
\begin{align}
\label{eq.coLebsub}
  &\|\mathcal{O}_{Q,g}u\|+\|(\pm
  \mathcal{Y}_{\mathcal{E}}-i\lambda)u\|+\langle
  \lambda\rangle^{\frac{1}{2}}\|u\|+\|u\|_{H^{2/3}(Q;\mathcal{H}^{0})}\leq
  C\|(K_{\pm,Q,g}-i\lambda)u\|\,,
\\
\label{eq.coLebacc}
&\frac{\dim Q}{2}\|u\|^{2}\leq
\|u\|_{L^{2}(Q;\mathcal{H}^{1})}^{2}\leq \Real \langle u\,,\, (K_{\pm,Q,g}-i\lambda)u\rangle
\end{align}
holds all $u\in D(K_{\pm,Q,g})$  and all $\lambda\in\rz$\,.\\
Moreover for any $\lambda\in\rz$\,, 
the resolvent $(K_{\pm,Q,g}-i\lambda)^{-1}$
extends by continuity 
\begin{itemize}
\item
 as an element of $\mathcal{L}(L^{2}(Q;\mathcal{H}^{-1});L^{2}(Q;\mathcal{H}^{1}))$
  with 
$$
\|(K_{\pm,Q,g}-i\lambda)^{-1}u\|_{L^{2}(Q;\mathcal{H}^{1})}\leq
\|u\|_{L^{2}(Q;\mathcal{H}^{-1})}\,,
$$
\item and as an element of
  $\mathcal{L}(H^{s}(Q;\mathcal{H}^{0});H^{s+\frac{2}{3}}(Q;\mathcal{H}^{0}))$
  for any $s\in [-\frac{2}{3},0]$ with 
$$
\|(K_{\pm,Q,g}-i\lambda)^{-1}u\|_{H^{s+\frac{2}{3}}(Q;\mathcal{H}^{0})}\leq C_{s}
\|u\|_{H^{s}(Q;\mathcal{H}^{0})}\,.
$$ 
\end{itemize}
\end{corollary}
\begin{proof}
The operator $P_{\pm,Q,g}-\frac{\dim Q}{2}$ initially defined with the domain
$\mathcal{S}(X)$\,, is accretive owing to the simple integration by
parts
$$
\forall u\in \mathcal{S}(X)\,, \quad
\Real\langle u\,,\, P_{\pm,Q,g}u\rangle =
\|u\|_{L^{2}(Q;\mathcal{H}^{1})}^{2}\geq \frac{\dim Q}{2}\|u\|^{2}\,.
$$
According to \cite{ReSi}, its closure denoted by $K_{\pm,
  Q,g}=\overline{P_{\pm,Q,g}}$
is accretive. It is actually maximal accretive because the equation
$P_{\mp,Q,g}=P_{\pm,Q,g}^{*}u=0$
 in $\mathcal{S}'(X)$ and the regularity result of
 Theorem~1.2 in \cite{Leb2} implies $u\in \mathcal{S}(X)$ and
 $u=0$\,. Therefore the range $\Ran P_{\pm,Q,g}$ is dense in
 $L^{2}(X)$ (use $\overline{\Ran P_{\pm,
     Q,g}}=\ker(P_{\pm,Q,g}^{*})=\left\{0\right\}$). This suffices to
 prove the maximal accretivity of $K_{\pm,Q,g}$\,.\\
From this definition of $K_{\pm,Q,g}$\,, $\mathcal{S}(X)$ (and
therefore $\mathcal{C}^{\infty}_{0}(X)$ because $P_{\pm,Q,g}$ has
polynomially controlled coefficients) is a core for $K_{\pm,Q,g}$\,.
The subelliptic estimates, the identification of the domain
$D(K_{\pm,Q,g})$ as it is stated, the equality
$D(K_{\pm,Q,g})=D(K_{\mp,Q,g})=D(K_{\pm,Q,g}^{*})$ are direct consequences of
\eqref{eq.subellLeb} applied with $s=0$\,. We use the upper bound
$$
\|u\|_{H^{\frac{2}{3}}(Q;\mathcal{H}^{0})}\leq C\|u\|_{\mathcal{W}^{\frac{2}{3}}_{g,i\lambda}}\,.\\
$$
The extension of
$(K_{\pm,Q,g}-i\lambda)^{-1}\in
\mathcal{L}(L^{2}(Q;\mathcal{H}^{-1});L^{2}(Q;\mathcal{H}^{-1}))$ is
contained in \eqref{eq.coLebacc}.\\
The extension of  $(K_{\pm,Q,g}-i\lambda)^{-1}\in
\mathcal{L}(H^{s}(Q;\mathcal{H}^{0});
H^{s+\frac{2}{3}}(Q;\mathcal{H}^{0}))$ comes from
$$
\|u\|_{H^{s+\frac{2}{3}}(Q;\mathcal{H}^{0})}\leq
C_{s}\|u\|_{\mathcal{W}^{s+\frac{2}{3}}_{g,i\lambda}}
\leq
C_{s}'\|u\|_{\mathcal{W}^{s}}\leq C''_{s}\|u\|_{H^{s}(Q;\mathcal{H}^{0})}\,,
$$
due to \eqref{eq.subellLeb} with $s\leq 0$ and $s+\frac{2}{3}\geq 0$\,.
\end{proof}
\subsection{Partitions of unity}
\label{se.partunitGKFP}
Like in \cite{Leb1}\cite{Leb2}, we introduce
 two partitions of unity: a spatial partition of unity
which allows to use coordinate systems, a dyadic partition of unity in
the momentum variable which replaces the analysis as $p\to \infty$ by
some parameter dependent problem.
In \cite{Leb1}\cite{Leb2}, all formulas are first proved for elements of
$\mathcal{S}(X)$ and extended by density arguments.
Here, in order to extend those decompositions to boundary value problems
where a possibly non regularity-preserving operator $A$ prevents from
identifying a core of regular functions, various quantities are
directly calculated for $u\in L^{2}(Q;\mathcal{H}^{1})$ such that
$P_{\pm,Q,g}u\in L^{2}(Q;\mathcal{H}^{-1})$ by using the $L^{2}(X)$ scalar
product and the 
$L^{2}(Q;\mathcal{H}^{1})-L^{2}(Q;\mathcal{H}^{-1})$ duality product
(both denoted by $\langle~,~\rangle$).
\subsubsection{Spatial partition of unity}
\label{se.spaPU}
The riemannian manifold $(\overline{Q}=Q\sqcup \partial Q,g)$  can be
\begin{itemize}
\item a compact manifold without or with a boundary\,,
\item a compact perturbation of the euclidean space $\rz^{d}$ or the
  half-space $\overline{\rz^{d}}_{-}=(-\infty,0]\times \rz^{d-1}$ with
  $g_{ij}-\delta_{ij}\in \mathcal{C}^{\infty}_{0}(\rz^{d})$\,,
\item a cylinder $\rz\times Q'$ or a half cylinder $(-\infty,0]\times
  Q'$ with $Q'$ compact, $\partial Q'=\emptyset$\,, 
$g=1\oplus m(q^{1},q')$\,, $m_{ij}-m_{ij}(-\infty,q')\in \mathcal{C}^{\infty}_{0}(\overline{Q})$\,.
\end{itemize}
Such a manifold $\overline{Q}$ admits finite
and countable
locally finite coverings $\overline{Q}=\ccup_{\ell\in \mathcal{L}}\Omega_{\ell}$ where
$\Omega_{\ell}$ is diffeomorphic to a regular domain of $\rz^{d}$ and
$g_{ij}(q)\in \mathcal{C}^{\infty}_{b}(\Omega_{\ell})$ in the
corresponding coordinate system\,, with a
bounded intersection number
$$
\sup_{\ell'\in \mathcal{L}}\sharp\left\{\ell\in\mathcal{L},\;
  \Omega_{\ell}\cap \Omega_{\ell'}\neq \emptyset \right\}\leq N_{\Omega}\in\nz\,.
$$
Furthermore the covering can be chosen such that there exists a
partition of unity
 $\chi=(\chi_{\ell})_{\ell\in \mathcal{L}}$\,, $\chi_{\ell}\in 
\mathcal{C}^{\infty}_{b}(\Omega_{\ell})$ with uniform seminorms w.r.t
$\ell\in \mathcal{L}$ \,, $\supp \chi_{\ell}\subset
\Omega_{\ell}$ and 
$\sum_{\ell\in \mathcal{L}}\chi_{\ell}^{2}=1$\,. The simplest way
consists in taking $\mathcal{L}$ finite.
\\
Due to the vertical nature
of $\mathcal{O}_{Q,g}$\,, the condition $u\in
L^{2}(Q;\mathcal{H}^{s'})$ is equivalent to
$$
\left\{
  \begin{array}[c]{l}
\forall \ell\in \mathcal{L}\,, \quad \chi_{\ell}u \in
L^{2}(\Omega_{\ell};\mathcal{H}^{s'})\,,\\
\sum_{\ell\in \mathcal{L}}\|\chi_{\ell}u\|_{L^{2}(\Omega_{\ell};\mathcal{H}^{s'})}^{2}< +\infty\,,
\end{array}
\right.
$$
with the equality of squared norms
$\|u\|_{L^{2}(Q;\mathcal{H}^{s'})}^{2}=\sum_{\ell\in
  \mathcal{L}}\|\chi_{\ell}u\|_{L^{2}(\Omega_{\ell};\mathcal{H}^{s'})}^{2}$\,.
We recall also (see for instance \cite{ChPi}), that when $\partial Q=\emptyset$ the usual Sobolev
norms in $q$ defined by
$\|u\|_{H^{s}(\overline{Q};\mathcal{H}^{0})}=\|(1-\Delta_{q})^{s/2}u\|$ with
$\Delta_{q}=\partial_{q^{i}}g^{ij}(q)\partial_{q^{j}}$ are equivalent
via a change of the metric and satisfy 
$$
\left(\frac{\sum_{\ell\in
      \mathcal{L}}\|\chi_{\ell}u\|_{H^{s}(Q;\mathcal{H}^{0})}^{2}}{\|u\|_{H^{s}(Q;\mathcal{H}^{0})}^{2}}\right)^{\pm
1}\leq C_{s,g,\chi}\,.
$$
In the case when $\partial{Q}\neq \emptyset$\,, Sobolev spaces are defined by using locally
the reflection principle (see e.g. \cite{ChPi}).
\begin{proposition}
\label{pr.partq}
 Take the above partition of unity and assume $z\in\cz$\,, $s'\in [-1,0]$\,.
 For $u\in L^{2}(Q;\mathcal{H}^{1})$\,,
 the condition $(P_{\pm,Q,g}-z)u\in L^{2}(Q;\mathcal{H}^{s'})$ is
    equivalent to 
\begin{equation}
   \label{eq.localPQg}
\forall \ell\in \mathcal{L}\,,\quad P_{\pm, Q,g}(\chi_{\ell}u)\in
L^{2}(Q;\mathcal{H}^{s'})\,.
\end{equation}
There exists a constant $C_{\chi}>0$ such that 
\begin{equation}
\label{eq.encaPQgloc2}
\left(\frac{\sum_{\ell\in
  \mathcal{L}}\|(P_{\pm,Q,g}-z)\chi_{\ell}u\|_{L^{2}(Q;\mathcal{H}^{s'})}^{2}+\|\chi_{\ell}u\|_{L^{2}(Q;\mathcal{H}^{1})}^{2}}{\|(P_{\pm,Q,g}-z)u\|_{L^{2}(Q;\mathcal{H}^{s'})}^{2}+
\|u\|_{L^{2}(Q;\mathcal{H}^{1})}^{2}}\right)^{\pm 1}
\leq C_{\chi} \,,
\end{equation}
holds uniformly with respect to $(u,z,s')$\,.\\
Two functions $u_{1},u_{2}\in L^{2}(Q;\mathcal{H}^{1})$ such that
$P_{\pm,Q,g}u_{1}, P_{\pm,Q,g}u_{2}\in L^{2}(Q;\mathcal{H}^{-1})$\,, satisfy the identity
\begin{equation}
  \label{eq.eqpartscal}
  \langle u_{1}\,,\, P_{\pm,Q,g}u_{2}\rangle=\sum_{\ell\in \mathcal{L}}\langle
  \chi_{\ell}u_{1}\,,\, P_{\pm,Q,g}\chi_{\ell}u_{2}\rangle\,.
\end{equation}
\end{proposition}
\begin{proof}
When $u\in L^{2}(Q;\mathcal{H}^{1})$ satisfies $(P_{\pm,Q,g}-z)u\in
L^{2}(Q;\mathcal{H}^{s'})$\,, \eqref{eq.localPQg} comes from
$$
(P_{\pm,Q,g}-z)\chi_{\ell}u=\chi_{\ell}(P_{\pm,Q,g}-z)u
\pm g^{ij}(q)(\partial_{q_{j}}\chi_{\ell}(q))p_{j}u\,.
$$
With $s'\leq 0$\,, it implies
\begin{eqnarray*}
  \sum_{\ell\in
    \mathcal{L}}\|(P_{\pm,Q,g}-z)\chi_{\ell}u\|_{L^{2}(Q;\mathcal{H}^{s'})}^{2}
&\leq & 2\sum_{\ell\in
    \mathcal{L}}\|\chi_{\ell}(P_{\pm,Q,g}-z)u\|_{L^{2}(Q;\mathcal{H}^{s'})}^{2}\\
&&\hspace{3cm}
+\|g^{ij}(q)(\partial_{q_{j}}\chi_{\ell}(q))p_{j}u\|^{2}\\
\\
&\leq & C_{\chi}^{1}\left[\|(P_{\pm,Q,g}-z)u\|_{L^{2}(Q;\mathcal{H}^{s'})}^{2}+\|u\|_{L^{2}(Q;\mathcal{H}^{1})}^{2}\right]\,,
\end{eqnarray*}
because $\sum_{\ell'\in
  \mathcal{L}}|\partial_{q}\chi_{\ell}(q)|^{2}\leq
N_{\Omega}\sup_{\ell\in \mathcal{L}}\|\partial_{q}\chi_{\ell}\|_{L^{\infty}}^{2}<\infty$\,.\\
The relation~\eqref{eq.partunitRe} applied with $P=P_{\pm,Q,g }$ and
$\ad_{\chi_{\ell}}^{2}P_{\pm,Q,g }=0$\,, gives here
$$
P_{\pm,Q,g}=\sum_{\ell\in \mathcal{L}}\chi_{\ell}P_{\pm,Q,g}\chi_{\ell}\,.
$$
Hence \eqref{eq.localPQg} implies $\chi_{\ell'}(P_{\pm,Q,g}-z)u\in
L^{2}(Q;\mathcal{H}^{s'})$ for $\ell'\in \mathcal{L}$\,, with the estimate
\begin{eqnarray*}
\|\chi_{\ell'}(P_{\pm,Q,g}-z)u\|_{L^{2}(Q;\mathcal{H}^{s'})}^{2}
&\leq&\|\sum_{\ell\in\mathcal{L}}\chi_{\ell'}\chi_{\ell}
P_{\pm,Q,g}\chi_{\ell}u\|_{L^{2}(Q;\mathcal{H}^{s'})}^{2}
\\
&\leq& N_{\Omega}\sum_{\chi_{\ell}\chi_{\ell'}\neq
  0}\|(P_{\pm,Q,g}-i\lambda)\chi_{\ell}u\|_{L^{2}(Q;\mathcal{H}^{s'})}^{2}\,.
\end{eqnarray*}
Taking the sum over $\ell'\in \mathcal{L}$ yields
$$
\|(P_{\pm,Q,g}-i\lambda)u\|_{L^{2}(Q,\mathcal{H}^{s'})}^{2}\leq
N_{\Omega}^{2}\sum_{\ell\in \mathcal{L}}
\|(P_{\pm,Q,g}-i\lambda)\chi_{\ell}u\|_{L^{2}(Q;\mathcal{H}^{s'})}^{2}\,.
$$
The above identity also proves\eqref{eq.eqpartscal}.
\end{proof}

\subsubsection{Dyadic momentum partition of unity}
\label{se.momentpart}
The  previous partition of unity in the $q$-variable allows to
focus on the case when $\overline{Q}$ is endowed with a global
coordinate system, namely:
\begin{itemize}
\item $\overline{Q}$ is a compact perturbation of the euclidean space
  $\rz^{d}$ or half-space $\overline{\rz^{d}}_{-}=(-\infty,0]\times \rz^{d-1}$\,,
\item $\overline{Q}$ is the torus $\tz^{d}$\,, the cylinder $\rz\times
  \tz^{d-1}$ or half-cylinder $(-\infty,0]\times \tz^{d-1}$
 with the metric $1\oplus m(q^{1},q')$\,, $m_{ij}-m_{ij}(-\infty, q')\in
 \mathcal{C}^{\infty}_{0}(\overline{Q})$\,. 
\end{itemize}
Take two functions $\tilde{\chi}_{0}\in \mathcal{C}^{\infty}_{0}(\rz)$ and
$\tilde{\chi}_{1}\in \mathcal{C}^{\infty}_{0}((0,+\infty))$ such that
$\tilde{\chi}_{0}^{2}(t)+\sum_{\ell\in \nz^{*}}\tilde\chi_{1}^{2}(2^{-\ell}t)\equiv 1$ on
$[0,+\infty)$\,. Set $\chi_{\ell}(q,p)=\tilde{\chi}_{1}(2^{-\ell}|p|_{q})$
for $\ell\in\nz^{*}$\,, $\chi_{0}(q,p)=\tilde{\chi}_{0}(|p|_{q})$
and
$\chi=(\chi_{\ell})_{\ell\in \mathcal{L}}$\,, $\mathcal{L}=\nz$\,.\\
The two following lemmas recall the equivalence of norms for
 the functional spaces
$L^{2}(Q;\mathcal{H}^{s'})$ and $H^{s}(Q;\mathcal{H}^{0})$\,.
A proposition similar to Proposition~\ref{pr.partq} is proved afterwards.
\begin{lemma}
\label{le.partp} 
With the above partition of unity
$\chi=\left(\chi_{\ell}\right)_{\ell\in\mathcal{L}}$\,,
$\chi_{\ell}=\tilde{\chi}_{1,0}\left(\frac{|p|_{q}}{2^{\ell}}\right)$\,,
$\sum_{\ell\in \mathcal{L}}\chi_{\ell}^{2}\equiv 1$ and
$\mathcal{L}=\nz$\,, the norms $\left(\sum_{\ell\in
  \mathcal{L}}\|\chi_{\ell}u\|_{L^{2}(Q;\mathcal{H}^{s'})}^{2}\right)^{1/2}$
 and 
$\|u\|_{L^{2}(Q;\mathcal{H}^{s'})}$
 are equivalent
for any $s'\in\rz$\,.
\end{lemma}
\begin{proof}
  For the equivalence of the norms in $L^{2}(Q;\mathcal{H}^{s'})$\,, we can
assume $Q=\rz^{d}$\,,  $g_{ij}(q)=\delta_{ij}$ and consider first $s'=n\in\nz$\,. 
The
result for a general $s'\in \rz$ follows by duality and interpolation (see the
application of Lemma~\ref{le.duainterp}).
Moreover density arguments are allowed here and all the calculations
can be made with $u\in \mathcal{S}(X)=\mathcal{S}(\rz^{2d})$\,.\\
We have to compare $\sum_{|\alpha|+m\leq n}\|\langle
p\rangle^{m}\partial_{p}^{\alpha}u\|^{2}$ and 
$\sum_{\ell\in\nz, |\alpha|+m\leq n}\|\langle
p\rangle^{m}\partial_{p}^{\alpha}
\chi_{\ell}u\|^{2}$\,.
Introduce for $m\in \nz$ and $m'\in\rz$\,, the class $\DiffP^{m,m'}$ of differential operators 
$a(p,D_{p})=\sum_{|\alpha|\leq m}a_{\alpha}(p)D_{p}^{\alpha}$ 
$$
\forall \beta\in\nz^{\dim Q}\,,
\quad
N^{m,m'}_{\alpha,\beta}(a)=\sup_{p\in\rz^{d}}\frac{|\partial_{p}^{\beta}a_{\alpha}(p)|}{\langle
  p\rangle^{m'-|\beta|}} <+\infty\,.
$$
For every $(m,m')\in\nz\times\rz$\,, $\DiffP^{m,m'}$ endowed with the seminorms
$N^{m,m'}_{\alpha,\beta}$, $|\alpha|\leq m, \beta\in\nz^{\dim Q}$ is
a Fr{\'e}chet space. The union
$\DiffP=\ccup_{(m,m')\in\nz\times\rz}\DiffP^{m,m'}$ is a bigraded
algebra
\footnote{The aware reader will recognize a subalgebra of the
  Weyl-H{\"o}rmander pseudodifferential class associated with the metric 
$\frac{dp^{2}}{\langle p\rangle^{2}}+\frac{d\eta^{2}}{\langle
  \eta\rangle^{2}}$ and the gain function $\langle p\rangle\langle
\eta\rangle$ (see \cite{HormIII}-Chap~18 or \cite{Helgl}). The restriction to
differential operators, which contains differential operators with
polynomial coefficients, allows an easier handling of supports without
requiring the more sophisticated notion of confinement (see \cite{BoLe}).} 
such that
\begin{eqnarray*}
&&\DiffP^{m_{1},m_{1}'}\circ \DiffP^{m_{2},m_{2}'}\subset
\DiffP^{m_{1}+m_{2},m_{1}'+m_{2}'}\,,\\
&&
\left[\DiffP^{m_{1},m_{1}'}\,,\, \DiffP^{m_{2},m_{2}'}\right]\subset
\DiffP^{m_{1}+m_{2}-1,m_{1}'+m_{2}'-1} \,,
\end{eqnarray*}
the corresponding mappings $(A,B)\mapsto A\circ B$ and $(A,B)\mapsto
[A,B]$ being bilinearly continuous. Note also that the
$p$-support
$$
p-\supp a(p,D_{p})=\ccup_{|\alpha|\leq m}\supp a_{\alpha}\,,
$$
satifies
$$
p-\supp \left[a(p,D_{p})\circ b(p,D_{p})\right]\subset p-\supp
a(p,D_{p}) \cap p-\supp b(p,D_{p})\,.
$$
The $\chi_{\ell}$'s, with
$\chi_{\ell}(p)=\tilde{\chi}_{0,1}\left(\frac{|p|}{2^{\ell}}\right)$\,,
make a uniformly bounded family 
of $\DiffP^{0,0}$ with $p-\supp \nabla\chi_{\ell}\subset
\left\{\frac{2^{\ell}}{C}\leq |p|\leq C2^{\ell}\right\}$\,.\\
The proof of the equivalence of norms is now done by induction:
\begin{itemize}
\item For $n=0$ it is obvious: $L^{2}(Q;\mathcal{H}^{0})=L^{2}(X)$\,.
\item If it is true for all $n'\leq N\in\nz$\,,  apply the formula
  \eqref{eq.partunit2} with $P=\langle p\rangle^{m}D_{p}^{\alpha}$ and
  $|\alpha|+m=n+1$\,. It gives for any $u\in \mathcal{S}(X)=\mathcal{S}(\rz^{2d})$\,,
$$
\|\langle
p\rangle^{m}D_{p}^{\alpha}u\|^{2}-\sum_{\ell\in \mathcal{L}}\|\langle
p\rangle^{m} D_{p}^{\alpha}\chi_{\ell}u\|^{2}
=-\sum_{\ell\in\mathcal{L}}\|A_{\ell}u\|^{2}+\Real\langle
u\,,\, B_{\ell}u\rangle\,,
$$
where
$A_{\ell}=\ad_{\chi_{\ell}}(\langle p\rangle^{m}D_{p}^{\alpha})$
is uniformly bounded in $\DiffP^{(|\alpha|-1)_{+},(m-1)_{+}}$\,, 
$B_{\ell}=D_{p}^{\alpha}\langle
p\rangle^{m}\ad_{\chi_{\ell}}^{2}(\langle p\rangle^{m}D_{p}^{\alpha})$
is uniformly bounded in $\DiffP^{(2|\alpha|-2)_{+},(2m-2)_{+}}$\,, with
$p-$supports contained in $\left\{\frac{2^{\ell}}{C}\leq |p|\leq
  C2^{\ell}\right\}$ and covered by a fixed number of
$\chi_{\ell'}$\,. The equivalence of norms for $n'\leq n$ implies
that
$$
\left|
\|\langle
p\rangle^{m}D_{p}^{\alpha}u\|^{2}-\sum_{\ell\in \mathcal{L}}\|\langle
p\rangle^{m} D_{p}^{\alpha}\chi_{\ell}u\|^{2}
\right|
$$
is estimated by
$$
 C\|u\|_{L^{2}(Q;\mathcal{H}^{n})}^{2}
\quad\text{or}\quad C'\sum_{\ell\in \mathcal{L}}\|\chi_{\ell}u\|_{L^{2}(Q;\mathcal{H}^{n})}^{2}\,.
$$
Since $|a_{n+1}-b_{n+1}|\leq C\sum_{p\leq n}b_{p}$ and $\sum_{p\leq
  n}a_{p}\leq C \sum_{p\leq n}b_{b}$ imply 
$\sum_{p\leq n+1}a_{p}\leq C' \sum_{p\leq n+1}b_{p}$\,, this proves
the equivalence for $n'=n+1$\,.
The density of $\mathcal{S}(X)=\mathcal{S}(\rz^{2d})$ in
$L^{2}(Q;\mathcal{H}^{n+1})$ ends the proof.
\end{itemize}
\end{proof}
\begin{lemma}
\label{le.partpq}
With the above partition of unity
$\chi=\left(\chi_{\ell}\right)_{\ell\in\mathcal{L}}$\,,
$\chi_{\ell}=\tilde{\chi}_{1,0}\left(\frac{|p|_{q}}{2^{\ell}}\right)$\,,
$\sum_{\ell\in \mathcal{L}}\chi_{\ell}^{2}\equiv 1$ and
$\mathcal{L}=\nz$\,.
For any $s\in [0,1]$ there exists $C_{\chi, s}$ such that 
$$
\left(
\frac{\sum_{\ell\in \mathcal{L}}
\|\chi_{\ell}u\|_{H^{s}(\overline{Q};\mathcal{H}^{0})}^{2}}{\|u\|_{H^{s}(\overline{Q};\mathcal{H}^{0})}^{2}}
\right)^{\pm 1} 
\leq C_{\chi,s}\,.
$$
\end{lemma}
\begin{proof}
Again the spatial partition of unity of Subsection~\ref{se.spaPU}, the
reflection principle and
the equivalence of norms 
for two different metrics restrict the analysis to the case 
$Q=\rz^{d}$\,, $g_{ij}=\delta_{ij}$\,.\\
 The result  is obvious when $s=0$\,,
 $H^{0}(\rz^{d};\mathcal{H}^{0})=L^{2}(\rz^{2d},dqdp)$\,.\\
Consider now the case $s=1$\,, with
$$
\|u\|_{H^{1}(Q;\mathcal{H}^{0})}^{2}=\|u\|^{2}+\sum_{j=1}^{d}\|\partial_{q^{j}}u\|^{2}\,.
$$
We use again the formula \eqref{eq.partunit2}
$$
\|\partial_{q^{j}}u\|^{2}-\sum_{\ell\in
  \mathcal{L}}\|\partial_{q^{j}}\chi_{\ell}u\|^{2}=
-\sum_{\ell\in
  \mathcal{L}}\|(\ad_{\chi_{\ell}}\partial_{q^{j}})u\|^{2}\,,
$$
with
\begin{eqnarray*}
  &&
\ad_{\chi_{\ell}}\partial_{q^{j}}=(\partial_{q^{j}}\chi_{\ell})=
\frac{\partial_{q^{j}}(|p|_{q})}{2^{\ell}}\tilde{\chi}_{0,1}'(\frac{|p|_{q}}{2^{\ell}})
\\
\text{and}&&
\ad_{\chi_{\ell}}^{2}\partial_{q^{j}}=0\,.
\end{eqnarray*}
The derivative $\partial_{q^{j}}|p|_{q}$ satisfies
$$
\frac{1}{2^{\ell}}\partial_{q^{j}}|p|_{q}=\frac{\partial_{q^{j}}g^{ik}(q)p_{i}p_{k}}{2\times
  2^{\ell}(g^{ik}(q)p_{i}p_{k})^{1/2}}=\mathcal{O}(1)\quad
\text{in}~\supp \tilde{\chi}_{0,1}'\left(\frac{|p|_{q}}{2^{\ell}}\right)\,.
$$
and this proves
$$
\left|\|u\|_{H^{1}(\rz^{d};\mathcal{H}^{0})}^{2}-\sum_{\ell\in
    \mathcal{L}}\|\chi_{\ell}u\|_{H^{1}(\rz^{d};\mathcal{H}^{0})}^{2}\right|\leq C_{\chi}^{1}\|u\|^{2}\,,
$$
which implies the equivalence of norms for $s=1$\,.\\
 The general case follows by
interpolation with
$\|u\|_{H^{r}(\rz^{d};\mathcal{H}^{0})}=\|(1-\Delta_{q})^{r/2}u\|$ (see
Lemma~\ref{le.duainterp} and its application).
\end{proof}
While working with $P_{\pm,Q,g}$\,, we avoid again density arguments in
order to allow the same line for general boundary value problems.
\begin{proposition}
\label{pr.partp}
Take the above dyadic partition of unity
$\chi=(\chi_{\ell})_{\ell\in\mathcal{L}}$\,,
$\chi_{\ell}=\tilde{\chi}_{1,0}\left(\frac{|p|_{q}}{2^{\ell}}\right)$\,,
$\sum_{\ell\in \mathcal{L}}\chi_{\ell}^{2}\equiv 1$ and
$\mathcal{L}=\nz$\,,
and let $u\in L^{2}(Q;\mathcal{H}^{1})$\,, $z\in\cz$ and
$s'\in[-1,0]$\,.\\
The condition $(P_{\pm,Q,g}-z)u\in
L^{2}(Q;\mathcal{H}^{s'})$ is equivalent 
to $\sum_{\ell\in\mathcal{L}}\|(P_{\pm,Q,g}-z)\chi_{\ell}u\|_{L^{2}(Q;\mathcal{H}^{s'})}^{2}<+\infty$\,,
and the equivalence of norms \eqref{eq.encaPQgloc2}
still holds (with the new choice 
$\chi_{\ell}=\chi_{\ell}(p)$).\\
For any pair $u_{1},u_{2}\in L^{2}(Q;\mathcal{H}^{1})$ such that
$P_{\pm,Q,g}u_{1}, P_{\pm,Q,g}u_{2}\in
L^{2}(Q;\mathcal{H}^{-1})$\,, the identity~\eqref{eq.eqpartscal} is now
replaced by
\begin{equation}
  \label{eq.ineqpartpscal}
  \left| \langle u_{1}\,,\, P_{\pm,Q,g}u_{2}\rangle-\sum_{\ell\in \mathcal{L}}\langle
  \chi_{\ell}u_{1}\,,\, P_{\pm,Q,g}\chi_{\ell}u_{2}\rangle\right|
\leq C_{\chi}\|u_{1}\|\|u_{2}\|\,.
\end{equation}
\end{proposition}
\begin{proof}
  The proof is essentially the same as for
  Proposition~\ref{pr.partq}. We shall use the notation
  $\DiffP^{m,m'}$ and $p-\supp$ introduced in the proof of
  Lemma~\ref{le.partp}.\\
With
$\left[\mathcal{Y}_{\mathcal{E}},\chi_{\ell}\right]=(\mathcal{Y}_{\mathcal{E}}\chi_{\ell})=0$\,,
the commutator $\ad_{\chi_{\ell}}(P_{\pm,Q,g}-z)$ equals
$$
\ad_{\chi_{\ell}}P_{\pm,Q,g}=\left[\chi_{\ell},-\frac{\Delta_{p}}{2}\right]=(\nabla_{p}\chi_{\ell}).\nabla_{p}+\frac{1}{2}(\Delta_{p}\chi_{\ell})\;
\in \DiffP^{1,-1}\,,
$$
where the gradients in the $p$-variable, the Laplace operator
$\Delta_{p}$ and the scalar product denoted by $x.y$ are given by the
scalar product $g(q)$\,. In particular
$\chi_{\ell}(p)=\tilde{\chi}_{0,1}\left(\frac{|p|_{q}}{2^{\ell}}\right)$
gives:
$$
\nabla_{p}\chi_{\ell}(p)=\frac{p}{2^{\ell}|p|_{q}}
\tilde{\chi}_{0,1}'\left(\frac{|p|_{q}}{2^{\ell}}\right)
\quad\text{and}\quad
\Delta_{p}\chi_{\ell}(p)=\frac{1}{2^{2\ell}}\tilde{\chi}_{0,1}''
\left(\frac{|p|_{q}}{2^{\ell}}\right)\,.
$$
This also implies $p-\supp (\ad_{\chi_{\ell}}P_{\pm,Q,g})\subset
\left\{\frac{2^{\ell}}{C}\leq |p|_{q}\leq C2^{\ell}\right\}$\,.\\
\noindent\textbf{a)} Assume $u\in L^{2}(Q;\mathcal{H}^{1})$ and $(P_{\pm,Q,g}-z)u\in
L^{2}(Q;\mathcal{H}^{s'})$ with $s'\in [-1,0]$\,.
The equality 
$$
\chi_{\ell}(P_{\pm,Q,g}-z)u=(P_{\pm,Q,g}-z)\chi_{\ell}u +(\ad_{\chi_{\ell}}P_{\pm,Q,g})u
$$
leads to
\begin{eqnarray*}
  \sum_{\ell\in
    \mathcal{L}}\|(P_{\pm,Q,g}-z)\chi_{\ell}u\|_{L^{2}(Q;\mathcal{H}^{s'})}^{2}
&\leq& 2\sum_{\ell\in
  \mathcal{L}}\|\chi_{\ell}(P_{\pm,Q,g}-z)u\|_{L^{2}(Q;\mathcal{H}^{s'})}^{2}
\\
&&
\hspace{3cm}
+
2\sum_{\ell\in \mathcal{L}}\|(\ad_{\chi_{\ell}}P_{\pm,Q,g})u\|^{2}\,.
\end{eqnarray*}
Lemma~\ref{le.partp} says that the first term is less than
$C_{\chi}^{1}\|(P_{\pm,Q,g}-z)u\|_{L^{2}(Q;\mathcal{H}^{s'})}^{2}$\,. 
The listed properties of $\ad_{\chi_{\ell}}P_{\pm,Q,g}$ imply
$$
\sum_{\ell\in \mathcal{L}}\|(\ad_{\chi_{\ell}}P_{\pm,Q,g})u\|^{2}
\leq C_{\chi}^{2}\sum_{\ell\in \mathcal{L}}\|\sum_{\chi_{\ell}\chi_{\ell'}\neq
  0}\chi_{\ell'}u\|_{L^{2}(Q;\mathcal{H}^{1})}^{2}\leq
C_{\chi}^{3}\sum_{\ell'\in \mathcal{L}}\|\chi_{\ell'}u\|_{L^{2}(Q;\mathcal{H}^{1})}^{2}
$$
because $\sharp\left\{\ell'\in \mathcal{L}\,,
  \chi_{\ell}\chi_{\ell'}\neq 0\right\}$ is uniformly bounded w.r.t
$\ell\in\mathcal{L}$\,.
The right-hand side $\sum_{\ell'\in
  \mathcal{L}}\|\chi_{\ell'}u\|_{L^{2}(Q;\mathcal{H}^{1})}^{2}$ is
equivalent to $\|u\|_{L^{2}(Q;\mathcal{H}^{1})}^{2}$ by
Lemma~\ref{le.partp}.
We have proved
\begin{multline*}
\sum_{\ell\in
  \mathcal{L}}\|(P_{\pm,Q,g}-z)\chi_{\ell}u\|_{L^{2}(Q;\mathcal{H}^{s'})}^{2}
+\|\chi_{\ell}u\|_{L^{2}(Q;\mathcal{H}^{1})}^{2}
\\
\leq C_{\chi}\left[
(P_{\pm,Q,g}-z)u\|_{L^{2}(Q;\mathcal{H}^{s'})}^{2}
+\|u\|_{L^{2}(Q;\mathcal{H}^{1})}^{2}
\right]
<+\infty.
\end{multline*}
\noindent \textbf{b)} The formula \eqref{eq.partunitRe} applied with
$P=P_{\pm,Q,g}$ and
$\ad_{\chi_{\ell}}^{2}P_{\pm,Q,g}=|\nabla_{p}\chi_{\ell}|^{2}_{q}=-\left[\frac{1}{2^{\ell}}\tilde{\chi}_{0,1}'\left(\frac{|p|_{q}}{2^{\ell}}\right)\right]^{2}$\,,
gives
\begin{equation}
  \label{eq.partscal}
P_{\pm,Q,g}-\sum_{\ell\in
  \mathcal{L}}\chi_{\ell}P_{\pm,Q,g}\chi_{\ell}
=
+\frac{1}{2}\sum_{\ell\in
  \mathcal{L}}\left[\frac{1}{2^{\ell}}\tilde{\chi}_{0,1}'\left(\frac{|p|_{q}}{2^{\ell}}\right)\right]^{2}\quad\in \DiffP^{0,-2}\,.
\end{equation}
This yields the  inequality \eqref{eq.ineqpartpscal}.\\
\noindent\textbf{c)} Let us finish the proof of the first
statement. Conversely to a) assume now $u\in L^{2}(Q;\mathcal{H}^{1})$
and $\sum_{\ell\in
  \mathcal{L}}\|(P_{\pm,Q,g}-z)\chi_{\ell}u\|_{L^{2}(Q;\mathcal{H}^{s'})}^{2}<+\infty$\,. The
relation \eqref{eq.partscal} implies for any $\ell'\in \mathcal{L}$\,,
\begin{eqnarray*}
  \chi_{\ell'}(P_{\pm,Q,g}-z)u&=&\chi_{\ell'}\sum_{\ell\in\mathcal{L}}\chi_{\ell}(P_{\pm,Q,g}-z)\chi_{\ell}u
+\frac{1}{2}\left[\frac{1}{2^{\ell}}\tilde{\chi}_{0,1}'\left(\frac{|p|_{q}}{2^{\ell}}\right)\right]^{2}u\\
&=&\sum_{\chi_{\ell}\chi_{\ell'}\neq 0}
\chi_{\ell'}\chi_{\ell}(P_{\pm,Q,g}-z)\chi_{\ell}u
+\frac{1}{2}\left[\frac{1}{2^{\ell}}\tilde{\chi}_{0,1}'\left(\frac{|p|_{q}}{2^{\ell}}\right)\right]^{2}\chi_{\ell'}u\,,
\end{eqnarray*}
and the right-hand side is a finite sum with a number of terms
uniformly bounded with respect to $\ell'\in \mathcal{L}$\,. Thus
$$
\|\chi_{\ell'}(P_{\pm,Q,g}-z)u\|_{L^{2}(Q;\mathcal{H}^{s'})}^{2}\leq
C_{\chi}^{1}\sum_{\chi_{\ell}\chi_{\ell'}\neq 0}
\|(P_{\pm,Q,g}-z)\chi_{\ell}u\|_{L^{2}(Q;\mathcal{H}^{s'})}^{2}+\|\chi_{\ell'}u\|^{2}
$$
holds for all $\ell'\in \mathcal{L}$ and summing over $\ell'\in
\mathcal{L}$ gives
\begin{eqnarray*}
  \|(P_{\pm,Q,g}-z)u\|_{L^{2}(Q;\mathcal{H}^{s'})}^{2}
&\leq&
  C_{\chi}^{2}\sum_{\ell'\in
    \mathcal{L}}\|\chi_{\ell'}(P_{\pm,Q,g}-z)u\|_{L^{2}(Q;\mathcal{H}^{s'})}^{2}\\
&\leq& C_{\chi}^{3}\left[\sum_{\ell'\in
    \mathcal{L}}\|(P_{\pm,Q,g}-z)\chi_{\ell'}u\|_{L^{2}(Q;\mathcal{H}^{s'})}^{2}+
\|\chi_{\ell'}u\|^{2}\right]\,.
\end{eqnarray*}
We have proved 
\begin{multline*}
  \|(P_{\pm,Q,g}-z)u\|_{L^{2}(Q;\mathcal{H}^{s'})}^{2}+\|u\|_{L^{2}(Q;\mathcal{H}^{1})}^{2}
\\
\leq C_{\chi}\left[\sum_{\ell\in
    \mathcal{L}}\|(P_{\pm,Q,g}-z)\chi_{\ell}u\|_{L^{2}(Q;\mathcal{H}^{s'})}^{2}
+\|\chi_{\ell}u\|_{L^{2}(Q;\mathcal{H}^{1})}^{2}\right]<+\infty\,.
\end{multline*}
\end{proof}

\subsection{Geometric KFP-operator on cylinders}
\label{se.cyl}
We consider here the case when $Q=\rz\times Q'$\,, with $Q'$ compact
without boundary, is endowed with the metric $g=1\oplus m(q')$\,,
$\partial_{q^{1}}m\equiv 0$\,. The contangent bundles are respectively
denoted by
$X=T^{*}Q$ and $X'=T^{*}Q'$\,.
The space $\mathcal{S}(X)$ is defined as the space of rapidly decaying
(in terms of $(q_{1},p_{1},p')$) $\mathcal{C}^{\infty}$-functions on
$X$\,. Its dual is denoted by $\mathcal{S}'(X)$\,.\\ 
On $Q$ we shall use a partition of unity $\sum_{\ell\in
  \zz}\chi^{2}(q^{1}-\ell)=1$
with $\chi\in \mathcal{C}^{\infty}_{0}((-\frac{1}{2},\frac{1}{2}))$\,.\\
On $Q$\,, the geometric Kramers-Fokker-Planck operator has the form
\begin{eqnarray*}
&&P_{\pm,Q,g}=\pm p_{1}\partial_{q^{1}}+\frac{-\partial_{p_{1}}^{2}+p_{1}^{2}+1}{2}+L=P_{\pm}^{L}
\\
\text{with}&&
L=K_{\pm,Q',m}=P_{\pm,Q',m}-\frac{1}{2}\\
\text{and}&&\mathfrak{L}=L^{2}(X')\quad,\quad D(L)=\left\{u\in
  L^{2}(Q',\mathcal{H}^{1})\,, P_{\pm,Q',m}u\in L^{2}(X')\right\}\,.
\end{eqnarray*}
described in Section~\ref{se.insvar} and in Section~\ref{se.genBC}.
It satisfies the Hypothesis~\ref{hyp.L1} because
$D(L)=D(L^{*})=D(P_{\pm,Q',m})$ contains 
$\mathcal{S}(X')$ and $\mathcal{C}^{\infty}_{0}(X')$ which are cores
for the self-adjoint operator $\Real(L)=\mathcal{O}_{Q',m}$\,.
The space
$L^{2}(\rz,dq^{1};\mathcal{H}^{1,L})=L^{2}(\rz,dq^{1};\mathcal{H}^{1,L^{*}})$
is nothing but $L^{2}(Q;\mathcal{H}^{1})$\,.\\
The Proposition~\ref{pr.idenDK} provides the maximal accretivity and
the corresponding  domain of $K_{\pm,Q,g}=P_{\pm,Q,g}=K_{\pm}^{L}$ while
Corollary~\ref{co.Leb} combined with the partition of unity
$\sum_{\ell\in \zz}\chi^{2}(q^{1}-\ell)$ allows to check the
subelliptic estimates of Hypothesis~\ref{hyp.whole}.
\begin{proposition}
\label{pr.cylacc} On the cylinder $(Q=\rz\times Q', g= 1\oplus m(q'))$\,,
the operator $K_{\pm,Q,g}-\frac{d}{2}$ defined by
\begin{eqnarray*}
  D(K_{\pm,Q,g})=\left\{u\in
    L^{2}(Q;\mathcal{H}^{1})\,,~P_{\pm,Q,g}u\in L^{2}(X)\right\}
\end{eqnarray*}
is maximal accretive with $K_{\pm,Q,g}^{*}=K_{\mp,Q,g}$\,.\\
It is the unique maximal accretive extension of $P_{\pm,Q,g}$ defined
on $\mathcal{C}^{\infty}_{0}(X)$ (or $\mathcal{S}(X)$) and it satisfies
$$
\forall u\in D(K_{\pm,Q,g})\,, \quad d\|u\|^{2}\leq
\|u\|_{L^{2}(Q;\mathcal{H}^{1})}^{2}
=\Real \langle u\,,\, (K_{\pm,Q,g}+\frac{d}{2})u\rangle\,.
$$
The estimate 
$$
\|(\pm
\mathcal{Y}_{\mathcal{E}}-i\lambda)u\|+\|\mathcal{O}_{Q,g}u\|+\langle
\lambda\rangle^{\frac{1}{2}}\|u\|+\|u\|_{H^{\frac{2}{3}}(Q;\mathcal{H}^{0})}\leq C\|(K_{\pm,Q,g}-i\lambda)u\|
$$
holds for all $u\in D(K_{\pm,Q,g})$ and all $\lambda\in \rz$\,.
In particular,
$$
D(K_{\pm,Q,g}^{*})=D(K_{\pm,Q,g})=\left\{u\in
  L^{2}(Q;\mathcal{H}^{1})\,,\, \mathcal{Y}_{\mathcal{E}}u\,,\,
  \mathcal{O}_{Q,g}u\in L^{2}(X)\right\}\,.
$$
\end{proposition}
\begin{remark}
  This result means that Hypothesis~\ref{hyp.whole} is verified with
  $L_{\pm}=K_{\pm,Q',m}$\,, $L^{2}(\rz,dq;\mathcal{H}^{1,L})=L^{2}(Q;\mathcal{H}^{1})$
and $\mathcal{Q}_{0}=H^{\frac{2}{3}}(Q;\mathcal{H}^{0})$\,.
\end{remark}
\begin{proof}
  The first part, including the identity for $\Real\langle u\,,\,
  K_{\pm,Q,g}u\rangle$\,,  is a direct application of
  Proposition~\ref{pr.idenDK}. Since
  $D(K_{\pm,\rz,1})\otimes^{alg}D(K_{\pm,Q',m})$ is dense in
  $D(K_{\pm,Q,g})=D(K_{\pm}^{L})$\,, with $L_{\pm}=K_{\pm,Q',m}$\,, and
$\mathcal{C}^{\infty}_{0}(\rz^{2})$
(resp. $\mathcal{C}^{\infty}_{0}(X')$) is dense in $D(K_{\pm,\rz,1})$
(resp. $D(L_{\pm})=D(K_{\pm,Q',m})$)\,, $\mathcal{C}^{\infty}_{0}(X)$ is a
core for $K_{\pm,Q,g}$\,.\\
For the subelliptic estimate, we use the partition of unity
$\sum_{\ell\in \zz}\chi^{2}_{\ell}\equiv 1$\,, with $\chi_{\ell}(q)=\chi(q^{1}-\ell)$\,.
For $u\in D(K_{\pm,Q,g})$ and $\lambda\in\rz$\,,
Proposition~\ref{pr.partq} applied with $s'=0$ gives
\begin{eqnarray*}
 C_{d}\|(K_{\pm,Q,g}-i\lambda)u\|^{2}
&\geq&
\|(P_{\pm,Q,g}-i\lambda)u\|^{2}+\|u\|_{L^{2}(Q;\mathcal{H}^{1})}^{2}
\\
&\geq&
C^{-1}\sum_{\ell\in \zz}\|(P_{\pm,Q,g}-i\lambda)\chi_{\ell}u\|^{2}
+\|\chi_{\ell}u\|_{L^{2}(Q;\mathcal{H}^{1})}^{2}\,.
\end{eqnarray*}
With $\supp \chi_{\ell}\subset \left\{\ell-\frac{1}{2}< q^{1}<
  \ell+\frac{1}{2}\right\}$ and the $q^{1}$-translation invariance, 
$(Q,g)$ can be replaced in the right-hand side by the compact manifold $(\sz^{1}\times Q',
g=1\oplus m)$\,. Thus Corollary~\ref{co.Leb} provides $\chi_{\ell}u\in
D(K_{\pm,\sz^{1}\times Q',g})$ and the uniform lower bound
\begin{multline*}
C'\|(P_{\pm,\sz^{1}\times Q',g}-i\lambda)\chi_{\ell}u\|^{2}
\geq 
\|\mathcal{O}_{\sz^{1}\times Q',g}\chi_{\ell}u\|^{2}+\|(\pm
\mathcal{Y}_{\mathcal{E}}-i\lambda)\chi_{\ell}u\|^{2}
\\+\langle \lambda\rangle\|\chi_{\ell}u\|^{2}
+
\|\chi_{\ell}u\|^{2}_{H^{\frac{2}{3}}(\sz^{1}\times Q';\mathcal{H}^{1})}\,.
\end{multline*}
Summing over $\ell\in\zz$ ends the proof.
\end{proof}
We end this section by checking Hypothesis~\ref{hyp.Q} with
$$
\mathcal{Q}=H^{t}(Q;\mathcal{H}^{0}))\,,\quad t\in[0,\frac{1}{9})\,. 
$$
The estimate \eqref{eq.hypQ2} is a consequence of
Proposition~\ref{pr.cylacc} with
\begin{multline*}
\forall u\in D(K_{\pm,Q,g})\,,\quad \langle \lambda\rangle^{\frac{1}{4}}
\|u\|_{\mathcal{Q}}\leq \langle
\lambda\rangle^{\frac{1}{4}}\|u\|_{H^{\frac{1}{3}}(Q;\mathcal{H}^{0})}\leq
(\langle\lambda\rangle^{\frac{1}{2}}\|u\|)^{\frac{1}{2}}\|u\|_{H^{\frac{2}{3}}(Q;\mathcal{H}^{0})}^{\frac{1}{2}}
\\
\leq C\|(K_{\pm,Q,g}-i\lambda)u\|\,.
\end{multline*}
We still have to check the estimate \eqref{eq.hypQ1}.
After identifying $Q'$ with the submanifold  $\left\{q^{1}=0\right\}$
of $Q$\,, the fiber bundle $\partial X=T_{Q'}^{*}Q$ is decomposed into
$$
\partial X=T^{*}_{Q'}Q=
T^{*}_{Q'}Q'
\oplus_{Q'}N^{*}_{Q'}Q=
T^{*}_{Q'}Q'\oplus_{Q'}(Q'\times \rz_{p_{1}})\,,
$$
The measures $|p_{1}|^{\nu}dq'dp=|p_{1}|^{\nu}dp_{1}dq'dp'$\,,
$\nu\in\rz$\,,  are
naturally defined on $T_{Q'}^{*}Q$ and the space
$L^{2}(\rz,\frac{dp_{1}}{|p_{1}|};\mathfrak{L})$
is nothing but $L^{2}(\partial X,\frac{dq'dp}{|p_{1}|})$\,.
Hypothesis~\ref{hyp.whole} for $K_{\pm,L}=K_{\pm,Q,g}$ with
$L_{\pm}=K_{\pm,Q',m}$ is validated by 
Proposition~\ref{pr.cylacc}. Therefore Proposition~\ref{pr.trL}
applies and we know that  a solution to the equation \eqref{eq.PLgd}, namely
\begin{eqnarray*}
  &&
  (P_{\pm,Q,g}-i\lambda)v=\gamma(q',p_{1},p')\delta_{0}(q^{1})\quad\text{in}~
\mathcal{S}'(\rz^{2};D(L^{*})')\subset
  \mathcal{S}'(X)\\
\text{with}&& \gamma\in L^{2}(\partial X,\frac{dq'dp}{|p_{1}|})\,,
\end{eqnarray*}
belongs to 
$L^{2}(\rz^{2},dqdp;\mathfrak{L})=L^{2}(X,dqdp)$\,. By adapting
Theorem~1.2 of \cite{Leb2} to the cylinder case, the above
equation admits a unique solution in $\mathcal{S}'(X)$ but we do not
really need this.
\begin{proposition}
\label{pr.Qverif} The Hypothesis~\ref{hyp.Q} is satisfied by 
the operator $P_{\pm}^{L}=P_{\pm,Q,g}$ and the space 
$\mathcal{Q}=\mathcal{H}^{t}(Q;L^{2}(\rz^{d-1}dp'))$ when
$t\in[0,\frac{1}{9})$\,.
\end{proposition}
\begin{proof}
The estimate \eqref{eq.hypQ2} has already been checked.\\
We must prove that a solution $v\in
L^{2}(\rz,dq;\mathcal{H}^{1,L})=L^{2}(Q;\mathcal{H}^{1})$
 to the
equation $(P_{\pm,Q,g}-i\lambda)v=\gamma
\delta_{0}(q^{1})$ satisfies
\eqref{eq.hypQ1}:
\begin{equation*}
\|v\|_{\mathcal{Q}}\leq
C_{\mathcal{Q},L}\left[\|\gamma\|_{L^{2}(\partial X,\frac{dq'dp}{|p_{1}|})}+
\|v\|_{L^{2}(Q;\mathcal{H}^{1})}\right]\,.
\end{equation*}
First note 
$$
\|\langle p_{1}\rangle v\|\leq \|v\|_{L^{2}(Q;\mathcal{H}^{1})}\,.
$$
Set $\tilde{v}=\langle p_{1}\rangle^{-\frac{1}{2}}v$\,. It still
belongs $L^{2}(Q;\mathcal{H}^{1})$ and solves 
\begin{eqnarray*}
  && (P_{\pm,Q,g}-i\lambda)\tilde{v}=\langle
p_{1}\rangle^{-\frac{1}{2}}\gamma \delta_{0}(q^{1})+ f\,,\\
\text{with}&&
f=\frac{p_{1}}{4\langle p_{1}\rangle^{\frac{5}{2}}}\partial_{p_{1}}v
+\frac{1-3p_{1}^{2}/2}{4\langle p_{1}\rangle^{\frac{9}{2}}}v\in L^{2}(X)\,.
\end{eqnarray*}
With the cut-off function $\theta\in
\mathcal{C}^{\infty}_{0}((-\frac{1}{2},\frac{1}{2}))$\,, one gets
\begin{eqnarray}
\label{eq.thtildv0}
  && (P_{\pm,Q,g}-i\lambda)(\theta(q^{1}) \tilde{v})= \langle
p_{1}\rangle^{-\frac{1}{2}}\gamma \delta_{0}(q^{1})+ \theta(q^{1})f +
p_{1}\theta'(q^{1})\tilde{v}\,,\\
\label{eq.thtildv1}
&& (P_{\pm,Q,g}-i\lambda)((1-\theta(q^{1})) \tilde{v})=
(1-\theta(q^{1}))f
- p_{1}\theta'(q^{1})\tilde{v}\,.
\end{eqnarray}
The right-hand side of \eqref{eq.thtildv1} belongs to $L^{2}(X,dqdp)$
while $(1-\theta(q^{1}))\tilde{v}$ belongs to
$L^{2}(Q;\mathcal{H}^{1})$\,. Hence Proposition~\ref{pr.cylacc} gives
\begin{multline*}
\|(1-\theta(q^{1}))\tilde{v}\|_{H^{\frac{2}{3}}(Q;\mathcal{H}^{0})}\leq
C
\left[\|\tilde{v}\|_{L^{2}(Q;\mathcal{H}^{1})}
+\|f\|\right]
\leq C'\|v\|_{L^{2}(Q;\mathcal{H}^{1})}\,.  
\end{multline*}
In \eqref{eq.thtildv0}, $(Q;g)$ can be replaced by the compact
manifold $(\sz^{1}\times Q',g)$ and Corollary~\ref{co.Leb} gives
\begin{align*}
\|\theta(q^{1})\tilde{v}\|_{H^{s+\frac{2}{3}}(Q;\mathcal{H}^{0})}
&\leq
C_{s}\left[\|\langle
  p_{1}\rangle^{-\frac{1}{2}}\gamma\|_{L^{2}(\partial X,dq'dp)}
+\|v\|_{L^{2}(Q;\mathcal{H}^{1})}\right]
\\
&\leq
C_{s}\left[\|\gamma\|_{L^{2}(\partial X,\frac{dq'dp}{|p_{1}|})}
+\|v\|_{L^{2}(Q;\mathcal{H}^{1})}\right]\,,
\end{align*}
as soon as $s<-\frac{1}{2}$ and $s\geq -\frac{2}{3}$\,.\\
We have proved  for $s\in [-\frac{2}{3},-\frac{1}{2})$
\begin{align*}
\|\langle
p_{1}\rangle^{-\frac{1}{2}}v\|_{H^{s+\frac{2}{3}}(Q;\mathcal{H}^{0})}
&=
\|\tilde{v}\|_{H^{s+\frac{2}{3}}(Q;\mathcal{H}^{0})}
\\
&\leq C_{s}'
\left[\|\gamma\|_{L^{2}(\partial X,\frac{dq'dp}{|p_{1}|})}+\|v\|_{L^{2}(Q,\mathcal{H}^{1})}\right]\,,
\end{align*}
while the first estimate is $\|\langle p_{1}\rangle
v\|_{H^{0}(Q;\mathcal{H}^{0})}\leq \|v\|_{L^{2}(Q;\mathcal{H}^{1})}$\,.
The interpolation inequality
$$
\|u\|_{H^{\frac{2\sigma}{3}}(Q;\mathcal{H}^{0})}\leq C_{\sigma}\|\langle
p_{1}\rangle^{-\frac{1}{2}}u\|_{H^{\sigma}(Q;\mathcal{H}^{0})}^{\frac{2}{3}}\|\langle
p_{1}\rangle u\|_{H^{0}(Q;\mathcal{H}^{0})}^{\frac{1}{3}}\,,
$$
applied with $\sigma=s+\frac{2}{3}\in [0,\frac{1}{6})$ yields the result.
\end{proof}
\subsection{Comments}
\label{se.comments}
There are essentially two strategies to prove maximal, i.e. with
optimal exponents, subelliptic estimates:
\begin{itemize}
\item a geometric approach following H{\"o}rmander
in 
\cite{HormIV}-Chap~27  (or Lerner in \cite{Ler}) 
based on a \underline{microlocal} reduction;
\item a more
algebraic approach after Rotschild-Stein in \cite{RoSt} or
Helffer-Nourrigat in \cite{HeNo} based on \underline{local} reduction
via Taylor approximations and solving algebraic models.
\end{itemize}
In \cite{Leb2} (and even in \cite{Leb1} for non maximal estimates)
 Lebeau follows the first approach.
Actually, another proof of Theorem~1.2 of \cite{Leb2} is possible via
a local point of view on $T^{*}Q$ which involves a canonical
transformation associated with a solution to the Hamilton-Jacobi
equation $|\partial_{q}\varphi(q)|_{q}^{2}=Cte$\,. In the case with a
boundary, as it is well known in the analysis of propagation of
singularities for the wave equation (see
\cite{AnMe}\cite{Tay1}\cite{Tay2}\cite{MeSj1}\cite{MeSj2}), a solution $\varphi$
shows singularities in the presence of glancing (i.e. gliding or
grazing) rays.
A microlocal point of view seems  out of reach for general boundary value
problems and we will follow the second local approach. 
Our strategy, is to reach some
subelliptic estimates for the boundary value problem
with the most flexible \underline{local} approach.\\
Nevertheless,  verifying
Hypothesis~\ref{hyp.Q} requires strong enough subelliptic estimates
for whole space problems, in order to absorb the
$-\frac{1}{2}-\varepsilon$-Sobolev singularity of
$\delta_{0}(q^{1})$\,.
Referring to the result by Lebeau was the most efficient way. Since it
was written only for compact riemannian manifolds, the most direct
application including a translation invariance in the $q^{1}$-variable
is the one proposed in Subsection~\ref{se.cyl}. The local model of
boundary manifolds will be $(-\infty,0]\times \tz^{d-1}$
 rather the standard half-space $(-\infty,0]\times \rz^{d-1}$\,.

\section{Geometric KFP-operators on manifolds with boundary}
\label{se.multi}

The proof of Theorem~\ref{th.main0} and Theorem~\ref{th.mainA} will be done in several steps
relying on a careful analysis of local models for boundary manifolds.
\subsection{Review of notations and outline}
\label{se.revout}
A neighborhood of $q_{0}\in \partial Q$ in 
the manifold with boundary $\overline{Q}=Q\sqcup \partial Q$
(resp. $\overline{X}=T^{*}\overline{Q}$)
can be identified with a domain of  
$(-\infty,0]\times \tz^{d-1}$ (resp. 
$(-\infty,0]\times\tz^{2d-1}$), diffeomorphic to a domain of
$(-\infty,0]\times\rz^{d-1}$
(resp. $(-\infty,0]\times\rz^{2d-1}$),
 with the corresponding global
coordinates $(q^{1},q')$
(resp. $(q^{1},p_{1},q',p')$)\,.  The coordinate $q'\in\tz^{d-1}$ means
$q'\in\rz^{d-1}$ plus periodicity conditions and we assume
$q^{1}(q_{0})=0$ and $q'(q_{0})=0$\,.\\
The coordinate system $(q^{1},q')$ can be chosen so that the metric equals
\begin{equation}
\label{eq.metbord1}
g(q)=
\begin{pmatrix}
  1&0\\
0& m(q^{1},q')
\end{pmatrix}\,,
\end{equation}
and since it is a local description, we can assume
\begin{equation}
\label{eq.metbord2}
 m(q^{1},q')-m(0,q')\in
 \mathcal{C}^{\infty}_{0}((-\infty,0]\times\tz^{d-1};\mathcal{M}_{d-1}(\rz))\,.
\end{equation}
The kinetic energy, the vertical (in $p$) harmonic oscillator operator
and the Hamilton vector field on $\rz^{2d}$ are given by 
\begin{eqnarray*}
  && 2\mathcal{E}=|p|_{g(q)}^{2}=|p_{1}|^{2}+m^{ij}(q^{1},q')p'_{i}p'_{j}\,,\\
&&
\Delta_{p}=\partial_{p_{1}}^{2}+\partial_{p'_{i}}m_{ij}(q^{1},q')\partial_{p'_{j}}\quad,\quad
\mathcal{O}_{Q,g}=\frac{-\Delta_{p}+|p|^{2}_{q}}{2}\,,
\\
&&
\mathcal{Y}_{\mathcal{E}}=p_{1}\partial_{q^{1}}
+
m^{ij}(q^{1},q')p_{i}'\partial_{q^{'j}}
-\frac{1}{2}\partial_{q^{'k}}m^{ij}(q^{1},q')p_{i}'p_{j'}\partial_{p_{k}'}\\
&&\hspace{4cm}-\frac{1}{2}\partial_{q^{1}}m^{ij}(q^{1},q')p_{i}'p_{j}'\partial_{p_{1}}\,.
\end{eqnarray*}
Since we will work with various metrics which fulfill
\eqref{eq.metbord1}\eqref{eq.metbord2} with the same $m(0,q')$\,,
 the previous notation
$|p|_{q}$ is replaced by $|p|_{g(q)}$\,. When $q=(0,q')\in \partial
Q=Q'$\,, the equalities $g(0,q')=m(0,q')=g_{0}(q)$\,,
$|p|_{g(q)}=|p|_{g_{0}(q)}$ allow to use $|p|_{q}$\,.\\
Following the Definition~\ref{de.geKFP}, the corresponding
Kramers-Fokker-Planck operator equals 
\begin{align}
\nonumber
P_{\pm,Q,g}&= \pm\mathcal{Y}_{\mathcal{E}}+\mathcal{O}_{Q,g}\,,\\
 \label{eq.gKFPhalf}
&=\pm p_{1}\partial_{q^{1}}+\frac{-\partial_{p_{1}}^{2}+|p_{1}|^{2}}{2}+
P_{\pm,Q', m(q^{1})}-\frac{1}{2}\partial_{q^{1}}m^{ij}(q^{1},q')p_{i}'p_{j}'\partial_{p_{1}}\,,
\end{align}
with here $Q'=\tz^{d-1}$\,, $Q=(-\infty,0]\times Q'$\,.
The operator $P_{\pm, Q',m(q^{1})}$ is the geometric Kramers-Fokker-Planck
operator on $X'=T^{*}Q'$ associated with the metric $m(q^{1})$ on
$Q'$\,.\\
The phase-space $X=T^{*}Q$ is endowed with the metric $g\oplus g^{-1}$
and the symplectic volume $dqdp$\,.
We shall use the notation $H^{s}(\overline{Q};\mathcal{H}^{s'})$\,,
$\mathcal{S}(X)$\,, $\mathcal{S}'(X)$ introduced in
Section~\ref{se.geoKFP}, extended to $\mathfrak{f}$-valued functions
or distributions, where $\mathfrak{f}$ is a complex Hilbert space.\\
When $Q'=\tz^{d-1}$\,, $Q=(-\infty,0]\times Q'$ (resp.
$Q=\rz\times Q'$) and $\partial_{q^{1}}m\equiv 0$\,, we shall use the results of
Section~\ref{se.insvar} and Section~\ref{se.genBC} with the
identification
\begin{eqnarray*}
  && P_{\pm,Q,g}= P_{\pm}^{L} =\pm
  p_{1}\partial_{q^{1}}+\frac{-\partial_{p_{1}}^{2}+|p_{1}|^{2}+1}{2}+L_{\pm}\\
&& L_{\pm}=K_{\pm,Q',m}-\frac{1}{2}\quad,\quad
\mathfrak{L}=L^{2}(X',dq'dp';\mathfrak{f})\,,\\
&&L^{2}(\rz_{-},dq^{1};\mathcal{H}^{1,L})=L^{2}(Q;\mathcal{H}
^{1})\\
\text{resp.}&&
L^{2}(\rz,dq^{1};\mathcal{H}^{1,L})
=L^{2}(Q;\mathcal{H}^{1})\,,
\end{eqnarray*}
where $K_{\pm,Q',m}=K_{\mp,Q',m}^{*}$ is the maximal accretive realization of
$P_{\pm,Q',m}$ provided by Proposition~\ref{pr.cylacc}. We checked in
Proposition~\ref{pr.cylacc} and Proposition~\ref{pr.Qverif} that
Hypothesis~\ref{hyp.whole} and Hypothesis~\ref{hyp.Q} are true
respectively with $\mathcal{Q}_{0}=H^{\frac{2}{3}}(\rz\times Q';\mathcal{H}^{0})$
and $\mathcal{Q}=H^{t}(\rz\times Q';\mathcal{H}^{0})$\,, $t\in [0,\frac{1}{9})$\,.\\
The boundary $\partial Q$ is identified with $Q'$ and the measure
$|p_{1}|dq'dp$ is well defined on  $\partial X=T^{*}_{Q'}Q$\,.
We recall that the trace of $u$ at $q^{1}=0$ is denoted by $\gamma
u$\,, that  $j$ is a unitary involution on $\mathfrak{f}$\,. 
The boundary conditions are written
\begin{eqnarray*}
&& \gamma_{odd}u=\pm\sign(p_{1})A\gamma_{ev}u\quad,\quad
  \gamma_{ev,odd}u=\Pi_{ev,odd}u\,,\\
\text{with}&&
\Pi_{ev}\gamma(q',p_{1},p')=\frac{\gamma(q',p_{1},p')+j\gamma(q',-p_{1},p')}{2}\,,\\
\text{and}&&
\Pi_{odd}\gamma(q',p_{1},p')=\frac{\gamma(q',p_{1},p')-j\gamma(q',-p_{1},p')}{2}\,;\\
\text{or}
&& \Pi_{\mp}\gamma u=\frac{1-A}{1+A}\Pi_{\pm}\gamma u\,,\\
\text{with}&& \Pi_{+}=\Pi_{ev}+\sign(p)\Pi_{odd}\quad,\quad \Pi_{-}=\Pi_{ev}-\sign(p)\Pi_{odd}\,.
\end{eqnarray*}
The operator $(A,D(A))$ commutes with  $\Pi_{ev,odd}$ and it is
 bounded and accretive in
$L^{2}(\partial X,|p_{1}|dq'dp;\mathfrak{f})$\,.
The norm of $A$ in $\mathcal{L}(L^{2}(\partial
X,|p_{1}|dq'dp;\mathfrak{f}))$ will be denoted $\|A\|$ and
the framework of Hypothesis~\ref{hyp.Abdd} contains the alternative
\begin{eqnarray*}
\textbf{either} &&c_{A}=\min \sigma(\Real A)>0\,;\\
\textbf{or}&& A=0\,.
\end{eqnarray*}
The treatment of general metrics on half-cylinders, $Q=\rz_{-}\times
Q'$ with $Q'=\partial Q$ compact, we will assume the
commutation
\begin{equation}
  \label{eq.commutAp}
\left[A,e^{it|p|_{q}^{2}}\right]=0\quad, \forall t\in\rz\,.  
\end{equation}
By writing the momentum $p\in T^{*}_{q}Q$\,, $p=r\omega$ with
$r=|p|_{q}$ and $\omega \in S^{*}_{q}Q$\,, the space $L^{2}(\partial
X, |p_{1}|dq'dp;\mathfrak{f})$ equals
$$
L^{2}(T^{*}_{\partial Q}Q,|p_{1}|dq'dp;\mathfrak{f})=L^{2}((0,+\infty),r^{d-1}dr;
L^{2}(S^{*}_{\partial Q}Q, |\omega_{1}|dq'd\omega;\mathfrak{f}))\,.
$$
The commutation \eqref{eq.commutAp} combined with
Hypothesis~\ref{hyp.Abdd} means that $A=A(|p|_{q})$ acts as a
multiplication in the radial coordinate $r=|p|_{q}$ with
\begin{eqnarray}
\label{eq.Ar1}
  && \|A(r)\|_{\mathcal{L}(L^{2}(S^{*}_{\partial Q}Q,
    |\omega_{1}|dq'd\omega;\mathfrak{f}))}\leq
  \|A\|\quad\text{for~a.e.}~r>0\,,\\
\label{eq.Ar2}
\text{with~either}&& \min \sigma(\Real A(r))\geq c_{A}>0\quad\text{for~a.e.}~r>0\,,\\
\label{eq.Ar3}\text{or}&& A(r)=0 \quad\text{for~a.e.}~r>0\,.
\end{eqnarray}
Finally the analysis of a general manifold with boundary, will be
studied with a spatial partition of unity. For this final analysis, we
will assume
$A=A(q,|p|_{q})$ with 
\begin{eqnarray}
\label{eq.Arq1}
  && \hspace{-1cm}\|A(q,r)\|_{\mathcal{L}(L^{2}(S^{*}_{q}Q,
    |\omega_{1}|d\omega ;\mathfrak{f}))}\leq
  \|A\|\quad\text{for~a.e.}~(q,r)\in \partial Q\times \rz_{+}\,,\\
\label{eq.Arq2}
\text{with~either}\hspace{-0.5cm}&& \min \sigma(\Real A(q,r))\geq
c_{A}>0\quad\text{for~a.e.}~(q,r)\in\partial Q\times \rz_{+}\,,\\ 
\label{eq.Arq3}
\text{or}&& A(q,r)=0 \quad\text{for~a.e.}~(q,r)\in\partial Q\times \rz_{+}\,.
\end{eqnarray}
\\
Here is the outline of the proof:
We will carefully study the case of half-cylinders
$Q=\rz_{-}\times Q'$ with $Q'=\tz^{d-1}$\,.
The case 
 when $\partial_{q^{1}}m\equiv 0$ 
will be an application of Section~\ref{se.genBC} and the result will
hold under Hypothesis~\ref{hyp.Abdd} which is more general than the
one of Theorem~\ref{th.mainA}.  
By assuming $A=A(|p|_{q})$ with
\eqref{eq.Ar1}\eqref{eq.Ar2}\eqref{eq.Ar3} and with a dyadic partition
of unity in $p$\,, the corresponding subelliptic estimates will be
written in a parameter dependent form. This allows an accurate
parameter dependent analysis  of some relatively bounded perturbations.
In a second step, the dyadic partition of unity for a general metric
$m$ on a half-cylinder, $\rz_{-}\times \tz^{d-1}$\,, 
and a non symplectic change of variable  in $X=T^{*}Q$ near
$\partial X$ will relate the general problem to the previous
pertubative analysis. Finally the case of a general manifold will be
treated by gluing the local models $Q=\rz_{-}\times \tz^{d-1}$ 
with a spatial partition of unity, by assuming
$A=A(q,|p|_{q})$ with \eqref{eq.Arq1}\eqref{eq.Arq2}\eqref{eq.Arq3}.

\subsection{Half-cylinders with $\partial_{q^{1}}m\equiv 0$}
\label{se.specmet}

Let us consider the case when $Q=(-\infty,0]\times Q'$ with
$Q'=\tz^{d-1}$ and
a specific metric $g_{0}=1\oplus^{\perp}m(q)$ satisfying
 $\partial_{q^{1}}m\equiv 0$\,.
\begin{proposition}
\label{pr.plat}
Let $Q=(-\infty,0]\times Q'$ with $Q'=\tz^{d-1}$ be endowed with the
metric $g_{0}=1\oplus m$ with $\partial_{q^{1}}m\equiv 0$\,.\\
 Assume Hypothesis~\ref{hyp.Abdd} for the operator $A$\,.\\
The operator $K_{\pm,A,g_{0}}-\frac{d}{2}$ defined by
\begin{eqnarray*}
  &D(K_{\pm,A,g_{0}})=
\left\{\begin{array}[c]{ll}
u\in L^{2}(Q;\mathcal{H}^{1})\,,
&P_{\pm,Q,g_{0}}u\in
    L^{2}(X,dqdp;\mathfrak{f}),\\
&\gamma u\in L^{2}(\partial X, |p_{1}|dq'dp;\mathfrak{f})\,,\\
&\gamma_{odd}u=\pm\sign(p_{1}) A\gamma_{ev}u
    \end{array}
\right\}\,,\\
&\forall u\in D(K_{\pm,A,g_{0}})\,,~K_{\pm,A,g_{0}}u=P_{\pm,Q,g_{0}}u=(\pm\mathcal{Y}_{\mathcal{E}}+\mathcal{O}_{Q,g_{0}})u\,,
\end{eqnarray*}
is maximal accretive, with
$$
\forall u\in D(K_{\pm,A,g_{0}})\,,\quad
\|u\|_{L^{2}(Q;\mathcal{H}^{1})}^{2}+\Real\langle
\gamma_{ev}u\,,\, A\gamma_{ev}u\rangle = \Real\langle u\,,\, (K_{\pm,A,g_{0}}+\frac{d}{2})u\rangle\,.
$$
The adjoint $K_{\pm,A,g_{0}}^{*}$ equals $K_{\mp,A^{*},g_{0}}$\,.\\
Respectively to the cases $A\neq 0$ and $A=0$\,, fix $t\in
[0,\frac{1}{9})$ and $\nu=\frac{1}{4}$ (resp. $t=\frac{2}{3}$ and $\nu=0$)\,.
There exist $C_{t}\geq 1$ and $C\geq 1$ independent of $t$\,, such that
\begin{multline}
\label{eq.subellg0}
\langle
\lambda\rangle^{\frac{1}{4}}\|\gamma u\|_{L^{2}(\partial X,|p_{1}|dq'dp;\mathfrak{f})}
+
\langle \lambda\rangle^{\frac 1 2}\|u\|
+
\langle \lambda\rangle^{\frac{1}{4}}\|u\|_{L^{2}(Q;\mathcal{H}^{1})}
+
C_{t}^{-1}\langle\lambda\rangle^{\nu}\|u\|_{H^{t}_{q}(\overline{Q};\mathcal{H}^{0})}
\\
\leq C\|(K_{\pm,A,g_{0}}-i\lambda)u\|
\end{multline}
for all $\lambda\in\rz$ and all $u\in D(K_{\pm,A,g_{0}})$\,.\\
For any $\Phi\in
\mathcal{C}^{\infty}_{b}((-\infty,0])$ such that $\Phi(0)=0$ there
exists $C_{\Phi}>0$ and $C'>0$ independent of $\Phi$\,, such that
$$
\|\Phi(q^{1})\mathcal{O}_{Q,g_{0}}u\|\leq C'\|\Phi\|_{L^{\infty}}\|(K_{\pm,A,g_{0}}-i\lambda)u\|+C_{\Phi}\|u\|\,,
$$
for all $\lambda\in\rz$ and all $u\in D(K_{\pm,A,g_{0}})$\,.\\
Finally in the case $A=0$\,, the set $\left\{u\in
  \mathcal{C}^{\infty}_{0}(\overline{X};\mathfrak{f})\,,
  \gamma_{odd}u=0\right\}$ is dense in $D(K_{\pm,0,g_{0}})$ endowed
with the graph norm.
\end{proposition}
\begin{proof}
Tensorizing with $\mathfrak{f}$ does not change the scalar results of Subsection~\ref{se.cyl}.\\
In Proposition~\ref{pr.cylacc} and Proposition~\ref{pr.Qverif},
Hypothesis~\ref{hyp.whole} and Hypothesis~\ref{hyp.Q} have been
checked with $L_{\pm}=K_{\pm,Q',m}$\,,
$\mathcal{Q}_{0}=H^{\frac{2}{3}}(\rz\times Q';\mathcal{H}^{0})$\,, 
$\mathcal{Q}=H^{t}(\rz\times Q';\mathcal{H}^{0})$\,,
 while Hypothesis~\ref{hyp.Abdd}
is assumed. We can refer to the results of  Section~\ref{se.insvar} and
Section~\ref{se.genBC}.
The maximal accretivity, the integration by part identity and the
identification of $K_{\pm,A,g_{0}}$ are provided by
Proposition~\ref{pr.maxaccKAL} when $A\neq 0$ and by
Proposition~\ref{pr.idenDKAId} when $A=0$\,.
The definition of the Sobolev spaces
$H^{s}(\overline{Q};\mathcal{H}^{0})$ says
$$
\|u\|_{H^{s}(\overline{Q};\mathcal{H}^{0})}\leq C_{s}\|\Sigma
u\|_{H^{s}(\rz\times Q';\mathcal{H}^{0})}\,.
$$
Hence the subelliptic estimate \eqref{eq.subellg0} and   the upper
bound
for $\|\Phi(q^{1})\mathcal{O}_{Q,g_{0}}\|$ are proved in Proposition~\ref{pr.regL1}
when $A\neq 0$ and in Proposition~\ref{pr.trKL0} when $A=0$ (actually
it works
with $\Phi(q^{1})=1$ in this case).\\
Because $\mathcal{C}^{\infty}_{0}(\rz\times Q';\mathfrak{f})$ is a
core for $K_{\pm,g_{0}}$\,,
 Proposition~\ref{pr.trKL0} also implies that the set 
$$
\mathcal{D}_{0}=\left\{u\in L^{2}(X,dqdp;\mathfrak{f})\,, \Sigma u\in
  \mathcal{C}^{\infty}_{0}(\rz\times Q';\mathfrak{f})\right\}
$$
is dense in $D(K_{\pm,0,g_{0}})$ endowed with its graph norm. This set is contained in
$$
\mathcal{D}_{1}=
\left\{u\in \mathcal{C}^{\infty}_{0}(\overline{X};\mathfrak{f})\,,
  \gamma_{odd}u=0\right\}\,,
$$
and
$\mathcal{D}_{0}\subset \mathcal{D}_{1}\subset D(K_{\pm,0,g_{0}})$
implies that $\mathcal{D}_{1}$ is dense in $D(K_{\pm,0,g_{0}})$\,.
\end{proof}

\subsection{Dyadic partition of unity and rescaled estimates}
\label{se.dyadresc}
We still work on $\overline{Q}=(-\infty,0]\times Q'$ with
$Q'=\tz^{d-1}$\,, with global coordinates $(q,p)\in
\overline{\rz_{-}}\times \tz^{d-1}\times \rz^{d}$\,, 
with a metric $g=1\oplus m$ which satisfies
\eqref{eq.metbord1}\eqref{eq.metbord2}. The notation $g_{0}$ is
specific to the case $\partial_{q^{1}}m\equiv 0$\,.
The operator $A$ is assumed to satisfy Hypothesis~\ref{hyp.Abdd} and
the commutation 
\eqref{eq.commutAp} with $|p|_{q}^{2}$\,, which can be
written as \eqref{eq.Ar1}\eqref{eq.Ar2}\eqref{eq.Ar3}\,.

\subsubsection{Rescaling}
\label{se.resca}
When $\partial_{q^{1}}m\equiv 0$\,, the maximal accretive realization
$K_{\pm,A,g_{0}}$ is given by Proposition~\ref{pr.plat}. The strengthened assumptions on
$A$ in  \eqref{eq.commutAp} yields
$$
(u\in D(K_{\pm,A,g_{0}}))\rightarrow\left(\chi(|p|_{g_{0}(q)}^{2})u\in D(K_{\pm,A,g_{0}})\right)
$$
for all $\chi\in \mathcal{C}^{\infty}_{0}(\rz)$\,.
Therefore we can use a dyadic partition of unity in the momentum
variable, $\chi=(\chi_{\ell})_{\ell\in\nz}$\,, like in
Subsection~\ref{se.momentpart}. This will be used for  the two metrics
$\tilde{g}=g$ and $\tilde{g}=g_{0}$:
\begin{eqnarray*}
  && \chi_{\ell}(q,p)=
  \tilde{\chi}_{1}(2^{-\ell}|p|_{\tilde{g}(q)})~\text{for}~\ell\in\nz^{*}
\quad,\quad\chi_{0}(q,p)=\tilde{\chi}_{0}(|p|_{\tilde{g}(q)})\,,\\
&&\tilde{\chi}_{0}\in \mathcal{C}^{\infty}_{0}(\rz)\quad,\quad 
\tilde{\chi}_{1}\in \mathcal{C}^{\infty}_{0}((0,+\infty))\,,\\
&&\sum_{\ell\in\mathcal{L}}\chi_{\ell}^{2}(q,p)\equiv 1\quad,\quad \mathcal{L}=\nz\,.
\end{eqnarray*}
The equivalence of norms 
\begin{align*}
&\left(\frac{\sum_{\ell\in\mathcal{L}}\|\chi_{\ell}u\|_{L^{2}(Q;\mathcal{H}^{s'})}^{2}}{\|u\|_{L^{2}(Q;\mathcal{H}^{s'})}^{2}}\right)^{\pm
1}\leq C_{g,\chi,s'}\quad, s'\in\rz\,,
\\
&
\left(
\frac{\sum_{\ell\in \mathcal{L}}
\|\chi_{\ell}u\|_{H^{s}(\overline{Q};\mathcal{H}^{0})}^{2}}{\|u\|_{H^{s}(\overline{Q};\mathcal{H}^{0})}^{2}}
\right)^{\pm 1} 
\leq C_{g,\chi,s}\,,\quad s\in [-1,1]\,,
\end{align*}
are given by Lemma~\ref{le.partp} and Lemma~\ref{le.partpq}.\\
With $\partial X=\tz^{d-1}\times \rz^{d}$\,, the two metrics $g$ and
$g_{0}$
coincide and the traces satisfy
$$
\|\gamma\|_{L^{2}(\partial X,|p_{1}|dq'dp;\mathfrak{f})}^{2}=\sum_{\ell\in\nz}\|\chi_{\ell}u\|_{L^{2}(\partial X,|p_{1}|dq'dp;\mathfrak{f})}^{2}\,.
$$
By working with $\tilde{g}=g_{0}$\,,  Proposition~\ref{pr.plat} ensures
$$
\forall \lambda\in\rz\,, 
\forall w\in D(K_{\pm,A,g_{0}})\,,\quad
\|w\|_{L^{2}(Q;\mathcal{H}^{1})}\leq C\|(K_{\pm,A,g_{0}}-i\lambda)w\|\,.
$$
Therefore the result of Proposition~\ref{pr.partp} reads
$$
\forall u\in D(K_{\pm,A,g_{0}})\,,\quad
\left(\frac{\sum_{\ell\in\rz}\|(K_{\pm,A,g_{0}}-i\lambda)\chi_{\ell}u\|^{2}
}{\|(K_{\pm,A,g_{0}}-i\lambda)u\|^{2}}\right)^{\pm 1}\leq C_{\chi}\,.
$$
After the unitary change of variables 
$$
u(q,p)=2^{-\ell d/2}v(q,\frac{p}{2^{\ell}})\quad,\quad f(q,p)=2^{-\ell
d/2}\varphi(q,\frac{p}{2^{\ell}})\,,
$$
the equation $(K_{\pm,A,g_{0}}-i\lambda)u=f$\,, i.e. the boundary value problem
$$
\left\{
  \begin{array}{l}
(P_{\pm,Q,g_{0}}-i\lambda)u=f\\
\gamma_{ev}u=\pm\sign(p_{1}) A(|p|_{q})\gamma_{ev}u\,,
\end{array}
\right.
$$
becomes
$$
\left\{
  \begin{array}{l}
h^{-1}(P_{\pm,Q,g_{0}}^{h}-i\lambda h)v=\varphi\\
\gamma_{ev}v=\pm\sign(p_{1}) A^{h}(|p|_{q})\gamma_{ev}v\,,
\end{array}
\right.
$$
where $P_{\pm,Q,g_{0}}^{h}$\,, $\mathcal{O}_{Q,g_{0}}^{h}$ and $A^{h}$
are defined according to the following general definition.
\begin{definition}
\label{de.resc}
For a metric $\tilde{g}$ on $(-\infty,0]\times Q'$ or $\rz\times Q'$ with
$Q'=\tz^{d}$ which satisfies \eqref{eq.metbord1}\eqref{eq.metbord2},
a bounded operator  $A$ on
$L^{2}(\partial X,|p_{1}|dq'dp;\mathfrak{f})$ and $h>0$\,, we set
\begin{eqnarray}
  \label{eq.rescaledP}
  &&P_{\pm,Q,\tilde{g}}^{h}=\pm\sqrt{h}\mathcal{Y}_{\mathcal{E}}+\mathcal{O}_{Q,\tilde{g}}^{h}\quad,\quad
  h=2^{-2\ell}\,,\\
\label{eq.rescaledO}
&&
\mathcal{O}_{Q,\tilde{g}}^{h}=
\frac{1}{2}\left[-h^{2}\Delta_{p}+|p|_{\tilde{g}(q)}^{2}\right]=
\frac{1}{2}\left[-h^{2}\partial_{p}^{T}\tilde{g}(q)\partial_{p}+p^{T}\tilde{g}^{-1}(q)p\right]\,,
\\
\label{eq.rescaledA}
  \text{and}&&
A^{h}(|p|_{\tilde{g}(q)})=A(h^{-\frac{1}{2}}|p|_{\tilde{g}(q)})\,.
\end{eqnarray}
\end{definition}
The  operator $K_{\pm,A,g_{0}}^{h}=P_{\pm,Q,g_{0}}^{h}$\,, 
with the domain $D(K_{\pm,A,g_{0}}^{h})$ 
characterized by
$$
\begin{array}[c]{ll}
u\in L^{2}(Q;\mathcal{H}^{1})\,,
&P_{\pm,Q,g_{0}}^{h}u\in
    L^{2}(X,dqdp;\mathfrak{f}),\\
&\gamma u\in L^{2}(\partial X, |p_{1}|dq'dp;\mathfrak{f})\,,\\
&\gamma_{odd}u=\pm\sign(p_{1}) A^{h}\gamma_{ev}u
    \end{array}
$$
is maximal accretive
with
\begin{multline*}
\frac{\left[\|g_{0}^{1/2}(q)(h\partial_{p})v\|^{2}+\|g_{0}^{-1/2}(q)pv\|^{2}\right]}{2}+\sqrt{h}\Real\langle
\gamma_{ev}v\,,\, A^{h}\gamma_{ev}v\rangle_{L^{2}(\partial X,|p_{1}|dq'dp;\mathfrak{f})}
\\
= \Real\langle
u\,,\, K_{\pm,A,g_{0}}^{h}u\rangle\,.
\end{multline*}
Actually, $h^{-1}K_{\pm,A,g_{0}}^{h}$ is by construction unitarily
equivalent to $K_{\pm,A,g_{0}}$\,, studied in Proposition~\ref{pr.plat}.\\
The other estimates of Proposition~\ref{pr.plat} now say
\begin{multline}
\label{eq.N2}
\langle
\lambda\rangle^{\frac{1}{4}}h^{3/4}\|\gamma v\|_{L^{2}(\partial X,|p_{1}|dpdq';\mathfrak{f})}
+\langle \lambda\rangle^{\frac 1 2}h\|v\|
+
\langle
\lambda\rangle^{\frac{1}{4}}\sqrt{h}\left[\|h\partial_{p}v\|+\|pv\|\right]
\\
+
C_{t}^{-1}\langle\lambda\rangle^{\nu}h\|u\|_{H^{t}(\overline{Q};\mathcal{H}^{0})}
\leq 
C\|(K_{\pm, A,g_{0}}^{h}-i\lambda h)u\|
\end{multline}
for all $\lambda\in \rz$ and all $u\in D(K_{\pm,A,g_{0}}^{h})$\,, with
$(\nu,t)\in \left\{\frac{1}{4}\right\}\times [0,\frac{1}{9})$ when
$A\neq 0$ and $(\nu,t)=(0,\frac{2}{3})$ when $A=0$\,.\\
 Similarly for
$\Phi\in \mathcal{C}^{\infty}_{b}((-\infty,0])$ such that
$\Phi(0)=0$\,, the estimates
\begin{equation}
  \label{eq.rescaledPhi}
\|\Phi(q^{1})\mathcal{O}_{Q,g_{0}}^{h}u\|\leq C\|\Phi\|_{L^{\infty}}\|(K_{\pm,A,g_{0}}^{h}-i\lambda h)u\|+C_{\Phi}h\|u\|
\end{equation}
holds for all $\lambda\in\rz$ and all $u\in D(K_{\pm,A,g_{0}}^{h})$\,.

\subsubsection{A perturbative result for $K_{\pm,A,g_0}^{h}$}
\label{se.perturbhalf}
Remember $\overline{Q}=(-\infty,0]\times\tz^{d}$ and $g_{0}=1\oplus m(q')$\,.\\
We shall consider now a perturbation of the operator $P_{\pm,Q,g_{0}}^{h}$:
\begin{eqnarray*}
  \label{eq.PbSTbord}
 P_{\pm,Q,g_{0}}^{h}(\mathcal{T})&=&\pm\sqrt{h}\left[\mathcal{Y}_{\mathcal{E}}+p_{i}b^{i}(q)
\theta(|p|_{g_{0}(q)})+S^{ij}_{k}(q)\theta(|p|_{g_{0}(q)})p_{i}p_{j}\partial_{p_{k}}\right]\\
&&\hspace{1cm}
+\frac{-(h\partial_{p})^{T}(g_{0}(q)+T_{1}(q))(h\partial_{p})+p^{T}(g_{0}^{-1}(q)+T_{2}(q)p}{2}\,,\\
\text{with}~\mathcal{Y}_{\mathcal{E}}&=&g_{0}^{ij}(q)p_{i}\partial_{q^{j}}-\frac{1}{2}\partial_{q^{k}}g_{0}^{ij}(q)p_{i}p_{j}\partial_{p_{k}}\\
&=&p_{1}\partial_{q^{1}}+m^{ij}(q')p_{i}'\partial_{{q'}^{j}}-\frac{1}{2}\partial_{{q'}^{k}}m^{ij}(q')p'_{i}p'_{j}\partial_{p'_{k}}\,,
\end{eqnarray*}
indexed by 
$$
\mathcal{T}=(b,S,\theta,T_{1},T_{2})\,.
$$
The collection $\mathcal{T}$ fullfills the following assumptions
\begin{equation}
\label{eq.condpert}
\left\{
  \begin{array}[c]{l}
b\in L^{\infty}(Q;\rz^{d})\quad,\quad S\in
  L^{\infty}(Q;\rz^{d^{3}})\,,\\
\theta\in \mathcal{C}^{\infty}_{0}([0,2R_{\theta}))\quad,\quad 0\leq
\theta\leq 1\quad,\quad \theta\equiv
1~\text{in}~\left[0,R_{\theta}\right]\,,  R_{\theta}\geq 1\\
 T_{1,2}=T_{1,2}^{T}\in
 L^{\infty}(Q;\mathcal{M}_{d}(\rz))\,,\\
 \|T_{1,2}\|_{L^{\infty}}\leq \varepsilon_{T}\quad,\quad 
|T_{1,2}(q)|\leq C_{T}|q^{1}|\,.
\end{array}
\right.
\end{equation}
\begin{proposition}
\label{pr.perturbflat}
Let $\mathcal{T}=(b,S,\theta,T_{1},T_{2})$ satisfy \eqref{eq.condpert}.\\
There exists $\varepsilon_{m,A}>0$ fixed by the metric
$g_{0}=1\oplus^{\perp} m(q')$ and the operator $A$\,,
such that the operator $K_{\pm,A,g_{0}}^{h}(\mathcal{T})=P_{\pm,Q,g_{0}}^{h}(\mathcal{T})$\,, with
the domain $D(K_{\pm,A,g_{0}}^{h}(\mathcal{T}))=D(K_{\pm,A,g_{0}}^{h})$\,,
 satisfies the following properties
uniformly w.r.t $h\in (0,1]$ for some  constant $C(\mathcal{T})$\,, 
as soon as $\varepsilon_{T}\leq
\varepsilon_{m,A}$ \,.\\
The operator $C(\mathcal{T})\sqrt{h}+K_{\pm,A,h}^{h}(\mathcal{T})$ is maximal accretive with 
\begin{multline*}
  \Real\langle v\,,\, K_{\pm,A,g_{0}}^{h}(\mathcal{T})v\rangle\geq \sqrt{h}\Real
  \langle\gamma_{ev}v\,,\,
  A^{h}\gamma_{ev}v\rangle_{L^{2}(\partial X,|p_{1}|dq'dp;\mathfrak{f})}
\\ \pm \sqrt{h}\Real
\langle v\,,\,
p_{i}b^{i}(q)\theta(|p|_{g_{0}(q)})v\rangle
\mp \frac{\sqrt{h}}{2}\langle v\,,\, S^{ij}_{k}(q)[\partial_{p_{k}}(p_{i}p_{j}\theta(|p|_{g_{0}(q)})]
v\rangle
\\
+\frac{\|(g_{0}+T_{1})^{1/2}(q)(h\partial_{p})v\|^{2}+\|(g_{0}^{-1}+T_{2})^{1/2}(q)pv\|^{2}}{2}\,.
\end{multline*}
By taking $(\nu_{1},\nu_{2},\nu_{3},t)\in
\left\{(\frac{1}{8},\frac{1}{8}, \frac{3}{4})\right\}\times[0,\frac{1}{18})$
when $A\neq 0$ and
$(\nu_{1},\nu_{2},\nu_{3},t)=(0,\frac{1}{4},\frac{5}{4},\frac{1}{3})$
when $A=0$\,, 
there exists a constant
${C'}_{t}>0$ such that
\begin{align*}
N(v,\lambda,t,h)=
\langle\lambda\rangle^{\frac 1 4} h\|v\| + 
&
\langle\lambda\rangle^{\frac 1 8}\sqrt{h}\left[\|h\partial_{p}v\|
+\|pv\|\right] +
{C'}_{t}^{-1}\langle \lambda\rangle^{\nu_{1}}h\|u\|_{H^{t}(\overline{Q};\mathcal{H}^{0})}\\
&
+\langle \lambda\rangle^{\nu_{2}}h^{\nu_{3}}\|\gamma v\|_{L^{2}(\partial X,|p_{1}|dq'dp;\mathfrak{f})}
\end{align*}
 is estimated as follows.\\
When  $h=1$ and $D(K_{\pm, A,g_{0}}^{1}(\mathcal{T}))=D(K_{\pm,A,g_{0}})$\,,
\begin{multline*}
\forall v\in D(K_{\pm,A,g_{0}})\,,~\forall \lambda\in\rz\,,\\
N(v,\lambda,t,1)\leq C(\mathcal{T})\left[\|(K_{\pm,A,g_{0}}^{1}(\mathcal{T})-i\lambda)v\|+\|v\|
\right]\,.
\end{multline*}
When  $h\leq \frac{1}{C(\mathcal{T})}$\,, the operator
$K_{\pm, A, g_{0}}^{h}(\mathcal{T})+\theta(4R_{\theta}^{2}|p|_{g_{0}(q)})-\frac{1}{C(\mathcal{T})}$ is maximal accretive with
\begin{multline*}
\forall v\in D(K_{A,h})\,,~\forall \lambda\in\rz\,,\\
N(v,\lambda,t,h)\leq
C(\mathcal{T})\|(K_{\pm,A,g_{0}}^{h}(\mathcal{T})+
\theta(4R_{\theta}^{2}|p|_{g_{0}(q)})-i\lambda h)v\|\,.
\end{multline*}
\end{proposition}
\begin{proof}
  When $\|T_{1,2}\|_{L^{\infty}}\leq \varepsilon_{T}$ is small enough
  $(g_{0}+T_{1})(q)$ and $(g_{0}^{-1}+T_{2})(q)$ are uniformly positive.
  The lower bound for $\Real\langle v\,,\, K_{\pm,A,g_{0}}^{h}(\mathcal{T})v\rangle$
  comes from the integration by part identity for $K_{\pm,A,g_{0}}^{h}$ (the
  rescaled version of $K_{\pm,A,g_{0}}$) after
  \begin{multline*}
K_{\pm,A,g_{0}}^{h}(\mathcal{T})-K_{\pm,A,g_{0}}^{h}
=
\pm \sqrt{h}p_{i}b^{i}(q)\theta(|p|_{g_{0}(q)})\pm
\sqrt{h}S^{ij}_{k}(q)p_{i}p_{j}\theta(|p|_{g_{0}(q)})\partial_{p_{k}} 
\\+\frac{-(h\partial_{p})^{T}T_{1}(q)(h\partial_{p})+p^{T}T_{2}(q)p}{2}\,,
\end{multline*}
where the right-hand side is a vertical operator sending
$L^{2}(Q;\mathcal{H}^{1})$ into its dual
$L^{2}(Q;\mathcal{H}^{-1})$ with
$$
[S^{ij}_{k}(q)p_{i}p_{j}\theta(|p|_{g_{0}(q)})\partial_{p}]^{*}=-
S^{ij}_{k}(q)p_{i}p_{j}\theta(|p|_{g_{0}(q)})\partial_{p}
-S^{ij}_{k}(q)\partial_{p_{k}}[p_{i}p_{j}\theta(|p|_{g_{0}(q)})]\,.
$$ 
This proves the accretivity of $C_{1}(\mathcal{T})\sqrt{h}+K_{\pm,A,g_{0}}^{h}(\mathcal{T})$
owing to the cut-off function $\theta$\,, when
$C_{1}(\mathcal{T})$ is chosen large enough.\\
More precisely, this inequality also implies
\begin{eqnarray*}
  \Real\langle v\,,\, 
\big[\theta(4R_{\theta}^{2}|p|_{g_{0}(q)})+K_{\pm,A,g_{0}}^{h}(\mathcal{T})\big]
  v\rangle
&\geq &
\frac{1}{2}\langle v\,,\,
\big[
\theta(4R_{\theta}^{2}|p|_{g_{0}(q)})
+
p^{T}(g_{0}^{-1}+T_{2})pv\big]\rangle
\\
&&\qquad\qquad
-C_{1}(\mathcal{T})\sqrt{h}\|v\|^{2}
\\
&\geq& \frac{1}{C_{2}(\mathcal{T})}\|v\|^{2}\,,
\end{eqnarray*}
for $C_{2}(\mathcal{T})$-large enough, when $h\leq h(\mathcal{T})$\,.\\
With the cut-off $\theta$\,, with the estimate \eqref{eq.rescaledPhi} and our
assumption on $T_{1}$ and $T_{2}$\,, $K_{\pm,A,g_{0}}^{h}(\mathcal{T})$ is a relatively
bounded perturbation of $K_{\pm,A,g_{0}}^{h}$ with relative bound $a<1$ as soon as 
$\varepsilon_{T}\leq \varepsilon_{m,A}$ and $\varepsilon_{m,A}$ is chosen small enough (the
constant $C$ in \eqref{eq.rescaledPhi} is fixed by the metric
$g_{0}=1\oplus m$ and $A$)\,.
Hence $C_{1}(\mathcal{T})\sqrt{h}+K_{\pm,A,g_{0}}^{h}(\mathcal{T})$ is a maximal accretive
operator.\\
We still have to check the upper bound for $N(v,\lambda,t,h)$\,.\\
\noindent\textbf{Case 1 ($h=1$):}
We already know that
$$
\forall u\in D(K_{\pm,A,g_{0}})\,,\quad
N(u,\lambda,t,1)\leq C_{m,A}\|(K_{\pm,A,g_{0}}-i\lambda)u\|
$$
while 
$$
\|(K_{\pm,A,g_{0}}^{1}(\mathcal{T})-K_{\pm,A,g_{0}})u\|\leq
C_{3}(\mathcal{T})\|u\|_{L^{2}(Q;\mathcal{H}^{1})}
+\|\Phi_{T}(q^{1})\mathcal{O}_{Q,g_{0}}u\|
$$
for some function $\Phi_{T}\in C^{\infty}_{0}((-\infty,0])$ such that
$\Phi_{T}(0)=0$ and $\|\Phi_{T}\|_{L^{\infty}}\leq
C_{m}\varepsilon_{T}$\,.
Again the estimate \eqref{eq.rescaledPhi} provides
$$
\|\Phi_{T}(q^{1})\mathcal{O}_{Q,g_{0}}u\|\leq
\frac{1}{2}\|(K_{\pm,A,g_{0}}-i\lambda)u\|+C'_{T,m,A}\|u\|\,, 
$$
when $\varepsilon_{T}\leq \varepsilon_{m,A}$ if $\varepsilon_{m,A}$ is chosen small enough.
Meanwhile the lower bound for $\Real\langle u\,,\,
K_{\pm,A,g_{0}}^{1}(\mathcal{T})u\rangle$ provides
$$
\|u\|_{L^{2}(Q;\mathcal{H}^{1})}^{2}\leq
C_{3}(\mathcal{T})\left[\|u\|+\|(K_{\pm, A, g_{0}}^{1}(\mathcal{T})-i\lambda)u\|\right]
$$
when $C_{3}(\mathcal{T})$ is chosen large enough.
This yields 
$$
N(u,\lambda,t,1)\leq
C_{4}(\mathcal{T})\left[\|u\|+\|(K_{\pm,A,g_{0}}^{1}(\mathcal{T})-i\lambda)u\|\right]\,. 
$$
\noindent\textbf{Case 2 ($h\leq h(\mathcal{T})$):}
We shall distinguish the case $(b,S)=(0,0)$ from the general case,
with the notations $\mathcal{T}_{0}=(0,0,\theta,T_{1},T_{2})$\,,
$\mathcal{T}=(b,S,\theta,T_{1},T_{2})$\,.
 For $C_{5}(\mathcal{T})$ large enough we know that
$(K_{\pm,A,g_{0}}^{h}(\mathcal{T})-z)^{-1}$ and $(K_{\pm,A,g_{0}}^{h}(\mathcal{T}_{0})-z)^{-1}$ are
well defined for $\Real z\leq-C_{5}(\mathcal{T})\sqrt{h}$\,.  The
proof relies 
 on the second resolvent formula
\begin{multline}
  \label{eq.secres}
(K_{\pm,A,g_{0}}^{h}(\mathcal{T}_{0})-z)^{-1}-(K_{\pm,A,g_{0}}^{h}(\mathcal{T})-z)^{-1}
\\
=
(K_{\pm,A,g_{0}}^{h}(\mathcal{T}_{0})-z)^{-1}
B
(K_{\pm,A,g_{0}}^{h}(\mathcal{T})-z)^{-1}
\end{multline}
with
$$
B=\pm\left[\sqrt{h}p_{i}b^{i}(q)\theta(|p|_{g_{0}(q)})+\frac{[S^{ij}_{k}(q)\theta(|p|_{g_{0}(q)})p_{i}p_{j}(h\partial_{p_{k}})]}{\sqrt{h}}\right]\,.
$$
after considering the case $(b,S)=(0,0)$\,.\\
\noindent\textbf{For $b=0$\,, $S=0$:} Take the quantity \eqref{eq.N2},
\begin{multline*}
N_{2}(v,\lambda,t,h)=
\langle\lambda\rangle^{\frac 1 2} h\|v\| + 
\langle\lambda\rangle^{\frac 1
  4}\sqrt{h}\left[\|h\partial_{p}v\|+\|pv\|\right] +
C_{2t}^{-1}\langle \lambda\rangle^{2\nu_{1}}h\|u\|_{H^{2t}(Q;\mathcal{H}^{1})}\\
+\langle \lambda\rangle^{\frac 1
  4}h^{3/4}\|\gamma v\|
_{L^{2}(\partial X,|p_{1}|dq'dp;\mathfrak{f})}\,,
\end{multline*}
where $(2\nu_{1},2t)\in \left\{\frac{1}{4}\right\}\times[0,\frac{1}{9})$
when $A\neq 0$ and $(2\nu_{1},2t)=(0,\frac{2}{3})$ when $A=0$\,.\\
We know
$$
\forall v\in D(K_{\pm, A,g_{0}}^{h})\,,\quad
N_{2}(v,\lambda,t,h)\leq C\|(K_{\pm,A,g_{0}}^{h}-i\lambda h)v\|\leq C\|(K_{\pm,A,g_{0}}^{h}-z)v\|
$$
for $z=-2C_{5}(\mathcal{T})+i\lambda$\,, $\lambda\in\rz$\,.
Like in \textbf{Case 1}\,, our assumptions on $T_{1},T_{2}$ imply that
the difference
$$
(K_{\pm,A,g_{0}}^{h}(\mathcal{T}_{0})-K_{\pm,A,g_{0}}^{h})v
=\frac{-(h\partial_{p})^{T}T_{1}(q)(h\partial_{p})+p^{T}T_{2}(q)p}{2}v
$$
satisfy
$$
\|(K_{\pm,A,g_{0}}^{h}(\mathcal{T}_{0})-K_{\pm,A,g_{0}}^{h})v\|
\leq C_{g}\|\Phi_{T}(q^{1})\mathcal{O}_{Q,g_{0}}^{h}v\|
$$
for some function $\Phi_{T}\in \mathcal{C}^{\infty}_{0}((-\infty,0])$
such that $\Phi_{T}(0)=0$ and $\|\Phi_{T}\|_{L^{\infty}}\leq
C_{m}\varepsilon_{T}$\,.
The estimate \eqref{eq.rescaledPhi} implies
$$
\|(K_{\pm,A,g_{0}}^{h}(\mathcal{T}_{0})-K_{\pm,A,g_{0}}^{h})v\|\leq \frac{1}{2}\|(K_{\pm,A,g_{0}}^{h}-i\lambda)v\|+C_{T,m,A}\|v\|
$$
if $\varepsilon_{m,A}>0$ is chosen small enough ($\varepsilon_{T}\leq
\varepsilon_{m,A}$)\,.
We infer
\begin{eqnarray*}
N_{2}(v,\lambda,t,h)
&\leq &
2C\|(K_{\pm,A,g_{0}}^{h}(\mathcal{T}_{0})-i\lambda
h)v\|+2C_{T,m,A}\|v\|\\
&\leq &C_{6}(\mathcal{T})\|(K_{\pm,A,g_{0}}^{h}(\mathcal{T}_{0})-z)v\|
\end{eqnarray*}
when $\Real z\leq -2C_{5}(\mathcal{T})$\,. \\
\noindent\textbf{$L^{2}$-estimates for general $(b,S,v)$:} 
The
condition $\Real z\leq -2C_{5}(\mathcal{T})$\,,implies
$$
z\not\in \sigma(K_{\pm,A,g_{0}}^{h}(\mathcal{T}_{0})\cup
\sigma(K_{\pm,A,g_{0}}^{h}(\mathcal{T})\,,
$$
and
\begin{eqnarray*}
\langle\lambda\rangle^{\frac 1 4}h^{1/2}\|v\|
&\leq&
(C_{5}(\mathcal{T})+\langle \lambda\rangle^{\frac 1 2}h)
\|v\|
\\
&\leq &C_{5}(\mathcal{T})\|v\|+N_{2}(v,\lambda,t,h)\\
&\leq&
C_{7}(\mathcal{T})\|(K_{\pm,A,g_{0}}^{h}(\mathcal{T}_{0})-z)v\|\,,
\end{eqnarray*}
while 
the lower bound for $\Real \langle v\,,\, K_{\pm,A, g_{0}}^{h}(\mathcal{T})v\rangle$
also gives
$$
\|h\partial_{p}v\|+\|pv\|\leq
C_{7}(\mathcal{T})\|(K_{\pm,A,g_{0}}^{h}(\mathcal{T})-z)v\|\,,
$$
for some $C_{7}(\mathcal{T})$ large enough and all $v\in D(K_{\pm,A,g_{0}}^{h})$\,.\\
Put in the second resolvent formula \eqref{eq.secres}, where the worst
term in $B$ is $h^{-\frac{1}{2}}(h\partial_{p})$\,, this implies
\begin{eqnarray*}
\|(K_{\pm,A,g_{0}}^{h}(\mathcal{T})-z)^{-1}\|
&\leq
&C_{7}(\mathcal{T})\|(K_{\pm,A,g_{0}}^{h}(\mathcal{T}_{0})-z)^{-1}\|
\times
\\
&&\qquad
\big[
h^{-1/2}\|(h\partial_{p})\circ
(K_{\pm,A,g_{0}}^{h}(\mathcal{T})-z)^{-1}\|+1
\big]
\\
&\leq&
\frac{C_{8}(\mathcal{T})}{2}\left[h^{-1}\langle\lambda\rangle^{-1/4} +
  h^{-1/2}\langle \lambda\rangle^{-1/4}\right]
\\
&\leq &C_{8}(\mathcal{T})h^{-1}\langle\lambda\rangle^{-1/4}\,,
\end{eqnarray*}
 for all $z$ such that $\Real z \leq -2C_{5}(\mathcal{T})$\,.\\
\noindent\textbf{$H^{t}(\overline{Q};\mathcal{H}^{0})$-estimate for general
  $(b,S,v)$:}
The estimate
$$
C_{5}(\mathcal{T})\|v\|+N_{2}(v,\lambda,t,h)\leq
C_{7}(\mathcal{T})\|(K_{\pm,A, g_{0}}^{h}(\mathcal{T}_{0})-z)v\|\quad,\quad \Real z\leq -2C_{5}(\mathcal{T})\,,
$$
implies, with the interpolation inequality $\|u\|_{H^{t}(\overline{Q};\mathcal{H}^{0})}\leq
\kappa_{t}\|u\|^{1/2}\|u\|_{H^{2t}(Q;\mathcal{H}^{0})}^{1/2}$\,,
\begin{eqnarray*}
\langle\lambda\rangle^{\nu_{1}}h^{1/2}\|v\|_{H^{t}(\overline{Q};\mathcal{H}^{0})}
&\leq &
\kappa_{t}C_{2t}\left[C_{5}(\mathcal{T})\|v\|+
C_{2t}^{-1}\langle
\lambda\rangle^{2\nu_{1}}
h\|u\|_{H^{2t}(Q;\mathcal{H}^{0})}\right]
\\
&\leq & \kappa_{t}C_{2t}
\left[C_{5}(\mathcal{T})\|v\|+
C_{2t}^{-1}\langle
\lambda\rangle^{2\nu_{1}}
h\|u\|_{H^{2t}(Q;\mathcal{H}^{0})}\right]\\
&\leq & \kappa_{t}C_{2t}C_{7}(\mathcal{T})\|(K_{\pm,A,g_{0}}^{h}(\mathcal{T}_{0})-z)v\|\,.
\end{eqnarray*}
Applying again the second resolvent formula \eqref{eq.secres} leads to
$$
C'_{t}\|(K_{\pm,A,g_{0}}^{h}(\mathcal{T})-z)^{-1}f\|_{H^{t}(\overline{Q};\mathcal{H}^{0})}\leq 
C_{9}(\mathcal{T})\langle \lambda\rangle^{-\nu_{1}}h^{-1}\|f\|\,,
$$
for some new function $t\in [0,\frac{1}{18})\to C_{t}'$\,, when $\Real
z\leq -2C_{5}(\mathcal{T})$ and $f\in L^{2}(X,dqdp;\mathfrak{f})$\,.\\
\noindent\textbf{Rough trace estimates for general $(b,S,v)$:}
The following argument holds for both cases $A=0$ and $A\neq 0$\,. The inequality
$$
\langle\lambda\rangle^{\frac{1}{4}}h^{\frac{3}{4}}
\|\gamma_{ev}v\|_{L^{2}(\partial X,|p_{1}|dq'dp;\mathfrak{f})}
\leq
N_{2}(v,\lambda,t,h)
\leq C_{7}(\mathcal{T})\|(K_{\pm,A,g_{0}}^{h}(\mathcal{T}_{0}-z)v\|
$$
for $v\in D(K_{\pm,A,g_{0}}^{h})$ and $\Real z\leq
-2C_{5}(\mathcal{T})$ 
can be written
$$
\|\gamma (K_{\pm,A,g_{0}}^{h}(\mathcal{T}_{0})-z)^{-1}f
\|_{L^{2}(\partial X,|p_{1}|dq'dp;\mathfrak{f})}
\leq C_{7}(\mathcal{T})\langle
\lambda\rangle^{-\frac{1}{4}}h^{-\frac{3}{4}}\|f\|\,,
$$
with $f\in L^{2}(X,dqdp;\mathfrak{f})$\,.
Inserted in the second resolvent formula \eqref{eq.secres}, this implies
$$
\|\gamma (K_{\pm,A,g_{0}}^{h}(\mathcal{T})-z)^{-1}f\|_{L^{2}(\partial X,|p_{1}|dq'dp;\mathfrak{f})}
\leq
C_{10}(\mathcal{T})
\langle \lambda\rangle^{-\frac{1}{4}}h^{-\frac{5}{4}}\|f\|
\,,
$$
when $\Real z\leq -2C_{5}(\mathcal{T})$\,.
Especially in the case $A=0$\,,
 this proves the upper bound of $\langle
 \lambda\rangle^{\nu_{2}}h^{\nu_{3}}\|\gamma v\|_{L^{2}(\partial
   X,|p_{1}|dq'dp;\mathfrak{f})}$
with $(\nu_{2},\nu_{3})=(\frac{1}{4},\frac{5}{4})$ when $v\in
D(K_{\pm,A,g_{0}}^{h}(\mathcal{T}))=D(K_{\pm,A,g_{0}}^{h})$\,.\\
\noindent\textbf{Using the integration by part inequality for general
  $(b,S,v)$:}
For $\Real z\leq -2C_{5}(\mathcal{T})$\,, the lower bound for $\Real\langle v\,,\,
K_{\pm,A,g_{0}}^{h}(\mathcal{T})v\rangle$ leads to
\begin{eqnarray*}
\|v\|\|(K_{\pm,A,g_{0}}^{h}(\mathcal{T})-z)v\|
&\geq&
\Real\langle v\,,\, (K_{\pm,A,g_{0}}^{h}(\mathcal{T})+2C_{5}(\mathcal{T}))v\rangle\\
&\geq& \sqrt{h}\Real\langle \gamma_{ev}v\,,\,
A^{h}\gamma_{ev}v\rangle_{L^{2}(\partial X,|p_{1}|dq'dp);\mathfrak{f}}\\
&&\quad
+\frac{1}{C_{11}(\mathcal{T})}\left[\|h\partial_{p}v\|^{2}+\|pv\|^{2}\right]\,.
\end{eqnarray*}
From 
$$
\|v\|\leq C_{8}(\mathcal{T})\langle\lambda\rangle^{-\frac 1 4}h^{-1}
\|(K_{\pm,A,g_{0}}^{h}(\mathcal{T})-z)v\|\,,
$$
we deduce
$$
\|h\partial_{p}v\|+\|p v\|\leq
\sqrt{2C_{11}(\mathcal{T})C_{8}(\mathcal{T})}\langle
\lambda\rangle^{-\frac{1}{8}}h^{-\frac{1}{2}}
\|(K_{\pm,A,g_{0}}^{h}(\mathcal{T})-z)v\|
$$
for all $v\in
D(K_{\pm,A,g_{0}}^{h}(\mathcal{T}))=D(K_{\pm,A,g_{0}}^{h})$\,.\\
When $A\neq 0$\,, this improves the estimate for $\|\gamma
v\|_{L^{2}(\partial X,|p_{1}|dq'dp;\mathfrak{f})}$ with
\begin{eqnarray*}
\frac{c_{A}}{2(1+\|A\|^{2})}\|\gamma v\|_{L^{2}(\partial
    X,|p_{1}|dq'dp;\mathfrak{f})}^{2}
&\leq& \langle \gamma_{ev}v\,,\, A^{h}\gamma_{ev}v\rangle
\\
&\leq & 
C_{8}(\mathcal{T})\langle \lambda\rangle^{-\frac{1}{4}}h^{-\frac{3}{2}}\|(K_{\pm,A,g_{0}}^{h}(\mathcal{T})-z)v\|^{2}\,,
\end{eqnarray*}
and we can take the pair of exponents
$(\nu_{2},\nu_{3})=(\frac{1}{8},\frac{3}{4})$\,.\\
\noindent\textbf{Summary:}
We have proved
$$
N(v,\lambda,t,h)\leq C_{12}(\mathcal{T})\|(K_{\pm,A,g_{0}}^{h}(\mathcal{T})-z)v\|
$$
when
$v\in D(K_{\pm,A,g_{0}}^{h})$ and $\Real
z\leq-2C_{5}(\mathcal{T})$
and $C_{12}(\mathcal{T})$ is  large enough. Take for $C_{12}(\mathcal{T})$\,, the
maximum value of all the intermediate $C_{k}(\mathcal{T})$ and
$\frac{2(1+\|A\|^{2})}{c_{A}}$ when $A\neq 0$\,.\\ 
\noindent\textbf{Taking $z=i\lambda$ for
  $K_{\pm,A,g_{0}}^{h}(\mathcal{T})
+
\theta(4R_{\theta}^{2}|p|_{g_{0}(q)})$:}
The operator $\theta(4R_{\theta}^{2}|p|_{g_{0}(q)})$ is a
bounded perturbation of $K_{\pm,A,g_{0}}^{h}(\mathcal{T})$ such that 
$$
\Real \langle v\,,\,
\big[K_{\pm,A,g_{0}}^{h}(\mathcal{T})+
\theta(4R_{\theta}^{2}|p|_{g_{0}(q)})
\big]
v\rangle\geq
\frac{1}{C_{2}(\mathcal{T})}\|v\|^{2}\,,
$$
when $h\leq h(\mathcal{T})$\,.
This implies
$$
\|(K_{\pm,A,g_{0}}^{h}(\mathcal{T})+
\theta(4R_{\theta}^{2}|p|_{g_{0}(q)})
-i\lambda)v\|\geq \frac{1}{C_{2}(\mathcal{T})}\|v\|
$$
Therefore there exists $C_{13}(\mathcal{T})$
such that 
\begin{multline*}
C_{13}(\mathcal{T})\|(K_{\pm,A,g_{0}}^{h}(\mathcal{T})+
\theta(4R_{\theta}^{2}|p|_{g_{0}(q)})-i\lambda)v\|\\
\geq
C_{12}(\mathcal{T})\|(K_{\pm,A,g_{0}}^{h}(\mathcal{T})
+2C_{5}(\mathcal{T})-i\lambda)v\|
\geq N(v,\lambda,t,h)\,.
\end{multline*}
We finally take
$$
C(\mathcal{T})=\max\left\{\frac{1}{h(\mathcal{T})}\,,
  C_{12}(\mathcal{T})\,, C_{13}(\mathcal{T})\right\}\,,
$$
where $C_{12}(\mathcal{T})$  depends on $A$ in the case $A\neq 0$\,.
\end{proof}
\subsection{General local metric on half-cylinders}
\label{se.basicpert}
We prove the maximal accretivity and subelliptic estimates for
Kramers-Fokker-Planck operators on $\overline{Q}=(-\infty,0]\times
Q'$\,, $Q'=\tz^{d-1}$\,, endowed with  metric $g=1\oplus m(q^{1},q')$
which fulfill \eqref{eq.metbord1}\eqref{eq.metbord2} without the
condition $\partial_{q^{1}}m\equiv 0$\,. Actually we shall work
locally near $Q'=\left\{q^{1}=0\right\}$ and introduce a parameter
$\varepsilon>0$ to be fixed within the proof.
The partial metric $m(q^{1},q')$
can be written
\begin{eqnarray*}
  &&m(q^{1},q')=m_{0}(q')+q^{1}\tilde{m}(q^{1},q')\,,\\
\text{with}&& \tilde{m}\in \mathcal{C}^{\infty}_{0}(\overline{Q};\mathcal{M}_{d-1}(\rz))\,.
\end{eqnarray*}
Since we are interested in the problem near the boundary, we can
replace $m(q^{1},q')$ by
$$
m_{\varepsilon}(q^{1},q')=m_{0}(q')+\chi_{m}(\frac{q^{1}}{\varepsilon})q^{1}\tilde{m}(q^{1},q')\,,
$$
with $\chi_{m}\in \mathcal{C}^{\infty}_{0}((-1,0])$\,, $\chi_{m}\equiv 1$
in a neighborhood of $\supp \tilde{m}$\,, so that $g_{1}=g$ while
$g_{0}=1\oplus m_{0}(q')$ satisfies $\partial_{q^{1}}m_{0}\equiv 0$\,.
The metrics $m$ and $m_{\varepsilon}$ ( resp. $g$ and
$g_{\varepsilon}$) coincide
 in $\left\{|q^{1}|\leq
  C_{g}\varepsilon\right\}$\,. 
In all our estimates, the constants 
determined by the metric $g$ and independent of $\varepsilon\in [0,1]$
 will be denoted with a subscript
$_{g}$ while the dependence with respect to $\varepsilon$ for the
metric $g_{\varepsilon}$ will be traced carefully.\\
The usefull properties of the metric $g_{\varepsilon}$ for $\varepsilon\in
[0,1]$
are summarized by
\begin{eqnarray*}
  && g_{\varepsilon}=
  \begin{pmatrix}
    1&0\\0&m_{\varepsilon}(q^{1},q')
  \end{pmatrix}\,,\quad m_{\varepsilon}(0,q')=m_{0}(q')\,,
\\
&& m_{\varepsilon}(q^{1},q')-m_{\varepsilon}(0,q')\in
\mathcal{C}^{\infty}_{0}(\overline{Q};\mathcal{M}_{d-1}(\rz))\,,
\\
&&\supp\left[m_{\varepsilon}(q^{1},q')-m_{\varepsilon}(0,q')\right]\subset\left\{|q^{1}|\leq C_{g}\varepsilon\right\}\\
&&|\partial_{q}m_{\varepsilon}(q)|\leq C_{g}\quad,\quad |m_{\varepsilon}(q^{1},q')-m_{\varepsilon}(0,q')|\leq
C_{g}|q^{1}|\leq C_{g}\varepsilon\,.
\end{eqnarray*}
The spaces $H^{s}(\overline{Q};\mathcal{H}^{s'})$ do not depend on
$g_{\varepsilon}$ but their norms
$\|u\|_{H^{s}(\overline{Q};\mathcal{H}^{s'}),g_{\varepsilon}}$ do. 
The above estimates ensures
the uniform equivalence of norms
$$
\left(
\frac{\|u\|_{H^{s}(\overline{Q};\mathcal{H}^{s'}),g_{\varepsilon}}}{\|u\|_{H^{s}(\overline{Q};\mathcal{H}^{s'}),g_{0}}}\right)^{\pm
1}\leq C_{s,s'} \quad\text{when}~(s,s')\in [-1,1]\times \rz\,,
$$
and since only the case $s\in [-1,1]$ occurs, we keep the notation
$\|u\|_{H^{s}(\overline{Q};\mathcal{H}^{s'})}$ without specifying the
metric. The choice of the norm $\|u\|_{L^{2}(Q;\mathcal{H}^{1})}$ in
the integration by parts inequality (or identity)
 will be clear from the context. With
 $g_{\varepsilon}(0,q')=m_{0}(q')$\,, the scalar product $\langle
 \gamma\,,\,
 \gamma'\rangle_{L^{2}(\partial X,|p_{1}|dq'dp;\mathfrak{f})}$ does
 not depend on $\varepsilon$\,.\\
The operator $A$ involved in the boundary condition satisfies
Hypothesis~\ref{hyp.Abdd} and \eqref{eq.commutAp} and is denoted
$A(|p|_{q})$\,.
The notations $A^{h}$\,, $\mathcal{O}_{Q,g_{\varepsilon}}^{h}$ and
$P_{\pm,Q,g_{\varepsilon}}^{h}$ are the ones of
Definition~\ref{de.resc} with $A=A^{1}$\,,
$\mathcal{O}_{Q,g_{\varepsilon}}^{1}=\mathcal{O}_{Q,g_{\varepsilon}}$
and 
$P_{\pm,Q,g_{\varepsilon}}^{1}=P_{\pm,Q,g_{\varepsilon}}$\,.
The result of this section is
\begin{proposition}
\label{pr.cylgen} Assume that $A$ fulfills Hypothesis~\ref{hyp.Abdd}
and \eqref{eq.commutAp}. There exists $\varepsilon_{g}>0$ 
such that for $\varepsilon\leq\varepsilon_{g}$ the conclusions of
Theorem~\ref{th.main0} and Theorem~\ref{th.mainA} 
are true when $\overline{Q}=(-\infty,0]\times
\tz^{d-1}$ is endowed with the metric $g_{\varepsilon}$\,.
\end{proposition}
This will be done in several steps: We first prove that for
$\varepsilon\leq \varepsilon_{g}$ small enough,
$K_{\pm,A,g_{\varepsilon}}$ is maximal accretive 
with most of the subelliptic estimates of Theorem~\ref{th.main0} and
Theorem~\ref{th.mainA}.
Then we prove the estimate of
$\|\Phi(q)\mathcal{O}_{Q,g_{\varepsilon}}u\|$\,, the equality
$K_{\pm,A,g_{\varepsilon}}^{*}=K_{\mp,A^{*},g_{\varepsilon}}$ and the
density of $\mathcal{D}(\overline{X},j)$ in $D(K_{\pm,0,g_{\varepsilon}})$\,.
\subsubsection{Maximal accretivity and first subelliptic estimates}
\label{se.mafirstsubell}
Assuming Hypothesis~\ref{hyp.Abdd} and \eqref{eq.commutAp} for $A$\,,
the domain of $K_{\pm,A,g_{\varepsilon}}$ is now defined without
prescribing a global estimate for $\gamma u$:
\begin{equation}
  \label{eq.DKAtemp}
  D(K_{\pm,A,g_{\varepsilon}})
=\left\{
 \begin{array}[c]{ll}
u\in L^{2}(Q;\mathcal{H}^{1})\,,\quad&
 P_{\pm,Q,g}u\in
  L^{2}(X,dqdp; \mathfrak{f})\,,\\
&
\gamma u \in L^{2}_{loc}(\partial X,|p_{1}|dq'dp; \mathfrak{f})\,,\\
&\gamma_{odd}u=\pm\sign(p_{1})A\gamma_{ev}u\,.    
  \end{array}
\right\}
\end{equation}
\begin{proposition}
\label{pr.dyadbord}
Assume that $A$ fulfills Hypothesis~\ref{hyp.Abdd}
and \eqref{eq.commutAp} and that
$\overline{Q}=(-\infty,0]\times\tz^{d-1}$ is
endowed with the metric $g_{\varepsilon}$\,.\\
There exists $\varepsilon_{g}>0$ such that the following properties
hold when $\varepsilon\leq \varepsilon_{g}$\,.\\
The operator $K_{\pm,A,g_{\varepsilon}}-\frac{d}{2}$\,, with
$K_{\pm,A,g_{\varepsilon}}u=P_{\pm,Q,g_{\varepsilon}}u$ and the domain
$D(K_{\pm,A,g_{\varepsilon}})$ 
given according to \eqref{eq.DKAtemp} is maximal accretive.\\
When $\supp
u\subset\left\{g^{ij}_{\varepsilon}(q)p_{i}p_{j}\leq R_{u}^{2}\right\}$
for some  $R_{u}\in (0,+\infty)$\,, $u\in
D(K_{\pm,A,g_{\varepsilon}})$ is equivalent to 
$u\in D(K_{\pm,A,g_{0}})$\,. The set of such $u$'s is dense in
$D(K_{\pm,A,g_{\varepsilon}})$ endowed with its graph norm.\\
The identity
\begin{equation}
  \label{eq.loweREKAg}
\|u\|_{L^{2}(Q;\mathcal{H}^{1})}^{2}+\Real \langle \gamma_{ev}u\,,\,
A\gamma_{ev}u\rangle_{L^{2}(\partial X,|p_{1}|dq'dp;\mathfrak{f})}
=\Real\langle
u\,,\, (\frac{d}{2}+K_{\pm,A,g_{\varepsilon}})u\rangle\,,
\end{equation}
where the left-hand side is larger thant $d\|u\|^{2}$\,,
holds for all $u\in D(K_{\pm,A,g_{\varepsilon}})$\,.\\
The adjoint of $K_{\pm,A,g_{\varepsilon}}^{*}$ equals
$K_{\mp,A^{*},g_{\varepsilon}}$\,.\\
By taking $(\nu_{1},\nu_{2},\tilde{\nu}_{3},t)\in
\left\{(\frac{1}{8},\frac{1}{8},0)\right\}\times [0,\frac{1}{18})$
when $A\neq 0$\,, and
$(\nu_{1},\nu_{2},\tilde{\nu}_{3},t)=(0,\frac{1}{4},-1,\frac{1}{3})$
when $A=0$\,, there exist a  constant $C_{t}>0$ and a constant
$C>0$ independent of $t$\,, such that
 the quantity
\begin{multline}
\label{eq.defN}
 N_{t,g_{\varepsilon}}(u,\lambda)= \langle \lambda\rangle^{\frac 1
    4}\|u\|+\langle\lambda\rangle^{\frac 1
    8}\|u\|_{L^{2}(Q;\mathcal{H}^{1})}+C_{t}^{-1}\langle\lambda\rangle^{\nu_{1}}\|u\|_{H^{t}(\overline{Q};\mathcal{H}^{0})}\\
+\langle
\lambda\rangle^{\nu_{2}}\|(1+|p|_{q})^{\tilde{\nu}_{3}}\gamma
u\|_{L^{2}(\partial X,|p_{1}|dq'dp;\mathfrak{f})}
\end{multline}
is less than $C\|(K_{\pm,A,g_{\varepsilon}}-i\lambda)u\|$\,, for all
 $u\in D(K_{\pm,A,g_{\varepsilon}})$ and all $\lambda\in\rz$\,.
\end{proposition}
When $\mathcal{T}=(b,S,T_{1},T_{2},\theta)$ fulfills the conditions
\eqref{eq.condpert},
$P_{\pm,Q,g_{0}}^{h}(\mathcal{T})$ (resp. $K_{\pm,A,g_{0}}^{h}(\mathcal{T})$)
denotes the perturbation of $P_{\pm,Q,g_{0}}^{h}$
(resp. $K_{\pm,A,g_{0}}^{h}$)
studied in Subsection~\ref{se.perturbhalf}.
\begin{lemma}
\label{le.supppborne}
Consider the ball
$$
B_{R}=\left\{g^{ij}_{\varepsilon}(q)p_{i}p_{j}\leq
  R^{2}\right\}\quad,\quad R\in (0,+\infty)\,,
$$
and the unitary change of variable
\begin{eqnarray}
  \label{eq.unitch}
  &&(U_{g_{\varepsilon}}v)(q,p)=\det(\Psi(q))^{-1/2}v(q,\Psi(q)^{-1}p)
\,,\\
\label{eq.unitch2}
&&\Psi(q)=g_{\varepsilon}(q)g_{0}(q)^{-1}\,.
\end{eqnarray}

 There exists
$\mathcal{T}=(b,S,T_{1},T_{2},\theta)$ fulfilling the conditions
\eqref{eq.condpert}, with $R_{\theta}=C_{g}R$ and
$\varepsilon_{T}\leq C_{g}\varepsilon$\,,
such that 
$$
K_{\pm,A,g_{\varepsilon}}u=U_{g_{\varepsilon}}K_{\pm,A,g_{0}}^{1}(\mathcal{T})U^{*}_{g_{\varepsilon}}u\,,
$$
for all $u\in D(K_{\pm,A,g_{\varepsilon}})$ with $\supp u\subset B_{R}$\,.\\
There exists a constant $C_{R,g_{\varepsilon}}$ 
independent of $(\lambda,t,u)$ such that the quantity
\eqref{eq.defN} is estimated by
$$
N_{t,g_{\varepsilon}}(u,\lambda)\leq C_{R,g_{\varepsilon}}
\left[\|(K_{\pm,A,g_{\varepsilon}}-i\lambda)u\|+\|u\|\right]\,.
$$
for all $\lambda \in \rz$ and all $u\in
D(K_{\pm,A,g_{\varepsilon}})$ such that $\supp u\subset
B_{R}$\,. Moreover the integration by parts identity
\eqref{eq.loweREKAg} holds for such $u$'s.
\end{lemma}
\begin{proof}
The unitary change of variable 
\begin{eqnarray*}
  u(q,p)&=&(U_{g_{\varepsilon}}v)(q,p)=\det(\Psi(q))^{-1/2}v(q,\Psi(q)^{-1}p)\\
&=&\det(\mu(q))^{-1/2}v(q,p_{1},\mu(q)^{-1}p')\,,
\\
 \Psi(q)&=&g_{\varepsilon}(q)g_{0}(q)^{-1}=g_{\varepsilon}(q^{1},q')g(0,q')^{-1}
=\begin{pmatrix}
  1&0\\0&\mu(q)
\end{pmatrix}\\
\mu(q)&=&
m_{\varepsilon}(q)m_{0}(q)^{-1}=m_{\varepsilon}(q^{1},q')m_{0}(q')^{-1}\,,
\end{eqnarray*}
does not change the functional spaces
$L^{2}(Q;\mathcal{H}^{s'})$\,, $s'\in\rz$\,,
$L^{2}(\partial X,|p_{1}|dq'dp)$\,, nor the traces of $\gamma_{ev}u$
and $\gamma_{odd}u$ because $\mu(q)=\Id_{\rz^{d-1}}$ when
$q^{1}=0$\,.  
The norms of $U_{g_{\varepsilon}}$ and $U_{g_{\varepsilon}}^{-1}$ as a bounded operator in
$L^{2}(Q;\mathcal{H}^{s'})$ are uniformly bounded by
$C_{s',g}\geq 1$ when $\varepsilon\in (0,\varepsilon_{g}]$\,.
The support condition $\supp u\subset\left\{g^{ij}_{\varepsilon}(q)p_{i}p_{j}\leq R^{2}\right\}$
tranformed into $\supp v\subset \left\{g^{ij}_{0}(q)p_{i}p_{j}\leq R_{\theta}^{2}\right\}$
for $R_{\theta}\leq C_{g}R$\,. A direct calculation gives
\begin{eqnarray*}
U^{-1}_{g_{\varepsilon}}P_{\pm,Q,g_{\varepsilon}}U_{g_{\varepsilon}}&=&P_{\pm,Q,g_{0}}\pm
p_{i}b^{i}(q)
\pm S^{ij}_{k}(q)p_{i}p_{j}\partial_{p_{k}}\\
&&\hspace{4cm}
+\frac{-(\partial_{p})^{T}T_{1}(q)(\partial_{p})+p^{T}T_{2}(q)p}{2}\,,\\
\text{with}&&
T_{1}(q)=g_{0}(q)[g_{\varepsilon}(q)^{-1}-g_{0}(q)^{-1}]g_{0}(q)\,,
\\
&&T_{2}(q)=g_{0}(q)^{-1}[g_{\varepsilon}(q)-g_{0}(q)]g_{0}(q)^{-1}\,,\\
&&\|b\|_{L^{\infty}}+\|S\|_{L^{\infty}}\leq C_{g}\,.
\end{eqnarray*}
The support condition on $u$ (and $v$) allows to replace $p_{i}b^{i}(q)$
by $\theta(|p|_{g_{0}(q)})p_{i}b^{i}(q)$ and
$S^{ij}(q)p_{i}p_{j}\partial_{p_{k}}$ by
$S^{ij}(q)\theta(|p|_{g_{0}(q)})p_{i}p_{j}\partial_{p_{k}}$\,.\\
We also deduce
$$
|T_{1}(q)|+|T_{2}(q)|\leq C_{g}|q^{1}|\leq C_{g}\varepsilon\,.
$$
Since $\Psi(q)=\Id_{\rz^{d}}$ along $\left\{q^{1}=0\right\}$\,, the
traces are not changed and the boundary condition for $u\in
D(K_{\pm,A,g_{\varepsilon}})$ becomes $v\in D(K_{\pm,A,g_{0}})=D(K_{\pm,A,g_{0}}^{1}(\mathcal{T}))$\,.
We are exactly in the case of Proposition~\ref{pr.perturbflat} with $h=1$\,, as soon
as $\varepsilon\leq \varepsilon_{g}$ for $\varepsilon_{g}$ small
enough. With
$$
U_{g}^{-1}K_{\pm,A,g}U_{g_{\varepsilon}}v=K_{\pm,A,g_{0}}^{1}(\mathcal{T})v\,,
$$
the estimate of $N_{t,g_{\varepsilon}}(u,\lambda)\leq C_{R,g_{\varepsilon}}N(v,\lambda,t,h=1)$
 is given by Proposition~\ref{pr.perturbflat}.\\
For the identity \eqref{eq.loweREKAg}, assume $\supp u\subset B_{R}$
and compute
\begin{eqnarray*}
&& \Real\langle
u\,,\, [\frac{d}{2}+K_{\pm,A,g}]u\rangle -\|u\|_{L^{2}(Q;\mathcal{H}^{1})}^{2}-
\Real \langle \gamma_{ev}u\,,\, A\gamma_{ev}u\rangle_{L^{2}(\partial X,|p_{1}|dq'dp;\mathfrak{f})}
\\
&&=\Real \langle v\,,\, [\frac{d}{2}+K_{\pm,A,g_{0}}^{1}(\mathcal{T})]v\rangle
-\|u\|_{L^{2}(\rz^{d}_{-};\mathcal{H}^{1})}^{2}
\\
&&
\hspace{6cm}-
\Real \langle \gamma_{ev}v\,,\,
A\gamma_{ev}v\rangle_{L^{2}(\partial X,|p_{1}|dq'dp;\mathfrak{f})}\\
&&=\pm \langle v\,,\, p_{i}b^{i}(q)v\rangle
\mp\frac{1}{2}\langle v\,,\, (S^{kj}_{k}(q)p_{j}+S^{ik}_{k}(q)p_{i})v\rangle\,.
\end{eqnarray*}
The same computation done for $u\in \mathcal{C}^{\infty}_{0}(X;\mathfrak{f})$ 
leads to a vanishing right-hand side. Since this right-hand side is
continuous on $L^{2}(Q;\mathcal{H}^{1})$\,, it always vanishes when
$u\in D(K_{\pm,A,g_{\varepsilon}})\subset L^{2}(Q;\mathcal{H}^{1})$
and $\supp u\subset B_{R}$\,.
\end{proof}
\begin{lemma}
\label{le.supppdyad}
Consider the shell
$$
S_{R,\ell}=
\left\{R^{-2}2^{2\ell}\leq g_{\varepsilon}^{ij}(q)p_{i}p_{j}\leq
  R^{2}2^{2\ell}\right\}\,, \quad R\in (0,+\infty)\,,\quad \ell\in\nz\,,
$$
and the unitary change of variable
$U_{g_{\varepsilon},\ell}=U_{g_{\varepsilon}}V_{\ell}$ where
$U_{g_{\varepsilon}}$ is defined by
\eqref{eq.unitch}\eqref{eq.unitch2} and $V_{\ell}w=2^{-\ell
  d/2}w(q,2^{-\ell}p)$\,.\\
There exists
$\mathcal{T}=(b,S,T_{1},T_{2},\theta)$ fulfilling the conditions
\eqref{eq.condpert}, with $R_{\theta}=C_{g}R$ and 
$\varepsilon_{T}\leq C_{g}\varepsilon$\,, such that 
$$
K_{\pm,A,g_{\varepsilon}}u=\frac{1}{h}U_{g_{\varepsilon},\ell}
\left[\theta(4R_{\theta}^{2}|p|_{g_{0}(q)})+K_{\pm,A,g_{0}}^{h}(\mathcal{T})\right]
U_{g_{\varepsilon},\ell}^{*}u\quad,\quad h=2^{-2\ell}\,.
$$
for all $u\in D(K_{\pm,A,g_{\varepsilon}})$ with $\supp u\subset S_{R,\ell}$\,.\\
There exists two constants $C_{R,g_{\varepsilon}}$ and $\ell_{R,g_{\varepsilon}}$\,,  independent of
$(\lambda,t,u,\ell)$\,, such that the quantity
\eqref{eq.defN} is estimated by
$$
N_{t,g_{\varepsilon}}(u,\lambda)\leq C_{R,g_{\varepsilon}}\|(K_{\pm,A,g_{\varepsilon}}-i\lambda)u\|\,,
$$
for all $\lambda\in\rz$
for all $u\in D(K_{\pm,A,g_{\varepsilon}})$ such that $\supp u\subset
S_{R,\ell}$ as soon as $\ell\geq
\ell_{R,g_{\varepsilon}}$\,.\\
Moreover the operator
$\theta(4R_{\theta}^{2}|p|_{g_{0}(q)})+K_{\pm,A,g_{0}}^{h}(\mathcal{T})$\,,
with $h=2^{-2\ell}$\,, is maximal accretive when
$\ell\geq \ell_{R,g_{\varepsilon}}$\,.
\end{lemma}
\begin{proof}
  We start with the same unitary change of variable $U_{g_{\varepsilon}}$ as in
  Lemma~\ref{le.supppborne}.
With $u=U_{g_{\varepsilon}}v$ one gets
$$
U_{g_{\varepsilon}}^{-1}K_{\pm,A,g_{\varepsilon}}U_{g_{\varepsilon}}v= 
K_{\pm,A,g_{0}}^{1}(\mathcal{T}_{\ell})v
$$
where $\mathcal{T}_{\ell}$ equals
$(b,S,T_{1},T_{2},\theta(2^{-\ell}.))$  
while $\mathcal{T}=(b,S,T_{1},T_{2},\theta)$ satisfies the same
conditions as in Lemma~\ref{le.supppborne}. Additionally one can
choose
$\varepsilon_{g}>0$ small enough  and
$R_{\theta}\leq C_{g}R$  such that $\supp
u\subset S_{R,\ell}$
implies
$$
\supp v \subset \left\{R_{\theta}^{-2}2^{2\ell}\leq
 g_{0}^{ij}(q)p_{i}p_{j}\leq R_{\theta}^{2}2^{2\ell}\right\}\,.
$$
The function $v$ is
supported in $\left\{1\leq R_{\theta}2^{-\ell}|p|_{g_{0}(q)}\leq
  R_{\theta}^{2}\right\}$ while
$\theta(4R_{\theta}^{2}2^{-\ell}|p|_{g_{0}(q)})$ is supported in 
$\left\{R_{\theta}2^{\ell}|p|_{g_{0}(q)}\leq \frac{1}{2}\right\}$\,.
We obtain
$$
U_{g_{\varepsilon}}^{-1}K_{\pm,A,g_{\varepsilon}}U_{g_{\varepsilon}}v
= 
\left[\theta(4R_{\theta}^{2}2^{-\ell}|p|_{g_{0}(q)})
+K_{\pm,A,g_{0}}^{1}(\mathcal{T}_{\ell})\right]v
$$
The unitary change of variable $v=V_{\ell}w=2^{-\ell d/2}w(q,2^{-\ell}p)$
gives
$$
V_{\ell}^{-1}U_{g_{\varepsilon}}^{-1}K_{\pm,A,g_{\varepsilon}}U_{g_{\varepsilon}}V_{\ell}w=\frac{1}{h}
\left[\theta(4R_{\theta}^{2}|p|_{g_{0}(q)})+K_{\pm,A,g_{0}}^{h}(\mathcal{T})\right]w
$$
with $h=2^{-2\ell}$\,.  After noticing that
$N_{t,g_{\varepsilon}}(u,\lambda)\leq
\frac{C_{g_{\varepsilon},R}N(v,\lambda,t, h)}{h}$ with
$h=2^{-2\ell}$\,, it 
  suffices to take $U_{g_{\varepsilon},\ell}=U_{g_{\varepsilon}}V_{\ell}$
and to apply the subelliptic estimates of
Proposition~\ref{pr.perturbflat} for
$K_{\pm,A,g_{0}}^{h}(\mathcal{T})-i\lambda h$ which are valid for $h\leq
\frac{1}{C(\mathcal{T})}$ small enough.\\
The maximal accretivity of
$\left[\theta(4R_{\theta}^{2}|p|_{g_{0}(q)})+K_{\pm,A,g_{0}}^{h}(\mathcal{T})\right]$
was also checked in Proposition~\ref{pr.perturbflat} for $h\leq
\frac{1}{C(\mathcal{T})}$\,, which means simply $\ell\geq \ell_{R,g_{\varepsilon}}$\,.
\end{proof}
\begin{proof}[Proof of Proposition~\ref{pr.dyadbord}:]~\\
We fix $\varepsilon\leq \varepsilon_{g}$ with $\varepsilon_{g}$ small
enough so that Lemma~\ref{le.supppborne} and Lemma~\ref{le.supppdyad}
apply\,.
The cylinder $\overline{Q}=(-\infty,Ø]\times\tz^{d}$ is then endowed
with the metric $g_{\varepsilon}$\,. The quantities
$|p|_{g_{\varepsilon}(q)}$ and $|p|_{g_{0}(q)}$ satisfy uniformly with
respect to $(q,\varepsilon)\in Q'\times [0,\varepsilon_{g}]$\,,
$\left(\frac{|p|_{g_{\varepsilon}(q)}}{|p|_{g_{0}(q)}}\right)^{\pm
  1}\leq (1+C_{g}\varepsilon)$\,.\\
\noindent\textbf{Dyadic partition of unity:}
The dyadic partition of unity $\chi=(\chi_{\ell})_{\ell\in \nz}$
with $\nz$ is given like in Paragraph~\ref{se.momentpart}
by  $\chi_{0}(q,p)=\tilde{\chi}_{0}(|p|_{g_{\varepsilon}(q)})$ and
$\chi_{\ell}(q,p)=\tilde{\chi}_{1}(2^{-\ell}|p|_{g_{\varepsilon}(q)})$
with $\sum_{\ell\in \nz}\chi_{\ell}^{2}\equiv 1$\,.
The
cut-off functions $\tilde{\chi}_{0}$ and $\tilde{\chi}_{1}$ are
 assumed to be supported respectively in $\left\{|p|_{g_{\varepsilon}(q)}\leq
   R_{0}\right\}$ and
$\left\{R_{1}^{-1}\leq |p|_{g_{\varepsilon}(q)}\leq R_{1}\right\}$ for some $R_{1}\geq 1$\,.\\
By applying Lemma~\ref{le.supppdyad} with $R=2R_{1}$\,, there exists
$\ell_{\chi,g_{\varepsilon}}=\ell_{2R_{1},g_{\varepsilon}}\geq 1$ and
$\mathcal{T}_{1}=(b,S,T_{1},T_{2},\theta)$ such that for $\ell\geq
\ell_{\chi,g_{\varepsilon}}$ the operator
$$
K_{\pm,A,g_{\varepsilon},\ell}=\frac{1}{h}U_{g_{\varepsilon},\ell}
\left[\theta(4R_{\theta}^{2}|p|_{g_{0}(q)})+K_{\pm,A,g_{0}}^{h}(\mathcal{T}_{1})\right]U_{g,\ell}^{*}\quad,\quad h=2^{-2\ell}\,,
$$
is maximal accretive while the condition \eqref{eq.commutAp} ensures
\begin{eqnarray*}
 \forall \ell\leq \ell_{\chi,g_{\varepsilon}}\,, &&
  \left(u\in
    D(K_{\pm,A,g_{\varepsilon}})\right)\Rightarrow\left(\chi_{\ell}u\in
  D(K_{\pm,A,g_{\varepsilon}})\subset
  D(K_{\pm,A,g_{\varepsilon},\ell})\right)\,,\\
&& \forall u\in D(K_{\pm,A,g_{\varepsilon}})\,,\quad 
K_{\pm,A,g_{\varepsilon}}\chi_{\ell}u=K_{\pm,A,g_{\varepsilon},\ell}\chi_{\ell}u\,.\\
\end{eqnarray*}
Moreover, there exists a constant $C_{\chi,g_{\varepsilon}}^{1}$
independent of $\ell\geq \ell_{\chi,g_{\varepsilon}}$ such that
$$
N_{t,g_{\varepsilon}}(\chi_{\ell}u\,,\lambda)\leq
C_{\chi,g_{\varepsilon}}^{1}\|(K_{\pm,A,g_{\varepsilon},\ell}-i\lambda)\chi_{\ell}u\|
$$
holds for all $(u,\lambda)\in D(K_{\pm,A,g_{\varepsilon}})\times\rz
$ and all $\ell\geq \ell_{\chi,g_{\varepsilon}}$\,.\\
We now replace $\nz$ by $\mathcal{L}=\left\{-1\right\}\cup \left\{\ell\geq
  \ell_{\chi,g_{\varepsilon}}\right\}$ and set 
$\chi_{-1}^{2}=\sum_{\ell=0}^{\ell_{2R_{1},g_{\varepsilon}}-1}\chi_{\ell}^{2}$\,.
All the results of Paragraph~\ref{se.momentpart}, namely
Lemma~\ref{le.partp}, Lemma~\ref{le.partpq} and
Proposition~\ref{pr.partp}, still hold with 
$\mathcal{L}=\nz$ replaced by $\mathcal{L}=\left\{-1\right\}\cup \left\{\ell\geq
  \ell_{\chi,g_{\varepsilon}}\right\}$\,.\\
Again when $u\in D(K_{\pm,A,g_{\varepsilon}})$ the
condition~\eqref{eq.commutAp} ensures $\chi_{-1}u\in
D(K_{\pm,A,g_{\varepsilon}})$ while $\supp \chi_{-1}u\subset
  \left\{|p|_{g_{\varepsilon}(q)}\leq
    R_{1}2^{\ell_{\chi,g_{\varepsilon}}-1}\right\}$\,.
By  Lemma~\ref{le.supppborne} applied with
$R=R_{1}2^{\ell_{\chi,g_{\varepsilon}}}$\,, there exists
$\mathcal{T}_{-1}=(b,S,T_{1},T_{2},\theta_{-1})$  and a constant
$C_{\chi,g_{\varepsilon}}^{2}>0$ such that
$$
C_{\chi,g_{\varepsilon}}^{2}+K_{\pm,A,g_{\varepsilon},-1}=
U_{g_{\varepsilon}}(C_{\chi,g_{\varepsilon}}^{2}+K_{\pm,A,g_{0}}(\mathcal{T}_{-1}))U_{g_{\varepsilon}}^{*}
$$
is
maximal accretive, and the relations
\begin{eqnarray*}
&& K_{\pm,A,g_{\varepsilon}}\chi_{-1}u=
K_{\pm,A,g_{\varepsilon},-1}
\chi_{-1}u\,,\\
&& N_{t,g_{\varepsilon}}(\chi_{-1}u,\lambda)
\leq
C_{\chi,g_{\varepsilon}}^{2}\|(K_{\pm,A,g_{\varepsilon},-1}+C_{\chi,g_{\varepsilon}}^{2}-i\lambda)\chi_{-1}u\|\,,
\end{eqnarray*}
hold for all $(u,\lambda,t)\in D(K_{\pm,A,g_{\varepsilon}})\times\rz\times
[0,\frac{1}{18})$ when $A\neq 0$ (take $t=\frac{1}{3}$ for the case $A=0$)\,.\\
\noindent\textbf{Domain and subelliptic estimate:} Take $C\geq
C_{\chi,g_{\varepsilon}}^{2}$ and $\lambda\in\rz$\,.
 By definition $u\in
D(K_{\pm,A,g_{\varepsilon}})$ means 
\begin{eqnarray*}
  &&
  \|u\|_{L^{2}(Q;\mathcal{H}^{1})}^{2}+\|(P_{\pm,Q,g_{\varepsilon}}+C-i\lambda)u\|^{2}+\|\gamma_{odd}u\|_{L^{2}(\partial X;|p_{1}|dq'dp;\mathfrak{f})}^{2}<
  +\infty\,,\\
&& \gamma_{ev}u\in L^{2}_{loc}(\partial X,|p_{1}|dq'dp;\mathfrak{f})
\quad \gamma_{odd}u=\sign(p_{1})A(|p|_{q}) \gamma_{odd}u~\text{for~a.e.~}|p|_q\,.
\end{eqnarray*}
With Proposition~\ref{pr.partp} applied with 
$\sum_{\ell\in\mathcal{L}}\chi_{\ell}^{2}\equiv 1$\,,
$\mathcal{L}=\left\{-1\right\}\cup \left\{\ell\geq
  \ell_{\chi,g_{\varepsilon}}\right\}$\,, and
$P_{\pm,Q,g}\chi_{\ell}u= K_{\pm,Q,g,\ell}u$\,, it means simply
\begin{align*}
  &
\forall \ell\in \mathcal{L}\,,\quad \chi_{\ell}u\in
  D(K_{\pm,A,g_{\varepsilon},\ell})\,,\\
&\sum_{\ell\in \mathcal{L}}
\|\chi_{\ell}u\|_{L^{2}(Q;\mathcal{H}^{1})}^{2}+\|(K_{\pm,Q,g_{\varepsilon}}+C-i\lambda)\chi_{\ell}u\|^{2}\\
&\hspace{5cm}+\|\gamma_{odd}\chi_{\ell
u}\|_{L^{2}(\partial X;|p_{1}|dq'dp;\mathfrak{f})}^{2}<
  +\infty\,.
\end{align*}
Summing the lower bound $\sum_{\ell\in \mathcal{L}}N_{t,g_{\varepsilon}}(\chi_{\ell}u,
\lambda)^{2}$ of
$\|(K_{\pm,A,g_{\varepsilon},\ell}+C-i\lambda)\chi_{\ell}u\|^{2}$
after setting, $R_{\ell}=2R_{1}2^{\ell}$ for $\ell\geq
\ell_{\chi,g_{\varepsilon}}$ and 
$R_{\ell}=R_{1}2^{\ell_{\chi,g_{\varepsilon}}}$ for $\ell=-1$\,,
$$
\sum_{\ell\in\mathcal{L}}
N_{t,g_{\varepsilon}}(\chi_{\ell}u,
\lambda)^{2}
\leq C_{\chi,g_{\varepsilon}}^{3}\left[\|(K_{\pm,Q,g_{\varepsilon}}+C-i\lambda)u\|^{2}+
\|u\|_{L^{2}(Q;\mathcal{H}^{1})}^{2}\right]\,.
$$
It contains $(1+|p|)^{\tilde{\nu}_{3}}\gamma u\in
L^{2}(\partial X;|p_{1}|dq'dp;\mathfrak{f})$  with
$\tilde{\nu}_{3}=-1$ (resp. $\tilde{\nu}_{3}=0$)
when $A=0$ (resp. $A\neq 0$). The summation of
$\|\chi_{\ell}u\|_{L^{2}(Q;\mathcal{H}^{1})}^{2}$ 
and $\|\chi_{\ell}u\|_{H^{t}(\overline{Q};\mathcal{H}^{0})}^{2}$ are
estimated from below by Lemma~\ref{le.partp} and
Lemma~\ref{le.partpq} so that
$$
\langle\lambda\rangle^{\frac{1}{4}}\|u\|_{L^{2}(Q;\mathcal{H}^{1})}^{2}\leq
N_{t,g_{\varepsilon}}(u,\lambda)^{2}
\leq C_{\chi,g_{\varepsilon}}^{3}\left[\|(K_{\pm,Q,g_{\varepsilon}}+C-i\lambda)u\|^{2}+
\|u\|_{L^{2}(Q;\mathcal{H}^{1})}^{2}\right]\,.
$$
Taking $\lambda_{0}$
large enough implies
$$
\forall u\in D(K_{\pm,A,g_{\varepsilon}})\,,\quad
N_{t,g_{\varepsilon}}(u, \lambda)\leq C_{\chi,g_{\varepsilon}}^{4}
\|(K_{\pm,A,g_{\varepsilon}}+C-i\lambda)u\|\,,
$$
for all
$u\in D(K_{\pm,A,g_{\varepsilon}})$ and all $\lambda\in\rz$ such that $|\lambda|\geq
\lambda_{0}$\,.\\
\noindent\textbf{Approximation in $D(K_{\pm,A,g_{\varepsilon}})$:} The
squared graph norm
$\|u\|^{2}+\|K_{\pm,A,g_{\varepsilon}}u\|^{2}$
is equivalent to the series
$$
\sum_{\ell\in \mathcal{L}}\|(K_{\pm,Q,g_{\varepsilon},\ell}+C-i\lambda)\chi_{\ell}u\|^{2}
$$
which is the limit of finite sums. Thus the set of $u\in
D(K_{\pm,A,g_{\varepsilon}})$ such that 
$$
\supp u\subset\left\{g_{\varepsilon}^{ij}(q)p_{i}p_{j}\leq R_{u}^{2}\right\}
$$
for some $R_{u}\in (0,+\infty)$\,, is dense in
$D(K_{\pm,A,g_{\varepsilon}})$ endowed with the graph norm.\\
\noindent\textbf{Accretivity of
  $K_{\pm,A,g_{\varepsilon}}-\frac{d}{2}$:} When $u\in
D(K_{\pm,A,g_{\varepsilon}})$ satisfies $\supp u\subset
\left\{g_{\varepsilon}^{ij}(q)p_{i}p_{j}\leq R_{u}^{2}\right\}$\,,
Lemma~\ref{le.supppborne} with $R=R_{u}$ implies that $u$ fulfills 
the inequality
\eqref{eq.loweREKAg}\,. Since all the terms of \eqref{eq.loweREKAg}
are continuous on $D(K_{\pm,A,g_{\varepsilon}})$ endowed with the
graph norm, it is extended to all $u\in
D(K_{\pm,A,g_{\varepsilon}})$\,. In particular the subelliptic
estimate
can be written now with $C=0$ and any $\lambda\in\rz$\,.\\
\noindent\textbf{Maximal accretivity:} For $C\geq
C_{\chi}^{2}$ and $\lambda\in \rz$\,, all the operators
$K_{\pm,A,g_{\varepsilon},\ell}+C-i\lambda$ are invertible. 
For $f\in L^{2}(X;\mathfrak{f})$ set
$v_{\ell}=\chi_{\ell}(K_{\pm,A,g_{\varepsilon},\ell}+C-i\lambda)^{-1}\chi_{\ell}f$\,. 
The subelliptic estimates of Lemma~\ref{le.supppborne} and
Lemma~\ref{le.supppdyad} imply
$$
\|(K_{\pm,A,g_{\varepsilon},\ell}+C-i\lambda)^{-1}\chi_{\ell}f\|_{L^{2}(Q;\mathcal{H}^{1})}
\leq
  C_{\chi,g_{\varepsilon}}\langle \lambda
  \rangle^{-\frac{1}{8}}\|\chi_{\ell}f\|\,.
$$
With $\sup_{\ell'\in
  \mathcal{L}}\sharp\left\{\ell\in
  \mathcal{L}\,,\chi_{\ell}\chi_{\ell'}\neq 0\right\}\leq N_{\chi}$\,,
we deduce that $v=\sum_{\ell\in \mathcal{L}}v_{\ell}$ satisfies
\begin{eqnarray*}
\|v\|_{L^{2}(Q;\mathcal{H}^{1})}^{2}
&\leq& C_{\chi}\sum_{\ell\in
  \mathcal{L}}\|\chi_{\ell}(K_{\pm,A,g_{\varepsilon},\ell}+C-i\lambda)^{-1}\chi_{\ell}f\|_{L^{2}(Q;\mathcal{H}^{1})}^{2}\\
&\leq& C_{\chi,g_{\varepsilon}}'\sum_{\ell\in \mathcal{L}}\|\chi_{\ell}f\|^{2}=C_{\chi,g_{\varepsilon}}\|f\|^{2}\,.
\\
\gamma v=\sum_{\ell\in \mathcal{L}}\gamma v_{\ell}&\in & L^{2}_{loc}(\partial X,
|p_{1}|dq'dp;\mathfrak{f})\,,\\
\gamma_{odd}v&=&\sum_{\ell}\chi_{\ell}\gamma_{odd}(K_{\pm,A,g_{\varepsilon},\ell}+C-i\lambda)^{-1}\chi_{\ell}f
\\
&=&
\pm\sum_{\ell\in\mathcal{L}}\chi_{\ell}\sign(p_{1})A(|p|_{q})\gamma_{odd}(K_{\pm,A,g_{\varepsilon},\ell}+C-i\lambda)^{-1}\chi_{\ell}f\\
&=&\pm\sign(p_{1})A(|p|_{q})\gamma_{ev}v\,,
\end{eqnarray*}
by using the commutation $[\chi_{\ell},A]=0$ for the last line.\\
For any $\ell\in \mathcal{L}$\,, the function $v_{\ell}\in L^{2}(Q;\mathcal{H}^{1})$ satisfies
\begin{eqnarray*}
(P_{\pm,A,g_{\varepsilon}}+C-i\lambda)v_{\ell}&=&
(K_{\pm,A,g_{\varepsilon},\ell}+C-i\lambda)v_{\ell}
\\
&=&
\chi_{\ell}^{2}f-\frac{1}{2}\left[\Delta_{p},\chi_{\ell}\right](K_{\pm,A,g_{\varepsilon},\ell}+C-i\lambda)^{-1}\chi_{\ell}f\\
&=&
\chi_{\ell}^{2}f+B_{\ell}(K_{\pm,A,g_{\varepsilon},\ell}+C-i\lambda)^{-1}\chi_{\ell}f\,,
\end{eqnarray*}
after setting
$B_{\ell}=-\frac{1}{2}\left[\Delta_{p},\chi_{\ell}\right]=-\left[(\nabla_{p}\chi_{\ell}).\nabla_{p}+\frac{(\Delta_{p}\chi_{\ell})}{2}\right]$\,.
With the notations of Paragraph~\ref{se.momentpart}, the
operators $B_{\ell}$ are  bounded in $\DiffP^{1,-1}$\,, uniformly
w.r.t $\ell\in \mathcal{L}$\,,  with
$p-\supp B_{\ell}\subset \left\{R_{\chi,g_{\varepsilon}}^{-1}2^{\ell}\leq
  |p|_{g_{\varepsilon}(q)}\leq R_{\chi,g_{\varepsilon}}2^{\ell}\right\}$\,.\\
We get like in Proposition~\ref{pr.partp}
$$
\|\sum_{\ell\in
  \mathcal{L}}B_{\ell}(K_{\pm,A,g_{\varepsilon},\ell}+C-i\lambda)\chi_{\ell}f\|^{2}
\leq
C_{\chi,g_{\varepsilon}}^{7}\langle \lambda\rangle^{-\frac{1}{4}}
\|f\|^{2}\,.
$$
For any $f\in L^{2}(X;\mathfrak{f})$ we have found $v=\sum_{\ell\in
  \mathcal{L}}v_{\ell}\in D(K_{\pm,A,g_{\varepsilon}})$ such that 
$$
(K_{\pm,A,g_{\varepsilon}}+C-i\lambda)v=(\Id+B)f
$$
with $\|B\|_{\mathcal{L}(L^{2}(X;\mathfrak{f}))}\leq
\sqrt{C_{\chi,g_{\varepsilon}}^{7}}\langle\lambda\rangle^{-\frac{1}{8}}$\,. For
$\lambda$ large enough $(\Id+B)$ is invertible and
$(K_{\pm,A,g_{\varepsilon}}+C-i\lambda)$ is invertible.
\end{proof}

\subsubsection{Estimate of $\|\Phi(q^{1})\mathcal{O}_{Q,g_{\varepsilon}}u\|$ }
\label{se.estimO}
In the definition \eqref{eq.defN} estimated in
Proposition~\ref{pr.dyadbord}, the quantity
$\|\Phi(q^{1})\mathcal{O}_{Q,g_{\varepsilon}}u\|$ is missing. Actually
our approach which handles non maximal subelliptic estimates with
exponent divided by $2$ in the end does not provide directly an upper
bound for it. It can be obtained in a second step by adapting the
prof of Proposition~\ref{pr.regL1}.
\begin{proposition}
\label{pr.estimO} Let $K_{\pm,A,g_{\varepsilon}}$ be the maximal
accretive realization of $P_{\pm,Q,g_{\varepsilon}}$ defined in
Proposition~\ref{pr.dyadbord}
with $\overline{Q}=(-\infty,0]\times\tz^{d-1}$\,, $\varepsilon\leq
\varepsilon_{g}$ small enough and $A$ satisfying
Hypothesis~\ref{hyp.Abdd} and \eqref{eq.commutAp}.
For any $\Phi\in \mathcal{C}^{\infty}_{b}((-\infty,0])$ such that
$\Phi(0)=0$\,. There exists a constant $C_{g_{\varepsilon}}$
independent of $\Phi$ and a constant $C_{g_{\varepsilon},\Phi}$ such that 
$$
\|\Phi(q^{1})\mathcal{O}_{Q,g_{\varepsilon}}u\|\leq
C_{g_{\varepsilon}}\|\Phi\|_{L^{\infty}}\|(K_{\pm,A,g_{\varepsilon}}-i\lambda)u\|+
C_{g_{\varepsilon},\Phi} \|u\|\,,
$$
holds for all $u\in D(K_{\pm,A,g_{\varepsilon}})$ and all $\lambda\in\rz$\,.
\end{proposition}
\begin{proof}
  We embed the manifold $\overline{Q}=(-\infty,0]\times Q'$\,,
  $Q'=\tz^{d-1}$\,,
 into $\tilde{Q}=\rz\times Q'$\,. The metric $\tilde{g}_{\varepsilon}$ on
$\tilde{Q}$  is assumed to be
$\mathcal{C}^{\infty}$ with
$\tilde{g}_{\varepsilon}\big|_{\overline{Q}}=g_{\varepsilon}$ and
$\tilde{g}_{\varepsilon}-g_{0}\in
\mathcal{C}^{\infty}_{0}(\tilde{Q})$\,. Due to the curvature of
  $\partial Q$ the metric $\tilde{g}_{\varepsilon}$ is not given by a
  simple symmetry argument but a $\mathcal{C}^{\infty}$ extension is
  always possible, first locally (see \cite{ChPi}) and then globally
  with a partition of unity.\\
When $u\in D(K_{\pm,A,g_{\varepsilon}})$ and $\Phi(0)=0$\,, the function
$\Phi(q^{1})u$ belongs to $
D(K_{\pm,\tilde{Q},\tilde{g}_{\varepsilon}})$ and
Proposition~\ref{pr.cylacc} about whole cylinders provides
$$
\|\mathcal{O}_{\tilde{Q},\tilde{g}_{\varepsilon}}\Phi(q^{1})u\|\leq
C_{g_{\varepsilon}}\|
(K_{\pm,\tilde{Q},\tilde{g}_{\varepsilon}}-i\lambda)\Phi(q^{1})u\|\,,
$$
equivalently written 
$$
\|\Phi(q^{1})\mathcal{O}_{Q,g_{\varepsilon}}u\|\leq
C_{g_{\varepsilon}}\|
(K_{\pm,Q,g_{\varepsilon}}-i\lambda)\Phi(q^{1})u\|\,.
$$
We compute
$$
(K_{\pm,Q,g_{\varepsilon}}-i\lambda)(\Phi(q^{1}))u)=
\Phi(q^{1})(K_{Q,A,g_{\varepsilon}}-i\lambda)u+p_{1}\Phi'(q^{1})u\,,
$$
in order to get
$$
\|\Phi(q^{1})\mathcal{O}_{Q,g_{\varepsilon}}u\|\leq
C_{g_{\varepsilon}}^{1}\|
\Phi(q^{1})(K_{\pm,Q,g_{\varepsilon}}-i\lambda)u\|
+ C_{\Phi,g_{\varepsilon}}^{1}\|u\|_{L^{2}(Q;\mathcal{H}^{1})}\,.
$$
But the integration by part identity \eqref{eq.loweREKAg} and the
accretivity of $K_{\pm,A,g_{\varepsilon}}-\frac{d}{2}$
\begin{align*}
\|u\|_{L^{2}(Q;\mathcal{H}^{1})}^{2}&\leq \|u\|\|(K_{\pm,A,g_{\varepsilon}}+\frac{d}{2}-i\lambda)u\|\\
&2\|u\|\|(K_{\pm,A,g_{\varepsilon}}-i\lambda)u\|
\leq \left[\delta\|(K_{\pm,A,g_{\varepsilon}}-i\lambda)u\|
+ \delta^{-1}\|u\|\right]^{2}\,,
\end{align*}
for all $\delta>0$\,.
Choosing
$\delta=\frac{C_{g_{\varepsilon}}^{1}\|\Phi\|_{L^{\infty}}}{C_{\Phi,g_{\varepsilon}}^{1}}$
 leads
to 
$$
\|\Phi(q^{1})\mathcal{O}_{Q,g_{\varepsilon}}u\|\leq
2C_{g_{\varepsilon}}^{1}\|\Phi\|_{L^{\infty}}\|(K_{\pm,A,g_{\varepsilon}}-i\lambda)u\|+ C_{\Phi,g_{\varepsilon}}^{2}\|u\|\,,
$$
with $C_{\Phi,g_{\varepsilon}}^{2}=\frac{(C_{\Phi,g_{\varepsilon}}^{1})^{2}}{C_{g_{\varepsilon}}^{1}\|\Phi\|_{L^{\infty}}}$\,.
\end{proof}
\subsubsection{Adjoint}
\label{se.adjoincyl}
\begin{proposition}
\label{pr.hcyladj} Let $K_{\pm,A,g_{\varepsilon}}$ be the maximal
accretive realization of $P_{\pm,Q,g_{\varepsilon}}$ defined in
Proposition~\ref{pr.dyadbord}
with $\overline{Q}=(-\infty,0]\times\tz^{d-1}$\,, $\varepsilon\leq
\varepsilon_{g}$ small enough and $A$ satisfying
Hypothesis~\ref{hyp.Abdd} and \eqref{eq.commutAp}.
The adjoint of $K_{\pm,A,g_{\varepsilon}}$ is
$K_{\mp,A^{*},g_{\varepsilon}}$\,.
\end{proposition}
\begin{proof}
  The adjoint $A^{*}$ of $A$ fulfills Hypothesis~\ref{hyp.Abdd} and
  \eqref{eq.commutAp}. Thus $K_{\mp,A^{*},g_{\varepsilon}}$ is maximal
  accretive and shares the same properties as
  $K_{\pm,A,g_{\varepsilon}}$ while replacing $A$ with $A^{*}$\,,
  $\pm$ with $\mp$\,. 
Take first $u\in D(K_{\mp, A^{*},g_{\varepsilon}})$ and $v\in
D(K_{\pm,A,g_{\varepsilon}})$ such that $\supp u\subset B_{R_{u}}$ and
$\supp v\subset B_{R_{v}}$\,, where we recall
$$
B_{R}=\left\{g_{\varepsilon}^{ij}(q)p_{i}p_{j}\leq R^{2}\right\}\,.
$$
By applying Lemma~\ref{le.supppborne} with
$R=\max\left\{R_{u},R_{v}\right\}$ there exists
$\mathcal{T}=(b,S,T_{1},T_{2},\theta)$ such that 
\begin{eqnarray*}
K_{\mp,A^{*},g_{\varepsilon}}
&=&
U_{g_{\varepsilon}}K_{\mp,A^{*},g_{0}}^{1}(\mathcal{T})U_{g_{\varepsilon}}^{*}
\\
K_{\mp,A^{*},g_{0}}^{1}(\mathcal{T})&=& K_{\mp,A^{*},g_{0}}\mp
p_{i}b^{i}(q)\theta(|p|_{g_{0}(q)})\mp
S^{ij}_{k}(q)\theta(|p|_{g_{0}(q)})p_{i}p_{j}\partial_{p_{k}}
\\
&&\hspace{4cm}
+\frac{-\partial_{p}^{T}T_{1}(q)\partial_{p}+p^{T}T_{2}(q)p}{2}\,,
\\
K_{\pm,A,g_{\varepsilon}}
&=&
U_{g_{\varepsilon}}K_{\pm,A,g_{0}}^{1}(\mathcal{T})U_{g_{\varepsilon}}^{*}\\
K_{\pm,A,g_{0}}^{1}(\mathcal{T})&=& K_{\pm,A,g_{0}}\pm
p_{i}b^{i}(q)\theta(|p|_{g_{0}(q)})\pm
S^{ij}_{k}(q)\theta(|p|_{g_{0}(q)})p_{i}p_{j}\partial_{p_{k}}
\\
&&\hspace{4cm}
+\frac{-\partial_{p}^{T}T_{1}(q)\partial_{p}+p^{T}T_{2}(q)p}{2}\,.
\end{eqnarray*}
By setting $\tilde{u}=U_{g_{\varepsilon}}u\in D(K_{\mp,A^{*},g_{0}}^{1}(\mathcal{T}))=D(K_{\mp,A^{*},g_{0}})$ and
$\tilde{v}=U_{g_{\varepsilon}}v\in D(K_{\pm,A,g_{0}}^{1}(\mathcal{T})=D(K_{\pm,A,g_{0}}))$\,, the equality
$K_{\pm,A,g_{0}}^{*}=K_{\mp,A^{*},g_{0}}$ given in
Proposition~\ref{pr.plat} leads to
\begin{align*}
  \langle K_{\mp,A^{*},g_{\varepsilon}}u\,, v\rangle-\langle u\,,\,
  K_{\pm,A,g_{\varepsilon}}v\rangle
&=
\langle K_{\mp,A^{*},g_{0}}^{1}(\mathcal{T})\tilde{u}\,,\,
\tilde{v}\rangle
-
\langle \tilde{u}\,,\,
K_{\pm,A,g_{\varepsilon}}^{1}(\mathcal{T})\tilde{v}\rangle\\
&\hspace{-1cm}=\mp 2\langle \tilde{u}\,,\, p_{i}b^{i}(q)\tilde{v} \rangle
\pm \langle \tilde{u}\,,\, [p_{i}S^{ij}_{j}(q)+p_{j}S^{ij}_{i}(q)] \tilde{v}\rangle\,,
\end{align*}
where we used $\theta(|p|_{g_{0}(q)})\equiv 1$ in a neighborhood of $\supp
\tilde{u}\cup \supp\tilde{v}$\,.\\
The right-hand side is continuous on $L^{2}(Q;\mathcal{H}^{1})$\,. It
vanishes when $u, v\in \mathcal{C}^{\infty}_{0}(X;\mathfrak{f})$ by
direct calculations\,, and therefore when 
$\tilde{u}\,, \tilde{v}\in \mathcal{C}^{\infty}_{0}(X;\mathfrak{f})$
with $\theta\equiv 1$ in a neighborhood of $\supp
\tilde{u}\cup \supp\tilde{v}$\,. By density it always vanishes.
For $u\in D(K_{\mp,A^{*},g_{\varepsilon}})$ such that $\supp u\subset
B_{R_{u}}$\,, $R_{u}<+\infty$\,, we have proved
$$
\langle K_{\mp,A^{*},g_{\varepsilon}}u\,,\, v\rangle= \langle u\,,\, K_{\pm,A,g_{\varepsilon}}v\rangle
$$
for all $v\in D(K_{\pm,A,g_{\varepsilon}})$   which fulfill $\supp
v\subset B_{R_{v}}$ for some $R_{v}<+\infty$\,.
Those $v$'s are dense in $D(K_{\pm,A,g_{\varepsilon}})$ endowed with
the graph norm. Therefore $u\in D(K_{\pm,A,g_{\varepsilon}}^{*})$\,. 
All those $u$'s are dense in $D(K_{\mp,A^{*},g_{\varepsilon}})$ and we
obtain
$K_{\mp,A^{*},g_{\varepsilon}}\subset
K_{\pm,A,g_{\varepsilon}}^{*}$\,. Both operators are maximal accretive
and this yields the equality.
\end{proof}
\subsubsection{Density of $\mathcal{D}(\overline{X},j)$ in
  $D(K_{\pm,0,g_{\varepsilon}})$}
When
$\overline{Q}=(-\infty,0]\times \tz^{d-1}$
we recall that $u\in\mathcal{D}(\overline{X},j)$  is characterized by
\begin{eqnarray*}
  && u\in
  \mathcal{C}^{\infty}_{0}(\overline{X};\mathfrak{f})\quad,\quad
\gamma_{odd}u=0\,, i.e.~\gamma u(q',-p_{1},p')=j\gamma u(q',p_{1},p')\\
&& \partial_{q^{1}}u=\mathcal{O}(|q^{1}|^{\infty})\,.
\end{eqnarray*}
\begin{proposition}
\label{pr.core} Within the framework of Proposition~\ref{pr.dyadbord}
and in the case $A=0$\,, $\mathcal{D}(\overline{X},j)$ is dense in
$D(K_{\pm,0,g_{\varepsilon}})$ endowed with its graph norm.
\end{proposition}
\begin{proof}
Remember the notation
$$
B_{R}=\left\{(q,p)\in \overline{X}\,, g^{ij}_{\varepsilon}(q)p_{i}p_{j}\leq R\right\}
$$
By Proposition~\ref{pr.dyadbord} the set 
$$
\left\{u\in D(K_{\pm,0,g_{\varepsilon}})\,,\, \exists R_{u}>0\,,\, 
\supp u \subset B_{R_{u}}\right\}
$$
is dense in $D(K_{\pm,0,g_{\varepsilon}})$\,.\\
By Lemma~\ref{le.supppborne} the unitary change of variables \eqref{eq.unitch}\eqref{eq.unitch2},
\begin{eqnarray*}
  &&(U_{g_{\varepsilon}}v)(q,p)=\det(\Psi(q))^{-1/2}v(q,\Psi(q)^{-1}p)
\,,\\
&&\Psi(q)=g_{\varepsilon}(q)g_{0}(q)^{-1}\,,
\end{eqnarray*}
allows to write
$$
K_{\pm,0,g_{\varepsilon}}u= U_{g_{\varepsilon}}K_{\pm,0,g_{0}}^{1}(\mathcal{T})U_{g\varepsilon}^{*}u
$$
for some $\mathcal{T}=(b,S,T_{1},T_{2},\theta)$ fixed by $R>0$ as soon as $\supp
u\subset B_{R}$\,. But $\mathcal{D}(\overline{X},j)$ is dense in
$D(K_{\pm,0,g_{0}})=D(K_{\pm,0,g_{0}}(\mathcal{T}))$\,.
Since the unitary transform preserves  the set of $u\in
\mathcal{C}^{\infty}_{0}(\overline{X};\mathfrak{f})$ such that
$\gamma_{odd}u=0$\,, we deduce that
this set is dense in $D(K_{\pm,0,g_{\varepsilon}})$ endowed with its
graph norm.\\
Let $u\in \mathcal{C}^{\infty}_{0}(\overline{X};\mathfrak{f})$ satisfy 
$\gamma_{odd}u=0$\,. We now prove that $u$ can be approximated in
$D(K_{\pm, 0,g_{\varepsilon}})$ by elements of
$\mathcal{D}(\overline{X},j)$\,.
By Hadarmard's lemma (or Taylor expansion with integral remainder),
there exists a cut-off function $\chi\in
\mathcal{C}^{\infty}_{0}((-\infty,0])$\,, $\chi\equiv 1$ in a
neighborhood of $0$\,, and $v\in
\mathcal{C}^{\infty}_{0}(\overline{X},\mathfrak{f})$
such that
$$
u(q,p)=\chi(q^{1})u(0,q',p) + q^{1}v(q,p)\,.
$$ 
We take for $n\in\nz^{*}$\,, 
$$
u_{n}= u+(1-\chi(nq^{1}))q^{1}v(q,p)=u-\chi(nq^{1})q^{1}v(q,p)\,.
$$
Clearly $u_{n}\in \mathcal{D}(\overline{X},j)$ for all $n\in\nz^{*}$\,
$\lim_{n\to\infty}\|u_{n}-u\|=0$ and
$$
P_{\pm,g_{\varepsilon}}u_{n}=P_{\pm,g_{\varepsilon}}u
-\chi(nq^{1})q^{1}P_{\pm,g_{\varepsilon}}v 
-p_{1}\left[\chi(nq^{1})+\chi'(nq^{1})nq^{1}\right]v
$$
converges to $P_{\pm,g_{\varepsilon}}u$ in $L^{2}(X,dqdp;\mathfrak{f})$\,.
\end{proof}
\subsection{Global result}
\label{se.withougev}
We now end the  proof of  Theorem~\ref{th.main0} and
Theorem~\ref{th.mainA}
 by using a partition of unity in $q$\,. For this we assume
 $A=A(q,|p|_{q})$
like in
 \eqref{eq.intArq1}\eqref{eq.intArq0}\eqref{eq.intArq2}\eqref{eq.intArq3}
 or \eqref{eq.Arq1}\eqref{eq.Arq2}\eqref{eq.Arq3}. 
When $\overline{Q}=Q\sqcup \partial Q$ is a compact manifold or a
compact perturbation of the euclidean half-space
$\overline{\rz^{d}_{-}}$\,. Every interior chart domain diffeomorphic to a
bounded domain of $(-1,1)^{d}$ can be embedded in $\tz^{d}$\,, while
every boundary chart domain, diffeomorphic to  a bounded domain of
$(-\infty,0]\times (-1,1)^{d-1}$ can be embedded in the half-cylinder 
$(-\infty,0]\times \tz^{d-1}$\,. For every $q\in \partial Q$ one can
find a neighborhood $\mathcal{U}$ of $q$, embedded in $(-\infty,0]\times
\tz^{d-1}$\,, a metric $g_{\mathcal{U},\varepsilon}$ of the form 
\begin{eqnarray*}
g_{\mathcal{U},\varepsilon}=
\begin{pmatrix}
  1&0\\
0& m_{\mathcal{U},\varepsilon}(q^{1},q')
\end{pmatrix}\\
m_{,\mathcal{U}, \varepsilon}(q^{1},q')=m_{\mathcal{U},0}(q')+\chi_{m}(\frac{q^{1}}{\varepsilon})q^{1}\tilde{m}_{\mathcal{U}}(q^{1},q')\,,  
\end{eqnarray*}
coincide with $g$ on $\mathcal{U}$\,. Furthermore $\mathcal{U}$ can be
chosen small enough so that the condition $\varepsilon\leq
\varepsilon_{g_{\mathcal{U}}}$ of Proposition~\ref{pr.dyadbord} is
satisfied.
By compactness, there exists an (at most countable) locally finite
covering $\overline{Q}=\ccup_{\ell\in \mathcal{L}}\mathcal{U}_{\ell}$
so that the condition $\varepsilon\leq \varepsilon_{g_{\mathcal{U}}}$
is satisfies for all $\mathcal{U}$ such that $\mathcal{U}\cap \partial
Q\neq \emptyset$\,. The partition of unity
$(\chi_{\ell}(q))_{\ell\in \mathcal{L}}$ such that $\sum_{\ell\in
  \mathcal{L}}\chi_{\ell}^{2}(q)\equiv 1$ is subordinate to this
covering.
This partition of unity can be chosen so that
$$
\supp_{\ell'\in \mathcal{L}}\sharp \left\{\ell\in
  \mathcal{L}, \chi_{\ell}\chi_{\ell'}\neq 0\right\}\leq N_{\chi}<+\infty\,.
$$
\begin{proof}[Proof of  Theorem~\ref{th.main0} and Theorem~\ref{th.mainA}:]
Let $(\chi_{\ell}(q))_{\ell\in \mathcal{L}}$ be the above 
partition of unity.\\
We recall that  $u\in D(K_{\pm,A,g})$ is characterized by
\begin{eqnarray*}
  && u\in L^{2}(Q;\mathcal{H}^{1})\quad,\quad P_{\pm,Q,g}u\in
  L^{2}(X,dqdp; \mathfrak{f})\,,\\
&& \forall R>0\,, 1_{[0,R]}(|p|_{q})\gamma u \in L^{2}(\partial X,|p_{1}|dq'dp; \mathfrak{f})\,,\\
&& \gamma_{odd}u=\pm\sign(p_{1})A\gamma_{ev}u\,.
\end{eqnarray*}
The additional condition $A=A(q,|p|_{q})$ and 
this definition~of $D(K_{\pm,A,g})$  implies that for any $\ell\in \mathcal{L}$\,, 
$$
\left(u\in D(K_{A,g})\right)
\rightarrow
\left(\chi_{\ell}u\in D(K_{A,g})\right)\,.
$$
Moreover with $\chi_{\ell}=\chi_{\ell}(q)$\,, one gets
$$
\ad_{\chi_{\ell}}K_{\pm,A,g}=g^{ij}(q)p_{i}\partial_{q^{j}}\chi_{\ell}\quad,\quad
\ad_{\chi_{\ell}}^{2}K_{\pm,A,g}=0\,. 
$$
Like in Subsection~\ref{se.spaPU}, we infer the equivalences
\begin{eqnarray*}
  &&\left(
\frac{\|u\|_{H^{t}(\overline{Q};\mathcal{H}^{0})}}{\sum_{\ell\in \mathcal{L}}\|\chi_{\ell}u\|_{H^{t}(\overline{Q};\mathcal{H}^{0})}}
\right)\leq C_{t}\\
\text{and}&&
\left(
\frac{\|(K_{\pm,A,g}-i\lambda)u\|^{2}+\|u\|_{L^{2}(Q;\mathcal{H}^{1})}^{2}}{\sum_{\ell=1}^{L_{\varepsilon}}\|(K_{\pm,A,g}-i\lambda)\chi_{\ell}u\|^{2}+\|\chi_{\ell}u\|_{L^{2}(Q;\mathcal{H}^{1})}^{2}}
\right)^{\pm 1}\leq C_{\chi}\,,
\end{eqnarray*}
while the squared norms 
$\|u\|_{L^{2}(Q;\mathcal{H}^{s})}^{2}$ and $\|\gamma
u\|_{L^{2}(\partial X, |p_{1}|dq'dp;\mathfrak{f})}^{2}$ are equal to
\begin{eqnarray*}
\|u\|_{L^{2}(Q;\mathcal{H}^{s})}^{2}&=&\sum_{\ell\in
  \mathcal{L}}\|\chi_{\ell}u\|_{L^{2}(Q;\mathcal{H}^{s})}^{2}
\\
\|\gamma
u\|_{L^{2}(\partial X, |p_{1}|dq'dp;\mathfrak{f})}^{2}
&=&
\sum_{\ell\in \mathcal{L}}\|\gamma \chi_{\ell}u\|_{L^{2}(\partial X,
  |p_{1}|dq'dp;\mathfrak{f})}^{2}\,.
\end{eqnarray*}
Hence $u$ belongs to $D(K_{\pm,A,g})$ if and only if all the
$\chi_{\ell}u$'s belong to $D(K_{\pm,A,g})$ and 
$$
\sum_{\ell\in
  \mathcal{L}}\|K_{\pm,A,g}\chi_{\ell}u\|^{2}+\|\chi_{\ell}u\|^{2}_{L^{2}(Q,\mathcal{H}^{1})}< +\infty\,.
$$
This proves that the set $u\in D(K_{\pm,A,g})$ with a  compact
$q$-support are dense in $D(K_{\pm,A,g})$ endowed with its graph
norm. In particular when $A=0$\,, the results of
Proposition~\ref{pr.core} and Corollary~\ref{co.Leb} imply that
$\mathcal{D}(\overline{X},j)$ is dense in $D(K_{\pm,0,g})$ when $A=0$\,.\\
According to \eqref{eq.partunitRe}, 
the identity $\ad_{\chi_{\ell}}^{2}K_{\pm,A,g}=0$ implies
$$
\Real\langle u\,,\,
K_{\pm,A,g}u\rangle=\sum_{\ell\in \mathcal{L}}
\Real\langle\chi_{\ell}u\,,\,K_{\pm,A,g}\chi_{\ell}u\rangle\,,
$$
for all $u\in D(K_{\pm,A,g})$\,. This proves the accretivity of
$K_{\pm,A,g}-\frac{d}{2}$ and the integration by part identity
$$
\Real \langle u\, K_{\pm,A,g}u\rangle=
\|u\|_{L^{2}(Q;\mathcal{H}^{1})}^{2}+\Real\langle
\gamma_{ev}u\,,\,A\gamma_{ev}u\rangle_{L^{2}(\partial X,|p_{1}|dq'dp;\mathfrak{f})}\,,
$$
for all $u\in D(K_{\pm,A,g})$\,.\\
Since the subelliptic estimates of Proposition~\ref{pr.dyadbord} when
$\supp \chi_{\ell}\cap \partial Q\neq \emptyset$ and
Corollary~\ref{co.Leb} are uniformly satisfied for the local models of
$K_{\pm,A,g}\chi$\,, the above equivalence of norms imply
the subelliptic estimates of Theorem~\ref{th.main0} and
Theorem~\ref{th.mainA} by simple summation over $\mathcal{L}$\,.\\
The maximal accretivity is checked like in the final proof of
Proposition~\ref{pr.dyadbord}.
For $f\in L^{2}(X,dqdp;\mathfrak{f})$ take $u=\sum_{\ell\in
  \mathcal{L}}v_{\ell}$ with 
$$
v_{\ell}=\chi_{\ell}(K_{\pm,A,g}-i\lambda)^{-1}\chi_{\ell}f\,.
$$
The function $u$ belongs to $D(K_{\pm,A,g})$ and satisfies
$$
(K_{\pm,A,g}-i\lambda)u=f+\sum_{\ell\in
  \mathcal{L}}[g^{ij}(q)p_{i}\partial_{q^{j}}\chi_{\ell}(q)](K_{\pm,A,g}-i\lambda)^{-1}\chi_{\ell}f=(\Id
+B)f\,.
$$
Owing to $\supp_{\ell'\in \mathcal{L}}\sharp \left\{\ell\in
  \mathcal{L}, \chi_{\ell}\chi_{\ell'}\neq 0\right\}\leq N_{\chi}$\,,
we infer from the subelliptic estimates for $(K_{\pm,A,g}-i\lambda)$
$$
\|B\|_{\mathcal{L}(L^{2}(X,dqdp;\mathfrak{f}))}\leq C\langle \lambda\rangle^{-\frac{1}{8}}\,.
$$
Hence $\Id+B$ is invertible when $\lambda$ is chosen large enough and
$(K_{\pm,A,g}-i\lambda)$ is invertible. This proves the maximal
accretivity.\\
For $u\in D(K_{\pm,A,g})$ and $v\in D(K_{\mp,A^{*},g})$ we compute
\begin{align*}
&\langle v\,, K_{\pm,A,g}u\rangle-\langle K_{\mp,A^{*},g}v\,,
  u\rangle
-\sum_{\ell\in \mathcal{L}}
 \langle \chi_{\ell}v\,, K_{\pm,A,g}\chi_{\ell}u\rangle-\langle \chi_{\ell}K_{\mp,A^{*},g}v\,,
  \chi_{\ell}u\rangle\\
&\qquad=\langle v\,,\, (P_{\pm,Q,g}-\sum_{\ell\in
  \mathcal{L}}\chi_{\ell}P_{\pm,Q,g}\chi_{\ell})u\rangle
-\langle
(P_{\mp,Q,g}-\sum_{\ell\mathcal{L}}\chi_{\ell}P_{\mp,Q,g}\chi_{\ell})v\,,\,
u\rangle\\
&\qquad=\sum_{\ell\in \mathcal{L}}\langle v\,,\,
\chi_{\ell}g^{ij}(q)p_{i}(\partial_{q^{j}}\chi_{\ell}) u\rangle
+\langle \chi_{\ell}g^{ij}(q)p_{i}(\partial_{q^{j}}\chi_{\ell})v\,,\,
 u\rangle
=0\,.
\end{align*}
But every term $\langle \chi_{\ell}v\,,\,
K_{\pm,A,g}\chi_{\ell}u\rangle-\langle
K_{\mp,A^{*},g}\chi_{\ell}v\,,\, \chi_{\ell}u\rangle$ vanishes.
We have proved
$$
\forall u\in D(K_{\pm,A,g})\,,\, \forall v\in D(K_{\mp,A^{*},g})\,,
\langle v\,,\, K_{\pm,A,g}u\rangle=\langle K_{\mp,A^{*},g}v\,,\,u\rangle
$$
and the adjoint of $K_{\pm,A,g}$ equals $K_{\mp,A^{*},g}$\,.
\end{proof}

\section{Variations on a Theorem}
\label{se.variation}
In this section, some  straightforward consequences and variants of
Theorem~\ref{th.main0} and Theorem~\ref{th.mainA} are listed.
\subsection{Corollaries}
\label{se.coro}
We refer to the definitions and results of Section~\ref{se.cusp}.
\begin{corollary}
\label{co.KAcusp} 
Within the framework of Theorem~\ref{th.main0} and
Theorem~\ref{th.mainA} the operator $K_{\pm,A,g}$ is
$\frac{1}{4}$-pseudospectral. Its spectrum is contained in 
$S\cap \left\{\Real z\geq \frac{d}{2}\right\}$ with
$$
S=\left\{z\in\cz\,,|z+1|\leq C(\Real z +1)^{4}\,, \Real
  z\geq -1\right\}\,,
$$
and the resolvent estimate
$$
\forall z\not\in S\,, \quad \|(z-K_{\pm,A,g})^{-1}\|\leq
C\langle z\rangle^{-\frac{1}{4}}\,.
$$
The constant $C>0$ can be chosen so that the semigroup
$(e^{-tK_{\pm,A,g}})_{t\geq 0}$ is given by the convergent contour
integral
$$
e^{-tK_{\pm,A,g}}=\frac{1}{2\pi i}\int_{\partial
  S}e^{-tz}(z-K_{\pm,A,g})~dz\quad,\quad t>0\,,
$$
with $\partial S$ oriented from $+i\infty$ to $-i\infty$\,.\\
It satisfies the estimate
$$
\sup_{t>0}\|t^{7}K_{\pm,A,g}e^{-tK_{\pm,A,g}}\|< +\infty\,.
$$
\end{corollary}
The above contour integral easily implies exponential decay estimate
under a spectral gap condition.
\begin{corollary}
\label{co.expo} Within the framework of Theorem~\ref{th.main0} and
Theorem~\ref{th.mainA}, assume that the spectrum $\sigma
(K_{\pm,A,g})$ is partitioned into to parts
\begin{eqnarray*}
  &&\sigma(K_{\pm,A,g})\subset \sigma_{0}\cup \sigma_{\infty}\\
\text{with}
&& \sigma_{0}\subset \left\{z\in\cz\,,~\Real z \leq
 \mu_{0} \right\}\quad,\quad
\sigma_{\infty}\subset \left\{z\in\cz\,, \Real z \geq \mu\right\}\quad,\quad \mu>\mu_{0}\,.
\end{eqnarray*}
Let $\Pi_{0}$ be the spectral projection associated with
$\sigma_{0}$\,. For any $\tau <\mu$ there exists a constant
$C_{\tau}>0$ such that
$$
\forall t\geq 1\,,\quad
\|e^{-tK_{\pm,A,g}}-e^{-tK_{\pm,A,g}}\Pi_{0}\|\leq C_{\tau}e^{-\tau t}\,.
$$
\end{corollary}
\begin{corollary}
\label{co.KAcomp} Within the framework of Theorem~\ref{th.main0} and
Theorem~\ref{th.mainA} and when $\overline{Q}$ is compact, the
resolvent of $K_{\pm,A,g}$ is compact and its spectrum
$\sigma(K_{\pm,A,g})$ is discrete.
\end{corollary}
\subsection{PT-symmetry}
\label{se.PTsymm}
While studying accurately the spectrum of scalar Kramers-Fokker-Planck
operators on $Q=\rz^{d}$\,, H{\'e}rau-Hitrik-Sj{\"o}strand in
\cite{HHS2} 
used the following version of
PT-symmetry: When $Q=\rz^{d}$\,, the operator $K$ satisfies
\begin{equation}
  \label{eq.PT}
  UKU^{*}=K^{*}\quad\text{with}\quad Uu(q,p)=u(q,-p) \quad (U*=U)\,.
\end{equation}
We keep the same notation $U$ for the unitary action on $L^{2}(\partial
X,|p_{1}|dq'dp;\mathfrak{f})$\,.
\begin{proposition}
\label{pr.PT} Within the framework of Theorem~\ref{th.main0} and
Theorem~\ref{th.mainA} assume additionally that $A$ satisfies
$UAU^{*}=A^{*}$\,. Then $K_{\pm,A,g}$ is PT-symmetric according to
\eqref{eq.PT} and its spectrum is symmetric with respect to the real
axis, 
$$
\overline{\sigma(K_{\pm,A,g})}=\sigma(K_{\mp,A^{*},g})=\sigma(K_{\pm,A,g})\,.
$$
\end{proposition}
\begin{proof}
The relation $UP_{\pm,Q,g}U^{*}=P_{\mp,Q,g}$ is straightforward and
$U$ preserves the conditions
\begin{eqnarray*}
 && u\in L^{2}(Q;\mathcal{H}^{1})\quad,\quad P_{\pm,Q,g}u\in
  L^{2}(X,dqdp; \mathfrak{f})\,,\\
&& \forall R>0\,, 1_{[0,R]}(|p|_{q})\gamma u \in L^{2}(\partial X,|p_{1}|dq'dp; \mathfrak{f})\,,\\
\end{eqnarray*}
which occurs the definition of $D(K_{\pm,A,g})$ and of $D(K_{\mp,A^{*},g})$\,.
With 
\begin{align*}
\gamma_{odd}(U^{*}u)(q',p)&=\frac{\gamma
    (U^{*}u)(q',p_{1},p')-j\gamma (U^{*}u)(q',-p_{1},p')}{2}\\
&=\frac{\gamma u(q',-p_{1},-p')-j\gamma u(q',p_{1},-p')}{2}
=(U^{*}\gamma_{odd}u)(q',p)\,\\
\gamma_{ev}(U^{*}u)(q',p)&=\frac{\gamma
    (U^{*}u)(q',p_{1},p')+j\gamma (U^{*}u)(q',-p_{1},p')}{2}
\\
&=\frac{\gamma u(q',-p_{1},-p')+j\gamma u(q',p_{1},-p')}{2}
=(U^{*}\gamma_{ev}u)(q',p)\,\\
\text{and}\hspace{2cm}& U\text{sign}(p_{1})U^{*}=-\text{sign}(p_{1})
\end{align*}
the condition $UAU^{*}=A^{*}$ leads to
\begin{equation*}
  [\gamma_{odd}(U^{*}u)]
=\mp
\sign(p_{1})A^{*}\gamma_{ev}(U^{*}u)\,,
\end{equation*}
when $u\in D(K_{\pm,A,g})$\,.\\
Obviously,
$u\in D(K_{\pm, A,g})$  is equivalent to $U^{*}u\in
D(K_{\mp,A^{*},g})=D(K_{\pm,A,g}^{*})$ and this ends the proof.
\end{proof}
\subsection{Adding a potential}
\label{se.poten}
When we add a potential the energy is 
$\mathcal{E}_{V}(q,p)=\frac{|p|^{2}_{q}}{2}+V(q)$ and the Hamiltonian
vector field is
$$
\mathcal{Y}_{\mathcal{E}_{V}}=\mathcal{Y}_{\mathcal{E}}-\partial_{q^{i}}V(q)\partial_{p_{i}}\,.
$$
The corresponding Kramers-Fokker-Planck operator is
$$
P_{\pm,Q,g}(V)=P_{\pm,Q,g}\mp \partial_{q^{i}}V(q)\partial_{p_{i}}\,.
$$
\begin{proposition}
\label{pr.poten}
When $V$ is a globally Lipschitz function on $\overline{Q}$\,,
the operator
$K_{\pm,A,g}(V)-\frac{d}{2}=K_{\pm,A,g}-\frac{d}{2}\mp\partial_{q^{i}}V(q)\partial_{p_{i}}$
with the domain $D(K_{\pm,A,g}(V))=D(K_{\pm,A,g})$ is maximal
accretive and shares the same properties as $K_{\pm,A,g}$ by simply
changing the constants in the subelliptic estimates.\\
Its adjoint is $K_{\mp,A^{*},g}(V)$\,.\\
When $\overline{Q}$ is compact, the resolvent of $K_{\pm,A,g}$ is
compact and its spectrum is discrete.\\
If $UAU^{*}=A^{*}$ with $Uu(q,p)=u(q,-p)$\,, then $K_{\pm,A,g}(V)$ satisfies the PT-symmetry
property \eqref{eq.PT} and 
$$
\sigma(K_{\pm,A,g}(V))=\overline{\sigma(K_{\mp,A^{*},g}(V))}=\overline{\sigma(K_{\pm,A,g}(V))}\,.
$$
\end{proposition}
\begin{proof}
The subelliptic estimates of Theorem~\ref{th.main0} and
Theorem~\ref{th.mainA} include
$$
\forall u\in D(K_{\pm,A,g})\,,
\langle
\lambda\rangle^{\frac{1}{8}}\|u\|_{L^{2}(Q;\mathcal{H}^{1})}\leq C\|(K_{\pm,A,g}-i\lambda)u\|\,.
$$
Combined with Proposition~\ref{pr.perturb} applied with
$K=C+K_{\pm,A,g}$ and $C>0$ large enough and
Corollary~\ref{co.perturb}, the inequality
$$
\|\partial_{q^{i}}V(q)\partial_{p_{i}}u\|\leq
\|\partial_{q}V\|_{L^{\infty}}\|u\|_{L^{2}(Q;\mathcal{H}^{1})}
$$
yields the result.
\end{proof}
\begin{remark}
\label{re.co} Of course all the consequences listed in
Subsection~\ref{se.coro} are valid.
\end{remark}
\begin{remark}
\label{re.flow} The hamiltonian flow is not well defined under the
sole assumption that $V$ is Lipschitz continuous. It is not a surprise
that the dynamics (the semigroup) is well defined as soon as the
diffusion term $\mathcal{O}_{Q,g}$ is added. It was already observed in
\cite{HelNi} that the Kramers-Fokker-Planck operator on $\rz^{2d}$ with any
$\mathcal{C}^{\infty}$ potential $V\in \mathcal{C}^{\infty}(\rz^{d})$ is
essentially maximal accretive on
$\mathcal{C}^{\infty}_{0}(\rz^{2d})$\,, although the hamiltonian flow
is not always well defined.
\end{remark}
\subsection{Fiber bundle version}
\label{se.connec}
The final proof of Theorem~\ref{th.main0} and Theorem~\ref{th.mainA}
reduces the problem to local ones via spatial partition of unities.
There is therefore no difficulties to replace $\overline{Q}\times
\mathfrak{f}$ by some Hermitian bundle.
More precisely we assume that $\pi_{F}:F\to \overline{Q}$ is a smooth finite
dimensional Hermitian bundle, with typical fiber $F_{q}\sim
\mathfrak{f}$\,,  endowed with a connection $\nabla^{F}$ which is an
$\text{End}(F)$-valued 1-form on $TQ$\,. Locally a section of $F$, $s\in
\mathcal{C}^{\infty}(Q;F)$ may be written 
$$
s=\sum_{k=1}^{d_{F}}s_{k}(q)f_{k}\quad\text{with}~\mathfrak{f}=\oplus_{k=1}^{d_{F}}\cz f_{k}
$$
and the covariant derivative $\nabla^{F}_{\partial_{q^{j}}}s$ equals
$$
\nabla^{F}_{\partial_{q^{j}}}[s(q)]=\sum_{k=1}^{d_{F}}(\partial_{q^{j}}s_{k})(q)f_{k}+s(q)(\nabla^{F}_{\partial_{q^{j}}}f_{k})\,.
$$
This connection is
 compatible
with the hermitian structure of $F$ when
$$
\partial_{q^{j}}\langle f_{k}\,,\,f_{k'}\rangle_{g^{F}}=\langle
\nabla^{F}_{\partial_{q^{j}}}f_{k}\,,\,
f_{k'}\rangle_{g^{F}}
+
\langle f_{k}\,,\,
\nabla^{F}_{\partial_{q^{j}}} f_{k'}\rangle_{g^{F}}\,,
$$
for all $(j,k,k')\in \left\{1,\ldots,d\right\}\times\left\{1,\ldots,d_{F}\right\}^{2}$\,.\\
When $\pi:X=T^{*}Q\to Q$ is the natural projection, we shall work with
the fiber bundle $\pi_{F_{X}}: F_{X}=\pi^{*}F\to X$ which also equals
$T^{*}Q\otimes_{Q} F=X\otimes_{Q}F$\,. The tangent bundle $TX$ is
decomposed into $TX=(TX)^{H}\oplus (TX)^{V}$ with 
\begin{eqnarray*}
  && (T_{x}X)^{H}=\sum_{j=1}^{d}\rz e_{i}\sim T_{q}Q\quad\text{with}\quad 
e_{j}=\partial_{q^{j}}+\Gamma^{\ell}_{ij}p_{\ell}\partial_{p_{i}}\\
\text{and}
&&
(T_{x}X)^{V}=\sum_{j=1}^{d}\rz \partial_{p_{j}}\sim T_{q}^{*}Q\,,
\end{eqnarray*}
where $x=(q,p)$ and the $\Gamma^{\ell}_{ij}$'s are the Christoffel
symbol for the metric $g$ on $Q$ (see Remark~\ref{re.sign} for the
sign convention).
Then the connection on $F_{X}$ is given by 
\begin{equation}
  \label{eq.connFX}
\nabla^{F_{X}}_{e_{j}}=\nabla_{\partial_{q^{j}}}^{F}\quad,\quad \nabla^{F_{X}}_{\partial_{p_{j}}}=0\,.
\end{equation}
For $1\leq j\leq \dim Q$\,,
the covariant derivative $\nabla^{F_{X}}_{e_{j}}(sf)=\nabla^{F_{X}}_{e_{j}}[s(q,p)f]$ then equals
$$
\nabla^{F_{X}}_{e_{j}}(sf)=(e_{j}s)f+s\nabla^{F}_{\partial_{q^{j}}}f=(\partial_{q^{j}}s)f+(\Gamma^{\ell}_{ij}p_{\ell}\partial_{p_{i}}s)+
s\nabla^{F}_{\partial_{q^{j}}}f
$$
while
$$
\nabla^{F_{X}}_{\partial_{p_{j}}}(sf)=(\partial_{p_{j}}s)f\,.
$$
When $\nabla^{F}$ is compatible with the metric $g^{F}$\,,
$\nabla^{F_{X}}$ is compatible with $g^{F_{X}}=\pi^{*}g^{F}$\,.
\begin{definition}
\label{de.genGKFP} 
Assume that $(F,\nabla^{F},g^{F})$ is a hermitian fiber
bundle on $Q$ with the connection $\nabla^{F}$ and the metric
$g^{F}$\,. Let $(F_{X},g^{F_{X}},\nabla^{F_{X}})$ be the pull-back,
$F_{X}=\pi^{*}F$\,, by the projection $\pi:X=T^{*}Q\to Q$ with
$g^{F_{X}}=\pi^{*}g^{F}$ and $\nabla^{F_{X}}$ defined by \eqref{eq.connFX}. 
A geometric Kramers-Fokker-Planck operator is a differential operator
on $\mathcal{C}^{\infty}(X;F^{X})$ of the form
$$
P_{\pm,Q,g}^{F}+M(q,p,\partial_{p})=
P_{\pm,Q,g}^{F}+M_{j}^{0}(q,p)\nabla^{F_{X}}_{\partial_{p_{j}}}
+M^{1}(q,p)\,,
$$
where 
\begin{eqnarray*}
  && P_{\pm,Q,g}^{F}=\pm
g^{ij}(q)p_{i}\nabla^{F_{X}}_{e_{j}}+\mathcal{O}_{Q,g}\,,\\
\text{with}&&\mathcal{O}_{Q,g}=\frac{-\Delta_{p}+|p|^{2}_{q}}{2}\,,\\
&& e_{j}=\partial_{q^{j}}+\Gamma^{\ell}_{ij}p_{\ell}\partial_{p_{i}}\,,
\end{eqnarray*}
and where
$M_{j}^{0}\,, M^{1}\in \mathcal{C}^{\infty}(X;\text{End}(\pi^{*}F))$
satisfy
\begin{eqnarray}
\label{eq.M0symb}&&\|\partial_{p}^{\alpha}\partial_{q}^{\beta}M^{0}_{j}(q,p)\|\leq
  C_{\alpha,\beta}\langle p\rangle_{q}^{-|\alpha|}\,,\\
\label{eq.M1symb}
\text{and}
&&\|\partial_{p}^{\alpha}\partial_{q}^{\beta}M^{1}(q,p)\|\leq
  C_{\alpha,\beta}\langle p\rangle_{q}^{1-|\alpha|}\,,
\end{eqnarray}
for all multi-indices $(\alpha,\beta)\in\nz^{2d}$\,.
\end{definition}
The spaces $L^{2}(Q;\mathcal{H}^{s'})$\,, $s'\in\rz$\,, and
$H^{s}(\overline{Q};\mathcal{H}^{s'})$
  are defined locally as spaces of sections of $F_{X}$\,.\\
Let us specify the framework for boundary conditions. The fiber bundle
$
\pi_{F_{\partial X}}: F_{\partial X}=\pi^{*} F_{\partial Q}\to \partial X$ is the pull-back of
$F_{\partial Q}$ by the projection $\pi\big|_{\partial X}:\partial
X\to \partial Q$ and equals $\partial X\otimes_{\partial Q}F_{\partial
Q}$\,. Traces will lie in $L^{2}(\partial X,|p_{1}|dq'dp;F_{\partial
X})$\,. 
We assume that $F_{\partial Q}$ is endowed with an
involution  $\mathbf{j}\in \mathcal{C}^{\infty}(\partial
Q;\text{End}(F_{\partial Q}))$ such that for the metric $g^{F}$\,,
$\mathbf{j}^{*}=\mathbf{j}=\mathbf{j}^{-1}$\,. The regularity of
$\mathbf{j}$ ensures that for all
$q\in \partial Q$\,, the fiber $F_{q}$ can be decomposed into the
orthogonal direct sum
$F_{q}=\ker(\mathbf{j}(q)-\Id)\oplus^{\perp}\ker(\mathbf{j}(q)+\Id)$
where both parts have constant dimensions. Equivalently there exists
 a unitary mapping $U\in
C^{\infty}(\text{End}(F_{\partial Q};\partial Q\times \mathfrak{f}))$
such that $U(q)\mathbf{j}U(q)^{*}=j$  for a.e. $q\in \partial Q$\,,
where $j=j^{*}=j^{-1}\in \text{End}(\mathfrak{f})$ is constant.\\
The mapping which associates to  $\gamma(q,p)\in L^{2}(\partial X,|p_{1}|dq'dp;F_{\partial X})$
the function $U(q)\gamma(q,p)\in L^{2}(\partial X,
|p_{1}|dq'dp;\mathfrak{f})$ is a unitary isomorphism.
\begin{definition}
  \label{de.BCfiber}
Let $(F,g^{F})$ be endowed with a unitary involution $\mathbf{j}$
which belongs to $\mathcal{C}^{\infty}(\partial Q;\text{End}(F_{\partial Q}))$\,. For
$\gamma\in L^{2}_{loc}(\partial X,|p_{1}|dq'dp;F_{\partial X})$ we
define
\begin{eqnarray*}
&&
\gamma_{ev}(q,p_{1},p')=\Pi_{ev}\gamma(q,p_{1},p')=\frac{\gamma(q,p_{1},p')+\mathbf{j}(q')\gamma(q,-p_{1},p')}{2}\\
\text{and}
&&\gamma_{odd}(q,p_{1},p')=
\Pi_{odd}\gamma(q,p_{1},p')=\frac{\gamma(q,p_{1},p')-\mathbf{j}(q')\gamma(q,-p_{1},p')}{2}\,.
\end{eqnarray*}
A bounded operator $\mathbf{A}$ on $L^{2}(\partial X, |p_{1}|dq'dp;F_{\partial
  X})$ is admissible if $U(q)\mathbf{A}U^{*}(q)$ on $L^{2}(\partial X,
|p_{1}|dq'dp;\mathfrak{f})$ has the form $A(q,|p_{q}|)$ and fulfills
the conditions \eqref{eq.intArq1}\eqref{eq.intArq0} and either
\eqref{eq.intArq2} or \eqref{eq.intArq3}.
\end{definition}
\begin{proposition}
\label{pr.bundle}
Assume that $\overline{Q}$ is compact or a compact perturbation of the
euclidean half-space $\overline{\rz}^{d}_{-}$ and set
$\overline{X}=T^{*}\overline{Q}$\,. In the second case, the fiber bundle
$(F,g^{F},\nabla^{F})$ is assumed to coincide with $(Q\times
\mathfrak{f},g^{F}_{0},0)$\,, $\partial_{q}g_{0}^{F}\equiv 0$\,, while
the pair $(F_{\partial Q},\mathbf{j})$ is trivial outside a compact domain of
$\overline{\rz^{d}_{-}}$\\
Let the geometric Kramers-Fokker-Planck operator
$P_{\pm,Q,g}^{F}+M(q,p,\partial_{p})$
satisfy the conditions of Definition~\ref{de.genGKFP} and assume that the
bounded operator
$\mathbf{A}$ on $L^{2}(\partial X,|p_{1}|dq'dp;F_{\partial X})$ is
admissible according to Definition~\ref{de.BCfiber}.\\
There exists a constant $C>0$ such that the operator
$C+K_{\pm,\mathbf{A},g,M}^{F}=C+P_{\pm,Q,g}^{F}+M(q,p,\partial_{p})$
defined with
 the domain
$$
D(K_{\pm,\mathbf{A},g,M}^{F})=\left\{u\in L^{1}(Q;\mathcal{H}^{1})\,,\quad
  \begin{array}[c]{l}
   [P_{\pm,Q,g}^{F}+M(q,p,\partial_{p})]u\in L^{2}(X,dqdp;F_{X})\\
    \gamma_{odd}u=\pm \sign(p_{1})\mathbf{A}\gamma_{ev}u
  \end{array}
  \right\}
$$
 is maximal accretive and satisfies
the same subelliptic estimates as in Theorem~\ref{th.main0} when $\mathbf{A}=0$
and as Theorem~\ref{th.mainA} when $\mathbf{A}\neq 0$\,.\\
When $\overline{Q}$ is compact, $K_{\pm,\mathbf{A},g}^{F}$ has a
 compact resolvent and its spectrum is discrete.\\
The domain $D(K_{\pm,\mathbf{A},g;M}^{F})=D(K_{\pm,\mathbf{A},g;0}^{F})$ does not
depend on $M=M(q,p,\partial_{p})$\,. In particular when $A=0$ and for
any such $M$\,, 
$$
D(K_{\pm,\mathbf{A},g;M}^{F})\cap
\left\{u\in \mathcal{C}^{\infty}_{0}(\overline{X};F_{X})\,, \partial_{q^{1}}u=\mathcal{O}(|q^{1}|^{\infty})
 ~\text{near}~\partial X\right\}
$$
 is dense in
$D(K_{\pm,\mathbf{A},g;M}^{F})$ endowed with its graph norm.\\
If  additionally  
the connection $\nabla^{F}$ is compatible with the metric $g^{F}$\,,
then the following properties are true:
The integration by part identity
$$
\|u\|_{L^{2}(Q;\mathcal{H}^{1})}^{2}+\Real\langle \gamma_{ev}u\,,\,
\mathbf{A}\gamma_{ev}(u)\rangle_{L^{2}(\partial
  X,|p_{1}|dq'dp;F_{\partial X})}
=\Real\langle u\,,\, (K_{\pm, \mathbf{A},g;0}^{F}-i\lambda)u\rangle
$$
holds for all $\lambda\in\rz$ and all $u\in D(K_{\pm,\mathbf{A},g;0}^{F})$\,.\\
The adjoint of
$K_{\pm,\mathbf{A},g;M}^{F}=K_{\pm,\mathbf{A},g;0}^{F}+M^{0}_{j}(q,p)\nabla^{F_{X}}_{\partial_{p_{j}}}+M^{1}(q,p)$
equals
$$
K_{\mp,\mathbf{A}^{*},g;M^{*}}^{F}=K_{\mp,\mathbf{A}^{*},g;0}^{F}-\nabla^{F_{X}}_{\partial_{p_{j}}}\circ M^{0}_{j}(q,p)^{*}+M^{1}(q,p)^{*}\,.
$$
When $U\mathbf{A}U^{*}=\mathbf{A}^{*}$\,, with $Uu(q,p)=u(q,-p)$\,, and
$M(q,-p,-\partial_{p})$ equals the formal adjoint $M(q,p,\partial_{p})^{*}$\,,
 $K_{\pm,\mathbf{A},g,M}^{F}$ satisfies the PT-symmetry \eqref{eq.PT} and
$$
\sigma(K_{\pm,\mathbf{A},g;M}^{F})=\overline{\sigma(K_{\mp,\mathbf{A}^{*},g;M^{*}}^{F})}=\overline{\sigma(K_{\pm,\mathbf{A},g;M}^{F})}\,.
$$
\end{proposition}
\begin{remark}
\label{re.corofiber} All the consequences listed in
Subsection~\ref{se.coro} are valid.
\end{remark}
\begin{proof}
\textbf{a)}
After introducing the suitable finite partition of unity, we can assume $F=Q\times
\mathfrak{f}$ and the problem is reduced to the comparison
\begin{align*}
\|(\hat{U}(q)[P_{\pm,Q,g}^{F}&+M(q,p,\partial_{p})]\hat{U}(q)^{*}-P_{\pm,Q,g}\otimes
\Id_{\mathfrak{f}})u\|
\\
&\leq \|M(q,p,\partial_{p})u\|+
\|g^{ij}(q)p_{i}(\hat{U}(q)(\partial_{q^{i}}\hat{U}^{*}(q))+\nabla_{\partial_{q^{j}}}^{F}u\|
\\
&\leq C\|u\|_{L^{2}(Q;\mathcal{H}^{1})}\,,
\end{align*}
where the local unitary transform $\hat{U}(q)$: a) is  nothing but $\Id$ 
when $u$ in supported in an interior chart; b) 
trivializes $\mathbf{j}$ (and $\hat{U}\mathbf{A}\hat{U}^{*}=A(q,|p_{q}|)$) when the support of $u$ meets $\partial
Q$\,.\\ Consider the two realization
$K_{\pm,A,g}^{g^{F}}$ and $K_{\pm,A,g}^{g^{\mathfrak{f}}}$ of
$P_{\pm,Q,g}\otimes \Id$ when $Q\times \mathfrak{f}$ is endowed with the variable
metric $g^{F}$ and when it is endowed with the constant metric
$g^{\mathfrak{f}}$\,. A simple conjugation shows that
$K_{\pm,A,g}^{g^{\mathfrak{f}}}$ is unitarily equivalent to
$K_{\pm,A,g}^{g^{\mathfrak{f}}}+M^{1}(q,p)$ where $M^{1}(q,p)$
satisfies \eqref{eq.M1symb}. 
We have found a local unitary transform $\tilde{U}(q)$ from
$L^{2}(X,dqdp;(\mathfrak{f},g^{\mathfrak{f}}))$ to
$L^{2}(X,dqdp;F_{X})$ such that 
$$
\|P_{\pm,Q,g;M}^{F}u-\tilde{U}(q)^{*}K_{\pm,A,g}^{g^{\mathfrak{f}}}
\tilde{U}(q)u\|\leq C\|u\|_{L^{2}(Q;\mathcal{H}^{1})}\,,
$$
for any $u\in D(K_{\pm,\mathbf{A},g}^{F})$ with $\supp u$ contained in
a chart of the partition of unity.
The subelliptic estimates hold for
$K_{\pm,A,g}^{g^{\mathfrak{f}}}$ because it corresponds exactly to the
situation of Theorem~\ref{th.main0} and Theorem~\ref{th.mainA}.
Therefore
$K_{\pm,\mathbf,g}^{F}$ is locally a relatively bounded perturbation
of $\tilde{U}^{*}K_{\pm,A,g}^{g^{\mathfrak{f}}}\tilde{U}$ and can be
treated like the perturbation by a potential in
Proposition~\ref{pr.poten}, i.e. by applying the general perturbation
result of Proposition~\ref{pr.perturb}. Putting together all the
pieces of the spatial partition of unity follows the line of
Subsection~\ref{se.withougev}. This proves the maximal accretivity of
$C+K_{\pm,\mathbf{A},g;M}^{F}$\,, the extension of
the  subelliptic estimates of Theorem~\ref{th.main0} and
Theorem~\ref{th.mainA} and as well as the density statement when $A=0$\,.\\
The global comparison 
$$
\forall u\in D(K_{\pm,\mathbf{A},g;0})\,,\quad
\|(P_{\pm,Q,g}^{F}+M(q,p,\partial_{p}))u-K_{\pm,\mathbf{A},g;0}^{F}u\|\leq C\|u\|_{L^{2}(Q;\mathcal{H}^{1})}\,,
$$
implies that $K_{\pm,\mathbf{A},g;M}^{F}$ can be treated as a
perturbation of $K_{\pm,\mathbf{A},g;0}^{F}$ like in
Proposition~\ref{pr.perturb} and $D(K_{\pm,\mathbf{A},g;M}^{F})$ does
not depend on $M=M(q,p,\partial_{p})$\,.\\
\noindent\textbf{b)} Assume now that the connection $\nabla^{F}$ is
compatible with the metric $g^{F}$\,.
Then the connection $\nabla^{F_{X}}$ is compatible with the metric
$g^{F_{X}}$\,.
After using a partition of unity in $q$ and the integration by parts
formula for the models $K_{\pm,A,g}^{g^{\mathfrak{f}}}$ for which all
traces are well defined\,, one gets
 for all $u\in D(K_{\pm,A,g;0}^{F})$
\begin{eqnarray*}
\Real\langle u\,,\, K^{F}_{\pm,\mathbf{A},g;0}u\rangle&=&
\|u\|_{L^{2}(Q;\mathcal{H}^{1})}^{2}\pm \int_{X}\mathcal{L}_{g^{ij}(q)p_{i}e_{j}}\left(\langle
  u\,,\,u\rangle_{g^{F}}\right)~dqdp\\
&=&
\|u\|_{L^{2}(Q;\mathcal{H}^{1})}^{2} \pm \int_{\partial X}\langle
  \gamma u\,,\,\sign(p_{1})\gamma u\rangle_{g^{F}}~|p_{1}|dq'dp
\\
&=&
\|u\|_{L^{2}(Q;\mathcal{H}^{1})}^{2}+ \Real \langle \gamma_{ev}u\,,\,
A\gamma_{ev}u\rangle_{L^{2}(\partial X, |p_{1}|dq'dp;F)}\,.
\end{eqnarray*}
 The equality
$(K_{\pm,A,g}^{F})=K_{\mp,A^{*},g}^{F}$ and the PT-symmetry when
$UAU^{*}=A^{*}$  are deduced with the same argument from
the compatibility of the connection with the hermitian
structure.
\end{proof}
\section{Applications}
\label{se.appli}
In this section we introduce several kinds of boundary conditions and
we check that they enter in our formalism. They are motivated by the
stochastic process formulation or by the relationship between the
geometric Kramers-Fokker-Planck equation (hypoelliptic Laplacian
introduced by Bismut) and the Witten Laplacian. These 
boundary conditions will be introduced by doing formal
calculations or by considering simple models. Complete justifications
from the analysis or probabilistic point of view require additional work, for
which our regularity results on the semigroup may be useful.
\subsection{Scalar Kramers-Fokker-Planck equations in a domain of $\rz^d$}
\label{se.scal}
Although the adjoint is obtained by changing the signed index $\pm$\,, we focus on
the case of 
\begin{equation}
  \label{eq.KFP+}
P_{+,Q}(V)=p.\partial_{q}-\partial_{q}V(q).\partial_{p}+\frac{-\Delta_{p}+|p|^{2}}{2}\,.
\end{equation}
We recall without proofs standard results about the relationship
between stochastic processes and parabolic PDE's on the euclidean
space, $\rz^{d}$\,. In particular we fix the interpretations in terms of
(Kramers)-Fokker-Planck equations and Einstein-Smoluchowski
equations. With this respect, it is simpler to discuss here the
case when the force field is the gradient of a potential. This could
be extended to the case of more general force fields while paying
attention to the computation of adjoint operators. We refer the reader to
\cite{Ris}\cite{Nel} for a more detailed introduction and to \cite{BisLNM}
and \cite{IkWa} for the extension to riemannian manifolds which is
more involved. Then we introduce a way to
think of boundary conditions as jump processes. Specific examples are
discussed afterwards.

\subsubsection{Einstein-Smoluchowski case}
\label{se.einsmo}
The pure spatial description of Brownian motion in a gradient field is provided by
the stochastic differential equation
\begin{equation}
  \label{eq.einsmo}
\left\{
  \begin{array}[c]{l}
  dQ=-\nabla V(Q)dt + \sqrt{\frac{2}{\beta}}dW_{t}\\
Q_{0}=q\,,
\end{array}
\right.
\end{equation}
where $dW$ denotes the $d$-dimensional white noise with covariance
matrix $E(dW(t)dW(t'))=\Id_{\rz^{d}}\delta(t-t')$ while $\beta>0$ is
the inverse temperature\,. The function $V$
is assumed to describe a confining $\mathcal{C}^{\infty}$ potential, the
confinement being replaced by boundary conditions on a bounded domain.
For $\tilde{u}_{0}\in \mathcal{C}^{\infty}_{0}(\rz^{d})$ the conditional
expectation 
$$
\tilde{u}(q,t)=E(\tilde{u}_{0}(Q_{t})|Q_{0}=q)
$$ 
solves the backward Kolmogorov equation
$$
\left\{
  \begin{array}[c]{l}
    \partial_{t}\tilde{u}=-\nabla
    V(q).\nabla_{q}\tilde{u}+\beta^{-1}\Delta_{q}\tilde{u}=-\beta^{-1}(-\partial_{q}+\beta\partial_{q}V(q)).\partial_{q}\tilde{u}=-\beta^{-1}L_{\beta
      V}\tilde{u}\\ 
   \tilde{u}(q,0)=\tilde{u}_{0}(q)\,.
  \end{array}
\right.
$$
When $Q_{0}$ is randomly distributed with the law $\tilde{\mu}_{0}=\tilde{\varrho}_{0}(q)dq$\,, the
expectation $E(\tilde{u}_{0}(Q_{t}))$ equals
\begin{align*}
E(\tilde{u}_{0}(Q_{t}))&=\int_{\rz^{d}}\tilde{u}(q,t)~d\tilde{\mu}_{0}(q)=
\int_{\rz^{d}}(e^{-t\beta^{-1}L_{\beta V}}\tilde{u}_{0})~d\tilde{\mu}_{0}
\\
&=\int_{\rz^{d}}\tilde{u}_{0}(q)~d\tilde{\mu}_{t}(q)
=\int_{\rz^{d}}\tilde{u}_{0}(q) \tilde{\varrho}(q,t)~dq\,,
\end{align*}
where $\tilde{\mu}_{t}=e^{-t\beta^{-1}L_{\beta V}^{*}}\tilde{\mu}_{0}$\,, 
or the density $\tilde{\varrho}(t)=e^{-t\beta^{-1} L_{\beta V}^{*}}\tilde{\varrho}_{0}$\,, solves
the Fokker-Planck equation
$$
\left\{
  \begin{array}[c]{l}
    \partial_{t}\tilde{\varrho}=\nabla_{q}.(\nabla
    V(q)\tilde{\varrho})+\beta^{-1}\Delta_{q}\tilde{\varrho}=\beta^{-1}\partial_{q}.(\partial_{q}+\beta \partial_{q}V(q))\tilde{\varrho}=-\beta^{-1}L_{\beta
      V}^{*}\tilde{\varrho}\\ 
   \tilde{\varrho}(q,0)=\tilde{\varrho}_{0}(q)\,.
  \end{array}
\right.
$$
With rather general confining assumptions
(e.g. $|\partial_{q}V(q)|\geq C^{-1}|q|^{\delta}-C$ for $\delta>0$\,, $\lim_{q\to\infty}V(q)=+\infty$ and $|\Hess
V(q)|=o(|\nabla_{q}V(q)|)$ as $q\to \infty$)
$\mu_{\infty}=\frac{e^{-\beta V(q)}dq}{\int_{\rz^{d}}e^{-\beta V(q)}~dq}=\varrho_{\infty}(q)~dq$ is the
unique invariant measure and the backward Kolmogorov equation can be
studied in $L^{p}(\rz^{d};d\mu_{\infty})$ for $p\in [1,+\infty]$\,.
Especially when $p=2$\,,
setting $u=e^{-\frac{\beta V(q)}{2}}\tilde{u}\in L^{2}(\rz^{d},dq)$ the backward
Kolmogorov equation is transformed into
$$
\left\{
  \begin{array}[c]{l}
    \partial_{t}u=-\Delta_{\beta V/2}^{(0)}u\\
    u(q,0)=u_{0}(q)=e^{-\frac{\beta V(q)}{2}}\tilde{u}_{0}(q)\,,
  \end{array}
\right.
$$
where $\Delta_{\beta V/2}^{(0)}$ is the Witten Laplacian acting on functions
$$
\Delta_{\beta V/2}^{(0)}=(-\partial_{q}+\frac{\beta}{2}\partial_{q}V(q)).
(\partial_{q}+\frac{\beta}{2}\partial_{q}V(q))
=-\Delta_{q}+\frac{\beta^{2}}{4}|\nabla_{q}V(q)|^{2}-\frac{\beta}{2}\Delta V(q)\,.
$$
With the suitable confining assumptions, $\Delta_{\beta V/2}^{(0)}$ is
self-adjoint on $L^{2}(\rz^{d},dq)$\,, its resolvent is compact and
$\ker(\Delta_{\beta V/2}^{(0)})=\cz e^{-\frac{\beta V(q)}{2}}$\,.\\
By setting $\varrho(q,t)=e^{\frac{\beta V(q)}{2}}\tilde{\varrho}(q,t)$\,, the
Fokker-Planck equation is also transformed into
$$
\left\{
  \begin{array}[c]{l}
    \partial_{t}\varrho=-\Delta_{\beta V/2}^{(0)}\varrho\\
    \varrho(q,0)=\varrho_{0}(q)=e^{\frac{\beta V(q)}{2}}\tilde{\varrho}_{0}(q)\,.
  \end{array}
\right.
$$
This property and the self-adjointness of $\Delta_{\beta V/2}^{(0)}$ is
summarized by saying that the process~\eqref{eq.einsmo} is
reversible.
Note the relation
$$
E(\tilde{u}_{0}(Q_{t}))=\int_{\rz^{d}}\tilde{u}(q,t)\tilde{\varrho}_{0}(q)~dq
=\int_{\rz^{d}}u(q,t) \varrho_{0}(q)~dq=\int_{\rz^{d}}u_{0}(q)\varrho(q,t)~dq\,,
$$
when $Q_{0}$ is randomly
distributed according to $\tilde{\varrho}_{0}(q)~dq=e^{-\frac{\beta V(q)}{2}}\varrho_{0}(q)~dq$\,.
\subsubsection{Langevin stochastic process}
\label{se.lan}
The Langevin theory of Brownian motion takes place in the phase-space
$x=(q,p)\in \rz^{2d}$ and the stochastic differential equation is
written
\begin{equation}
  \label{eq.langevin}
\left\{
  \begin{array}[c]{l}
    dq=pdt\\
    dp=-\partial_{q}V(q)dt -\nu pdt+\sqrt{\frac{2\nu}{\beta}}dW
  \end{array}
\right.
\end{equation}
where $dW$ denotes again the $d$-dimensional white noise with covariance
matrix $E(dW(t)dW(t'))=\Id_{\rz^{d}}\delta(t-t')$\,, $\beta>0$ is the
inverse temperature and $\nu>0$ is the
friction coefficient (the mass is set to $1$ here).
It is a non reversible process which has nevertheless 
very strong relations
with the reversible Einstein-Smoluchowski case. Actually the
 process \eqref{eq.einsmo} can
be proven to be some weak limit, in the large friction regime $\nu\to +\infty$\,,  of the
Langevin process \eqref{eq.langevin}. We refer the reader to
\cite{Nel} for a simple approach. Within the $L^{2}$-framework, a more
general approach proving the relationship
between the hypoelliptic Laplacian introduced by Bismut and the Witten
Laplacian  has been considered in \cite{BiLe}, using a Schur complement technique in the
spirit of \cite{SjZw}.\\
We shall use the notation $X_{t}=(q_{t},p_{t})$\,.  For
$\tilde{u}_{0}\in \mathcal{C}^{\infty}_{0}(\rz^{2d})$\,, the
expectation 
$$
\tilde{u}(x,t)=E(\tilde{u}_{0}(X_{t})|X_{0}=x)\,,
$$
solves 
$$
\left\{
  \begin{array}[c]{l}
\partial_{t}\tilde{u}=p.\partial_{q}\tilde{u}-\partial_{q}V(q).\partial_{p}\tilde{u}
-\nu p.\partial_{p}\tilde{u}+\frac{\nu}{\beta}\Delta_{p}\tilde{u}=-B_{V}\tilde{u}\\
\tilde{u}(x,0)=\tilde{u}_{0}(x)\,.
\end{array}
\right.
$$
When $X_{0}$ is randomly distributed with the law
$\tilde{\mu}_{0}=\tilde{\varrho}_{0}dqdp$\,, the expectation
$E(\tilde{u}_{0}(X))$ equals
\begin{align*}
E(\tilde{u}_{0}(X_{t}))&=\int_{\rz^{2d}}\tilde{u}(x,t)~d\tilde{\mu}_{0}(x)=
\int_{\rz^{d}}(e^{-tB_{V}}\tilde{u}_{0})~d\tilde{\mu}_{0}
\\
&=\int_{\rz^{2d}}\tilde{u}_{0}(x)~d\tilde{\mu}_{t}(x)
=\int_{\rz^{2d}}\tilde{u}_{0}(x) \tilde{\varrho}(x,t)~dx\,,
\end{align*}
where $\tilde{\mu}_{t}=e^{-tB_{V}^{*}}\tilde{\mu}_{0}$ or the density
$\tilde{\varrho}(t)=e^{-tB_{V}^{*}}\tilde{\varrho}_{0}$\,, solves the
equation
$$
\left\{
  \begin{array}[c]{l}
\partial_{t}\tilde{\varrho}=
-p.\partial_{q}\tilde{\varrho}+\partial_{q}V(q).\partial_{p}\tilde{\varrho}
+\nu \partial_{p}.(p\tilde{\varrho})+\frac{\nu}{\beta}\Delta_{p}\tilde{\varrho}=-\tilde{B}_{V}^{*}\tilde{\varrho}\\
\tilde{\varrho}(x,0)=\tilde{\varrho}_{0}(x)\,.
\end{array}
\right.
$$
With a confining potential the unique invariant measure is the
Maxwellian probability
$\tilde{\mu}_{\infty}(q,p)=\frac{e^{-\beta(\frac{|p|^{2}}{2}+V(q))}dqdp}{
\int_{\rz^{2d}}e^{-\beta(\frac{|p|^{2}}{2}+V(q))}dqdp}$\,.
As we did for the Einstein-Smoluchowski case, the PDE approach for
the $L^{2}$-theory is better understood after taking
$$
  u(q,p,t)=e^{-\frac{\beta}{2}(\frac{|p|^{2}}{2}+V(q))}\tilde{u}(q,p,t)\quad\text{and}\quad
\varrho(q,p,t)=e^{\frac{\beta}{2}(\frac{|p|^{2}}{2}+V(q))}\tilde{\varrho}(q,p,t)\,.
$$
The conjugation relations
\begin{eqnarray*}
 &&
 e^{\pm\frac{\beta}{2}(\frac{|p|^{2}}{2}+V(q))}\partial_{p}e^{\mp\frac{\beta}{2}(\frac{|p|^{2}}{2}+V(q))}=\partial_{p}\mp\frac{\beta}{2}p\\
\text{and}&&
 e^{\pm\frac{\beta}{2}(\frac{|p|^{2}}{2}+V(q))}(p.\partial_{q}-\partial_{q}V(q).\partial_{p})e^{\mp\frac{\beta}{2}(\frac{|p|^{2}}{2}+V(q))}=
(p.\partial_{q}-\partial_{q}V(q).\partial_{p})\,,
\end{eqnarray*}
imply that $\varrho(x,t)$ and $u(x,t)$ solve respectively
$$
\left\{
  \begin{array}[c]{l}
\partial_{t}\varrho=
-p.\partial_{q}\varrho+\partial_{q}V(q).\partial_{p}\varrho
+\frac{\nu}{\beta}(\partial_{p}-\frac{\beta}{2}p)(\partial_{p}+\frac{\beta}{2}p)\varrho\\
\varrho(x,0)=\varrho_{0}(x)=e^{\frac{\beta}{2}(\frac{|p|^{2}}{2}+V(q))}\tilde{\varrho}_{0}(x)\,.
\end{array}
\right.
$$
and
$$
\left\{
  \begin{array}[c]{l}
\partial_{t}u=
p.\partial_{q}u-\partial_{q}V(q).\partial_{p}u
+\frac{\nu}{\beta}(\partial_{p}-\frac{\beta}{2}p)(\partial_{p}+\frac{\beta}{2}p)u\\
u(x,0)=u_{0}(x)=e^{-\frac{\beta}{2}(\frac{|p|^{2}}{2}+V(q))}\tilde{u}_{0}(x)\,.
\end{array}
\right.
$$
Like in the Einstein-Smoluchowski case, we have the relation
\begin{align*}
  E(\tilde{u}_{0}(X_{t}))&=\int_{\rz^{2d}}\tilde{u}(q,p,t)\tilde{\varrho}_{0}(q,p)~dqdp
\\
&=\int_{\rz^{2d}}u(q,p,t)\varrho_{0}(q,p)~dqdp
=\int_{\rz^{2d}}u_{0}(q,p)\varrho(q,p,t)~dqdp
\end{align*}
when $X_{0}$ is randomly distributed according to
$\tilde{\varrho}_{0}~dqdp= e^{\frac{\beta}{2}(\frac{|p|^{2}}{2}+V(q))}\varrho_{0}(q,p)~dqdp$\,.\\
Taking $\beta=2$ and $\nu=1$ gives
$$
\partial_{t}\varrho = -(P_{+,\rz^{d}}-\frac{d}{2})\varrho\quad\text{and}\quad \partial_{t}u=-(P_{-,\rz^{d}}-\frac{d}{2})u\,.
$$
The operator \eqref{eq.KFP+} corresponds to the Kramers-Fokker-Planck
equation in the sense that it provides the evolution of probability
measures (divided by the square root of the Maxwellian density). 
In \cite{HelNi} it was proved that for any $V\in
\mathcal{C}^{\infty}(\rz^{d})$\,, $P_{\pm,\rz^{d}}$ initially defined
on $\mathcal{C}^{\infty}_{0}(\rz^{2d})$ has a unique
maximal accretive extension $K_{\pm,\rz^{d}}$\,. We can write
$\varrho(t)=e^{-t(K_{+,\rz^{d}}-\frac{d}{2})}\varrho_{0}$ and
$u(t)=e^{-t(K_{-,\rz^{d}}-\frac{d}{2})}u_{0}$\,.
The process is no more reversible but up to the conjugation by the
exponential factors, the
adjoint evolution is obtained by changing $p$ into $-p$ (see
Subsection~\ref{se.PTsymm}).
This works only for gradient force fields.\\
Finally keep in mind that the choice $\beta=2$  corresponds to our
normalization of the Kramers-Fokker-Planck operator.

\subsubsection{Jump process at the boundary}
\label{se.jump}
This paragraph is a heuristic introduction of boundary conditions for the
Kramers-Fokker-Planck equation in terms of stochastic processes.
We start with the Kramers-Fokker-Planck equation written for
$\tilde{\varrho}(x,t)=e^{\frac{|p|^{2}}{2}+V(q)}\varrho(x,t)$
(remember $\beta=2$) which is naturally interpreted as a probability
density when there is no loss of mass. The domain $Q$ is a regular
bounded 
open domain of the euclidean
space $\rz^{d}$ and we set $\overline{Q}=Q\sqcup \partial Q$\,,
$X=T^{*}Q$ and $\overline{X}=X\sqcup \partial X=T^{*}Q\sqcup
T^{*}_{\partial Q}Q$\,. The Kramers-Fokker-Planck equation in the
interior $X$ reads
\begin{equation}
  \label{eq.kramfp}
\left\{
  \begin{array}[c]{l}
\partial_{t}\tilde{\varrho}=
-p.\partial_{q}\tilde{\varrho}+\partial_{q}V(q).\partial_{p}\tilde{\varrho}
+\nu \partial_{p}.(p\tilde{\varrho})+\frac{\nu}{\beta}\Delta_{p}\tilde{\varrho}\\
\tilde{\varrho}(x,0)=\tilde{\varrho}_{0}(x)\,,
\end{array}
\right.
\end{equation}
and it has to be completed with boundary conditions along $\partial
X$\,. We keep the conventions of the Introduction and take in a
neighborhood $U\subset \overline{Q}$ of
point $q_{0}\in \partial Q$ a coordinate system such that
$U\cap \partial Q=\left\{q^{1}=0\right\}$ and the
euclidean metric on $\rz^{d}$ follows \eqref{eq.defgm}:
$$
|dq|^{2}=(dq^{1})^{2}+m_{ij}(q^{1},q')dq'^{i}dq'^{j}\,.
$$
The corresponding symplectic coordinates on $T^{*}U\subset
\overline{X}$ are written $(q,p)=(q^{1},q',p_{1},p')$\,, $p_{1}>0$
corresponding to the outgoing conormal component. The measure
$|p_{1}|dq'dp$ on $\partial X=T^{*}_{\partial Q}Q$ does not depend on
such a choice of coordinates and by introducing polar coordinates
$p=|p|\omega$ with $\omega \in S^{*}_{q}Q\sim\sz^{d-1}$\,, the same holds for
$|\omega_{1}|dq'd\omega=|p|^{-d}|p_{1}|dq'dp$\,. 
The boundary $\partial X$ is partitionned into
$$
\partial X =\partial X^{+}\sqcup \partial X^{-}\sqcup \partial X^{0}\,,
$$
where $\partial X^{+}$ is the set of strictly outgoing rays
$$
(0,q',p_{1},p')\in \partial X\Leftrightarrow p_{1}>0\,,
$$
$\partial X^{-}$ is the set of strictly incoming rays
$$
(0,q',p_{1},p')\in \partial X\Leftrightarrow p_{1}< 0\,,
$$
and $\partial X^{0}$ is the set of glancing rays
$$
(0,q',p_{1},p')\in \partial X\Leftrightarrow p_{1}=0\,.
$$
By assuming that it makes sense, the trace of $\tilde{\varrho}$
along $\partial X$ is denoted $\gamma \tilde{\varrho}$ and 
$\gamma_{+}\tilde{\varrho}=\gamma \tilde{\varrho}\big|_{\partial
  X^{+}}$ and $\gamma_{-}\tilde{\varrho}=\gamma \tilde{\varrho}\big|_{\partial
  X^{-}}$\,. 
When $\tilde{\varrho}(x,t)$ is smooth enough, a simple integration by
parts gives
$$
\frac{d}{dt}\int_{X}\tilde{\varrho}(q,p,t)~dqdp
=\int_{\partial X^{-}}\gamma_{-}\tilde{\varrho}(q,p)~|p_{1}|dq'dp - 
\int_{\partial X^{+}}\gamma_{+}\tilde{\varrho}(q,p)~|p_{1}|dq'dp\,.
$$
Here and in the sequel, the glancing region $\partial X^{0}$ is
assumed to have a null measure contribution and  the possible
singularities around this region are absorbed by the weight
$|p_{1}|$\,. The system is dissipative when the incoming flow 
$\int_{\partial X^{-}}\gamma_{-}\tilde{\varrho}(q,p)~|p_{1}|dq'dp$ is
less than the outgoing flow $\int_{\partial
  X^{+}}\gamma_{+}\tilde{\varrho}(q,p)~|p_{1}|dq'dp$\,.
The difference of fluxes may not vanish and it is convenient to send
the missing mass to an artificial exterior point $\mathfrak{e}$\,. The
relation between $\gamma_{+}\tilde{\varrho}|p_{1}|dq'dp$ and
$\gamma_{-}\tilde{\varrho}|p_{1}|dq'dp$ can then be modelled by a
Markov kernel $M(x^{-},x^{+})$ from $\partial X^{+}$ to $\partial
X^{-}\sqcup \left\{\mathfrak{e}\right\}$\,: $M(x^{-},x^{+})$ is a
$\partial X^{-}\sqcup \left\{\mathfrak{e}\right\}$-probability measure
valued function of $x^{+}\in \partial X^{+}$ such that
\begin{multline*}
\left[\int_{\partial X^{+}}\gamma_{+}\tilde{\varrho}(x^{+})|p_{1}|dq'dp-\int_{\partial X^{-}}\gamma_{-}\tilde{\varrho}(x^{-})|p_{1}|dq'dp\right]
\delta_{\mathfrak{e}}+ (\gamma_{-}\tilde{\varrho})(x^{-})|p_{1}|dq'dp\\
=\int_{\partial X^{+}}M(x^{-},x^{+})\gamma_{+} \tilde{\varrho}(x^{+})|p_{1}|dq'dp\,.
\end{multline*}
While considering the stochastic process \eqref{eq.langevin} with
initial data $X_{0}\in X\sqcup \partial X^{-}$\,, the quantity
$$
\tau_{x}=\inf\left\{t\in [0,+\infty)\,, X_{t}(x)\in \partial
  X^{+}\right\}\quad,\quad x\in X\sqcup \partial X^{-}\,,
$$
is a stopping time. The Langevin stochastic process
\eqref{eq.langevin} can thus be completed as a set of cadlag
trajectories in $\overline{X}\sqcup \left\{\mathfrak{e}\right\}$ by
introducing a jump process $J_{M}$ at times $\tau$ when
$X_{\tau^{-}}\in \partial X^{+}$\,:
\begin{equation}
  \label{eq.langevincompl}
\left\{
  \begin{array}[c]{l}
dX_{t}=
\begin{pmatrix}
  pdt\\
-\partial_{q}V(q)dt-pdt+dW
\end{pmatrix}
+dJ_{M}
\quad\text{when}~X_{t}=(q_{t},p_{t})\in \overline{X}\setminus \partial X^{0}\,\\
   dX_{t}=0\quad \text{for~all}~t\geq
   t_{0}~\text{if}~X_{t_{0}}=\partial X^{0}\sqcup\mathfrak{e}\,.
  \end{array}
\right.
\end{equation}
We chose the coefficients $\nu=1$ and $\beta=2$ according to the
discussion of the previous paragraph. The transition matrix of the
jump process $J_{M}$ is given by the Markov kernel
$M(x^{-},x^{+})$\,.\\
In order to be consistent with our framework, the jump process sending
$x^{+}$ to $x^{-}$ is
assumed to be local with respect to $q\in \partial Q$\,, which means
$x^{-}=(q,p^{-})$ or $x^{-}=\mathfrak{e}$ when $x^{+}=(q,p^{+})$\,,
and elastic\,, which means $|p_{-}|^{2}=|p_{+}|^{2}$ when $x^{-}\neq
\mathfrak{e}$\,. 
More precisely
the Markov kernel is assumed to satisfy
$$
M(x^{-},x^{+})\big|_{\partial X^{-}\times \partial
  X^{+}}=\delta(q^{+}-q^{-})\delta(|p^{-}|-|p^{+}|) \tilde{R}(q,r,\omega^{-},\omega^{+})
$$
where we used $q^{+}=q^{-}=q$ and $p^{\pm}=r\omega^{\pm}$ with
$r=|p^{+}|=|p^{-}|$ and $\omega^{\pm}\in S^{\pm}=\left\{\omega\in
  S^{d-1}\,,\, \pm \omega_{1}>0\right\}$\,. Allowing $\omega^{-}\in
S^{-}\sqcup \mathfrak{e}$\,, $\tilde{R}$ can be considered as a Markov
kernel from $S^{+}$ to $S^{-}\sqcup\left\{\mathfrak{e}\right\}$ such
that
\begin{multline*}
\left[\int_{S^{+}}\gamma_{+}\tilde{\varrho}(r\omega^{+})|\omega_{1}^{+}|d\omega^{+}-\int_{S^{-}}\gamma_{-}\tilde{\varrho}(r\omega^{-})|\omega_{1}^{-}|d\omega^{-}\right]\delta_{\mathfrak{e}}
+ \gamma_{-}\tilde{\varrho}(r\omega^{-})|\omega_{1}^{-}|d\omega^{-}
\\
=\int_{S^{+}}\tilde{R}(q,r,\omega^{-},\omega^{+})\gamma_{+}\tilde{\varrho}(q,r\omega^{+})|\omega_{1}^{+}|d\omega^{+}\,,
\end{multline*}
for all $(q,r)\in \partial Q \times (0,+\infty)$\,.
For the regularity we assume first that for all $(q,r)\in \partial Q
\times (0,+\infty)$\,, the Markov kernel sends
$L^{1}(S^{+},d\omega^{+})$ into $L^{1}(S^{-}\sqcup \mathfrak{e},
d\omega^{-}\oplus \delta_{\mathfrak{e}})$\,.\\
Set
$$
R(q,r,\omega^{-},\omega^{+})=|\omega_{1}^{-}|^{-1}\tilde{R}(q,r,\omega^{-},\omega^{+})|\omega_{1}^{+}|
\quad\text{for}~(\omega^{-},\omega^{+})\in S^{-}\times S^{+}\,,
$$
and for $\gamma \in L^{1}(S^{+},|\omega_{1}|d\omega)$\,,
$\gamma'=R(q,r)\gamma$ when
$$
\gamma'(\omega^{-})=
\int_{S^{+}}R(q,r,\omega^{-},\omega^{+})\gamma_{+}(\omega^{+})~d\omega^{+}
\quad \text{for}~a.e.~\omega^{-}\in S^{-}\,.
$$
The two conditions
\begin{eqnarray}
\label{eq.condcontr1}
&&  \tilde{R}(q,r,\{\mathfrak{e}\},\omega^{+})\geq
  \alpha\quad \text{for}~a.e.~\omega^{+}\in S^{+}\\
\label{eq.condcontr2}
\text{and}
&&
R(q,p)1_{S^{+}}\leq 1_{S^{-}}\,,
\end{eqnarray}
imply for any $s\in [1,+\infty]$ (interpolate between $s=1$ and $s=\infty$)
\begin{equation}
  \label{eq.contrLp}
\forall \gamma\in L^{s}(S^{+},|\omega_{1}|d\omega)\,,\quad
  \|R(q,p)\gamma\|_{L^{s}(S^{-},|\omega_{1}|d\omega)}\leq
  (1-\alpha)^{\frac{1}{s}}\|\gamma\|_{L^{s}(S^{+}, |\omega_{1}|d\omega)}\,.
\end{equation}
When $\mu_{t}$ is the law of $X_{t}\in\overline{X}\sqcup
\left\{\mathfrak{e}\right\}$ solving \eqref{eq.langevincompl} with a
suitable initial data $\mu_{0}$\,, we claim that
$\mu_{t}=\tilde{\varrho}(q,p,t)dqdp \oplus
(1-\int_{X}\tilde{\varrho}(q,p,t)~dqdp)\delta_{\mathfrak{e}}$ where
$\tilde{\varrho}$ solves the Kramers-Fokker-Planck equation
\eqref{eq.kramfp} with the boundary condition
$$
\gamma_{-}\tilde{\varrho}(q,r.)= R(q,r) \gamma_{+}\tilde{\varrho}(q,
r.)\quad\text{for}~a.e.~(q,r)\in \partial Q\times (0,+\infty)\,,
$$
shortly written
$\gamma_{-}\tilde{\varrho}=R\gamma_{+}\tilde{\varrho}$\,.
Since the operator $R$ preserves the energy (the jump process is local
and elastic), replacing $\tilde{\varrho}$ by
$\varrho(x,t)=e^{-(\frac{|p|^{2}}{2}+V(q))}\tilde{\varrho}(x,t)$ is
straightforward:
\begin{equation}
  \label{eq.kramfpbc}
\left\{
  \begin{array}[c]{l}
    \partial_{t}\varrho=-P_{+,Q}\varrho\,,\\
    \gamma_{-}\varrho=R \gamma_{+}\varrho \,,\\
     \varrho(x,t=0)=\varrho_{0}(x)\,.
  \end{array}
\right.
\end{equation}
Moreover the conditions \eqref{eq.condcontr2} and
\eqref{eq.condcontr1} with a uniform lower bound $\alpha\geq 0$
imply that $R$ is a contraction from $L^{2}(\partial X^{+},
|p_{1}|dq'dp)$ to $L^{2}(\partial X^{-}, |p_{1}|dq'dp)$ with norm
$(1-\alpha)^{\frac{1}{2}}$\,.
By following the arguments of Lemma~\ref{le.eqA+-} (with
$\mathfrak{f}=\cz$ and $j=1$) the boundary
condition can be rewritten $\Pi_{odd}\gamma
\varrho=\sign(p_{1})A\Pi_{ev}\gamma$ where $A$ is a maximal accretive
operator so that $\varrho(t)=e^{-tK_{+,A}}\varrho_{0}$\,. The
identification of $A$ with the corresponding assumptions will be
specified in the examples listed below. 
\subsubsection{Comments}
\label{se.commentsappl}
The first work to justify the stochastic trajectorial interpretation
of boundary conditions for the Einstein-Smoluchowski equation and the
Kramers-Fokker-Planck equation date back to Skorohod who introduced the
so called Skorohod reflection map in \cite{Sko1}\cite{Sko2}. 
The Einstein-Smoluchowski case, even when the force field is not a
gradient field, has been widely studied
because it relies on the standard regularity theory for elliptic
boundary value problems. We refer the reader for example to
\cite{IkWa}\cite{StVa} and to \cite{LiSz} where discontinuous boundary value
problems and corner problems were also taken into account.
For the Langevin process and the Kramers-Fokker-Planck equation little
seems to be known. The weak formulation developed for kinetic theory
in \cite{Luc}\cite{Car} does not provide any information on the
operator domain, nor any sufficient regularity. The one dimensional
case has been studied with specular reflection in \cite{Lap}  and with
some exotic non elastic case in \cite{Ber}. More recently the
half-space problem with specular reflection has been considered in \cite{BoJa}.

\subsubsection{Specular reflection}
\label{se.ref}
The specular reflection is a deterministic jump process which
transforms $(p_{1},p')$ into $(-p_{1},p')$ when the particle hits the
boundary $X_{\tau^{-}}\in \partial X^{+}$\,. Within the presentation
of Paragraph~\ref{se.jump}\,, the Markov kernel
$\tilde{R}(q,r,\omega^{-},\omega^{+})$ is simply given by
$$
\tilde{R}(q,r,\omega^{-},\omega^{+})=\delta(\omega_{1}^{-}+\omega_{1}^{+})
\delta({\omega'}^{-}-{\omega'}^{+})\,.
$$
The stochastic process \eqref{eq.langevincompl} takes place in
$\overline{X}$ and all the mass of the probability measure
$\tilde{\varrho}(q,p,t)dqdp=e^{-(\frac{|p|^{2}}{2}+V(q))}\varrho(q,p,t)dqdp$
lies in $\overline{X}$\,. We can forget the exterior $\mathfrak{e}$\,.
The corresponding boundary condition for $\varrho(x,t)$
will be
$$
\varrho(0,q',p_{1},p)=\varrho(0,q',-p_{1},p')\quad\text{for}~p_{1}<0\,.
$$
Within our formalism for boundary value problem for the
Kramers-Fokker-Planck operator $P_{+,Q}$, it can be written 
\begin{eqnarray*}
  &&\gamma_{odd}u(x,t)=0\,\\
\text{with}&& \mathfrak{f}=\cz\quad\text{and}\quad j=1\,.
\end{eqnarray*}
We can apply Theorem~\ref{th.main0} (case $A=0$) and its Corollaries of
   Subsection~\ref{se.coro} (see Subsection~\ref{se.poten} for adding a
potential). The PT-symmetry \eqref{eq.PT} is also satisfied and
Proposition~\ref{pr.PT} also applies.\\
The kernel $\Ker(K_{+,0})$ equals $\cz \varrho_{\infty}$ with
$\varrho_{\infty}(q,p)=e^{-(\frac{|p|^{2}}{2}+V(q))}$  and
$\tilde{\mu}_{\infty}=\frac{e^{-2(\frac{|p|^{2}}{2}+V(q))}dqdp}{\int_{\overline{X}}e^{-2(\frac{|p|^{2}}{2}+V(q))}dqdp}$
is the 
equilibrium probability measure. Corollary~\ref{co.expo} ensures the exponential
return to the equilibrium in the $L^{2}$-sense for the Langevin
process
with specular reflection in a bounded domain.\\
We expect that in the large friction limit (introduce the parameter
$\nu$ and consider the limit as $\nu\to +\infty$ while $\beta=2$), the
solution $\varrho(q,p,t)=e^{-tK_{+,0}}\varrho_{0}$ converges to
$f(q,t)e^{-\frac{|p|^{2}}{2}}$ with
$f(q,t)=e^{-t\Delta_{V}^{(0),N}}f_{0}$\,, where $\Delta_{V}^{(0),N}$ is the
Neumann realization of the Witten Laplacian
$\Delta_{V}^{(0)}=-\Delta +|\nabla V(q)|^{2}-\Delta V(q)$\,. This realization is
characterized by
$$
u\in D(\Delta_{V}^{(0),N})\Leftrightarrow  \left\{
\begin{array}[c]{ll}
u \in L^{2}(Q,dq)\quad,\quad &
\Delta_{V}^{(0)}u\in L^{2}(Q,dq)\,,\\
&   \partial_{n}u+\partial_{n}V(q)u\big|_{\partial Q}=0\,.
  \end{array}
\right\}
$$
where $\partial_{n}$ is the outgoing normal derivative. Intuitively it
corresponds to putting a $+\infty$ repulsive potential outside
$\overline{Q}$
 and to the specular reflection within the phase-space dynamical
 picture.
\subsubsection{Full absorption}
\label{se.absor}
Within the Einstein-Smoluchowski approach another natural boundary
condition is the homogeneous Dirichlet boundary condition for
$\tilde{\varrho}(q,t)$ or $\varrho(q,t)$\,. This means that the
particles are absorbed by the exterior as soon as they hit the
boundary $\partial Q$ and it can be intuitively presented by putting a
$-\infty$ attractive potential outside $\overline{Q}$\,. The evolution
of $\varrho(q,t)$ is given by
$\varrho(t)=e^{-t\Delta_{V}^{(0),D}}\varrho_{0}$ where
$\Delta_{V}^{(0),D}$ is the Dirichlet realization of the Witten
Laplacian $\Delta_{V}^{(0)}=-\Delta +|\nabla V(q)|^{2}-\Delta V(q)$\,.
Within the Langevin description this can be modelled by sending the
particle to the exterior $\mathfrak{e}$ as soon as it hits $\partial
X^{+}$\,.
Within the presentation of Paragraph~\ref{se.jump}, this deterministic
jump process is simply given by the Markov kernel
$$
\tilde{R}(q,r,\omega^{-},\omega^{+})=\delta_{\mathfrak{e}}\,.
$$
The corresponding boundary condition in \eqref{eq.kramfpbc} is simply
$$
\varrho(0,q',p_{1},p')=0\quad\text{for}~p_{1}<0\,.
$$
With our notations they can be simply written as
\begin{eqnarray*}
  &&\gamma_{odd}u(q,p)=\sign(p_{1})\gamma_{ev}u(q,p)\quad,\quad
q\in\partial Q\,,\\
\text{with}&& \mathfrak{f}=\cz\quad\text{and}\quad j=1\,.
\end{eqnarray*}
We are in the case when $A=1$ and Theorem~\ref{th.mainA} and its
Corollaries of Subsection~\ref{se.coro} apply (see
Subsection~\ref{se.poten} for the case with the potential $V$). Since
$j1j=1$\,, the PT-symmetry \eqref{eq.PT} is fulfilled and
Proposition~\ref{pr.PT} is valid.\\
We also expect some kind of convergence to the Einstein-Smoluchowski equation,
with the generator $\Delta_{V}^{(0),D}$ after the conjugation with the
exponential weight, in the large friction limit
$\nu\to \infty$\,. Note nevertheless that the density
$\int_{T^{*}_{q}Q}\varrho(q,p,t)~dp$ does not vanish in general when
$q\in \partial Q$\,. Some numerical simulations provided by
T.~Leli{\`e}vre in the one-dimensional case
 show that in the large friction limit 
exponentially decaying boundary layers persist in the neighborhood of $\partial Q$\,,
while comparing the local quantity 
 $\int_{\rz^{d}}e^{-(\frac{|p|^{2}}{2}+V(q))}\varrho(q,p,t)~dp$
 and the solution to the Einstein-Smoluchowski equation with Dirichlet
 boundary conditions.\\
A reason for studying such boundary conditions is that they are
related with quasi-stationnary distributions. We refer to
\cite{LbLLP} for the presentation for the Einstein-Smoluchowski equation.
 A detailed analysis in the
low temperature limit has
been performed  recently with the help of semiclassical tools for
Witten Laplacians, with T.~Leli{\`e}vre in \cite{LeNi}. 
The Langevin phase-space approach is a more natural framework for molecular dynamics.

\subsubsection{More general boundary conditions} 
\label{se.moregen}
In \cite{Car} the weak formulation of the Kramers-Fokker-Planck
equation is studied for more general conditions
$\gamma_{-}\varrho(q,r.)=R(q,r)\gamma_{+}\varrho(q,r.)$ similar to
\eqref{eq.kramfp}. Within the assumptions introduced in
Paragraph~\ref{se.jump}, let us check that they enter in our
formalism.
We assume in particular that the lower bound \eqref{eq.condcontr1}
holds with $\alpha>0$ uniform w.r.t to $(q,r)\in \partial Q\times
(0,+\infty)$\,.
We shall work here with $\mathfrak{f}=\cz$ and $j=1$ so that the
projections $\Pi_{+}$ and 
$\Pi_{-}$ defined on $L^{2}(\partial X, |p_{1}|dq'dp;\cz)$ (see Definition~\ref{de.invoL}) are
characterized by 
$$
\Pi_{+}\gamma (q,p_{1},p')=\gamma(q,|p_{1}|,p')\quad,\quad 
\Pi_{-}\gamma(q,p_{1},p')=\gamma(q,-|p_{1}|,p')\,.
$$
The relation $\gamma_{-}\varrho (q,r.)=R(q,r)\gamma_{+}\varrho(q,r.)$ can be
written $\Pi_{-}\gamma
\varrho(q,r.)=R'(q,r)\Pi_{+}\gamma\varrho(q,r.)$ with
$$
R'(q,r,\omega^{*},\omega)=
1_{\rz_{-}}(\omega_{1}^{*})R(q,r,\omega^{*},\omega)1_{\rz_{+}}(\omega_{1})
+ 1_{\rz_{+}}(\omega_{1}^{*})R(q,r, \widehat{\omega^{*}},\widehat{\omega} )1_{\rz_{-}}(\omega_{1})\,,
$$
with $\widehat{\omega}=(-\omega_{1},\omega')$ when
$\omega=(\omega_{1},\omega')$\,.
The estimate \eqref{eq.contrLp} imply that $R'$ is a contraction of
$\Pi_{ev}L^{2}(\partial X, |p_{1}|dq'dp)$ with the norm
$(1-\alpha)^{\frac{1}{2}}$\,.
The operator $R'$ commutes with $\Pi_{ev}$ and it is local in the
variables $(q,r=|p|)$\,. Since $R'$ is a strict contraction with norm
$(1-\alpha)^{1/2}< 1$ we can define the bounded operator
$$
A=\frac{1-R'}{1+R'}
$$
which is local in $(q,r=|p|)$ and bounded in $L^{2}(\partial
X,|p_{1}|dq'dp)$ with the norm $\|A\|\leq
\frac{1+(1-\alpha)^{1/2}}{1-(1-\alpha)^{1/2}}$\,.
Moreover the inequality
\begin{align*}
\|(1-A)u\|_{L^{2}(\partial X, |p_{1}|dq'dp)}^{2}
&=\|R'(1+A)u\|_{L^{2}(\partial X, |p_{1}|dq'dp)}^{2}\\
&\leq (1-\alpha)\|(1+A)u\|_{L^{2}(\partial X, |p_{1}|dq'dp)}^{2}
\end{align*}
implies
$$
\forall u\in L^{2}(\partial X,|p_{1}|dq'dp)\,,\quad
\alpha \|u\|_{L^{2}(\partial X, |p_{1}|dq'dp)}^{2}\leq 4\Real \langle
u\,,\, A u \rangle_{L^{2}(\partial X, |p_{1}|dq'dp)}\,.
$$
Hence  the operator $A=A(q,r)$ fulfills all the conditions
\eqref{eq.intArq1}\eqref{eq.intArq0}\eqref{eq.intArq2}\eqref{eq.intArq3}. We
can apply Theorem~\ref{th.mainA} and the corollaries of
Subsection~\ref{se.coro} and Subsection~\ref{se.poten}.
Note in particular that $A=jAj=A^{*}$ ($j=1$), when 
 $R'(q,r,\omega^{*},\omega)$ is symmetric in
$(\omega^{*},\omega)$\,.\\
A specific case is when a particle hitting $\partial X^{+}$ at
$x=(q,p_{1},p')$\,, $p_{1}>0$\,, jumps to $(q,-p_{1},p)$ with
probability $\varepsilon(q,p)$ and to $\mathfrak{e}$ with probability
$1-\varepsilon(q,p)$\,.
The Markov kernel is then
$$
R(q,r,\omega^{-},\omega^{+})=\varepsilon(q,p)\delta(\omega_{1}^{-}+\omega_{1}^{+})\delta(\omega'^{-}-\omega'^{+})
+(1-\varepsilon(q,p))\delta_{\mathfrak{e}}\,.
$$
We assume $1-\varepsilon(q,p)\geq \alpha>0$ for some $\alpha >0$
independent of $(q,p)\in \partial X^{+}$\,.
Then the boundary condition reads simply
\begin{equation}
  \label{eq.vareps}
\gamma_{-}\varrho(q,-p_{1},p')=\varepsilon(q,p_{1},p')
\gamma_{+}\varrho(q,p_{1},p')\quad\text{with}~0\leq \varepsilon(q,p)\leq 1-\alpha\,,
\end{equation}
or
$$
\gamma_{odd}\varrho(q,p)=\sign(p_{1})\frac{1-\varepsilon(q,|p_{1}|,p')}{1+\varepsilon(q,|p_{1}|,p')}
\gamma_{ev}\varrho(q,p)\quad\forall
(q,p)\in \partial X\,.
$$
\begin{remark}
For boundary conditions of the form \eqref{eq.vareps},
Theorem~\ref{th.main0} applies to the case $\varepsilon\equiv 1$ which
corresponds to specular reflection (see Paragraph~\ref{se.ref})  while
Theorem~\ref{th.mainA} applies to the case \eqref{eq.vareps} with the
uniform upper bound $\varepsilon(q,p)\leq 1-\alpha$ with
$\alpha>0$\,. We are unfortunately not able to treat the general case when
$\varepsilon(q,p)\leq 1$\,. A norm smaller than $1$ for the contraction
 $C^{L}_{\pm}(\lambda)$ in Proposition~\ref{pr.imfermeecontr}
would suffice to handle the more general case $\varepsilon(q,p)\leq
1$\,. A result that we are not able to achieve in the general abstract
setting of Section~\ref{se.insvar}. Whether taking for $L_{\pm}$ a tangential
Kramers-Fokker-Planck operator would help, is an open question.\\
Another type of boundary conditions which occur sometimes in kinetic
theory and that we are not able to treat with Theorem~\ref{th.main0}
and Theorem~\ref{th.mainA} is the case of bounce-back boundary
conditions. It corresponds to the case when a particle hitting the
boundary $\partial X^{+}$ is sent back with an opposite velocity:
$$
\gamma_{-}\varrho(q,-p)=\gamma_{+}\varrho(q,p)\quad\text{when}~(q,p)\in \partial X^{+}\,.
$$
 The Markov kernel is then
$$
R(q,r,\omega^{-},\omega^{+})=\delta(\omega^{-}+\omega^{+})
$$
and it does not lead to a strict contraction in $\Pi_{ev}L^{2}(\partial X,|p_{1}|dq'dp)$\,.
\end{remark}
\subsubsection{Change of sign}
\label{se.chan}
 Forget for a while the potential
$V$\,. Take $V=0$ so that the Witten Laplacian $\Delta_{V}^{(0)}$ 
is the opposite standard Laplacian $-\Delta$ and consider 
a half-space problem $Q=\rz^{d}_{-}=\left\{(q^{1},q')\,,\,
  q^{1}<0\right\}$\,.  In Paragraph~\ref{se.ref}
(resp. Paragraph~\ref{se.absor})
 the Neumann (resp. Dirichlet) boundary
conditions for $-\Delta$ where interpreted as a reflecting
(resp. absorbing) boundary $\partial Q$ for the stochastic process and
this led us to natural boundary conditions for the
Kramers-Fokker-Planck equation. 
Here is another phase-space version of the Dirichlet case.
By embedding the half-space problem into a whole-space problem with
the symmetry $q^{1}\to -q^{1}$\,, the Neumann (resp. Dirichlet) boundary
conditions can be introduced by considering even (resp. odd) elements
of $L^{2}(\rz^{d},dq)$\,. More precisely the boundary value problem
$$
(1-\Delta)u=f\quad\text{with}~\partial_{q^{1}}u\big|_{\partial
  Q}=0~(\text{resp.~}u\big|_{\partial Q}=0)\,,
$$
for $f\in L^{2}(Q,dq)$ is equivalent to 
$$
(1-\Delta)\tilde{u}=\tilde{f}
$$
where $\tilde{u}, \tilde{f}\in L^{2}(\rz^{d},dq)$ are the even
(resp. odd) extensions of $u$ and $f$\,.\\ 
In the phase-space $\rz^{2d}$ with 
$X=T^{*}Q=\rz^{2d}_{-}=\left\{(q^{1},q',p_{1},p')\,,\,
  q^{1}<0\right\}$\,, the symmetry which preserves the
Kramers-Fokker-Planck operator is
$$
(q^{1},q',p_{1},p')\to (-q^{1},p',-p_{1},p')
$$
and the even (resp. odd) extension of $u\in L^{2}(\rz^{2d}_{-},dqdp)$
is given by the operator $\Sigma$ of Definition~\ref{de.sigma}
$$
\Sigma u(q^{1},q',p_{1},p')=\left\{
  \begin{array}[c]{ll}
    u(q^{1},q',p_{1},p')&\text{if}~q^{1}<0\,,\\
   ju(-q^{1},q',-p_{1},p')&\text{if}~q^{1}>0\,,\\
  \end{array}
\right.
$$
with $j=1$ (resp. $j=-1$)\,.
A solution to $P_{+,Q}u=f$ with the specular reflection boundary
condition is equivalent to $P_{+,\rz^{d}}(\Sigma u)=\Sigma f$ with
$j=1$\,. When we take $j=-1$ the boundary condition
$$
\gamma_{odd}\varrho=0
$$
is equivalent to 
\begin{equation}
  \label{eq.chsign}
\gamma_{-}\varrho(q,-p_{1},p')=-\gamma_{+}(q,p_{1},p')\quad,\quad
\forall (q,p)\in \partial X^{+}\,.
\end{equation}
Now when $Q\subset\rz^{d}$ is a bounded regular domain of $\rz^{d}$\,,
Theorem~\ref{th.main0} with $\mathfrak{f}=\cz$ and $j=-1$ applies, as
well at its corollaries of Subsection~\ref{se.coro} and
Subsection~\ref{se.PTsymm}. Adding a potential $V$ is treated in
Subsection~\ref{se.poten}.\\
We expect again that in the large friction limit (add the parameter
$\nu>0$ and let $\nu\to \infty$), the solution of the
Kramers-Fokker-Planck equation converges to
$f(q,t)e^{-\frac{|p|^{2}}{2}}$ with
$f(t)=e^{-t\Delta_{V}^{(0),D}}f_{0}$\,. The good point of the boundary
condition \eqref{eq.chsign} is that the density
$\int_{T^{*}_{q}Q}\varrho(q,p)~dp$ vanishes for all $q\in\partial
Q$\,. The drawback is that the Kramers-Fokker-Planck equation with
such boundary conditions does not preserve the positivity. There is no
straightforward interpretation in terms of stochastic processes.
\subsection{Hypoelliptic Laplacian}
\label{se.hypolapl}
In a series of works J.M.~Bismut (see~for example
\cite{Bis05}\cite{BisSurv}\cite{BiLe}) introduced and analyzed what he
called the hypoelliptic Laplacian. It is the phase-space version of
Witten's deformation of Hodge theory and it leads to an hypoelliptic
non self-adjoint second order operator. The main part of it is
actually the scalar Kramers-Fokker-Planck operator. As there are
natural boundary conditions for Witten or Hodge Laplacians, which
extends the scalar Dirichlet an Neumann boundary conditions to
$p$-forms and correspond respectively to the relative homology and
absolute homology (see
\cite{ChLi}\cite{HeNi}\cite{Lau}\cite{Lep3}\cite{LNV}), we propose a
phase-space version for the hypoelliptic Laplacian. After recalling
the writing of the Witten Laplacian and of the hypoelliptic Laplacian,
a simple half-space problem will first be considered with symmetry
arguments like in Paragraph~\ref{se.chan}. We finally propose a
general version of these boundary conditions and check that they enter
in our formalism. Especially in this section, we will see the interest of
a fiber bundle presentation of Subsection~\ref{se.connec} 
with fiber $\mathfrak{f}$ and a general
involution $j$ (or $\mathbf{j}$)\,.

\subsubsection{Witten Laplacian and Bismut's hypoelliptic Laplacian}
\label{se.witbis}

We consider here the case when $\overline{Q}=Q$ is a riemannian manifold
without boundary\,. The exterior fiber bundle  is denoted by
$\bigwedge T^{*}Q=\oplus_{p=0}^{d}\bigwedge^{p}T^{*}Q$ and we shall
consider its flat complexified  version $(\bigwedge
T^{*}Q)\otimes_{Q} (Q\times\cz)=\sqcup_{q\in Q}(\bigwedge T^{*}_{q}Q)\otimes
\cz $\,. We shall use the shorter notation $\bigwedge T^{*}Q\otimes
\cz$ when there is no ambiguity.
The set of $\mathcal{C}^{\infty}_{0}$
differential forms is $\mathcal{C}^{\infty}_{0}(Q;\bigwedge
T^{*}Q\otimes \cz)=\bigoplus_{p=0}^{d}\mathcal{C}^{\infty}_{0}(Q;\bigwedge^{p}
T^{*}Q\otimes_{Q}(Q\times\cz))$ (with $\mathcal{C}^{\infty}_{0}=\mathcal{C}^{\infty}$ when
$\overline{Q}=Q$ is compact). The metric is denoted by
$g=g_{ij}(q)dq^{i}dq^{j}$ and it provides a natural Hermitian bundle
structure on $\bigwedge T^{*}Q\otimes\cz$\,. When $d\text{Vol}_{g}(q)$ is the
riemannian volume, the $L^{2}$-scalar product of two differential
forms equals
$$
\langle \omega,\eta\rangle=\int_{Q}\langle \omega(q)\,,\, \eta(q)\rangle_{g(q)}~d\text{Vol}_{g}(q)\,.
$$
The differential acting on $\mathcal{C}^{\infty}_{0}(Q;\bigwedge
T^{*}Q\otimes \cz))$ is denoted by $d$: For $\omega=\sum_{\sharp
  I=p}\omega_{I}(q)dq^{I}\in \mathcal{C}^{\infty}_{0}(Q;\bigwedge^{p}
T^{*}Q\otimes\cz)$
$$
d\left(\sum_{\sharp I=p}\omega_{I}(q)dq^{I}\right)=\sum_{i=1}^{d}\sum_{\sharp
  I=p}\partial_{q^{i}}\omega_{I}(q)~dq^{i}\wedge dq^{I}~\in \mathcal{C}^{\infty}_{0}(Q;\bigwedge^{p+1}
T^{*}Q\otimes \cz)\,.
$$
The codifferential $d^{*}$ is its formal adjoint for the above
$L^{2}$-scalar product.
For $V\in \mathcal{C}^{\infty}(Q;\rz)$\,, Witten's deformations of the
differential and codifferential (see \cite{CFKS}\cite{Wit}\cite{Zha}) are respectively given by
$$
d_{V}=e^{-V}de^{V}=d+dV\wedge\quad,\quad d_{V}^{*}=e^{V}d^{*}e^{-V}=d^{*}+\mathbf{i}_{\nabla_{V}}\,.
$$
Owing to $d\circ d=0$\,, they also satisfy
$$
d_{V}\circ d_{V}=0\quad,\quad d_{V}^{*}\circ d_{V}^{*}=0\,.
$$
The Witten Laplacian equals
$$
\Delta_{V}=(d_{V}+d_{V}^{*})^{2}=d_{V}^{*}d_{V}+d_{V}d_{V}^{*}=\bigoplus_{p=0}^{d}\Delta_{V}^{(p)}
$$
with
$$
\Delta_{V}^{(0)}=-\Delta_{q}+\left|\nabla V(q)\right|^{2}-\Delta_{V}(q)
$$
and more generally
$$
\Delta_{V}=(d+d^{*})^{2}+|\nabla V(q)|^{2}+\mathcal{L}_{\nabla
  V}+\mathcal{L}_{\nabla V(q)}^{*}\,,
$$
where $\mathcal{L}_{X}$ is the Lie derivative along the vector field
$X$\,. When $V=0$\,, the Witten Laplacian is nothing but the Hodge
Laplacian. An important property, which requires a specific attention
when $\partial
Q\neq \emptyset$ and boundary
conditions are added, is
$$
\Delta_{V}\circ d_{V}=d_{V}\circ \Delta_{V}\quad,\quad \Delta_{V}\circ
d_{V}^{*}=d_{V}^{*}\circ \Delta_{V}\,.
$$
Bismut's hypoelliptic Laplacian is constructed like Witten's Laplacian, now on
the phase-space $X=T^{*}Q$\,. The differential $d^{X}$ is defined as usual on
$\mathcal{C}_{0}^{\infty}(X;\bigwedge T^{*}X\otimes\cz)$\,, with
$\bigwedge T^{*}X\otimes \cz=(\bigwedge T^{*}X)\otimes_{X}(X\times \cz)$\,. But the
codifferential $\overline{d}_{\phi_{b}}^{X}$ is now defined with respect to a non
degenerate but non hermitian sesquilinear form and the weight $e^{V(q)}$
has to be replaced by $e^{\mathcal{H}(x)}=e^{\mathcal{H}(q,p)}$\,,
where $\mathcal{H}(q,p)$ is some energy functional on the phase-space.
The sesquilinear form is given by
$$
\langle s\,,\, s'\rangle_{\phi_{b}}=\int_{X}\langle s(x)\,,\,s'(x)\rangle_{\eta_{b}^{*}}~dqdp\,,
$$
where $dqdp$ is the symplectic volume on $X$ and $\eta_{b}^{*}$ is the
dual form, extended to $\bigwedge T^{*}X$  and a trivial hermitian
version on $\bigwedge T^{*}X\otimes\cz$, of
$$
\eta_{b}(U,V)=\langle
\pi_{*}U\,,\,\pi_{*}V\rangle_{g}+b\omega(U,V)\,,\quad U, V\in TX=T(T^{*}Q)\,.
$$
In the above definition $\pi:X=T^{*}Q\to Q$ is the natural projection
and $\omega=dp_{j}\wedge dq^{j}$ is the symplectic form (an element
$p\in T^{*}_{q}Q$ is written $p=p_{j}dq^{j}$). 
The mapping $\phi_{b}:TX\to T^{*}X$ is given by $\eta_{b}(U,V)=\langle
U\,,\, \phi_{b}V\rangle$\,. For a section $s\in
\mathcal{C}^{\infty}(X;\bigwedge T^{*}X  \otimes\cz)$\,, 
$\overline{d}^{X}_{\phi_{b}} s$ is then defined by
$$
\forall s'\in
\mathcal{C}^{\infty}_{0}(X;T^{*}X\otimes \cz)\,,\quad
\langle s\,,\,d^{X}s'\rangle_{\phi_{b}}=\langle \overline{d}^{X}_{\phi_{b}}s\,,\,s'\rangle_{\phi_{b}}\,.
$$
The function
$\mathcal{H}(q,p)=\frac{1}{2}|p|_{q}^{2}+V(q)=\mathcal{E}(q,p)+V(q)$
and the corresponding Hamiltonian vector field on $X=T^{*}Q$ endowed
with the symplectic form $\omega$ is 
$$
\mathcal{Y}_{\mathcal{H}}=\mathcal{Y}_{\mathcal{E}}+\mathcal{Y}_{V}=g^{ij}(q)p_{i}e_{j}-\partial_{q^{j}}V(q)\partial_{p_{j}}\,,
$$
after introducing the vector field
$e_{j}=\partial_{q^{j}}+\Gamma^{\ell}_{kj}p_{\ell}\partial_{p_{k}}$
(see Subsection~\ref{se.notgFKP} and Remark~\ref{re.sign}).\\
The deformed differential and codifferential are  then given by
$$
d_{\mathcal{H}}^{X}=e^{-\mathcal{H}}d^{X}e^{\mathcal{H}}\quad\text{and}\quad
\overline{d}^{X}_{\phi_{b},\mathcal{H}}=e^{\mathcal{H}}\overline{d}^{X}_{\phi_{b}}e^{-\mathcal{H}}\,,
$$
and the hypoelliptic Laplacian by the square
\begin{equation}
  \label{eq.defhyplapl}
\mathcal{U}_{\phi_{b},\mathcal{H}}^{2}=\frac{1}{4}(\overline{d}_{\phi_{b},\mathcal{H}}^{X}+d^{X}_{\mathcal{H}})^{2}=\frac{1}{4}(\overline{d}_{\phi_{b},\mathcal{H}}^{X}\circ
d_{\mathcal{H}}^{X}+d_{\mathcal{H}}^{X}\circ\overline{d}_{\phi_{b},\mathcal{H}}^{X})\,.
\end{equation}
The Weitzenbock type formula of \cite{Bis05} expresses
$\mathcal{U}_{\phi_{b},\mathcal{H}}^{2}$
in a coordinate system or as a sum of elementary geometric operators.
It relies on the identification
$T_{x}T^{*}X\sim T_{q}Q\oplus
T_{q}^{*}Q$ related with the adjoint Levi-Civita connection on
$X=T^{*}Q$ associated with the riemannian metric $g$\,. 
For $x=(q,p)\in X=T^{*}Q$\,, the adjoint Levi-Civita connection
on $T^{*}Q$ provides the  vertical-horizontal decomposition $T_{x}X=(T_{x}X)^{H}\oplus
(T_{x}X)^{V}$ where $(T_{x}X)^{V}\sim T^{*}_{q}Q$ is spanned by 
$(\hat{e}^{j}=\partial_{p_{j}})_{j=1,\ldots,d}$ and $(T_{x}X)^{H}\sim
T_{q}Q$ is spanned by
$(e_{j}=\partial_{q^{j}}+\Gamma^{\ell}_{ij}p_{\ell}\partial_{p_{i}})$\,.
The horizontal tangent vector $\partial_{q^{j}}+\Gamma^{\ell}_{ij}p_{\ell}\partial_{p_{i}}$ is
sent to $\partial_{q^{j}}$ by the isomorphism $\pi_{*}\big|_{(T_{x}X)^{H}}:(T_{x}X)^{H}\to
T_{q}Q$\,.
The dual basis of $(e_{j},\hat{e}^{j})_{j=1,\ldots, d}$ is denoted by
$(e^{j},\hat{e}_{j})_{j=1,\ldots, d}$ with $e^{j}=dq^{j}$ and
$\hat{e}_{j}=dp_{j}-\Gamma^{\ell}_{ij}p_{\ell}dq^{j}$\,. Note that
$(e^{j})_{j=1,\ldots, d}$  (resp. $(\hat{e}_{j})_{j=1,\ldots,d}$) is
nothing but another copy of $(\hat{e}^{j})_{j=1,\ldots,d}$ (resp. of
$(e_{j})_{j=1,\ldots,d}$) after
identifying $(T^{*}_{x}X)^{H}$ with $T^{*}_{q}Q$
(resp. $(T^{*}_{x}X)^{V}$ with $T_{q}Q$).
Hence $\bigwedge^{1}T^{*}_{x}X\sim \bigwedge^{1}T_{q}^{*}Q\oplus
\bigwedge^{1}T_{q}Q$ and the exterior algebra $\bigwedge
T^{*}_{x}X$ is
identified with $(\bigwedge T_{q}^{*}Q)\otimes(\bigwedge
T_{q}Q)$\,. Some details are given below.\\
According to Theorem~3.7 of \cite{Bis05} with $\beta=\frac{1}{b}$ (see also
Theorem~3.8 of \cite{Bis05} and Theorem~3.7 in \cite{BisSurv} for
some simplifications), the hypoelliptic Laplacian
equals
\begin{align}
\nonumber
\mathcal{U}_{b,\mathcal{H}}^{2}=&
\frac{1}{4b^{2}}\big[-\Delta_{p}+|p|^{2}+2(\hat{e}_{i}\wedge)\mathbf{i}_{\hat{e}^{i}}-d
\\
  \label{eq.weitzhyp}
&\qquad-\frac{1}{2}\langle R^{TX}(e_{i},e_{j})e_{k}\,,\, e_{\ell}\rangle
(e^{i}\wedge)(e^{j}\wedge)\mathbf{i}_{\hat{e}^{k}}\mathbf{i}_{\hat{e}^{\ell}}\big]
-\frac{1}{2b}\mathcal{L}_{\mathcal{Y}^{\mathcal{H}}}\,,
\end{align}
where $R^{TX}$ is the Riemann curvature tensor on $TQ$\,.
\begin{remark}
\label{re.weitz}
More precisely the Weitzenbock formula \eqref{eq.weitzhyp} is really
written in this way in Theorem~3.7 and Theorem~3.8 of \cite{Bis05} 
for $V\equiv 0$ (take $\omega(F,\nabla F)=0$ with the notations of
\cite{Bis05}). The way to include a non zero potential $V(q)$ is
explained in Remark~2.37  of \cite{Bis05}: It suffices to endow the
complexification bundle $F_{1}=Q\times\cz$  with the
metric $g^{F_{1}}=e^{-2V(q)}$ then to apply the Weitzenbock formula of
Theorem~3.7 and Theorem~3.8 in \cite{Bis05} with
$\omega(\nabla^{F_{1}},g^{F_{1}})$ (denoted by $\omega(\nabla^{F},g^{F})$ in
\cite{Bis05}) equal to $-2\nabla_{q}V(q)$ and finally to compute
$e^{-V(q)}\mathcal{U}_{b,\mathcal{E}}^{2}e^{V(q)}$ from the
general formula for $\mathcal{U}_{b,\mathcal{E}}$\,.
\end{remark}
The hypoelliptic Laplacian $\mathcal{U}_{b,\mathcal{H}}$ is not yet in
the form of a geometric Kramers-Fokker-Planck operator (see
Definition~\ref{de.genGKFP}). For this we need some precisions about
the identification $\bigwedge T^{*}X\sim (\bigwedge
T^{*}Q)\otimes_{Q}(\bigwedge TQ)$\,. These details may also help the
neophytes to grasp the more involved formalism of \cite{Bis05}.
Remember that the
$\hat{e}^{j}=\partial_{p_{j}}$ and
$e_{j}=\partial_{q^{j}}+\Gamma^{\ell}_{ij}p_{\ell}\partial_{p_{i}}$
form a basis of $T_{x}X$ while  the $e^{j}=dq^{j}$ and
$\hat{e}_{j}=dp_{j}-\Gamma^{\ell}_{ij}p_{\ell}dq^{i}$ form a basis of
$\bigwedge^{1}T_{x}^{*}X$\,. When we work on $Q$\,,
$(e_{j}=\partial_{q^{j}})_{1\leq j\leq}$ is a basis of $T_{q}Q$\,,
$(e^{j}=dq^{j})_{1\leq j\leq d}$ a basis of $T^{*}_{q}Q$ and
$(\hat{e}_{j}=dp_{j})_{1\leq j\leq d}$ another copy of
$(\partial_{q^{j}})_{1\leq j\leq d}$\,.
The fiber bundle $F=
(\bigwedge
T^{*}Q)\otimes_{Q}(\bigwedge TQ)\otimes_{Q}(Q\times \cz)$  is endowed with
the hermitian metric $g^{F}=(\bigwedge g^{-1})\otimes (\bigwedge g)$ 
and the Levi-Civita connection characterized by
$$
\nabla^{F}_{\partial
  q^{i}}(dq^{j})=-\Gamma^{j}_{i\ell}dq^{\ell}\quad
\text{and}\quad\nabla^{F}_{\partial_{q^{i}}}(dp_{j})=\Gamma^{\ell}_{ij}dp_{\ell}\,.
$$
After identifying $(T_{x}X)^{H}$ with $T_{q}Q$ and $(T_{x}X)^{V}$ with $T^{*}_{q}Q$
for $x=(q,p)$\,, the complexified Grassman bundle $\bigwedge
T^{*}X\otimes_{X}(X\times\cz)$ is nothing but $F_{X}=\pi^{*}F$
introduced in Subsection~\ref{se.connec} with the natural projection
$\pi:X=T^{*}Q\to Q$\,.
A basis of $\bigwedge T^{*}_{x}X\otimes
\cz$ (or $F$) is given by
$e^{I}\hat{e}_{J}=e^{i_{1}}\wedge\ldots\wedge e^{i_{p}}\wedge
\hat{e}_{j_{1}}\wedge\ldots \wedge \hat{e}_{j_{p'}}$ with
$I=\left\{i_{1},\ldots,i_{p}\right\}$ and
$J=\left\{j_{1},\ldots,j_{p'}\right\}$\,, the $i_{k}$'s and $j_{k}$'s
being written in the increasing order. The hermitian metric
$g^{F_{X}}=\pi^{*}g^{F}$ is extended to the exterior algebra from
$$
g^{F_{X}}(x;e^{i},e^{j})=g^{ij}(q)\quad,\quad
g^{F_{X}}(x;\hat{e}_{i},\hat{e}_{j})=g_{ij}(q)\quad,\quad
g^{F_{X}}(x;e^{i},\hat{e}_{j})=0\,,
$$
where $e^{j}=dq^{j}$ and
$\hat{e}_{j}=dp_{j}-\Gamma^{\ell}_{ij}p_{\ell}dq^{i}$ at
$x=(q,p)$\,.\\
According to \eqref{eq.connFX} the connection $\nabla^{FX}$ is characterized by 
\begin{eqnarray*}
  &&
  \nabla^{F_{X}}_{e_{i}}\hat{e}_{j}=\Gamma^{\ell}_{ij}\hat{e}_{\ell}\quad,\quad
  \nabla^{F_{X}}_{e_{i}}e^{j}=-\Gamma^{j}_{i\ell}e^{\ell}\,,
\\
\text{and}
&& \nabla^{F_{X}}_{\hat{e}^{i}}\hat{e}_{j}=0\quad,\quad \nabla_{\hat{e}^{i}}e^{j}=0\,,
\end{eqnarray*}
and it is compatible with the metric $g^{F_{X}}$ because $\nabla^{F}$
is compatible with $g^{F}$\,.\\

We take $b=\mp 1$ and we note that the hypoelliptic Laplacian
\begin{multline}
\label{eq.weitz1}
2\mathcal{U}_{\mp 1,\mathcal{H}}^{2}=
\frac{1}{2}\big[-\Delta_{p}+|p|^{2}+2(\hat{e}_{i}\wedge)\mathbf{i}_{\hat{e}^{i}}-\dim Q
\\
\qquad-\frac{1}{2}\langle R^{TX}(e_{i},e_{j})e_{k}\,,\, e_{\ell}\rangle
(e^{i}\wedge)(e^{j}\wedge)\mathbf{i}_{\hat{e}^{k}}\mathbf{i}_{\hat{e}^{\ell}}\big]
\pm \mathcal{L}_{\mathcal{Y}^{\mathcal{H}}}\,,
\end{multline}
and we follow the presentation of \cite{Leb1} in order to write it in
the form of Definition~\ref{de.genGKFP}. The Lie derivative
$\mathcal{L}_{\mathcal{Y}^{\mathcal{H}}}$ applied to the form
$\omega=\omega_{I}^{J}e^{I}\wedge \hat{e}_{J}$ gives
\begin{equation}
  \label{eq.liesec}
\mathcal{L}_{\mathcal{Y}^{\mathcal{H}}}\omega
=(\mathcal{Y}^{\mathcal{H}}\omega_{I}^{J})e^{I}\wedge\hat{e}_{J}
+\omega_{I}^{J} \mathcal{L}_{\mathcal{Y}^{\mathcal{H}}}(e^{I}\wedge \hat{e}_{J})\,.
\end{equation}
We compute
\begin{eqnarray*}
  &&\mathcal{L}_{\mathcal{Y}^{\mathcal{H}}}(dq^{j})=d(\mathcal{Y}^{\mathcal{H}}q^{j})=
  d(\partial_{p_{j}}\mathcal{H})=(\partial^{2}_{q^{k}p_{j}}\mathcal{H})dq^{k}
+(\partial^{2}_{p_{k}p_{j}}\mathcal{H})
dp_{k}\,,\\
&& \mathcal{L}_{\mathcal{Y}^{\mathcal{H}}}(dp_{j})=d(\mathcal{Y}^{\mathcal{H}}p_{j})=-d(\partial_{q^{j}}\mathcal{H})=-(\partial^{2}_{q^{k}q^{j}}\mathcal{H})dq^{k}-(\partial^{2}_{p_{k}q^{j}}\mathcal{H})dp_{k}\,.
\end{eqnarray*}
For $\mathcal{H}=\mathcal{E}$\,, let us compute 
$$
\mathcal{L}_{\mathcal{Y}^{\mathcal{E}}}(e^{j})-\nabla^{F_{X}}_{\mathcal{Y}^{\mathcal{E}}}e^{j}\quad\text{and}\quad
\mathcal{L}_{\mathcal{Y}^{\mathcal{E}}}(\hat{e}_{j})-\nabla^{F_{X}}_{\mathcal{Y}^{\mathcal{E}}}\hat{e}_{j}\,.
$$
From $e^{j}=dq^{j}$\,, $\hat{e}_{j}
=dp_{j}-\Gamma^{i}_{kj}p_{i}dq^{k}$\,, $\Gamma_{ij}^{k}=\Gamma_{ji}^{k}$ (the
Levi-Civita is torsion free) and 
$$
\partial_{q^{\ell}}g^{ij}=-g^{jk}\Gamma^{i}_{\ell
  k}-g^{i\ell}\Gamma^{j}_{\ell k}\,,
$$
we deduce 
\begin{align}
\nonumber
\mathcal{L}_{\mathcal{Y}^{\mathcal{E}}}e^{j}-\nabla^{F_{X}}_{\mathcal{Y}^{\mathcal{E}}}e^{j}
&= (\partial_{q^{k}}g^{ij})p_{i}dq^{k}+g^{jk} dp_{k}
-g^{ik}p_{i}\nabla^{F_{X}}_{e_{k}}e^{j}\\
\nonumber
&=(\partial_{q^{k}}g^{ij})p_{i}e^{k}+g^{jk}dp_{k}+g^{ik}p_{i}\Gamma^{j}_{k\ell}e^{\ell}\\
\label{eq.diffej}
&= g^{jk}\hat{e}_{k}\,.
\end{align}
The computation of
$\mathcal{L}_{\mathcal{Y}^{\mathcal{E}}}\hat{e}_{j}-\nabla^{F_{X}}_{\mathcal{Y}^{\mathcal{E}}}\hat{e}_{j}$
with $\hat{e}_{j}=dp_{j}-\Gamma^{i}_{jk}p_{i}dq^{k}$ relying
on the same argument is a bit more involved:
\begin{align}
\nonumber
  \mathcal{L}_{\mathcal{Y}^{\mathcal{E}}}\hat{e}_{j}-\nabla^{F_{X}}_{\mathcal{Y}_{\mathcal{E}}}\hat{e}_{j}
\nonumber
=&
-(\partial^{2}_{q^{j}q^{k}}g^{i\ell})p_{i}p_{\ell}dq^{k}-(\partial_{q^{j}}g^{ik})p_{i}dp_{k}\\
\nonumber
&\quad
-\Gamma_{jm}^{n}p_{n}(\partial_{q^{k}}g^{im})p_{i}dq^{k}-\Gamma_{jm}^{n}p_{n}g^{mk}dp_{k}
-g^{i\ell}p_{i}\nabla^{F_{X}}_{e_{\ell}}\hat{e}_{j}
\\
\nonumber
=&\left((\partial^{2}_{q^{j}q^{k}}g^{i\ell})+\Gamma^{\ell}_{jm}(\partial_{q^{k}}g^{im})\right)p_{i}p_{\ell}dq^{k}\\
\nonumber
&\quad
+g^{im}\Gamma^{k}_{mj}p_{i}dp_{k}-g^{i\ell}p_{i}\Gamma^{k}_{\ell
  j}\hat{e}_{k}
\\
\nonumber
&=\left((\partial^{2}_{q^{j}q^{k}}g^{i\ell})+\Gamma^{\ell}_{jm}(\partial_{q^{k}}g^{im})
\right)p_{i}p_{\ell}dq^{k}\\
\nonumber&\quad + g^{i\ell}p_{i}\Gamma^{k}_{\ell j}\Gamma^{n}_{k
  n'}p_{n}dq^{n'}\\
\label{eq.diffhatej}
&=\left((\partial^{2}_{q^{j}q^{k}}g^{i\ell})+\Gamma^{\ell}_{jm}(\partial_{q^{k}}g^{im})+
g^{in}\Gamma^{m}_{nj}\Gamma^{\ell}_{mk}\right)p_{i}p_{\ell} e^{k}\,,
\end{align}
where the last factor of $e^{k}$ can be related with the Riemann
curvature tensor.
The computation of
$\mathcal{L}_{\mathcal{Y}^{V}}-\nabla^{F_{X}}_{\mathcal{Y}^{V}}$ is even
simpler because
$\mathcal{Y}^{\mathcal{V}}=-\partial_{q^{j}}V(q)\partial_{p_{j}}$ and
$\nabla^{F_{X}}_{\mathcal{Y}^{V}}=0$ while
$$
 \mathcal{L}_{\mathcal{Y}^{V}}(dq^{j})=0
\quad\text{and}\quad
\mathcal{L}_{\mathcal{Y}^{V}}(dp_{j})=
-(\partial^{2}_{q^{j}q^{k}}V)dq^{k}
$$
imply
\begin{eqnarray}
\label{eq.lievej}
  && (\mathcal{L}_{\mathcal{Y}^{V}}-\nabla^{F_{X}}_{\mathcal{Y}^{V}})(e^{j})=0\\
\label{eq.lievhej}
\text{and}&&
(\mathcal{L}_{\mathcal{Y}^{V}}-\nabla^{F_{X}}_{\mathcal{Y}^{V}})(\hat{e}_{j})
=\left(-(\partial^{2}_{q^{j}q^{k}}V)+\Gamma^{\ell}_{jk}(\partial_{q^{\ell}}V)\right)e^{k}\,.
\end{eqnarray}
The relations \eqref{eq.liesec}\eqref{eq.diffej}\eqref{eq.diffhatej}\eqref{eq.lievej}\eqref{eq.liesec},
with
$\mathcal{Y}^{\mathcal{H}}=g^{ij}p_{i}\partial_{q^{j}}-(\partial_{q^{j}}V)\partial_{p_{j}}$
for $\mathcal{H}(q,p)=\frac{|p|^{2}_{q}}{2}+V(q)$\,,
specify the action of the derivation
$\mathcal{L}_{\mathcal{Y}^{\mathcal{H}}}-\nabla_{\mathcal{Y}^{\mathcal{H}}}$
on $\mathcal{C}^{\infty}(X;\bigwedge
T^{*}X\otimes\cz)=\mathcal{C}^{\infty}(X;\pi^{*}F)$\,. 
If one keeps the metric $g^{F_{X}}$ for the definition of
$L^{2}$-sections, it is not a geometric Kramers-Fokker-Planck operator
because the right-hand side \eqref{eq.diffhatej} is homogeneous of
degre $2$ w.r.t $p$ and cannot be included in the remainder term
$M(q,p,\partial_{p})$ of Definition~\ref{de.genGKFP}\,. This is corrected by
considering the $|p|_{q}$-homogeneity of
$e^{j}=dq^{j}$ with degree $0$ and of
$\hat{e}_{j}=dp_{j}-\Gamma^{i}_{jk}p_{i}dq^{k}$ with degree $1$\,. 
 Consider the
weighted metric $\langle p\rangle^{2\widehat{\text{deg}}}g^{F_{X}}$
given by 
$$
[
\langle p\rangle^{2\widehat{\text{deg}}}
g^{F_{X}}](x; e^{I}\hat{e}_{J}, e^{I'}\hat{e}_{J'})=
\langle p\rangle_{q}^{2|J|}g^{F_{X}}(x; e^{I}\hat{e}_{J}, e^{I'}\hat{e}_{J'})\,,
$$
which is non zero iff $|I|=|I'|$ and $|J|=|J'|$\,.
The local operator $\langle p\rangle^{\pm \widehat{\text{deg}}}$ is
defined by $\langle
p\rangle^{\pm\widehat{\text{deg}}}e^{I}\hat{e}_{J}=\langle
p\rangle_{q}^{\pm|J|}e^{I}\hat{e}_{J}$ and $\langle
p\rangle_{q}^{-\widehat{\text{deg}}}$ is a unitary operator from
$L^{2}(X,dqdp;(F_{X},g^{F_{X}}))$ onto $L^{2}(X,dqdp;(F_{X},\langle
p\rangle^{2\widehat{\text{deg}}}g^{F_{X}}))$\,. Thus studying the
hypoelliptic Laplacian $2\mathcal{U}_{\mp 1,\mathcal{H}}^{2}$ in
$L^{2}(X,dqdp;(F_{X},\langle p\rangle^{2\widehat{\text{deg}}}g^{F_{X}}))$ is
the same as studying
\begin{equation}
  \label{eq.hypoweight}
\mathcal{K}_{\pm,\mathcal{H}}=\langle p\rangle^{+\widehat{\text{deg}}}\circ
2\mathcal{U}_{\mp 1,\mathcal{H}}^{2}\circ \langle p\rangle^{-\widehat{\text{deg}}}
\end{equation}
in $L^{2}(X,dqdp;(F_{X},g^{F_{X}}))$\,. 
After the conjugation with $\langle p\rangle^{\widehat{\text{deg}}}$
all the right-hand sides of
\eqref{eq.diffej}\eqref{eq.diffhatej}\eqref{eq.lievej}\eqref{eq.lievhej}
take the form of the remainder term $M(q,p,\partial_{p})$ (remember
that
$\mathcal{L}_{\mathcal{Y}^{H}}-\nabla^{F_{X}}_{\mathcal{Y}^{\mathcal{H}}}$
is a derivation).\\
Hence $\mathcal{K}_{\pm,\mathcal{H}}$ is a geometric
Kramers-Fokker-Planck operator according to
Definition~\ref{de.genGKFP}.
General boundary problems for such operators have been studied in Subsection~\ref{se.connec}.

\subsubsection{Review of natural boundary Witten Laplacians}
\label{se.wittbc}
Consider now a compact manifold with boundary
$\overline{Q}=Q\sqcup \partial Q$\,.
The so-called Dirichlet and Neumann boundary conditions for
the
Witten Laplacian acting on $p$-forms correspond after de~Rham
duality to the determination of relative homology and
absolute homology groups respectively (see
\cite{ChLi}\cite{HeNi}\cite{Lep3}
\cite{Lau}).\\
After introducing coordinates for which the metric $g$ has the form
$g=dq^{1}+m_{ij}(q^{1},q')dq'^{i}dq'^{j}$ according to
\eqref{eq.defgm} along the boundary, a general $p$-form can be written
$$
\omega=\sum_{|I|=p\,,\, 1\not\in
  I}\omega_{I}(q)dq^{I}+\sum_{|I'|=p-1\,,\, 1\not\in
  I'}\omega_{\left\{1\right\}\cup I'}(q)dq^{1}\wedge dq^{I'}\,.
$$
The tangential  component and the normal component  along $\partial Q$
are respectively defined by 
\begin{eqnarray*}
  &&\mathbf{t}\omega=\sum_{|I|=p\,,\, 1\not\in
    I}\omega_{I}(0,q')dq^{I}\,,\\
&& \mathbf{n}\omega=\sum_{|I'|=p-1\,,\, 1\not\in
  I'}\omega_{\left\{1\right\}\cup I'}(0,q')dq^{1}\wedge dq^{I'}.
\end{eqnarray*}
The domain of  Neumann Witten Laplacian is 
$$
D(\Delta_{V}^{N})=\left\{u\in H^{2}(Q;\bigwedge
  T^{*}Q\otimes\cz)\,,\quad\left[
    \begin{array}[c]{l}
      \mathbf{n}\omega=\omega\\
      \mathbf{n}d_{V}\omega=0
    \end{array}
\right.\right\}\,,
$$
and the domain of the Dirichlet Laplacian is
$$
D(\Delta_{V}^{D})=
\left\{u\in H^{2}(Q;\bigwedge
  T^{*}Q\otimes\cz)\,,\quad\left[
    \begin{array}[c]{l}
      \mathbf{t}\omega=\omega\\
      \mathbf{t}d_{V}^{*}\omega=0
    \end{array}
\right.\right\}\,.
$$
Note that the potential appears only in the boundary condition which
involves the second trace, the trace of the (co-)differential,
i.e. 
\begin{eqnarray*}
  && \mathbf{n}d\omega
+dV\wedge\omega\big|_{\partial Q}=0\quad\text{for}~\Delta_{V}^{N}\,,\\
\text{and}&&
\mathbf{t}d^{*}\omega +\mathbf{i}_{\nabla f}\omega\big|_{\partial Q}=0\quad\text{for}~\Delta_{V}^{D}\,.
\end{eqnarray*}
Since the boundary conditions for the Kramers-Fokker-Planck equation
involve only the first trace, natural boundary conditions should not
depend on the potential (see~a.e. the case of the specular
reflection in Paragraph~\ref{se.ref}).

\subsubsection{Neumann  and Dirichlet boundary conditions for
  the hypoelliptic Laplacian on flat cylinders}
\label{se.neumdir}
We now work in the phase-space $\overline{X}=X\sqcup \partial
X=X\sqcup T^{*}_{\partial Q}Q$ where $\overline{Q}$ is  the half-space
$\overline{\rz^{d}}_{-}$ or possibly a half-cylinder
$(-\infty,0]\times Q'$\,. The metric is assumed to be $g=1\oplus m$
with $\partial_{q^{1}}m=0$\,. Within this framework, boundary
conditions for 
 the hypoelliptic Laplacian, which should hopefully be related with
 Neumann and Dirichlet realization of the Witten Laplacian,
 are introduced with the help of a simple
reflection principle.
 Remember that the functional analysis of the
hypoelliptic Laplacian $\mathcal{U}_{\mp 1,\mathcal{H}}^{2}$ 
given in  \eqref{eq.defhyplapl} and \eqref{eq.weitzhyp}, is made
in standard $L^{2}$-spaces by using the conjugated form
$\mathcal{K}_{\pm,\mathcal{H}}$ defined in \eqref{eq.hypoweight}.
\\
According to the previous discussion we consider boundary conditions
which do not depend on the potential $V(q)$\,. Since its contribution
in the conjugated hypoelliptic Laplacian, 
$\mathcal{K}_{\pm,\mathcal{H}}$\,, is reduced to a relatively bounded
perturbation of the case $V=0$\,, we can forget it at first.\\
The most natural introduction of 
boundary conditions for the hypoelliptic Laplacian, corresponding to 
Neumann and Dirichlet realization of the Hodge Laplacian ($V\equiv 0$), relies on
the extension by symmetry or anti-symmetry like in Paragraph~\ref{se.chan}.
In the position variable the extension by symmetry
(resp. anti-symmetry), $\sigma_{k}(\omega)(q^{1},q')$ defined for
$q^{1}\in \rz$\,, of  a form
$\omega(q^{1},q')=\omega_{I}(q^{1},q')dq^{I}$\,, initially defined for
$q^{1}<0$\,, equals
$$
\sigma_{k}(\omega_{I}dq^{I})
=\left\{
  \begin{array}[c]{ll}
    \omega(q^{1},q') &\text{if}~q^{1}<0\,,\\
    (-1)^{k}\left[(-1)^{|\left\{1\right\}\cap I|}\omega_{I}(-q^{1},q')dq^{I}\right]
&\text{if}~q^{1}>0\,,
  \end{array}
\right.
$$
with $k=0$ (resp. $k=1$)\,.
On the whole space the mapping
$\tilde{\sigma}_{k}:L^{2}(Q'\times\rz;\bigwedge
T^{*}(Q'\times\rz)\otimes\cz)$ is defined
by 
\begin{equation*}
\tilde{\sigma}_{k}\left[\omega_{I}
dq^{I}\right](q^{1},q')
=(-1)^{k}\left[(-1)^{|\left\{1\right\}\cap I|}\omega_{I}(-q^{1},q')dq^{I}
\right]\,.
\end{equation*}
A regular $\omega\in
\mathcal{C}^{\infty}_{0}(\overline{Q},\bigwedge
T^{*}Q\otimes \cz)$ satisfies the Neumann
(resp. Dirichlet) boundary conditions, $\mathbf{n}\omega=0$ and
$\mathbf{n}d\omega=0$
(resp. $\mathbf{t}\omega=0$ and $\mathbf{t}d^{*}\omega =0$) iff 
the extended form
$\tilde{\eta}=\sigma_{k}(\eta)$ 
for $\eta\in
  \left\{\omega,d\omega\right\}$\,, (resp. $\eta\in
  \left\{\omega,d^{*}\omega\right\}$) has a
  trace along  $\partial Q=\left\{0\right\}\times Q'$ (no
  discontinuity):
$$
\sigma_{k}(\eta)(0^{+},q')=\eta(0^{-},q')\,.
$$
This condition implies that the extended form
$\tilde{\omega}=\sigma_{k}\omega$ satifies in addition to 
$\tilde{\sigma}_{k}\tilde{\omega}=\tilde{\omega}$\,, the continuity
of $\tilde{\omega}$ and $d\tilde{\omega}$
(resp. $d^{*}\tilde{\omega}$) along $\partial Q$\,.\\
While keeping $\overline{Q}=Q'\times (-\infty,0]$ with
$\partial_{q^{1}}m\equiv 0$ and $V\equiv 0$\,, the phase-space
reflection extensively used in Subsection~\ref{se.halfwhL},
Subsection~\ref{se.specmet}
and Paragraph~\ref{se.chan} is given by 
$$
(q^{1},q',p_{1},p')\to (-q^{1},q',-p_{1},p')\,.
$$
The differentials $dq^{1}$ and $dp_{1}$ are pushed forward to
$-dq^{1}$ and $-dp_{1}$ by this mapping. Hence the extension
$\Sigma_{k}$\,, with the general 
Definition~\ref{de.sigma}\,, becomes here
$$
\Sigma_{k}(\omega)(q^{1},q',p_{1},p')=\omega(q^{1},q',p_{1},p')\qquad\text{if}~q^{1}<0\,,
$$
and
$$
\Sigma_{k}(\omega)(q^{1},q',p_{1},p')=
    (-1)^{k}\left[(-1)^{|\left\{1\right\}\cap
        I|+|\left\{1\right\}\cap
        J|}\hat{\omega}^{J}_{I}(-q^{1},q',-p_{1},p')dq^{I}\wedge dp_{J}\right]
$$
if $q^{1}>0$\,, when 
$$
\omega=\hat{\omega}^{J}_{I}dq^{I}\wedge dp_{J}\,.
$$
The Neumann case corresponds to $k=0$ and the Dirichlet case to
$k=1$\,. Accordingly the mapping $\tilde{\Sigma}_{k}$ is given by
$$
\tilde{\Sigma}_{k}\left[\hat{\omega}^{J}_{I}dq^{I}\wedge dp_{J}\right]=(-1)^{k}\left[(-1)^{|\left\{1\right\}\cap
        I|+|\left\{1\right\}\cap
        J|}\hat{\omega}^{J}_{I}(-q^{1},q',-p_{1},p')dq^{I}\wedge dp_{J}\right]\,.
$$
The corresponding boundary conditions for the conjugated hypoelliptic
Laplacian $\mathcal{K}_{\pm,\mathcal{H}}$\,, can be written
$$
\hat{\omega}^{J}_{I}(0,q',p_{1},p')\big|_{p_{1}<0}=
\left\{
  \begin{array}[c]{ll}
    \hat{\omega}^{J}_{I}(0,q',-p_{1},p')&\text{if}~|\left\{1\right\}\cap
    I|+|\left\{1\right\}\cap J|~\text{is~even}\,,\\
-\hat{\omega}^{J}_{I}(0,q',-p_{1},p')&\text{if}~|\left\{1\right\}\cap
    I|+|\left\{1\right\}\cap J|~\text{is~odd}\,,
  \end{array}
\right.
$$
for the Neumann boundary conditions ($k=0$), and by 
$$
\hat{\omega}^{J}_{I}(0,q',p_{1},p')\big|_{p_{1}<0}=
\left\{
  \begin{array}[c]{ll}
    \hat{\omega}^{J}_{I}(0,q',-p_{1},p')&\text{if}~|\left\{1\right\}\cap
    I|+|\left\{1\right\}\cap J|~\text{is~odd}\,,\\
-\hat{\omega}^{J}_{I}(0,q',-p_{1},p')&\text{if}~\left\{1\right\}\cap
    I|+|\left\{1\right\}\cap J|~\text{is~even}\,,
  \end{array}
\right.
$$
for the Dirichlet boundary conditions ($k=1$). Those boundary
conditions extend to the case of $p$-forms the specular boundary
condition for the scalar case reviewed in Paragraph~\ref{se.ref} and
the boundary condition with a change of sign presented in
Paragraph~\ref{se.chan}.\\
Nevertheless they are not yet written in the proper form considered in
Subsection~\ref{se.connec} for $(dq^{I}\wedge
dp_{J})_{|I|\leq d, |J|\leq d}$ is not the basis $(e^{I}\hat{e}_{J})_{|I|\leq d, |J|\leq d}$ corresponding to the identification $\bigwedge T^{*}X=F_{X}$ with
$F=(\bigwedge T^{*}Q)\otimes_{Q}(\bigwedge TQ)$\,.
Fortunately, our assumption $g=1\oplus m$ with $\partial_{q^{1}}m=0$
imply
$$
\Gamma_{1j}^{i}=\Gamma_{j1}^{i}=\Gamma^{1}_{ij}=0\,,
$$
which means $e^{1}=dq^{1}$\,, $\hat{e}_{1}=dp_{1}$ and
$\mathbf{i}_{\partial_{q^{1}}}e^{I}\hat{e}_{J}=\mathbf{i}_{\partial_{p_{1}}}e^{I}\hat{e}_{J}=0$
when $1\not\in I$ and $1\not \in J$\,. 
Therefore the forms $\Sigma_{k}(\omega)1_{\rz_{+}}(q_{1})=\Sigma_{k}
(\omega^{J}_{I}e^{I}\hat{e}_{J})1_{\rz_{+}}(q_{1})$ 
and
$\tilde{\Sigma}_{k}(\omega)=\tilde{\Sigma}_{k}(\omega^{J}_{I}e^{I}\hat{e}_{J})$
are respectively equal to 
\begin{eqnarray*}
  && (-1)^{k}\left[(-1)^{|\left\{1\right\}\cap
        I|+|\left\{1\right\}\cap
        J|}\omega^{J}_{I}(-q^{1},q',-p_{1},p')e^{I}\hat{e}_{J}\right]1_{\rz_{+}}(q^{1})\\
\text{and}
&&
(-1)^{k}\left[(-1)^{|\left\{1\right\}\cap
        I|+|\left\{1\right\}\cap
        J|}\omega^{J}_{I}(-q^{1},q',-p_{1},p')e^{I}\hat{e}_{J}\right]\,.
\end{eqnarray*}
Hence the Neumann ($k=0$) boundary conditions now written as
\begin{equation}
\label{eq.NeuHyp}
\omega^{J}_{I}(0,q',p_{1},p')\big|_{p_{1}<0}=
\left\{
  \begin{array}[c]{ll}
    \omega^{J}_{I}(0,q',-p_{1},p')&\text{if}~|\left\{1\right\}\cap
    I|+|\left\{1\right\}\cap J|~\text{is~even}\,,\\
-\omega^{J}_{I}(0,q',-p_{1},p')&\text{if}~\left\{1\right\}\cap
    I|+|\left\{1\right\}\cap J|~\text{is~odd}\,,
  \end{array}
\right.
\end{equation}
and the Dirichlet ($k=1$) boundary conditions now written as
\begin{equation}
  \label{eq.DirHyp}
\omega^{J}_{I}(0,q',p_{1},p')\big|_{p_{1}<0}=
\left\{
  \begin{array}[c]{ll}
    \omega^{J}_{I}(0,q',-p_{1},p')&\text{if}~|\left\{1\right\}\cap
    I|+|\left\{1\right\}\cap J|~\text{is~odd}\,,\\
-\omega^{J}_{I}(0,q',-p_{1},p')&\text{if}~\left\{1\right\}\cap
    I|+|\left\{1\right\}\cap J|~\text{is~even}\,,
  \end{array}
\right.
\end{equation}
correspond  to the case 
\begin{eqnarray*}
  && \mathbf{A}=0\,,\\
\text{and}&&
\mathbf{j}(e^{I}\hat{e}_{J})=(-1)^{k}(-1)^{|\left\{1\right\}\cap
  I|+|\left\{1\right\}\cap J|}e^{I}\hat{e}_{J}
\end{eqnarray*}
of Proposition~\ref{pr.bundle} with respectively $k=0$ for  and
$k=1$\,.
We shall denote $\mathcal{K}_{\pm,\mathcal{H}}^{N}$
(resp. $\mathcal{K}_{\pm,\mathcal{H}}^{D}$) the Neumann (resp. Dirichlet)
realization corresponding to $k=0$ (resp. $k=1$) of
$\mathcal{K}_{\pm,\mathcal{H}}$\,.
Proposition~\ref{pr.bundle}  applies to
$\mathcal{K}_{\pm,\mathcal{H}}^{N}$ and
$\mathcal{K}_{\pm,\mathcal{H}}^{D}$ also when
$V(q)\neq 0$ is a globally Lipschitz potential.\\
We end this section by checking the commutation with the differential
on a dense set of the domain. The hypoelliptic Laplacian
$\mathcal{U}_{\mp,\mathcal{H}}$ as a differential operator acting on
$\mathcal{C}^{\infty}(\overline{X};F_{X})$ commutes with
$d^{X}_{\mathcal{H}}$\,. Equivalently $\mathcal{K}_{\pm,\mathcal{H}}$
commutes with 
$$
\langle p\rangle^{\widehat{\text{deg}}}d^{X}_{\mathcal{H}}\langle
p\rangle^{-\widehat{\text{deg}}}
=
\langle p\rangle^{\widehat{\text{deg}}}(e^{-\mathcal{H}}d^{X}e^{\mathcal{H}})\langle
p\rangle^{-\widehat{\text{deg}}}
=
\langle p\rangle^{\widehat{\text{deg}}}(d^{X}+d\mathcal{H}\wedge)\langle
p\rangle^{-\widehat{\text{deg}}}\,.
$$
\begin{proposition}
\label{pr.commut}
Let $\mathcal{K}^{N}_{\pm,\mathcal{H}}$
(resp. $\mathcal{K}_{\pm,\mathcal{H}}^{D}$) be the Neumann
(resp. Dirichlet) realization of $\mathcal{K}_{\pm,\mathcal{H}}$ when
$\overline{Q}=(-\infty,0]\times Q'$ and $\partial_{q^{1}}m\equiv 0$
($Q'$ is either compact or a compact perturbation of the euclidean
space $\rz^{d-1}$) and
$\mathcal{H}(q,p)=\frac{|p|_{q}^{2}}{2}+V(q)$ and $\mathcal{V}$
globally Lipschitz on $\overline{Q}$\,.
Then the set $\mathcal{D}$ of sections $\omega=\omega^{J}_{I}e^{I}\hat{e}_{J}\in
\mathcal{C}^{\infty}_{0}(\overline{X};F_{X})$ which satisfy
\eqref{eq.NeuHyp} (resp. \eqref{eq.DirHyp}) and
$$
|\partial_{q^{1}}[e^{\mathcal{H}(q,p)}
\langle p\rangle_{q}^{-\widehat{\text{deg}}}\omega](q,p)|=\mathcal{O}(|q^{1}|^{\infty})
$$
is dense in $D(K_{\pm,\mathcal{H}}^{N})$
(resp. $D(K_{\pm,\mathcal{H}}^{D})$) endowed with its graph norm.\\
For any $\omega \in \mathcal{D}$\,, $\langle p\rangle^{\widehat{\text{deg}}}d^{X}_{\mathcal{H}}\langle
p\rangle^{-\widehat{\text{deg}}}\omega$ belongs to
$D(K_{\pm,\mathcal{H}}^{N})$ (resp. $D(K_{\pm,\mathcal{H}}^{D})$) and 
\begin{eqnarray}
\label{eq.commutN}
  && K_{\pm,\mathcal{H}}^{N}
\langle p\rangle^{\widehat{\text{deg}}}d^{X}_{\mathcal{H}}\langle
p\rangle^{-\widehat{\text{deg}}}\omega = \langle p\rangle^{\widehat{\text{deg}}}d^{X}_{\mathcal{H}}\langle
p\rangle^{-\widehat{\text{deg}}} \mathcal{K}_{\pm,
  \mathcal{H}}^{N}\omega\,,\\
\label{eq.commutD}
\text{resp.}&&
K_{\pm,\mathcal{H}}^{D}
\langle p\rangle^{\widehat{\text{deg}}}d^{X}_{\mathcal{H}}\langle
p\rangle^{-\widehat{\text{deg}}}\omega = \langle p\rangle^{\widehat{\text{deg}}}d^{X}_{\mathcal{H}}\langle
p\rangle^{-\widehat{\text{deg}}} \mathcal{K}_{\pm,
  \mathcal{H}}^{D}\omega\,.
\end{eqnarray}
\end{proposition}
\begin{proof}
First of all, multiplying $\omega$ by
$e^{\pm\mathcal{H}(q,p)}=e^{\pm\frac{|p|_{q}^{2}\pm V(q)}{2}}$ or
$\langle p\rangle^{\pm\widehat{\text{deg}}}$ does not affect the
boundary conditions \eqref{eq.NeuHyp} and \eqref{eq.DirHyp} because
$|(-p_{1},p')|_{q}=|(p_{1},p')|_{q}$\,.\\
The analysis of both cases, Dirichlet and Neumann, follow the same
lines and we focus on the Neumann case.\\
Let us consider the density of $\mathcal{D}$ in
$D(\mathcal{K}_{\pm,\mathcal{H}}^{N})$\,.
We already know that the set of $\omega\in
\mathcal{C}^{\infty}_{0}(\overline{X};F_{X})$ which satisfy the 
boundary conditions \eqref{eq.NeuHyp} is dense in
$D(K_{\pm,\mathcal{H}}^{N})$\,.
Hence it suffices to approximate compactly supported elements of
$D(K_{\pm,\mathcal{H}}^{N})$ by elements of $\mathcal{D}$\,. When
$\chi\in \mathcal{C}^{\infty}_{0}(\overline{X};\rz)$ is a cut-off
function such that $\chi(0,q',-p_{1},p')=\chi(0,q',p_{1},p')$\,, one
has
$$
e^{+\mathcal{H}}\langle p\rangle^{-\widehat{deg}}
\mathcal{K}_{\pm,\mathcal{H}}^{N}e^{-\mathcal{H}}\langle
p\rangle^{\widehat{deg}}\chi
=
\left[
\mathcal{K}_{\pm,\mathcal{H}}^{N}+M_{\chi}(q,p,\partial_{p})\right]\chi\,,
$$ 
where $M_{\chi}(q,p,\partial_{p})$ satisfies the conditions
\eqref{eq.M0symb} and \eqref{eq.M1symb} of
Definition~\ref{de.genGKFP}.
Like in the proof of Proposition~\ref{pr.bundle},
$\mathcal{K}_{\pm,\mathcal{H}}^{N}$ can be identified locally and  up to an
additional correcting term $M(q,p,\partial_{p})$\,, with
$\tilde{U}(q)K_{\pm,0,g}^{g^{\mathfrak{f}}}\tilde{U}(q)$ where
$\tilde{U}$ is a unitary transform $\tilde{U}(q)$ from
$L^{2}(X,dqdp;(\mathfrak{f},g^{\mathfrak{f}}))$ to
$L^{2}(X,dqdp;F_{X})$\,. We can additionally assume
$\partial_{q^{1}}\tilde{U}(q)\equiv 0$ in our case where
$\partial_{q^{1}}m\equiv 0$\,.
With those local transformations (introduce also finite partition in
the position variable of
unity in order to reduce the problem to local ones), the density
actually amounts to the density of $\mathcal{D}(\overline{X},j)$  in
$D(K_{\pm, 0,g}^{g^{\mathfrak{f}}})$ stated in
Theorem~\ref{th.main0}\,.\\
Assume now that $\omega$ belongs to $\mathcal{D}$\,. Then
$\eta =e^{\mathcal{H}(q,p)}\langle
p\rangle^{-\widehat{\text{deg}}}\omega$ belongs to
$\mathcal{C}^{\infty}_{0}(\overline{X};F_{X})$\,, satisfies the
boundary conditions \eqref{eq.NeuHyp} and 
$$
|\partial_{q^{1}}\eta (q,p)|=\mathcal{O}(|q^{1}|^{\infty})\,.
$$

The extension by symmetry $\tilde{\eta}=\Sigma_{1}(\eta)$
satisfies $\tilde{\Sigma}_{1}(\tilde{\eta})=\tilde{\eta}$\,. 
We deduce
$\tilde{\Sigma}_{1}(d^{X}\tilde{\eta})=d^{X}\tilde{\Sigma}(\tilde{\eta})=d^{X}\tilde{\eta}$
outside $\left\{q^{1}=0\right\}$ while the above vanishing
condition ensures that $d^{X}\tilde{\eta}$ has no discontinuity at
$\left\{q^{1}=0\right\}$\,.
We deduce that $d^{X}\eta$ and therefore 
$\langle p\rangle^{\widehat{\text{deg}}}d^{X}_{\mathcal{H}}\langle
p\rangle^{-\widehat{\text{deg}}}\omega=\langle
p\rangle^{\widehat{\text{deg}}}d^{X}\eta$ fulfill
 the boundary conditions \eqref{eq.NeuHyp}. 
We have proved that $\langle p\rangle^{\widehat{\text{deg}}}d^{X}_{\mathcal{H}}\langle
p\rangle^{-\widehat{\text{deg}}}\omega$ belongs to
$D(K_{\pm,\mathcal{H}}^{N})$\,. The commutation 
\eqref{eq.commutN} is now nothing but the known commutation of
differential operators.
\end{proof}

\subsubsection{Neumann and Dirichlet realizations of the hypoelliptic Laplacian}
\label{se.NeuDirHyp}
We consider now the general case when $\overline{Q}$ is either a
compact manifold with boundary or a compact perturbation of the
euclidean half-space $\overline{\rz^{d}_{-}}$\,.
The hamiltonian function $\mathcal{H}$ is given by 
$$
\mathcal{H}(q,p)=\mathcal{E}(q,p)+V(q)=\frac{|p|_{q}^{2}}{2}+V(q)
$$
where $V(q)$ is a globally Lipschitz function.\\
Remember that $\mathcal{K}_{\pm,\mathcal{H}}$ is a geometric
Kramers-Fokker-Planck operator according to
Definition~\ref{de.genGKFP}.
After introducing the basis $e^{I}\hat{e}_{J}$ of $F_{X}=\bigwedge
T^{*}X$ like in Paragraph~\ref{se.witbis} with $e^{i}=dq^{i}$ and
$\hat{e}_{j}=dp_{j}-\Gamma^{\ell}_{ij}p_{\ell}dq^{i}$\,, our proposal
of  Neumann
(resp. Dirichlet) boundary condtions for
$\mathcal{K}_{\pm,\mathcal{H}}$ is the same as in the flat case:
\begin{eqnarray}
  \label{eq.NeuDirHyp}
  &&\omega^{J}_{I}(0,q',p_{1},p')\big|_{p_{1}<0}=
 (-1)^{k+|\left\{1\right\}\cap
    I|+|\left\{1\right\}\cap J|}\omega^{J}_{I}(0,q',-p_{1},p')
\\
\nonumber
\text{when}&& \omega =\omega^{J}_{I}e^{I}\hat{e}_{J}\,,
\end{eqnarray}
with $k=0$ (resp. $k=1$).\\
It still corresponds to the case
\begin{eqnarray*}
  && \mathbf{A}=0\,,\\
\text{and}&&
\mathbf{j}(e^{I}\hat{e}_{J})=(-1)^{k}(-1)^{|\left\{1\right\}\cap
  I|+|\left\{1\right\}\cap J|}e^{I}\hat{e}_{J}
\end{eqnarray*}
of Definition~\ref{de.BCfiber} and all the results of
Proposition~\ref{pr.bundle} with $\mathbf{A}=0$ apply to the
corresponding realization $\mathcal{K}_{\pm,\mathcal{H}}^{N}$
(resp. $\mathcal{K}_{\pm,\mathcal{H}}^{D}$) of
$\mathcal{K}_{\pm,\mathcal{H}}$\,.
This can be translated directly in terms of the hypoelliptic Laplacian
 $\mathcal{U}_{\mp,\mathcal{H}}^{2}$ because applying the weight
 $\langle p\rangle^{\pm\widehat{\text{deg}}}$ does not affect the
 boundary condition~\eqref{eq.NeuDirHyp}.
 \begin{proposition}
 \label{pr.DirNeuHyp}
Let $[\mathcal{U}_{\mp 1,\mathcal{H}}^{2}]^{N}$
(resp. $[\mathcal{U}_{\mp,\mathcal{H}}^{2}]^{D}$) be the closure in
the Hilbert space
$L^{2}(X,dqdp;(F_{X},\langle p\rangle^{2\widehat{\text{deg}}}g^{F_{X}}))$ of
$\mathcal{U}_{\mp 1,\mathcal{H}}^{2}=(\overline{d}_{\mp 1,
  \mathcal{H}}^{X}+d^{X}_{\mathcal{H}})^{2}$
 initially defined with the domain
$$
\mathcal{D}_{0}=\left\{\omega=\omega^{J}_{I}e^{I}\hat{e}_{J}\in\mathcal{C}^{\infty}_{0}(\overline{X};F^{X})\,,
  \eqref{eq.NeuDirHyp}\text{~holds}
 \right\}
$$
with $k=0$ (resp. $k=1$)\,. Then there exists $C\in\rz$ such that
$C+[\mathcal{U}_{\mp 1, \mathcal{H}}^{2}]^{N}$
(resp. $C+[\mathcal{U}^{2}]^{D}$) is maximal accretive. Moreover all
the results and estimates of Proposition~\ref{pr.bundle} with
$\mathbf{A}=0$
 apply to
$\mathcal{K}_{\pm,\mathcal{H}}^{N,D}=\langle p\rangle^{\widehat{\text{deg}}}2[\mathcal{U}_{\mp
  1, \mathcal{H}}^{2}]^{N,D}\langle p\rangle^{-\widehat{\text{deg}}}$
and are easily translated in $L^{2}(X,dqdp;(F_{X},\langle
p\rangle^{2\widehat{\text{deg}}}g^{F_{X}}))$
by conjugating with the unitary weight $\langle p\rangle^{\widehat{\text{deg}}}$\,.
 \end{proposition}
 \begin{proof}
   The fact that $\langle
   p\rangle^{-\widehat{\text{deg}}}\mathcal{K}_{\pm,\mathcal{H}}^{N,D}\langle
   p\rangle^{\widehat{\text{deg}}}$ is exactly the closure of
   $2\mathcal{U}_{\mp 1,\mathcal{H}}^{2}$ with the domain
   $\mathcal{D}_{0}$ comes from the density of
   $\mathcal{C}^{\infty}_{0}(\overline{X};F_{X})$ in
   $D(\mathcal{K}_{\pm,\mathcal{H}}^{N,D})$ endowed with its graph norm.
 \end{proof}
We end this paragraph with two questions:
\begin{enumerate}
\item The accurate spectral analysis of boundary Witten Laplacians
  developed in \cite{ChLi}\cite{HelNi}\cite{Lep3} relies strongly on
  the commutations
  \begin{eqnarray*}
&&(z-\Delta^{D,N}_{V})^{-1}d_{V}\omega =
d_{V}(z-\Delta_{V}^{D,N})^{-1}\quad,\\
&&
(z-\Delta_{V}^{D,N})^{-1}d_{V}^{*}\omega =
d_{V}^{*}(z-\Delta_{V}^{D,N})^{-1}\quad\,,
\end{eqnarray*}
for all $z\not\in \sigma (\Delta_{V}^{D,N})$ and all $\omega$ belonging to the
form domain of $\Delta_{V}^{D,N}$\,.\\
Is there a dense set in $L^{2}(X,dqdp;(F_{X},\langle
p\rangle^{\widehat{\text{deg}}}g^{F_{X}})$ in which the commutation
$$
(z-\mathcal{U}_{\mp
  1,\mathcal{H}}^{2})^{-1}d^{X}_{\mathcal{H}}=d_{\mathcal{H}}^{X}(z-\mathcal{U}_{\mp
1,\mathcal{H}}^{2})^{-1}\,,
$$
holds ?
\item Are the names Neumann and Dirichlet relevant~? In particular by
  introducing the parameter $b$ in
  $\mathcal{U}_{b,\mathcal{H}}^{2}$\,, does the operator
  $\mathcal{U}_{b,\mathcal{H}}^{N}$
  (resp. $\mathcal{U}_{b,\mathcal{H}}^{D}$) converge in some sense to
  $\Delta_{V}^{N}$ (resp. $\Delta_{V}^{D}$) as $b\to 0$~(large
  friction limit)~?
\end{enumerate}

\appendix
\section{Translation invariant model problems}
\label{se.kfpline}
In this appendix we review a few properties of the operator
$$
P_{\pm}=\pm
p.\partial_{q}+\frac{-\Delta_{p}+|p|^{2}}{2}=\sum_{j=1}^{d}\pm
p_{j}\partial_{q^{j}}+\frac{-\partial_{p_{j}}^{2}+p_{j}^{2}}{2}\quad\text{on~}
T^{*}Q\,,
$$
where $Q=\rz^{d'}\times \tz^{d''}$\,, $d',d''\in\nz$\,,
$\tz^{d''}=\rz^{d''}/\zz^{d''}$\,,
and its cotangent bundle $X=T^{*}Q=\rz^{d'}\times \tz^{d''}\times
\rz^{d}$\,, $d=d'+d''$\,,
 are endowed
with euclidean
metrics\,.\\
The operator $P_{\pm}$ is a typical example of a type~II hypoelliptic
operator (see \cite{Hor67, HelNi,Kol}) in the sense that $P_{\pm}u\in
\mathcal{C}^{\infty}(X)$ implies $u\in
\mathcal{C}^{\infty}(X)$\,.
More explicit subelliptic and pseudospectral estimates can be
given. This appendix collects a few of them and we refer to
\cite{HerNi} and
\cite{HelNi} for results with a potential, to  \cite{HiPr}\cite{Pra} for more
general results concerned with differential operators
with quadratic symbols and to \cite{Leb1, Leb2} for geometric
Kramers-Fokker-Planck operators (see also Section~\ref{se.geoKFP}).
\\
Remember the notation
$
\mathcal{H}^{s'}=
\left\{
u\in \mathcal{S}'(\rz^{d})\,,
(\frac{d}{2}+\mathcal{O})^{\frac{s'}{2}}u
\in L^{2}(\rz^{d})
\right\}
$
 with
$\mathcal{O}=\frac{-\Delta_{p}+|p|^{2}}{2}$\,, and remember that
$$
\|\big(\frac{-\Delta_{p}+|p|^{2}+d}{2}\big)^{N}u\|_{L^{2}}\quad\text{and}\quad 
\sum_{|\alpha|+|\beta|=2N}\|p^{\alpha}\partial_{p}^{\beta}u\|_{L^{2}}
$$
are equivalent norms on $\mathcal{H}^{2N}$ when $s'=2N\in 2\nz$\,,
owing to the global ellipticity of $\mathcal{O}+\frac{d}{2}$ (see
\cite{Helgl}\cite{HelNi} or \cite{HormIII}-Chap~18). The
operator $\mathcal{O}$ is diagonal in the Hermite basis
$(\varphi_{\nu})_{\nu\in\nz^{d}}$ given, with the multi-index notation, by
\begin{eqnarray*}
  &&
  \varphi_{\nu}=\frac{1}{\sqrt{\nu!}}(a^{*})^{\nu}\varphi_{0}\quad,\quad
  \varphi_{0}=
\frac{e^{-\frac{|p|^{2}}{2}}}{\pi^{\frac{d}{4}}}\,,\\
&& a=\frac{\partial_{p}+p}{\sqrt{2}}\quad,\quad
a^{*}=\frac{-\partial_{p}+p}{\sqrt{2}}\,,\\
&& (\mathcal{O}-\frac{d}{2})\varphi_{\nu}=a^{*}a\varphi_{\nu}=|\nu|\varphi_{\nu}\,.
\end{eqnarray*}
One easily checks with this Hermite basis,
$\ccap_{s'\in\rz}\mathcal{H}^{s'}=\mathcal{S}(\rz^{d})$ and
$\ccup_{s'\in\rz}\mathcal{H}^{s'}=\mathcal{S}'(\rz^{d})$\,.
We shall also use the vertical weighted $L^{2}$-spaces 
$$
L^{2}_{s'}=\left\{u\in
  \mathcal{S}'(\rz^{d})\,, \langle p\rangle^{s'}u\in
  L^{2}(\rz^{d})\right\}\,.
$$
In all variables $(q,p)=(q,q',p)$\,, the Schwartz space
$\mathcal{S}(X)$ is the set 
$\mathcal{C}^{\infty}$-functions $u$ such that
$$
\forall (\alpha,\beta)\in
\nz^{d_{1}+d}\times\nz^{2d}\,, 
\sup_{(q,p)\in X}
\left|(q',p)^{\alpha}\partial_{(q,p)}^{\beta}u(q,p)\right|< +\infty\,,
$$
and its dual is denoted by $\mathcal{S}'(X)$\,.
By using  the
Fourier transform in $q'\in\rz^{d}$ and Fourier series in the variable $q''\in\tz^{d''}$\,,
$$
\hat{u}(\xi',\xi'',p)=\int_{\rz^{d'}\times\tz^{d''}}e^{-i\xi.
  q}u(q,p)~dq\,,\quad \xi=(\xi',\xi'')\in \rz^{d'}\times(2\pi\zz)^{d''}\,,
$$
the norms $(p_{N})_{N\in\nz}$ defined by
\begin{equation}
  \label{eq.semiS}
p_{N}(u)=\sqrt{\sum_{|\alpha|+|\beta|+m'\leq
  N}\|\xi^{\alpha}\partial_{\xi'}^{\beta}\hat{u}\|_{L^{2}(\rz^{d'}\times
  (2\pi\zz)^{d''};\mathcal{H}^{m'})}^{2}}\,,
\end{equation}
provide to $\mathcal{S}(X)$ its Fr{\'e}chet-space structure.\\
For $\mathcal{E}=\mathcal{H}^{s'}$ or
$\mathcal{E}=L^{2}_{s'}$\,, the Sobolev spaces of $\mathcal{E}$-valued
functions are defined by
$$
H^{s}(Q;\mathcal{E})=\left\{u\in \mathcal{S}'(X)\,,
  (1-\Delta_{q})^{\frac{s}{2}}u\in L^{2}(Q,dq;\mathcal{E})\right\}\,.
$$
When $d'=0$ (and only in this case), we know
 $\mathcal{S}(X)=\ccap_{s,s'\in\rz}H^{s}(Q;\mathcal{H}^{s'})$\,, 
$\mathcal{S}'(X)=\ccup_{s,s'\in\rz}H^{s}(Q;\mathcal{H}^{s'})$ and
the embedding $H^{s_{1}}(Q;\mathcal{H}^{s_{1}'})\subset
H^{s_{2}}(Q;\mathcal{H}^{s_{2}'})$ is compact as soon as $s_{1}>s_{2}$
and $s_{1}'> s_{2}'$\,.
\begin{theorem}
\label{th.app} 
The operator $K_{\pm}-\frac{d}{2}$ defined in $L^{2}(X,dqdp)$ by
\begin{eqnarray*}
  && D(K_{\pm})=\left\{u\in L^{2}(X,dqdp),  P_{\pm}u\in
    L^{2}(X,dqdp)\right\}\\
&&\forall u\in D(K_{\pm})\,,\quad K_{\pm}u=P_{\pm}u\,,
\end{eqnarray*}
is maximal accretive. The spaces $\mathcal{S}(X)$ and
$\mathcal{C}^{\infty}_{0}(X)$ are dense in $D(K_{\pm})$ endowed
with its graph norm.\\
For any $\lambda\in\rz$ the operator $(P_{\pm}-i\lambda)$ is
an automorphism of $\mathcal{S}(X)$ and
of $\mathcal{S}'(X)$\,. 
For any $(\lambda,s,s')\in\rz^{3}$\,, the quantities
\begin{equation}
  \label{eq.subell1}
\|\big(\frac{-\Delta_{p}+|p|^{2}+d}{2}\big)u\|_{H^{s}(Q;L^{2}_{s'})}\quad,\quad
\|u\|_{H^{s+\frac 2 3}(Q;L^{2}_{s'})}\quad\text{and}\quad
\langle \lambda\rangle^{\frac 1 2}\|u\|_{H^{s}(Q;L^{2}_{s'})}\,,
\end{equation}
are bounded by $C_{s'}\|(P_{\pm}-i\lambda)u\|_{H^{s}(Q;L^{2}_{s'})}$
with $C_{s'}$ independent of $(s,\lambda)$\,.\\
In particular with $s,s'=0$\,, this implies
$$
D(K_{\pm})=D(K_{\mp})=\left\{u\in L^{2}(X,dqdp)\,,\quad
  p.\partial_{q}u\,,\,
  \mathcal{O}u\in L^{2}(X,dqdp)\right\}\,.
$$
The adjoint of $K_{\pm}$ is $K_{\pm}^{*}=K_{\mp}$\,.
\end{theorem}
\begin{proof}
The accretivity of $P_{\pm}-\frac{d}{2}$ on $\mathcal{S}(X)$
comes from the direct integration by parts:
$$
\Real \langle u\,,\, P_{\pm}u\rangle=\langle u\,,\, \mathcal{O}
u\rangle \geq \frac{d}{2}\|u\|^{2}\,,
$$
for all $u\in \mathcal{S}(X)$\,. The maximal accretivity of $K_{\pm}$
comes from the fact that $(P_{\pm}-i\lambda)$ is an automorphism of
$\mathcal{S}'(X)$\,. The space $\mathcal{S}(X)$ 
is dense in $D(K_{\pm})$ because $(P_{\pm}-i\lambda)$ is an
automorphism of $\mathcal{S}(X)$\,. The density of
$\mathcal{C}^{\infty}_{0}(X)$ is obtained after truncating elements of
$\mathcal{S}(X)$ while $P_{\pm}$ is a differential operator with
polynomial coefficients.\\
Thus the problem amounts to solving
$(P_{\pm}-i\lambda)u=f$ in $\mathcal{S}'(X)$ and to study the
regularity properties of $u$ according to $f$\,.
Due to the translational invariance in $q$ and after a Fourier transform
(series in $\xi''$), this is reduced to
solving and finding parameter dependent estimates for the equation
$$
(P_{\pm\xi}-i\lambda)\hat{u}(\xi,p)=\hat{f}(\xi,p),\quad
\xi=(\xi',\xi'')\in\rz^{d'}\times(2\pi\zz)^{d''}\,,
$$
with 
$$
P_{\xi}=ip.\xi -\frac{\Delta_{p}+|p|^{2}}{2}\,.
$$
Since the topology of $\mathcal{S}(X)$ is given by the
countable family of
norms $(p_{N})_{N\in\nz}$ given by \eqref{eq.semiS}  
the invertibility of $P_{\pm}-i\lambda$ in $\mathcal{S}(X)$ is a
consequence of Lemma~\ref{le.Pxila}. The invertibility in
$\mathcal{S}'(X)$ is deduced by duality with $P_{\pm}^{*}=P_{\mp}$\,.
The estimates of the quantities \eqref{eq.subell1}  are
$\xi$-integrated versions of Lemma~\ref{le.Pxila2}.\\
For the adjoint the identity
$$
\langle u\,,\, K_{\pm}v\rangle=\langle K_{\mp} u\,,\, v\rangle
$$
holds for any $u,v\in \mathcal{S}(X)$\,. By the density of
$\mathcal{S}(X)$ in $D(K_{\pm})$ this means that $u\in \mathcal{S}(X)$
belongs to $D(K_{\pm}^{*})$ with $K_{\pm}^{*}u=K_{\mp}u$\,. Since
$\mathcal{S}(X)$ is dense in $D(K_{\mp})$\,, this implies
$K_{\mp}\subset K_{\pm}^{*}$\,. Both operators $K_{\mp}$ and
$K_{\pm}^{*}$ are maximal accretive. The inclusion is therefore an equality.
\end{proof}
Let us first consider the resolvent $(P_{\xi}-i\lambda)^{-1}$ in
the $\mathcal{H}^{s'}$-space, $s'\in\rz$\,.
\begin{lemma}
\label{le.Pxila}
For $\xi\in\rz^{d'}\times(2\pi\zz)^{d''}$ and $\lambda\in\rz$\,, the operator
$P_{\xi}-i\lambda=i(p.\xi-\lambda)+\frac{-\Delta_{p}+p^{2}}{2}$ is an
automorphism of $\mathcal{S}(\rz^{d})$ and $\mathcal{S}'(\rz^{d})$\,.
For any $s'\in\rz$ and $\alpha\in\nz^{d_{1}}$ there exists $C_{s',\alpha}>0$
such that
\begin{multline*}
(1+|\xi|^{\frac 2 3}+|\lambda|^{\frac 1
  2})\|\partial_{\xi'}^{\alpha}(P_{\xi}-i\lambda)^{-1}\|_{\mathcal{L}(\mathcal{H}^{s'})}+
\|\partial_{\xi'}^{\alpha}(P_{\xi}-i\lambda)^{-1}\|_{\mathcal{L}(\mathcal{H}^{s'};\mathcal{H}^{s'+2+|\alpha|})}
\\
\leq C_{\alpha,s}\langle \xi\rangle^{|s'|(1+|\alpha|)}\,,
\end{multline*}
holds for all $(\xi,\lambda)\in\rz^{d'}\times(2\pi \zz)^{d''}\times \rz$\,.
\end{lemma}
\begin{proof}~
\noindent\textbf{a)} Consider first the case $\alpha=0$ and $s'=0$:
The operator $P_{\xi}-\frac{d}{2}$ with domain $D(P_{\xi})=\mathcal{H}^{2}$ is
maximal accretive and $(P_{\xi}-i\lambda)^{-1}\in
\mathcal{L}(L^{2}(\rz^{d},dp);\mathcal{H}^{2})$ for all $\xi\in\rz^{d'}\times(2\pi\zz)^{d''}$\,. In
order to prove the resolvent estimates, take $u\in D(P_{\xi})$ and compute
\begin{equation}
\label{eq.subellFou}
\|(P_{\xi}-i\lambda)u\|^{2}
=\|(p.\xi-\lambda)u\|^{2}+\|(\frac{-\Delta_{p}+|p|^{2}}{2})u\|^{2}+\langle
u\,,\, \xi.D_{p}u\rangle\,,
\end{equation}
with $D_{p}=\frac 1 i \partial_{p}$ and $\|~\|=\|~\|_{L^{2}(\rz^{d},dp)}$\,.
With a rotation in the variable $p$\,, the problem is reduced to
the case $\xi=(\xi_{1},0,\ldots,0)\in\rz^{d}$\,.
The separation of variables $p=(p_{1},p')\in \rz\times\rz^{d-1}$ gives
\begin{eqnarray*}
\frac{-\Delta_{p}+|p|^{2}}{2}
&=&\frac{D_{p_{1}}^{2}+p_{1}^{2}}{2}+\frac{-\Delta_{p'}+|p'|^{2}}{2}
\\
&=&
\sum_{\nu\in\nz^{d-1}}\left[\frac{D_{p_{1}}^{2}+p_{1}^{2}}{2}+\frac{d-1}{2}+|\nu|\right]|\varphi_{\nu}\rangle
\langle \varphi_{\nu}|\,.
\end{eqnarray*}
where $(\varphi_{\nu})_{\nu\in\nz^{d-1}}$ is the Hermite
basis in $L^{2}(\rz^{d-1},dp)$\,.
After writing
$u=\sum_{\nu\in\nz^{d-1}}u_{\nu}(\xi_{1},p_{1})\varphi_{\nu}$
the problem is reduced to finding lower bounds of
$$
\|(p_{1}\xi_{1}-\lambda)u_{\nu}\|^{2}+\|(\frac{D_{p_{1}}^{2}+p_{1}^{2}}{2}+\frac{d-1}{2}+|\nu|)u_{\nu}\|^{2}
+\langle u_{\nu}\,,\, \xi_{1}D_{p_{1}}u_{\nu}\rangle\,,
$$
which is a one dimensional problem parametrized by $(\xi_{1},\lambda,\nu)$\,.
The identity
$$
(D_{p_{1}}^{2}+p_{1}^{2}+a)^{2}=D_{p_{1}}^{4}+p_{1}^{4}+2a(D_{p_{1}}^{2}+p_{1}^{2})+a^{2}
+2D_{p_{1}}p^{2}_{1}D_{p_{1}}-2\quad,\quad a\geq 0\,,
$$
gives
\begin{multline*}
\|\big(\frac{D_{p_{1}}+p_{1}^{2}}{2}+\frac{d-1}{2}+|\nu|\big)u_{\nu}\|^{2}\geq 
\|\frac{D_{p_{1}}^{2}}{2}u_{\nu}\|^{2}+\|\frac{p_{1}^{2}}{2}u_{\nu}\|^{2}
\\+\left[(\frac{d-1}{2}+|\nu|)^{2}-\frac{1}{2}\right]\|u_{\nu}\|\,.
\end{multline*}
We are thus looking for a lower bound of
\begin{multline}
  \label{eq.estimfour}
\|(p_{1}\xi_{1}-\lambda)u_{\nu}\|^{2}+\|\frac{D_{p_{1}}^{2}}{2}u_{\nu}\|^{2}
+\xi_{1}\langle
u_{\nu}\,,\, D_{p_{1}}u_{\nu}\rangle
+
\|\frac{p_{1}^{2}}{2}u_{\nu}\|^{2}\\
+\left[(\frac{d-1}{2}+|\nu|)^{2}-\frac{1}{2}\right]\|u_{\nu}\|^{2}\,.
\end{multline}
For $\xi_{1}\neq 0$ (the case $\xi_{1}=0$ is obvious)\,, set 
$$
u_{\nu}(\xi_{1},p_{1})=|\xi_{1}|^{\frac 1
  6}\phi(\xi_{1},\sign(\xi_{1})|\xi_{1}|^{\frac 1 3}p_{1}-|\xi_{1}|^{-\frac 2 3}\lambda)\,,
$$
corresponding to the change of variable $t=\sign(\xi_{1})|\xi_{1}|^{\frac 1 3}p_{1}-|\xi_{1}|^{-\frac 2 3}\lambda$\,,
and compute the first three terms
\begin{multline*}
\|(p_{1}\xi_{1}-\lambda)u_{\nu}\|^{2}
+\|\frac{D_{p_{1}}^{2}}{2}u_{\nu}\|^{2}+\xi_{1}\langle
u_{\nu}\,,\, D_{p_{1}}u_{\nu}\rangle
\\
= |\xi_{1}|^{\frac 4 3}\left[
\|t\phi\|^{2}+\|\frac{D_{t}^{2}}{2}\phi\|^{2}
+\sign(\xi_{1})\langle \phi\,, D_{t} \phi\rangle
\right]
\\
=
\frac{|\xi_{1}|^{\frac 4 3}}{4}\|(\sign(\xi_{1}) 2it+D_{t}^{2})\phi\|^{2}\,.
\end{multline*}
On $L^{2}(\rz_{t},dt)$\,, the complex Airy operator
$D_{t}^{2}+\sign(\xi_{1})2it$ has a compact resolvent because
$4\|f\|^{2}+\|(D_{t}^{2}+\sign(\xi_{1})2it)f\|^{2}\geq \|2t
f\|_{L^{2}}^{2}+\|D_{t}^{2}f\|^{2}$
and it is injective $(D_{t}^{2}+\sign(\xi_{1})2it)f=0$ (or equivalently
$(\tau^{2}+\sign(\xi_{1})2\partial_{\tau})\hat{f}=0$) has no solution in
$L^{2}(\rz,dt)$\,.
We deduce
$$
\|(D_{t}^{2}+\sign(\xi_{1})2it)\phi\|^{2}\geq
C^{-1}\left[\|\phi\|^{2}+\|2t\phi\|^{2}
+\|D_{t}^{2}\phi\|^{2}\right]\,.
$$
 and
\begin{multline}
  \label{eq.estimxi}
\|(p_{1}\xi_{1}-\lambda)u_{\nu}\|^{2}
+\|\frac{D_{p_{1}}^{2}}{2}u_{\nu}\|_{L^{2}}^{2}+\xi_{1}\langle
u_{\nu}\,,\, D_{p_{1}}u_{\nu}\rangle
\\
\geq
\frac{1}{C'}\left[|\xi_{1}|^{\frac{4}{3}}\|u_{\nu}\|^{2}
+\|(p_{1}\xi_{1}-\lambda)u_{\nu}\|^{2}
+\|D_{p_{1}}^{2}u_{\nu}\|^{2}
\right]
\,,
\end{multline}
 with a uniform constant
$C'\geq 1$\,.  By summing over $\nu\in\nz^{d-1}$\,,
\eqref{eq.subellFou}, \eqref{eq.estimfour} and \eqref{eq.estimxi}
lead to
\begin{multline*}
C'\|(P_{\xi}-i\lambda)u\|^{2}\geq
\|(p.\xi-\lambda)u\|_{L^{2}}^{2}
+\|D_{p_{1}}^{2}u\|_{L^{2}}^{2}+\|p_{1}^{2}u\|^{2}
\\
+\||\xi|^{\frac 2 3}u\|_{L^{2}}^{2}
+\|(\frac{-\Delta_{p'}+|p'|^{2}+d-1}{2})u\|^{2}
-\frac{1}{2}\|u\|^{2}\,.
\end{multline*}
With
$$
\|u\|\|(P_{\xi}-i\lambda)u\|\geq \Real \langle u\,,\,
(P_{\xi}-i\lambda)u\rangle\geq \frac{d}{2}\|u\|^{2}\,,
$$
which gives $\|(P_{\xi}-i\lambda)u\|^{2}\geq
\frac{d^{2}}{4}\|u\|_{L^{2}}^{2}$\,,  we obtain
$$
C''\|(P_{\xi}-i\lambda)u\|_{L^{2}}^{2}\geq
(1+|\xi|^{4/3})\|u\|_{L^{2}}^{2}+\|u\|_{\mathcal{H}^{2}}^{2}+\|(p.\xi-\lambda)u\|_{L^{2}}^{2}\,.
$$
For  $|\lambda|\|u\|_{L^{2}}^{2}$\,, simply use
\begin{eqnarray*}
|\lambda|\|u\|_{L^{2}}^{2}&=&i\sign(\lambda)\langle u\,,\,
i(p.\xi-\lambda)u\rangle + \sign(\lambda)\langle u\,,\,
p.\xi u\rangle\\
&\leq
&\|u\|\|(p.\xi-\lambda)u\|+\int_{\rz}|\xi||p||u(p)|^{2}~dp
\end{eqnarray*}
with $|\xi||p|\leq C_{1}(|\xi|^{4/3}+|p|^{4})$ in the integral.
This finishes the proof for $\alpha=0$ and $s'=0$\,.\\
\noindent\textbf{b)} Consider the case $\alpha=0$ and $s'=2N\in
2\nz\setminus\left\{0\right\}$:
By the global ellipticity of $P_{\xi}$ and the uniqueness in $L^{2}(\rz^{d},dp)$\,, the equation 
$(P_{\xi}-i\lambda)u=f$ has a unique solution in
$\mathcal{S}'(\rz^{d})$ which belongs necessarily to
$\mathcal{H}^{2N+2}$ when $f\in \mathcal{H}^{2N}\subset
L^{2}(\rz^{d},dp)$\,. In $\mathcal{S}'(\rz^{d})$ and for
$\varepsilon>0$\,,
$(P_{\xi}-i\lambda)\left[1+\varepsilon^{2}(-\Delta_{p}+p|^{2})\right]^{N}u$ equals
$$
\left[1+\varepsilon^{2}(-\Delta_{p}+|p|^{2})\right]^{N}f
+ i
\sum_{k=1}^{d}\sum_{M=1}^{N}\xi_{k}\left(
  \begin{array}[c]{c}
N\\M
  \end{array}
\right)
\left[p_{k}\,, \varepsilon^{2M}(-\Delta_{p}+|p|^{2})^{M}\right]u\,.
$$
The expression  of the commutator 
$\left[p_{k}\,,  \varepsilon^{2M}(-\Delta_{p}+|p|^{2})^{M}\right]$
occurring in the second term as
$$
\sum_{0\leq n\leq M-1, |\alpha|+|\beta|\leq 1}
\varepsilon^{2(M-n)-1}c_{M,n,\alpha,\beta}(\varepsilon
p)^{\alpha}(\varepsilon D_{p})^{\beta}[\varepsilon^{2}|D_{p}|^{2}+\varepsilon^{2}|p|^{2}]^{n}\,,
$$
is easily obtained by induction on $M\geq 1$\,. With $2(M-n)-1\geq 1$
and the operator
$$
(\varepsilon p)^{\alpha}(\varepsilon D_{p})^{\beta}[\varepsilon^{2}
|D_{p}|^{2}+\varepsilon^{2}|p|^{2})]^{n}
\left[1+\varepsilon^{2}(|D_{p}|^{2}+|p|^{2})\right]^{-N}
$$
uniformly bounded in $\mathcal{L}(L^{2}(\rz^{d},dp))$\,.
The function
$v=\left[1+\varepsilon^{2}(-\Delta_{p}+|p|^{2})\right]^{N}u\in L^{2}(\rz^{d},dp)$ solves
$$
\left[(P_{\xi}-i\lambda)+\sum_{k=1}^{d}\varepsilon
  \xi_{k}R_{N,k}(\varepsilon)\right]
v=\tilde{f}=\left[1+\varepsilon^{2}(-\Delta_{p}+|p|^{2})\right]^{N}f\,,
$$ 
with $\|R_{N,k}\|_{\mathcal{L}(L^{2})}\leq C_{N}$\,.
Meanwhile we know
$\|(P_{\xi}-i\lambda)^{-1}\|_{\mathcal{L}(L^{2})}\leq
\frac{2}{d}$\,. By taking $\varepsilon \leq \min\left\{\frac{1}{4C_{N}\langle
  \xi\rangle}, \sqrt{\frac{2}{d}}\right\}$\,, one gets
$$
v=(P_{\xi}-i\lambda)^{-1}
\left[\Id +\sum_{k=1}^{d}\varepsilon
  \xi_{k}R_{k}(P_{\xi}-i\lambda)^{-1}\right]^{-1}
\tilde{f}\,,
$$
and 
\begin{eqnarray*}
\varepsilon^{2N}\|u\|_{\mathcal{H}^{2N}}
&\leq&
\left\|\left[1+\varepsilon^{2}(-\Delta_{p}+|p|^{2})\right]^{N}u
\right\|=\|v\|\\
&\leq& 
C(1+|\xi|^{\frac 2 3}+|\lambda|^{\frac 1 2})^{-1}\left\| \left[\Id +
  \sum_{k=1}^{d}\varepsilon \xi_{k}R_{k}(P_{\xi}-i\lambda)^{-1}\right]^{-1}\tilde{f}\right\|\,,\\
&\leq &2C(1+|\xi|^{\frac 2 3}+|\lambda|^{\frac 1 2})^{-1}\|\tilde{f}\|_{L^{2}}
\leq C'(1+|\xi|^{\frac 2 3}+|\lambda|^{\frac 1 2})^{-1}\|f\|_{\mathcal{H}^{2N}}\,.
\end{eqnarray*}
Similarly we obtain
\begin{eqnarray*}
  \varepsilon^{2N}\|u\|_{\mathcal{H}^{2N+2}}
&\leq&
\left\|\left[1+\varepsilon^{2}(-\Delta_{p}+|p|^{2})\right]^{N}u
\right\|_{\mathcal{H}^{2}}
\\
&\leq &
C''\|\tilde{f}\|
\leq C^{(3)}\|f\|_{\mathcal{H}^{2N}}\,.
\end{eqnarray*}
We have proved
$$
(1+|\xi|^{\frac 2 3}+|\lambda|^{\frac 1 2})
\|(P_{\xi}-i\lambda)^{-1}\|_{\mathcal{L}(\mathcal{H}^{2N})}+
\|(P_{\xi}-i\lambda)^{-1}\|_{\mathcal{L}(\mathcal{H}^{2N},\mathcal{H}^{2N+2})}\leq
C_{N}'\langle \xi\rangle^{2N}\,.
$$
\noindent\textbf{c)} The general result for $s'\in\rz$ and $\alpha=0$
follows by duality and interpolation after noticing that the formal
adjoint of $P_{\xi}-i\lambda$ is $P_{-\xi}+i\lambda$\,.\\
\noindent\textbf{d)}  The mapping $\xi'\to (P_{\xi}-i\lambda)^{-1}$ is
differentiable as an $\mathcal{S}$-continuous (or
$\mathcal{S}'$-continuous)
 operator valued function of $\xi'\in \rz^{d'}$\,. Its $\alpha$-th derivative has
 the form
$$
\partial_{\xi'}^{\alpha}(P_{\xi}-i\lambda)^{-1}=
\sum_{\substack{
(\beta_{1},\ldots, \beta_{|\alpha|})\in (\nz^{d})^{|\alpha|}\,,\\
|\beta_{j}|=1
} 
}
c_{\beta_{1},\ldots,\beta_{|\alpha|}}(P_{\xi}-i\lambda)^{-1}
\left[
\prod_{j=1}^{|\alpha|} (p')^{\beta_{j}}(P_{\xi}-i\lambda)^{-1}\right]\,.
$$
From \textbf{b)} and \textbf{c)}, we deduce
$$
\|(p')^{\beta_{j}}(P_{\xi}-i\lambda)^{-1}\|_{\mathcal{L}(\mathcal{H}^{s'})}\leq
C_{s}\langle \xi\rangle^{|s'|}\,.
$$
The result then follows from the estimates for $\alpha=0$ applied to
the first factor 
$(P_{\xi}-i\lambda)^{-1}$\,.
\end{proof}

\begin{lemma}
\label{le.Pxila2}
For any $s'\in\rz$\,, there exists $C_{s'}>0$
such that
\begin{equation*}
(1+|\xi|^{\frac 2 3}+|\lambda|^{\frac 1 2})\|(P_{\xi}-i\lambda)^{-1}\|_{\mathcal{L}(L^{2}_{s'})}+
\|(-\Delta_{p}+|p|^{2})
(P_{\xi}-i\lambda)^{-1}\|_{\mathcal{L}(L^{2}_{s'})}
\leq C_{s'}\,,
\end{equation*}
holds for all $(\xi,\lambda)\in\rz^{d'}\times(2\pi\zz)^{d''}\times \rz$\,.
\end{lemma}
\begin{proof}
  The case $s'=0$ is already proved in Lemma~\ref{le.Pxila}.
Assume $(P_{\xi}-i\lambda)u= f$
with $f\in L^{2}_{s'}(\rz^{d})$
and write in $\mathcal{S}'(\rz^{d})$:
$$
(P_{\xi}-i\lambda)\langle \varepsilon p\rangle^{s'}u=\langle
\varepsilon p\rangle^{s'}f - \sum_{k=1}^{d}\partial_{p_{k}}\left[\varepsilon g_{s',k}(\varepsilon
p)\langle
\varepsilon p\rangle^{s'}u\right]+\frac{\varepsilon^{2}}{2}g_{s'}(\varepsilon
p)\langle \varepsilon p\rangle^{s'}u\,,
$$
for $\varepsilon>0$ to be fixed, 
with 
$$
g_{s',k}(p)=\frac{\partial_{p_{k}}\langle p\rangle^{s'}}{\langle p\rangle^{s'}}=s'\langle
p\rangle^{-2}p_{k}\quad\text{and}\quad
g_{s'}(p)=\frac{\Delta_{p}\langle p\rangle^{s'}}{\langle p\rangle^{s'}}=(s'^{2}+(d-2)s')\langle p\rangle^{-2}.
$$
We obtain
$$
\left[(P_{\xi}-i\lambda)+\varepsilon R(s',\varepsilon)\right]\langle
\varepsilon p\rangle^{s'}u=\langle \varepsilon p\rangle^{s'}f\,,
$$
with $\|R(s',\varepsilon)\|_{\mathcal{L}(\mathcal{H}^{1}; L^{2}(\rz^{d},dp))}\leq
C_{s',1}$\,. We know
$\|(P_{\xi}-i\lambda)^{-1}\|_{\mathcal{L}(L^{2}(\rz^{d},dp);\mathcal{H}^{1})}\leq
C_{2}$ and by taking $\varepsilon_{s'}\leq \frac{1}{2 C_{s',1}C_{2}}$ we get
$$
\langle \varepsilon_{s'} p\rangle^{s}u
=(P_{\xi}-i\lambda)^{-1}\left[\Id + \varepsilon_{s'}
  R(s',\varepsilon_{s'})(P_{\xi}-i\lambda)^{-1}\right]^{-1}\langle
\varepsilon_{s'} p\rangle^{s}f\,.
$$
This proves
$$
(1+|\xi|^{\frac 2 3}+|\lambda|^{\frac 1 2})\|\langle \varepsilon_{s'}
p\rangle^{s'}u\|+
\|(-\Delta_{p}+\frac{|p^{2}|}{4}+\frac{d}{2})\langle \varepsilon_{s'}
p\rangle^{s'}u\|\leq C_{s'}\|\langle
\varepsilon_{s'}p\rangle^{s'}f\|\,.
$$
The equivalences of norms between $\|u\|_{L^{2}_{s'}}$ and
$\|\langle \varepsilon_{s'}p\rangle^{s'}u\|$\,,
and between $\|(-\Delta_{p}+|p|^{2})\langle
\varepsilon_{s'}p\rangle^{s'}u\|$ and 
$
\|(\frac{-\Delta_{p}+|p|^{2}+d}{2})u\|_{L^{2}_{s'}}=\|\langle
p\rangle^{s'}(\frac{-\Delta_{p}+|p|^{2}+d}{2})u\|$
finally yields the result.
\end{proof}
Theorem~\ref{th.app} can be used to solve $(P_{\pm}-i\lambda)u=f$ in the
distributional sense.
In particular when $Q=\rz\times Q'$ with $Q'=\rz^{d_{1}-1}\times
\tz^{d_{2}}$\,, $d=d_{1}+d_{2}$\,, with the coordinates $q=(q^{1},q',p_{1},p')$ on $X=T^{*}Q$\,,
and $\gamma\in L^{2}(Q'\times\rz^{d},\frac{dq'dp}{|p_{1}|})\subset
L^{2}(Q'\times\rz^{d}, \frac{dq'dp}{\langle p\rangle})$ the distribution
$\gamma(q',p)\delta_{0}(q_{1})$
 belongs to $H^{s_{0}}(Q;L^{2}_{-\frac 1
  2})$ for any $s_{0}< -\frac 1 2$\,. The equation
 $(P_{\pm}-i\lambda)u=\gamma(q',p)\delta_{0}(q_{1})$ then admits a unique
 solution in $H^{s_{0}}(Q;L^{2}_{-\frac 1 2})$ with
 \begin{multline}
   \label{eq.estimsing}
\||D_{p}|^{2}u\|_{H^{s_{0}}(Q;L^{2}_{-\frac 1
    2})} + \||p|^{2}u\|_{H^{s_{0}}(Q;L^{2}_{-\frac 1
    2})}+\langle \lambda\rangle^{\frac{1}{2}}\|u \|_{H^{s_{0}}(Q;L^{2}_{-\frac 1 2})}
\\
+
\|u\|_{H^{s_{0}+\frac 2 3}(Q;L^{2}_{-\frac 1 2})}\leq C_{s_{0}}
\|\gamma\|_{L^{2}(Q'\times\rz^{d},\frac{dq'dp}{|p_{1}|})}\,.
\end{multline}
This estimate can be given a more precise form with a simple
integration by parts.
\begin{proposition}
\label{pr.trace1D}
Assume $Q=\rz\times Q'$ with $Q'=\rz^{d_{1}-1}\times
\tz^{d_{2}}$\,, $d_{1}\geq 1$\,, $d=d_{1}+d_{2}$\,,
 with the coordinates $q=(q^{1},q',p_{1},p')$ on $X=T^{*}Q$\,.\\
For any $u\in D(K_{\pm})$\,, the following trace estimate holds
\begin{equation}
  \label{eq.subell2}
\langle \lambda\rangle^{\frac 1 4}\|u(q_{1}=0,q',p)\|_{L^{2}(Q'\times \rz^{d},|p_{1}|dq'dp)}\leq
C'\|(K_{\pm}-i\lambda)u\|\,.
\end{equation}
The dual version says that for any  $\gamma\in L^{2}(Q'\times\rz^{d},
\frac{dq'dp}{|p_{1}|})$\,, the solution $u$ to
$(P_{\pm}-i\lambda)u=\gamma(q',p)\delta_{0}(q_{1})$ belongs to 
$ L^{2}(X,dqdp)$ with the estimate
\begin{equation}
  \label{eq.subell3}
\langle \lambda\rangle^{\frac 1 4}\|u\|\leq
C''\|\gamma\|_{L^{2}(Q'\times \rz^{d},\frac{dq'dp}{|p_{1}|})}\,.
\end{equation}
\end{proposition}
\begin{proof}Set $X_{-}=(-\infty,0)\times Q'\times
\rz^{d}$\,.
For $u\in \mathcal{C}^{\infty}_{0}(X)$ and $\lambda\in\rz$\,,  write 
\begin{eqnarray*}
\int_{Q'\times\rz^{d}}|p_{1}||u(0,q',p)|^{2}~dq'dp
&&
=
2\Real\int_{X_{-}}|p_{1}|\partial_{q_{1}}\overline{u}(q,p)
u(q,p)~dq dp
\\
&&\hspace{-2cm}=
\pm 2\Real\int_{X}(\pm p.\partial_{q}+i\lambda)\overline{u}(q,p)
\sign(p_{1})u(q,p)~dq dp
\\
&&\hspace{-2cm}=
\pm
2\Real\langle (\pm p\partial_{q}-i\lambda)u\,,\, 1_{\rz_{-}}(q)\sign(p)u \rangle
\\
&&\hspace{-2cm}\leq 2\|(\pm p\partial_{q}-i\lambda)u\|\|u\|\leq
C\langle\lambda\rangle^{-\frac 1 2}\|(K_{\pm }-i\lambda)u\|^{2}\,,
\end{eqnarray*}
by applying \eqref{eq.subell1} with $s,s'=0$\,.
The extension of the estimate to any $u\in D(K_{\pm})$ is a
consequence of the density of $\mathcal{C}^{\infty}_{0}(X)$ in $D(K_{\pm})$\,.\\
The existence and uniqueness of a solution $u\in L^{2}(X,dqdp)$ to
$(P_{\pm}-i\lambda)u=\gamma(q',p)\delta_{0}(q_{1})$ comes from
Theorem~\ref{th.app} as explained for
\eqref{eq.estimsing}.
 For any $\phi\in L^{2}(Q'\times\rz^{d},|p_{1}|dq'dp)$ and any
$f\in L^{2}(X,dqdp)$ the previous result applied to $K_{\pm}^{*}=K_{\mp}$  gives
\begin{multline*}
\left|\int_{Q'\times\rz^{d}}
\overline{\phi(q',p)}[(K_{\pm}^{*}+i\lambda)^{-1}f](0,q',p)|p_{1}|dq'dp\right|
\\
\leq
C'\langle\lambda\rangle^{-\frac 1 4}\|\phi\|_{L^{2}(Q'\times\rz^{d},|p_{1}|dq'dp)}\|f\|\,.
\end{multline*}
But the above integral equals
\begin{multline*}
\int_{\rz}\overline{\phi(q',p)}[(K^{*}+i\lambda)^{-1}f](0,q',p)|p_{1}|dq'dp
\\
=\langle \phi(q',p)|p_{1}|\delta_{0}(q_{1})\,,\, (K_{\pm}^{*}+i\lambda)^{-1}f\rangle
\\
=\langle (P_{\pm}-i\lambda)^{-1}[\phi(q',p)|p_{1}|\delta_{0}(q_{1})]\,,\, f\rangle\,.
\end{multline*}
After setting $\gamma=|p_{1}|\phi(q',p)$\,, we obtain
\begin{multline*}
\left|
\langle f\,,\, (P_{\pm}-i\lambda)^{-1}[\gamma(q',p)\delta_{0}(q_{1})]\rangle
\right|
\leq C''\langle \lambda\rangle^{-\frac 1
  4}\|f\|\||p_{1}|^{-1}\gamma\|_{L^{2}(Q'\times\rz^{d},|p_{1}|dq'dp)}
\\=
C''\langle \lambda\rangle^{-\frac 1
  4}\|f\|\|\gamma\|_{L^{2}(Q'\times\rz^{d},\frac{dq'dp}{|p_{1}|})}\,,
\end{multline*}
which ends the proof.
\end{proof}
We end with another application of \eqref{eq.estimsing}.
\begin{proposition}
\label{pr.regQcomp}
Assume $Q=\rz\times Q'$ with $Q'=\tz^{d'}\times \rz^{d''-1}$\,, 
$d''\geq 1$\,,$d=d'+d''$\,,  and take the coordinates
$(q^{1},q',p_{1},p')\in \rz\times Q'\times \rz\times\rz^{d-1}$ in $X=T^{*}Q$\,.\\
If $u\in L^{2}(Q;\mathcal{H}^{1})$ solves 
$$
(P_{\pm}-i\lambda)u=\gamma(q',p)\delta_{0}(q^{1})
$$
with $\gamma\in L^{2}(Q'\times\rz^{d},\frac{dq'dp}{|p_{1}|})$\,. 
Then for any $s\in [0,\frac{1}{9})$\,, $u$ belongs to
$H^{s}(Q;\mathcal{H}^{0})$ with the estimate
$$
\|u\|_{H^{s}(Q;\mathcal{H}^{0})}\leq
C_{s}\left[\|\gamma\|_{L^{2}(Q'\times\rz^{d},\frac{dq'dp}{|p_{1}|})}
+\|u\|_{L^{2}(Q;\mathcal{H}^{1})}\,.
\right]
$$
When $d''=1$\,, $R\in (0,+\infty)$ and $s\in (0,\frac{1}{9})$\,, the
embedding $H^{s}(Q;\mathcal{H}^{0})\cap
L^{2}(Q;\mathcal{H}^{1})\subset L^{2}([-R,R]\times Q'\times \rz^{d};
dqdp)$ is compact.
\end{proposition}
\begin{proof}
  From \eqref{eq.estimsing} we deduce
$$
\|\langle
p\rangle^{-\frac{1}{2}}u\|_{H^{\nu}(Q;\mathcal{H}^{0})}\leq 
\|u\|_{H^{\nu}(Q;L^{2}_{-\frac{1}{2}})}\leq
C_{\nu}\|\gamma\|_{L^{2}(Q'\times\rz^{d},\frac{dq'dp}{|p_{1}|})}\,,
$$
for all $\nu\in [0,\frac{1}{6})$\,. The assumption $u\in
L^{2}(Q;\mathcal{H}^{1})$ implies
$$
\|\langle p\rangle u\|_{L^{2}(X,dqdp)}\leq C_{d}\|u\|_{L^{2}(Q;\mathcal{H}^{1})}\,.
$$
With 
$$
\|u\|_{H^{\frac{2}{3}\nu}(Q;\mathcal{H}^{0})}\leq C_{\nu}'\|\langle
p\rangle^{-\frac{1}{2}}u\|_{H^{\nu}(Q;\mathcal{H}^{0})}^{\frac{2}{3}}\|\langle
p\rangle u\|_{L^{2}(X,dqdp)}^{\frac{1}{3}}
$$
we choose $s=\frac{2}{3}\nu\subset [0,\frac{1}{9})$\,. We obtain
\begin{eqnarray*}
\|u\|_{H^{s}(Q;\mathcal{H}^{0})}
&\leq&
C_{s}\|\gamma\|_{L^{2}(Q'\times\rz^{d},\frac{dq'dp}{|p_{1}|})}^{\frac{2}{3}}\|u\|_{L^{2}(Q;\mathcal{H}^{1})}^{\frac{1}{3}}
\\
&\leq&
C_{s}'\left[\|\gamma\|_{L^{2}(Q'\times\rz^{d},\frac{dq'dp}{|p_{1}|})}+\|u\|_{L^{2}(Q;\mathcal{H}^{1})}
\right]\,.
\end{eqnarray*}
With $[-R,R]\times Q'=[-R,R]\times \tz^{d-1}$\,, the compactness statement is obvious. 
\end{proof}
\section{Partitions of unity}
\label{se.partunit}

We briefly review a few formulas with partitions of unities.
We shall work  on a riemannian manifold $M$ possibly with boundary and
$L^{2}(M)$ is endowed with the scalar product $\langle
u\,,\,v\rangle=\int_{M}\overline{u}v~d\text{Vol}_{g}$  and the norm
$\|u\|=\sqrt{\langle u\,,\,u\rangle}$\,, associated with the metric $g$\,.
The operator
$P=\sum_{|\alpha|\leq m}a_{\alpha}(x)\partial_{x}^{\alpha}$
 is a differential operator with $\mathcal{C}^{\infty}$-coefficients and its
 formal adjoint is denoted by $P^{*}$\,. The
family $(\chi_{\ell})_{\ell\in \mathcal{\mathcal{L}}}$ is a locally
finite family $\mathcal{C}^{\infty}_{0}(M;\rz)$ such that
$\sum_{\ell\in \mathcal{L}}\chi_{\ell}^{2}\equiv 1$\,.
 For every $\ell\in
\mathcal{L}$\,, $P_{\ell}$ will denote a local model of $P$\,, that is
$P\chi_{\ell}=P_{\ell}\chi_{\ell}$\,. We also assume that possible boundary
conditions behave well with the partition of unity. The following
formula are written for $u\in \mathcal{C}^{\infty}_{0}(M)$ (or possibly $u\in
  \mathcal{C}^{\infty}_{0}(\overline{M})$ fulfilling some boundary
  conditions) or $u\in L^{2}(M)$ provided that all the terms make sense in
  $\mathcal{D}'(M)$ or in $L^{2}(M)$\,.
Finally the notation $\ad^{k}_{A}B$ is defined by
 $\ad_{A}B=\left[A,B\right]$ and $\ad_{A}^{k}B=\ad_{A}\ad_{A}^{k-1}B$\,.
\begin{eqnarray}
&&\hspace{-1cm}
\label{eq.partunit2}
\|Pu\|^{2}-\sum_{\ell\in
  \mathcal{L}}\|P_{\ell}\chi_{\ell}u\|^{2}=-\sum_{\ell\in
  \mathcal{L}}\|(\ad_{\chi_{\ell}}P) u\|^{2}
+\Real \langle
Pu\,,(\ad_{\chi_{\ell}}^{2}P)u\rangle\,,\\
&&
\hspace{-1cm}
\label{eq.partunitRe}
 P-\sum_{\ell\in \mathcal{L}}
\chi_{\ell} P_{\ell}\chi_{\ell}
=-\frac{1}{2}\sum_{\ell\in
  \mathcal{L}}\ad_{\chi_{\ell}}^{2}P \,.
\end{eqnarray}
For \eqref{eq.partunitRe} write
$$
P-\sum_{\ell\in \mathcal{L}}\chi_{\ell}P\chi_{\ell}=\sum_{\ell\in
  \mathcal{L}}\left[P,\chi_{\ell}\right]\chi_{\ell}
=\sum_{\ell\in \mathcal{L}}\chi_{\ell}\left[\chi_{\ell},P\right]
$$
and take half the sum of the two last expressions with
$
\left[P,\chi_{\ell}\right]\chi_{\ell}+\chi_{\ell}[\chi_{\ell},P]=-\ad_{\chi_{\ell}}^{2}P$\,.\\
For \eqref{eq.partunit2} introduce the formal adjoints $P^{*}$ of $P$ (resp
$P_{\ell}^{*}$ of $P_{\ell}$) and write
$$
\|Pu\|^{2}-\sum_{\ell\in
  \mathcal{L}}\|P_{\ell}\chi_{\ell}u\|^{2}
=\sum_{\ell\in\mathcal{L}}\langle u\,,\, [P^{*}\chi_{\ell}^{2}P-\chi_{\ell}P^{*}P\chi_{\ell}]u\rangle\,.
$$
For one given $\ell\in \mathcal{L}$ compute
\begin{eqnarray*}
  P^{*}\chi_{\ell}^{2}P-\chi_{\ell}P^{*}P\chi_{\ell}
&=&P^{*}\chi_{\ell}\left[\chi_{\ell},P\right]+\left[P^{*},\chi_{\ell}\right]P\chi_{\ell}\\
&=&
P^{*}\chi_{\ell}\left[\chi_{\ell},P\right]+\left[P^{*},\chi_{\ell}\right]\chi_{\ell}P
+\left[P^{*},\chi_{\ell}\right]\left[P,\chi_{\ell}\right]\,.
\end{eqnarray*}
For the second and third term, use
$$
\left[P^{*},\chi_{\ell}\right]=[\chi_{\ell},P]^{*}=-[P,\chi_{\ell}]^{*}\,.
$$
For the first term, use
$$
2\chi_{\ell}[\chi_{\ell},P]=\left[\chi_{\ell}^{2},P\right]-\left[\left[\chi_{\ell},P\right],\chi_{\ell}\right]
=\left[\chi_{\ell}^{2},P\right]+\ad_{\chi_{\ell}}^{2}P\,.
$$
This gives
$$
\|\chi_{\ell}P\|^{2}-\|P\chi_{\ell}u\|^{2}=\Real\langle
Pu\,,\, \left[\chi_{\ell}^{2},P\right]u\rangle+\Real\langle Pu\,,\,
(\ad_{\chi_{\ell}}^{2}P)u\rangle
-\|(\ad_{\chi_{\ell}}P)u\|_{L^{2}}^{2}\,,
$$
and summing over $\ell\in \mathcal{L}$ proves \eqref{eq.partunit2}.
We conclude with a duality and interpolation argument.
\begin{lemma}
\label{le.duainterp}
Let $E_{1}\subset E_{0}\subset E_{-1}$ and $F_{1}\subset F_{0}\subset
F_{-1}$ be two Hilbert triples (the inclusion are continuous and dense
embeddings, the ${}_{-1}$ space being the dual of the ${}_{+1}$
space). The interpolated Hilbert spaces are denoted by $E_{r}, F_{r}$
for any $r\in [-1,1]$\,.
 Let $J:E_{0}\to F_{0}$ be an isometry, $J^{*0}J=\Id_{E_{0}}$
(non necessarily surjective), such that
$$
\exists C>0\,,
\forall u\in E_{1}\,,\quad C_{1}^{-1}\|u\|_{E_{1}}\leq \|Ju\|_{F_{1}}\leq C_{1}\|u\|_{E_{1}}\,.
$$
Then for any $r\in [-1,1]$\,,  $J$ defines a bounded operator from
$E_{r}$ to $F_{r}$ and the equivalence
$$
\forall u\in E_{r}\,,\quad C_{r}^{-1}\|u\|_{E_{r}}\leq \|Ju\|_{F_{r}}\leq C_{r}\|u\|_{E_{r}}\,,
$$
holds with $C_{r}=C_{1}^{|r|}$\,.
\end{lemma}
\begin{proof}
The equivalence of $\|Ju\|_{F_{1}}$ and $\|u\|_{E_{1}}$ implies that
$J\big|_{E_{1}}$ has a closed range $G_{1}$ in $F_{1}$ with the
inverse $J^{*0}\big|_{G_{1}}$\,. Let
$\Pi_{G_{1}}$ a continuous projection and note
$$
|\langle Jv\,,\, u\rangle|=|\langle v\,,\,
J^{*0}(\Pi_{G_{1}}u)\rangle|\leq C\|v\|_{E_{-1}}\|u\|_{F_{1}}\,.
$$
Therefore $J$ extends as a bounded operator from $E_{-1}\to F_{-1}$\,.
By interpolation $J\in \mathcal{L}(E_{r},F_{r})$ and $J^{*0}\in
\mathcal{L}(F_{r},E_{r})$ with $J^{*0}J=\Id_{E_{r}}$\,.
\end{proof}
\noindent\textbf{Application:} If $(\chi_{\ell})_{\ell\in \mathcal{L}}$ is a
locally finite partition of unity such that $\sum_{\ell\in
  \mathcal{L}}\chi_{\ell}^{2}=1$ and $(A,D(A))$ is positive
 self-adjoint operator  in $L^{2}(M)$ such that
$$
\left(\frac{\sum_{\ell\in
      \mathcal{L}}\|A\chi_{\ell}u\|^{2}}{\|Au\|^{2}}\right)^{\pm
1}\leq C\,,
$$
then  the equivalence
$$
\left(\frac{\sum_{\ell\in
      \mathcal{L}}\|A^{r}\chi_{\ell}u\|^{2}}{\|A^{r}u\|^{2}}\right)^{\pm
1}\leq C_{r}
$$
holds for any $r\in [-1,1]$\,. Simply apply the previous lemma with
\begin{eqnarray*}
  && E_{0}=L^{2}(M)\quad,\quad E_{1}=D(A)\quad,\quad E_{-1}=D(A)'\,,\\
&& F_{0,\pm 1}=\bigoplus_{\ell\in \mathcal{L}}E_{0,\pm 1}\,,\\
\text{and}&&
Ju=(\chi_{\ell}u)_{\ell\in \mathcal{L}}\,.
\end{eqnarray*}

\vspace{2cm}
\noindent\textbf{Acknowledgments:} This work was developed mainly while the
author had a ``D{\'e}l{\'e}gation INRIA'' at CERMICS in Ecole Nationale des
Ponts et Chauss{\'e}es. The author acknowledges the support of INRIA and thanks
the people of CERMICS for their hospitality. This work is issued and
has benefited from various, sometimes short, discussions with
J.M.~Bismut, C.~G{\'e}rard, F.~H{\'e}rau, T.~Hmidi, B.~Lapeyre, G.~Lebeau
T.~Leli{\`e}vre, D.~Le~Peutrec, L.~Michel, M.~Rousset,  H.~Stephan, G.~Stoltz.
The author takes the opportunity to thank all of them.


\begin{thebibliography}{}

\end{thebibliography}


\begin{thebibliography}{99}
\bibitem[AbMa]{AbMa} R.~Abraham, J.E.~Marsden.
\newblock \textit{Foundation of Mechanics, Second edition revised and
  enlarged.}
\newblock Benjamin/Cummings Publishing Co, Advanced Book Program,
(1978).

\bibitem[AnMe]{AnMe} K.G.~Anderson, R.B.~Melrose.
\newblock The propagation of singularities along gliding rays.
\newblock Inventiones Mathematicae \textbf{41}, pp.~197-232 (1977).

\bibitem[Bar]{Bar} C.~Bardos.
\newblock Probl{\`e}mes aux limites pour les {\'e}quations aux d{\'e}riv{\'e}es
partielles du premier ordre {\`a} coefficients r{\'e}els; th{\'e}or{\`e}mes
d'approximation; application {\`a} l'{\'e}quation de transport.
\newblock Annales Scientifiques de l'E.N.S 4e s{\'e}rie \textbf{3}-2,
pp.~185-233 (1970).

\bibitem[BeLo]{BeLo} J.~Bergh, J.~L{\"o}fstr{\"o}m.
\newblock \textit{Interpolation spaces. An introduction.}
\newblock Grundlehren der Mathematischen Wissenschaften,
Springer-Verlag (1976).

\bibitem[Ber]{Ber} J.~Bertoin.
\newblock Reflecting a Langevin process at an absorbing boundary.
\newblock  The Annals of Probability \textbf{35}-6, pp.~2021--2037
(2007).

\bibitem[BisLNM]{BisLNM} J.M.~Bismut.
\newblock \textit{M{\'e}canique al{\'e}atoire.}
\newblock Lect. Notes in Mathematics \textbf{866}, Springer-Verlag
(1981).

\bibitem[{BiLe}]{BiLe} J.M.~Bismut, G.~Lebeau.
\newblock{\sl The hypoelliptic Laplacian and Ray-Singer metrics.}
\newblock Annals of Mathematics Studies, 167, Princeton University Press (2008).

\bibitem[Bis1]{BisSurv} J.M.~Bismut.
\newblock A survey of the hypoelliptic Laplacian.
\newblock G{\'e}om{\'e}trie diff{\'e}rentielle, physique math{\'e}matique,
math{\'e}matiques et soci{\'e}t{\'e}~II. Ast{\'e}risque \textbf{322}, pp.~39--69 (2008). 

\bibitem[{Bis2}]{Bis2} J.M.~Bismut.
\newblock Hypoelliptic Laplacian and Bott-Chern cohomology.
\newblock Comptes Rendus Mathematique, \textbf{349}--2, pp.~75--80 (2011).

\bibitem[Bis05]{Bis05} J.M.~Bismut.
\newblock The hypoelliptic Laplacian on the cotangent bundle.
\newblock J.~Amer. Math. Soc., \textbf{18}-2, pp. 379--476 (2005).

\bibitem[Bisorb]{Bisorb} J.M.~Bismut
\newblock Hypoelliptic Laplacian and orbital integrals.
\newblock Annals of Mathematics Studies \textbf{177}, Princeton
University Press (2011).

\bibitem[BoLe]{BoLe} J.M.~Bony, N.~Lerner.
\newblock Quantification asymptotique et microlocalisation d'ordre
sup{\'e}rieur I.
\newblock Ann. Scient. Ec. Norm. Sup. $4^{e}$ s{\'e}rie \textbf{22},
pp.~377-433 (1989).
  
\bibitem[BdM]{BdM} L.~Boutet~de~Monvel.
\newblock Boundary problems for pseudodifferential operators.
\newblock Acta Math., \textbf{126}--1-2,  pp.~11--51 (1971).

\bibitem[BoJa]{BoJa} M.~Bossy, J.F.~Jabir.
\newblock On confined McKean Langevin processes satisfying the mean
no-permeability boundary condition.
\newblock Stochastic Processes and their Applications,
\textbf{121}-12, pp.2751--2775 (2011).

\bibitem[{Bre83}]{Bre} H.~Brezis,
\newblock {\it Analyse Fonctionnelle. Th{\'e}orie et applications.}
\newblock Masson, (1983).

\bibitem[Car]{Car} J.~Carrillo.
\newblock{Global weak solutions for the initial boundary value
  problems to the Vlasov-Poisson-Fokker-Planck system.}
\newblock Math. Meth. Appl. Sci. \textbf{21}, pp.~907--938 (1998).

 \bibitem[{ChLi}]{ChLi} K.~C.~Chang, J.~Liu.
 \newblock A cohomology complex for manifolds with boundary.
 \newblock Topological Methods in Non Linear Analysis,
 Vol.~5, pp.~325--340 (1995).

\bibitem[ChPi]{ChPi} J.~Chazarain, A.~Piriou.
\newblock{\textit{Introduction {\`a} la th{\'e}orie des {\'e}quations aux d{\'e}riv{\'e}es
  partielles lin{\'e}aires.}}
\newblock Gauthiers-Villars, (1981).

\bibitem[{CFKS}]{CFKS}
H.L.~Cycon, R.G.~Froese, W.~Kirsch, and B.~Simon.
\newblock{\it Schr{\"o}dinger operators with application to quantum
  mechanics and global geometry}.
\newblock Text and Monographs in Physics. Springer Verlag (1987).

\bibitem[DeMG]{DeMG} P.~Degond, S.~Mas-Gallic.
\newblock Existence of solutions and diffusion approximation for a
model Fokker-Planck equation.
\newblock Transp. Theory and Stat. Phys., \textbf{16} pp.~589-636
(1987).

\bibitem[EcHa]{EcHa} J.P.~Eckmann, M.~Hairer.
\newblock Spectral properties of hypoelliptic operators.
\newblock Comm. Math. Phys. \textbf{223}-2, pp. 233-253 (2003).

\bibitem[EnNa]{EnNa} K.J.~Engel, R.~Nagel.
\newblock One-parameter semigroups for linear evolution equations.
\newblock Graduate Texts in Mathematics 194, Springer-Verlag (2000).

\bibitem[Hel1]{Helgl} B.~Helffer.
\newblock Th{\'e}orie spectrale pour des op{\'e}rateurs globalement
elliptiques.
\newblock Soc. Math. de France, Ast{\'e}risque \textbf{112} (1984).
 
\bibitem[{HKN}]{HKN} B.~Helffer, M.~Klein, F.~Nier.
 \newblock Quantitative analysis of metastability in reversible
 diffusion processes via  a Witten complex approach.
 \newblock Matematica Contemporanea, 26, pp.~41--85 (2004).

\bibitem[{HeNi1}]{HeNi} B.~Helffer, F.~Nier.
 \newblock \textit{Quantitative analysis of metastability in 
reversible diffusion processes via a Witten complex approach: the case with boundary}.
 \newblock M{\'e}moire 105, Soci{\'e}t{\'e} Math{\'e}matique de France (2006).

 \bibitem[{HelNi}]{HelNi} B.~Helffer, F.~Nier.
\newblock \textit{Hypoelliptic estimates and spectral theory for Fokker-Planck operators and Witten Laplacians}.
Lecture Notes in Mathematics 1862, Springer-Verlag (2005).

\bibitem[HeNo]{HeNo} B.~Helffer, J.~Nourrigat.
\newblock \textit{Hypoellipticit{\'e} maximale pour des op{\'e}rateurs
  polyn{\^o}mes de champs de vecteur.}
\newblock Progress in Mathematics 58, Birkha{\"u}ser (1985).

 \bibitem[{HerNi}]{HerNi}
F.~H{\'e}rau, F.~Nier.
\newblock Isotropic hypoellipticity and trend to the equilibrium
  for the Fokker-Planck equation with high degree potential. 
\newblock Archive for Rational Mechanics and Analysis 171 (2), pp.~151--218
(2004).
 
\bibitem[{HHS2}]{HHS2}
F.~H{\'e}rau, M. Hitrik, J.~Sj{\"o}strand. 
\newblock Tunnel effect and symmetries for Kramers-Fokker-Planck type operators.
\newblock J. Inst. Math. Jussieu 10, no. 3, pp.~567--634 (2011).

\bibitem[HiPr09]{HiPr} M.~Hitrik, K.~Pravda-Starov.
\newblock Spectra and semigroup smoothing for non-elliptic quadratic operators.
\newblock  Mathematische Annalen, 344 no.4  (2009) pp.~801--846.

\bibitem[HormIII]{HormIII} L.~H{\"o}rmander.
\newblock\textit{The analysis of linear partial differential
  operators, III.}
\newblock Springer-Verlag (1985).

\bibitem[HormIII]{HormIV} L.~H{\"o}rmander.
\newblock\textit{The analysis of linear partial differential
  operators, IV.}
\newblock Springer-Verlag (1985).

\bibitem[Hor67]{Hor67} L.~H{\"o}rmander.
\newblock Hypoelliptic second order differential equations.
\newblock Acta Math., \textbf{119} pp. 147--171 (1967).

\bibitem[IkWa]{IkWa} N.~Ikeda, S.~Watanabe.
\newblock \textit{Stochastic differential equations and diffusion
  processes.}
\newblock 2nd Ed., North-Holland Mathematical Library \textbf{24}
(1989).

\bibitem[Kol]{Kol} A.N.~Kolmogorov.
\newblock Zuf{\"a}llige Bewegungen.
\newblock Ann. of Math. \textbf{2}-35, pp.~116-117 (1934).

\bibitem[{Lau}]{Lau} F.~Laudenbach.
 \newblock A Morse complex on manifolds with boundary.
 \newblock Geometrica Dedicata \textbf{153}-1, pp.~47--57 (2011).

\bibitem[Lap]{Lap} B.~Lapeyre.
\newblock Une application de la th{\'e}orie des excursions {\`a} une diffusion
r{\'e}fl{\'e}chie d{\'e}g{\'e}n{\'e}r{\'e}e.
\newblock Probab. Th. Rel. Fields \textbf{87} (1990) pp.~189--207.

\bibitem[Leb1]{Leb1} G.~Lebeau.
\newblock Geometric Fokker-Planck equations.
\newblock Port. Math. \textbf{62}-4, pp. 469--530 (2005).

\bibitem[Leb2]{Leb2} G.~Lebeau.
\newblock Equations de Fokker-Planck g{\'e}om{\'e}triques, II. Estimations
maximales.
\newblock Ann. Inst. Fourier \textbf{57}-4, pp.~1285--1314 (2007).

\bibitem[LeNi]{LeNi} T.~Leli{\`e}vre, F.~Nier.
\newblock Low temperature asymptotics for Quasi-Stationary
Distributions in a bounded domain.
\newblock arXiv~1309.3898

\bibitem[{Lep3}]{Lep3} D.~Le~Peutrec.
 \newblock Small eigenvalues of the Neumann realization of the semiclassical Witten Laplacian.
 \newblock Ann. Fac. Sci. Toulouse Math. (6),
 Vol.~19, no 3--4, pp.~735--809 (2010).
 
\bibitem[Ler]{Ler} N.~Lerner.
\newblock Energy methods via coherent states and advanced
pseudo-differential calculus,
\newblock in \textit{Multidimensional complex analysis and partial
  differential equations (Sao-Carlos 1995)?}
\newblock Contemp. Math. \textbf{205}, Amer. Math. Soc., pp~117--201
(1997).

\bibitem[LiSz]{LiSz} P.L.~Lions, A.S.~Sznitman.
\newblock Stochastic differential equations with reflecting boundary
conditions.
\newblock C.P.A.M. \textbf{37}-4, pp~511--537 (1984).

\bibitem[Lio]{Lio}J.L.~Lions.
\newblock{\textit{Equations diff{\'e}rentielles op{\'e}rationelles et
    probl{\`e}mes aux limites.}}
\newblock Springer (1961).

\bibitem[LNV]{LNV} D.~Le~Peutrec, F.~Nier, C.~Viterbo.
\newblock Precise Arrhenius law for p-forms: The Witten Laplacian and
Morse-Barannikov complex.
\newblock Ann. Henri Poincar{\'e}, \textbf{14} pp.~567--610 (2013).

\bibitem[LBLLP]{LbLLP} C.~Le~Bris, T.~ Leli{\`e}vre, M.~Luskin, D.~Perez.
\newblock A mathematical formalization of the parallel replica
algorithm.
\newblock Monte Carlo methods and Applications \textbf{18}-2,
pp.~119--146 (2012).

\bibitem[Luc]{Luc} B.~Lucquin-Desreux.
\newblock The Milne problem for the linear Fokker-Planck operator with
a force term.
\newblock Math. Meth. Appl. Sci., \textbf{26}, pp.~389--442 (2003).

\bibitem[MeSj1]{MeSj1} R.B.~Melrose, J.~Sj{\"o}strand.
\newblock Singularities of boundary value problems.~I.
\newblock CPAM \textbf{31}-5, pp.~593--617 (1978).

\bibitem[MeSj2]{MeSj2} R.B.~Melrose, J.~Sj{\"o}strand.
\newblock Singularities of boundary value problems.~II.
\newblock CPAM \textbf{35}-2, pp.~129--168 (1982).

\bibitem[Nie]{Nie} F.~Nier.
\newblock Remarques sur les algorithmes de d{\'e}composition de domaines.
\newblock S{\'e}minaire Equations aux D{\'e}riv{\'e}es Partielles,
Ec. Polytechnique \textbf{IX}-26 (1999).

\bibitem[{Nel}]{Nel} E.~Nelson.
\newblock {\it Dynamical Theories of Brownian Motion}.
\newblock 2nd Ed., Princeton University Press (2002).

\bibitem[Pra]{Pra} K.~Pravda-Starov.
\newblock Subelliptic quadratic differential operators.
\newblock Amer. J. Math. \textbf{133}-1, pp.~39--89 (2011).

\bibitem[{ReSi75}]{ReSi} M.~Reed, B.~Simon.
\newblock {\it Method of Modern Mathematical Physics}.
\newblock Academic press, (1975).

\bibitem[Ris]{Ris} H.~Risken.
\newblock\textit{The Fokker-Planck equation. Methods of solution and
  applications.}
\newblock 2nd Ed., Springer-Verlag (1989).

\bibitem[RoSt]{RoSt} L.P.~Rothschild, E.M.~Stein.
\newblock Hypoelliptic differential operators and nilpotent groups.
\newblock Acta Mathematica \textbf{137}, pp~248--315 (1977).

\bibitem[Sko1]{Sko1} A.V.~Skorohod.
\newblock Stochastic differential equations for diffusions in a
bounded region~1.
\newblock Theor. Veroyatnost. i Primenen \text{6}, pp.~264--274 (1961).

\bibitem[Sko2]{Sko2} A.V.~Skorohod.
\newblock Stochastic differential equations for diffusions in a
bounded region~2.
\newblock Theor. Veroyatnost. i Primenen \text{7},  pp.~3--23 (1962).

\bibitem[SjZw]{SjZw} J.~Sj{\"o}strand, M.~ Zworski.
\newblock Elementary linear algebra for advanced spectral problems.
\newblock Annales Inst. Fourier \textbf{57}-7, pp.~2095-2141 (2007).

\bibitem[StVa]{StVa} D.~Stroock, S.R.S~Varadhan.
\newblock Diffusion processes with boundary condtions.
\newblock Comm. Pure Appl. Math. \textbf{24}  pp.~147--225 (1971).

\bibitem[{Tay}]{Tay} M.E.~Taylor.
\newblock{\it Partial Differential Equations~1, Basic Theory}.
\newblock Applied Mathematical Sciences 115, Springer-Verlag (1997).

\bibitem[Tay1]{Tay1} M.E.~Taylor.
\newblock Grazing rays and reflection of singularities of solutions to
wave equations.
\newblock CPAM, \textbf{29}-5, pp.~1--38 (1976).
 
\bibitem[Tay2]{Tay2} M.E.~Taylor.
\newblock Grazing rays and reflection of singularities of solutions to
wave equations. Part~II (Systems).
\newblock \newblock CPAM, \textbf{29}-5, pp.~463--481 (1976).

\bibitem[Vil]{Vil} C.~Villani.
\newblock Hypocoercive diffusion operators.
\newblock International Congress of Mathematicians, Eur. Math. Soc.,
\textbf{3} pp.~473--498 (2006).

\bibitem[{Wit}]{Wit} E.~Witten.
\newblock Supersymmetry and Morse inequalities.
\newblock J. Diff. Geom. 17, no. 4, pp.~661--692 (1982).

\bibitem[{Zha}]{Zha} W.~Zhang.
\newblock \textit{Lectures on Chern-Weil theory and Witten deformations.}
\newblock Nankai Tracts in Mathematics, 4, World Scientific Publishing Co. (2001).
\end{thebibliography}
\end{document}